# IL FONDO **GIAN ANTONIO MAGGI**
# (1885-1937)
# CONSERVATO PRESSO
# LA BIBLIOTECA "G. RICCI"
# DEL DIPARTIMENTO DI MATEMATICA
# DELL'UNIVERSITÀ DEGLI STUDI DI MILANO

**a cura di Paola Testi Saltini**

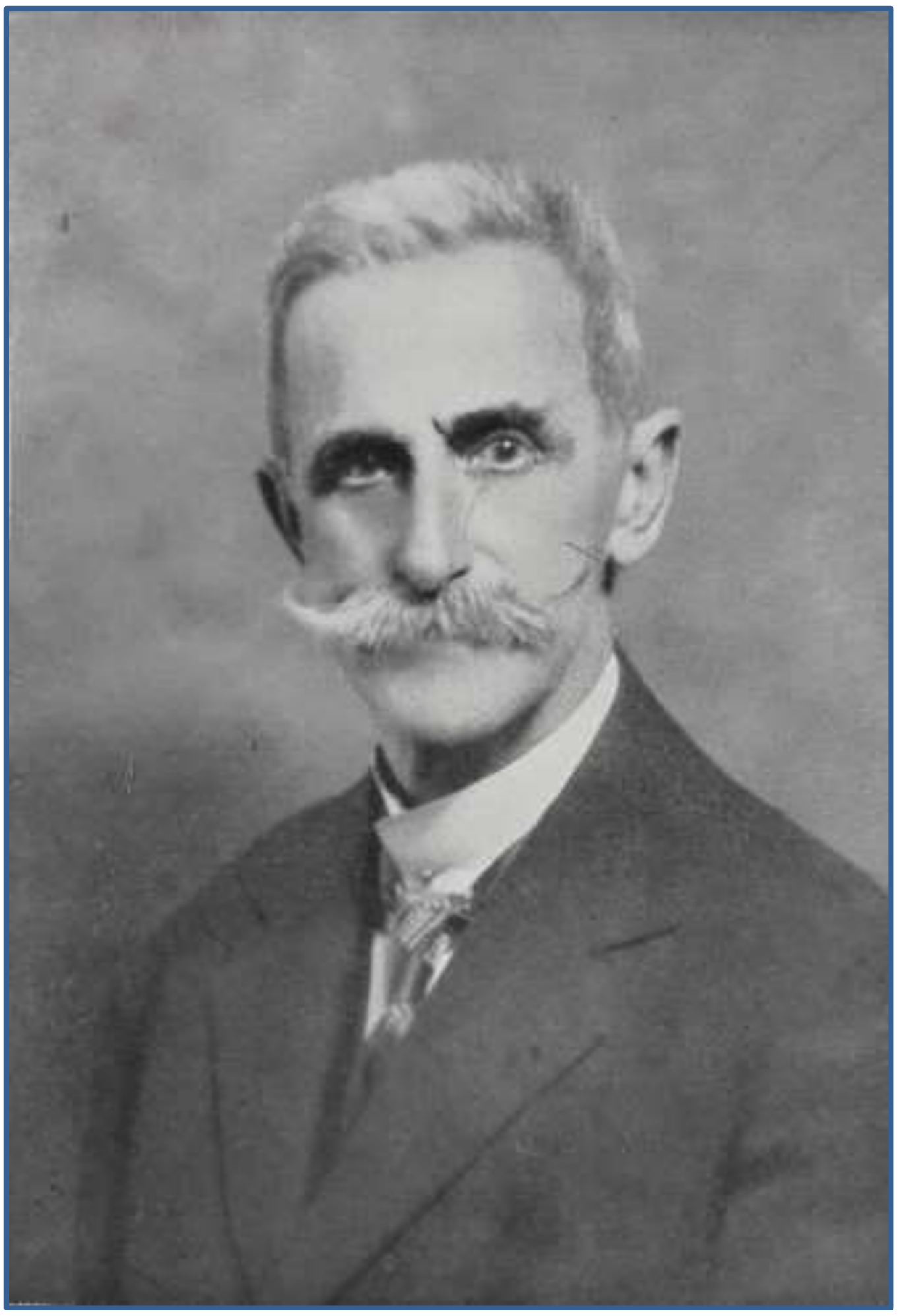

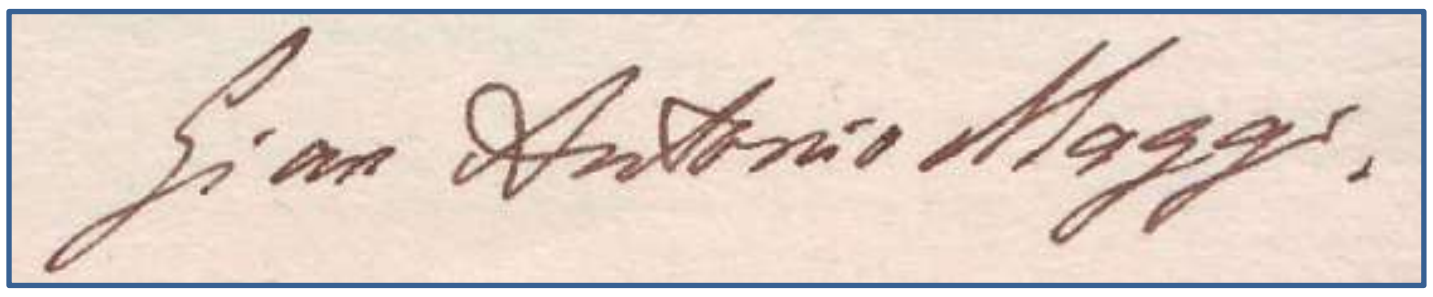

## Indice

All'incirca nel 1990 il "Fondo Maggi" fu reperito in un sopralzo della Biblioteca del Dipartimento di Matematica dell'Università degli Studi di Milano. Conteneva lettere scambiate tra il Professor Gian Antonio Maggi, vissuto dal 1856 al 1937, e noti scienziati italiani - quali Tullio Levi-Civita, Carlo Somigliana, Roberto Marcolongo, Salvatore Pincherle, Vito Volterra - nonché vari altri tipi di manoscritti.

Poiché Gian Antonio Maggi ha rappresentato nella cultura scientifica italiana degli anni 1890-1937 un punto d'incontro significativo per molti scienziati del Paese, il suo archivio presenta interesse non solo dal punto di vista istituzionale e del costume, ma anche strettamente dal punto di vista del dibattito scientifico.

Questa pubblicazione, che conclude il lavoro da me iniziato nel 1993 con la tesi *Alcuni problemi di Teoria della Relatività nella corrispondenza di Gian Antonio Maggi*, discussa presso il Dipartimento di Matematica, intende mettere a disposizione degli studiosi documenti inediti connessi alla storia del pensiero scientifico italiano. Ciò sarà utile non solo per impedire la dispersione o la perdita di documenti rilevanti, ma anche per evitare che la difficoltà ad accedervi favorisca giudizi non adeguatamente argomentati.

Ho scelto di rendere pubblico questo lavoro sul sito dedicato a Luigi Cremona - anche se egli non compare tra i corrispondenti di Maggi - perché la sua formazione scientifica ha origine nell'ambiente plasmato dall'attività della generazione di Cremona e ne rappresenta la continuazione nel panorama nazionale. Maggi, infatti, pur non essendo direttamente uno degli allievi di Cremona, è loro coetaneo ed è anche in ovvio contatto con i colleghi dei suoi maestri che, molto probabilmente, frequentava anche al di fuori dell'ambiente accademico. Per esempio, Cremona conosce Maggi, come si legge dalla lettera a lui indirizzata da Felice Casorati, n. 050-10769 del 7/7/1883: "Mia figlia Eugenia è fidanzata con il dott. G.A. Maggi, che mi assicura non soltanto di esserti noto, ma di aver anche ricevuto qualche segno della tua simpatia".

Paola Testi Saltini

## Sulla figura di Gian Antonio Maggi

L'intensa vita di Gian Antonio Maggi comincia a Milano il 19 febbraio 1856, in una nobile famiglia di origine piacentina: il padre, Giuseppe Pietro, è un noto orientalista, membro dell'Istituto Lombardo, la madre è Clara Anelli. Egli cresce in un ambiente di scienza, di letteratura, di studi classici per i quali eredita la passione, tanto che, al termine degli studi medi, esita prima di iscriversi alla Facoltà di Scienze dell'Università di Pavia. Si laurea il 7 dicembre 1877 in Fisica (e un mese più tardi in Matematica, l'8 gennaio 1878), dopo aver seguito i corsi di insigni maestri quali Eugenio Beltrami, Giovanni Schiaparelli, Giovanni Cantoni,[1] Felice Casorati; di quest'ultimo sposa la figlia Eugenia, il 23 settembre 1883.[2]

La sua attività nel campo scientifico, dopo qualche lavoro iniziale di fisica sperimentale, sotto la guida di Beltrami e Casorati si indirizza definitivamente agli studi matematici.

La sua prima sede di insegnamento è la stessa Università di Pavia, dove è inizialmente assistente alla cattedra di Fisica sperimentale, poi incaricato del corso di Fisica per i farmacisti, libero docente in Fisica matematica dal 1882 al 1885 e quindi professore aggiunto alla Scuola di Magistero. Nel 1882 compie un anno di perfezionamento a Berlino e frequenta i corsi di Hermann von Helmoltz, Gustav Kirchhoff e Karl Weierstrass.[3]

Nel 1885 viene nominato professore straordinario di Analisi Infinitesimale presso l'Università di Modena. L'anno successivo si trasferisce a Messina mantenendo, come professore ordinario, l'insegnamento di Analisi Infinitesimale fino al 1895. Nel 1891 viene eletto Rettore, fatto che, data la sua giovane età, testimonia la stima di cui gode presso i colleghi.

Nell'anno accademico 1895-96 viene chiamato a coprire la cattedra di ordinario di Meccanica razionale e incaricato di Fisica matematica all'Università di Pisa, succedendo a

---

[1] Il 13/3/1881 Beltrami così scriveva a Enrico Betti: "Caro Betti. Ti mando il manoscritto di una Memoria presentata all'Accademia dei Lincei dal D.r Maggi, assistente qui al Cantoni, e che fu anche mio scolaro nei primi anni che fui qui. Ti unisco anche uno schema di relazione che tu del resto puoi sostituire liberamente con altra, o modificare come credi. Mi è sembrato inutile di inserirci una più minuta indicazione delle varie parti del lavoro, perché l'Autore la dà egli stesso nelle prime pagine. Il lavoro in origine era molto più esteso, ed io, rivedendolo coll'Autore, l'ho persuaso a restringerlo più che poteva, il che ha fatto, assicurandomi che non si poteva far di più, stante la natura della questione.
Io desidererei che il lavoro fosse stampato, specialmente perché ritengo che il Maggi dovrà quanto prima domandare la libera docenza per mettersi in regola coll'insegnamento di fisica che dà, già da parecchi anni, ai farmacisti. E siccome egli è d'altronde un signore (ed anche un nobile) che fa tutto questo per passione allo studio, così mi rincrescerebbe che qualche circostanza potesse scoraggiarlo o svogliarlo, mentre credo che abbia attitudine e capacità di far bene, tanto più che possiede anche molta abilità nella parte sperimentale, ed ha quindi tutti i requisiti che possono servire a diventare un buon fisico. [...]" (L. Giacardi, R. Tazzioli (a cura di), *Le lettere di Eugenio Beltrami a Betti, Tardy e Gherardi. Pel lustro della Scienza italiana e pel progresso dell'alto insegnamento*, Collana Materiali per costruzione delle biografie di matematici italiani dopo l'Unità, Mimesis, Milano, 2012, lettera n. 66)

[2] I loro figli sono: Maria, Felice, Clara, Emma e Bianca.

[3] Il 24/3/1882, ancora per mano di Beltrami: "Caro Betti. Ieri, trovandomi a Milano, ho incontrato il D.r Maggi, reduce (per le vacanze di Pasqua) da Berlino, dove studia sotto Helmholtz e Kirchhoff. Parlandomi egli del suo desiderio, di trovare qualche occupazione pel venturo anno, ho pensato che egli potrebbe esser molto opportuno pel corso di Meccanica razionale a Pisa (naturalmente senza fargli cenno di questo mio pensiero). Come insegnante egli ha già tre o quattro anni di tirocinio, poiché insegnava la fisica ai farmacisti, e faceva un corso speciale sulla teoria degli errori di osservazione agli allievi della scuola di magistero.
Veramente la sua intenzione sarebbe di insegnare, potendo, fisica matematica, ed ha già chiesta la libera docenza qui a Pavia, per tale materia. Ma io credo che sarebbe assai meglio, per lui, di restare, almeno per alcuni anni, nella meccanica razionale, sia per approfondirsi ben bene in questa dottrina fondamentale, sia per allargare le sue cognizioni prima di assumere un insegnamento elevato, che deve più o meno mutare d'anno in anno.
Ho creduto bene di scriverti questo, perché tu mi dica quel che ne pensi." (*ibidem*, lettera n. 74)

Vito Volterra. Rimane a Pisa per quasi trent'anni, collega di Ulisse Dini, Luigi Bianchi, Eugenio Bertini e Paolo Pizzetti.

Nel 1924 nasce l'Università di Milano e Maggi, ormai quasi settantenne, accoglie con entusiasmo l'invito a tornare nella sua città natale come primo Preside della Facoltà di Scienze. Svolge tale incarico per un quadriennio ed è prima professore di Fisica matematica, poi di Meccanica razionale, inizialmente con l'incarico di Meccanica celeste e in seguito di Fisica matematica.

Nel 1926, il 26 giugno, entra a far parte - a Roma - della commissione giudicatrice del concorso per la prima cattedra di fisica teorica in Italia, commissione eletta, con voto segreto, dalle facoltà di Scienze di tutte le università statali. Con lui sono Michele Cantone di Napoli, Quintino Majorana di Bologna, Vito Volterra di Roma, Carlo Somigliana di Torino, Orso Maria Corbino di Roma e Antonio Garbasso di Firenze.

Nel 1927, con il concorso dell'Università, del Politecnico e dell'Osservatorio astronomico, fonda il Seminario matematico e fisico di Milano del quale rimane Direttore fino alla morte.

È Socio nazionale della Reale Accademia dei Lincei, uno dei XL della Società italiana delle Scienze e membro di numerose accademie nazionali e straniere (Membro effettivo del Reale Istituto Lombardo, corrispondente della Reale Accademia delle Scienze di Torino, dell'Accademia Gioenia di Scienze naturali di Catania, della Reale Accademia Peloritana, della Società matematica di Kharkoff, Socio ordinario della Società matematica di Kasan).

È collocato a riposo nel 1931 per raggiunti limiti di età, ma ciò non pone fine alla sua attività scientifica. Così Giulio Vivanti ne tratteggia la figura il 10 marzo 1932 nell'Aula magna dell'Università di Milano,[4] poco dopo il suo collocamento a riposo, durante una cerimonia voluta da colleghi, discepoli, amici, ammiratori.

> "La tua esperienza della vita universitaria, la tua equanimità, la tua obbiettività serena hanno segnato ai colleghi la via giusta... Ma la tua figura non sarebbe completa se io tacessi che, anche al di là della Scienza a cui hai dedicato il meglio delle tue forze, nulla vi è di umano o di divino che sia a te estraneo: profondo nelle discipline storiche e filosofiche, ammiratore dei grandi classici dell'antichità, conoscitore di un numero infinito di lingue... È ben raro che si ricorra al tuo aiuto senza averne i lumi richiesti... Concedimi anche di dire della saldezza del tuo carattere e della profondità della tua fede, che ti hanno permesso di risollevare il tuo spirito dopo i più gravi dolori e di ritornare con immutata vigoria alla tua vita di lavoro e di dovere".

Muore il 12 giugno 1937 e viene sepolto al Cimitero Monumentale di Milano, accanto alla moglie e al figlio Felice.

La produzione di Maggi si articola in due filoni: quello più prettamente scientifico, che comprende circa settanta Memorie sparse in periodici scientifici e Atti di Accademie e quello didattico, con sette trattati.

I primi articoli riguardano questioni particolari di meccanica, già trattate in generale dai classici (le piccole oscillazioni dei fili flessibili, il moto del pendolo conico e circolare,

[4] Citazione tratta da: A. Signorini, "Commemorazione del Socio Gian Antonio Maggi", *Rendic. R. Accademia dei Lincei*, v. XXVII, 1938, p. 457

l'equilibrio delle superfici flessibili ed inestendibili, la teoria dell'attrito, il moto relativo), alle quali apporta perfezionamenti ed approfondimenti.

Per il primo di questi problemi egli offre una trattazione da un punto di vista generale scrivendo le equazioni del moto, supponendo variabili la densità del filo, e integrandole in casi particolari.

A proposito del problema del pendolo stabilisce le equazioni della dinamica dei sistemi anolonomi, dette appunto "equazioni di Maggi", nel suo volume *Principii della teoria matematica del movimento dei corpi*, pubblicato nel 1896.[5] Questa sua scoperta precede di tre anni le ricerche di Paul Appel (trovo una precisazione in tal senso dello stesso Maggi nella lettera a Giovanni Giorgi, in data 12/3/1937, #79).

Nella teoria dell'equilibrio delle superfici flessibili e inestendibili egli riesce a dar ragione di un'indeterminazione delle forze applicate che si era presentata a Volterra. Alla teoria dell'attrito, infine, porta un duplice contributo: un'acuta critica della teoria di Paul Painlevé e la proposta di annoverare l'attrito tra le forze attive invece che tra le reazioni vincolari.

Maggi si occupa anche di elettrologia, con studi di elettrostatica e con un volume sulla teoria del campo elettromagnetico che Bruno Finzi nella Commemorazione scritta per il *Periodico di Matematiche*[6] descrive come:

> "una personale esposizione, eminentemente sintetica e informata ai concetti più moderni, che porta dal campo vettoriale alla relatività einsteiniana"

La teoria della propagazione ondosa è stata arricchita da Maggi di risultati fondamentali, che testimoniano anche il rigore analitico e matematico della rappresentazione dei fatti fisici con cui egli è solito condurre i suoi lavori, rigore legato certamente all'influenza di Kirchhoff e di Beltrami. Nella Memoria pubblicata nel 1887[7] è data la dimostrazione rigorosa della formula di Kirchhoff fondata su di una formula integrale generale che la pone al sicuro da ogni obiezione. In questo modo egli riesce a mettere in maggior evidenza l'analogia tra la formula in questione e la classica formula di Green. Inoltre, in una Nota lincea del 1920[8] introduce quelle onde che sono dette appunto "onde di Maggi".

Anche per quanto riguarda l'Ottica fisica si trovano pubblicazioni originali, concernenti l'applicazione della suddetta formula di Kirchhoff alla teoria della formazione delle ombre e della diffrazione.

Ancora collegata al lavoro del fisico tedesco è una ricerca relativa alla definizione del raggio luminoso nell'ottica ondulatoria: Maggi specifica come

> "il raggio di Kirchhoff si possa ricavare dall'espressione dell'energia elastica, caratterizzando il raggio come linea di propagazione dell'energia. Si viene così a stabilire una perfetta

[5] G.A. Maggi, *Principii della Teoria Matematica del Movimento dei Corpi. Corso di Meccanica Razionale,* Hoepli Milano, 1896.

[6] B. Finzi, "Necrologio: Gian Antonio Maggi", *Per. di Mat.*, 1937, p. 186.

[7] G.A. Maggi, "Sulla propagazione libera e perturbata delle onde luminose in un mezzo isotropo", *Annali di Matematica*, s. II, v .XVI, 1887, pp. 21-48.

[8] G.A. Maggi, "Sulla propagazione delle onde di forma qualsivoglia nei mezzi isotropi", *Rendic. R. Accademia dei Lincei*, s. V, v. XXIX, 1920, pp. 371-378.

> analogia, anche analitica, fra il raggio luminoso dell'ottica ondulatoria elastica ed il raggio luminoso dell'ottica elettromagnetica, definito dal vettore di Poynting."[9]

Maggi ha dedicato molte Note e Memorie alla teoria della funzione potenziale riguardanti sia le proprietà di questa funzione nell'immediata vicinanza della massa agente sia la funzione potenziale di doppio strato e di superficie. Il suo ultimo lavoro in questo campo riguarda la rappresentazione matematica della tensione elettrostatica.[10]

Più che nelle ricerche particolari, però, la completezza della sua preparazione appare nell'esposizione e discussione dei principi fondamentali della meccanica. Egli domina questioni anche delicate avendo studiato approfonditamente i lavori di W.K. Clifford (del quale traduce nel 1886 l'opera: *The common Sense of the Exact Sciences*), di Ernst Mach, di G. Kirchhoff, di Heirich Hertz e, negli ultimi tempi, di Albert Einstein.

Frutto di questi studi è la chiara e obiettiva esposizione della teoria delle relatività nella Nota: "Esposizione compendiosa dei principi sostanziali della Relatività", pubblicata sul *Nuovo Cimento*, nel 1921.[11] Maggi

> "non fu di questa teoria né fautore, né avversario deciso, ma cercò di esserne un critico sereno, e metterne in luce, con forma rigorosa, i concetti fondamentali".

Così Somigliana definisce il rapporto che lega Maggi alla teoria della relatività. Ed è proprio in una lettera della corrispondenza con Somigliana (#197) che si trova la seguente frase di Maggi:

> "... Nessun dubbio intanto che io sia disposto ad esaminare le cose da un punto di vista perfettamente obbiettivo. Io non presto alla teoria della relatività che un'adesione relativa. Poiché sta il fatto che quella teoria, attraente per ampiezza di quadro e acutezza di vedute, complica estremamente la trattazione razionale dei problemi fisici, per arrivare a perturbazioni dei risultati della teoria classica che, nella quasi totalità dei casi, sono inaccessibili alla verificazione sperimentale. Gli stessi risultati verificabili si fondano sopra esperienze che l'estrema delicatezza rende discutibili...".

È da sottolineare, comunque, che Maggi espone la Relatività nel suo Corso di Fisica Matematica, come si legge nella lettera a Giuseppe Erede, inviata da Pisa il 26 gennaio 1924 (#64).

Egli si rivela, dunque, un attento osservatore della realtà scientifica che lo circonda; non ripudia i principi classici e la filosofia che li supporta, ma si pone come continuatore del pensiero scientifico nella storia.

Questo suo rispetto non solo per una teoria attualmente in discussione come poteva essere quella della Meccanica classica, ma per una teoria addirittura già sorpassata si rivela nella Nota di Ottica Fisica del 1926,[12] ed ancora in quella del 1933,[13] in cui, accanto alla

---

[9] Cfr. C. Somigliana, "Gian Antonio Maggi", *Atti della R. Accademia delle Scienze di Torino*, 1937/38, v. 73, pp. 525-526.

[10] G.A. Maggi, "Di un'estrinsecazione energetica della tensione elettrostatica", *Rendic. R. Istituto Lombardo,* v. LVIII, 1925.

[11] v. XXI, pp. 5-33.

[12] G.A. Maggi, "Sul raggio di luce nell'ottica fisica", *Rendic. R. Istituto Lombardo,* v. LIX, 1926.

[13] G.A. Maggi, "La questione della superficie d'onda", *Rendic. del Seminario Matematico e Fisico di Milano,* v. VII, 1933.

Teoria elettromagnetica, trova spazio la Teoria elastica della luce, perché, come riporta Antonio Signorini:[14]

> "Egli diceva: A malgrado della fortunata concorrenza di altre teorie, mantiene tuttavia un intrinseco significato".

Una motivazione di questa sua posizione può essere anche ricercata nella passione per l'indagine storica che lo accompagna fin dalla sua giovinezza. Scrive ancora il Signorini nella sua commemorazione:[15]

> "... proprio nei primi anni della Sua carriera scientifica, iniziando la teoria generale delle piccole oscillazioni dei fili flessibili e inestendibili, Egli volle fare uno studio diretto di tutto quanto precedentemente era stato scritto pel caso particolare dei fili gravi: e così ebbe occasione di ricavare dai più antichi atti dell'Accademia di Pietroburgo una notizia ancor oggi poco diffusa: le due funzioni cilindriche di ordine zero, prima che da Bessel e da Fourier, sono state adoperate da Bernoulli e da Eulero."

È da ricordare anche l'opera innovativa svolta da Maggi in campo didattico: nei sei Trattati sulla Meccanica e in quello sulla teoria dei fenomeni elettromagnetici, egli cerca di far conoscere i procedimenti didattici che ritiene più opportuni e di cui effettivamente si serve. Dice ancora Signorini nella commemorazione già citata:[16]

> "I Suoi postulati non sono gli ordinari. Sono più simili a quelli proposti dal Mach e collimano coi concetti di Clifford: massa e 'forza motrice' di una figura materiale vengono definite, conferendo opportuni attributi all'accelerazione media della figura. Successivamente Egli passa al 'sistema' di figure materiali, di densità generalmente diverse; da questo, con un procedimento di limite, al più generale 'corpo naturale'
> Punto cardinale della Dinamica del Maggi è la netta, decisa distinzione della 'Dinamica fisica' dalla 'Dinamica dei sistemi'."

Un'altra innovazione nella presentazione della Meccanica riguarda, come si è già sottolineato, l'idea di includere l'attrito tra le forze attive (o, come le chiama lo stesso Maggi, *forze impresse*). In questo modo egli riesce a mantenere validi i principi di Lagrange dei lavori virtuali e il conseguente principio di D'Alembert anche se i vincoli del sistema sono scabri. Scrive a tal proposito Cisotti nella commemorazione pubblicata nel 1838:[17]

> "L'idea, come dissi, è veramente geniale, perché dà una più vasta portata a quei principi immortali della Meccanica, anche tenuto conto delle difficoltà che possono presentarsi per stabilire nuove equazioni a cagione delle introdotte forze d'attrito, che a priori, in generale, sono incognite."

Maggi stesso espone completamente i criteri con cui sono ispirate le sue opere sulla Meccanica nell'articolo "Réflexions sur l'exposition des principes de la mécanique rationnelle"[18] che Cisotti così traduce:[19]

---

[14] Alla p. 459 della Commemorazione citata in nota 1.
[15] *Ibidem*
[16] *Ibidem*, p. 461.
[17] U. Cisotti, "Gli scritti scientifici di Gian Antonio Maggi", *Rendic. del Sem. Mat. e Fis. di Milano*, 1938, p. 174.
[18] *Enseignement Mathématique*, 1901.

> "La Meccanica razionale è, come la Geometria, il risultato della idealizzazione di una scienza fisica: essa è una continuazione della Geometria".
> "Noi consideriamo delle figure finite, e non già dei punti materiali".
> "Considerare direttamente il movimento, per dedurne il fatto e le modalità della forza, in luogo di costruire a priori la forza, per dedurne le leggi del movimento...".

Se ci si domandasse come un uomo di così grande cultura e che ha dato tali contributi alla Scienza possa essere tanto poco conosciuto, se ne troverebbe, forse, una ragione nelle parole dello stesso Cisotti:[20]

> "... Dobbiamo a questo punto esprimere il rincrescimento che il grande amore del Maggi per la forma letteraria, di cui era indubbiamente signore, l'abbia indotto ad usare un linguaggio che non sempre giova alla migliore chiarezza dell'esposizione scientifica. Questa circostanza ha di certo nociuto alla diffusione del suo pensiero colto e innovatore nell'introduzione e svolgimento dei principii della Meccanica razionale".

Potrebbe apparire una motivazione ingenua, ma certamente la lettura e la comprensione sia delle lettere di Maggi che dei suoi scritti scientifici è resa difficile dallo stile spesso involuto, ricco di avverbi, aggettivi, citazioni dotte, nonché da una costruzione complessa della fraseggiatura. E ciò appare in misura molto maggiore nei suoi scritti che in quelli dei suoi, pur contemporanei, corrispondenti.

Quella stessa cultura letteraria che poco si addice alla divulgazione scientifica dà luogo ad una curiosità: nelle aule universitarie egli cita Aristotele ed Archimede nella lingua originaria. È questo che si apprende dal necrologio comparso sul periodico letterario *Atene e Roma*.[21]

Per più di cinquant'anni Maggi svolge la sua attività nell'ambito delle scienze fisiche e matematiche con la sola "distrazione" dello studio delle lingue: classiche (latino e greco) e moderne (tedesco, inglese, francese, spagnolo, russo, rumeno, ebraico, ungherese). Riferisce Somigliana:[22]

> "... E la sua conoscenza di queste lingue non era superficiale, ma ne aveva approfondita la struttura grammaticale, e poteva leggere, senza difficoltà e senza dizionario, opere letterarie scritte in polacco, in russo, in ungherese."

Maggi non trascura neppure i dialetti: in milanese, del quale è un purista, traduce dal greco una satira di Simonide di Amorgo contro le donne. Riesce anche a coniugare la scienza al suo amore per i classici, impegnandosi nella traduzione di alcuni passi di Platone attinenti la matematica (si vedano i carteggi con Paolo Ubaldi e Francesco Zambaldi).

Per dare chiusura al ritratto di Maggi possiamo usare le parole con le quali lui stesso descrive il proprio ruolo di educatore e di docente nel discorso tenuto il giorno delle

---

[19] Si veda la nota n. 14, p. 172.
[20] *Ibidem*
[21] A.M. Pizzagalli, "Necrologio: Gian Antonio Maggi", *Atene e Roma*, Firenze, 1937, pp. 214-215.
[22] Si veda la nota n. 6, p. 521.

onoranze tributategli quando lasciò l'insegnamento per raggiunti limiti d'età a Milano, nel 1932:[23]

> "Gli studi, coltivati con amore, con zelo, con fiducia, sono destinati a promuovere quella educazione dell'animo, i cui benèfici frutti si devono aspettare sui più molteplici campi, e, supremo frutto, se ne deve aspettare la formazione di una generazione atta a prestare opera veramente valida alla conservazione, al progresso, alla difesa dell'amato nostro Paese. Ch'io possa restare colla convinzione che la mia Scuola ha contribuito a questo insegnamento, e veramente mi consolerò del riposo, col dire fra me e me, senza bisogno che altri mi senta non omnis morior."

[23] Si veda la. nota n. 1, p. 462

## Descrizione del "Fondo Maggi"

Il "Fondo Maggi" è conservato presso la Biblioteca "Giovanni Ricci" del Dipartimento di Matematica dell'Università degli Studi di Milano e faceva parte del materiale donato alla Biblioteca da Carlo Somigliana. Si presentava diviso in tre faldoni dei quali viene fornita la descrizione qui di seguito.

Il contenuto del primo e del secondo faldone appariva in disordine e non rivelava alcun tipo evidente di raggruppamento: né per autore né per argomento né cronologico.
Il primo faldone recava la scritta *Giudizii 1885-1920* e raccoglieva la corrispondenza ricevuta da Maggi e le minute delle risposte da lui scritte nel periodo 1885-1920.
Il secondo faldone recava la scritta *Giudizii 2? serie* e raccoglieva, in modo analogo alla precedente, la corrispondenza dal 1921 fino al 1937, anno della morte di Maggi.
Oltre alle lettere - per la cui consistenza si rimanda all'elenco dei nomi citati - nei due contenitori si trovava il seguente materiale (il numero tra parentesi rimanda alla trascrizione: se cardinale si riferisce alla sezione delle lettere, se romano a quella dei Giudizii):

1. 30 Giugno 1899 - "Objection au §30 des "Études Geometriques *[sic!]* sur la théorie des parallèles" de Lobatschevsky (traduction de l'Hoüel)" - "Mandata a M.C.A. Laisant Rédacteur de l'Interm. des Mathém. Avenue Victor Hugo 169, Paris" **(93)**
2. Milano, 13 Febbraio 1930 - Bigliettino ad un amico (probabilmente T. Levi-Civita), sul problema etimologico dell'uso dei termini "sinistrorso" e "destrorso" **(227)**
3. Modena, 15 Dicembre 1885 - Osservazioni ad uno scritto (di autore ignoto), ricevuto dal Maggi in data 31 Ottobre 1885 e ad altre Note del medesimo autore, riguardanti la termodinamica **(229)**
4. Minuta di lettera a destinatario ignoto riguardanti la definizione di massa **(230)**
5. Minuta di lettera su "Massa e forza motrice" a destinatario ignoto **(231)**
6. Giudizio su una Dissertazione riguardante "L'esperienza nelle Teorie di Maxwell e di Lorentz o l'interpretazione meccanica dei fenomeni elettrici" di autore ignoto **(I)**
7. Giudizio sulla Nota di A. Signorini: "Sulla teoria analitica dei fenomeni luminosi nei mezzi cristallini uniassici" **(II)**
8. 29 Luglio 1913 - Traduzione dal tedesco dell'articolo "Il Principio di Relatività dell'Elettrodinamica" di Otto Berg. - "Per il Prof. A. Forte, al quale mandata il 29 Luglio 1913" **(III)**
9. Relazione del Concorso di Pavia, 5 giugno 1913 - Giudizi su Ermenegildo Daniele e Ernesto Laura **(IV)**
10. Giudizio su E. Daniele dopo il concorso di Pavia **(V)**
11. Fascicolo "Concorso di Meccanica Razionale per Politecnico di Torino 1915" con giudizi sui Proff. Ernesto Laura e Giuseppe Armellini **(VI)**
12. Analisi dell'opera scientifica del Dott. A. Signorini **(VII)**
13. Analisi dell'opera scientifica del Dott. Lucio Silla **(VIII)**
14. Giudizi globali sui Proff. Lucio Silla, Matteo Bottasso e Goffredo Mancini **(IX)**
15. Relazione favorevole alla pubblicazione nei *Rendiconti dell'Istituto Lombardo* della Nota dell'Ing. Arnaldo Masotti: "Decomposizione intrinseca del vortice e sue applicazioni" **(X)**
16. Relazione favorevole alla pubblicazione nei *Rendiconti dell'Istituto Lombardo* della Nota del Prof. Aldo Pontremoli: "Sulla conducibilità elettrica e termica nei metalli" **(XI)**

17. Milano, 15 Febbraio 1928 - Relazione favorevole alla pubblicazione nei *Rendiconti dell'Istituto Lombardo* della Nota di Arnaldo Belluigi: “Sull'impiego delle isogamma” (in duplice copia) **(XII)**
18. Milano, 15 Febbraio 1928 - Giudizio per il Concorso Cagnola 1928 della Nota: “L'azione dei campi elettrici intermolecolari sulla emissione delle righe spettrali” il cui autore non viene nominato (in duplice copia). **(XIII)**
19. Milano, 30 Maggio 1928 - Relazione favorevole alla pubblicazione nei *Rendiconti dell'Istituto Lombardo* della Nota del Prof. Bruno Finzi: “Integrazione per successive approssimazioni delle equazioni di un liquido viscoso in moto stazionario” **(XIV)**
20. Analisi dell'opera scientifica del Dott. A. Signorini **(XV)**
21. “Relazione della Commissione Giudicatrice della Libera Docenza in Fisica Teorica del Dr. Gleb Wataghin” **(XVI)**
22. Milano, 22 Aprile 1930 - Relazione favorevole alla pubblicazione nei *Rendiconti dell'Istituto Lombardo* della Nota del Signor Guido Facciotti: “Sulla corrispondenza fra spostamenti e deformazioni in un solido elastico” **(XVII)**
23. Relazione sulla Nota della Dott. Maria Pacifico “Sopra alcuni problemi al contorno per le funzioni armoniche” per la pubblicazione nei *Rendiconti dell'Istituto Lombardo* **(XVIII)**
24. Profilo di Bruno Finzi **(XIX)**
25. “Copia del M.S. consegnato a Chisini pel prossimo numero del Periodico il 19 Nov. 1933” - Puntualizzazione alla propria Nota: “In quanto tempo un pianeta, fermato, cadrebbe sul Sole?” [*Periodico di Matematiche,* 1927] **(XX)**
26. Osservazioni sulla Nota di G. Vitali “Sulla forza centrifuga” - “Copia di relazione consegnata a Chisini su lettura del M.S. presentato pel Periodico di Matematica” **(XXI)**
27. Elenco delle pubblicazioni di A. Signorini sull’Elastostatica dal 1930 al 1932 **(XXII)**
28. Milano, 27 Aprile 1936 - Appunti dal testo di J. Bertrand, *Leçons sur la théorie mathématique de l'electricité,* (Paris,1890), p. 267 **(XXIII)**
29. Milano, 27 Maggio 1936 - Profilo di Maria Pastori **(XXIV)**
30. Milano, 9 Giugno 1937 - Appunti sulla “Statistica di Fermi” **(XXV)**

Frammenti vari con minute di calcoli, di traduzioni, di appunti e di parti di opere poi stampate.

Nel terzo faldone si trovavano, invece, bozze di libri, minute varie, giudizi per Concorsi ed alcune lettere e cartoline. Il titolo dato da Maggi a questa terza parte è *Rusticationes* ed i vari scritti sono a loro volta suddivisi in “fascicoli”, individuati dal luogo e dalla data di stesura del contenuto. Ogni fascicolo contiene materiale scritto da Maggi durante le vacanze estive dal 1920 al 1936 con poche eccezioni.

*Bagni di Lucca - 1920.*
- 23 luglio: studio sulle “Curve piane parallele” (due pagine);
- minuta della Nota: “Sulla propagazione delle onde di forma qualsivoglia nei mezzi isotropi” (forse non terminata; dodici pagine, con ripetizione della 4 e della 12, mancano le pagine 6, 7, 8, 11);
- 3 settembre: minuta “Sulla congruenze di rette e le superficie parallele” (trentotto pagine, con inserimento della 18');

- 3 settembre: appunti sullo studio della propagazione delle onde con bibliografia di riferimento.

*Cireglio - 1922*

Note Matematiche - Cireglio - (Estate 1922). [vuoto]

*Gavinana - 1923*

- frammenti per la stesura di *Elementi di Statica e Teoria dei Vettori applicati*, Litografie, Pisa, 1923;
- 6 agosto: minuta della prefazione del volume di cui sopra (da a ad h, con ripetizione delle pagine a, b, c).

*Gavinana - 1924*

- frammenti vari: da pag. I a pag. VIII minuta per una Nota sull'uso delle equazioni di Hertz "(Adattamento della Parte II della Tesi di Laurea della sig.na A.M. Simonatti)", più alcuni stralci di appunti a carattere letterario.

*Note. - Bormio - 1925*

- due minute della Nota "Recenti vicende della Relatività", pubblicata poi sul *Periodico di Matematiche*, s. IV,v. VI, pp. 1-19 (da pag. 1 a pag. 36, manca la 30, l'una; venti pagine la seconda avente per titolo: "Appunti per l'eventuale Discorso al Congresso di Mathesis").

*Note matematiche - Cheglio - 1926*

- appunti relativi ad un articolo di Kottler apparso sugli *Annalen der Physik*, 1923 (quattro pagine);
- "Teoria della radiazione termica. Sunto conforme a Planck." (trentacinque pagine);
- "Entropia e probabilità - Planck, Teoria della radiazione termica" riguardante il precedente "sunto conforme" (due pagine).

*Note Matematiche - Lanzo d'Intelvi - 1927*

- appunti vari: "Somma della progressione aritmetica", "Sulla derivata covariante. cfr. Marcolongo - Relatività - pagg. 35,36", "Schrödinger - Wellenmechanik. pag.110" (da pag. I a pag. V);
- 18 settembre: minuta di una lettera a Somigliana (da pag. VI a pag. VIII; si veda la lettera #**203**);
- "Dimostrazione che $\frac{z^n}{n!}$ riceve il massimo valore per n=z";
- appunti su lavori di Schrödinger (trentatré pagine).

*Porto Valtravaglia - 1928 - Hic est locus Marini Falerii (1929).* [vuoto]

*Valnegra - 1929*

- luglio: "Effetto del movimento del sistema solare nell'esperimento di Michelson";
- 29 luglio: "Minimo della somma delle distanze di un punto dai vertici di un triangolo", (undici pagine);
- sul "Poligono articolato"
- 9 agosto: "Distanza del punto del piano equatoriale dal centro della Terra in cui si eguaglia la gravitazione e forza centrifuga";
- 10 agosto: "Equilibrio di un filo flessibile e inestendibile...";
- 11 agosto: "Correzione della lunghezza del pendolo...";
- sul "Poligono funicolare chiuso", (da I a VIII);
- 7/8 settembre: frammenti letterari;
- 9/11 settembre: "Problema" su un "fascio di circoli" propostogli "dall'ing. Livraga di Treviglio, in villeggiatura a Piazza Brembana" (da A a D);
- 13 settembre: "Reminiscenze del tema di maturità" (da α a η);

- 21 settembre: “La forza centrifuga nelle voltate”; “Problema idrostatico”; “Velocità del perno di una ruota pesante”.

*Valnegra - 1930*

- “Tema agli esami di maturità scientifica - Sessione estiva - 1930” (da I a VII);
- appunti da una Nota di Gugino (da a. a c.);
- 9 agosto: minuta di una lettera a Pincherle (si veda la lettera #**145**);
- 19 agosto: cartolina di Pincherle da Modena (si veda la lettera #**146**);
- 5 settembre: lettera di Levi-Civita da Padova (si veda la lettera #**103**);
- 19 settembre: lettera di Levi-Civita da Padova (si veda la lettera #**104**);
- 23 settembre: cartolina di Levi-Civita da Padova (si veda la lettera #**105**);
- commento a un lavoro di Pfeiffer (da α a ε);
- “Bilancio energetico” (da I a IV);
- “L'espressione $\int_\tau Vkd\tau$ nella gravitazione e nell'elettromagnetismo” (pagine a. e b.);
- una pagina di commento ad un articolo sul Circolo polare apparso su *Débats*;
- “Problema al Concorso e Abilitazione pei Lincei - Estate 1930” (cinque pagine);
- “Relatività” (due pagine);
- “Problema”;
- 2 settembre: minuta di una lettera al Direttore del *Journal des Débats*, pubblicata nel numero del 27 settembre, col titolo “Le latin en Italie” (si veda la lettera #**60**);
- minuta della commemorazione del Senatore Scherillo; telegramma di condoglianze alla vedova.

*Valnegra - 1931*

- 30 luglio: stampato, “Da tener presente. De deflectione lucis in proximitate solis (Gianfranceschi).”, in latino: comunicato della Pontificia Academia Scientiarum Novi Lyncæi;
- 14 agosto: “Dell'uomo, della carrucola e della secchia”;
- “Movimento di un punto materiale assoggettato a mantenersi sopra un cerchio...”;
- “Godofredo Garcia” (sei pagine);
- “Pel Concorso”, analisi dei lavori di: Finzi, Gugino, Nobile, Sbrana, Krall e Teofilato).

*Valnegra (Mojo de' Calvi) - 1932*

- “Problema” di geometria (da 1 a 10);
- appunti da una Nota di Chisini: “La dimostrazione cinematica di un problema di minimo”;
- appunti dal testo di Levi-Civita: “Caratteristiche ecc. e propagazione ondosa, Bologna, 1931” (tredici pagine);
- frammenti dello studio su “Funzione potenziale di una sfera omogenea in un punto interno”.

*Lanzo d'Intelvi - 1933*

- “Problema” di ottica (sette pagine);
- elenco di correzioni alle “Poesie milanesi di Carlo Porta a cura di Decio e Capello - 1932”.

*Lanzo d'Intelvi - 1934*

- risoluzione di un problema apparso sul *Periodico di Matematiche*, 1 luglio 1934, pag. 251;
- 11 agosto: appunti della Nota “Polvani, Velocità della luce, Annuario del Seminario di Milano, 1933”;
- 16 agosto: appunti della Nota “G. De Marchi. Omogeneità, similitudine e modelli idraulici, Annuario del Seminario di Milano, 1933, pag. 49”;
- appunti sul “Metodo delle dimensioni nulle”, dal testo Lezioni di Meccanica Razionale di Levi-Civita e Amaldi;

- 29: frammenti letterari;
- stralci di quotidiani svedesi con le relative traduzioni del Maggi sul problema della regolazione della natalità in Italia sotto il regime fascista;
- minuta di una cartolina inviata a Finzi (si veda la lettera #**76**).

*Lanzo d'Intelvi - 1935*
- appunti dell'articolo "La chaîne de fortune (Débats, 10 Aout, 1931)" (tre pagine);
- agosto: "Sulle misure elettromagnetiche" (quattro pagine);
- "Teorie statistiche" (da A ad H);
- appunti di Meccanica (tre pagine);
- 17 settembre: minuta di lettera al prof. Alessandro Visconti della R. Università di Ferrara (in milanese; si veda la lettera #**220**).

*Lanzo d'Intelvi - 1936*
- "Problema delle erbe scagliate" (da 1 a 8);
- "Problema della Domenica del Corriere del 2 Agosto: I guerrieri";
- "Al prof. Puccianti per la sua Nota in Rendic. dei Lincei, vol. XXIII, 1936 - Lettera ideale";
- "Paradosso dimensionale";
- 3 settembre: biglietto di Facciotti (si veda la lettera #**65**);
- 8 settembre: minuta di lettera a Facciotti con trascrizioni della "Formola del Facciotti" medesimo (si veda la lettera #**66**);
- 8 settembre: minuta di lettera a Signorini (si veda la lettera #**178**);
- 10 settembre: minuta di telegramma a Facciotti con correzione rispetto alla precedente (si veda la lettera #**68**);
- "Il problema di Archimede";
- saggio di traduzione dall'ungherese;
- frammento di lettera.

In questa pubblicazione si trovano le trascrizioni di tutto il materiale contenuto nei primi due faldoni e delle lettere contenute nei fascicoli delle *Rusticationes*.

## Trascrizione delle lettere ordinate per corrispondente

*Criteri di edizione*

Nonostante le lettere di Maggi siano soltanto delle minute e la grafia non sia sempre molto comprensibile, il fatto che le lettere siano conservate in discreto stato e le grafie dei corrispondenti siano abbastanza chiare ha ridotti il numero degli interventi davvero necessari. Tali interventi sono stati resi graficamente con il corsivo, quindi, tutto ciò che non è scritto in questo stile, fa parte della trascrizione. Solo in caso di evidente errore di scrittura si è inserito un *[sic!]*; invece, nel caso di errata scrittura di nomi, si è preferito inserire una nota a piè di pagina nella quale è indicato il nome corretto. Nel caso di rimandi scritti in fondo alla pagina o alla lettera, si è scelto di inserirne il testo di fila. In alcuni casi, sono state trascritte parti di testo cancellate da Maggi (con opportune segnalazioni), perché ritenute comunque interessanti.

Si è scelto di inserire la scansione dei disegni che compaiono nelle lettere e sono state inserite le scansioni anche di qualche stampato, conservato insieme al carteggio, al quale si fa riferimento nel testo delle lettere, in modo da facilitare la comprensione.

Nel Fondo sono presenti lettere (o parti di esse) scritte in diversi idiomi e caratteri: quelle in francese, in dialetto milanese, in latino e alcune parole in greco sono state trascritte, mentre quelle in russo (cirillico) e in tedesco sono state inserite come immagini.

Alcune minute sono presenti in più copie quasi simili tra loro; ho scelto di trascrivere solo quelle che presentano differenze in qualche modo significative tra le varie stesure: sono numerate come “#BIS”. Per quanto riguarda l’indicazione del luogo e della data, ho utilizzato le usuali abbreviazioni “s.d.” e “s.l.” quando mancanti.

**1**
Gian Antonio Maggi ad **Apreda**

Pisa 9 luglio 1920

Caro Professor Apreda,

La informo di aver ricevuto il Suo manoscritto, appena ne fui in possesso, e avrei voluto rimandarglielo più sollecitamente con le osservazioni da Lei desiderate. Ma, per esaminarlo, coll'attenzione che richiede questo genere di scritti, mi occorreva un certo tempo, del quale non ho potuto disporre prima d'ora.

Non potrei che lodare il Suo proposito di approfondire codeste questioni, e l'impegno ch'Ella ha posto nella svolgimento del Suo assunto. E perciò tanto più mi spiace di non poter aderire che in assai ristretta misura alle Sue vedute e alle Sue considerazioni.

In primo luogo, io non reputo opportuno fondare i principi generali della Meccanica sulla Statica. Certo la teoria dell'equilibrio comporta e merita una trattazione speciale, mentre poi si può dedurre, come caso singolare, dalla teoria generale del movimento. Ma, per la singolarità del caso, alcune circostanze, attinenti alla teoria generale, vi si trovano ridotte ad una forma particolare, altre ne riescono senz'altro occultate, in modo da potersi giudicare non egualmente efficace il tradizionale procedimento inverso. Lo stesso principio dell'eguaglianza dell'azione o della reazione s'identifica, nel caso dell'equilibrio, colla eguaglianza ed opposizione della reazione dei sostegni e della forza impressa al mobile: la reazione di un filo reggente un grave in equilibrio è eguale ed opposto al peso del grave. E simile circostanza non si verifica più nel caso del movimento, in cui la reazione dei sostegni, dipende, in generale, dalla forza impressa e dalla velocità.

In secondo luogo, non potrei convenire nell'assumere lo sforzo muscolare a fondamento del concetto di forza. Lo sforzo muscolare si può ben invocare come esempio di manifestazione palpabile della forza, in alcuni casi speciali. Ma esso non è suscettibile, per sé stesso, che di una misura più o meno grossolana, e nella generalità dei casi, non si saprebbe cosa significhi ricondurre ad esso la forza. Che cosa dice lo sforzo muscolare nella valutazione dello sforzo sostenuto dai fondamenti di una casa? Ella stesso è obbligato a rinunciare ad affermare sensibile l'eguaglianza dei due sforzi muscolari che introduce per rappresentare l'azione di un peso nel deformare una molla e la reazione della molla deformata. Ora, messa in disparte la dimostrazione dell'eguaglianza, non possiamo in modi più ovvii che ricorrendo alla rappresentazione per mezzo dello sforzo muscolare metter in evidenza l'esistenza dell'azione e della reazione? Per esempio, la molla sostiene il peso e il peso schiaccia la molla. Ma poi dalla suddetta rinuncia Ella è tratto a concludere che bisogna ammettere per convenzione l'eguaglianza dell'azione e della reazione. E davvero non saprei come ci si potrebbe appagare di una convenzione stabilita in codesti termini. Io intendo che il principio dell'eguaglianza dell'azione e della reazione vada stabilito come postulato: ma dichiarando che si reputa appoggiato ad esperienze, atte a dimostrarlo colla desiderabile esattezza.

Il Poincaré, nel brano da Lei citato, mi sembra fare più che altro la critica dell'eguaglianza di due forze fondata sul criterio dell'elisione di uno stesso peso. Vero è, che subito dopo, aderendo a quello ch'egli chiama la definizione di Kirchhoff della forza, afferma esplicitamente che questa definizione va completata col principio dell'eguaglianza dell'azione e della reazione, riguardato, non come una legge sperimentale, ma come una definizione. Ma poi conclude che, definito con questo il rapporto di due masse, tocca all'esperienza verificarne la costanza. E non è questo rimettere all'esperienza di verificare l'eguaglianza dell'azione e della reazione?

Per quanto riguarda l'opera di Galileo, bisogna tener presente che molto fu scritto intorno ad essa, per evitare di tornare su cose note. Non saprei quanto di nuovo contengano le sue osservazioni sul concetto di Galileo della propensione alla caduta, commisurata colla quantità di materia del grave. Per riguardo ad intendere il peso in Galileo "nel senso di quantum di materia associato all'elemento costante del moto, cioè all'accelerazione" o l'espressione rimane vaga, o, se dobbiamo prenderla col significato che Galileo considerasse il peso col prodotto della massa per l'accelerazione, mi sembra un'asserzione troppo ardita.

Infine trovo da fare varie objezioni alle Sue conclusioni sull'intervento, nella stessa opera di Galileo, del principio dell'eguaglianza dell'azione e della reazione.

Nel paragrafo "La priorità del principio etc." mi sembra che varii passi di Galileo, a cui Ella accenna, si riferiscano a due forze eguali e contrarie applicate allo stesso corpo, piuttosto che ad azione e reazione, che sono applicate a corpi diversi. Nel caso poi del filo compresso fra due diti, io, per verità, vedo poco più che il meccanismo della stessa compressione, che implica la simultanea azione dei due diti. Certo, i fatti che un corpo compresso comprime, che un corpo a cui è applicata una forza oppone una resistenza ad un ostacolo sono così ovvii che è facile trovarli implicati in varii ragionamenti. Ma, per discorrere di priorità, bisogna ben misurare il passo che occorre fare per enunciare il principio dell'eguaglianza dell'azione e della reazione come una delle tre leggi fondamentali del movimento.

Più oltre, quando Ella ricorda la parte di Galileo sul moto circolare uniforme, mi sembra che di nuovo si tratti piuttosto di opposte tendenze relative ad uno stesso corpo. Per quanto poi ai brani che si riferiscono al principio delle velocità virtuali, la questione mi pare assolutamente un'altra.

Com'Ella riconoscerà, i miei appunti sono, almeno in buona parte, questione di modo di vedere. Guardi Lei in quanto ne resta persuaso, e non pensi ch'io mi avrei menomamente a male che sentisse altri, e accogliesse un diverso parere, per pubblicare il Suo scritto. Creda mi spiace assai di non poter consigliarlo io stesso, tanto più che c'è del vero in quello che Ella dice sul giudizio dell'attività degli insegnanti.

Mi conservi sempre la Sua buona memoria che io Le ricambio di cuore, e mi tenga sempre disposto a giovarLe per quanto io posso. Intanto Le stringo cordialmente la mano, e mi confermo

Affmo Suo

**2**
Gian Antonio Maggi a **Moisè Ascoli**

Pisa 5 aprile 1918

Carissimo Moisè,

La controversia tra il nostro Guido Grassi e il Guglielmo[24] ha dato occasione di ripetuti discorsi sull'argomento tra me e il collega Puccianti, per cui avendo questo ricevuto la tua Nota sul Nuovo Cimento e trovata assai interessante si apprestò a mettermene a parte. D'accordo col Puccianti sull'attenzione che merita il tuo scritto trovo pure da fare qualche riflessione sul principio da cui prendo le mosse e penso che non ti dispiacerà che te ne discorro. Tanto più che un discorso fra me e te ci richiama ai bei tempi, oramai fatti troppo lontani, e le cui riprese sono diventate troppo scarse e troppo fugaci.

Dici che:

$$dq=C_p\frac{\partial T}{\partial v}+C_v\frac{\partial T}{\partial p}dp \qquad (1)$$

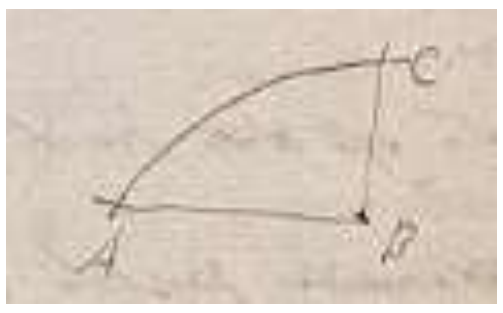

è in contraddizione col primo principio della termodinamica, perché esprime che un corpo riceve la stessa quantità di calore passando direttamente da A a C o passandovi per la via ABC.

Certo, col primo principio della termodinamica è in contraddizione l'ipotesi che il secondo membro di (1) sia un differenziale esatto: perché ipotesi equivalente a quella che la quantità di calore che riceve un corpo per passare da uno stato ad un altro sia primariamente dipendente dagli stati estremi. Ma colla precedente interpretazione della (1), si discorre di quantità di calore infinitesimali relative al passaggio da uno stato ad un altro infinitamente vicino. La relazione è stabilita tra infinitesimi del primo ordine. E poiché una relazione così fatta è sempre conciliabile coll'ommissione dell'aggiunta di infinitesimi di ordine superiore, esso stesso conciliabile con la prima legge della Termodinamica, allo stesso modo che l'espressione del lavoro infinitesimale Xdx+Ydy+Zdz si concilia colla mancanza di potenziale quantunque si presti all'interposizione che il lavoro della supposta forza è lo stesso del passaggio dal punto mobile direttamente da A a B e per la spezzata terminata agli stati precedenti avendo dati di grandezza dx, dy e dz, rispettivamente paralleli all'asse delle x, all'asse delle y e all'asse delle z.

Ammesso lo stesso primo principio della Termodinamica è pur certo che la (1) si deduce da

$$\Delta q_{AC}-\Delta q_{AB}-\Delta q_{BC}=\frac{1}{E}\text{Area ABC}$$

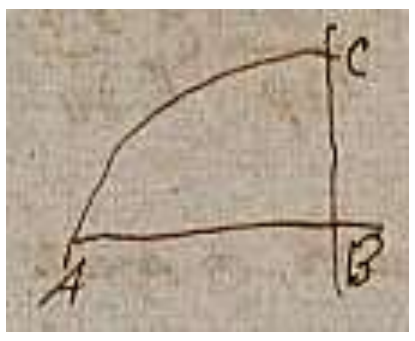

[24] Nel 1914, sui *Rendiconti della Reale Accademia dei Lincei*, il prof. G. Guglielmo pubblica l'articolo "Sull'esperienza di Clément e Desormes e sulla determinazione dell'equivalente meccanico della caloria" (pp. 698-703). Questo scritto, e quello del 1916, "Sulle leggi di Poisson e dello stato aeriforme in relazione al primo principio di termodinamica", innescano una controversia con Guido Grassi che pubblica negli anni successivi tre articoli a commento (si trovano tutti in rete).

trascurando l'infinitesimale di ordine superiore rappresentato dal secondo membro e, in generale, le differenze, pure infinitesimali di ordine superiore tra $\Delta q_{AB}$ e $C_p \frac{\delta T}{\delta v} dv$, $\Delta q_{BC}$ e $C_v \frac{\delta T}{\delta p} dp$. Quando invece s'intendesse prescindere da quel principio, la ricordata interpretazione dalla (1) si potrebbe sempre ricondurre al principio della indipendente coesistenza delle variazioni infinitesimali, la quale significa appunto l'ammissione dell'aggiunta di infinitesimi di ordine superiore.
Così se capita che, nell'epoca pretermodinamica, la (1) fu stabilita nel supposto che il secondo membro sia un differenziale esatto, non mi sembra che il principio dell'equivalenza sia necessario per scriverla senza fare tale ipotesi che si traduce nel principio della conservazione del calore.

In conclusione, a mio modo di vedere, la differenza fra le due ipotesi relativa al secondo membro della (1) - l'una che afferma, l'altra che nega la conservazione del calore - anziché nella considerazione di infinitesimi di ordine superiore, la cui traccia modificherebbe il significato dell'espressione - va tutta nel verificarsi o no la

$$\frac{\partial}{\partial p} C_p \frac{\partial T}{\partial v} = \frac{\partial}{\partial v} C_v \frac{\partial T}{\partial p},$$

da te appurato, in seguito, molto opportunamente rievocate.

S'intende che farai il caso che crederai di queste riflessioni, ne' occorre dire che non perturbano menomamente il corso del tuo manoscritto, che ho veduto con piacere anche come segno di ristabilimento pienamente conseguito.

Noi stiamo tutti bene, e abbiamo il guerriero[25] con noi nei suoi quindici giorni di licenza. Non lo vedevamo da un anno. Ma siamo in pena per il caro amico Pizzetti ammalato di polmonite. Oggi è la nona giornata e sembra annunciarsi un principio di crisi, per cui la prognosi si fa migliore, e speriamo di essere presto rassicurati, quand'anche non dovesse verificarsi quella soluzione discontinua, colla quale ogni pericolo cessa improvvisamente. Avrei voluto risparmiarti simile notizia, tardando anche a scriverti. Ma penso che possa avertela recata il prof Baravelli che, venuto ieri da Roma, ne ripartì in giornata.

Alla tua gentile Signora, alla Signora Lia e famiglia le più cordiali e distinte cose, e un abbraccio a te dal

Tuo affmo
Antonio.[26]

[25] Potrebbe trattarsi del figlio, Felice, che aveva 25 anni.
[26] Probabilmente Maggi scrive ad Ascoli anche la lettera #228 sullo stesso argomento, ma non la spedisce.

**3**
Gian Antonio Maggi a **Carlo Barzanò**

Pisa 6 Novembre 1908

*[Nota a margine]* All'Ing C. Barzanò non mandata V. la lettera al Prof Colombo.

Pregiatissimo signor Ingegnere,

Ho letto con molto interesse il suo articolo "Sul riordinamento degli studi nell'istituto tecnico superiore", e mi duole di non averne potuto ascoltare l'esposizione al Congresso di Firenze, al quale, quantunque iscritto, fui da una gita a Milano impedito d'intervenire. Ho particolarmente rilevato le Sue considerazioni sull'insegnamento della Meccanica Razionale; e sentendomi, con questo, quasi chiamato in causa, mi permetto di sopperire all'assenza collo scriverLe, e, chiedendoLe scusa della libertà, pregarLa d'indicarmi; con più estesi particolari, quali Ella giudica che debbano essere il programma e lo spirito di quell'insegnamento, per meglio soddisfare le esigenze della Scuola d'Applicazione degli Ingegneri.

Ch'io debba o no insegnare una volta o l'altra in una Scuola d'Applicazione il problema m'interessa molto per sé stesso: ed io sarei obbligatissimo a chi, con competenza nella Tecnica pari alla Sua mi fornisca gli elementi necessari per risolverlo.

Io mi permetto tanto più di farLe questa domanda perché ho veduto con piacere ch'Ella esprime dei desiderii che nel mio caso, da parecchio tempo procuro già di soddisfare. Io soglio far precedere l'impostazione del problema reale, discorrendo di projettili, di ruote ecc. a quello del problema ideale: mi studio costantemente di richiamare l'obbiettività fisica ad illustrazione della rappresentazione matematica: mi preservo di risolvere al minimo il macchinario analitico: soglio ancora sviluppare i principi della simmetria e della dinamica dei corpi variabili ponendo con ciò le basi della teoria dei sistemi elastici. Con tutto ciò io non mi dissimulo che il mio corso resta quel che si suol chiamare un corso teorico, per la doppia ragione ch'io mi sento pure tenuto a soffermarmi sulla fondazione di concetti che nella Pratica diventano poi grandezze come oggetto d'intuizione, e che tal corso, nella sua generalità, per essere utilizzato dalla Pratica richiede il tramite d'insegnamenti specialmente destinati a questo scopo. Né, per ben ragionevoli adattamenti ch'io vi facessi, mi sentirei in grado di mutarne tale carattere. Debbo ritenere che è, in massima, un corso teorico in questo senso che non si reputa conciliabile coi desiderii della scuola di Applicazione? A me pare che, se non si preferisce smembrare la Meccanica Razionale in tanti capitoli da porsi come introduzione ai varii insegnamenti tecnici - ciò che potrebbe essere una soluzione - ma si vuole invece conservarne l'integrità, essa non possa altro proporsi che l'esposizione dei metodi generali per la traduzione in formole matematiche dei problemi presentati dalla Fisica o dalla Pratica, e, con tale scopo, resti uno di quegli insegnamenti di cui i pratici bisogna che portino pazienza ad aspettare col tempo i frutti.

Sono pienamente d'accordo con Lei nella critica degli attuali ordinamenti del primo biennio di Facoltà, per quanto concerne l'insegnamento degli allievi ingegneri, o d'accordo nel reclamare pel Cotesto un posto privilegiato, in modo ch'esso possa riuscire un vero strumento alla mano. Ma non posso convenire sull'opportunità dello sdoppiamento in un medesimo luogo del primo biennio di Facoltà e di Politecnico, parendomi più conveniente riformare il primo biennio di Facoltà, anziché creargli a fianco un concorrente, che, levandogli la quasi totalità degli studenti, ne ridurrà all'estremo l'efficacia e l'intensità?

Di nuovo mi scuso della libertà che mi sono preso, e aggradisca i miei più distinti saluti, coi quali La prego di credermi

Dev$^{mo}$ Suo
Gian Antonio Maggi.

**4**
Gian Antonio Maggi a **Eugenio Bertini**

Pisa 11 Maggio 1897

Caro Bertini,

[Non mi consta che sia stata trattata qualche questione simile a quella proposta dal prof Uzielli; poiché non potrei chiamare analoghi, malgrado la coincidenza di certi dati, alcuni problemi di distribuzione di temperatura.]

Non mi pare che il problema proposto [D'altra parte], in quei termini, e, con quei soli dati, [non credo che] si possa risolvere.

Non vedo Innanzi tutto, [non vedo] perché cercare la quantità di calore, che si deve sottrarre al miscuglio per unità di tempo, mentre, mentre per ridurre la temperatura della sezione da $t_2$ a $t_3$, in un modo determinato (dato, questo modo, che non si trova fra quelli del proposto problema), occorrerà sottrarre una quantità di calore di fissa.

Ammesso che si cerchi quest'ultima, la sua grandezza assoluta sarà chiaramente

$$P_c(t_2 - t_3)+Q;$$

dove P e c indicano il peso del miscuglio che riempie la sezione, e il suo calore specifico: e Q la quantità di calore ceduto per conduttività dal contorno.

Ed è appunto pel calcolo di Q che bisogna stabilire la legge del raffreddamento: a seconda della quale Q avrà un valore diverso. Sia T la durata del raffreddamento: e supponiamo, per esempio, che, ad ogni istante, corrispondente alla fine del tempo $\tau$, contato dal principio dell'operazione, la temperatura del miscuglio sia uniforme, in tutto il cilindro, e di misura t. In tal caso, rappresentando con S l'area della superficie laterale della sezione, e con a il coefficiente di conduttività termica fra l'involucro solido e il miscuglio sarà:

$$Q=aS\int_0^T (t_1 - t)d\tau;$$

nella qual formula a ed S sono dati dal nostro problema, ma T e la funzione t di $\tau$ nell'intervallo (0T) dipendono dal processo di raffreddamento, col quale potranno variare in infiniti modi, salvo la condizione che t riceva i valori $t_2$ e $t_3$ per $\tau$=0 e $\tau$=T.

Su questa influenza del modo di raffreddamento ho discusso con Bardelli. Non posso rinunciare a ricordare a proposito la celebre esperienza degli Accademici del Cimento che una massa d'acqua tiepida circondata di ghiaccio fondente si raffredda più rapidamente circondando il ghiaccio d'acqua calda.

Se non che la questione pratica su cui la lettera insiste [il prof Uzielli] mi pare suggerisca un problema diverso, assai semplice. E cioè il calcolo della quantità di calore, che bisogna sottrarre nell'unità di tempo (allora sta bene) da una sezione determinata del cilindro, per mantenere il miscuglio a $t_3$, mentre l'involucro solido si mantiene a temperatura $t_1$ (senza preoccuparsi della quantità di calore, che sarà occorsa una volta tanto per abbassare la temperatura da $t_2$ a $t_3$).
In tal caso, la quantità di calore cercata non è altro che $aS(t_1 - t_3)$, se si suppone che $t_3$ si estenda al miscuglio fuori della sezione, o questa aumentata di [$2bA(t_2 - t_3)$, indicando con A l'area delle basi, e con b il relativo coefficiente di conduttività termica presso il miscuglio, se fuori della sezione la temperatura si suppone mantenersi $t_2$] di una quantità, che si rappresenta in modo analogo, se si suppone che fuori della sezione la temperatura sia diversa, eguale cioè alla quantità di calore ceduta nello stesso tempo per conduttività dell'involucro solido, o da questo e dalle basi [espressioni calcolabili immediatamente.]
Questo è quanto io so dire intorno al quesito proposto dal prof Uzielli, senza pretendere di esaurire la questione, ché sicuramente potrà dare una più completa risposta chi abbia in fatto di Fisica Tecnica la competenza che mi manca.

[Sono] sempre a tua disposizione, se hai qualche schiarimento da domandarmi, o da farmi qualche osservazione [.], [Intanto] mi è grato intanto mandarti i più cordiali saluti

Affmo collega
GAMaggi.

## 4 BIS[27]

Gian Antonio Maggi a Eugenio Bertini

Problema proposto dal prof Uzielli al prof Bertini con lettera in data di Firenze 21 aprile 97.

"Dato un cilindro circolare infinito pieno di un miscuglio di aria e di vapore d'acqua, e che trapassa una massa infinita solida a temperatura $t_1$, e dati i coefficienti termici della sostanza solida costituente la massa traversata, qual'è la quantità di calore che bisogna sottrarre, nell'unità di tempo, al miscuglio gassoso, in una sezione del cilindro compresa fra due piani paralleli e normali all'asse del cilindro, perché la temperatura $t_2$ del miscuglio (temperatura comunicata dalla massa solida a $t_1$, e quindi poco differente da $t_1$) si abbassi a $t_3$."

———————————————————————

Pisa 11 maggio 97

Caro Bertini,

Non mi consta che sia stata trattata qualche questione del genere di quella proposta dal prof. U.; poiché non potrei chiamare analoghi, malgrado la coincidenza di certi dati, alcuni problemi di distribuzione di temperatura.

D'altra parte, in quei termini, e, con quei soli dati, non credo che si possa risolvere.

Innanzi tutto, non vedo perché cercare la quantità di calore, che si deve sottrarre al miscuglio <u>per unità di tempo</u>, mentre, per ridurre la temperatura della sezione del cilindro da

[27] Si tratta di una seconda stesura della lettera precedente.

$t_2$ a $t_3$, nelle supposte condizioni, aggiunto il modo in cui s'intende che la sottrazione si faccia, occorrerà sottrarre una quantità di calore fissa.

Ammesso che si cerchi quest'ultima, la sua grandezza assoluta sarà chiaramente

$$P_c(t_2 - t_3)+Q,$$

dove P e c indicano il peso del miscuglio che riempie la sezione, e il suo calore specifico: e Q la quantità di calore ceduto per conduttività dell'involucro contorno - cioè dall'involucro solido mantenuto a $t_1$. Ed è appunto pel calcolo di Q, che bisogna fissare la legge del raffreddamento, a seconda della quale, avrà un valore diverso. Sia T il tempo in cui si compie il raffreddamento: e supponiamo, per esempio, che ad ogni istante, corrispondente alla fine del tempo $\tau$, contato dal principio dell'operazione, la temperatura del miscuglio sia uniforme, e di misura t. In tal caso, rappresentando con S la superficie della sezione, e con a il coefficiente di conduttività termica fra l'involucro solido e il miscuglio sarà

$$Q=aS\int_0^T(t_1 - t)d\tau;$$

nella qual formula a ed S sono dati dal problema proposto, ma T e la funzione t di $\tau$ nell'intervallo (0T) dipendono dal processo di raffreddamento col quale potranno variare in infiniti modi, salvo la condizione che t riceva i valori $t_2$ e $t_3$ per $\tau=0$ e $\tau=T$.

[Su questa influenza del modo di raffreddamento ho discorso anche col Bardelli. Non posso rinunciare a ricordare a proposito l'esperienza degli Accademici del Cimento che una massa d'acqua tiepida circondata di ghiaccio fondente si raffredda più rapidamente circondando il ghiaccio d'acqua calda.]

Se non che la questione pratica su cui insiste il prof U. mi pare suggerire un problema diverso, e assai semplice. E cioè il calcolo della quantità di calore, che bisogna sottrarre nell'unità di tempo (allora sta bene) da una sezione determinata del cilindro, per mantenere il miscuglio alla temperatura $t_3$, mentre l'involucro solido si conserva alla temperatura $t_1$ (senza preoccuparsi della quantità di calore che sarà occorsa una volta tanto per abbassare la temperatura da $t_2$ a $t_3$.) In tal caso, la quantità di calore cercata non è altro che $aS(t_1 - t_3)$,se si suppone che $t_3$ regni in tutto il cilindro, o questa aumentata di $2bA(t_2 - t_3)$, indicando, con A l'area delle basi della sezione, e con b il relativo coefficiente di conduttività termica verso il miscuglio, se fuori della sezione la temperatura si suppone mantenersi $t_2$ - eguale cioè alla quantità ceduta nello stesso tempo per conduttività dall'involucro solido, o da questo e dalle basi - espressioni calcolabili immediatamente.

Sono sempre a tua disposizione, se hai qualche schiarimento da domandarmi, o da farmi qualche osservazione. Intanto m'è grato mandarti i più cordiali saluti.

Affmo collega
GAM

**5**

Gian Antonio Maggi a **Giuseppina Biggiogero Masotti**
[cartolina postale indirizzata a:
Gent.[a] Sig.[a] Dott. Giuseppina Biggiogero - Melegnano - (Milano);
intestata R. Università di Milano - Istituto Matematico
Via C. Saldini, 50 - Milano - Città degli Studi]

Milano 19 Dicembre 1930

Gentilissima Signorina,

Baricentro non si ricava da Archimede il quale usa κέντρον τοῦ βάρεος, che, leggendo βάρεος per βαρέος, significa centro del grave. Per quanto allo spostamento dell'accento, vedo che Eutocio Ascalonita, commentatore di Archimede, scrive βαρέος, per cui suppongo la differenza da attribuirsi al dialetto dorico, usato dal Siracusano.
La traduzione latina reca centrum gravitatis

Con questo Le rinnovo augurii e saluti, e mi confermo Suo

Aff.mo G.A. Maggi.

**6**

Gian Antonio Maggi a **Vittorio Emanuele Boccara**
[s.d., ma probabilmente del 1894]

Egregio Sigr Boccara,

Ricevo con qualche ritardo il suo manoscritto, mandatomi qui da Milano, donde manco da parecchio tempo: e glielo rendo colle mie osservazioni.

Mi limito all'essenziale. E noto, innanzitutto, che due volte, per formar l'integral generale dell'equazione del movimento, ricorre all'espressione:

$$A_1 y_1 + A_2 y_2,$$

dove $y_1$, $y_2$ sono due integrali particolari, senza avvertire che quelli che Ella prende sono legati da una semplice relazione lineare, per modo che non soddisfanno alla condizione di costituire un sistema fondamentale. Il quale è un grosso guajo. Ma ancora di più grossi ne troveremo in seguito.

Poiché Ella, cercando l'integral generale dell'equazione

$$\left.\begin{array}{c} l\dfrac{d^2\Theta}{dt^2}+2\dfrac{dl}{dt}\dfrac{d\Theta}{dt}+g\Theta=0 \\ (l=gt^2+Bt+C) \end{array}\right\} \quad (1)$$

per mezzo della formula

$$y=y_1\int\frac{1}{y_1^2}e^{-\int\frac{B}{A}dt}dt$$

$$y_1=\frac{g}{B+gt\pm\sqrt{B^2-2gC}}$$

dimentica di moltiplicare per $y_1$ la costante additiva della seconda relazione (da Lei designata con E): e conclude, sostituendo, che non può aver altro valore all'infuori di O - il solo caso infatti in cui si elimina quell'errore. E questo doveva ben metterla in guardia: poiché la sua conclusione si riduce a questo; che l'equazione in discorso non assumerebbe un integral generale con due costanti arbitrarie, vale a dire integral generale, che se guardiamo il risultato dal punto di vista meccanico, vorrebbe dire che si ha il movimento un punto, al quale non si può assegnare ad un istante posto e velocità a piacere.

Ciò ch'Ella non mostra d'avvertire: ma richiamando l'espressione precedentemente trovata - erronea, per altro ordine di ragioni, come s'è veduto - finisce per identificarla colla nuova, e trovare una condizione per le costanti $\alpha_2$, $\beta_2$, che il procedimento di deduzione di quell'integrale dimostra chiaramente come non possa essere richiesta.

In tutto ciò c'è persino troppo! debolezze di fondamenti, e, indipendentemente dalle cognizioni matematiche, mancanza di ragionamento: e bisogna ch'Ella si rinsaldi in quelli, e curi ben questo, per arrivare a qualche attendibile conclusione.

Effettivamente l'integral generale delle (1) è dato da:

$$\Theta l=\alpha+\beta t \quad (2)$$

($\alpha$, $\beta$ costanti arbitrarie). Risultato che si può trovare immediatamente, osservando che, pel teorema di Leibniz, la (1) si può porre sotto la forma:

$$\frac{d^2(\Theta l)}{dt^2}=0.$$

E si vede anche che l'ipotesi $A=\frac{g}{2}$ è costituisce, in certo qual modo, un caso singolare: di semplicità presumibilmente eccezionale; per modo che, a parte qualche considerazione intorno alla (2), bisognerà cercare altrove materia di lavoro.

Al qual fine proverà che riprenda il caso generale, e segua la ricerca del Lecornu, che Le dovrebbe servire come una falsa-riga, secondo il consiglio che Le ho dato ripetutamente. Ma prima riveda le equazioni differenziali lineari specialmente del 2° ordine, e le teorie ipergeometriche

Per queste vedrà con profitto la monografia del prof. Pincherle nel giornale di Battaglini[28] (ultimo o penultimo volume); e per quelle, oltre i trattati, il capitolo dedicatovi nella Physique Mathématique del Mathieu.[29]

Per quanto al pendolo conico, ne trova la trattazione nei limiti che Le possono occorrere, in qualunque trattato di Meccanica Razionale: e può anche cercarlo nel Durege - Elliptische Functionen[30] dove, se la memoria non mi tradisce, v'è un'esposizione completa.

Mi lasci credere che questi avvertimenti Le riescano di vantaggio, e riceva i mie saluti

Aff Suo

[28] S. Pincherle, "Delle funzioni ipergeometriche e di varie questioni ad esse attinenti", *Giornale di matematiche di Battaglini*, Napoli, 32, 1894, p. 411.
[29] É. Mathieu, *Traité de physique mathématique*, Gauthier-Villars, Paris, 1873-1890.
[30] H. Durege, *Theorie der Elliptischen Functionen*, E.G. Teubner, Leipzig, 1861.

**7**

Gian Antonio Maggi a **Tommaso Boggio**

Pisa 17 Aprile 1924

Caro Professore Boggio,

La ringrazio vivamente del dono, e La prego di voler assumersi l'espresso compito di ringraziare da parte mia il prof. Burali-Forti: più espliciti ringraziamenti che aggiungo a quelli che riceveranno, colla mia recente commemorazione ai Lincei del Van der Waals.[31]

Vi ho unito, oltre ad una recensione di un libretto di Nørlund,[32] la mia Nota ai Lincei "Sulle varie interpretazioni della trasformazione di Lorentz".[33] Apparisce da questo breve scritto come io non divida l'entusiasmo di altri per la teoria della Relatività. Tuttavia, mettendomi questa volta dalla parte dei relativisti, piuttosto che dalla contraria, io ho inteso, con quella Nota, di opporre una ragione pregiudiziale alle conclusioni del mio ottimo amico Somigliana, da Loro adottate e riprodotte: quella che non sono affatto paragonabili tra loro gli usi che della trasformazione in discorso fa il Voigt, senza decampare dai principii newtoniani, e fa l'Einstein, in base ai principii della relatività, l'uno e l'altro, per la teoria matematica del fenomeno di Doppler. La differenza mi sembra abbastanza chiaramente spiegata. Invece ho lasciato quasi interamente al lettore di cavare la conseguenza principale del mio discorso, che non si può quindi accettare senza riserve l'affermazione del Somigliana, ricordata in principio. E, per esempio, io domanderei che l'affermazione che la trasformazione di Lorentz può essere applicata coll'interpretazione newtoniana, usata pel fenomeno di Doppler, alla spiegazione dell'esperimento di Michelson, fosse comprovata dalla relativa dimostrazione. O vogliamo proprio tornare alla <u>concreta</u> contrazione di Lorentz? Colpa la forma... dolce, che Loro converranno di non aver adoperato coi relativisti, credo che la mia critica, sostanzialmente concettuale, sia passata inavvertita.

Io poi mi sono limitato alla parte essenziale, ma, a proposito della priorità del Voigt, avrei potuto osservare che il Lorentz applicò la nota trasformazione, non già all'equazione di d'Alembert, ma alle equazioni di Hertz, che, col calcolo ipervettoriale del Minkowski, vuol dire, com'Ella ben sa, alla divergenza di un sestivettore, che ne riesce trasformata in sé stessa. Ma ancora, osservazione assai più importante, le molteplici sostituzioni che il Somigliana deduce dal suo principio, e pone accanto alla sostituzione di Lorentz, diventata di Voigt-Lorentz, servono come questa, finché si tratta dell'equazione $\frac{1}{a^2}\frac{\delta^2\varphi}{\delta t^2} - \Delta_2\varphi = 0$, ma non egualmente per la $\frac{1}{a^2}\frac{\delta^2\varphi}{\delta t^2} - \Delta_2\varphi = \psi$, che occorre, più in generale, di considerare; poiché quelle sostituzioni, ad eccezione della sostituzione di Lorentz, riproducono il dalambertiano $\frac{1}{a^2}\frac{\delta^2\varphi}{\delta t^2} - \Delta_2\varphi$, moltiplicato per un fattore. <u>Ergo</u>, non sta precisamente che, accanto alla teoria di Einstein, se ne possa formare una moltitudine d'altre equivalenti, almeno per quella via. Di tutto questo ho parlato e scritto all'amico Somigliana, senza che mai siamo riusciti a intenderci reciprocamente.

[31] "Commemorazione del Socio straniero Van der Waals", *Rendiconti della R. Accademia dei Lincei*, v. XXXIII, pp. 152-159.

[32] "Recensione dell'opera: N.E. Nörlund, Videnskabelige Causerier (Conferenze scientifiche)", *Bollettino della Unione Matematica Italiana*, v. III, pp. 32-35.

[33] *Rendiconti della R. Accademia dei Lincei*, v. XXXII, pp. 196-197.

Con questo, che a Loro potrà parere la parte dell'Avvocato del Diavolo, Le presento le osservazioni che ho potuto fare ad una prima visione del volume, in conseguenza delle analoghe discussioni col Somigliana. Mi riserbo poi di dedicare al resto la dovuta attenzione, disposto anche a trarne la conclusione di non aver prima capito quello che, ad ogni modo, ero persuaso di aver capito.

Aggradisca, con nuovi ringraziamenti, per sé e pel prof. Burali-Forti, i miei migliori saluti e augurii di Buona Pasqua, e mi creda sempre

Suo affmo collega
Gian Antonio Maggi.

**7 BIS**[34]
Gian Antonio Maggi a Tommaso Boggio

Pisa 17 Aprile 1924

Caro Professore Boggio,

La ringrazio vivamente del dono, e La prego di voler assumersi l'espresso compito di ringraziare da parte mia il prof. Burali-Forti: più espliciti ringraziamenti, che aggiungo a quelli che riceveranno, colla mia recente commemorazione ai Lincei del Van der Waals.

Vi ho unito, oltre ad una recensione di un volumetto del Nørlund, la mia Nota ai Lincei "Sulle varie interpretazioni della trasformazione di Lorentz". Apparisce da questo breve scritto come io non divido l'entusiasmo di altri per la teoria della Relatività. Tuttavia, mettendomi, questa volta, dalla parte dei relativisti, piuttosto che dalla contraria, io ho inteso, con quella Nota, di opporre una ragione pregiudiziale alle conclusioni del mio ottimo amico Somigliana, da Loro riprodotte e adottate: quella che l'uso che della trasformazione in discorso fa il Voigt, senza decampare dai principii newtoniani per la teoria del fenomeno di Doppler, e l'uso che, allo stesso scopo, ne fa l'Einstein, in base ai principii della Relatività, non sono affatto paragonabili fra loro. La differenza mi sembra abbastanza chiaramente spiegata. Invece ho lasciato quasi interamente al lettore di cavare la conseguenza principale del mio discorso, che non si può quindi accettare senza riserve l'affermazione del Somigliana, ricordata in principio. E, per esempio, io domanderei che l'affermazione che la trasformazione di Lorentz, coll'interpretazione newtoniana, usata pel fenomeno di Doppler, può essere applicata alla spiegazione dell'esperimento di Michelson, fosse comprovata dalla relativa dimostrazione. O vogliamo proprio tornare alla contrazione di Lorentz concreta? Colpa la forma... dolce, che Loro converranno di non aver adoperato coi relativisti, credo che la mia critica, sostanzialmente concettuale, sia passata inavvertita.

Io poi mi sono limitato alla parte essenziale, ma, a proposito della priorità del Voigt, avrei potuto osservare che il Lorentz applicò la nota sostituzione, non alla trasformazione dell'equazione di d'Alembert, ma a quella delle equazioni di Hertz, che, col calcolo ipervettoriale del Minkowski, vuol dire, com'Ella ben sa, alla divergenza di un sestivettore, che ne riesce trasformato in sé stesso. Ma ancora, osservazione assai più importante, le molteplici sostituzioni che il Somigliana deduce dal suo principio, e pone accanto alla sostituzione di Lorentz, diventata di Voigt-Lorentz, servono come questa, finché si tratta dell'equazione $\square\varphi=0$, ma non egualmente per la $\square\varphi=\psi$, che occorre, più generalmente, di considerare: poiché quelle sostituzioni, ad eccezione di quella di Lorentz, riproducono il dalambertiano $\square$, moltiplicato per un fattore. Ergo, non sta precisamente che, accanto alla teoria di Einstein, se ne possa formare una moltitudine d'altre equivalenti, almeno per quella via. Di tutto ciò ho parlato e scritto all'amico Somigliana, senza che mai siamo riusciti a intenderci reciprocamente.

[34] Si tratta della prima stesura della lettera precedente.

Con questo, che a Loro potrà parere la parte dell'Avvocato del Diavolo, Le presento le osservazioni che ho potuto fare ad una prima visione del volume, in conseguenza delle analoghe discussioni col Somigliana. Mi riserbo poi di dedicare al resto la dovuta attenzione, disposto anche a trarne la conclusione di non aver prima capito quello che, ad ogni modo, ero persuaso di aver capito.

Aggradisca, con nuovi ringraziamenti, per sé e pel prof. Burali-Forti, i miei migliori saluti e augurii di Buona Pasqua, e mi creda sempre

Suo affmo collega
Gian Antonio Maggi.

**8**

**Ettore Bortolotti** a Gian Antonio Maggi
[busta intestata Unione Matematica Italiana e
indirizzata Al Ch.mo Signor Prof G. A. Maggi della R. Università Pisa;
lettera su carta intestata Regia Università di Bologna -
Facoltà di scienze fisiche, matematiche e naturali]

21.V 1924

Carissimo,

il famoso Ivaldi manda per il nostro Bollettino una piccola Nota dove dice di aver fatto non so che esperienza con l'apparecchio del tuo amico Cirinei.[35]

Mando a te la nota, pregandoti, anche da parte del presidente, che tu ci dica se si può lasciar passare, e, nel caso negativo perché ci suggerisca un grazioso eufemismo con cui si possa esprimere senza parer scortesi, il nostro rifiuto all'Illmo Signor prof Ing Ivaldi.
Ti prego di salutarmi il nostro prof Bianchi, cui dirai che ho consegnato subito i manoscritti di Pincherle; e di ricordarmi con affetto ai colleghi.

Tuo aff[mo]
Ettore Bortolotti

[35] Ivaldi aveva pubblicato un articolo su questo argomento: "Il rendimento delle macchine", *La Scienza per tutti*, n. 6, marzo 1921. In tale articolo - che si trova in rete - descrive alcuni esperimenti del Cav. Egisto Cirinei e conclude così: "Questo aumento *[del coefficiente di rendimento]* è stato ammesso e riconosciuto dal Senatore Blaserna, in una lettera diretta al Cav. Cirinei. Esso è incompatibile con la seconda legge della termodinamica, legge che noi abbiamo sempre affermato, come recisissimamente affermiamo, essere falsa, perché è falso che l'energia di un fluido sia indipendente dalla pressione."

**9**

Gian Antonio Maggi a Ettore Bortolotti

Pisa 25 Maggio 1924

Carissimo,

L'Ill.mo Sig.r Prof. Ing. Ivaldi, supposto un punto materiale di massa m, fissato ad un estremità *[sic!]* di una leva, rappresentata da una retta, girevole intorno al proprio fulcro, a distanza r da questo, e assoggettato ad una forza motrice di grandezza assoluta f, e retta d'applicazione, la cui distanza dal fulcro della leva è b, indicando poi con a la grandezza assoluta dell'accelerazione tangenziale (da lui chiamata lineare) del punto, al tempo generico, deduce conv.mente dal teorema della forza viva (da lui chiamato principio delle energie di moto) la relazione:

(1) $$fb=mra.$$

Ma lo stesso Ill.mo Sig.r ecc. ecc. regala poi un proprio grosso sproposito alla Meccanica, facendole dedurre, a sentir lui, dal secondo principio newtoniano, nel caso del movimento in discorso, la relazione

(2) $$f=ma.$$

La Meccanica, nello stesso caso, ne deduce invece la relazione

(3) $$f\cos\alpha=ma,$$

dove α indica l'angolo acuto formato dalla retta di applicazione della forza motrice e dalla perpendicolare alla retta rappresentante la leva.

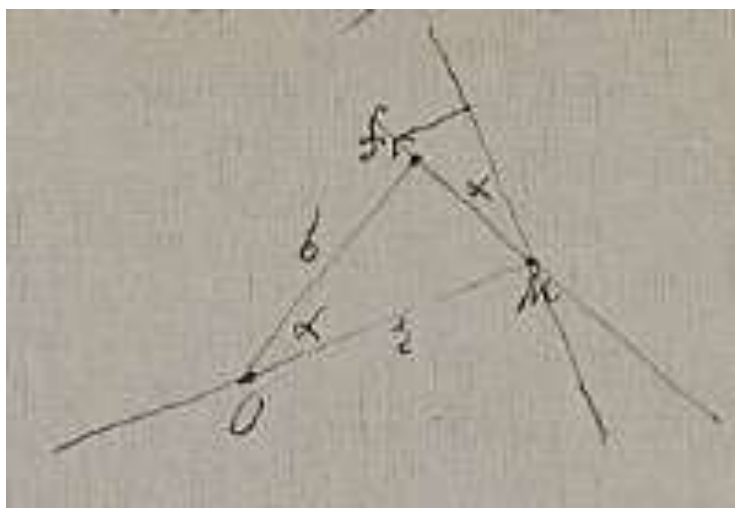

E poiché, manifestamente, si ha

$$\frac{b}{a}=\cos\alpha,$$

Le equazioni (1) e (3), si deducono, a piacere, l'una dall'altra.

Fin qua ci arrivo, ma davvero non arrivo a trovare un grazioso eufemismo, per rifiutare quel gioiello.

Per conto mio, mando il Prof. Ing. Ivaldi a imparare che il secondo principio newtoniano della Dinamica essendo "vettore forza motrice (incluso, nel caso di punto vincolato, la reazione del vincolo, parallelo, nel supposto movimento, alla retta-leva) = prodotto della massa pel vettore accelerazione (totale)", se ne ricava la (3) prendendo la componente, secondo la tangente alla trajettoria di ambedue i membri della corrispondente equazione, e non del solo secondo membro.

Imparare ho poi detto, per modo di dire. Perché so, per esperienza, come questa specie di paranoia sia incurabile; e dalla discussione con loro non ci sia da ricavare che una rottura di tasche, da dover poi ricorrere, per ripararle, …a Voronoff![36]

Le nostre migliori alla tua cara famiglia e a te, cordiali saluti al Pincherle, al Tonelli, al Burgatti, col quale ho avuto il piacere di trovarmi, gli scorsi giorni, a Roma ancora un abbraccio a te dal

Tuo affmo collega
Gian Antonio Maggi.

[36] Probabilmente Maggi si riferisce a Serge Voronoff, chirurgo e scienziato molto noto in quegli anni.

**10**
Gian Antonio Maggi a **Mario Canavari**

Pisa 15 Aprile 1903

Caro Canavari,

Una volta, quando riportai il noto M.S. non ti trovai in Laboratorio, un’altra volta non riuscii a farmi aprire la porta, né, in questi giorni, ebbi occasione di vederti altrimenti. Ti mando quindi brevemente per iscritto le mie osservazioni.

L’autore non impiega che la somma dei vettori, vulg. composizione delle forze: e se ne serve per formare un nuovo ente, che chiama combinante di due vettori. Egli parte perciò da un’analogia meccanica, dove, ad ogni modo, occorre sostituire quantità di moto a velocità, ma ancora il concetto vuol essere meglio chiarito, e i fenomeni mancano di precisione. Definisce poi il combinante con un doppio significato, e sia ora di un solo vettore ed ora di insieme di due vettori eguali ed opposti. Così dice espressamente di fare, ma io troverei necessario di giustificare in qualche modo tale procedimento. Infine, quando si verificano certe circostanze, cambia completamente la definizione, affermando come evidente cosa di cui non mi pare apparisca la ragione. In conclusione, il nuovo ente non mi risulta fondato con sufficiente precisione e chiarezza.

Mi sono fermato, secondo la nostra intelligenza, alle prime pagine. Se renderne ragione è riserbato alla parte chimica, che segue, non è argomento di mia competenza.

Aggradisci i miei migliori saluti, e credimi

Tuo aff. collega
Gian Antonio Maggi.

**11**
**Ignazio Canestrelli** a Gian Antonio Maggi
[busta indirizzata a: Illustre Prof. G.A. Maggi di meccanica razionale nella R. università a riposo Corso plebiscito N.° 3 Milano]

Roma 30/5/1936[37]

Multi fuerunt a me vocati, pauci vero electi.
Illustre professore
Roma Lungo Tevere Sanzio N.°2 (Ponte Sisto)

L’espressione in calce è quella del mio interno affanno con il quale depongo la mia spregiata lacrima a coloro che si distinsero sempre per nobiltà d’animo ed altruismo sotto ogni rapporto.
E ne comprenderà la ragione dopo compresa la ragione della presente.
Un poco di storia innanzi tutto. Dopo d’aver io insegnato fisica per 40 anni, come ordinario nei R. Licei ed anche come incaricato nei R istituti tecnici, ho inteso il bisogno spirituale di cercare di comprendere i mirabili principi della meccanica. Tale bisogno cominciò in me dall’86 (mezzo secolo fa)! quando il tema dato agli esami dal ministero per la licenza liceale fu risolto e con varietà di forma dai miei scolari al Liceo d’Aquila. Il compianto prof. Giulio

[37] Dal timbro postale.

Pittarelli che in quell'anno insegnava geom. descrittiva in quell'istituto tecnico e col quale io mi consigliai per alcune classificazioni, volle che io avessi dato 9/10 allo Alunno Stromei, e 10/10 all'alunno Oreste Ranelletti (ora prof. di diritto ammin. alla univ. di Milano; e parecchie altre classificazioni elevate. Vi furono soluzioni di forma varia, anche ingegnose del tema, che conservo. Io avevo esercitato per l'intero anno molto gli scolari, ai quali più che Fisica svolgeva esercizii di meccanica nei vari capitoli della meccanica e della Fisica, tratti da opere di problemi, fattici venire specialm. dalla Francia come lo *[*Chevellier*]*, il Jacquier[38] il *[*Bonant*]* etc.
Ella ricorderà, come generalmente nei vari Licei del Regno fu data agli alunni carta bianca. In Napoli il tema fu risolto solo nel Liceo dove insegnava il Pinto, professore di fisica matematica in quell'università; a Roma solo nel Liceo dove insegnava il Perotti già prof di S.M. il Re d'Italia, cui fu proposto dal mio maestro il Blaserna. In seguito, il mio desiderio di esporre con chiarezza i principii della meccanica fu da me coltivato e si fece più vivo quando il Levi-Civita mi fece avere la relazione al convegno didattico tenuto nel Giugno 1922 dai professori Levi Civita-Marcolongo-Volterra all'istituto nazionale d'educazione profession. in Roma per invito del direttore di quell'istituto ing. Andreoni perché avessero dato le norme per un trattato di meccanica per gli alunni delle scuole professionali.
Così il Levi Civita trattò della cinematica; il Marcolongo della statica, il Volterra della dinamica. Poi dopo andato io in riposo nel 23, nei miei 12 anni di pensione già trascorsi io mi occupai di ritirare le sarte[39] su ciò che avevo rilevato dai vari testi, che tutti più o meno furono uno dopo l'altro usati da me nella scuola. Cito tra i più severi il Roiti[40] (grande edizione) il Pinto,[41] il Pitoni[42] (che venne dalla scuola di Pisa), il Vanni e Monti.[43] Ma non escluderei Felice Marco,[44] il quale sebbene prolisso, puerile ed anche inesatto in alcuni punti, fu sempre un libro didattico (Il prof Palazzo mi diceva che sul Marco aveva capito la meccani. elementare), mentre il Roiti per quanto mirabile sintesi dei principi della meccanica non è affatto didattico. Non parliamo del Battelli[45] (mancante di unità), del Amaduzzi[46] parolaio, di libri o testi usciti da poco come il Fermi,[47] di altri come il Corbino,[48] poco sviluppati e pochissimo usati: anche per poca chiarezza.
In questi anni di riposo, io mi illusi di avere degli ausili dagli esperti nella meccanica. Veramente il Levicivita per preghiera del suo collega compianto Pittarelli cominciò a darmi chiarimenti. Ma mi concesse un unico convegno in sua casa. Poi si scusò che era troppo occupato.
Col Volterra, sebbene raccomandato a lui dal Prof Dilegge nulla conclusi.
Il Marcolongo finì per consigliarmi di piantarla.
Il Giorgi (ora prof. alla scuola di Ingegneria di Roma mi disse che mi era messo a pelare una gatta, ma che la gatta occorre pelarla da se. *[sic!]*

---

[38] F. Jacquier, *Elementi di perspettiva*, Salomoni, Giovanni Generoso, 1755.
[39] Forse nel senso di "ritirare le vele".
[40] A. Roiti, *Elementi di Fisica*, Le Monnier, 1898.
[41] L. Pinto, *Trattato elementare di fisica*, Morano, Napoli, 1906.
[42] R. Pitoni, *Storia della fisica*, Società tipografico–editrice nazionale, 1913.
[43] G. Vanni, V. Monti, *Corso di Fisica ad uso dei Licei*, Vallardi, Milano, 1906.
[44] Felice Marco, *Elementi di Fisica per Licei, Istituti tecnici, Scuole tecniche e magistrali*, Paravia, Torino, 1887.
[45] A. Battelli, P. Cardani, *Trattato di fisica sperimentale dei proff. Battelli e Cardani*, Milano, Vallardi, 1932.
[46] L. Amaduzzi, *Elementi di Fisica*, Zanichelli, 1921.
[47] Ad esempio: E. Fermi, *Fisica ad uso dei Licei*, Zanichelli, 1929.
[48] O.M. Corbino, *Nozioni di Fisica, per le scuole secondarie*, R. Sandron, Palermo, 1908.

Il Silla che non aveva tempo di sentirmi. Badi che a queste degne pressioni fui vivamente raccomandato e per lettera ed oralmente dal prof. Dilegge.
Il Marcolongo anziché limitarsi a consigliarmi di esprimere con chiarezza i suddetti principii vorrebbe si portassero in questi delle novità così p. es. dedurre la legge del parallel. dalle forze dal principio delle delle forze parallele (come anche dice in quei 2 modesti!!! volumetti Hoepli *[la parentesi non è chiusa: sic!]* (che il Dilegge assolutamente non voleva darmi dicendomi che non ne avrei mai capito il contenuto). In questi il Marcolongo, *[*propone*]* come cosa facile a *[*giungere*]* dedurre il principio del parallelogr. forze ortogonali, poi mano mano si passerebbe su forze oblique. A me disse a voce che con sole cinque lezioni ci si arriva. Ma... domando io dove si prende il tempo di 5 lezioni per spiegare il parallel. delle forze quando appena nelle scuole medie vi si può dare una sola lezione!
Un trattato completo di meccanica sia per scuole classiche sia per scuole industriali in Italia è assolutamente ricercato non solo dagli scolari, ma anche degli insegnanti i quali usciti dalle università, piena la testa di mat. superiori cogli orari onerosi non hanno tempo a consultare opere per farsi un corso di lezioni.
Veramente esiste un piccolo trattato del Burgatti,[49] il quale nella prefazione dice di avere risoluto la questione. Ma tutti gli insegnanti di Roma in coro si sono espressi che se il testo è tra i più rigorosi è tra i più difficili perché si esprime troppo analiticamente.
Io avevo pensato di scriverne al Burgatti, ma un suo collega (di Geom. all'università di Bologna) che vidi a Roma, mi disse che il Burgatti non risponde neppure ai convenevoli che gli fanno gli amici.
Dopo di che comprenderà se io ho per avere ausili espatriato da Roma, in oggi quanto mai più imperiale! e dopo di avere tentato inutilmente altre vie (come p. es. il Somigliana) mi sia rivolto a Lei, che il figliuolo del prof Mauri (Gian Carlo) mi disse essere persona molto buona e cortese, e mi consigliò di scriverle direttamente e senza bisogno di intermediario.
Io mi auguro che Ella stia ancora bene. Ma certo si sentirà stanco sebbene voglia tenersi ancora sulla breccia, dappoiché mi fu detto recentemente che sta riordinando le Sue memorie. Del resto anche per me non si scherza. Il 12 del prossimo agosto valico l'ottantesimo di età pur sentendomi in gamba organicamente e di spirito. Dunque io mi sono proposto di rivolgermi a Lei; ma col pianino cioè solo per avere da Lei a mia richiesta qualche chiarimento o consiglio, ora che vado mettendo assieme i varii capitoli del testo quasi ultimato. Così p. es. io lessi sopra un trattato di una Sua nota «cosa è la forza centripeta». Se ne avesse una copia potrebbe mandarmela, o dirmi almeno dove potrei leggerla?
E così per altre cose del genere, sempre su tema elementare, che valgono a chiarire certi punti di non facile intuizione. Poi in seguito potrei chiederle qualche chiarimento su questioni speciali sull'ordine migliore di trattazione dei vari capitoli.
Il prof Di Legge che entra nell'89^mo^ di età il 2 luglio prossimo è allettato su poltrona dopo una caduta in casa che gli ha rotto il collo del femore. Del resto stà abbastanza bene ed ha lucidità di mente. Se Ella capitasse in Roma, come per l'occasione della adunanza reale all'accademia dei Lincei, mi sentirei onorato di poter venire ad ossequiarla.
Perdoni la libertà, con profondo rispetto ossequiandola.
Dev.mo Prof Ignazio Canestrelli, dottore in Fisica, Roma Lungo Tevere Sanzio 2.

[49] P. Burgatti, *Lezioni di meccanica razionale*, 1918.

*[aggiunte a lato nelle varie pagine:]*

A completamento del testo di meccanica mi manca il potenziale gravitazionale. Il Roiti mi diceva non conoscere chi lo abbia trattato meglio di *[*lui*]* nel suo testo. È fatto bene ma non troppo facile. Ella saprebbe dove elementarmente fosse trattato in modo più facile?

Scrivere un trattato elementare è più difficile che un trattato Superiore. Trattasi di studio meticoloso che non finisce mai. Come la perfezione è assintotica *[sic!]*.

Mi sono servito anche del trattato Levicivita *[sic!]* Amaldi.

Nel mio testo il moto armonico, trascurato negli altri testi che precedono la Fisica è trattato ampiamente.

Il prof. Somigliana mi consigliava il Perucca[50] come trattato elementare di mecc. e fisica. Ma è fatto col calcolo superiore, è testo per ingegneria.

Mi sono servito anche delle opere francesi, come la meccanica del Gabriel tanto lodata dal Marcolongo.

**12**

Gian Antonio Maggi a Ignazio Canestrelli

Milano 2 Giugno 1936

Pregiatissimo Signor Professore,

"Vocatus" non mancherò certo di rispondere, ma spiacemi di dover molto dubitare di riuscire l'"electus". Poiché Ella mi nomina, in fatto di Meccanica, tutti i miei conoscenti, che sembrano non siano riusciti ad esserlo, e fra questi parecchi, che, per l'insegnamento medio - che, dai trattati che ricorda e dal resto, apparisce ciò che particolarmente Le interessa - hanno competenza ben altrimenti superiore alla mia. Io sarei ben lieto di corrispondere adeguatamente alla raccomandazione del mio illustre collega Prof. Di Legge[51] (al quale invio augurii di presto e intero ristabilimento), ma credo di dover in precedenza metterla in guardia che non mi reputo persona atta a servirLa più e meglio degli altri che mi scrive di essersi rivolto. Un saggio delle risposte che potrei darLe le fornisce il mio articolo di cui, secondo l'accennato Suo desiderio, Le invio una copia, "Che cos'è la forza centrifuga?"

Con questo gradirò pur sempre che mi scriva. In tanto La prego di aggradire i miei saluti e di credermi

Dev.mo Suo
Gian Antonio Maggi.

[50] Potrebbe trattarsi di: E. Perucca, *Fisica Generale e Sperimentale*, Avalle, Torino, 1931.
[51] Cfr. lettera #59.

**13**
Gian Antonio Maggi a **[Tommaso] Cannizzaro**

Casa 10 Luglio 94

Chiarissimo Signore,

Non lascerò Messina[52] senza mandarLe una risposta intorno alla soluzione della trisezione dell'angolo su cui Ella mi fece l'onore di domandare il mio giudizio; che se non Le ho scritto prima, Ella ne sa le tristi ragioni.

Mi spiace di non poterla approvare. L'autore che, cominciando dalla forma dell'espressione, non rivela troppa famigliarità col metodo matematico, fa un seguito d'asserzioni che dovrebbe, prima di tutto, dimostrare.

Se non che si può agevolmente riconoscere che alcune non si potrebbero in alcun modo giustificare perché conducono ad una conseguenza assurda. Mi limiterò all'esempio seguente. L'autore asserisce (V l'annessa figura)
KF=PC,
per modo che, levando la parte comune PK, si ha per legittima conseguenza:
FP=KC (1);
e così pure asserisce che
HC=PL (2).
Ora, almeno una di queste due asserzioni non può esser vera. Difatti, si ha
PG=LG–PL,
donde, essendo LG=KF e HC=HK+KC e KF–HK=FH.

2 1

Se fosse essere vera (2), dovrebb'essere
PG=FH–KC.
E allora, ammessa la (1), si trova:
FP+PG=FH,
donde, essendo FH=GF, anche:
FP+PG=GF,
cioè la somma di due lati del triangolo FPG eguale al terzo lato.

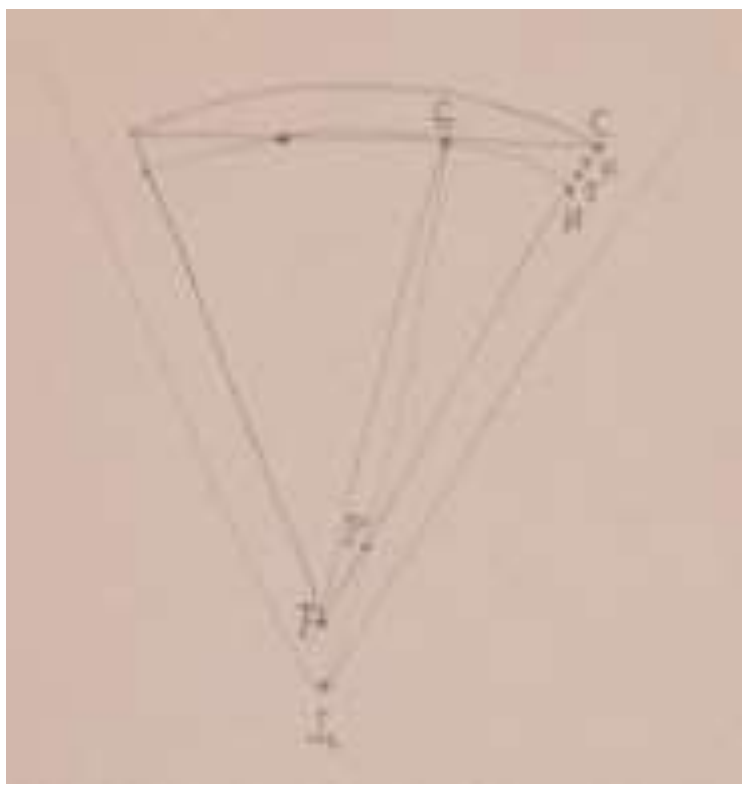

Mi è grato etc.

All'Illustr. Prof. T Cannizzaro
Messina

(Soluzione proposta da un incognito finlandese).

[52] Maggi insegnò Analisi Infinitesimale presso l'Università di Messina dal 1886 al 1895. Tommaso Cannizzaro e Giuseppe Zuppone-Strani (si veda la lettera #226) erano due poeti di Messina. Il nonno di Maggi aveva curato nel 1862 l'edizione italiana dei *Lusiadi* di Camoens e anche i due poeti si interessavano di quell'autore; inoltre erano tutti colti linguisti e orientalisti, quindi probabilmente si conobbero mentre Maggi si trovava lì.

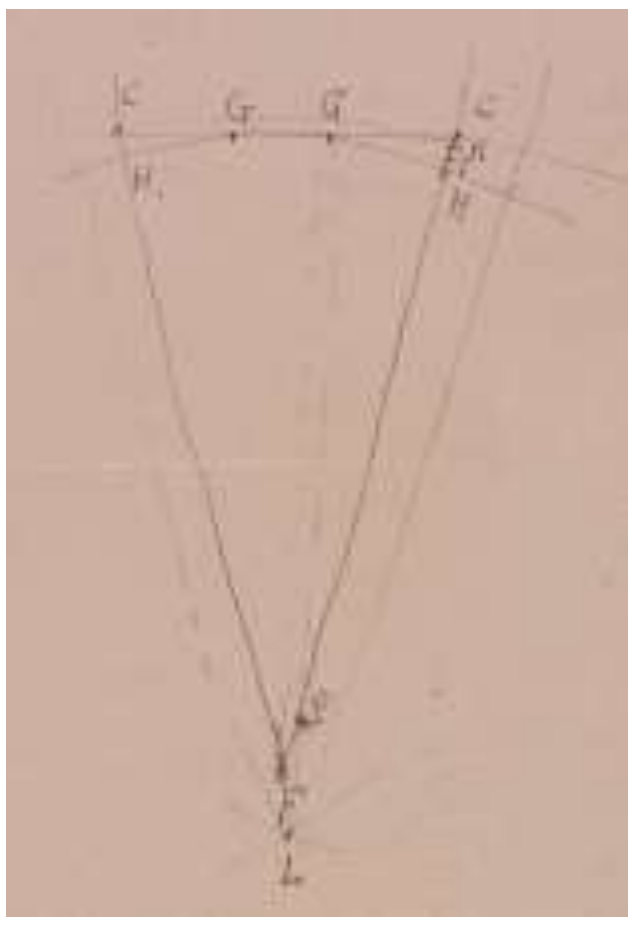

*[appunti allegati]*
L'a. ammette le due eguaglianze LG=PC e HC=PL. Ne viene l'assurdo FP+PG=GF

GF=FH
FP=HI
FP+PG=FH
PG=LG–LP=KF–HK–HI=FH–HI

*[sul retro di questo foglio probabilmente appunti per un articolo]*

**14**
**Cosimo Canovetti** a Gian Antonio Maggi
[busta e carta intestate: Ing. C. Canovetti - CAV. LÉGION D'HONNEUR - Ing. Consulente - MILANO (17) - Via Sebeto, 3]

Milano, 4 Nov 1926
Preg Sig Prof Maggi Preside Facoltà Scienze
Università Milano

Fino dal 1921 con conferenze e pubblicazioni inviate a tutte le Accademie mondiali, ho dimostrato che le formule di EINSTEIN nulla dimostrano e che sono scientemente alterate. Tale dimostrazione comporta solo la conoscenza dell'algebra elementare. Ma debbo confessarLe che nessuna accademia, nemmen quelle che mi hanno assegnato premi, non si pronuncia.

Recentemente al Congresso di Bologna (Scienze) il Prof Gianfranceschi lo attaccò con ragione.[53]

Il Lallemand Vice presidente Académie des Sciences ha pubblicato un attacco Astronomico.[54]

Io fin dal 1921 chiesi al Senatore Mangiagalli di tenere una conferenza all'Università. Gli ripetei la domanda recentemente. Per chiederLe di voler appoggiare la mia domanda sono a pregarLa di volermi fissare un convegno ove meglio Le piacerà per farle vedere le mie pubblicazioni con grafici che dimostrano che le formule di Einstein non solo non diminuiscano il tempo di un sistema in moto, ma che PER ERRORE lo accrescono. Tale verifica è alla portata di tutti.

In attesa con profonda stima di Lei

OBBMO
C. Canovetti
con numerosi
premi scientifici.

---

[53] Si tratta probabilmente della XV Riunione della SIPS, tenutasi a Bologna nel 1926 (30/10-5/11). G. Gianfranceschi in quell'occasione, tenne due conferenze: "Sulla legge di distribuzione dell'energia dello spettro del corpo nero" (*Atti*, p. 551) e "Sulle attuali teorie della fisica" (*Atti*, pp. 32-43). Presidente era C. Somigliana.
[54] Potrebbe trattarsi di: C. Lallemand, "La théorie de la relativité et les expériences du prof. D.C. Miller", 1926.

**15**

Gian Antonio Maggi a Cosimo Canovetti

Milano 8 Dicembre 1926

Pregiatissimo Signor Ingegnere,

Se Ella desidera parlarmi, potrà cercare di me, per essere sicuro di trovarmi, alla Università, il Mercoledì alle 10 e il Martedì, Giovedì e Sabbato alle 11. Volentieri io farò la Sua conoscenza personale, ma per la domanda a cui Ella accenna, io non potrò che attenermi, per ragione di massima, al consiglio ripetutamente espresso in merito dalla Facoltà.

Aggradisca intanto i miei distinti saluti, coi quali La prego di credermi

Dev.mo Suo
G.A. Maggi.

Ill.mo Sig.r
Ing. Cav. C. Canovetti
Via Sebeto, 3 Città (17).

**16**[55]

Gian Antonio Maggi a Cosimo Canovetti
[busta indirizzata a: Ill.mo Sig.r Ing. Cav. C. Canovetti
Via Sebeto 3 - Città]

Milano 19 Gennaio 1927

Pregiatissimo Signor Ingegnere,

Sull'argomento de' Suoi scritti, che mi ha favorito a suo tempo, e di cui La prego scusarmi di non avermi potuto occupare prima d'ora, le equazioni

(1) $x-ct=0, \qquad x'-ct'=0$

non sono verificate che nel caso particolare della propagazione di un segnale luminoso, con velocità c, lungo l'asse delle x, in senso positivo, mentre le x,t devono essere considerate, in generale, come variabili indipendenti, atte a individuare un punto qualsivoglia dell'asse delle x, ad un istante pur qualsivoglia: punto e istante, individuati, col secondo riferimento, da x',t'. Ma, per l'ipotesi fondamentale della Relatività, da Lei pure ricordata, se si verifica l'una di esse, deve allora verificarsi anche l'altra, e questa necessaria simultaneità (non necessaria separata verificazione) si traduce nell'equazione

(2) $x'-ct'=\lambda(x-ct),$

con $\lambda$ indipendente dalle coordinate, dove x,t serbano il loro significato di variabili indipendenti: la quale serve a stabilire una relazione tra x,t e x',t'. Per la stessa ipotesi le

(3) $x+ct=0, \qquad x'+ct'=0$

sono tali che, come le precedenti, se si verifica l'una di esse (segnale che si propaga lungo l'asse delle x, alla velocità c, in senso negativo), deve allora verificarsi anche l'altra, e questa necessaria simultaneità si traduce nella

[55] Questa minuta è conservata, con qualche variazione, in quattro copie, ma è dubbio che sia stata effettivamente spedita.

(4) $$x'+ct'=\mu(x+ct),$$
con $\mu$ indipendente dalle coordinate, che stabilisce un'altra relazione tra x,t e x',t'.
Le (2), (4) traducono in tal modo la suddetta ipotesi. Non così le

$$x-ct=0, \qquad x'-ct'=0, \qquad x'=x-vt,$$

che manifestamente non sono ad esse equivalenti. Per modo che, se da queste ultime equazioni si ricavano conseguenze in contraddizione coi risultati dell'Einstein, non sono però esse stesse d'accordo coi principii della Relatività.

Dalle (2), (4) poi si ricavano le note formule

(5) $$x'=\frac{x-vt}{\sqrt{1-\left(\frac{v}{c}\right)^2}}, \qquad t'=\frac{t-\frac{v}{c^2}x}{\sqrt{1-\left(\frac{v}{c}\right)^2}}$$

La formazione delle (2), (4) e la deduzione da esse delle (5) sono dimostrate dall'Einstein nell'Appendice al suo opuscolo "Sulla Relatività Speciale e Generale - Volgarizzazione - Traduzione dell'Ing. Calisse" (Bologna, Zanichelli), e a quei ragionamenti io non trovo da muovere alcuna obbiezione.
Né ho mai trovato da muovere obbiezioni, al procedimento più antico e più comune fondato sulla

$$x'^2+y'^2+z'^2-c^2t'^2=x^2+y^2+z^2-c^2t^2,$$

che si stabilisce con ragionamento simile a quello che serve per stabilire (2) e (4), o sulla

$$x'^2-c^2t'^2=x^2-c^2t^2,$$

che se ne deduce, facendosi y'=y, z'=z (non necessariamente y=0, z=0 etc.), conformemente alla nota scelta (in nostro arbitrio) degli assi coordinati.

Questa la mia sincera opinione, che, poiché Ella mi ha usato la cortesia di chiedermi, credo doverLe manifestare.

Mi è grato presentarLe insieme i miei distinti saluti, e pregarLa di credermi

Dev.mo Suo
G.A.Maggi.

(*). *[non c'è il richiamo corrispondente all'interno del testo]* Resta naturalmente a parte la discussione del valore delle ipotesi, sotto l'aspetto della verificazione sperimentale - obbiezioni di fisici, di astronomi - e - obbiezioni di filosofi - sotto l'aspetto della dottrina della conoscenza: obbiezioni alle quali ho io pure accennato nella mia conferenza al Congresso di Mathesis nell'Ottobre del 1925, a Milano (pubblicata nel "Periodico di Matematica", 1926, p. 1).

**17**
**Guglielmo Capitelli** a Gian Antonio Maggi
[Biglietto da visita: Conte Guglielmo Capitelli - Prefetto di Messina.]

Messina - 4-5-93[56]

Egregio professore non rida, e possibilmente mi dia una risposta.

Pare a Lei siovi qualche legge matematica, che stabilisca la verità di questo assunto "essere negli eserciti più numerosi maggiori le probabilità di più rapide promozioni nei varii gradi, massima negli alti".

Desidero una risposta in iscritto essendone stato richiesto da alcuni amici. Mi scusi e mi creda

Aff Suo
G. Capitelli

**18**
Gian Antonio Maggi a Guglielmo Capitelli

Illustrissimo signor Conte,

Non mi consta che siano stati fatti studii sulla questione di cui Ella mi scrive, almeno non mi risulta dai trattati sul Calcolo delle Probabilità ch'io conosco, alcuni dei quali sono scritti da Militari.

Posta la questione a priori, bisogna innanzi tutto stabilire come eserciti diversamente numerosi s'intendono composti: e fin qui ammetteremo che il numero dei posti di un certo grado sia proporzionale al numero totale dei soldati. Ma bisogna definire precisamente la legge con cui s'intendono fatte le promozioni.

Ammesso che le promozioni si facciano per pura anzianità, la probabilità che ha ciascun individuo di conseguire in un certo tempo un certo grado sembra indipendente dal numero totale; perché essa sarà il prodotto delle probabilità che, in quel tempo, si renda vacante un posto di quel grado, per la probabilità che, rendendosi vacante quel posto, ha l'individuo di avere l'anzianità necessaria per conseguirlo. Ora la prima probabilità è direttamente proporzionale al numero dei posti di quel grado, e quindi al numero totale, mentre la seconda è inversamente proporzionale a questo numero.

Se poi si ammette la promozione a scelta, ciò che, almeno in parte, si verifica in pratica, le circostanze determinatrici sono troppe e troppo complicate perché si possa fare a priori un'ipotesi che traduca con sufficiente approssimazione la realtà. Ma pare però di poter fare le riflessioni seguenti:

Ammesso che un individuo abbia certi titoli segnalati pei quali può aspirare ad un dato grado, la probabilità che in un certo tempo si renda vacante un posto di quel grado sarà, per la nostra supposizione, proporzionale al numero totale. Invece sembra ipotesi conforme al fatto che non sia proporzionale a questo numero quello che gl'individui che avranno titoli eguali, ma assai minore di quello che vorrebbe questa proporzionalità. Per esempio, se in un esercito di 100.000 uomini vi è un colonnello insigne per opere militari, ciò non vuol dire che ve ne saranno 10 egualmente insigni in un esercito di un milione, ma presumibilmente ve ne saranno meno. Sembra quindi che il supposto individuo in un esercito più numeroso potrà conseguire quel dato grado, in un certo tempo, più facilmente. E poiché si richiederanno titoli tanto più eccezionali quanto più alto è il grado potrà realmente darsi che il fatto si verifichi, coll'elevarsi del grado, con crescenti sensibilità.

[56] Probabilmente si tratta del 1892 sia perché Maggi così scrive nella risposta, sia perché, nel maggio 1893, Capitelli si trovava a Firenze.

Non so se questa risposta, che le do senza pensarci di più per non farla di più attendere, Le potrà bastare. Mi sarà sempre grato, s'Ella crede, d'aggiungere tutto quello che posso. In massima poi (senza far torto alla diva Urania) bisogna accogliere il "matematicamente dimostrato" colle debite riserve: perché spesso, e specialmente quando si tratta di questioni pratiche complicate, non è che la logica conseguenza d'ipotesi che traducono solo parzialmente la realtà, e non possono avere un valore pratico diverso da esse.

Aggradisca i miei ossequii, e mi creda sempre

Devotissimo aff$^{mo}$
Gian Antonio Maggi.

Di casa 7 maggio 1892.

**19**
Gian Antonio Maggi a **Michele Casamassima**

Messina 27 Febbrajo 95

Egregio Professor Casamassima,

Mi sono ricordato senza dubbio del suo scritto; soltanto che ho dovuto aspettare le vacanze per trovare il tempo d'occuparmene. L'ho esaminato; e poiché Ella mi ha chiesto il mio parere, Le debbo dire che mi rincresce di non trovarlo sufficiente, secondo i miei criterii, per una dissertazione dottorale. Non pretendo molto, e ammetto anche una semplice compilazione; purché la maturità di studii del candidato emerga per qualche pregio dal lato della misura, dell'ordine, del metodo e della forma oltre di che da un certo studio della letteratura relativa all'argomento; dal qual complesso di circostanze riesca dimostrato che il candidato ha saputo approfondire una teoria, e ricavare dalle proprie idee una buona esposizione. E ciò non so trovare nel Suo scritto.

Non vedrei che cosa Ella possa aver aggiunto alla ricerca delle note equazioni differenziali della propagazione del calore, che ne forma la massima parte; bensì alcuni punti richiederebbero qualche commento; e talvolta pare ch'Ella non abbia perfettamente compreso. Per esempio, insiste sopra un terzo principio della teoria di Fourier, e mostra, di non avvertire che, com'è da Lei posto, riesce una diretta conseguenza del secondo, il quale alla sua volta ne scaturisce coll'ammissione di termini sostanzialmente trascurabili, mentre Ella parla di questi a proposito del secondo principio, e come compariscono non si vede.

Il problema speciale, che considera infine, è posto in forma assai imperfetta; talché, chi non lo sappia, difficilmente credo io potrebbe farsene un concetto preciso. Vi si parla promiscuamente d'un punto riscaldato, d'un'areola infinitamente piccola e d'un area qualunque, in modo da non potersi rilevare dallo scritto la parte che si attribuisce a ciascun ipotesi. La formula di risoluzione, è appunto ciò che, concentrando parecchie nozioni, fra le più importanti in simili problemi, meriterebbe d'essere più largamente discusso e illustrato;

e, invece Ella si limita a rimandare al Mathieu,[57] e, al poco che ne dice, c'è da domandarsi come la u, che è una funzione di t e delle coordinate del punto qualunque della lamina, diventa per t=0 una funzione delle variabili d'integrazione del doppio integrale definito.

Per quanto alla parte tratta dalla Nota del sig Jannetaz,[58] - che è un'introduzione ad un lavoro sperimentale, per se stessa d'assai limitato interesse - va osservato che questo autore ammette a priori la costanza della temperatura sopra ogni circonferenza di circolo avente per centro il punto riscaldato, e riproduce così caso più semplice del noto equilibrio della temperatura in un cilindro circolare. L'integrazione per serie è anche più nota di quella per integrali definiti di cui si vale l'autore; solo sarebbe desiderabile che fosse avvertito come comparisce la funzione besseliana $J_o$, che in simili questioni ha tanta importanza. È poi evidente che l'eguaglianza stabilita fra le due serie è soddisfatta da $hx^2=h'x'^2$, e l'autore dice d'esser ricorso al sig. Picard per la dimostrazione che questa è la sola soluzione possibile.

S'Ella crede di restare in questo argomento, Le suggerisco di consultare Riemann - Partielle Differentialgleichungen, e Kirchhoff - Theorie der Wärme; e procurar di trattare qualche problema, che, per esempio, si presti, all'uso di coordinate curvilinee e non sia esaurito.

Riceva i miei saluti, e mi creda

Dev^mo^Suo
G.A. Maggi.

## 20

**Gino Cassinis** a Gian Antonio Maggi
[su carta intestata: R. Istituto Superiore di Ingegneria - (R. Politecnico) - Milano
Istituto di Topografia e Geodesia]

Milano 28 gennaio 1935 XIII

Ch/mo Professore,

Le restituisco il dattiloscritto dell'ing. Paroli "Sul calcolo numerico delle coordinate geodetiche rettangolari", che ho esaminato in vista della eventuale pubblicazione nei Rendiconti dell'Istituto Lombardo.

Supponendo di dover calcolare le coordinate rettangolari, riferite ad un'unica origine, di un notevole numero di vertici situati in una regione non troppo estesa, il Pàroli si propone di ridurre la mole dei computi, introducendo un metodo di interpolazione.

A tale scopo, egli suddivide la zona in un certo numero di trapezi ellissoidici, che in un caso concreto assume dell'ampiezza di 5' in latitudine e in longitudine, e calcola le coordinate geodetiche rettangolari X, Y in alcuni dei loro vertici con le consuete formule, negli altri con formule di interpolazione che egli stesso ricava e che risultano di applicazione abbastanza semplice tabellandone i valori.

Una volta in possesso delle coordinate di tutti i vertici dei trapezi di 5', questi possono suddividersi in trapezi di minore estensione, p. es. corrispondenti ad ampiezze di 1' in latitudine e longitudine, ed è facile costruire delle tabellette ausiliarie che, in unione a quella principale consentono di ottenere le coordinate rettangolari di qualsiasi punto della zona, mediante semplici operazioni di addizione e sottrazione, con errore al più di qualche cm.

[57] Si veda la nota 26.
[58] Potrebbe trattarsi del chimico francese Édouard e la Nota in questione potrebbe essere *Sur la propagation de la chaleur dans les corps cristallisés*, Gauthier-Villars, Paris, 1873.

Il procedimento è conveniente solo quando il numero dei vertici della zona, di cui sian note le coordinate geografiche, e dei quali si desiderino quelle ortogonali, è piuttosto ingente.

Il lavoro del Pàroli, anche se non sempre il suo metodo possa essere di impiego praticamente opportuno, presenta interesse, e quindi ritengo debba essere presentato all'Istituto e pubblicato nei Rendiconti.

Però occorre pregare l'A. di rivedere accuratamente il manoscritto, pieno di errori e di dimenticanze (ho segnato qualche cosa per ricordo), e, in certi punti non troppo chiaro, sia come scrittura (perché eccessivamente fitto), sia anche come forma.

Bisognerebbe anche pregarlo di sostituire nelle formule arc 1'' in luogo di sen 1''. È vero che i valori numerici di queste due quantità differiscono di poco e la differenza è trascurabile nel calcolo delle formule applicate dal Pàroli; e, purtroppo, è anche vero che in molti trattati si trova scritto sen 1'' al posto di arc 1''; ma l'errore concettuale deve essere evitato qui come negli altri casi analoghi.

Voglia perdonarmi, Illustre Professore, se ho tardato, ma Lei conosce le molte occupazioni che non sempre mi consentono di fare ciò che vorrei.

Con molti devoti saluti

G. Cassinis

Ill/mo sig. prof. G.A. Maggi
Corso Plebisciti 3
MILANO

**21**[59]

Gian Antonio Maggi a **Oscar Chisini**

[busta indirizzata a: Ch.mo Prof. Oscar Chisini - S.P.M. con la dicitura "Non mandata"; busta e carta intestate: R. Università di Milano - Istituto Matematico - Via C. Saldini, 50 (Città degli studi) - Il Direttore]

Letta ma non consegnata

Milano 21 Aprile 1930

Caro collega,

Le rendo la lettera circolare del Sig.r Brodovitzky.
Non reputo necessario ricorrere agli iperspazii, per trovare la doppia relazione (pag. 5), che forma la capitale obbiezione del Brodovitzky alla teoria della Relatività, bastante valersi del piano e di una superficie cilindrica rientrante. Indicando con L e L' le lunghezze del circuito, rappresentato da una sezione normale della superficie cilindrica, per gli osservatori appartenenti ai due riferimenti, in moto (inerziale), l'uno relativamente all'altro, da

$$x' = \beta(x - vt), \qquad x' + L' = \beta(x + L - vt)$$

si deduce

(1) $$L' = \beta L,$$

e da

$$x = \beta(x' + vt'), \qquad x + L = \beta(x' + L' + vt')$$

si deduce

(2) $$L = \beta L'.$$

Da (1) e (2) il Brodovitzky deduce ß=1, invece del supposto $\beta = \left(1 - \left(\frac{v}{c}\right)^2\right)^{-\frac{1}{2}}$.

Ora, le (1), (2) non sono che le notissime e incontrastate relazioni tra i rapporti delle lunghezze di due segmenti, misurati, una volta, col tempo t, appartenente al riferimento S, e, un'altra volta, col tempo t', appartenente al riferimento S'.

Con cordiali saluti

Affmo Suo
Gian Antonio Maggi.

---

[59] Questa versione è conservata insieme a una minuta.

**22**
**Umberto Cisotti** a Gian Antonio Maggi
[busta indirizzata a: Chiar.mo Sig.re Prof.re Gian Antonio Maggi
della R. Università - Via Risorgimento, 1 - Pisa;
busta e carta intestate: R. Istituto tecnico superiore - Milano]

7 Dicembre 1917

Chiarissimo Professore,

Le sono vivamente grato della lettera cortese che Ella mi ha inviato in seguito all'omaggio, che ebbi il piacere di farle, delle mie Lezioni di Analisi matematica. La ringrazio altrettanto vivamente di avermi manifestata la sua impressione generale, arricchita di riflessioni su punti particolari.

Dato il tipo di scuola dove l'analisi matematica deve svolgersi era un compito non troppo facile lo svolgerlo degnamente. Esso rappresenta un tentativo e come tale desidero venga giudicato. La mia personale esperienza e il benevolo giudizio di altri eminenti colleghi (basti citare il Levi-Civita e il Somigliana) mi confortano e mi inducono a ritenere di non avere speso male il tempo (non breve) destinato allo scopo.

Il titolo stesso "Lezioni…" e non "Trattato" mette in rilievo che l'A. non ha avuto la pretesa di avere toccato tutte, anche fra le più importanti, le innumerevoli questioni di analisi algebrica e infinitesimale; mi sono imposto dei limiti precisi che credo sarebbe pericoloso oltrepassare senza turbare l'equilibrio delle singole parti e dell'assieme.

Perciò ben lieto di accogliere osservazioni e consigli, quelli tra questi ultimi che riguardano aggiunte da farsi mi lasciano alquanto trepidante. Tuttavia intanto ne prendo atto e rimando la questione a quando, se mi sarà concesso, dovrò fare una ristampa del libro.

Ho letto con molto interesse le due Note che Ella ha avuto la cortesia di inviarmi. La loro lettura mi ha fatto riesumare alcuni appunti manoscritti che feci in riguardo alle Note del Somigliana sulle discontinuità delle derivate prime e seconde delle funzioni potenziali di semplice e di doppio strato. L'applicazione delle formole di Somigliana mi indusse a mettere in rilievo le discontinuità in superficie delle derivate seconde e terze di una funzione potenziale di volume.

A questo si perviene con tutta facilità ricorrendo alla circostanza ben nota che le derivate prime di detta funzione potenziale si possono esprimere ciascuna quale somma di una funzione potenziale di superficie e di una funzione potenziale di volume. La

$$\varphi = \int_S \frac{k dS}{r}$$

è la funzione potenziale in discorso, e si immagina che S sia limitata da una superficie <u>regolare</u> σ (o da un insieme di tali superficie), seguendo le Sue notazioni si trova

$$D \frac{\delta^2 \varphi}{\delta z^2} = -4\pi k,$$

mentre nell'attraversare σ tutte le altre derivate seconde si mantengono continue. Per le derivate terze si ottiene

$$D\frac{\delta^3\varphi}{\delta x\delta z^2}=-4k\frac{\delta k}{\delta s_u},\qquad D\frac{\delta^3\varphi}{\delta y\delta z^2}=-4k\frac{\delta k}{\delta s_v},$$

$$D\frac{\delta^3\varphi}{\delta x^2\delta z}=-\frac{4\pi k}{r_1},\qquad D\frac{\delta^3\varphi}{\delta y^2\delta z}=-\frac{4\pi k}{r_2},$$

$$D\frac{\delta^3\varphi}{\delta z^3}=4\pi k\left(\frac{1}{r_1}+\frac{1}{r_2}\right)-4\pi\frac{\delta k}{\delta z};$$

le rimanenti derivate terze rimangono continue nell'attraversare la superficie σ nel punto considerato.

Poiché penso che questi risultati - per quanto facilmente giustificabili noti i precedenti - possano interessarla, così mi sono permesso di sottoporli alla sua attenzione. Essi potrebbero fare seguito alle formole che Ella mette in rilievo nel #:2 della Nota lincea.

Già da tempo ho promesso al Loria di recensire il Suo ultimo volume sulla Dinamica dei sistemi, incarico che ben volentieri mi sono assunto non appena avessi ultimato la stampa delle mie lezioni. Ora debbo mantenere l'impegno assunto: avrò così occasione di intrattenermi a lungo spiritualmente coll'A. la cui opera mi ha lasciato un'impressione così gradita.

Con i più distinti e cordiali ossequi mi creda

Suo dev.mo
M Cisotti.

P.S. Distinto saluto al Prof Pizzetti.

**23**
Gian Antonio Maggi a **Giuseppe Colombo**[60]

Al Prof. G. Colombo.
Pisa 6 Novembre1908

Illustre e Gentilissimo Professore,

Che mi sia o no riserbato d'insegnare, una volta o l'altra, in una Scuola d'Applicazione (ciò che, secondo certe nascenti aspirazioni, mi potrà accadere anche fissando stabile dimora a Pisa), mi sta pur sempre a cuore la questione dei rapporti della Meccanica Razionale colla Tecnica, che io reputo vitale per l'avvenire di codesto insegnamento. Perciò ho veduto con molto interesse la lettura fatta a Firenze dall'Ing Barzanò:[61] spiacente anche per non aver potuto partecipare a tale discussione che una gita a Milano mi abbia impedito d'intervenire, quantunque iscritto, al Congresso.
L'Ing. Barzanò reclama un radicale rimaneggiamento dell'insegnamento della Meccanica Razionale, ed espone, in proposito, varii suoi desiderii. Io, con compiacenza, rilevo che a codesti varii desiderii io procuro di soddisfare nel corso, che, oramai da tredici anni, impartisco in questa Università. Io soglio premettere il problema reale, concernente projettili, ruote etc. al problema idealizzato o teorico: insisto più particolarmente sulla meccanica di sistemi, perché più agevole da essi che dal punto è il passaggio ai mobili del mondo reale: sviluppa i principii della cinematica e della dinamica dei sistemi deformabili, ponendo così le basi della teoria dei sistemi elastici: infine riduco allo stretto necessario l'uso degli assi coordinati, attenendomi di preferenza al ragionamento sintetico e al metodo vettoriale. Con tutto questo, in confronto dello spirito pratico invocato dall'Ing Barzanò, non mi dissimulo che il mio corso rimane ciò che si chiama un corso teorico, cioè tale che, nel suo complesso, non si presta ad essere dalla Pratica utilizzato direttamente. Ma io mi chiedo, a questo proposito, se, volendo conservare l'integrità del corso di Meccanica Razionale (poiché una riforma potrebbe consistere nel decomporlo in capitoli, e premetterli ai varii corsi tecnici), questo insegnamento possa mai altro prefiggersi che lo svolgimento dei metodi generali per tradurre in espressioni matematiche i problemi presentati dalla Fisica e dalla Pratica. Col quale scopo, mi sembra destinato o rimanere un insegnamento di cui gli allievi debbano portar pazienza ad aspettare i frutti attraverso al tramite d'insegnamenti tecnici speciali ad esso efficacemente coordinati.

La prego di scusare, colla Sua consueta benevolenza, la libertà che mi prendo di metterla a parte di queste mie impressioni, ché io non saprei davvero a chi meglio rivolgermi che a Lei, che tanto posto tiene nel nostro Insegnamento. All'Ing Barzanò non mi permetto di scrivere, non avendo il piacere d'essere con lui in relazione: e dal renderle di pubblica ragione mi trattiene il dubbio di passar per prendere la parola per fatto personale.

La prego ritenere di aggradire i miei più vivi e più distinti rispetti, e di credermi sempre

Dev$^{mo}$ Suo
Gian Antonio Maggi.

---

[60] La corrispondenza con Giuseppe Colombo è contenuta nel fascicolo denominato da Maggi: *Sull'articolo dell'ing. Barzanò "Alcune considerazioni pel riordinamento dell'insegnamento tecnico superiore" (v.* <u>Industria</u> *1 Novembre 1908).*
[61] Si vedano la lettera di Camillo Guidi (#84) e all'Ing. Barzanò (#3).

**24**
Giuseppe Colombo a Gian Antonio Maggi
[busta indirizzata a: Prof. G.A. Maggi - R. Università Pisa;
busta e carta intestate: R. Istituto tecnico superiore - Milano - Direzione]

Milano 14 Nov 1908

Egregio Collega

Mi rincresce che l'ing Barzanò non mi abbia prevenuto che intendeva di trattare quell'argomento al Congresso, che poi fu tema, come lessi, di una discussione. A questa avrei preso parte io pure, e avrei cercato di chiarir meglio le necessità di una scuola d'ingegneria e il modo in cui io le sento e cerco di soddisfarle.

Non ho avuto ancora l'occasione di vedere il mio antico allievo ed amico Barzanò, per discutere con lui sulla questione. Io non credo (e non suppongo nemmeno che egli lo creda) che il prof. di Meccanica razionale debba trattare le questioni pratiche ed esclusivamente quelle che possono aver rapporti colla pratica. Io vorrei invece che il Corso non perdesse il suo carattere teorico; le applicazioni le faranno gli altri professori. Solamente il prof di M.R. non dovrebbe trattare la materia dal solo punto di vista matematico, come si faceva una volta e come lo faceva il compianto prof Bardelli; ma additare il senso fisico e meccanico delle questioni, soprattutto in dinamica, dando così interesse e vita agli argomenti, trattati, nel senso, per es., della Meccanica del Mach.[62] Dovrebbe anche abbandonare certi argomenti, che non interessino ne' direttamente, ne' indirettamente l'ingegneria, e tratti invece più largamente quelli che interessano, con nuovi metodi e nuovi punti di vista, l'elettrotecnica e la resistenza delle costruzioni. Ciò che Ella dice sullo svolgimento da Lei dato al suo corso, mi pare, quindi, corrispondere a questo mio modo di vedere.

Io non Le ho parlato più di quanto Ella mi disse un giorno a proposito della Cattedra di M.R. del Politecnico. Il fatto è che il Consiglio dei miei professori ha trattato più volte la questione, senza venire a una conclusione per quest'anno, e ciò per il contrasto di idee fra i fautori di un Corso di tendenze teoriche e quelli di un corso di tendenza spiccatamente pratica. Da questo dibattito è uscito il partito, unanimemente accettato, anche da me (per quanto dapprima vi fossi contrario) di aprire il concorso per un posto di straordinario. Avremo quindi probabilmente un giovane, cercheremo di avviarlo a trattare la materia nel modo che ci sembra migliore per la nostra Scuola, e se corrisponderà al nostro ideale, potremo domandarne la promozione a ordinario. Era la sola soluzione che metteva tutti d'accordo, ed io l'ho accettata senz'altro.

Aggradisca i miei cordiali saluti

Suo G.Colombo

[62] Potrebbe trattarsi di: E. Mach, *Die Mechanik in ihrer Entwicklung*, Leipzig, 1883.

**25**
Gian Antonio Maggi a Giuseppe Colombo

Pisa 20 Novembre 1908
Non mandata, e mandato invece il cartoncino[63]

Gentilissimo Professore,
(Sen. Colombo).

La ringrazio vivamente della pregiatissima Sua. Io ben sapevo quale doveva essere il Suo concetto dell’insegnamento della Meccanica Razionale. Solo mi sono permesso di far presente a persona tanto autorevole, e a me tanto benevola, come l’Ing Barzanò, del quale conosco per fama il merito, ne invocasse “un radicale rimaneggiamento”, di cui, stando ai suoi stessi desiderii, sembra che il mio corso non abbia veramente bisogno.

La ringrazio ancora della notizia che mi favorisce sulla decisione di codesto Consiglio dei Professori. Intendo perfettamente come, quando sorga disparità d’opinioni sopra la proposta di persona, convenga senz’altro abbandonarla. Giova alle stesse persone in questione, che arrischierebbero di trovarsi poi a disagio fra la persuasione di non accontentare gli oppositori e la preoccupazione di non render abbastanza ragione all’appoggio dei fautori.

Di nuovo La prego di aggradire, coi miei ringraziamenti, i più distinti e cordiali rispetti, e di credermi sempre

Dev$^{mo}$ Suo
Gian Antonio Maggi.

**26**
Gian Antonio Maggi a Giuseppe Colombo
[busta indirizzata a : Illustre Prof. Giuseppe Colombo - Senato del Regno - Direttore del R. Istituto Tecnico Superiore di Milano]

Pisa 20 Novembre 1908

Gentilissimo Professore,

Le sono vivamente obbligato per la Sua risposta, e per la notizia, che mi ha, allo stesso tempo, favorito, sulla decisione di codesto Consiglio dei Professori. Riconosco perfettamente che tale doveva essere, colla disparità, a cui Ella accenna, di opinioni.

Di nuovo La prego di aggradire, coi miei ringraziamenti, i più distinti e cordiali rispetti, e di credermi sempre

Dev$^{mo}$ Suo
Gian Antonio Maggi.

[63] Si riferisce alla lettera seguente.

**27**

Gian Antonio Maggi a **Gustavo Colonnetti**[64]

Viareggio 23 Luglio 1916

Egregio e caro Collega,

Rammentando il suo desiderio di conoscere la risposta di Somigliana sulla questione delle distorsioni elastiche di cui Le feci parola, la mia osservazione, a Somigliana, sembra perfettamente giusta.[65] Egli conviene che conseguenza dei suoi recenti risultati è che potrebbero esistere distorsioni determinate da spostamenti rigidi aventi tensioni con derivate discontinue: nei quali casi le distorsioni di Volterra non esisterebbero. Aggiunge che sarebbe interessante verificare con qualche esempio la possibilità di queste nuove discordanze con ciò che ci suggerisce la più semplice intuizione. Ma a me l'intuizione all'infuori del caso del *[*toro*]* per verità, su codesto punto dice poco. E credo piuttosto che occorrerebbe *[*dilagare*]* quale grado di generalità possiedono le distorsioni di Volterra definite da condizioni che appariscono sovrabbondanti. È una questione aperta, ch'io però io non ho intenzione di trattare, riputandomi soddisfatto che ne risulti provato come le mie obbiezioni allo spostamento polidromo non fossero vana logomachia, ma richiamo alla diritta via.

Al Prof Colonnetti.
Via Arcivescovado 6 Torino

Viareggio - 23 Luglio 1916.

**28**

Gustavo Colonnetti a Gian Antonio Maggi
[cartolina postale indirizzata: Al Ch^mo Sig Professore G.A. Maggi
Via Zanardelli 54 - Viareggio]

Torino 26 luglio 16

Egregio Professore,

Le son grato delle informazioni che Ella mi manda e che mi interessano vivamente. Per parte mia mi sono convinto che per bisogno delle applicazioni le distorsioni possano essere benissimo definite dalla discontinuità delle compon. di spostam. ottenuto mediante moto rigido relativo delle due faccie *[sic!]* del taglio, anche senza assicurarsi della continuità della derivata delle comp. di deformazione. Ma di questo spero di poterle presto parlare personalmente.

Intanto Le prego dei miei più vivi e rispettosi ossequii alla Sua Famiglia, anche con rinnovati ringraziamenti mi dico di Lei sempre dev^mo.

G Colonnetti

---

[64] La grafia della lettera non sembra di Maggi, il quale ha però scritto il destinatario e la data.
[65] Si veda la corrispondenza di Maggi con Somigliana relativa a questo periodo.

**29**

Gustavo Colonnetti a Gian Antonio Maggi
[busta indirizzata a: Ch$^{\underline{mo}}$ Prof. G. A. Maggi della R. Università di Pisa;
busta e carta intestate: Società anonima italiana per la fabbricazione dei proiettili]

Torino 4 agosto 1918

Egregio Professore,

Faccio assegnamento sulla sua pazienza e Le accludo alcune pagine di chiacchierata sulle quali desidererei il suo avviso. Si tratta di quei certi stati di coazione elastica in assenza di forze esterne di cui mi ero occupato l'anno scorso ed a proposito dei quali Ella mi aveva presentata una Nota ai Lincei il 19 giugno 1917.

Studiando alcuni casi particolari assai interessanti perché connessi colla teoria della resistenza delle artiglierie sono giunto ad un teorema di questo genere: la densità media di un solido elastico in equilibrio in assenza di forze esterne (e ciònonpertanto dotato di tensioni interne) è eguale alla densità propria del materiale allo stato naturale non deformato. O in altre parole: lo stato di tensione interna non muta il volume totale.

Questo teorema non pare fosse noto se non in casi particolarissimi (tubi cerchiati); non è stato rilevato, che io sappia, nel caso delle distorsioni. Tuttavia a me risulterebbe vero in un gran numero di casi; e mi sembra di essere riuscito a dimostrarlo in modo generale per qualunque solido omogeneo imitando il ragionamento fatto dal Betti per studiare la dilatazione cubica dei solidi elastici pesanti.

Vuol esaminare questa dimostrazione e dirmi se non ci trova nulla a ridire? Io vado a Roma dove Ella potrebbe scrivermi all'Hotel Santa Chiara. Se il Sen Volterra non mi tratterrà colà troppo a lungo, mi fermerò al ritorno a Pisa e verrò personalmente a chiederle scusa della libertà che mi prendo del disturbo che Le arreco.

Ad ogni modo accetti fin d'ora le mie scuse; e mi ritenga giustificato pensando che non sono abbastanza matematico per sentirmi sicuro del ragionamento qui esposto e che d'altra parte sarebbe vero peccato servirsi della proprietà in discorso nei singoli casi speciali (in cui è indubbiamente vera e molto utile) senza stabilirla prima colla maggior possibile generalità.

La ringrazio della sua ultima cartolina; La prego dei miei ossequi alla Sua signora e Signorine, mentre con anticipati ringraziamenti mi dico

di Lei sempre dev$^{\underline{mo}}$
G. Colonnetti

Consideriamo un solido elastico omogeneo (vale a dire: le cui costanti elastiche siano indipendenti dalle coordinate), libero da vincoli, il quale nel suo stato naturale - cioè nello stato di equilibrio stabile che spontaneamente assume quando non è soggetto ad alcuna forza esterna - si presenti in uno di quegli stati di coazione elastica di cui mi sono occupato nelle due note precedenti (°)

(°) Confr. ……………………...……………………………………..

Ed immaginiamo attraverso ad un tale solido praticati tutti quei tagli che sono necessarii perché i singoli elementi materiali che lo compongono - liberati dai vincoli che inizialmente essi si creavano mutuamente - possano tutti assumere il loro stato non deformato.

Nessuna restrizione vogliamo fare relativamente al numero o alla disposizione di quei tagli, i quali, nei casi più generali, potranno anche moltiplicarsi all'infinito fino a sconnettere, ove fosse necessario, il solido dato nei singoli elementi.

Né vogliamo che resti in alcun modo limitato il grado di libertà che coi tagli stessi si viene a conferire alle parti del solido che essi separano; perciò supporremo esplicitamente che, una volta praticato un taglio, i punti già affacciati, da una parte e dall'altra di esso, possano presentare spostamenti relativi affatto qualunque, tanto nel senso di scostare le due faccie *[sic!]* del taglio, quanto nel senso di dar luogo a sovrapposizioni di materia.

Ciò premesso non v'è dubbio che dallo stato non deformato a cui il solido è così ridotto, si potrà sempre ritornare al primitivo stato naturale non soltanto ripristinando quei vincoli interni che col descritto sistema di tagli noi avevamo soppressi, ma anche applicando alle due faccie di ciascun taglio due distribuzioni di forze, fra loro eguali e contrarie, equipollenti punto a punto alle azioni già sviluppate da quei vincoli.

Nei riguardi del solido tagliato queste forze debbono evidentemente considerarsi come delle forze esterne. Noi riterremo le loro tre componenti secondo i tre assi (riferite all'unità di area)

$$P_x \qquad P_y \qquad P_z$$

funzioni delle coordinate, completamente definite in ogni punto delle superficie lungo le quali sono stati praticati i tagli, e riferite per ciascun taglio ad una sua faccia ben determinata, se pure arbitraria, intendendo che sull'altra faccia dello stesso taglio le analoghe componenti siano rispettivamente

$$-P_x \qquad -P_y \qquad -P_z$$

Le componenti di spostamento

$$u \qquad v \qquad w$$

determinate dal solido tagliate da un tale sistema di forze, epperò per presupposto caratterizzanti il passaggio del solido stesso dallo stato non deformato allo stato naturale, saranno esse pure delle funzioni delle coordinate, che noi riterremo in ogni punto del solido tagliato soddisfacenti a tutte quelle restrizioni che sono abituali nella teoria classica dell'elasticità.

Ciò posto noi procederemo all'applicazione del principio di reciprocità di Betti scrivendo che il lavoro che le forze $P_x$, $P_y$, $P_z$ compirebbero qualora ai loro punti di applicazione venissero attribuiti certi spostamenti u', v', w' determinati da un secondo sistema arbitrario di forze esterne $P'_x$, $P'_y$, $P'_z$ deve essere eguale al lavoro che le forze di questo secondo sistema compirebbero nell'ipotesi che ai loro punti di applicazioni si attribuissero gli spostamenti u, v, w.

E ciò faremo assumendo il secondo sistema di forze in modo che le relative componenti di deformazione - e quindi anche, per la supposta omogeneità del solido, le componenti speciali di tensione - siano costanti in tutto il corpo cioè indipendenti dalle coordinate.

È facile vedere dove conduca questa ipotesi: le derivate parziali delle componenti speciali di tensione rispetto alle coordinate devono infatti risultar tutte nulle; perciò le equazioni indefinite per l'equilibrio richiedono che siano nulle tutte le forze di massa.

Quanto alle forze superficiali esse potranno essere diverse a seconda dei valori che noi vogliamo attribuire alle componenti speciali di tensione. Se per esempio noi assumiamo le tre componenti normali eguali dell'unità, supponendo nulle le tre componenti tangenziali, le equazioni ai limiti risultano identicamente verificate per

$$P'_x = -\cos(n,x) \qquad P'_y = -\cos(n,y) \qquad P'_z = -\cos(n,z)$$

n essendo come d'uso la normale in un punto generico della superficie, rivolta verso l'interno del corpo.

Ciò equivale a dire che il solido è soggetto ad una forza uniformemente distribuita su tutta la superficie, diretta in ogni punto secondo la normale alla superficie stessa, rivolta verso l'esterno, e di intensità ovunque eguale all'unità.

Il significato del primo membro dell'equazione di Betti risulta ora immediato. Si ha infatti, ove si indichi con V lo spazio occupato dal solido, e con S il complesso delle superficie che lo delimitano,

$$\int_S \left(P_x' u + P_y' v + P_z' w\right) dS =$$

$$-\int_S \left(u\cos(n,x) + v\cos(n,y) + w\cos(n,z)\right) dS = \int_V \left(\frac{\delta u}{\delta x} + \frac{\delta v}{\delta y} + \frac{\delta w}{\delta z}\right) dV$$

E questo integrale rappresenta notoriamente la dilatazione cubica totale che il solido tagliato subisce nel cambiamento di configurazione caratterizzato dalle u, v, w, vale a dire nel passaggio dallo <u>stato non deformato</u> allo <u>stato naturale</u>.

L'altro membro dell'equazione di Betti

$$\int_S \left(P_x u' + P_y v' + P_z w'\right) dS$$

risulta invece identicamente nullo.

Per convincersene basta osservare che la superficie S del solido tagliato a cui l'integrale va esteso, consta di due parti ben distinte: quella che propriamente limita lo spazio V occupato dal corpo, cioè che separa tale spazio dallo spazio circostante, e quella, situata all'interno di V, che ha origine dai tagli.

Sulla prima sono ovviamente nulle le $P_x$, $P_y$, $P_z$, nessun dubbio quindi che sia nullo anche il corrispondente contributo all'integrale.

Lo stesso non accade sull'altra: accade però, e lo abbiamo già rilevato a suo tempo, che per ogni elemento superficiale scelto su di una data faccia di un taglio, esiste sempre sull'altra faccia dello stesso taglio un elemento identico al primo e con esso coincidente, pel quale la forza differisce soltanto nel segno. Se si tien conto che, per due elementi siffatti, lo spostamento deve essere lo stesso - le u', v', w' essendo lineari nelle coordinate, epperò continue anche attraverso i tagli - si può concludere che ad ogni termine del tipo

$$P_x u' + P_y v' + P_z w'$$

deve sempre corrispondere un altro esattamente uguale a

$$-P_x u' - P_y v' - P_z w'$$

L'integrale deve dunque in definitiva annullarsi come avevamo annunciato.

Di qui il teorema

$$\int_V \left(\frac{\delta u}{\delta x} + \frac{\delta v}{\delta y} + \frac{\delta w}{\delta z}\right) dV = 0$$

**30**

Gian Antonio Maggi a Gustavo Colonnetti

Pisa 8 Agosto 1918

Caro Professor Colonnetti,

Il teorema di Betti, esteso all'ipotesi dell'esistenza di superficie di discontinuità, S*, per le due specie (in generale) di spostamenti (u,v,w) a cui corrispondono ($P'_x$, $P'_y$, $P'_z$), nell'assenza di forze di massa, diventa, indicando con S la superficie limitante il corpo

$$\int_S (P'_x u+P'_y v+P'_z w)dS+\int_{S^*}(P'_x\Delta u+P'_y\Delta v+P'_z\Delta w)dS^*$$

$$=\int_S(P_x u'+P_y v'+P_z w')dS+\int_{S^*}(P_x\Delta u'+P_y\Delta v'+P_z\Delta w')dS^*,$$

dove Δ indica, naturalmente, la discontinuità.

Per le ipotesi intrinseche del problema, sopra S, si ha $P_x=P_y=P_z=0$, e le Δu, Δv, Δw (su S*) non sono tutte nulle, ché diversamente le u, v, w sarebbero o nulle, o corrisponderebbero ad uno spostamento rigido.

Facevamo poi l'ipotesi, sopra S,

$$P'_x=-\cos\widehat{nx},\qquad P'_y=-\cos\widehat{ny},\qquad P'_z=-\cos\widehat{nz},$$

e sopra S*,

$$\Delta u'=\Delta v'=\Delta w'=0.$$

Le u', v', w' ne risultano determinate (sono uno spostamento rigido) in ogni punto del corpo: e quindi determinate, sopra S*, le $P'_x$, $P'_y$, $P'_z$.

Introduciamo queste ipotesi nella (1). Ne risulta (annullando il secondo membro)

$$-\int_S(u\cos\widehat{nx}+v\cos\widehat{ny}+w\cos\widehat{nz})dS+\int_{S^*}(P'_x\Delta u+P'_y\Delta v+P'_z\Delta w)dS^*=0,$$

ossia

$$\int_V\left(\frac{\delta u}{\delta x}+\frac{\delta v}{\delta y}+\frac{\delta w}{\delta z}\right)dV=\int_{S^*}(P'_x\Delta u+P'_y\Delta v+P'_z\Delta w)dS^*.$$

Così io trovo il secondo membro diverso da 0. E d'altra parte nella deduzione ch'Ella fa della sua formula, quando ne costruisce il nuovo membro, non vedo traccia della superficie S*.

Può darsi benissimo ch'io mi inganni, e perciò mi farà molto piacere a scriverne di nuovo, perché in questo come in altro, non potremmo far altrimenti che metterci completamente d'accordo.

**31**

Gian Antonio Maggi a Gustavo Colonnetti

Pisa 8 Agosto 1918

Caro Professor Colonnetti,

Faccio seguire, a distanza di un treno, questa alla mia lettera di stamane per la seguente cagione. Tenuto calcolo che le u, v, w presentano alla superficie S*, interna al campo τ, le discontinuità Δu, Δy, Δw, va scritto

$$\int_V \left(\frac{\delta u}{\delta x}+\frac{\delta v}{\delta y}+\frac{\delta w}{\delta z}\right) dV= -\int_S (u\cos\widehat{nx}+v\cos\widehat{ny}+w\cos\widehat{nz})dS -\int_{S^*} (\Delta u\cos\widehat{nx}+\Delta v\cos\widehat{ny}+\Delta w\cos\widehat{nz})dS^* .$$

Quindi la formola finale, secondo la mia deduzione, è

$$\int_V \left(\frac{\delta u}{\delta x}+\frac{\delta v}{\delta y}+\frac{\delta w}{\delta z}\right) dV= -\int_{S^*} \left[\Delta u(P_x'+\cos\widehat{nx})+\Delta v\left(P_y'+\cos\widehat{ny}\right)+\Delta w(\quad)\right]dS^*.$$

E si vede che il secondo si ridurrà a 0, quando sulla superficie S*, come sulla S, <u>risulti</u>

$$P_x'=-\cos\widehat{nx}, P_y'=-\cos\widehat{ny}, P_z'=-\cos\widehat{nz}.$$

Di nuovo le migliori cose, e a rivederci presto

Aff$^{mo}$ Suo

G.A. Maggi.

**32**

Gustavo Colonnetti a Gian Antonio Maggi

[cartolina postale indirizzata: Al Ch$^{mo}$ Sig Prof. Gian Antonio Maggi

Via Risorgimento 1 - Pisa]

Roma 9 ore 21[66]

Egregio Professore.

Trovo ora la sua seconda lettera, dopo la quale mi pare che l'accordo deve divenire completo: la mia di oggi le avrà infatti chiarito che la condizione a cui Ella giunse faceva parte integrante della mia ipotesi. Spero a voce di poter avere i suoi consigli sul miglior modo esporre la questione.

Con rinnovati ringraziamenti ed auguri

dev$^{mo}$

G. Colonnetti

[66] 1918, dal timbro postale.

**33**

Gustavo Colonnetti a Gian Antonio Maggi

[carta intestata: Ministero per le armi e munizioni]

Roma 9 agosto 18

Egregio Professore,

Grazie vivissime: quanto al disaccordo ecco dove ha origine:

Ella scrive: facciamo poi l'ipotesi sopra S:

$P'_x = -\cos nx,\ P'_y = -\cos ny,\ P'_z = -\cos nz$ e sopra S* Δ'u=Δ'v=Δ'w=0.

Ora la mia ipotesi era un poco diversa: io volevo che nel solido tagliato riuscissero costanti le componenti di tensione, e precisamente supponevo che di esse le tre normali fossero =1, e le tre tangenziali =0. Ora ciò implica che si abbia

$$P'_x = -\cos nx,\ P'_y = -\cos ny,\ P'_z = -\cos nz$$

su tutta la superficie del solido tagliato cioè non soltanto su S ma anche su S*.
Le Δ'u=Δ'v=Δ'w=0 vengono poi di conseguenza perché le u', v', w' riescono lineari nelle x y z e quindi continue anche attraverso i tagli.

Messe così le cose a me sembra che l'integr.

$-\int_S \big(u\cos(n,x)+v\cos(n,y)+w\cos(n,z)\big)dS$ rappresenta l'intiero primo membro dell'equazione di Betti purché lo si intenda esteso a tutta la superficie S + S* del solido tagliato. E d'altronde è necessario far così se si vuole che il suo trasformato

$$\int_V \left(\frac{\delta u}{\delta x}+\frac{\delta v}{\delta y}+\frac{\delta w}{\delta z}\right) dV$$

rappresenti la variazione totale del volume del solido tagliato nel passaggio dallo stato non deformato allo strato naturale (di equilibrio in assenza di forze esterne).

Se ha ancora dei dubbi pensi al caso classico delle distorsioni in un anello: dal momento che vogliamo applicare il teor. di Betti all'anello tagliato considerando le azioni attraverso i tagli come forze esterne, perché tener distinte le faccie AB A'B' dal resto della superficie? A me pare che bisogna considerare quelle faccie alla stessa stregua del resto: se no si calcolerà bensì la variazione di volume dovuta ai cambiamenti di forma della superficie laterale, ma si verrà a trascurare lo spazio inizialmente compreso fra AB e A'B' che viene occupato dal solido quando le due faccie sotto l'azione delle $P_x$ $P_y$ $P_z$ vengono ad accostarsi riproducendo la configurazione preesistente al taglio.

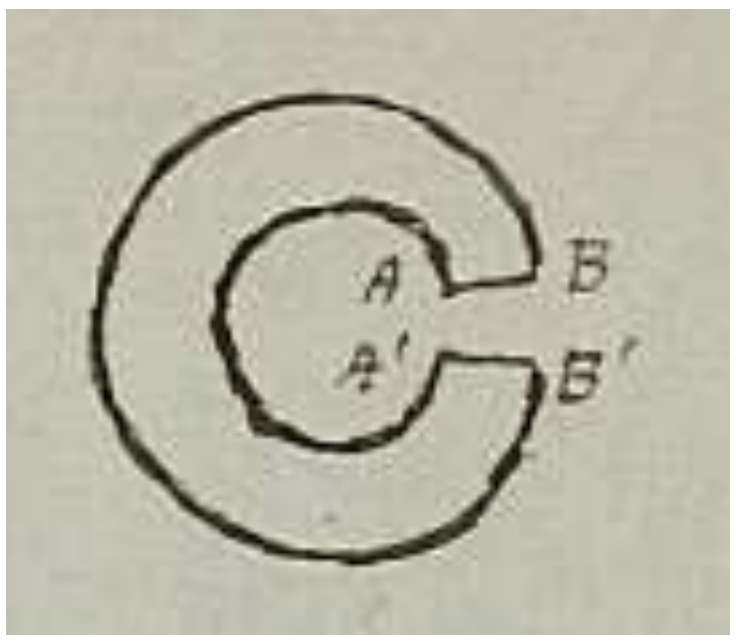

Mi sbaglio? O mi sono soltanto spiegato male nella mia prima lettera?

Io dovrò restare qui ancor parecchi giorni, ma spero di trovar egualmente modo di fermarmi a Pisa nel ritorno.

Intanto grazie ancora e distinti saluti anche per la Sua Famiglia.

di Lei dev^mo

G. Colonnetti

Poco fa girovagando con l'amico Vacca abbiamo trovati da un antiquario i varii libri suoi, insieme con quelli del Dini, del Bianchi ecc. in mezzo ad una collezione di libri russi. Portavano tutti la sigla E.S.B e pare provengano da quella certa signora *[*Bodareu*]* (dico bene?) di cui ho sentito a parlare tanto all'antro, la quale deve essere stata espulsa dal regno per ragioni di pubblica sicurezza. Lo sapeva Lei?

**34**
Gian Antonio Maggi a Gustavo Colonnetti

Pisa 15 Agosto 1918

Caro Professor Colonnetti,

Doppiamente mi duole ch'Ella abbia dovuto bruler Pisa, per la ragione che ve la determinò. Ma spero che il Suo disturbo sarà stato passeggero, e che sarà a quest'ora interamente ristabilito. Le avrei volontieri discorso della nota questione, destinata forse a ridestarsi anche prima del preveduto. Non però in questo stesso momento, per modo che se ne può rimandare il discorso senza conseguenza.

Per quanto al problema elastico, mi spiace che la mia seconda lettera non Le sia stata recapitata abbastanza presto, per dispensarla di trattenersi sopra una formula incompleta. Mi accorsi della necessaria correzione appena, tornato d'aver portato la mia prima lettera alla Posta, e la seconda avrebbe dovuto prendere la via di Roma, col treno successivo. Servì meglio allo scopo la seconda nave mandata a Mitilene da Atene, quando riconobbe l'errore del consiglio di Cleone. Ma è pur vero che furono presi i provvedimenti raccontati da Tucidide, e che s trattava di revocare la strage in massa dei Mitilenesi… La formula esatta, per tornare, abbastanza da lontano, alle distorsioni è, secondo la mia deduzione

$$\int_{\tau}\left(\frac{\delta u}{\delta x}+\frac{\delta v}{\delta y}+\frac{\delta w}{\delta z}\right)d\tau=\int_{S^*}\left[\Delta u(P_x'+\cos\widehat{nx})+\Delta v\left(P_y'+\cos\widehat{ny}\right)+\Delta w(P_z'+\cos\widehat{nz})\right]dS^*$$

Ho piacere ch'Ella ne tragga indizio di accordo, e aspetto, con interesse, di sentire il definitivo risultato. La via da Lei accennata nella Sua lettera mi sembra buona, tenuto calcolo che quando si suppone l'esistenza di una superficie di discontinuità, si aggiungono alle condizioni determinanti l'equazione elastica i valori delle stesse discontinuità nei punti della superficie, oppure quelli delle pressioni elastiche negli stessi punti. (cfr §13 della mia Dinamica dei Sistemi).

**35**
Gian Antonio Maggi a Gustavo Colonnetti

Pisa 18 Agosto 1918

Caro Professor Colonnetti,

Ho oramai verificato che il Suo risultato è perfettamente esatto, per le ragioni da Lei accennate. Difatti dalle condizioni che u', v',w', siano continue, e che sul contorno siano verificate le

(1) $P'_x = -\cos\widehat{nx}$, $P'_y = -\cos\widehat{ny}$, $P'_z = -\cos\widehat{nz}$,

segue che queste eguaglianze sono soddisfatte per ogni punto del corpo, e per ogni semiretta n. Quindi risulta nullo il secondo membro della

(2) $\int_\tau \left(\frac{\delta u}{\delta x} + \frac{\delta v}{\delta y} + \frac{\delta w}{\delta z}\right) d\tau = \int_{S^*} \left[\Delta u(P'_x + \cos\widehat{nx}) + \Delta v\left(P'_y + \cos\widehat{ny}\right) + \Delta w(P'_z + \cos\widehat{nz})\right] dS^*$

Per quanto alla miglior via per stabilire il risultato in discorso, mi parrebbe questa. Partire dal Teorema di Betti, nell'ipotesi di spostamenti discontinui

$$\int_\sigma (P'_x u + \ldots) d\sigma + \int_{S^*} (P'_x \Delta u + \ldots) dS^* = \int_\sigma (P_x + \ldots) d\sigma + \int_{S^*} (P_x \Delta u' + \ldots) dS^*.$$

Introdurvi l'ipotesi

$$\Delta u' = \Delta v' = \Delta w' = 0, \quad P'_x = -\cos\widehat{nx} \quad P'_y = -\cos\widehat{ny}, \quad P'_z = -\cos\widehat{nz}$$

In S* In σ

con che si ha

$$-\int_\sigma (u\cos\widehat{nx} + \ldots) d\sigma - \int_{S^*} (\Delta u \cos\widehat{nx} + \ldots) dS^* = \int_\sigma (P_x u' + \ldots) d\sigma$$

ossia, pel teorema di Gauss,

$$\int_V \left(\frac{\delta u}{\delta x} + \frac{\delta v}{\delta y} + \frac{\delta w}{\delta z}\right) dV = \int_\sigma (P_x u' + \ldots) d\sigma,$$

dove la u'=ax, v'=ay, w'=az.

Applicare la formola al caso della distorsione per cui $P_x \mp P_y \mp P_z$, su σ, sono nulli. E con questo

$$\int_V \left(\frac{\delta u}{\delta x} + \frac{\delta v}{\delta y} + \frac{\delta w}{\delta z}\right) dV = 0.$$

Non occorre insistere sulla soluzione lineare, ch'Ella avrà pure trovato (u'=ax, v'=−ay, w'=az), con che si ha la formula più generale

$$\int_V \left(\frac{\delta u}{\delta x} + \frac{\delta v}{\delta y} + \frac{\delta w}{\delta z}\right) dV = a \int_\sigma \left(P_x + P_y + P_z\right) d\sigma.$$

...

**36**

Gustavo Colonnetti a Gian Antonio Maggi
[carta intestata: R. Università di Pisa - Scuola di applicazione per gl'Ingegneri - Gabinetto di costruzioni]

Roma 22 sett 18 (Hotel Bologna)

Egregio Professore,

La prego scusarmi se ho tardato a rispondere alla Sua: il ritardo è dovuto al fatto che desideravo rileggere prima le bozze della mia nota, bozze che ho avuto solo in questi giorni.

Ho creduto bene mantenere la esplicita riserva già accennatale nella mia precedente perché la $\sigma_Y$ appare in principio del mio scritto come una funzione scelta ad arbitrio. Faccio poi osservare in nota che la supposta continuità non è (come ben dice Lei) una limitazione, ma una condizione necessariamente verificata nel problema fisico.

Su di un sol punto non siamo d'accordo ed è sulla continuità di $\sigma_t$ che a me non sembra necessaria perché si tratta di tensione che agisce su elementi radiali cioè attraversanti le supposte superficie cilindriche di discontinuità degli spostamenti.

Ma la cosa importa poco per ora perché ai fini del mio calcolo non occorre precisarla.

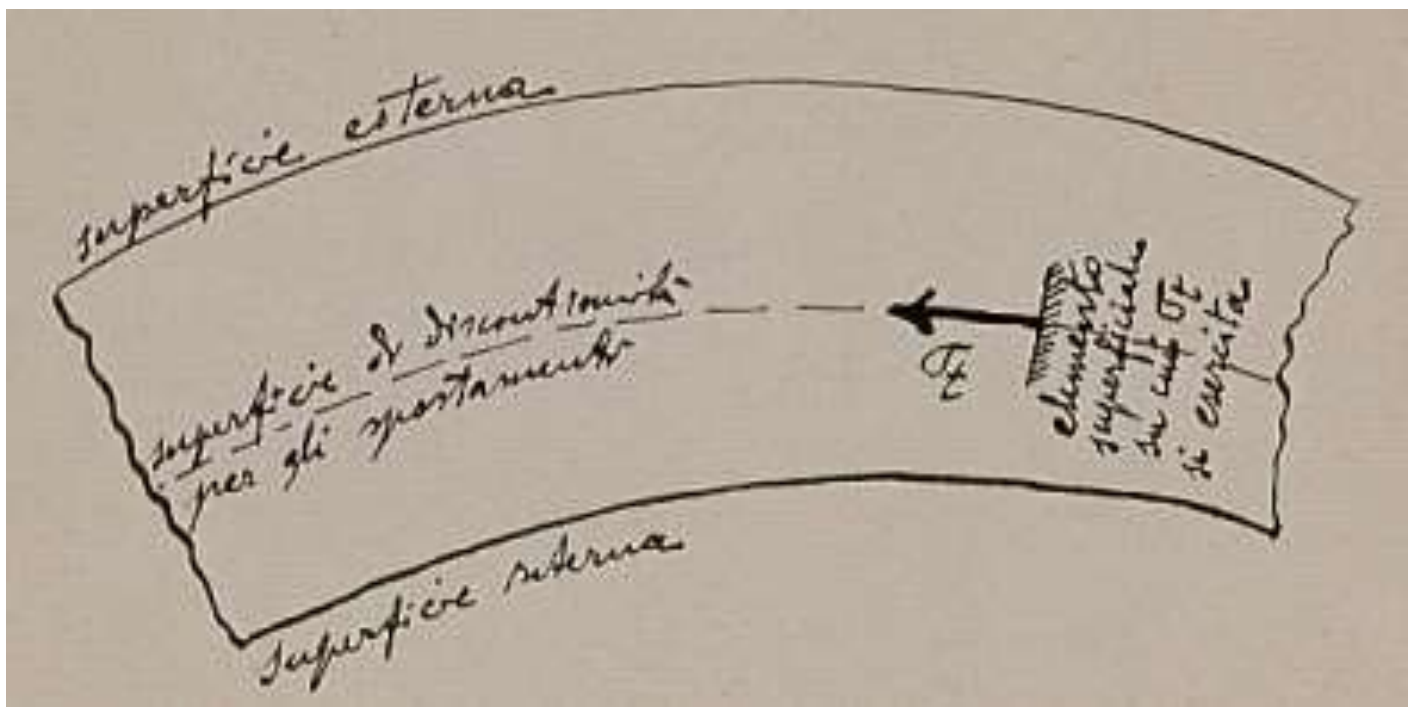

Le accludo il manoscritto della seconda nota che Ella potrà rileggere con tutto suo comodo ed inviare poi col suo visto al Comm. Mancini se le sembrerà meritevole. Resta sempre ben inteso che le sarò gratissimo di tutte le osservazioni e dei consigli di modificazioni che Ella avesse a farmi.

Come vedrà mi sono attenuto alla applicazione concrete del teorema di reciprocità <u>al solido tagliato</u>, perché credo sia più facilmente comprensibile dagli ingegneri che troverebbero forse ostica l'applicazione al solido intiero dotato di discontinuità a cui essi non sono abituati.

A quest'ora Ella avrà saputo che il Congresso è stato rinviato a data indeterminata (e credo che la notizia non Le avrà fatto molto dispiacere).

Il prof Nasini mi ha portato notizie molto gravi del sen. Dini; da lui io non ebbi più nulla e sarei molto grato a Lei se volesse tenermi informato della tristissima cosa. Il Rettore è sempre deciso a non adunar il Consiglio della Scuola? E non pensa che a questo modo arriveremo all'apertura dell'anno scolastico senza avere i nuovi incaricati?

Gradiremo anche vivamente notizie della Sua Famiglia, alla quale invio anche a nome di Mamma[67] e di Gemma i più cordiali voti ed i più distinti saluti. A Lei rinnovato grazie pei disturbi che Le arreco. Sempre dev.$^{mo}$

G. Colonnetti

[67] Paoletta Calligaris.

**37**
Gian Antonio Maggi a Gustavo Colonnetti
[carta intestata: R. Università di Pisa - Scuola di applicazione per gl'Ingegneri - Gabinetto di costruzioni]

Pisa 23 Dicembre 1921

*[in pastello rosso, scritto da Maggi:]* Copia conforme
Caro Professor Colonnetti,

Comincio col presentare a Lei, e alla Signora, che tanto lieta accoglienza fece alle mie viaggiatrici, i migliori e più cordiali augurii per le imminenti Feste di Natale e del Capo d'Anno. E per Lei, vi faccio seguire la retroguardia dei ringraziamenti, di cui Le mandai l'avanguardia colla mia Noticina dei Lincei, pel gradito dono delle sue Lezioni di Dinamica.

Ho percorso il volume dal principio alla fine, e vi ho rilevato i notevoli pregi dell'originalità della composizione, del rigore e della perspicacità dell'esposizione. Il suo metodo è abbastanza diverso dal mio, cominciando dalla circostanza che Ella ci presenta un buon saggio di metodo analitico, col premettere lo sviluppo storico di fondamentali argomenti. Io sono troppo convinto dei concetti a cui è informato quello che lo Schiaparelli, in una lettera al Pizzetti, chiamava la mia rivoluzionaria meccanica, (con quanta pace dell'Europa altri concetti rivoluzionari sarebbero stati parimente riservati agli scolari, che si preparano agli esami!) per non mantenermi fedele al mio modo di vedere, ma certo non sono così esclusivo da non riconoscere l'efficacia anche di un procedimento informato a concetti diversi.

Vi è però un punto, sul quale non posso accordarmi con Lei, ed è sul posto d'onore che Ella conferisce alla forza d'inerzia, per stabilirvi l'equazione fondamentale e il principio di d'Alembert. Io non so attribuire alla $F-ma=0$ altro significato che quella della $F=ma$, esprimente che la forza motrice è rappresentata dal prodotto della massa per l'accelerazione. Come equazione di Statica, essa esprime che al punto materiale considerato sono applicate due forze eguali e contrarie che danno per somma la forza nulla: e di qua io non so dedurre altro fuorché dove risultava nulla l'accelerazione, mentre, nella stessa ipotesi, compete al punto, al considerato istante, l'accelerazione A. Per la stessa ragione, non mi persuade simile applicazione della forza d'inerzia alla deduzione (quasi una verificazione a priori) del principio d'Alembert. Il ragionamento originario di d'Alembert, a parte le riserve a cui dà luogo, porta precisamente alla conclusione che le accelerazioni corrispondenti al sistema delle forze $F-ma$, subordinatamente all'azione dei vincoli, sono nulle.

Newton introduce la vis inertiae colla Definizione III dei Principia, ma né costì, né, ch'io sappia, altrove, ne indica la rappresentazione mediante $-ma$. La parola vis in questo caso, mi sembra significare piuttosto attitudine. E l'opposizione della vis impressa alla newtoniana vis inertiae credo debba paragonarsi con quella della buona volontà alla pigrizia, la quale deve dare un risultato positivo. Buona parte dei trattati ch'io conosco non accennano alla forza d'inerzia. Il Delaunay (Mécanique, pag. 221)[68] insiste perché s'intenda non applicata al punto, insieme colla F, ma dal punto all'altro, che vi esercita la F. Marcolongo e Burgatti hanno però ultimamente trasportato il nome di forza d'inerzia alla ma, quella che io chiamo forza motrice effettiva. Tutto quanto mi sembra mostrare che l'ostracismo della forza d'inerzia, come −ma, applicata al punto considerato, non è una mia particolarità.

---

[68] C.E. Delaunay, *Traité de mécanique rationnelle*, Paris, 1856.

Per quanto allo Scholium di Newton, avendo ricercato il relativo passo del Thomson e Tait, nella mia copia, l'ho trovato postillato in margine con un ?. E, dopo molti anni da quella postilla, resto sempre dall'opinione che, per riscontrare in quel brano di Newton il principio di d'Alembert, bisogna aggiungervi una quantità di cose, che, se Newton avesse pensato, avrebbe altrimenti sviluppato. Quel brano contempla, anche di più, l'equazione di d'Alembert e Lagrange. È ammissibile che Newton, quanto avesse scoperto un principio di così straordinaria importanza, l'avesse poi lasciato dormire, finché, non d'Alembert, non Lagrange, ma Thomson e Tait lo destassero, dopo due secoli di sonno? Infine Newton avrebbe dovuto valersi della vis inertiae, sotto le spoglie vires... ex acceleratione oriundæ,[69] con un concetto della medesima, che, come dicevo più sopra, non mi sembra risultare dalla relativa definizione.

Sentirò molto volontieri quanto Ella mi vorrà replicare in proposito. Chiudo ora, rinnovando gli augurii, e pregandoLa di aggradire insieme per sé e per la signora da parte di noi tutti, i più distinti saluti. Mi creda sempre

Aff.mo Collega
G.A. Maggi

## 38
## **Mario Crenna** a Gian Antonio Maggi

Illustre e venerato professore,

lo scrivente, meditando lungamente sulla concezione einsteiniana, è stato portato ad alcune considerazioni critiche che non gli sembrano prive di un qualche interesse. Memore della paterna benevolenza che la S.V. Ill.ma gli ha sempre dimostrato il sottoscritto ha pensato di rivolgersi a Lei, pregandoLa umilmente e vivamente di voler degnare di un'occhiata il manoscritto qui accluso e, nel caso di un giudizio favorevole, di volerlo benevolmente presentare presso il R. Istituto Lombardo-Veneto per la pubblicazione in quei Rendiconti.

Lo scrivente non sa trovare espressioni atte ad esprimere alla S.V. Ill.ma i sentimenti della sua obbligazione profondissima e nel porgerLe le più vive scuse, con infiniti ringraziamenti, La prega di voler gradire l'espressione della sua immutabile venerazione.

Devmo ed obblmo
Mario Crenna
Via Fiasella, 7-Sarzana (Spezia)

Sarzana, 16-2-30

[69] Da I. Newton, *Philosophiæ naturalis principia mathematica*, p. 27.

**39**

Gian Antonio Maggi a Mario Crenna
[carta intestata: R. Università di Milano - Istituto matematico -
Via C. Saldini, 50 (Città degli Studi) - Il Direttore]

Milano 1 Marzo 1930

Caro Dottor Crenna,

Per diverse occupazioni non ho potuto veder prima il Suo manoscritto, che mi spiace di non poter presentare all'Istituto Lombardo, per cui Glielo rimando, colle seguenti considerazioni, colle quali procuro di renderLe ragione della mia riserva.[70]

Per seguire la Sua esposizione, nel §1 l'enunciato del Suo proposito di dimostrare che l'interpretazione del problema relativista, fornita dall'Einstein è errata, mi fa veramente domandare che cos'è il problema relativista, all'infuori dell'Einstein. Segue poi una corretta impostazione di quello che si può chiamare il problema di Lorentz.

Nel §2, accennato al passaggio dal coup de pouce del Lorentz alla concezione einsteiniana, celebrata dall'Occhialini, Ella afferma esplicitamente il suo compito, che è di dimostrare che "quella concezione è assolutamente insostenibile, e che impostato il problema della relatività dei fenomeni fisici con le formole di Lorentz, l'unica interpretazione possibile di esse è quella che è fornita dalla concezione classica".

Io non vedo veramente interpretazione. Secondo la concezione classica, delle formole di Lorentz all'infuori della contrazione reale, a cui non tardò a rinunciare lo stesso Lorentz: ammessa la quale, le formole del Lorentz cessano di essere relativistiche. Ma poi, nel seguito del §, io non trovo che un'esposizione dei concetti einsteiniani, alla quale non trovo da obbiettare che la chiusa, dove dice che "la teoria einsteiniana dà, ad ogni osservatore, la conoscenza fisica del proprio sistema, ma circa la conoscenza fisica del complesso degli altri sistemi, essa non è, e non potrebbe essere, una teoria soggettiva, ma è evidentemente una teoria agnostica". Perché agnostica? L'osservatore appartenente al sistema galileiano K dalle misure che eseguisce sulle quantità appartenenti al sistema K', deduce, in base alle formole di Lorentz, le misure che alle stesse quantità appartengono al sistema K'. Ha però ragione se intende che queste misure non può ottenerle coll'osservazione diretta di un regolo o di un orologio appartenenti a K'.

Col §3 comparisce il tempo proprio, che Le fornisce l'oggetto principale della Sua critica. Ora, il tempo proprio si riferisce sempre ad un punto mobile che descrive una linea d'universo, e quando il movimento di questo punto sia rettilineo, uniforme, com'è supposto in questo §, il tempo proprio corrispondente si può introdurre semplicemente (cfr. Kopff, (53))[71] con

(1) $$t^{*}=\sqrt{1-\frac{x^2+y^2+z^2}{c^2}}\,t=\sqrt{1-\frac{q^2}{c^2}}\,t$$

(preferisco pel tempo proprio l'asterisco all'apice, riservando questo per un sistema galileiano qualsivoglia, K', diverso quando occorre intenderlo da K) dove t, x, y, z e q (grandezza della veloc. del punto) si riferiscono ad un certo sistema galileiano K, espressione invariante con questo, perché si ha egualmente

---

[70] Maggi scrive un fascicolo di 14 pagine di commenti puntuali al manoscritto ricevuto per poi riportare i punti salienti in questa lettera.

[71] Probabilmente si riferisce all'articolo A. Kopff, "Über eine Möglichkeit der Prüfung des speziellen Relativitätsprinzips auf astronomischem Wege", 1922.

$$t^{*}=\sqrt{1-\frac{x'^{2}+y'^{2}+z'^{2}}{c^{2}}}\,t'=\sqrt{1-\frac{q'^{2}}{c^{2}}}\,t'.$$

Ciò posto, se s'immagina il sistema galileiano K* che si muove colla stessa velocità del punto in discorso, questo è un punto invariabilmente unito con K*, vale a dire è per questo sistema q=0: per cui la formola (1) dà pel sistema medesimo t*=t. Cioè il tempo proprio corrispondente ad un punto, il cui moto è rettilineo uniforme, non è altro che il tempo appartenente al sistema galileiano K* che possiede la stessa velocità. Con questo esso è ancora il tempo segnato da un orologio appartenente a K*. Vedremo poi che questo discorso va limitato all'ipotesi del punto in movimento rettilineo uniforme.

Per fissare il concetto del tempo proprio nel presente caso non occorre altro. Ella riconosce la "ragione di questa denominazione sul fatto che la relativa espressione rappre*[senta]* realmente il tempo segnato fisicamente dall'orologio di K* legato al punto mobile."Il Minkowski (Raum und Zeit, pag. 9) introducendo questa variabile dice “lo chiamiamo Eigenzeit”, senza accennare a ragione, e io dubito che tempo proprio sia la miglior traduzione di Eigenzeit - io credo la traduzione giusta tempo speciale - (miracolo che non fu tradotto autotempo!). Comunque sia esso, fra i tempi a disposizione, quanti sono i K, è certo un tempo segnalato: ma non per questo di diversa natura del tempo t, appartenente ad un K qualsivoglia che è il tempo segnato da un orologio legato ad un punto qualunque di K. Per cui non vedo come il tempo proprio darebbe luogo ad un miscuglio di relativo e di assoluto, come se l'osservatore appartenente a K potesse leggere, sull'orologio legato ad un punto di K*, il tempo t*, e non il tempo t. Se qualcuno così dice, non c'è che da richiamare le definizioni, per riconoscere l'equivoco, al quale non è forse del tutto estraneo il vocabolo tempo proprio, traduzione, per me discutibile, di Eigenzeit.

Nel §4, premessa la simultaneità delle equazioni $x^2-c^2t^2=0$ e $x'^2-c^2t'^2=0$, ch'Essa vi connette colle formole di Lorentz ridotte (11), da cui segue $x^2-c^2t^2=\left(1-\frac{v^2}{c^2}\right)(x'^2-c^2t'^2)$, procede ad un confronto della propagazione della luce in linea retta nei due riferimenti A, M, B e A', M', B', il secondo mobile rispetto al primo, col quale, a istanti corrispondenti, attribuisce, in ambedue i riferimenti, lo stesso tempo. Per cui io non so veramente seguirlo in quel ragionamento, di cui l'ultima correlazione è poi che, per conseguire la nota inavvertibilità del movimento, bisogna attribuire alla veloc. di propagaz. della luce lo stesso valore in ambedue i riferimenti, rinunciando alle ipotesi precedenti. Che se fosse partito, senz'altro, da quella ipotesi avrebbe ritrovato le formole di Lorentz, come fa l'Einstein in Appendice al ricordato libro (Complemento del §11). Strada facendo, trovo poi l'affermazione che il principio di relatività non esige che i fenomeni fisici siano visti avvenire con la stessa velocità: a proposito della quale osservo che, se si tratta della relatività dinamica, la velocità non vi ha nessuna parte, la condizione d'invarianza limitandosi all'accelerazione, mentre nella relatività ottica o elettromagnetica, è essenziale l'invarianza delle veloc. di propagaz. della luce.

Infine nel §5 ricompare l'orologio connesso col punto mobile, che ne segna il tempo vero, oggetto da Lei particolarmente attaccato. Al qual proposito, osservo che mentre, come s'è detto, nel caso del punto in movimento rettilineo uniforme detto orologio ha un significato esente da ogni discussione, esso non ha più alcun significato, nel caso che il punto, a cui si riferisce il tempo vero t*, sia in movimento vario. Perché il tempo dt*, rappresentato da

$$c^2 dt^{*} = c^2 dt^2 - (dx^2 + dy^2 + dz^2) = c^2\left(1 + \frac{q^2}{c^2}\right) dt^2$$

corrispondente a dt, differenza tra i tempi t+dt e t, ambedue segnati dall'orologio appartenente al sistema galileiano K, connesso, se si vuole, con un punto qualunque di K, è bensì, alla sua volta, la differenza tra i tempi t*+dt* e t*, segnati da un orologio appartenente al sistema galileiano K* che si muove alla velocità del punto al tempo t, ma, variando questa velocità, di grandezza e d'orientazione, col tempo t, i dt* corrispondenti ai singoli tempi t, appartengono tutti a orologi diversi, e la loro somma, rappresentata dall'integrale, non più riferita, in nessun modo, ad un unico orologio connesso col punto. Questo non toglie nulla all'efficacia dell'applicazione del tempo proprio, al quale basta di essere una variabile atta a individuare i punti della linea d'universo, appartenente al considerato punto mobile, segnalata per la qualità di essere indipendente dal riferimento del moto del punto all'uno piuttosto che all'altro dei possibili sistemi galileiani K. Con questo, a mio modo di vedere, cade da sé ogni obbiezione attinente al suddetto orologio, e al tempo proprio concepito come tempo, vero o fisico, segnato da esso.

Così, per concludere, io non desumo dall'esame del Suo M.S. false interpretazioni che la Sua critica rilevi e corregga, imputabili all'Einstein, piuttosto che per avventura, ad altri, di cui mi è difficile riconoscere l'equivoco. Che se m'inganno, non è perché non mi sia studiato di capire. A voce potrei spiegarmi più ampiamente, e lo farei ben volentieri se qualche volta verrà da queste parti, procurandomi il piacere di rivederLa.

Intanto Le sono ben grato della Sua costante buona memoria, e con cordiali saluti mi confermo

Aff.mo Suo
G.A. Maggi.

**40**

Mario Crenna a Gian Antonio Maggi

[busta indirizzata a: Illustre e Chiarissimo Signore Sig. Prof. Gian Antonio Maggi - Corso Plebisciti, 3 - Milano; scritto a matita da Maggi: colla mia risposta]

Illustre e venerato professore,

la paterna benevolenza che la S.V. Ill.ma manifesta, in mille modi, allo scrivente ha posto il medesimo in una confusione somma; ma poiché lo scrivente ebbe già occasione di scoprire, durante il periodo della sua vita di studente, nella S.V. un cuore in tutto degno delle Sue doti intellettuali, così la sua confusione gli è riuscita meno imbarazzante, come a chi confida soltanto nell'indulgenza e benevolenza altrui, e tale da porgergli ancora il coraggio di chiarire alcuni punti sui quali la sua precedente esposizione è riuscita quanto mai infelice.

Il sottoscritto nel §1 ha chiamato problema relativista il problema del Lorentz, all'infuori dell'Einstein, nel senso che le formole del Lorentz traducono l'equivalenza di tutti i sistemi galileiani per la descrizione dei fenomeni fisici (e perciò il principio di relatività dei fenomeni fisici); di queste formole si ha l'interpretazione classica del Lorentz (con la reale contrazione) e quella einsteiniana; l'interpretazione einsteiniana era, secondo lo scrivente, contraddittoria e perciò insostenibile: quindi volendo risolvere il problema di relatività dei fenomeni fisici con le formole del Lorentz non restava che l'interpretazione classica (della contrazione reale); in paragrafi successivi a quelli che lo scrivente ha inviato alla S.V. egli prendeva poi ad esaminare questo punto: che cioè essendo poi a sua volta la contrazione reale alquanto artificiosa e non più sostenuta dallo stesso Lorentz, il principio di relatività dei fenomeni fisici avrebbe dovuto impostarsi partendo da un principio radicalmente diverso (postulato di Ritz, al quale si accenna in una nota del §2 e che ha ricevuto dal La Rosa[72] brillanti conferme astronomiche); questi successivi § lo scrivente non inviò alla S.V. sembrandogli eccessivo il disturbo e bastandogli di sottoporre all'ambito giudizio della S.V. Ill.ma la parte originale del lavoro. L'espressione «teoria agnostica» usata nel §2 ha appunto il significato che l'osservatore di K conosce le grandezze spaziali e temporali di K' soltanto attraverso la trasformazione di Lorentz, ma che essendo le entità spaziali e temporali di K' eterogenee con quelle di K, egli non potrà mai conoscere queste entità nella loro essenza fisica (mentre invece può conoscere quelle di K), perché per afferrare l'essenza spaziale e temporale di K' dovrebbe appartenere a K' (ma allora perderebbe la conoscenza fisica di tali entità del sistema K), ciò non si verifica nella concezione classica dell'universo: ogni sistema galileiano è perciò, dal lato della conoscenza fisica (non puramente formale, matematica) un mondo chiuso: cioè, concludendo con le parole della S.V. l'operatore di K non può avere la conoscenza delle entità spaziali e temporali di K' con l'osservazione diretta di un regolo e di un orologio appartenente a K'. Nel §3 ha insistito, oltre il bisogno, sul significato fisico del "tempo proprio" (per il moto rettilineo uniforme) perché lo Straneo contrastava un tale significato (anche nel caso del moto rettilineo uniforme): ora ha soppresso tutto quanto alla S.V. è apparso (giustissimamente) inutile. Circa il §4 si è accorto, dalle osservazioni della S.V., di non esser riuscito in nessun modo a tradurre il suo pensiero; ma costituendo esso una parentesi che, sull'argomento dibattuto, sta a sé (doveva essere, nell'intenzione iniziale un capitolo a sé di un volumetto di considerazioni sulla Relatività) e che perciò non interessa allo scopo, non osa abusare dell'indulgenza della S.V. con esporre tali considerazioni con maggior chiarezza, approfittando, per la prima occasione, di parlarne a voce alla S.V. Ill.ma che ha concesso allo scrivente un tanto piacere (gli pare che quel § contenga un nuovo punto di vista, nella questione dei principî essenziali della relatività): finalmente arrivando al §5, che contiene il punto centrale della discussione, lo scrivente esprime alla S.V. Ill.ma il pensiero che l'ha finora preoccupato: si consideri il sistema galileiano K con l'annesso orologio che segna il tempo t e il sistema K*, in moto vario, nel cui origine si trova un orologio: il moto vario percorso da K* si può scomporre in intervalli infinitesimi in ciascuno dei quali il moto è rettilineo uniforme (come una curva si scompone in tratti rettilinei infinitesimi); per ciascuno di questi intervalli l'espressione[73]

$$dt^{*}=dt^2-\frac{dx^2+dy^2+dz^2}{c^2}=\left(1-\frac{q^2}{c^2}\right)dt^2$$

dà il tempo proprio appartenente ad orologi diversi (variando la velocità di grandezza e orientazione col tempo t); ma ciascuno di questi orologi (appartenenti a differenti sistemi), nei singoli intervalli che ad ognuno compete, non possono differire dall'orologio connesso con l'origine di K*, perché vengono con esso ad identificarsi successivamente (perché nei singoli intervalli infinitesimi K* viene ad identificarsi con i singoli sistemi galileiani cui appartengono i differenti orologi dianzi accennati); quindi la somma, rappresentata dall'integrale, di tutti i successivi intervalli di tempo, appartenenti ai differenti orologi, fornisce il tempo segnato dall'orologio di K* (che si è successivamente identificato con quegli orologi). Lo scrivente ricorda anche d'aver letto un articolo dell'Einstein stesso su di una rivista redatta in tedesco «Die Naturwissenschaften» (anno 1916, se non erra; a Fiorenzuola avrebbe i dati esatti) nel quale, con l'esempio di un animaletto trasportato nello spazio con moto vario, egli interpretava il "tempo proprio" coll'identico significato dianzi accennato; tutti i trattati e gli articoli sulla Relatività che son capitati nelle mani del sottoscritto (e non sono pochi) danno questa interpretazione e lo stesso Marcolongo, che non tocca l'argomento, confutando, in una nota, il famoso esempio dei due gemelli, porta difficoltà che presuppongono l'interpretazione del "tempo proprio" data ora. Non parrebbe alla S.V. che l'articolo (con le opportune correzioni introdotte e con la soppressione del §4) possa meritare di essere pubblicato? ché nella peggiore delle ipotesi che l'interpretazione del "tempo proprio" non fosse conforme alle vedute dello scrivente (per qualche difetto inerente al ragionamento del medesimo) esso servirebbe, non foss'altro, a far chiarire tale concetto sul quale corre l'equivoco quasi universale. Se alla S.V. non apparisce opportuna la presentazione presso l'Istituto Lombardo-Veneto (richiedendo un lavoro di maggior consistenza) lo scrivente osa pregare vivamente la S.V. di volerlo raccomandare a qualche altra rivista ove non sia necessaria una presentazione ufficiale e che accolga il lavoretto, in merito alla buona parola interposta dalla S.V. Ill.ma; di che lo scrivente Le serberà gratitudine vivissima.[74] Il sottoscritto prega la S.V., con tutto il cuore, di volergli concedere venia per questa sua indiscrezione e per l'immenso disturbo che è venuto a recarle e mentre affretta l'occasione di poterLa rivedere e ringraziare a viva voce La prega anche di voler gradire l'espressione della sua obbligazione, devozione e venerazione immutabili.

[72] Potrebbe trattarsi di: M. La Rosa, "Il postulato di Ritz sulla velocità della luce ed i fenomeni delle stelle variabili", 1921.
[73] Manca un quadrato al primo membro.

Con i migliori voti del cuore

Dev.mo ed obbl.mo
Mario Crenna
Via Fiasella, 7 - Sarzana (Spezia)

Sarzana, 9-III-30

[74] Crenna pubblicherà nel 1931 l'articolo: "Considerazioni critiche sulla teoria della relatività einsteiniana", *Accademia di scienze, lettere e belle arti di Palermo. Atti*, v. XVI, pp. 131-144, nel 1948 "La teoria di relatività e la filosofia tradizionale" e nel 1949 "Essenza e contraddizioni della relatività".

**41**
Gian Antonio Maggi a Mario Crenna

Milano 23 Marzo 1930

Caro Dottor Crenna,

Le formole di Lorentz traducono l'equivalenza dei noti due sistemi galileiani coll'interpretazione dell'Einstein, ma non vedo come lo stesso possa dirsi coll'interpretazione classica. Con questa, se K è connesso coll'etere, è mobile bensì relativamente a K', ma il suo movimento è apparente, come è apparente il movimento di rotazione delle pareti, relativamente ad una persona che gira relativamente ad esse: e la contrazione, reale in K', è puramente apparente in K. Il termine relatività non mi sembra quindi più conveniente al caso.

L'ipotesi balistica cambia completamente i termini della questione. Rievocata dal Ritz e sostenuta dal La Rosa, non sembra reggere facilmente alle obbiezioni che le sono mosse, particolarmente dagli astronomi, sul confronto fra le sue conseguenze e l'esperienza. Ad ogni modo, nel Suo M.S., non era che accennato.

*[segue parte cancellata da Maggi. A fianco, in penna rossa si legge:* Soppresso*:*
Sull'espressione "teoria agnostica", sta bene, se s'intenda che dal sistema K non si possa ottenere <u>direttamente</u>, cioè <u>coll'osservazione diretta</u>, la misura della lunghezza e dei tempi appartenenti a K', quali risultano coll'osservazione diretta istituita nello stesso K'. Ma la parola <u>eterogenee</u> alle entità spaziali e temporali appartenenti ai due diversi riferimenti non mi sembra opportuna, perché sono tutte lunghezze e tutti tempi, e tutte misure di lunghezze e di tempi. La parola <u>eterogenee</u> è adoperata per <u>inconfrontabili</u>, che mi sembra tradurre il suo stesso concetto che le misure, quali risultano ad un osservatore nel proprio sistema, sono direttamente inaccessibili all'osservatore appartenente all'altro sistema.*]*

<u>Venendo</u> <u>infine</u>. Il <u>tempo</u> <u>proprio</u> non sottrae a di qualunque specie obbiezione, prendendolo per quello che risulta dalla sua definizione mediante $dt^{2*}=dt^2-\frac{dx^2+dy^2+dz^2}{c^2}=\left(1-\frac{q^2}{c^2}\right)dt^2$, che significa una variabile, avente le dimensioni del tempo, atta a individuare le posizioni di un punto mobile sopra una linea d'universo, invariante rispetto alla scelta di un riferimento piuttosto di un altro. Così si esplica il suo ufficio, nella teoria relativistica del campo elettromagnetico, e non occorre cercar altro. Per aver un orologio, si ricorrerà a K, e ci darà t, si ricorrerà a K', e ci darà t', con che t* resterà invariato. Ma t* non può esser dato, salvo il caso di un movimento del punto rettilineo uniforme, da un orologio qualsiasi; perché un così detto orologio, che coincide, a successivi istanti, con orologi diversi, il cui ritmo varia coll'istante, subordinatamente al movimento del supposto punto mobile, non può considerarsi come un orologio, implicando il concetto d'orologio quello di un assegnato ritmo costante.

Dopo questo, io non saprei più rendermi conto esattamente della sua dimostrazione dell'inconsistenza della teoria dell'Einstein. Ma, comunque sia, difficilmente io potrei presentare all'Istituto Lombardo, né altrimenti a qualche rivista, uno scritto avente per scopo tale dimostrazione. Perché di tale inconsistenza dovrei cominciare a persuadermi io stesso, mentre, come ho fatto le eccezioni che Le ho scritto alle Sue obbiezioni, non consento con altre, e, a parte possibili riserve di tutt'altro genere, io non riconosco, nella dottrina einsteiniana, nessuna menda, sotto l'aspetto della logica consistenza. Alcuni appunti, come la storia dei gemelli, non sono imputabili alla dottrina, ma ad una indebita interpretazione dei suoi concetti.

RinnovandoLe l'espressione dei più cordiali sentimenti, La prego di credermi sempre

Aff$^{mo}$ Suo

G.A. Maggi.

## 42
## Gian Antonio Maggi a **Francesco Crotti**

Sui principii del Calcolo degli Errori
Lettera all'Ing. F. Crotti a proposito del suo scritto
"Sulla compensazione degli errori nei rilievi geodetici"
(Atti del Collegio degli Ing. di Milano 1886) - Messina - Marzo 1887

Vedo che stai preparando per la stampa un Manuale sul Metodo dei Minimi Quadrati: e lo vedrò con molto interesse. Opera del tuo ingegno e del tuo fine spirito critico, questo riuscirà certamente un ottimo libro: e sono persuaso che, fondando la dimostrazione della maggior plausibilità delle note medie col "principio dell'imparzialità", senz'alcun dubbio assai persuasivo, la trattazione riunirà al rigore la semplicità indispensabile allo scopo che hai di mira. Ma in via assoluta, è realmente da ripudiarsi la teoria degli errori, secondo la quale la maggior plausibilità delle medie stesse riesca una conseguenza del fatto che a quei valori corrisponde la massima probabilità di essere la condizione per cui si ottiene il sistema di valori osservati effettivamente ottenuto? La maggior probabilità di essere il vero è certo il criterio ideale della maggior plausibilità di un'ipotesi!... E, innanzi tutto, siamo d'accordo che se l'ipotesi più probabile è necessariamente reputata come la più plausibile, la maggior plausibilità, come può essere stabilita dietro criterii diversi, non trae seco per nulla la maggior probabilità: circostanza che molti, e fra questi il Ferrero,[75] lasciano credere di non veder punto chiara, e parrebbe davvero impossibile!

Ora qual è il fondamento su cui posa la dimostrazione che, dal precedente punto di vista, la media aritmetica dei risultati, $a_1, a_2, ..., a_n$ di n misure d'una quantità è il valore a cui compete la massima probabilità di essere il vero? Stando a qualche tua proposizione, questa questione eccederebbe i limiti della matematica, per invadere il campo della metafisica. Pure non eccede senza dubbio quei limiti il teorema di Bayes, secondo il quale se un evento può essere prodotto da n cause, e la causa $q^{ma}$ ha la probabilità $Q_i$ di agire, mentre, supposto che agisca, $P_i$ è la probabilità che sorta l'evento, avvenuto che sia l'evento medesimo, la probabilità che sia avvenuto per opera della causa $C_r$ è data da

$$\frac{P_r Q_r}{\sum_i P_i Q_i}.$$

Da questo teorema segue immediatamente che se $\varphi(x)dx$ rappresenta la probabilità che un errore sia compreso tra x e x+dx, per modo che, la probabilità di commettere n errori compresi fra $x_1, x_2, ... x_n$ e $x_i+dx_i, x_2+dx_2, ... x_n+dx_n$ sia

$$\varphi(x_1)\varphi(x_2) ... \varphi(x_n)dx_1dx_2 ... dx_n$$

[75] A. Ferrero, *Esposizione dei minimi quadrati*, Barbera, Firenze, 1876.

ottenuto il sistema d'errori $x_1$, $x_2$, ... $x_n$, corrispondenti ai risultati $o_1$, $o_2$, ..., $o_n$, la probabilità che ne sia stata condizione l'esser il valore vero compreso fra v e v+dv (del qual evento la probabilità a priori è proporz. a dv) sarà

$$\frac{\varphi(v-o_1)\varphi(v-o_2)...\varphi(v-o_n)dv}{\int_{-\infty}^{\infty}\varphi(v-o_1)\varphi(v-o_2)...\varphi(v-o_n)dv}.$$

E supposto

$$\varphi(x)=\frac{h}{\sqrt{r}}e^{-h^2x^2},$$

si troverà subito che questa quantità diventa massima nell'ipotesi

$$v=\frac{o_1+o_2+...+o_n}{n}.$$

Tutto sta dunque nell'ammettere che la probabilità di commettere un errore compreso fra 0 e x si possa rappresentare con una certa funzione, di cui $\varphi(x)dx$ sarà il differenziale, e che abbiasi

$$\varphi(x)=\frac{h}{\sqrt{r}}e^{-h^2x^2}.$$

Qui sta il busillis; e convengo che in parecchi trattati la deduzione di questo risultato è fatta in modo da guadagnarsi ben poca fiducia.

A mio giudizio, ella va egualmente appoggiata al calcolo e all'esperienza: col calcolo, stabilite alcune ipotesi sulla natura degli errori considerati, fabbricando una formola che all'esperienza spetta di verificare, e coll'esperienza provando che la formola stessa si verifica entro limiti accettabili. L'espressione

$$\frac{h}{\sqrt{r}}\int_{-x}^{x}e^{-h^2x^2}dx$$

della probab. di commettere un errore in valore assoluto inf. ad x, riuscirà così qualcosa come una formula empirica, valida per una certa categoria d'errori, ed, entro questi limiti, non vedrei come si potrebbe ragionevolmente ripudiare.

Quella categoria d'errori è caratterizzata innanzitutto dalla proprietà che è eguale la probab. di commettere un errore compreso fra 0 e x e fra 0 e −x. Su questo postulato fondamentale tu trovi da ridire: ma badiamo che è la pura e semplice definizione degli errori di cui noi intendiamo occuparci. Riesce poi dimostrato a posteriori che almeno in molti casi, e con grandissima approssimazione, quest'ipotesi risponde alla realtà: poiché, prendiamo una delle tante verificaz. sperimentali, fatta la media delle n osservazioni, e sottratte le singole osservazioni da essa, è un fatto che il rapporto del numero di queste diff. il cui valore assoluto cade fra 0 e x, al numero totale è rappresentato assai prossimamente da $\frac{h}{\sqrt{r}}\int_{-x}^{x}e^{-h^2x^2}dx$: e questa coincidenza che, in quella ipotesi, si deve aspettare, non si saprebbe come spiegare con un'ipotesi diversa, in conseg. della quale la media di un numero grand. di oservaz. non sarebbe più con una probab. osservaz. al loro numero, il valor vero, o la probab. di commettere un errore compreso fra 0 e x non sarebbe più rappr. da quella funzione. Di queste verificazioni ne feci anch'io parecchie, con risultati di misure eseguite con diversi strumenti, e il risultato è invariabile.

Poi quanto ai ragionamenti che conducono a priori alla $\frac{h}{\sqrt{r}}\int_{-x}^{x}e^{-h^2x^2}dx$, il migliore mi sembra quello (meno conosciuto) di Encke, secondo il quale questa è una formola approssimata, valida quando l'errore è prodotto da un numero grandissimo di cause, che introducono altrettanti errori parziali, non eguali (come nel ragionamento di Hazen) ma dello stesso ordine. Ma colla riserva di confrontare il risultato coll'esperienza, credo accettabile anche la dimostraz. di Hazen, o, quando si voglia andare per le spicce, non mi pare nemmeno da condannare l'assumere addirittura per rappresentare la probab. di commettere un errore compreso tra x e x+dx la funz $ce^{-h^2x^2}$, che possiede, alcune esattamente, altre per approssimaz., le principali proprietà che giova concedere a priori alla funz. cercata.

A questo proposito non occorrerà partire da una supposta $e^{a_0+a_1x+a_2x^2+a_3x^3+a_4x^4+\dots}$: ma, nel caso, ammesso che stia la φ(x)=φ(−x), e che una stessa formola debba valere per x positive e negative, i coeff. delle potenze dispari di x dovranno eguagliarsi a 0.

Tu chiami erronea la deduz. che φ(x) sia pari dalla proprietà φ(x)=φ(−x). Ora, una funz. si dice pari per ciò solo che possiede questa proprietà. Nel caso da te citato la formola cambia nei due casi di x pos. e di x negativo, e perciò la funzione pari può essere rappresentata da una funzione lineare di x. Supposto invece che pei due casi debba valere la stessa formola, se questa sarà un pol. in x, o una serie di potenze di x, è patente che non potranno figurarvi che le potenze pari. Osservai anche che tu consideri del momento del trave soltanto la grandezza. Tenuto conto del verso della rotaz. quel momento sarà una funz. dispari della distanza x *[*della reg.cons.*]* dall'azione media, e la sua espressione, *[*per ogni lato*]*, λ−x, dove λ rappr. l'ascissa di quell'estremo del trave che si trova, in confronto della reg. media, dalla stessa parte della reg. considerata.

Mi permetterò di fare anche qualche altra osservazione.

Alla citazione della mia Nota soggiungi che, infine, dovetti io pure rassegnarmi a introdurre per la precisione l'ipotesi più probabile. Osservo a questo proposito che si è obbligati a prendere pel metodo di precisione, o pel così detto errore medio, legato con esso dalla nota relazione, in valore di ripiego, quando per calcolarlo non si dispone che dei risultati di un numero ristretto d'osservazioni. Allora, come $\frac{|0|}{n}$ è il valore più probab. della quantità misurata, $\frac{|x^2|}{n}$ è il valore più probabile di $\int_{-\infty}^{\infty} x^2e^{-h^2x^2}dx$. Ma gli errori x non son noti: solo si possono avere gli scostamenti s dal valor più probabile, e si trova pel quadr. dell'errore medio il valor $\frac{|s^2|}{n}$ assumendo per differenza fra il valor vero e la media aritmetica l'errore medio della media. Questo non è che uno spediente qualunque: e un valore attendibile dell'error medio non si avrà che documentando un numero molto grande d0osservazioni per modo che gli scartamenti con probabilità proporzionata al numero stesso rappresentino i corrisp. errori. Una volta poi calcolato, si potrà valersene per ogni altra osservaz. fatta collo stesso strumento, presumibilmente colla stessa precisione.

Finalmente, dici, a proposito di un'osservazione del Forti che un teorema di calcolo non può esser dimostrato che in un numero inf. di osservaz. il numero delle volte che un evento si riproduce è proporz. alla sua probab. Ora è si dimostra effettiv. che la probab. che in p prove, in ciascuna delle quali può sortire o un evento o l'evento contrario, il rapporto del numero delle volte in cui l'evento ha luogo al numero totale differisce dalla prob. che compete all'evento meno di un numero prefissato piccolo finché si vuole, si può andare mantenere prossimo finché si vuole all'unità prendendo il numero delle prove superiore ad un certo limite. Questo è il celebre teorema di Bernoulli, che fornisce una base razionale alle verificaz. aprimentali di cui sopra.

Così al Forti, in questo punto, si può rimproverare tutt'al più d'aver enunciato tale teorema alquanto alla carlona. Invece difficilmente si potrà difendere lo stesso Forti dall'essere caduto nel più vizioso dei circoli, colla sua verificazione a posteriori che la media è il valore più prossimo al vero. E censurabile è pure questa locuzione di valore più prossimo al vero, invero di valore più probabile, che è come se si dicesse che una persona che s'incontra per via a Messina è più siciliana di una che s'incontra a Milano. In quella locuzione vi è tutto un falso punto di vista: quello cioè che il Calcolo degli Errori possa servire in certo qual modo a correggere il risultato dell'esperienza, mentre solo può indicare la probabilità che sia affetto da un errore inferiore ad *[un]* certo limite.

**43**[76]

Francesco Crotti a Gian Antonio Maggi

[carta intestata: Società italiana per le strade ferrate del Mediterraneo
Società Anonima con Sede in Milano.
Capitale Sociale L. 180 milioni, interamente versato.
Servizio del Mantenimento, Sorveglianza e Lavori - Direzione]

Milano, il 28/11 1892

Caro Gian Antonio

Scusa se ti disturbo, ma bisogna farla finita con quella vera scempiaggine, tuttora trionfante nelle scuole e sui libri, che è la formola per il calcolo dell'area compresa tra una retta ed un perimetro di cui si conosce solo i valori di una serie di ordinate equidistanti. Trattasi della cosiddetta formola Simpson

$$A=\frac{h}{3}(y_1+4y_2+2y_3+4y_4+\ldots)$$

L'errore logico da cui si cava quella formola, è proprio atroce, poiché l'area del trapezio compreso fra $y_{2r}$ ed $y_{2r+1}$ si calcola ritenendo che il perimetro corrisponda a quello della parabola passante per i punti (2r−1), (2r), (2r+1). Invece è logico di ripetere il calcolo anche a mezzo della parabola passante per (2r) (2r+1) (2r+2); poi prendere la media dei due risultati.

La formola di Simpson è proprio enorme. Calcola l'area fra $y_1$ ed $y_{21}$ essendo

h=30 $y_1=y_3=y_5=\ldots=y_{21}=20$

$y_2=y_4=\ldots=y_{20}=14$

otterrai, colla Simpson, A=9600

Calcola ora l'area tra $y_1$ ed $y_{20}$ e risulta A=9720

Ossia <u>la parte maggiore del tutto</u>.

Che i nostri grandi matematici si occupino di spazi di n dimensioni, di circoli immaginari all'infinito, delle 27 rette della superficie del 3° ordine ecc. sta bene; ma che si chiudano in un altezzoso silenzio quando un galantuomo si lagna perché nelle scuole e sui libri si insegnano corbellerie come quella espressa dalla formola di Simpson, è cosa che forse passa il segno. A qualche cosa li obbliga la fama che si spappolano colle loro elucubrazioni che nessuno ha mai lette.

Tu dirai che sono un poco irritato. Lo sono ed a ragione. Fino dal 1885 ho scritto e dimostrato la verità sulla formola di Simpson. Che ne fu la conseguenza? Il più dignitoso silenzio da parte dei nostri grandi. Ma che, avranno detto, un profano che sa appena le prime quattro operazioni, deve egli osare di mettere in dubbio il valore di quelle perle che escono dalla nostra bocca da semidei!

[76] Le lettere che seguono sono state raccolte dallo stesso Maggi in un fascicolo dal titolo: *Corrispondenza coll'Ing. Crotti sulla Formola di Simpson. 1892-93.*

Ah semidei... a momenti vi mando a scuola si senso comune.

Tu invece sei e sei sempre stato valente ma onesto. Tu capirai che quello che dico sulla formola Simpson non è da spregiare affatto.

Ieri al Collegio Ingegneri di Milano ho letto l'adesione scritta, per la condanna della formola di Simpson, nientemeno che dei professori Colombo, Formenti, Jorini e di 20 Ingegneri.

Unisco foglietto di spiegazioni.

A parte lo sfogo fatto con te che devi perdonare desidero molto molto una tua risposta

Tuo aff. Crotti Francesco Milano Via Dante 14

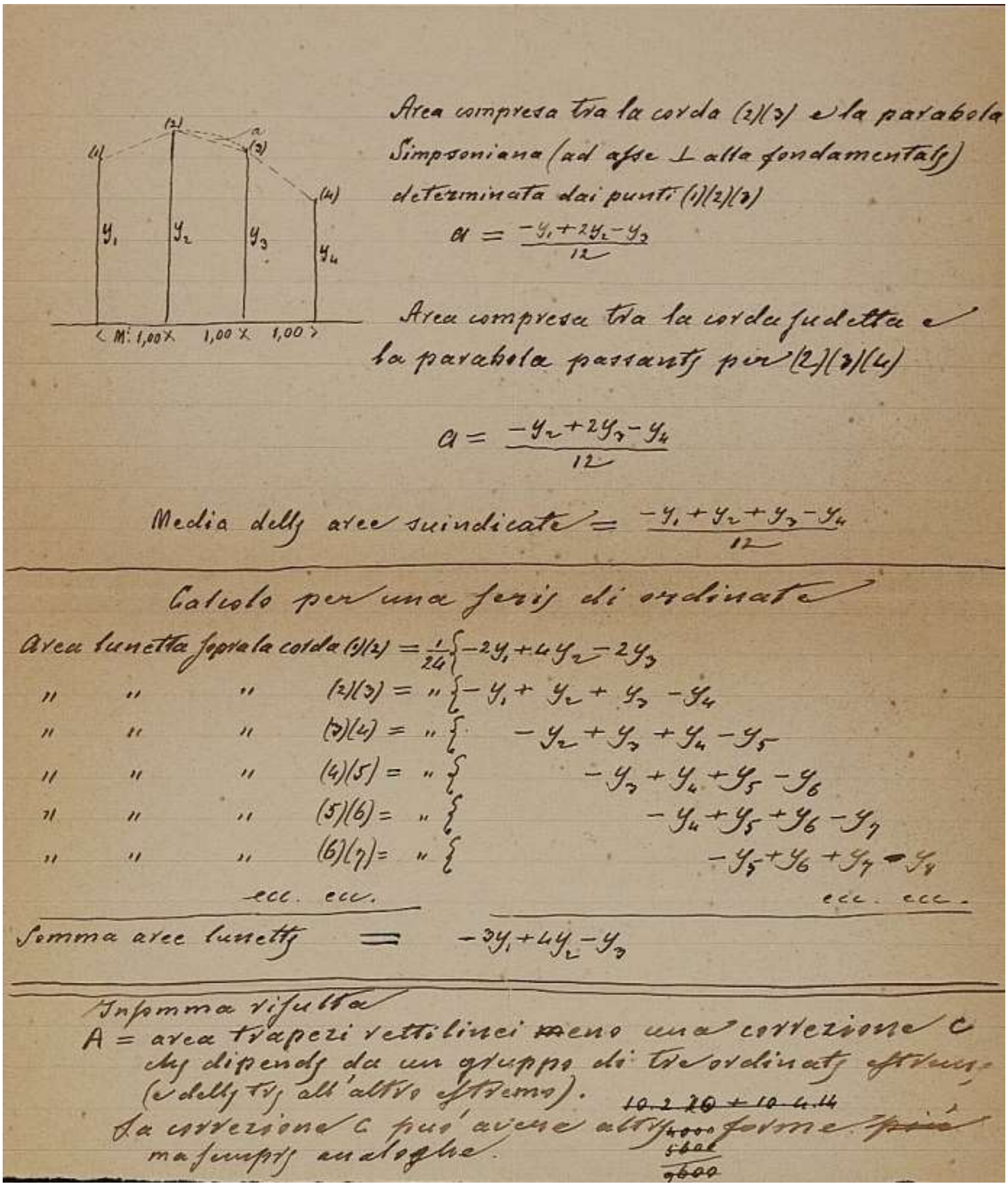

Area compresa tra la corda (2)(3) e la parabola Simpsoniana (ad asse ⊥ alla fondamentale) determinata dai punti (1)(2)(3)

$$a = \frac{-y_1+2y_2-y_3}{12}$$

Area compresa tra la corda suddetta e la parabola passante per (2)(3)(4)

$$a = \frac{-y_2+2y_3-y_4}{12}$$

Media delle aree suindicate $= \frac{-y_1+y_2+y_3-y_4}{12}$

Calcolo per una serie di ordinate

Area lunetta sopra la corda (1)(2) $= \frac{1}{24}\{-2y_1+4y_2-2y_3$

" " " (2)(3) = " $\{-y_1+y_2+y_3-y_4$

" " " (3)(4) = " $\{-y_2+y_3+y_4-y_5$

" " " (4)(5) = " $\{-y_3+y_4+y_5-y_6$

" " " (5)(6) = " $\{-y_4+y_5+y_6-y_7$

" " " (6)(7) = " $\{-y_5+y_6+y_7-y_8$

ecc. ecc. ecc. ecc.

Somma aree lunette $= -3y_1+4y_2-y_3$

Insomma risulta

A = area trapezi rettilinei meno una correzione c che dipende da un gruppo di tre ordinate estreme (e delle tre all'altro estremo).

La correzione c può avere altre forme ~~più~~ ma sempre analoghe.

**44**

Gian Antonio Maggi a Francesco Crotti

Messina 4/12 92

Carissimo amico

Sono d'accordo con te che allo spaziare troppo lontano in certi campi è preferibile restare modestamente a terra a curare i fondamenti: e che questo lavoro non è sempre debitamente rimunerato. Perciò volentieri rispondo al tuo appello: ma - non te l'avrai a male - per spezzare una lancia a favore di questa maltrattata formula di Simpson.

Non dimentichiamoci che la formula di Simpson è formula d'approssimazione: e non meravigliamoci che una formula d'approssimazione possa condurre anche all'assurdo, prendendola come esatta, senza preoccuparsi che si verifichino le circostanze occorrenti perché l'errore non oltrepassi certi limiti e diventi sproposito

Entro in medias res. Immaginiamo la successione degli archi di parabola coll'asse perpendicolare all'asse delle x passanti per le terne (2r−1), (2r), (2r+1) con r compreso tra 1 e 10, e la successione degli archi analoghi passanti invece per le terne (2r), (2r+1), (2r+2), con r compreso tra 1 e 9. Il numero 9608 che fornisce la formula di Simpson applicata ad un trapezio curvilineo compreso fra le ordinate $y_1$ e $y_{21}$ non è altro che l'area del trapezio la cui curva è la prima successione, e il numero 9720 che ti fornisce applicandola alla parte di quel trapezio compresa fra le ordinate non è altro che l'area del trapezio la cui curva è la seconda successione. Se il secondo numero è maggiore del primo, vuol dire che la somma delle aree delle figure *[tratteggiate in verticale]* è maggiore della somma delle aree dei due trapezii *[tratteggiate in orizzontale]* (V. fig.).

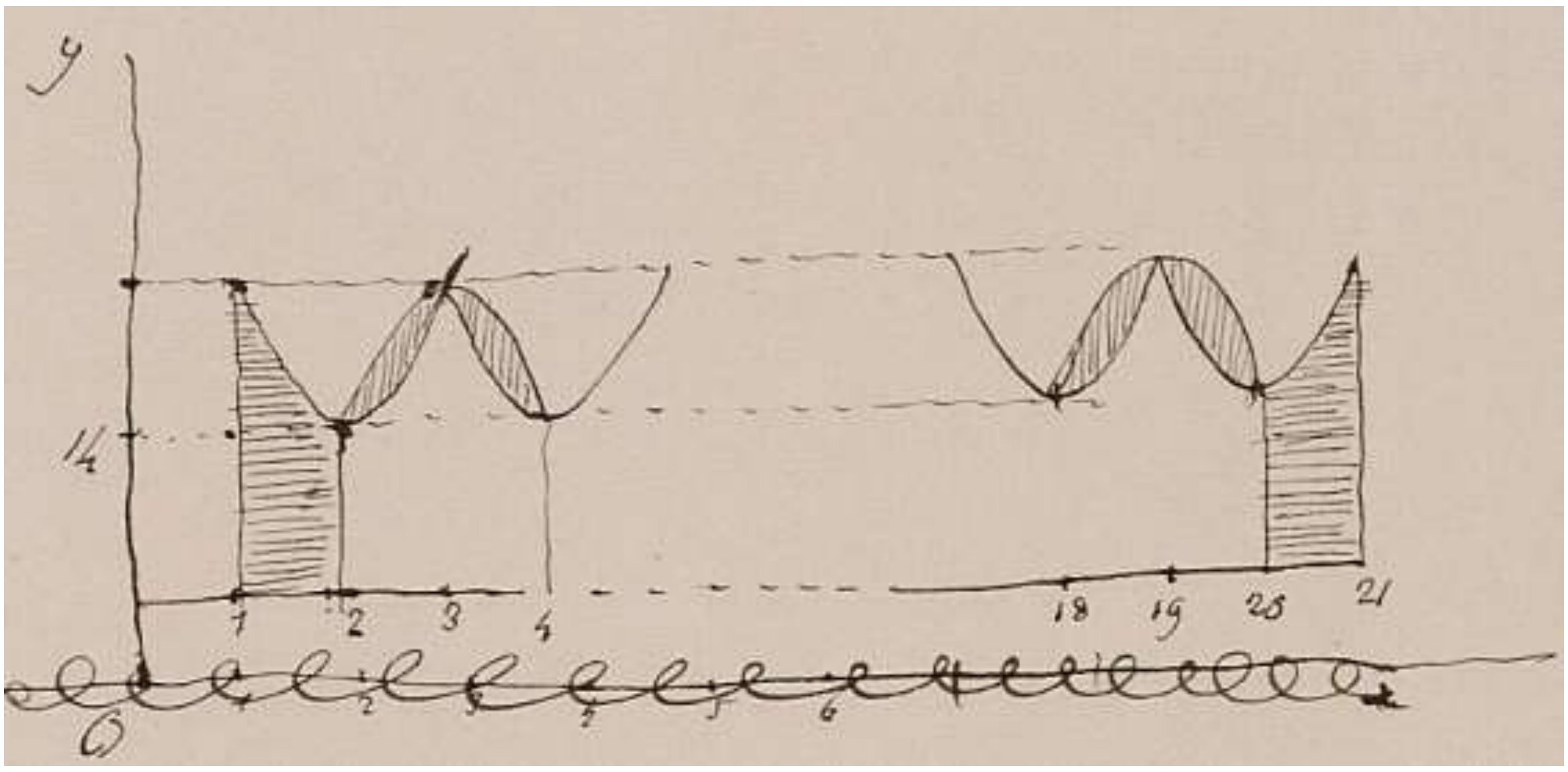

e di questo non c'è ragione di meravigliarsi; che se non ne inferiamo la conclusione meravigliosa che la parte è maggiore del tutto, la colpa non mi par proprio della formola di Simpson. Immaginando in fatti una curva qualunque alla quale appartengano i 21 punti del tuo esempio, e il relativo trapezio curvilineo, 9600 e 9720 non saranno in generale le misure dell'area né del tutto né della parte, ma valori più o meno approssimati a seconda della maggiore o minore adesione pel tutto, della curva alla prima successione d'archi parabolici, e, per la parte, del tratto di quella curva compresa fra i punti 2 e 20 alla seconda successione.

Non è difficile riconoscere che quanto più debolmente la curva aderirà alla prima successione, tanto meno aderirà quella sua parte alla seconda: lo vediamo nel caso che la curva sia la stessa prima successione: ce lo può mostrare ad *[*occhio e croce*]* la figura: ce lo dice la misura dell'area delle figure *[tratteggiate in verticale]*. E perciò, coi dati del tuo esempio, la formola di Simpson se si applica bene a dovere al tutto non si applica a dovere alla parte, e viceversa. Questa è la spiegazione del paradosso: o non credo più alla Matematica. Bisogna poi che il trapezio curvilineo, che si vuol misurare, sia effettivamente dato, per decidere se si dà un caso o l'altro: mentre potrebbe anche avvenire che quei dati fossero insufficienti per conseguire una discreta approssimazione in un caso e nell'altro.
Noto, per quanto al principio, che la deformazione approssimata di $\int_{\alpha}^{\beta} f(x)dx$ colla formula di Simpson si riduce puramente a concepire l'intervallo (αβ) diviso in un numero pari di intervalli eguali, e a prendere nelle successiva coppie di questi intervalli $a+bx+cx^2$ in luogo di f(x), determinando in ogni doppio intervallo le costanti a, b, c colla condizione che sia $a+bx+cx^2=f(x)$ negli estremi e nel centro: procedimento, in sostanza, anche per la rappresentazione delle leggi fisiche, così largamente applicato.

Non parlo d'opportunità pratica: questione nella quale mi dichiaro del tutto incompetente.
M'interessa veramente di sapere i capi d'accusa imputati alla nostra formola da tanti valentuomini. E non dubito che le tue osservazioni mi procureranno presto il piacere di ricevere altre tue lettere.

Comunque la formola di Simpson possa valere, le saremo obbligati d'aver riattivato la nostra corrispondenza.

Voglimi sempre e credimi sempre

Aff tuo

**45**

Francesco Crotti a Gian Antonio Maggi
[carta intestata: Società italiana per le strade ferrate del Mediterraneo
Società Anonima con Sede in Milano.
Capitale Sociale L. 180 milioni, interamente versato.
Servizio del Mantenimento, Sorveglianza e Lavori - Direzione]

Milano, il 28/12 1892

Caro Maggi

Ho ricevuto la carissima tua e te ne ringrazio vivamente.

Ti spedisco sottofascia il mio opuscolo del 1885.[77]

Ora il Collegio Ingegneri di Milano ha votato un ordine del giorno in cui delega al suo Comitato la nomina di una commissione di matematici, per esaminare la questione da me proposta. È dunque una cosa seria.

A proposito, il Peano di Torino mi ha dato ampia ragione.

E perché il mio Maggi, il mio valente, valentissimo Maggi, uno dei matematici che io stimo di più, mi tien broncio?

---

[77] F. Crotti, "Sopra alcune formole di planimetria e stereometria", *Il Politecnico-Giornale dell'ingegnere architetto civile ed industriale,* 1885, v. 17, pp. 193-207.

Pensaci su quieto 5 minuti e basteranno.

T'ho già detto che il termine di conversione può avere varie forme. Oltre alla

$$\frac{-y_1+y_2}{12} + \text{l'analoga all'ultimo estremo}$$

si può anco arrivare alla

$$\frac{-3y_1+4y_2-y_3}{24} + \text{l'analoga.}$$

Parlando d'altro chi sa che, se vien soppressa codesta Università, ti mandino nell'Alta Italia. Ci rivedremmo più di spesso.

E la tua gentil signora e i figli? Ti prego di presentare i miei omaggi.

Tornando al chiodo fisso, non venire a dirmi che di pratica non te ne interessi. Ma ti intendi di logica, e come!

Dunque pensaci e scrivi al tuo buon amico che ogni tanto dà un <u>bui</u> come diciamo a Milano.

Ciao, ciao

aff^mo Crotti

**46**

Gian Antonio Maggi a Francesco Crotti

Messina 26/1 93.

Carissimo!

Io non ti tengo certamente il broncio, tant'è vero che alla tua lettera del 28 Novembre ho risposto subito e a lungo, e ora rispondo all'altra tua del 28 Dicembre, appena mi dà tregua la bufera che mi portò alla Capitale. Ma, per quanto lo desideri, non posso darti ragione: neppure dopo letto il tuo Articolo.

Tu discorri come se la curva del trapezio fosse ignota, e dalle ordinate date dovesse dedursi un valore plausibile dell'area. Non saprei davvero come tale problema si potrebbe porre; ad ogni modo non è quello cui si riferisce la Formola di Simpson. La curva si deve intendere data, e misurate tante ordinate quante bastano per ottener l'approssimazione desiderata: ciò che equivale a non discostarsi oltre certi termini la successione degli archi parabolici dalla curva medesima. Perché ci preoccuperemo che le ordinate pari siano trattate diversamente dalle dispari?

Io sono disposto a continuare la discussione, purché tu risponda alle mie osservazioni. Che cosa dicono i matematici che mi hai nominato? Ho veduto che il Peano pubblicò il tuo quesito nella Rivista, ma si è pronunciato?[78] Per ora mi limito a notare che, coi dati di quel quesito *[*della*]* prima lettera, supposto che la curva sia precisamente la successione degli archi parabolici, 9688 è <u>esattamente</u> l'area del trapezio, mentre colla correzione da te proposta si troverebbe 10170.

*[in calce: calcoli del Maggi]*

[78] Il quesito di Crotti si trova in: *Rivista di Matematica*, v. II, 1892, p. 176. Sul volume successivo della *Rivista* c'è la risposta di Peano sul medesimo argomento alle pp. 17-18 e i commenti di Jadanza (pp. 15-16) e di Bardelli (pp. 16-17).

**47**
Francesco Crotti a Gian Antonio Maggi

Caro Gian A. Maggi

Milano 30/1 93

Con grande piacere ho ricevuto la tua. Nel subbuglio Messinese ho seguito, sui giornali, le tue vicende. A parte le simpatie che puoi avere per codesta nobile città, so quanto è in te elevato il concetto della dignità degli studi, e che per conseguenza ti piaceranno Università poche ma provvedute bene.

Venendo alla nostra tesi, tu mi scrivi: “Tu discorri come se la curva del trapezio fosse incognita, e dalle ordinate date dovesse dedursi un valore plausibile dell’area.” ecc. ecc.

Precisamente! È questo il caso. Quando si fa un rilievo sul terreno si misurano solo parecchie ordinate, e nulla si sa né si vuol sapere della natura delle linee di perimetro. Ci vorrebbe altro!

Ciò posto, siamo perfettamente d’accordo.

A parlare ci si intende.

Soffri dunque che io ti ponga fra gli egregi che si associano a me nel dire che la formola di Simpson è male applicata.

Ho fatto estratti di 20 autori che insegnano essere quella formola applicabile a toute courbe irrégulière.

Hanno torto.

Intanto ti saluto caramente e mi dico

Tuo aff
Ing. Crotti

**48**
Gian Antonio Maggi a Francesco Crotti

Messina 8 Febbrajo 93

Carissimo amico!

Io non desidero di meglio che intendermi con un caro amico come mi sei tu, e accrescerebbe, se fosse possibile, questo desiderio la gentilezza delle tue espressioni. Ma, per amore della verità, a costo di compromettere un’apparente *[sic!]* accordo, bisogna che ritorni sul nostro argomento.

Coloro che applicano la formula di Simpson senza preoccuparsi della natura della curva del trapezio, l’applicano male, verissimo! Ma non hanno torto per questo quelli che affermano che la formula di Simpson si applica ad ogni curva, perché c’è bisogno di dichiarare che bisogna applicarla a dovere?

Se gl’Ingegneri, quando fanno un rilievo, non si preoccupano della natura della linea di perimetro, ecco un fatto pratico ch’io imparo adesso! T’ho detto che di pratica non m’intendo! Ma, se m’intendo un po’ di logica, dovrebbero darsi quella briga, quanto basta (e non sarà gran cosa) per, scegliere la formola d’approssimazione più appropriata al caso, e per misurare poi le ordinate che occorrono affinché l’errore non oltrepassi certi limiti.

Per quanto poi al problema generale del trovare il valore più plausibile d'un'area del cui perimetro non sono dati che alcuni punti, come ti dicevo nell'ultima mia, non vedo come si possa porre perché quei punti potranno appartenere egualmente ad ogni linea, finita ed anche infinita, che non sia determinata da un numero di punti minore, ed io non so vedere un criterio razionale per sceglierne una e scartare le altre.

Avrai veduto nella Rivista del Peano la risposta dell'Jadanza. Così egli come l'altro ragionano nell'ipotesi che la curva sia data; e neppur io, senza il tuo Articolo, avrei immaginato che si pensasse di farne a meno.

Non mi hai mai detto che cosa affermano il Colombo, il Formenti e gli altri, procura di farne materia d'una prossima lettera.
Intanto ti stringo cordialmente la mano, e ti prego di credermi sempre

Aff$^{mo}$ tuo
G.A Maggi

**49**
Francesco Crotti a Gian Antonio Maggi

Milano 3/3 93

Caro Maggi,

Sottofascia avrai ricevuto un mio opuscoletto con un'Appendice... Sicuro, t'ho citato fra coloro che sostengono non me ma la verità.[79]

Tu hai fatto un poco l'avvocato del diavolo! Ma la tua convinzione, che la famosa formola a coefficienti alternati sia la più sgraziata cosa che i teorici hanno regalato alla pratica, traspira chiaramente dalle tue lettere.
Intanto, oltre le adesioni di già citate nell'opuscoletto, abbiamo già una stampa del Prof Bardelli che, pur mostrando di non avermi capito interamente, confessa essere deplorabile che in libri reputati si ammanisca *[sic!]* quella formola come valevole per i casi generali.

Da questo solo la mia campagna viene più che giustificata.
Quel buon Jadanza! Senza la falsariga della scienza succhiata nelle scuole, perde subito le staffe. Ma io ho la testa dura.
A proposito, Jorini mi disse che la Commissione nominata dal Collegio ha discusso tre ore senza andare d'accordo. Buon segno; la discordia è nel campo d'Agramante. Jorini mi disse anche che in Francia Levy ed un altro di cui non ricordo il nome, hanno pure combattuto la formola in discorso, giungendo alle conseguenze a cui sono giunto io.

Dunque spazziamo via l'errore.

Scrivi qualche cosa tu. Voi teorici di acuto ingegno ve ne state un poco nelle nubi e lasciate scorazzare parecchie corbellerie nel campo delle applicazioni.

Leggi il mio articolo sulle botti. È un argomento allegro anzi di...vino. Scusami della noja che ti avrò data e credimi

Tuo aff.
Ing. Crotti

[79] F. Crotti, "Alcune considerazioni sulla recente edizione (12$^a$) del Manuale del Prof. Colombo", *Il Politecnico - Giornale dell'ingegnere architetto civile ed industriale*, 1893, v. 25, pp. 229-242. Nella prima nota a piè di pagina c'è la citazione di cui parla Crotti.

**50**
Gian Antonio Maggi a Francesco Crotti

Messina 17 Marzo 93.

Carissimo!

Non so proprio capire in che modo che mi abbia potuto citare come pervenuto a conclusione analoga a quella che alla formola di Simpson manca ogni fondamento teorico logico. Le mie lettere (le ho rilette) sono pur chiare. Ma tu mi obblighi a dire che, per leggere fra le righe, non leggi le righe. Mi spiace anche perché ora non posso esimermi dal pubblicare la mia opinione, che a quella è direttamente opposta. A far conoscere com'io penso non tengo gran fatto: ma altra cosa è lasciar credere ch'io pensi il contrario.[80] Divisi, assolutamente divisi, su questo campo, restiamo sempre buoni amici: e con questo sentimento ti do una stretta di mano e mi affermo

**51**
Francesco Crotti a Gian Antonio Maggi

Milano 19/3 93

Caro G.A. Maggi

Alla tua gradita lettera rispondo che forse ho sbagliato citandoti, anzi senza forse vedo che ho sbagliato.

Smentisci pure l'opinione che ti ho attribuito, anche io da parte mia nella prossima Seduta del Collegio farò la rettifica.

Io annettevo e annetto ancora grande importanza al tuo voto; ma in linea obbiettiva lo scopo a cui tendevo è raggiunto lo stesso.

Già Poncelet, nientemeno che Poncelet, ha rilevato l'inesattezza assoluta teorica e pratica della formola che tu difendi.

Fu combattuta da Parmentier (Nouvelles annales de mathematique oct 1855)[81]; dal celebre Culmann (Traitè de Statique grafique, Dunod 1880 pag 103). Il nostro Colombo l'ha già portata nella nuova edizione del Manuale. La Commissione del Collegio ha già formulato un voto analogo a quanto io sostengo.

Che siano proprio i matematici puri i più refrattari a ragionare con un poco di buon senso? Non lo posso credere, specialmente quando fra i matematici puri vi è un Maggi, di cui tanto stimo il forte ingegno.

Sempre buoni amici del resto

Tuo aff.
Crotti

[80] Si confronti la corrispondenza con Giuseppe Peano: la lettera di Maggi a lui indirizzata - del 6/3/1893 - sarà pubblicata con il titolo "(Sulla formola di Simpson). Estratto da una lettera al Prof. Peano" sulla *Rivista di Matematica edita da G. Peano*, v. III, pp. 60-61.

[81] T. Parmentier, "Comparaison de quelques méthodes de quadrature et formule nouvelle", pp. 370-384.

**52**

Francesco Crotti a Gian Antonio Maggi

Milano 4/1 94

Caro GianAntonio

Non puoi imaginarti con quanto piacere ho ricevuto il tuo caro, carissimo biglietto. Quante volte ho pensato a te, alla tua gentil signora, alla tua famiglia sentendo le terribili notizie del terremoto.[82] Come avrai dovuto trepidare per i tuoi cari!

Basta! Ora almeno è finito lo spaventevole fenomeno? Pare di sì, sebbene si abbiano ogni tanto notizie di commozioni quà *[sic!]* e là. Auguriamoci dunque che non se ne parli più.

Tu come stai? Figli quanti? Desidero tanto tue notizie. Già la gran bontà de' cavalieri antiqui, come tu spiritualmente dicesti, fa che delle bizze scientifiche nulla riverberi sulla vera e forte amicizia che passa tra noi. Ricordo i bei giorni di Pavia in cui tu mi hai messo un poco in giornata colla crudele mia bella matematichetta. La metto al diminutivo, poiché né la testa né il tempo mi permisero di fare studi completi come tu hai potuto fare.

E posto che il terremoto come fa si tace non posso non esimermi, (abbi pazienza, tu col tuo gran cuore puoi compatire un tantin di riscossa) dal dirti che in quella tal tenzone ebbi piena ragione da un verdetto di una Commissione nominata dal Collegio Ingegneri e presieduta dal Jorini che si consigliò in proposito anche col Brioschi. Che la relazione della Commissione fu poi stampata e che sulla 13ª edizione del Manuale del Colombo, fu dichiarata la nuova formola preferibile a quella di Simpson. Dunque, non avevo poi tutti i torti.

Ma bando a ciò. Vedo di spesso a Roma il Prof. Beltrami che qualche volta mi volle a pranzo. Non so ora se la tua signora siasi pienamente ristabilita dalla sua ribaltata che fece in carrozza alcuni mesi fa.

Dunque tanti e tanti saluti: ti prego de' miei omaggi alla tua signora

Tuo aff. Amico
Ing. Crotti

[82] Probabilmente Crotti si riferisce al sisma che ebbe luogo il 16 novembre 1894 e ha scritto come anno il 94 invece del 95.

**53**
Gian Antonio Maggi a Francesco Crotti
[s.l. e s.d.]

*[in matita:]* Non mandata

Per quanto alla famosa questione della Formola di Simpson, io mantengo invariabilmente il parere che ho scritto e stampato. Piuttosto non mi sento troppo animato a discorrerne ancora; poiché tu, bisogna pur dirlo! le ragioni opposte alla tua tesi non le vuoi sentire; e questo è un modo - in grande onore presso i filosofi propriamente detti - d'aver sempre ragione, o sempre torto, secondo il punto di vista. Una sola eccezione fa un brano della tua lettera del 30/93[83] dove riproduci le mie parole: "Tu discorri come se la curva del trapezio fosse incognita, e dalle ordinate date dovesse dedursi un valore plausibile dell'area", e rispondi: "Precisamente! È questo il caso. Quando si fa un rilievo sul terreno si misurano solo parecchie ordinate, e nulla si sa né si vuol sapere delle linee di perimetro." Ebbene! Nella tua risposta al prof. Jadanza (il cui scritto merita tutta l'attenzione) [84] stampata nel Politecnico ( ) tu poni la restrizione che "condotte le successive corde rettilinee, il perimetro corre assai vicino ad esse". È proprio il passo dal molto al qualcosa, dal zero al non zero, fatto il quale ogni altro è un'inezia. Concedici dunque che ci facciamo un'idea un po' meno vaga del grado in cui il perimetro s'avvicina alla successione delle corde, e riconosciamo se la successione degli archi parabolici inclusi nella formola di Simpson non gli si addatta *[sic!]* meglio o invece se ne discosta di più e tanto basterà perché ci serviamo debitamente dell'una o dell'altra formola, a seconda del caso, senza bisogno di fonderla in una che va a cercare la probabilità d'un fatto che ci sta sotto gli occhi. Sono le precauzioni a cui accenna il prof. J. in fine al §20. Né è vero che sì fatte restrizioni sulla cui strada, pare senza accorgertene, ti sei messo tu pure, si accampano ora per la prima volta, perché sono incluse in più formole d'approssimazione. "E non meravigliamoci" etc. Copio queste parole dal principio della prima lettera del...[85] Tu attribuisci alla tua formola l'inestimabile vantaggio di non richiedere tali restrizioni; ma la verità è che non se ne possono fare, perché al criterio di un'approssimazione determinata, sostituisce l'arbitrio della media aritmetica e superflua la briga di preoccuparsi quanto ci si accosta o si discosta dal vero, perché diventa impossibile appurarlo.

Col prof. Colombo ha scambiato sull'argomento qualche parola, lo scorso Autunno; non di più, perché m'è parso che non avesse gran voglia di discorrerne. Con tutto il rispetto, non approvo ch'egli abbia escluso dal suo Manuale la formola di Simpson. Non mi so poi persuadere che né egli né altri aderisca alla tua sentenza che quella formola è un errore di logica; perché perciò bisognerebbe che contraddicesse ai principii donde si deduce, quelli sui quali insisto dalla mia prima lettera, quelli affermati dal prof. J. e non a principii diversi, lasciamo stare se ammissibili o no per sé stessi.

---

[83] #47.

[84] Potrebbe trattarsi di: "Una difesa della formola di Simpson ed alcune formole di quadratura poco note: contro-osservazioni dell'Ing. Crotti", *Il Politecnico - Giornale dell'ingegnere architetto civile ed industriale,* 1893, v. 25, pp. 354-362.

[85] I puntini sono di Maggi. La lettera richiamata è la #44 del 4/12/1892 scritta da Maggi stesso.

## 54

Gian Antonio Maggi a Francesco Crotti

[frammenti; s.l. e s.d.]

(dimostrazione non comunicata)

.......

Se qualcuno coi dati del quesito applica la formola di Simpson etc. vuol dire che nella valutazione del tutto o della parte s'accontenta d'un'approssimazione che implica errori la cui differenza supera 120:

Indichino T e P l'area vera del tutto e della parte: t e p gli errori che si commettono valutandole colla formola di Simpson valendosi dei dati del quesito. Sarà

$$T=9600+t$$
$$P=9720+p.$$

Di cui

$$T-P=-120+t-p.$$

E poiché dev' essere $T-P>0$, bisogna che sia $t-p>120$.

*[altro frammento contenuto nel fascicolo dedicato alla corispondenza con Crotti]*

perpendicolare all'asse delle x, avendo per equazione $y=a+bx+cx^2$ si può effettivamente obbligare a passare per tre punti presi ad arbitrio. La somma delle aree dei trapezi parabolici analoghi a quello che si considera corrispondenti a questi archi porta a fornire un valore ancora più prossimo a quello del trapezio medesimo.

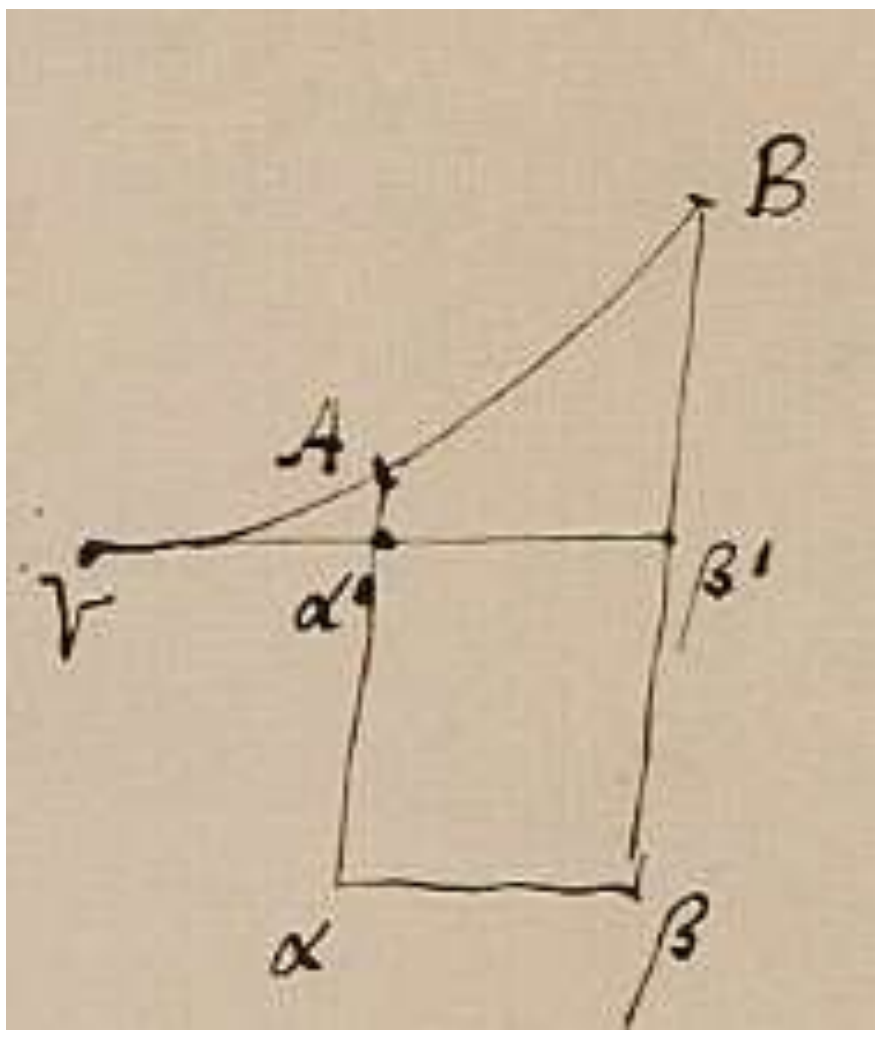

Calcoliamo l'area d'un trapezio parabolico, supposta la parabola coll'asse perpendicolare alla base, e convessa verso la base medesima. Supposto V il vertice della parabola a cui appartiene l'arco AB, la tangente ad essa in V sarà parallela alla base αβ del trapezio, la cui area sarà in primo luogo eguale alla somma di quelle del rettangolo rettilineo αβα'β' e del trapezio parabolico AB α'β'.

**55**
**Ermenegildo Daniele** a Gian Antonio Maggi
[busta indirizzata a: Ch.mo Prof. G.A. Maggi - Preside della Facoltà di Scienze
Corso Plebisciti 3 - Milano; busta e carta intestate R. Università di Pisa]

Pisa, 11.III.25.

Chiar.mo Professore,

La Sig.na Simonatti[86] mi presentò di questi giorni una tesi in Fisica mat. su un tema assegnato da Lei. Nella tesi non vi sono riferimenti né indicazioni bibliografiche, ma si comprende che la base della preparazione dev'essere costituita da un corso di Sue lezioni; inoltre, dal modo come la tesi è redatta, non risulta troppo bene quali siano stati i confini del tema da Lei proposto.

Perciò io Le sarei grato se Ella potesse darmi qualche indicazione in merito; come pure gradirei un Suo giudizio sommario sul valore della Signorina, cosa che nessuno meglio di Lei è in condizione di saper fare: e ciò non soltanto, per mia norma, ma anche per norma tutta la Commissione.

Mi creda, coi più cordiali saluti
Suo dev.mo

Daniele

**56**
Gian Antonio Maggi a Ermenegildo Daniele
[carta intestata: R. Università di Milano - Facoltà di Scienze - Il Preside]

Milano 15 Marzo 1925.

Gentilissimo Collega,

Per tesi di laurea, io proposi alla Sig.na Simonatti di fare, sulla base delle equazioni di Hertz, una ricerca simile a quello che, sulla base dell'equazione di d'Alembert, forma l'oggetto della mia Nota "Sulla propagazione delle onde di forma qualsivoglia nei mezzi isotropi" (Rendic. dell'Acc. dei Lincei, 19 Dicembre 1920): per estenderla possibilmente al problema della riflessione e della rifrazione di onde incidenti sferiche ad una superficie piana, sviluppando i concetti che si trovano accennati nel §6 della stessa Nota. E poiché la mia ricerca traeva origine, attraverso ad una Nota del Somigliana, dalla Nota del Laura "Sulla propagazione delle onde in un mezzo indefinito", le prestai questa Nota e alcune Note del Love sullo stesso argomento, per estendere se le riusciva, anche in quel senso le sue ricerche.

La Sig na Simonatti, per quanto mi consta dal manoscritto in preparazione, che mi mandò per visione, si limitò alla prima parte, nella cui esecuzione mi è parso che mostrasse una certa abilità nel maneggio del calcolo, e ottenesse anche qualche risultato non privo d'interesse; come, se ben ricordo, la riproduzione, colla condizione della trasversalità, delle equazioni da me ottenute, sulla base dell'equazione di d'Alembert, anziché delle equazioni di Hertz. A questo studio, nel suddetto manoscritto, si trovava poi premessa l'integrazione delle equazioni di Hertz, colle note condizioni che ne determinano univocamente, la soluzione, in un mezzo illimitato. La trattazione, non presa, ch'io sappia, da altra parte, non manca forse di una certa originalità. Altre parti, sempre col ricordo del manoscritto mostratomi, nelle scorse vacanze,[87] sono lavoro di compilazione, tratto particolarmente dalle mie lezioni di Ottica Fisica, che la Sig.$^{na}$ Simonatti seguì nell'anno 1921-22. Formalmente, d'altro, che conteneva il manoscritto, non resta sufficiente ricordo per accennarne un giudizio.

[86] Si vedano gli scritti da Maggi a Simonatti.

Ricordo che all'esame la Sig.$^{na}$ Simonatti ottenne la lode. Come, d'altronde, si è sempre distinta per diligenza e profitto, per modo che la credo meritevole di sincera raccomandazione.

Spero, con questo, di avere abbastanza soddisfatto il Suo desiderio, ed anche di arrivare in tempo utile, poiché la gradita Sua mi è giunta con qualche ritardo, in conseguenza del passaggio per l'Università. Il mio indirizzo, di cui confido che non manchi di valersi, è Corso Plebisciti 3.

Mi è grata intanto l'occasione per presentarLe i più cordiali saluti, colla preghiera di farne parte ai cari colleghi pisani, e per confermarmi

Suo aff.$^{mo}$ collega
G.A.M.

[87] Si veda la lettera #180 dell'agosto 1924.

**57**

**Carlo Del Lungo** a Gian Antonio Maggi

[busta intestata: Firenze - Gabinetto di Fisica e indirizzata a: prof. Gian Antonio Maggi - Carinana (Firenze); scritto dal Maggi in matita rossa: *Per "Giudizii"*]

Di campagna 7.Agosto 24.

Illustre professore; mi faccio lecito sottoporre alla sua profonda competenza in tale materia, questo che a me sembra un paradosso di meccanica razionale, del quale non riesco a trovare una spiegazione del tutto convincente.

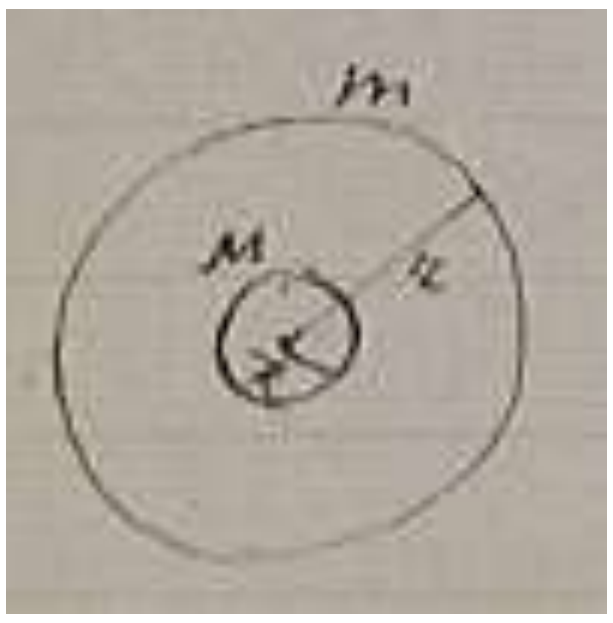

Uno strato sferico di massa m e di raggio r, contiene un altro simile strato di massa M e di raggio R. La sfera m non ha azione sui punti interni, quindi fra le due sfere non si esercita azione attrattiva e perciò il sistema non possiede energia potenziale relativa.
Ma questa conclusione che sembra legittima, mi apparisce non vera considerando le cose da un altro punto di vista.
Supponiamo le due masse riunite insieme nel raggio R: per distaccare dalla massa totale la massa m portandola al raggio r occorre certamente un lavoro e questo deve corrispondere all'energia potenziale reciproca dei due strati, ossia all'energia che vien restituita se lo strato esterno va in pezzi e precipita sulla massa centrale.

Per calcolare quest'energia potenziale, mi pare poter ragionare così.
Se la sfera unica (M+m) passa dal raggio R al raggio r il lavoro di dilatazione è dato da

$$\frac{(M+m)^2}{2}\left(\frac{1}{R}-\frac{1}{r}\right)=\frac{M^2+m^2+2Mm}{2}\left(\frac{1}{R}-\frac{1}{r}\right)$$

ed è composto di tre termini. Il 1° è il lavoro dovuto alla sola massa M; il 2° alla sola m; il 3° $\frac{Mm}{R}-\frac{Mm}{r}$ all'azione reciproca delle due masse nella dilatazione. Se r=∞ la separazione è completa, quindi $\frac{Mm}{R}$ misura il lavoro di distacco delle due masse.

Ora, se nella disposizione iniziale (quella della figura) le due masse non agiscono più l'una sull'altra, si può considerare avvenuto il distacco, quindi compiuto il lavoro $\frac{Mm}{R}$, che dovrebbe rappresentare l'energia potenziale del sistema, la quale sarebbe indipendente dal raggio r! Ma come resulta energia potenziale se fra le due masse non c'è attrazione?
Ecco, per me, il punto oscuro.

Può darsi che la questione sia già stata analizzata e discussa e che il dubbio mio dipenda soltanto dalla mia ignoranza.

In ogni modo, parendomi che il problema se pure facile ed elementare, presenti un certo interesse didattico, Le sarò veramente obbligato di quello che vorrà comunicarmi in proposito mettendolo nei termini veri quali forse io non ho saputo trovare.
Con profonda stima e memori rispettosi saluti mi abbia Suo dev[mo]

Carlo Del Lungo
R. Liceo Michelangiolo.
Firenze

P.S. Ragionando in altro modo l'espressione dell'energia potenziale resulta diversa. Partendo dalla disposizione della figura, e dilatando le due sfere fino ad un raggio infinito il lavoro resulterà

$$\frac{m^2}{2r}+\frac{M^2}{2R} \quad \{\text{escludendo ogni azione dell'una sull'altra}\}$$

Se le due masse sono riunite nel raggio R il lavoro nella dilatazione fino al raggio infinito sarà

$$\frac{(M+m)^2}{2R}$$

La differenza

$$\frac{(M+m)^2}{2R}-\frac{m^2}{2r}-\frac{M^2}{2R}$$

sarà il lavoro compiuto per distaccare la m dalla M portandola alla distanza r come nella figura, ossia sarà

$$\frac{Mm}{R}+\frac{m^2}{2}\left(\frac{1}{R}-\frac{1}{r}\right)$$

cioè dipendente dal raggio r.

## 58

Gian Antonio Maggi a Carlo Del Lungo

Gavinana (Firenze) 13 Agosto 1924.

Egregio Professor Del Lungo,

Alla gradita Sua, trasmessami da Pisa a Gavinana, dove mi sto io pure in campagna, trovo da rispondere, prima di tutto, ch'io non vedo perché l'energia potenziale relativa dei due strati sferici da Lei considerati dovesse essere nulla, per la ragione che è nulla la forza motrice di ciascuno sull'altro. Dovrebbe verificarsi per questo, fra l'energia potenziale relativa di due corpi, e la loro mutua forza motrice una relazione ch'io non so che esista.

Mentre l'energia potenziale relativa di due corpi si riconduce al lavoro delle forze esercitate dall'uno dei due corpi sui punti dell'altro, in conseguenza del passaggio di questi punti, dalla posizione supposta all'infinito: lavoro che non risulta nullo, e non ha ragione d'essere, perché, nella posizione iniziale dei punti mobili la forza motrice mutua dei due strati è nulla.

Nel caso dei due strati, è forse inutile osservare che, nella stessa posizione in cui è nulla la forza motrice di ciascun strato sull'altro, cioè, il risultante delle forze motrici di ciascuno strato sopra i singoli punti dell'altro, queste forze motrici sono nulle pei punti dello strato interno, ma non già pei punti dello strato esterno.

Giustamente, Ella invoca l'accennato lavoro, ma non mi rendo ragione dei calcoli.

Il potenziale mutuo dei due corpi, che rappresenta detto lavoro, col segno cambiato, è espresso, in generale, da

$$\int_{\tau_1} U_2 k_1 d\tau_1, \qquad U_2=\int_{\tau_2}\frac{k_2 d\tau_2}{\tau_2}, \qquad \Delta=\sqrt{(x_1-x_2)^2+(y_1-y_2)^2+(z_1-z_2)^2},$$

*[scritta aggiunta a lato:* V invece di U. *]*

dove $\tau_1$ e $\tau_2$ indicano i campi appartenenti a due corpi, $k_1 d\tau_1$, e $k_2 d\tau_2$ gli elementi di massa appartenenti al punto generico di $\tau_1$ e $\tau_2$, $\Delta$ la mutua distanza del punto generico $(x_1,y_1,z_1)$ dell'uno e del punto generico $(x_2,y_2,z_2)$ dell'altro.
Nel caso dei due strati sferici, attribuendo l'indice 2 allo strato interno, si ha

$$U_2=\frac{m_2}{\rho_1},$$

dove $\rho_1$ indica la distanza dal centro comune dei due strati dal punto generico dallo strato esterno. Quindi, pel potenziale mutuo,

$$m_2 \int_{\tau_1} \frac{k_1 d\tau_1}{\rho_1}=4\pi m_2 \int_R^r k_1 \rho_1 d\rho_1 .$$

Attribuendo invece l'indice 2 allo stato esterno si ha

$$U_2=4\pi \int_R^r k_2 \rho_2 d\rho_2 .$$

Quindi, pel potenziale mutuo,

$$4\pi \int_R^r k_2 \rho_2 d\rho_2 \int_{\tau_1} k_1 d\tau_1=4\pi m_1 \int_R^r k_2 \rho_2 d\rho_2 ,$$

che, scambiando 1 con 2, è il risultato di prima.
Mi farà sempre piacere a scrivermi se altro desidera. Intanto mi è grato presentarLe i mie migliori saluti, pregarLa di ricordarmi all'illustre Suo padre,[88] e confermarmi

Aff$^{mo}$ Suo
G.A.M.

**59**
**Alfonso Di Legge** a Gian Antonio Maggi
[busta indirizzata a: Illustre Preofessore G. A. Maggi
Milano - Corso Plebiscito 3]

Roma 27 maggio 1936
Piazza Borghese 91

Illustre professore,

Mi permetto di raccomandarle il mio amico prof Ignazio Canestrelli il quale avrebbe bisogno di qualche chiarimento dalla sua cortesia sui principi della meccanica.[89]
Saluti cordiali dal

Suo aff$^{mo}$ collega
Alfonso Di Legge

[88] Isidoro, letterato e politico.
[89] Si veda la corrispondenza con Canestrelli.

## 60

Gian Antonio Maggi al **Direttore del *Journal des Débats***[90]

Valnegra (Bergamo) 22 Sett. 1930

Monsieur le Directeur du Journal des Débats,

À propos de votre intéressant article, "Le latin en Italie" paru dans le Numéro du 20 Septembre,[91] permettez-moi de faire présent à votre correspondant que les vénitiens parlent communément le vénitien, qui était la langue de la Serenissima, et qu'il s'ensuit, à Venise et ailleurs dans la Vénète, une prononciation du latin à la vénitienne.

Quant aux tché et aux tchi, l'interlocuteur de votre correspondant a tout simplement perlé vénitien en disant vitsino pour vitchino; à Milan, en milanais, on aurait dit vesin, proche de votre voisin. À Milan (je suis milanais) on prononce quoui, coum spiritoum touo, Dominuos, sans appuyer, il est bien vrai, excessivement sur ou; mais, quand j'étais enfant on entendait bien souvent dominus et Spiritu tuo avec l'u prononcé à la lombarde, c'est-à-dire à la française. En Toscane votre correspondant aurait entendu Dominous-se vobiscoum-me, le toscan n'admettant pas la terminaison par consonne. Enfin, ce n'est pas par la prononciation du latin de la part de tous les fidèles et les officiants de la péninsule, à laquelle Rome n'a pas tardé à étendre la civitas romana, qu'on peut dire Roma locuta est. Agréez, Monsieur le Directeur, l'assurance de ma considération la plus distinguée.

Votre abonné
Gian Antonio Maggi

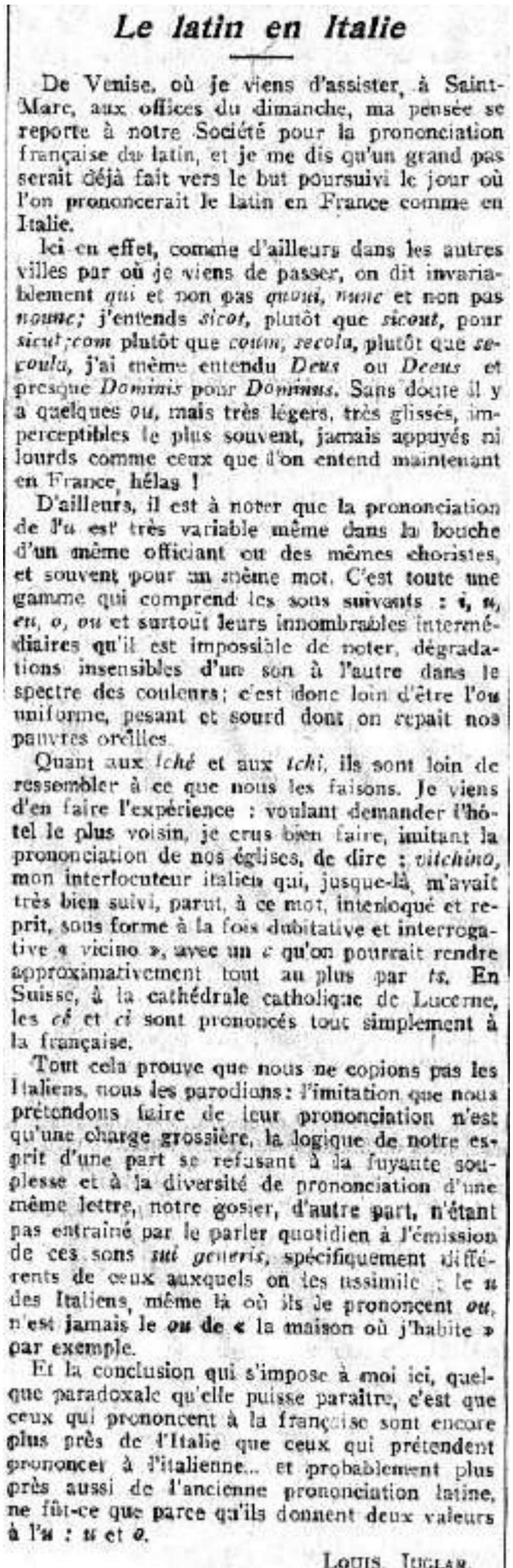

**Le latin en Italie**

De Venise, où je viens d'assister, à Saint-Marc, aux offices du dimanche, ma pensée se reporte à notre Société pour la prononciation française du latin, et je me dis qu'un grand pas serait déjà fait vers le but poursuivi le jour où l'on prononcerait le latin en France comme en Italie.

Ici en effet, comme d'ailleurs dans les autres villes par où je viens de passer, on dit invariablement *qui* et non pas *quoui*, *nunc* et non pas *noumc;* j'entends *sicot*, plutôt que *sicout*, pour *sicut*; *com* plutôt que *coum*, *secola*, plutôt que *secoula*, j'ai même entendu *Deus* ou *Deeus* et presque *Dominis* pour *Dominus*. Sans doute il y a quelques *ou*, mais très légers, très glissés, imperceptibles le plus souvent, jamais appuyés ni lourds comme ceux que l'on entend maintenant en France, hélas !

D'ailleurs, il est à noter que la prononciation de l'*u* est très variable même dans la bouche d'un même officiant ou des mêmes choristes, et souvent pour un même mot. C'est toute une gamme qui comprend les sons suivants : *i*, *u*, *eu*, *o*, *ou* et surtout leurs innombrables intermédiaires qu'il est impossible de noter, dégradations insensibles d'un son à l'autre dans le spectre des couleurs; c'est donc loin d'être l'*ou* uniforme, pesant et sourd dont on repait nos pauvres oreilles.

Quant aux *tché* et aux *tchi*, ils sont loin de ressembler à ce que nous les faisons. Je viens d'en faire l'expérience : voulant demander l'hôtel le plus voisin, je crus bien faire, imitant la prononciation de nos églises, de dire : *vitchino*, mon interlocuteur italien qui, jusque-là, m'avait très bien suivi, parut, à ce mot, interloqué et reprit, sous forme à la fois dubitative et interrogative « vicino », avec un *c* qu'on pourrait rendre approximativement tout au plus par *ts*. En Suisse, à la cathédrale catholique de Lucerne, les *cé* et *ci* sont prononcés tout simplement à la française.

Tout cela prouve que nous ne copions pas les Italiens, nous les parodions : l'imitation que nous prétendons faire de leur prononciation n'est qu'une charge grossière, la logique de notre esprit d'une part se refusant à la fuyante souplesse et à la diversité de prononciation d'une même lettre, notre gosier, d'autre part, n'étant pas entraîné par le parler quotidien à l'émission de ces sons *sui generis*, spécifiquement différents de ceux auxquels on les assimile ; le *u* des Italiens, même là où ils le prononcent *ou*, n'est jamais le *ou* de « la maison où j'habite » par exemple.

Et la conclusion qui s'impose à moi ici, quelque paradoxale qu'elle puisse paraître, c'est que ceux qui prononcent à la française sont encore plus près de l'Italie que ceux qui prétendent prononcer à l'italienne... et probablement plus près aussi de l'ancienne prononciation latine, ne fût-ce que parce qu'ils donnent deux valeurs à l'*u* : *u* et *o*.

LOUIS JUGLAR.

---

[90] Sono due le copie di minuta di questa lettera che si trovavano nelle *Rusticationes*. Sulla seconda Maggi riporta: *Pubblicata dai Débats nel n. del 27 Settembre 1930 col titolo* Le latin en Italie; nella prossima pagina viene riprodotta.

[91] Si trova in rete: http://gallica.bnf.fr/ark:/12148/bpt6k506359q/f2.item; a lato una riproduzione.

*Le latin en Italie.* — Nous recevons la lettre suivante :

« A propos de votre intéressant article « Le latin en Italie », paru dans le numéro du 20 septembre, permettez-moi de signaler à votre correspondant que les Vénitiens parlent communément le vénitien, qui était la langue de la Sérénissime, et, partant, qu'il existe, à Venise et dans la Vénétie, une prononciation du latin à la vénitienne. Quant aux *tché* et aux *tchi*, l'interlocuteur de votre correspondant a tout simplement parlé vénitien, disant *vitcino* pour *vitchino*; à Milan, en milanais, on aurait dit *vesin*, proche de votre *voisin*. A Milan (je suis Milanais), on prononce *quous, Dominous, coum spiritou touo*, sans appuyer, il est vrai, excessivement sur *ou* ; mais, quand j'étais enfant, on entendait bien souvent *dominus* et *spiritu tuo* avec l'*u* prononcé à la lombarde, c'est-à-dire à la française. En Toscane, votre correspondant aurait entendu *Dominous-se vobiscoum-me*, le toscan n'admettant pas la terminaison par consonne. Enfin, ce n'est pas par la prononciation du latin de la part de tous les fidèles et des officiants de la péninsule, à laquelle Rome n'a pas tardé à étendre la *civitas romana*, qu'on peut dire *Roma locuta est*.

» GIAN ANTONIO MAGGI. »

**61**

Gian Antonio Maggi al **Direttore de *la Meneghina***[92]

Spettabil Direzion de "la Meneghina"

Minga per fa lezion, tanto pu che sont in pension, ma per amor de la verginitaa de la nostra lingua meneghina, me permetti de fa l'osservazion che Scior Zucchetti, e sciori tai e tai (in la Circolar per l'Assemblea del 21) lighen nò cont "el me sur Lella", "sur professor Ronchett", "sur sergent maggior", ecc. ecc. E l'istes sura, e minga sciora, dennanz a nomm e cognomm d'ona donna. Mi poeu disi semper inscì, e me par che inscì disen semper i milanes compagn de mì. Scior e sciora, per sur e sura, se senten sì, ma in Brianza. E se po' legg scior Bernard, in la novella "Nostalgia di Milano" del "Corriere" d'incoeu. Ma andaremm minga a imparà el milanes del sur Panzini.[93]

Se poeu ciappò on gamber, sont chi a toeu su la lezion.

Con tanti salud,

El socio

Milan 17 Marz 1932.

Giovann Antoni Magg

[92] Nel 1931 viene nominato Presidente della Famiglia Meneghina (ai tempi si chiamava *Resgiò*, in milanese) il Cav. Pier Gaetano Venino, Senatore del Regno, che rimase in carica per molti anni, fino al 1945.

[93] Forse si riferisce a Alfredo Panzini, lessicografo.

**62**
**Giuseppe Erede** a Gian Antonio Maggi
[carta intestata: Ing. Giuseppe Erede - Via Assarotti 29$^{A}$ - Genova (2)]

Il 19 Gennaio 1924

Chiarissimo Professore

Mi presento da me, perché le buone conoscenze che avevo in Pisa sono tramontate da un pezzo, per esempio il prof Dini e Giovanni Cuppari (io sono molto vecchio, e da qualche anno: il prof. Italo Giglioli *[sic!]*. Sono stato professore nell'Istituto tecnico di Firenze, ai tempi di Emilio Villari e del Roiti! Poi sono stato nel nuovo Catasto, per uscirne quando vennero abolite le Direzioni Compartimentali, e quindi da molti anni sono a riposo, e faccio il dilettante de omnibus rebus. Ho fatto parecchie pubblicazioni, specialmente sulla Topografia e sul Catasto. La mia Topografia (Bemporad) ha 8 edizioni, e la Geometria pratica (Hoepli) ne ha sei.

Verso la fine del 1921 mi colpirono le stranezze che si andavano stampando sulla relatività, in libri, riviste e giornali; vidi che sul famoso esperimento si sbagliava stranamente, e pensai di potermene occupare senz'altro, con le poche cognizioni che avevo di fisica. Mandai una Nota ai Lincei per mezzo del Roiti, ma gli giunse durante l'ultima sua malattia; la feci dare al prof. Volterra, che la passò al prof. Levi-Civita, e questi me la rimandò osservando giustamente che io non avevo considerato l'effetto del moto degli specchi sulla riflessione, e avvertendomi che il Righi aveva studiato l'esperimento "in modo piuttosto difficile a seguire". Le parole citate in principio della mia pag. 3 sono sue.
Feci dei nuovi calcoli, anch'essi sbagliati, su molti progetti d'esperimento, e diversi esperimenti eseguii, con gran perdita di tempo e fatica. Ai professori di fisica in questa città ne parlai senza alcun frutto. Non si degnarono d'occuparsene, o mostrarono la loro ignoranza in materia. Finalmente, a furia di sbagliare, chiedendo qualche parere fuori di qui, sono giunto a conoscere tutto ciò che bisogna conoscere per giudicare dell'origine della relatività. Tanto che senza superbia posso dire di conoscerla più di ogni altro uomo al mondo.

Ora dirà, o ha già detto "quare conturbas me?". Ecco, io vorrei che Ella acconsentisse a presentare l'acclusa Nota ai Lincei per i Rendiconti.[94] Ella vi ha presentato una Nota sulla famosa trasformazione, segno che s'interessa alla relatività, e giudica "attraente" il fatto della dimostrazione del coefficiente di trascinamento del Fresnel, segno che non ne è perfettamente convinto. Si esprimerebbe in altro modo.

Ella non è dunque del gran numero di coloro che si sono fatti relativisti, come ora tanti si fanno fascisti, e non ricuserà d'aiutare la pubblicazione di una critica incontestabile. Forse Ella troverà soverchiamente polemiche le mie osservazioni a pag. 4bis. Ma, dopo letto un articolo, assai difettoso, del prof. Bouasse in Scientia,[95] io gli scrissi per uno schiarimento, ed egli alla risposta aggiunse che la questione in Francia era diventata una questione di religione. Gli israeliti si adoperarono per impedire la divulgazione delle sue critiche. "Toute la juïverie a donné". L'anno scorso, nell'estate, mandai un articolo alla Revue scientifique, e poi riscrissi per farvi un'aggiunta. Nessuna risposta. Oltre all'esprit de corps degli israeliti (temo molto dal Volterra) vi sono tutti coloro che si sono compromessi esaltando la relatività e il suo profeta, e la caduta della meccanica di Galileo e di Newton, e dicendo che capivano le teorie dell'Einstein.

[94] Non presente nel Fondo.
[95] H. Bouasse, "La question préalable contre la théorie d'Einstein", *Scientia*, v. XXXIII, 1923, pp. 13-24.

Veramente sarebbe indicato al mio scopo il prof. Somigliana, poiché nel parlare della trasformazione di Lorentz si disse lieto di portare un contributo alla guerra che opportunamente s'inizia contro la relatività. Ma l'anno scorso egli pubblicò in Scientia[96] uno scritto tanto sbagliato, che termina appunto mostrando la sua ignoranza dell'origine, che sono certo egli farebbe d'ogni erba fascio per impedire che il mio scritto venga pubblicato, e gli faccia fare una pessima figura.

Scusi tanto disturbo, e gradisca con i miei ringraziamenti anticipati il mio distinto ossequio

Dev.mo
Prof. Ing. G Erede

**63**

Giuseppe Erede a Gian Antonio Maggi
[cartolina postale indirizzata: Al Chiarmo Signor Commre Prof. Gio. Ant. Maggi
R. Università di Pisa]

Via Assarotti, 29 - Genova, 2
Il 21 Gennaio 1924

Chiarmo Professore

Confermando la mia lettera del 19 corrente, devo dirle che ho riguardato l'ultima delle varie minute della Nota che l'accompagnava, e ho veduto che distrattamente vi scrissi "meno di un centomillesimo" nella citazione dall'Enciclopedia britannica riguardante il bolometro, e certamente ricopiando in fretta ho messo le stesse parole nella Nota. Invece l'Enciclopedia citata dice: "un cento milionesimo", che va d'accordo colle parole precedenti "Nella sua forma più raffinata", e permette di contare sul risultato dell'esperimento. Questo non è facile a eseguire e penso di scriverne al Salmoiraghi,[97] col quale sono in buona relazione, dopo averne detto bene nella mia Topografia.

Con distinti ossequi
Devmo G Erede

---

[96] C. Somigliana, "I fondamenti della Relatività", *Scientia*, v. XXXIV, 1924, pp. 1-10.
[97] Potrebbe trattarsi di Angelo S., imprenditore, ottico e ingegnere italiano, produttore di strumenti di alta precisione per l'industria e la geodesia.

**64**[98]

Gian Antonio Maggi a Giuseppe Erede

[busta indirizzata a Ch. mo Sig. Prof. Ing. Giuseppe Erede - Via Assarotti 29[A] Genova (2); mittente: Raccomandata spedita Prof. G.A. Maggi - Via Risorgimento, 15 Pisa]

Pisa 26 Gennajo 1924.

Pregiatissimo Signor Ingegnere,

Col costante desiderio di prestarmi, dove posso, a chi crede di rivolgersi a me, sono sinceramente spiacentissimo di non poter appagare la Sua domanda, quando, a parte i pregiati volumi, Ella mi richiama i nomi del Dini, del Cuppari, del Giglioli, e manifesta tanto interesse per la ricerca della verità.

Ella mi concederà che, presentando un lavoro d'altri ad un'Accademia, per l'inserzione negli Atti, se non si subentra all'autore, si dichiara tuttavia solidarietà colle sue opinioni. Ella me lo concede tanto che dice di rivolgersi a me, perché da una mia espressione ha desunto ch'io non sono convinto della Relatività. Quell'espressione mi ha alquanto tradito. Certo io non professo per la Relatività l'adesione entusiasta di parecchi. E come sono io è il Volterra, a cui Ella accenna. Ma qualunque sia il mio grado di convinzione, io, per esempio, espongo la Relatività nel mio Corso di Fisica Matematica, e ho in corso di stampa, pel Bollettino del Prof. Gino Loria di codesta Università, una recensione favorevole sul recente volume relativistico del Vasiljef.[99] Non mi sento quindi adatto a presentare ai Lincei una Nota, dove le conclusioni del Michelson e de' suoi collaboratori si affermano involgere "due gravi errori, un errore di logica e un errore di Fisica", e sono attribuite ad "una specie di ebbrezza, data dall'idea dubbiosa di trovarsi davanti ad un grande avvenimento scientifico, ebbrezza che scacciò il dubbio...", ed è affermato prima come incontestabile il concetto dell'etere, sul quale io convengo colla Relatività nel fare molte riserve.

Potrei essermi convertito? Per questo, neppure mi posso arrendere alle Sue obbiezioni. Quando Ella domanda come non ci poteva essere una linea nera, dove c'era sconcordanza di n+½ lunghezze d'onda, trovo da rispondere che la Relatività, accettando dall'esperimento di Michelson la mancanza della linea nera, al luogo indicato dalla teoria classica, conclude colla esclusione della suddetta sconcordanza, e colla nota modificazione dei concetti classici, che deriva da quella esclusione. E quando Ella contesta la conclusione che non si trovò spostamento della linea nera, pur essendosene riscontrati piccoli spostamenti, io non ho modo d'impugnare l'affermazione dei fisici americani che le loro misure sono garantite tra tali limiti di approssimazione, da dover escludere che, per causa di errori d'osservazione, i valori osservati indicassero come risultato plausibile, lo spostamento nullo, piuttosto che lo spostamento, che, in base alla teoria classica, avrebbe dovuto presentare la linea nera.

Si può piuttosto impugnare il metodo, come fece il Righi, senza però venire, a quanto sembra, a conclusioni interamente chiare. Ad ogni modo, sono tutte discussioni di ordine sperimentale, da rimettere a chi, nella Fisica Sperimentale, possieda la competenza ch'io sono ben lontano d'attribuirmi. Per cui, quand'anche non La persuadessero le mie precedenti osservazioni, dovrei sempre chiederLe di dispensarmi di intervenire in una discussione che oltrepassa la mia competenza.

---

[98] In due copie di minuta: *Risposta all'Ing. Giuseppe Erede all'aclusa sua lettera in data 19 Gennajo 1924*. Su quella qui trascritta Maggi scrisse in lapis rosso: *Copia della lettera spedita*.

[99] Si tratta della recensione apparsa su *Il bollettino di Matematica* di Gino Loria (R. Università di Genova, 24 Maggio 1924, Fasc. I, anno III) del libro di A.V. Vasiljef, *Spazio, Tempo, Movimento*.

Simile osservazione troverei da fare alla Sua obbiezione al ragionamento del Lorentz, che, Ella dice, identifica due punti, che, nell'esperienza, sono distinti. Per quanto infine al Le Roux, mi sembra ch'egli, in base all'equazione notissima, e scarsamente usata, col nome di equazione di d'Alembert, non fa in sostanza che dimostrare che quello ch'egli chiama il pseudotempo di Lorentz si estende alle onde di qualsiasi specie. E questo è conforme alla Relatività, la quale reclama di estendere il suo dominio alla stessa Meccanica Razionale. Soltanto, quando non entrano in campo velocità dell'ordine di quella della luce, torna sensibilmente lo stesso valersi del pseudotempo o del tempo ordinario.

Le rimando quindi con questa, non senza chiederLe scusa di non prestarmi alla Sua domanda, e rinnovargliene l'espressione del mio dispiacere, i quattro fogli, ch'Ella mi mandò colla pregiata Sua del 19 corrente. Non so a che cosa si riferiscono le indicazioni pag. 3[bis] nel foglio 3, e pag. 4[bis] nel foglio 4, in lapis turchino. Io non ricevetti altro, all'infuori della Sua cartolina del 21, con cui corregge una citazione, di cui però io non trovo traccia nei quattro fogli suddetti, né nella lettera che li accompagnava. Dalla Enciclopedia Britannica Ella cita in Nota le parole "Non si ebbe la minima traccia di spostamento".

Aggradisca insieme i miei più distinti saluti, coi quali La prego di credermi sempre

Dev.[mo] Suo
Gian Antonio Maggi.

**65**
**Guido Facciotti** a Gian Antonio Maggi
[busta indirizzata: All'Ill.[mo] Signor Prof. Gian Antonio Maggi
Lanzo d'Intelvi - Como]

Ill.[mo] Signor Professore,

Sottopongo al Suo giudizio questo mio modesto lavoretto, grato per tutte le osservazioni ch'Ella crederà di addebitarmi.
Sono a Milano sino al 14 p.v. e forse l'ultima settimana di settembre.

Ringraziando vivamente pel disturbo, porgo auspici di ogni bene e prosperità, nonché ossequiosi saluti.

G Facciotti.

Milano, li 3 settembre 1936 b XIV° Via Soperga, 50

Formola del Facciotti[100]
(da sua lettera del 3 Sett. 1936)

Forza magnetica (unitaria) in O dal punto P in ds

$$dH = \frac{i \sin \varphi \, ds}{\rho^2}$$

[100] Lo scritto di Facciotti non è presente nel Fondo; ne restano solo i seguenti appunti scritti da Maggi che si trovavano nelle *Rusticationes* insieme al biglietto precedente, alla lettera seguente e al telegramma #68.

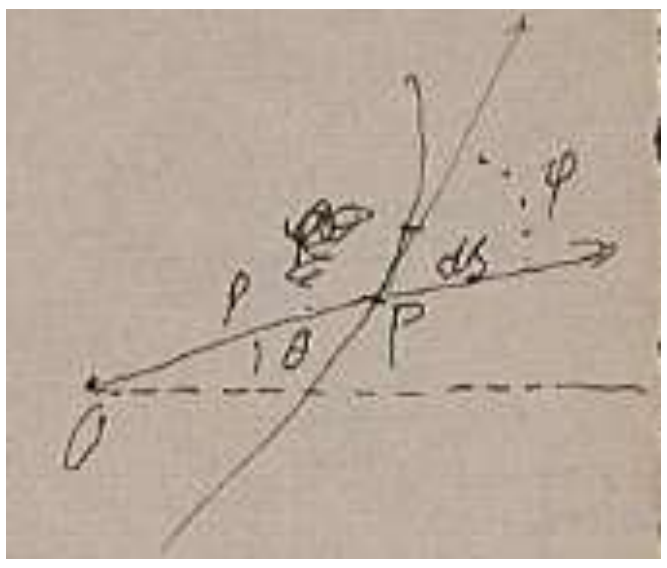

Assumiamo un sistema di coordinate polari col polo in O e l'asse polare comunque, con che competano al punto P le coord. ρ, θ.

Abbiamo

$$\tan\varphi = \frac{1}{\rho'}, \quad \sin\varphi = \frac{\rho}{\sqrt{\rho^2+\rho'^2}}, \quad ds=\sqrt{\rho^2+\rho'^2}d\theta$$

Quindi

$$dH=\frac{id\theta}{\rho}.$$

E, col reoforo, dato

$$\rho=\rho(\theta),$$

integrando,

$$H=i\int_{\theta_1}^{\theta_2}\frac{d\theta}{\rho(\theta)}.$$

## 66

Gian Antonio Maggi a Guido Facciotti

Lanzo d'Intelvi (Como) 8 Settembre 1936

Caro Facciotti,

Il Suo lavoretto si riduce in realtà ad un ben semplice esempio; ma dà un risultato, la (3), di notevole semplicità. Soltanto, appunto per la sua agevole deduzione, nasce il dubbio che sia già reso noto. Ella fatto ricerca a questo proposito?

Nell'Applicazione Ella dice "spire circondanti il polo soppresso per impedire che H diventi infinito": discorso, ad ogni modo, di colore oscuro. D'altra parte, com'è ben chiaro, tutte le spire, formate colla spirale, circondano il polo, tranne la prima $(0\leqq\vartheta\leqq 2\pi)$, che comincia col polo, dove è tangente all'asse polare; e questa sola (accogliendo la decomposizione in spire) è necessario escludere per evitare ρ=0, col qual valore (1) non ha più significato (meglio dire così di H che diventa infinito). Nulla impedisce di escludere $n_1>1$. Ma la dicitura, per le suddette ragioni, vuol essere mutata.

Per quanto alla costante α interpretata come diametro del filo, nell'ipotesi delle spire a contatto, bisogna regolarsi coll'osservazione che, col tendere del filo ad una linea, $\frac{1}{\alpha}$ tende a ∞.

Se non che il risultato dell'Applicazione non concerne di più che il valore di H nel polo della spirale da cui s'immagina ricavato il reoforo: risultato quindi troppo particolare per attribuirgli sufficiente importanza. Non si esce con questo dai confini di un esercizio.

Così, ho veduto il suo lavoretto con piacere, e lodo la coltura e abilità di cui vi dà saggio, ma se Ella ha inteso di sentire il mio giudizio anche sull'opportunità di pubblicarlo, mi sento tenuto a dire che, anche verificata la novità di (3), levati alcuni particolari da reputarsi superflui, e recato l'accennato ritocco nell'Applicazione, sempre non me ne sembrerebbe il caso.

Le rimando il M.S., perché Ella possa sentire il parere anche d'altri, affatto liberamente, sicuro ch'io non troverò nulla a ridire se, in materia di questo genere, non dovesse coincidere col mio.

Mi è grato intanto ricambiarLe i migliori saluti, e confermarmi

Aff$^{mo}$ Suo

P.S. La (3) richiama sicuramente una certa attenzione. E quando Ella potesse trovarne ad esporre un'applicazione di maggior importanza, il lavoro potrebbe acquistare consistenza così da valere anche giudizio favorevole per essere pubblicato. Perché non si metterebbe all'opera?

**67**

Gian Antonio Maggi a Guido Facciotti

Lanzo d'Intelvi (Como) 10 Sett. 1936

Caro Facciotti,

Mi sono accorto di un grave equivoco contenuto nel suo lavoro: quello di presentare φ - angolo del raggio vettore OP colla tangente nello stesso P - per anomalia di P, cioè angolo del raggio vettore OP con un asse fisso uscente da o. E tanto basta.

Voglia quindi considerare la mia lettera per annullata.

Con sinceri saluti,

Aff$^{mo}$ Suo

G.A. Maggi.

**68**

Gian Antonio Maggi a Guido Facciotti

[telegramma inviato al Prof. Facciotti - Via Soperga 50 - Milano]]

Prego annullare mia lettera oggi
Maggi.

**69**

Gian Antonio Maggi a Guido Facciotti

Lanzo d'Intelvi (Como) 10 Sett. 1936

Caro Facciotti,

Mi sono affrettato, col telegramma, a provvedere che, quando domani riceverà la mia prima lettera d'oggi, sappia che ho riconosciuto inesistente l'accennato equivoco. Fu invece un mio equivoco per non aver sotto gli occhi il M.S., che ho subito riconosciuto, ritrovando le formole.

Sono tornato sul Suo lavoro a proposito dell'espressione "per evitare ρ=0, pel qual valore la (1) non ha significato". Essa dà luogo a domandarsi che cosa avviene dell'integrale, quando vi si introduca ρ=0, come limite di ρ>0. Per cui meglio di tutto dire "per evitare ρ=0, introducendo il qual valore nell'integrale esso diventa infinito".

Con rinnovati saluti

Aff$^{mo}$ Suo
G.A. Maggi.

**70**

Guido Facciotti a Gian Antonio Maggi

Milano, li 24 settembre 1936 - XIV

Chiarissimo Signor Professore,

Ho proceduto oltre nel mio lavoretto e, incoraggiato dalle benevole parole dell'ultimo Suo gradito scritto, Le invio i risultati raggiunti.

Dopo essere stato a Trento per la sessione autunnale da ieri sera mi trovo a Milano e mi fermerò sino a tutto il 29, per raggiungere nuovamente la sede col giorno successivo .

Cordialissimamente ringraziando, porgo cari e distinti ossequî.

Dev.$^{mo}$ allievo

Guido Facciotti

Milano - Via Soperga, 50 - Trento - R Liceo Scientifico "G. Galilei"

**71**

Gian Antonio Maggi a Guido Facciotti

Milano 7 ottobre 1936

Caro Facciotti,

Per varie occupazioni, inerenti al ritorno dalla villeggiatura, non ho potuto scriverLe prima, e risponderLe a proposito del Suo M.S.

Ora, per questo, trovo da osservare, in primo luogo, che è superfluo il lungo calcolo a §2 per dedurre la (7), che coincide colla (1). Una volta che Ella trascura dh in confronto di a, non è necessario dedurre la dϑ' per concludere dϑ'=dϑ. Tanto vale, fin da principio, reputare confondibili, per gli effetti sul punto O, i due elementi a distanza dl.

Simile semplificazione comporta la deduzione del §3. Ma qua compare dI, che per

$$dH = dI\frac{ads}{\rho^3},$$

dI rappresenta intensità di corrente e d'altra parte, Ella tratta come un differenziale subordinato a dl per formare $j=\frac{dI}{dl}$ "densità di corrente". Io mi domando con questo quale significato possa avere questo dI, e debbo dirLe che non ne trovo la risposta. E poiché il dI compare nelle successive deduzioni e applicazioni, la mia riserva si estende al resto del suo M.S, che Le rimando, senza, pel momento, poterLe dire di più.

Alla ricerca di un dI, lo trovo in Perucca, Vol. II, §12, solenoide cilindrico indefinito.[101] Ma qua dI è subordinato a dx, elemento dell'asse del solenoide, mediante dI=indx, indicando con n il numero delle spire per unità di lunghezza.

Sempre disposto a darLe, dove posso, un buon consiglio, Le invio intanto cordiali saluti, coi quali mi confermo.

Aff.^mo Suo G.A. Maggi

P.S. A parte la suddetta riserva, l'argomento, come diretta e semplice applicazione della formula (1), non oltrepassa, secondo il mio modo di vedere, i confini dell'esercizio. Né credo che, con questo genere di problemi, sia facile uscirne.

**72**

Gian Antonio Maggi a Guido Facciotti

[biglietto]

Milano 16 Ottobre 1936

Caro Facciotti,

Io ho inteso puramente di dire che il Suo discorso, non si può connettere, per le accennate differenze, con quello del Perucca,[102] al quale Ella mi rimandava, e questo poi svolge una teoria, mentre il Suo discorso, a mio giudizio, introduce formole senza preparazione sufficiente, e non ne mostra sufficiente profitto.

Ad insistenza io non faccio, appunto, perché, non per pura forma, dico che mi fanno sempre piacere a rivolgersi a me. Ma è pur vero che l'insistenza, per sgombrare incertezze, finisce per obbligare ad esprimere la propria opinione con una chiarezza, che non fa sempre piacere ad usare.

Di nuovo, cordialmente La saluto e mi confermo

Aff^mo Suo

G.A. Maggi.

Egr Prof. Guido Facciotti

del R Liceo Scientifico G. Galilei di Trento

[101] Si veda la nota 47.

[102] Si veda la nota 47.

## 73

**[Paolo] Fassini Camossi** a Gian Antonio Maggi

[busta indirizzata a: Egregio Sig. Prof. Maggi - Via Risorgimento - Pisa]

Torino 4 aprile

Egregio Professore, non vorrei essere un importuno se mi rivolgo a Lei per pregarla di darmi uno schiarimento che, credo, non facile.

È possibile una soluzione dell'equazione: $x^4–ax–b=0$

E se possibile esiste un metodo piuttosto semplice, anche per approssimazione?

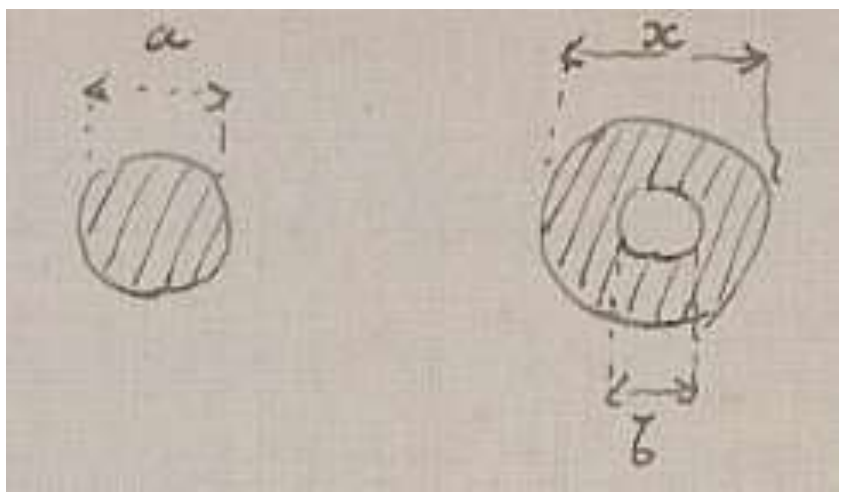

Tale risoluzione servirebbe per trovare facilmente il diametro esterno x di un albero cavo dati: il diametro interno b e il diametro a dell'albero pieno di egual resistenza.

La formola è ricavata dall'uguagliarne i momenti (*[...]* del Colombo).

Scusi, Professore la seccatura.

Colgo l'occasione per pregarla di presentare i miei omaggi alla sua Signora e di ricordarmi ai suoi figli.

Mia moglie manda loro i suoi amichevoli saluti.

Coi miei ossequi, voglia gradire la mia cordiale stretta di mano.

Dev$^{mo}$ P. Fassini Camossj[103] - Via Finanze 11$^{bis}$ - Torino

## 74

Gian Antonio Maggi a [Paolo] Fassini Camossi

Pisa 8 Aprile 1910

Egregio Capitano,

L'equazione $x^4–ax–b=0$, essendo del 4° grado è sempre risolubile algebricamente, per radicali quadratici e cubici. Ella può trovarne le formole di risoluzione nei trattati di Analisi Algebrica o Algebra Complementare: per esempio nel Trattato Elementare sulla Teoria delle Equazioni del Todhunter tradotto da Battaglini (Napoli, Pellerano). Ai §§180 e seguenti.

Queste formole però riescono, per la loro complicazione, di scarsa utilità pratica, e, nei problemi speciali, in cui a e T siano numeri dati, conviene piuttosto ricorrere alla risoluzione numerica per approssimazione. Anche per questo Ella può consultare i suddetti trattati e per esempio, il Todhunter-Batttaglini a §§211 e seguenti.

Dai Suoi pochi cenni non posso rilevare abbastanza il problema a cui si riferisce la equazione in discorso. Me ne scriva qualche parola di più, se crede ch'io possa aggiungere qualcosa di utile in proposito. Sarà sempre per me un vero piacere stare con Lei in corrispondenza.

Intanto, lieto di codesta prima occasione, Le ricambio i migliori saluti, da parte anche de' miei, e La prego di presentare alla gentile Baronessa i nostri rispetti.

Ella aggradisca in particolare da me una cordiale stretta di mano, e mi creda sempre

Aff$^{mo}$ Suo G.A. Maggi.

Al Cap. B$^{ne}$ Fassini Camossi - Via Finanze - Torino 11$^{bis}$

---

[103] Si tratta probabilmente del Barone Paolo Fassini-Camossi, marito di Adele Sella, parente di Quintino.

**75**

Gian Antonio Maggi a **Felice**[104]

Pisa 12 Maggio 1923

Carissimo Felice,

Senza aspettare la cara tua del 10,[105] avevo approfittato della Ascensione per riprendere la teleferica, e ti trascrivo quello che troverei per rispondere alle domande che mi fai colla lettera medesima; valendomi dei tuoi stessi simboli.

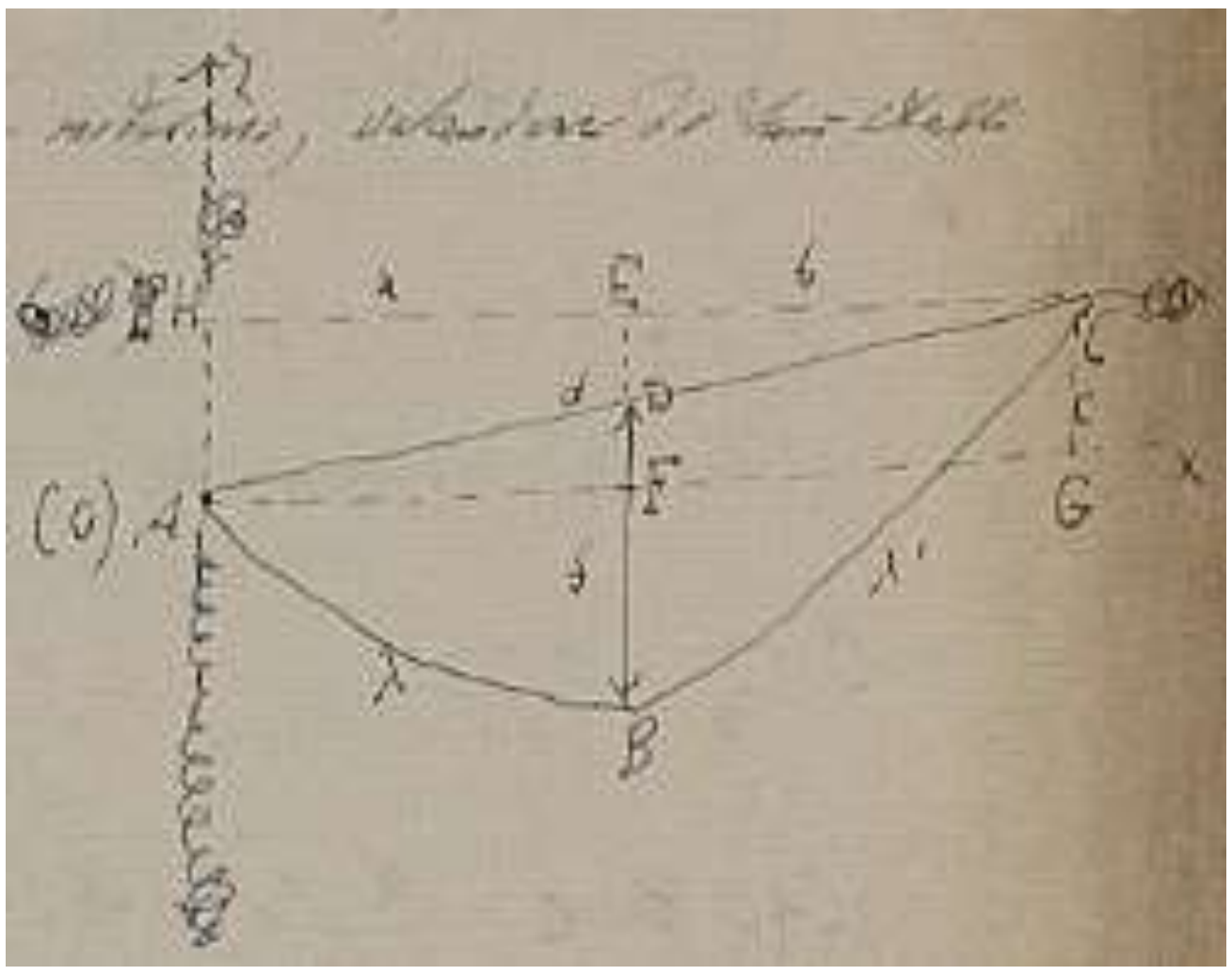

I dati da te indicati sono, in ordine alfabetico,

(1) a, b, c, f.

Con essi risulta dato

(2) $\alpha=a+b$,

ed anche

(3) $\alpha=\sqrt{\alpha^2+c^2}$.

Inoltre dal triangolo ACG si ricava DF:CG=AF:AG, ossia

$$DF=\frac{AF}{AG}CG=\frac{a}{\alpha}c$$

per cui, posto FB=e

(4) $e=\left(f-\frac{c}{\alpha}a\right)$.

Così, dai dati si ricava l'ordinata −e del punto B, rispetto alla coppia d'assi indicata nella figura, l'ascissa essendone la a direttamente data.

Le equazioni delle parabole, da sostituire colla nota approssimazione alla catenaria AB e BC, sono, posto peso specif. lineare p=1,

$$(5)\ z=z_0+\frac{(x-x_0)^2}{2H}, \qquad (6)\ z=z_0'+\frac{(x-x_0')^2}{2H},$$

dove $x_0$, $x'_0$ sono le ascisse dei vertici, H indica la grandezza della torsione orizzontale, e $z_0$, $z'_0$ sono incognite come i rimanenti parametri.

Da (4) ricaviamo manifestamente

$$0=z_0+\frac{x_0^2}{2H}, \qquad -e=z_0+\frac{(a-x_0)^2}{2H},$$

donde

$$+e=\frac{x_0^2}{2H}-\frac{(a-x_0)^2}{2H}$$

e sviluppando, poi riducendo,

$$+e=\frac{2ax_0-a^2}{2H}$$

ossia

$$(7)\ 2(-eH+ax_0)=a^2.$$

Da (6) analogamente si ha

[104] Non meglio identificato.

[105] Non si trova nel Fondo.

$$c = z'_0 + \frac{\left(\alpha - x'_0\right)^2}{2H}, \qquad - e = z'_0 + \frac{\left(a - x'_0\right)^2}{2H}.$$

donde

$$c + e = \frac{\left(\alpha - x'_0\right)^2}{2H} - \frac{\left(a - x'_0\right)^2}{2H} = \frac{(a+b)^2 - 2(a+b)x'_0 + x'^2_0 - a^2 + 2ax'_0 - x'^2_0}{2H}$$

da cui, posto

$$(8)\ x'_0 - a = x''_0 = \frac{b^2 + 2ab - 2bx'_0}{2H} = \frac{b^2 + 2b\left(a - x'_0\right)}{2H}$$

Si ricava, tenuto conto di (2),

$$(9)\ 2\left((c+e)H + bx''_0\right) = b^2$$

Per l'area, ricorreremo alla formola approssimata

$$ds = \left(1 + \frac{1}{2}\left(\frac{dz}{dx}\right)^2\right)dx$$

(ottenuta arrestando al secondo termine lo sviluppo di $\left(1 + \left(\frac{dz}{dx}\right)^2\right)^{\frac{1}{2}}$), donde, con

$$z = z_0 + \frac{(x - x_0)^2}{2H}, \qquad \frac{dz}{dx} = \frac{x - x_0}{H},$$

si ha, integrando, per la parabola (5), contando l'asse S dal vertice

$$(10)\ S = (x - x_0)\left(1 + \frac{(x - x_0)^2}{6H^2}\right),$$

e analogamente per la parabola (6), contando l'asse S' dal suo vertice

$$(11)\ S' = \left(x - x'_0\right)\left(1 + \frac{\left(x - x'_0\right)^2}{6H^2}\right).$$

Indichiamo con λ l'arco AB, e con λ' l'arco BC. Per (10)

$$(a - x_0)\left(1 + \frac{(a - x_0)^2}{6H^2}\right) + x_0\left(1 + \frac{x_0^2}{6H^2}\right) = \lambda,$$

e per (11)

$$\left(b - x''_0\right)\left(1 + \frac{\left(b - x''_0\right)^2}{6H^2}\right) + x''_0\left(1 + \frac{x''_0}{6H^2}\right) = \lambda'.$$

La prima, sviluppando e riducendo, diventa

$$(12)\ a\left(1 + \frac{1}{6H^2}\left(a^2 - 3ax_0 + 3x_0^2\right)\right) = \lambda$$

e la seconda

$$(13)\ b\left(1 + \frac{1}{6H^2}\left(b^2 - 3bx''_0 + 3x''^2_0\right)\right) = \lambda'.$$

$$\lambda' = l - \lambda$$

Le (7), (9), (12) e (13) sono quattro equazioni Le (7), (9), (12), (13) sono quattro equazioni fra le λ e λ', che sono la tua incognita, e le $x_0$, $x''_0$ e H. Esse possono servire a determinare λ e λ', se si suppone dato H, che tu non nomini tra i dati: oppure, se, contando H come incognita si suppone dato l=λ+λ'. E veramente a me sembra singolare che non si consideri come data la lunghezza totale (l) del filo. Ne sei ben sicuro?

Una quinta equazione si può aggiungere scrivendo che la somma del peso P sospeso a B e del filo, rappresentato da l (perché si è posto p=1) è uguale alla somma delle componenti della tensione nei due estremi A e C, secondo la verticale volta in alto (asse delle x). Da $T\frac{dz}{ds}=H_x$ si ha $T\frac{dz}{ds}=H\frac{dz}{dx}$. Quindi nel punto B $T\frac{dz}{ds}=(x-x_0')$, e nel punto A, mutando il segno per introdurre la tensione esercitata dal punto A sul filo, in luogo della tensione applicata dal filo allo stesso punto A, $-T\frac{dz}{ds}=-H(0-x_0)=+x_0$. Si trova così $[\alpha-x_0'+x_0]=P+l$.

Ma meglio ricorreremo alle formole esatte $\frac{dz}{dx}=\sin h\frac{x-x_0}{H}$ per la prima catenaria e $\frac{dz}{dx}=\sin h\frac{x-x_0'}{H}$, per la seconda. Sviluppando, e fermandoci al secondo termine, abbiamo, col mutamento di segno per la pressione in A

$$(14)\ x_0+\frac{x_0^3}{6H^2}+\alpha-x_0'+\frac{(\alpha-x_0')^3}{6H^2}=P+l.$$

Con questo si hanno tante equazioni quante sono le incognite, $\lambda$, $\lambda'$, $x_0$, $x'_0$, H: ma introducendo il carico P.

## 76

Gian Antonio Maggi a **Bruno Finzi**

Da cartolina postale al Prof. Finzi 16 Sett 1934[106]

Il prof S.[107] mi accompagnò cortesemente la Nota con una cartolina, alla quale risposi approvando le piccole varianti come conformi agli espressi miei desiderii. Avendola poi riletta, credo ormai che l'A si è dato pensiero di confinare nelle note a piede di pagina la ricerca della paternità del noto soggetto, per la ragione che la conclusione della sua "Osservazione"[108] è che la prole non ha corrisposto, come si suol dire, all'insegnamento che si faccia il genitore. Ad ogni modo, la nostra Nota è destinata a mettere tutte le cose interamente in chiaro. Non mi spiego come l'Istituto deve aver trattenuto gli estratti
...
Sono poche copie, ne rimetterò la diffusione a Lei, ché la parte Sua è tale e tanta che io mi sento autorizzato, senza offendere la modestia, a chiamar quello un assai bel lavoro[109]
...

[106] Si trovava nelle *Rusticationes*.
[107] Si tratta probabilmente del prof. Luigi Sona. Si veda la corrispondenza con Levi-Civita.
[108] Il titolo dell'articolo di L. Sona è: "Un'osservazione riguardante la propagazione delle onde elettromagnetiche armoniche".
[109] Potrebbe riferirsi alla Nota: "Condizioni d'esistenza delle onde elettromagnetiche armoniche" (in collaborazione con B. Finzi), *Rendiconti del R. Istituto Lombardo*, v. 67, 1934, pp. 331 e 363-370.

**77**
**Federigo Giordano** a Gian Antonio Maggi[110]
[biglietto intestato: Scuola di ingegneria (R Politecnico) Milano - Costruzione delle macchine]

Copia per il Not. Prof. Gr. Uff. G.A. Maggi - con preghiera di indicare le eventuali modificazioni - con deferenti cordiali saluti.[111]

F. Giordano.

Tel. 292.039

**78**
Gian Antonio Maggi a Federigo Giordano

Casa - Martedì
(18 Dicembre 1934)

Egregio e Gentilissimo Collega.

La ringrazio, prima di tutto, per la Sua Memoria sugli scritti del Colombo, che ha la bontà d'inviarmi, tanto più gradita e interessante col grato e riverente ricordo che serbo di lui.

Per la relazione, non potrei altrimenti che trovarla di mio perfetto aggradimento, salvo permettermi due piccole osservazioni.

Crederei opportuno, in primo luogo, riprodurre dalla mia relazione personale - non si sa come smarrita all'Istituto - che il lavoro, quantunque, di necessità, d'indole compilativa, per la composizione, l'impostazione, non comuni concetti, e non comuni metodi di trattazione, reca l'impronta personale.

In secondo luogo, pensando all'impressione che potrebbe fare all'Istituto, nostro giudice, l'appunto sui lavori italiani non abbastanza considerati, sulla portata del quale io certamente mi rimetto a chi possieda l'erudizione in materia di applicazioni che mi fa interamente difetto, Le chiederei se non credesse di poterlo attenuare colla scusa dei criterii a cui l'autore si è attenuto nel prescrivere i limiti del suo svolgimento del tema.

Mi rimetto, ben s'intende, del resto, a quanto Ella crederà meglio di fare, e coi migliori e più distinti saluti La prego di credermi sempre

Suo aff. collega
Gian Antonio Maggi

Ch. Sig. Prof. Ing. Gr. Uff. Federigo Giordano

---

Concorso Kramer, 1933[112] - Giudizio sulla Memoria
"La forza nasce dal movimento, il movimento dalla forza"

La poderosa Memoria, che comincia con notizie storiche, dove compare il paradosso di d'Alembert, germe della dottrina che ha preso ai giorni nostri così vasta estensione da riuscire ardua impresa dominarla nella sua integrità, e finisce con specificate applicazioni pratiche del materiale precedentemente esplorato, si dedica con successo a quell'arduo compito contemplando il problema proposto sotto i suoi molteplici aspetti, e rendendo conto dei risultati, sia nell'indirizzo razionale, che particolarmente interessa la Scienza pura, che in quello delle svariate pratiche applicazioni. Il rigore matematico, che informa la trattazione, e le conferisce tanta chiarezza e solidità, e i più progrediti e raffinati strumenti dell'Analisi, usati per svilupparlo, conferiscono al lavoro un pregio, che vuol essere particolarmente rilevato. La costruzione, concetti non comuni, e non comuni risultati imprimono al lavoro, che, a termini del tema proposto, non potrebbe avere altrimenti che carattere, per molta parte, compilativo, l'impronta della personalità.

[110] La corrispondenza con Giordano è contenuta in una busta dedicata a: *Premio Kramer dell'Ist. Lomb. 1934*.
[111] Probabilmente questo biglietto accompagnava il bando del Premio Kramer dell'Istituto Lombardo del 1934; sia Maggi che Giordano facevano parte della Commissione (insieme a Murani, Fantoli e Zunini).
[112] Probabilmente Maggi sbaglia la data, che è, invece, 1934.

La mia opinione è che a questa Memoria, per soddisfazione del tema proposto, e pel suo intrinseco pregio, sia conferito il premio intero.[113]

**79**

Gian Antonio Maggi a **Giovanni Giorgi**

Milano 12 Marzo 1937

Caro Professor Giorgi,

Io ho adottato la definizione di Maxwell del vettore, colla quantità qualsivoglia, nel mio vecchio trattato Principii della teoria matematica del movimento dei corpi, Milano, 1896 e riprodotto nei miei Elementi di statica, Bologna 1925 come dalla recensione del Finzi nel Periodico di Matem. del 1 Gennaio 1930, che testualmente ve la riporta. Mi ha fatto piacere, a suo tempo, vederlo adottato pure da Lei, nelle Sue pregevoli Lezioni di Meccanica,[114] unico, ch'io sappia, degli autori di trattati di Meccanica, per cui non è certo superfluo ch'Ella la rimetta in onore nel suo interessante articolo del volume per Berzolari.[115]

Per vettore applicato, io non potrei passare per seguace dell'Appell, perché introduco quel termine, unito con retta d'applicazione, nel suddetto mio libro del 1896, mentre l'Appell non fa la distinzione di vettore libero, applicato ad una retta (il suo cursore, che il Clifford chiama rotore) e applicato ad un punto che nelle edizioni posteriori alla seconda, che è del 1902. Facendo una sola distinzione dal vettore libero, io preferisco "applicato ad un punto", cioè la combinazione di un vettore e di un punto, perché mi sembra il caso poi che più spontaneamente presento o uso, e pel cursore mi serve la retta d'applicazione. Per aggettivo modificativo piuttosto che qualificativo, ricordo l'osservaz. del Burali Forti per quel mio termine, ma mi sembrava veramente esagerato.

Piuttosto l'unione di una orientazione con un numero mi sembra far una certa eccezione alla definizione del vettore, se si tiene conto della giusta sua distinzione tra numero e scalare. Per vettore omnibus io mi valgo del segmento rappresentativo, che è segmento orientato, e mi serve utilmente per le costruzioni con cui definisco la somma ecc. Ben inteso facendo attenzion di *[...]* descrivere da un punto un segmento avente la lunghezza, la direzione e il verso del segmento rappresentativo, e *[*non punto*]* lo stesso segmento rappresentativo che non saprei come prendere in mano.

---

113 Nella busta è presente anche lo stampato con la valutazione finale di conferimento del Premio

114 G. Giorgi, *Lezioni di Meccanica Razionale*: volume I, del 1931 e volume II del 1934, P. Cremonese, Roma.

115 G. Giorgi, "Spigolature di calcolo vettoriale", in: *Volume degli scritti matematici offerti a Luigi Berzolari*, Pavia, 1936, pp. 637-641.

Credo che non Le spiacerà questo mio richiamo del Suo articolo, che ho letto con vero interesse e piacere. Intanto lieto dell'occasione di ricordarmi, La prego di aggradire i miei migliori saluti e di credermi sempre

Suo aff$^{mo}$ collega
Gian Antonio Maggi.

*[sullo stesso foglio]*
Ch. Prof. Giovanni Giorgi della Facoltà d'Ingegneria della R. Università di Roma
Corso V.E. 39.

Cartolina Milano 13 Marzo 1937
...
Lei non manca di dire che il numero è una specie particolare di scalare e quindi nessuna eccezione alla definizione di vettore.
...

**80**

Giovanni Giorgi a Gian Antonio Maggi
[busta intestata: Prof. Ing. G. Giorgi - 39, Corso V.E. Roma e
indirizzata a: Illustre Prof.$^{\underline{re}}$ Gian Antonio Maggi - Corso Plebisciti 3 - Milano;
carta intestata: Prof. Ing. G. Giorgi - Telefono interurbano 63.406 - (Roma 17) -
39, Corso V.E. (Palazzo Ferretti)]

13 marzo 1937

Illustre Prof Maggi,

Sono molto onorato che quel mio modesto scritto abbia richiamato la Sua attenzione e mi abbia procurato il piacere di ricevere una Sua lettera.

Sì, le definizioni di scalare, vettore, orientazione, che si trovano nei Suoi Elementi di Statica sono bene in accordo con le idee mie: nel momento che ho scritto quelle "Spigolature" non le avevo sott'occhio: altrimenti le avrei citate per appoggiarmi all'autorità del Suo nome.

Quanto al vettore puro: - quando io dico che «scalare» non è sinonimo di «numero», intendo però che il numero sia un caso particolare dello scalare. Perché non dovrebbe esserlo? È una grandezza come un'altra, di dimensioni [L°M°,T°]. Quindi sta bene prendere come vettore astratto il numero puro associato a un'orientazione, di guisa che moltiplicando il vettore astratto per [L], per [$LT^{-1}$] se si hanno i diversi vettori concreti. Ciò non impedisce di prendere come vettore tipico, a scopo grafico, il vettore-linea, l'unico che si può disegnare sulla lavagna.

Venendo poi ai vettori applicati: -

Non ho mai pensato a assumere Lei fra i seguaci di Appell. I metodi che Lei segue, traggono se mai, origine dalle idee di Mach, ma Lei se ne è discostato sostanzialmente (come io ebbi già a rilevare nei miei articoli sulla Massa) fin dai primi Suoi libri; e più ancora se ne è discostato nelle redazioni successive, organizzando un'impostazione che è tutta sua propria. La maggior parte però degli altri autori di lingua latina sono "Appellisti".

I suoi sono vettori applicati a un punto? Come definizione sì, come trattazione no, perché Lei li tratta quasi sempre come se fossero forze applicate a un corpo rigido, si preoccupa della retta d'azione, trova i momenti, il sistema risultante, etc. Mi sembra che i vettori localizzati a un punto si presentino solo nella meccanica dei solidi deformabili e dei fluidi, nel qual caso non si compongono fra loro se non sono proprio applicati allo stesso punto, e non danno quindi luogo a studi di composizione o decomposizione.

Venendo alla parola "rotore", che ho imparata da Clifford: mi sembra debba essere applicata a quell'ente che rappresenta l'atto di un moto rotatorio di un solido rigido: quindi solo nell'$S_3$ equivale al cursore. Nell'$S_2$ il rotore è un punto caricato (un centroide) e nell'$S_n$ è un multivettore di ordine $n-2$, localizzato. Quindi è bene mantenere in uso entrambe le parole "cursore" e "rotore" nei loro significati propri.

Di nuovo, si abbia le mie grazie, Illustre Professore, e i miei ossequi deferenti.

Suo dev$^{mo}$
G. Giorgi

P.S.- Le interessa sapere una novità? Il Comit. Elettrotecnico francese, che pur avendo votato le mie unità, non voleva la razionalizzazione del $4\pi$ e aveva significato di volersi opporre fino alla morte a qualunque tentativo di razionalizzazione delle formule, si è ora convertito, e ha emesso una deliberazione per l'adozione delle eq. e unità razionalizzate. Con questo è tolto un grande ostacolo, e nulla impedirà che la decisione per la razional.$^{\underline{ne}}$ possa venir adottata internazionalmente.
A parte Le invio alcuni miei lavori ultimi

GG

**81**

Gian Antonio Maggi a **Francesco Giudice**

Sulla possibilità di funzioni di stesso valore, equivalenti, non identiche.
Nota del prof. F. Giudice nei Rendic. del Circolo Matem. di Palermo. Tomo III (1889).

Caro Giudice,

...

Alla domanda "Se vi possono essere funzioni dello stesso valore non identiche" io risponderei colla pregiudiziale che, essendo una funzione dai valori che acquista nei singoli posti sostanzialmente definita, dall'eguaglianza di quei valori deve necessariamente seguire l'eguaglianza di tutte le proprietà.

Da che dipende per esempio l'esistenza e il valore del coefficiente differenziale in un posto, se non dai valori che prende la funzione in un intorno di quel posto? E allora, se due funzioni, in un intorno del posto, ricevono gli stessi valori, come potrebbe accadere che l'una ammettesse il coefficiente differenziale e non l'altra, o l'avessero diverso?

Perciò dall'eguaglianza

$$\int_a^x \varphi(x)dx = \int_a^x f(x)dx$$

dove f(x), φ(x) sono le due funzioni da te considerate, io deduco, per ogni posto x dell'intervallo (a,b):

$$\frac{d}{dx}\int_a^x \varphi(x)dx = f(x).$$

Vuol dire che nei posti del gruppo di prima specie e suoi punti limiti il coeff. diff. di $\int_a^x \varphi(x)dx$ non sarà φ(x), ciò che s'accorda con la circostanza che φ(x) in quei posti è discontinua.

Perciò pare che tu escluda l'esistenza della derivata di $\int_a^x \varphi(x)dx$.

Ma, se una funzione φ(x) in un intervallo (a,b) è integrabile e discontinua, per questa circostanza, solo non si può asserire la necessaria esistenza della derivata di $\int_a^x \varphi(x)dx$, e, supposto che esista, la sua eguaglianza con φ(x).

...

Messina 6 Gennajo 1890.

**82**

Gian Antonio Maggi a **Eduardo Gugino**
[carta intestata: R. Università di Milano - Istituto Matematico -
Via C. Saldini, 50 (Città degli Studi) - Il Direttore]

Lettera al Prof. E. Gugino 6 Dicembre 1930

Sped. G.A. Maggi
Corso Plebisciti, 3, Milano.
Ch.$^{mo}$ Sig. Prof. Eduardo Gugino della R. Università di Messina
R. Università Facoltà di Scienze.

Milano 6 Dicembre 1930

Egregio Professor Gugino

Leggo la Sua Nota nel Numero oggi ricevuto dei Rendiconti dei Lincei,[116] e, prima di tutto, La ringrazio di ricordare, più di una volta, il mio nome, a proposito delle equazioni dinamiche, che col mio nome ha voluto chiamare la benevolenza del Levi-Civita. Per quanto al problema, che forma l'oggetto della Sua Nota, esso si trova posto e risoluto nel §82 della mia Dinamica dei Sistemi. Il metodo non è precisamente lo stesso, non però meno semplice, ma il principio finisce per essere lo stesso, perché io, assunte, per gli elementi delle reazioni vincolari, le note espressioni come funzioni lineari omogenee dei parametri λ di Lagrange, mostro come si possa stabilire, valendosi delle equazioni che traducono i vincoli, per sistemi olonomi o anolonomi, un sistema di equazioni lineari, dove le λ fungono da incognite, atto a fornirle in funzione del tempo, della posizione e dell'atto di moto.

Mi è grata l'occasione per inviarLe i migliori saluti, e per pregarLa di credermi sempre

Aff. Suo
Gian Antonio Maggi.

[116] E. Gugino, "Sul problema dinamico di un qualsivoglia sistema vincolato ridotto all'analogo problema relativo ad un sistema libero", *Rendiconti dei Lincei*, s. VI, v. XII, 1930, pp. 307-312, con presentazione di Levi-Civita (si veda qui la sua lettera del 2/1/1931).

**83**
Eduardo Gugino a Gian Antonio Maggi
[busta indirizzata a: Chiarissimo Signor Professore Gian Antonio Maggi
della R. Università - Corso Plebisciti 3 - Milano;
sul retro: Spedisce: E Gugino. - Via Dante 47 - Palermo.]

Palermo 14 Dic. 1930

Chiarissimo Prof. G.A. Maggi.
Milano.

Dalla Sua pregma del 6 c.m. apprendo che il problema, che forma oggetto di mia recente nota lincea, si trova posto e risoluto nel §82 della Sua "Dinamica dei sistemi". Effettivamente debbo riconoscere, avendo già presa visione del Suo trattato, che il problema della determinazione delle reazioni in funzione delle forze motrici impresse, della posizione, dell'atto di moto e del tempo, è stato da Lei nettamente risoluto con referenza ad un sistema di corpi rigidi, comunque vincolati. Appare quindi manifesta la priorità del risultato da Lei ottenuto, con procedimento diverso, e di tale priorità sono pronto dare pubblicamente atto alla prima occasione. Ho soltanto il rammarico di non avere consultato il Suo testo prima di avere pubblicato il mio lavoro; mi resta però la soddisfazione morale di essermi, per così dire, incontrato con Lei, avendo sentito, al pari di Lei, la necessità di quel significativo complemento ai fondamenti della Meccanica analitica.

La prego, Chiarissimo Professore, gradire i sensi della mia più deferente stima.

Suo devotissimo
E Gugino

Palermo. Via Dante 47

**84**

**Camillo Guidi** a Gian Antonio Maggi
[carta intestata: Associazione fra gli Utenti di Caldaie a vapore
avente sede a Milano - Il Direttore][117]

Milano 3 XI 1908

Carissimo,

Ti scrivo in treno. Abbi pazienza se scrivo male.

Riceverai il Giornale l'Industria dove leggerai un articolo dell'Ing. Barzanò sull'insegnamento nei Politecnici

Faresti cosa ottima scrivendo al Barzanò per esporgli i tuoi concetti in proposito all'insegnamento del calcolo e della meccanica razionale.

Allo stato attuale delle cose una tua lettera da presentare al Colombo, oppure scritta al Colombo stesso, e magari da pubblicare farebbe bene assai a noi che ti desideriamo.

Il Barzanò è effettivamente d'accordo col Colombo in quelle idee e il Colombo apprezza assai assai il Barzanò.

Nota la stangata al Jung, quella al modo di insegnamento attuale del calcolo a Milano. = Sviluppa te i concetti che mi hai esposti = Vedi *[Jorini]* nel cenno ai sistemi elastici.

Scrivi, se ci vuoi bene.

Saluti affettuosi a tutti

tuo Guidi

volta

Assai probabilmente in settimana verrò a Pisa per fermarmi poche ore; vorrò vederti

Guidi

**85**

Gian Antonio Maggi a Camillo Guidi

Lettera al Prof. C. Guidi sul Cirenei[118] (18 Maggio, 1930).

Milano 18 maggio 1930

A S.E. Prof. C. Guidi della R. Accademia d'Italia - Viale della Milizia 16. - Roma

Gentilissimo Professore,

Ho subito risposto alla gradita Sua in data 25 Marzo scorso,[119] e conto che Ella avrà potuto vedere ai Lincei la nostra Relazione suoi nuovi trovati della scienza di quel disgraziato di Cirenei # (v. in fine)

Forse Ella non sa che, qualche tempo dopo, egli fece stampare, in un giornale di Roma, che la Commissione, formata dai professori Levi-Civita e Maggi, aveva reso parere favorevole sulla sua Memoria, presentata ai Lincei, per l'inserzione negli Atti. In seguito a che, avendo il cancelliere Mancini fatto pubblicare, nello stesso giornale, a titolo di rettifica, che la Commissione aveva conchiuso colla semplice proposta dell'invio della Memoria agli Archivii dell'Accademia, il Cirenei, tornando alla carica, fece pubblicare la spiegazione di questo discorso, …per chi non sa di latino: e cioè che voleva dire che l'Accademia aveva trovato la Memoria così bella da volerla conservare nei proprii Archivii!

[117] Maggi scrive su un foglio di questa lettera: *Sull'articolo dell'ing. Barzanò "Alcune considerazioni pel riordinamento dell'insegnamento tecnico superiore" (v. Industria 1 Novembre 1908.).* Cfr. corrispondenze Barzanò e Colombo.

[118] Si tratta di Egisto Cirinei; si veda anche la nota alla lettera #8.

[119] Non presente nel Fondo.

Formata poi dal Ministero una Commissione, presieduta dal Volterra, che si destreggiò in modo da non riunirsi mai (io e Levi-Civita, nominativi, cominciammo per esimercene, allegando il pronunciato giudizio), si lesse, nei giornali dell'Urbe, che cotesta Commissione si mettesse alacremente all'opera, perché si avvicinava la sessione autunnale degli esami, e i giovani erano dai professori chiamati a rispondere sopra una Meccanica, che il Cirenei dimostrava falsa!...

Non so quali nuovi teoremi il Cirenei abbia potuto scoprire, dopo la nostra Commissione lincea. A quel tempo, se la prendeva col principio della leva, conciato da lui sotto la novissima forma: "ciò che si guadagna in forza, si perde in tempo". E lo demoliva, col mostrare che una leva metallica, girevole intorno ad un asse orizzontale, dove la resistenza era rappresentata da un contrappeso, e la potenza da pesetti, applicabili a crescenti distanze dal fulcro, in modo da guadagnare in forza, col crescere questa distanza, girava di un certo angolo, nel medesimo tempo, comunque si facesse crescere questa distanza, per diminuire il peso occorrente per l'equilibrio. L'esperienza avrebbe dovuto essere di fatto una confutazione della nota obbiezione al da mihi ubi consistam terramque movebo, consistente nell'osservazione che, concepita montata la Terra sopra una leva, per girare cotesta anche di un piccolissimo angolo, con una forza accessibile all'uomo, si trova occorrente un tempo infinitamente grande. E la ragione ne è chiaramente che, per fare equilibrio al peso ideale della Terra, con un peso accessibile all'uomo, occorre assegnare a questo un braccio estremamente grande, con che diventa estremamente grande il momento d'inerzia, rispetto all'asse, del mobile, che contiene già il momento d'inerzia dal Globo Terreno. Perciò l'esperienza del Cirenei avrebbe, di fatto, dovuto dare tempi sempre più lunghi col guadagnare in forza, coll'accrescere il braccio; e il motivo, per cui egli trovava sempre lo stesso tempo era semplicemente questo, che lo spostamento del piccolo peso aggiunto, sulla massiccia leva, possedente già, per proprio conto, un forte momento d'inerzia, non produceva variazioni di tempo rilevabili alla abbastanza grossolana osservazione. Naturalmente, io mi sono ben guardato di fare questo discorso al Cirenei, che del momento d'inerzia non sospettava neppure l'esistenza, e del prestarsi a fargli intendere la ragione non avrebbe inteso altro che approvazione delle sue incongruenze.

Come il candido principio della leva abbia potuto diventare oggetto di tanto scempio è cosa da fare strabiliare. Meno però, se si riflette che, col cervello scalcinato del Cirenei, bisogna portare in conto l'incoraggiamento e l'appoggio che gli sono prestati, colpa quell'ignoranza, che non resta, pur troppo, confinata ai più bassi strati della nostra Società, e a quella certa ribellione alla Scienza così detta ufficiale, che non è altro che una particolare manifestazione della stessa ignoranza.

Se non che mi accorgo di avviarmi a trattenerLa un tempo paragonabile con quello previsto per l'esperimento di Archimede. Per cui faccio punto, chiedendoGliene scusa, lieto del pensiero di vederLa presto a Roma, La prego di aggradire le migliori cose nostre, per sé e per la Sua famiglia, e di credermi sempre

Suo dev.$^{mo}$ collega
Gian Antonio Maggi

# Avendone veduto il tenore, Ella si meraviglierà certamente di sentire che il Cirenei, con sua lettera ricevuta jeri l'altro, mi prega di comunicare all'Officina Galilei, che ha incaricato della costruzione dei suoi istrumenti, il mio parere, conforme a quello del prof. Levi-Civita,[120] da me, in tale conformità, asserito, che i suoi strumenti sono idonei per l'insegnamento della Meccanica. E con questo tanto più m'interessa di conoscere le Sue conclusioni, sull'ultima offensiva del Cirenei all'Accademia d'Italia

**86**

**Nicolò Guzzardi** a Gian Antonio Maggi

[raccomandata con busta indirizzata: Al Comm Ufficiale Maggi Dottore Gian Antonio Professore di Fisica Matematica nella R. Università di Milano]

Palermo 20-1-928.

Gent$^{mo}$ Comm$^{re}$ Professore Maggi,

Leggendo una sua memoria "Sul raggio di luce nell'ottica fisica"[121] e trovandola molto piacevole, e di mio gusto mi sono accinto a studiarla specie nella parte fisica per la quale ho molto tendenza, ma ho dovuto però incontrare alla fine delle difficoltà e avendo il vivo desiderio di portarla al completo questa sua memoria vengo alla S.V. Illma perché voglia chiarire quanto a me sembra oscuro chiedendo umilmente perdono del disturbo.

Per quanto riguarda la teoria elastica e la nota definizione, definizione in virtù della quale, sotto la condizione $\cos\widehat{Sx}\cos\widehat{Vx}+\cos\widehat{Sy}\cos\widehat{Vy}+\cos Sz\cos Vz=0$ e così via... Sono arrivato alla conclusione basandomi sulla formola che mi da lo sforzo esercitato sull'elemento normale cioè

$$\Phi_V=X_Vi+Y_Vj+Z_Vk$$

dove i j k sono tre vettori; $X_V$, $Y_V$, $Z_V$ sono gli sforzi dovuti agli elementi normali secondo la direzione V.

Le componenti di questi sforzi mi sono dati da:

$X_V=X_X\alpha+Y_Y\beta+Z_Zj$

e così per la $Y_V$ e la $Z_V$; dove le $\alpha$, $\beta$, j sono le componenti degli spostamenti, cioè, sono le $\xi$, $\zeta$, $\varphi$.

Non so rendermi ora ragione nel caso generale e precisamente perché

$X_X=-2kx_X$ e $X_Y=-kx_Y$

Da dove si ricava quel $-2kx_X$? E perché il segno meno? E perché dovrà essere

$X_X+Y_Y+Z_Z=0$?

Desidererei sapere da quale formola potrei ricavare la $x_X$ e la $Y_Y$.

Voglia esaudire questi miei desideri per appagare il mio piacere nell'aver capita questa sua memoria fino alla fine.

Ringraziandola la riverisco

Devotissimo

Capitano Nicolò Guzzardi

[120] Si vedano, su tale argomento, le lettere #101 e 102.

[121] *Rendiconti del R. Istituto Lombardo,* s. II, V. LIX, 1926, pp. 688-656.

**87**

**Carlo Hoepli** a Gian Antonio Maggi
[timbro: raccomandato con 2 allegati;
Sapere - Amministrazione (Ulrico Hoepli - Milano);
l'angolo superiore sinistro è strappato]

Milano, 7 Marzo 1936-XIV.

Ill.mo Signor
Prof. Dott. G.A. Maggi
Milano Corso Plebisciti 3

Ho ricevuto qualche giorno fa da un anonimo un adattamento di un articoletto di volgarizzazione del Prof. russo Perelmann su quanto avverrebbe in un mondo senza gravità. Le unisco l'articoletto stesso.

L'argomento mi ha interessato e penso che avrebbe potuto interessare i lettori di "SAPERE". Ma rileggendo l'adattamento mi sono venuti alcuni dubbi, per chiarire i quali ho pregato tre diverse persone e precisamente il prof. Marcolongo, il Prof. Finzi e l'Ing. Carlo Rossi (a Lei noti) di darmi un loro giudizio personale.

Questi giudizi sono stati da me fedelmente raccolti nello accluso foglio.
Avrei voluto riassumere in SAPERE tutta questa simpatica discussione su di un argomento affascinante anche per un grande pubblico. Ma per sgravio di coscienza e desiderio di equanimità verso il Perelmann vorrei essere confortato dal di Lei autorevolissimo giudizio che ben s'intende resterà riservatissimo: sia sull'articolo stesso, sia sulla opportunità ed esattezza delle critiche stesse.

La prego di scusare la noia che Le dò e di gradire i miei più sentiti ringraziamenti.

Mi creda
dev.mo
Carlo Hoepli

**88**

Gian Antonio Maggi a Carlo Hoepli[122]

Minuta della lettera al Sig. C. Hoepli in risposta alla sua - acclusa - sopra un articolo "Una colazione in un mondo senza gravità"

Ill.mo Sig.r Comm. Carlo Hoepli.
S.M.

Milano, 11 Marzo 1936.

Pregiatissimo Signor Commendatore,

Volentieri ho letto ed esaminato l'articolo "Una colazione in un mondo senza gravità"; ed eccoLe, secondo il desiderio da Lei accennatomi, le mie osservazioni che, in buona parte, collimano con quelle de' miei colleghi.

[122] In due copie di minuta datate 10 e 11 marzo. Qui si trascrive la seconda, che reca la scritta in matita: *Ripresa dalle mani del Comm. Hoepli il 14 III e sostituita con lettera ridotta alle conclusioni finali*. Questa ultima lettera non è presente nel Fondo, ma probabilmente era rielaborata a partire dalle frasi sottolineate (in matita) nella presente.

E, prima di tutto, “mondo senza gravità”, per farvi la storia di una colazione, bisogna reputare, secondo me, oggetti su cui non agiscono sensibilmente forze esterne, tranne le supposte forze manuali dei personaggi chiamati in scena. Quindi da lasciar in disparte l’attrazione newtoniana, esercitata su ciascun oggetto dei rimanenti: la quale, per la sua piccolezza, non potrebbe produrre effetti percepibili nel tempo di una colazione.

Ciò premesso nella storia dei cioccolatini, che escono infuriati dal sacchettino, è solo da far eccezione alla direi sfuggita scossa leggiera. Con questa i cioccolatini non avrebbero potuto moversi che lentamente, in ragione della piccola velocità impressa. Sarebbe la volta dei “boni boni” del Marcolongo. Con una scossa vigorosa, non c’è nulla da opporre al descritto spettacolo.

Su Benedetto, a cui salta in gola il cioccolatino piovutovi, il mio egregio collega Finzi cita certi acrobati, che bevono colla testa, in giù sul palcoscenico. Ma, acrobati! Certo, la gravità non è indispensabile per la deglutizione. Mi sembra però che il suo concorso sia sufficiente per sentire da chiunque un imbarazzo, la volta che mancasse, dopo aver ingojato un cioccolatino.

Che, senza gravità, l’acqua non uscisse dalla bottiglia sta bene. Per quanto, cavata, a formare una palla attaccata al collo della bottiglia, direi che si può sfuggire a critiche, prendendo in parola palla, per non intendere necessariamente una sfera.

Il bicchiere che “si librava nell’aria qua e là” non può voler dir altro che stava sospeso dove lo poneva la mano.

Il rivestimento del bicchiere, da parte dell’acqua, per l’adesione, non contrastata dalla gravità, mi sembra, com’è esposto, non prestarsi a ragionevole interpretazione.

Col fornello, che non sta acceso, ricadiamo un po’ in Benedetto che si strozza col cioccolatino. Certo, in diversi modi si possono concepire liberati gli ingombranti gas inconsistenti, prodotti dalla combustione. Ma non potrebbe non farsi sentire la mancanza della ascesa di detti gas, in conseguenza della leggerezza conferita dall’elevata temperatura nativa, in seno ad una atmosfera possedente una pressione subordinata alla gravità. Badiamo che, nel mondo della nostra colazione, un turacciolo di sughero dal fondo di un bicchier d’acqua, non verrebbe a galla. Ci vorrebbe altro per tirarlo su come per spazzare i suddetti gas.

E simile osservazione mi sembra potersi fare, a proposito dell’acqua che non bolle, per la mancata ascesa della parte calda e leggiera nella parte fredda e pesante.

A questo punto si tronca, col foglio 3, la copia mandatami dell’articoletto.
Non so se c’è dell’altro. A noi intanto resta da domandarci come agiscano, senza gravità, gli organi del Benedetto, del cuoco, del professore e compagnia.

Ad ogni modo, quanto precede mi sembra sufficiente per giudicare lo scritto come leggiero, da considerazioni destinate a far colpo sul lettore, le quali reggono, in generale, all’ingrosso. Articolo da leggere per passatempo: se non è per sottoporlo ad un esame, come quello da Lei desiderato, che obbliga a rievocare i principii della Fisica, per cercarvi in quanto sono rispettati… o strapazzati.

Aggradisca con questo i miei migliori saluti, e mi creda sempre

Dev.mo Suo
Gian Antonio Maggi

**89**
**Giuseppe Jung** a Gian Antonio Maggi

Pro veritate - Osservazioni sull'opuscolo del sig. L.S. distribuito col fascicolo di marzo (1901) del "Politecnico".[123]

- I. -

1. Sia φ(x) una funzione di x; ciò significa che φ(x) è una quantità dipendente dalla variabile x - tale cioè che per ogni valor dato di x, assume un valore, determinato dalla natura (ossia dal significato, dalla forma) della funzione φ.
2. Se per valori differenti di x, la φ(x) assume sempre un identico valore, cosicché per x qualunque e per w qualunque, sia
$\varphi(x)=\varphi(x+w)$ (a)
si dice che φ(x) è una costante rispetto ad x e si scrive
$\varphi(x)=cost=C$ (poniamo) (b)
se C è il valor numerico di φ(x) per uno (e quindi per tutti) i valori di x. Dunque: dire che φ(x)=cost; oppure: dire che φ(x)=φ(x+w), per x e w qualunque, equivale a dire che per ogni valor dato di x, il corrispondente valore di φ(x) è sempre il medesimo; così per es. se per x=1, fosse φ(1)=3, per x=2, per x=3, ... per x=10 ... farebbe sempre φ(2)=3, φ(3)=3, ... φ(10)=3, ecc. ecc. – L'ordinata y=φ(x) di un punto M descrivente una retta AM parallela all'asse OX, e situata alla distanza OA=3, offre un esempio intuitivo di una funzione φ(x)=cost=C – In questo esempio la costante C è =3; se non fosse 3, ma un altro numero, e non si volesse precisar quale, si esprimerebbe la costante con un simbolo, per es. con C, come si è fatto sopra, e si scriverebbe

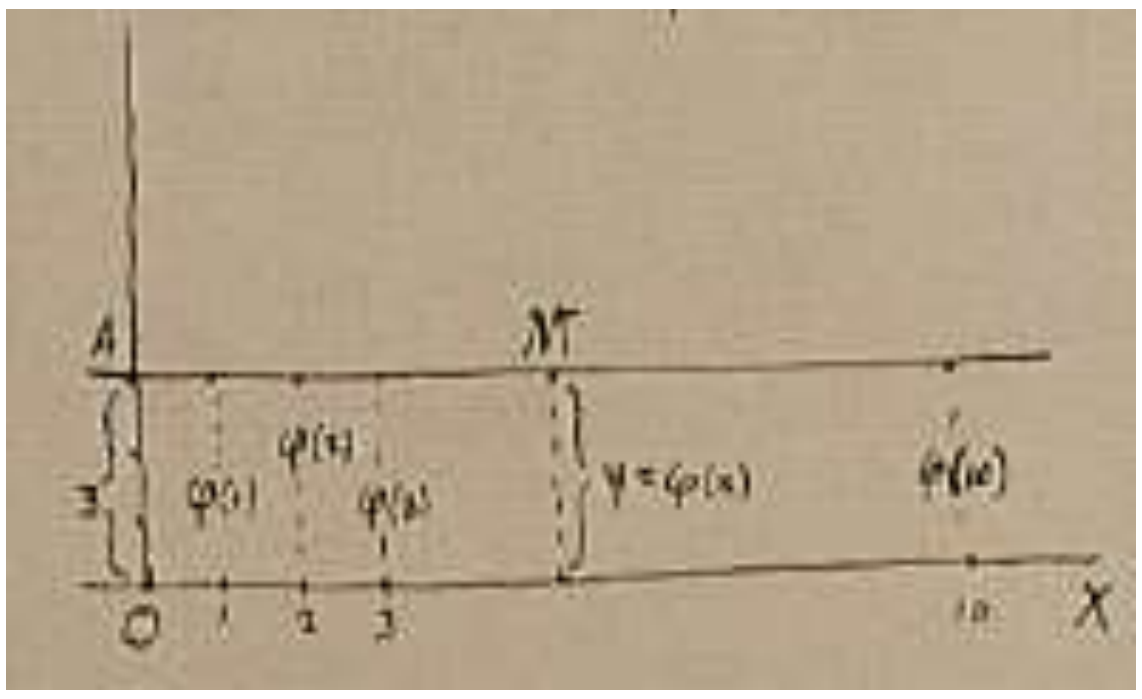

$$\varphi(x)=C \ (=\text{cost. risp. } x);$$
ove C è un parametro il cui valore non dipende dai valori di x.
3. Se per ogni valor dato di x, φ(x) assume un valore ad esso proporzionale, cosicché per x e w qualunque, sia
$\frac{\varphi(x)}{x}=\frac{\varphi(w)}{w},$ (c)
la forma o il significato di φ resta da questa condizione determinato: perché se $\frac{\varphi(x)}{x}$, per un valor particolare di w, risulta uguale a k, sarà k il valore numerico di $\frac{\varphi(x)}{x}$, per qualsiasi valore di x, e quindi la forma della funzione è
$\varphi(x)=kx$ (d)
ove k denota una costante indipendente da x, ossia un parametro.
Da questa forma (d) di funzione, che esprime la proporzionalità fra i valori assunti dalla funzione e i valori attribuiti alla variabile, si ricava la relazione
$\varphi(x+w)=\varphi(x)+\varphi(w)$ (e)
che sussiste per x e w qualunque [Infatti: da φ(x)=kx e φ(w)=kw risulta φ(x)+φ(w)=kx+kw=k(x+w); ma da (d) si ha φ(x+w)=k(x+w), dunque ecc. c.d.d.; oppure, altrimenti: dalla (c), facendo la composizione di rapporti, si ricava $\frac{\varphi(x)+\varphi(w)}{x+w}$=k; e da (d) si ha $\frac{\varphi(x+w)}{x+w}$=k; onde ecc. c.d.d.].

[123] Non è stato possibile individuare né l'autore né l'opuscolo.

Viceversa se si domanda : "di che specie è una funzione φ la quale sia definita dalla condizione (e)?" - Si trova (come egregiamente fa Bordoni, p. 8, 9, 10 opuscolo L.S.) che essa funzione è della forma (d).

Dunque da (d) si deduce la (e), e viceversa da (e) la (d); ossia di dette due relazioni ciascuna è conseguenza necessaria dell'altra.

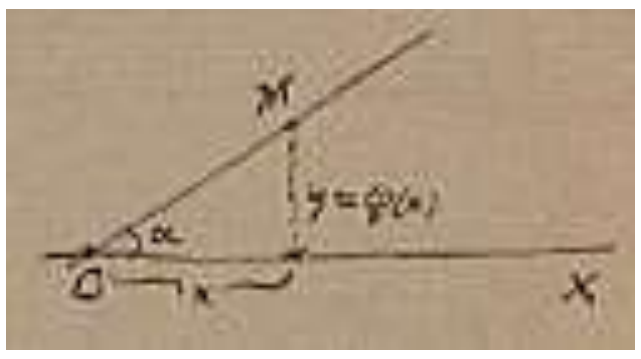

L'ordinata y=φ(x)=kx dal punto M descrivente una retta OM, che passa per l'origine O ed è inclinata all'asse OX di un angolo α tale che tangα=k, offre un esempio intuitivo della funzione considerata nell'esempio di questo n°.

- II. -

4. Siano ora f(x) una funzione qualunque di x ed f'(x) la sua derivata, cioè un'altra funzione di x, che si ricava (o che deriva) dalla prima, con regole di Calcolo determinate e conosciute e che qui non è il luogo di ripetere.

A sua volta il rapporto $\frac{f'(x)}{f(x)}$ è anch'esso evidentemente una funzione di x, il cui valore, per ogni valore assegnato di x, è dato dal quoziente dei corrispondenti valori di f'(x) e di f(x).

Per la f(x), per la f'(x), per $\frac{f'(x)}{f(x)}$ - come per ogni altra funzione imaginabile *[sic!]* di x - vale quanto al §I. si è detto per la φ(x). - Così, ad esempio, se la natura di f è tale che per x e w qualunque, coesistano le due relazioni

$$f(x)=f(x+w)=cost=c \text{ (poniamo)}$$

(f)

$$f'(x)=f'(x+w)=cost=z \text{ (poniamo)}$$

pel rapporto $\frac{f'(x)}{f(x)}$ sussisterà necessariamente (*) la relazione

$$\frac{f'(x)}{f(x)} = \frac{f'(x+w)}{f(x+w)} = \text{cost} = \text{a} \quad \left(\text{posto a} = \frac{z}{c}\right); \qquad \text{(g)}$$

ma le regole del Calcolo insegnano che, in questa ipotesi particolare, qualunque sia il valore del parametro costante c e nella 1ª (f), risulta z=o nella 2ª (f); e, per conseguenza, a=0 nella (g). Vale a dire che: "dato f(x)=cost=c, si deduce $\frac{f'(x)}{f(x)}$=cost=0" il qual risultato si accorda con la formola (2) pag. 8 dell'opuscolo di L.S. (che si trova trascritta qui sotto, al n.° 5).

Dunque: se φ(x)=cost. la (2) è vera per ogni valore di x.

======================

(*) Qui apro una parentesi per osservare che se dalla (f) si deduce necessariamente la (g), viceversa dalla (g) non si può in generale dedurre le (f): a persuadersene basti un esempio numerico: $\frac{3}{4}$ è = $\frac{75}{100}$; ma non per ciò 3 è =75, né 4=100. [Questa osservazione si metta a confronto con le ultime righe di pag. 14 (opuscolo di L.S)].

- III. -

5. Il sig. L.S. afferma (opuscolo, p. 8) che la formula

$$\frac{f'(x)}{f(x)} = \text{a} = \text{cost. risp. x} \qquad (2)$$

"non si riferisce punto, come pretende la Commissione del r. Istituto Lombardo di Scienze e lettere, ad un caso particolare, ma riguarda tutte le funzioni, derivate e integrali particolari, figuranti nella serie di Taylor: funzioni che possono essere potenziali, logaritmiche, esponenziali, trigonometriche."

Se questa affermazione non è erronea, la formola (2) sarà vera per qualsiasi forma di funzione; cosicché, se anche per una sola forma di funzione la (2) non fosse vera, l'affermazione sarebbe erronea.

Per la funzione f(x)=c=cost. risp. x si è veduto precedentemente (n.° 4) che la formola (2) è vera e che (qualunque sia il valore della costante c) risulta a=0; vediamo ora se la (2) sia vera per la funzione φ(x)=kx del n.°3.

Posto φ(x)=kx (k cost. indipendente da x; x valore qualunque) si ricava, come è noto, φ'(x)=k; e quindi, per ogni valore di x, si ha in questo caso:

$$\frac{\varphi'(x)}{\varphi(x)} = \frac{k}{kx} = a \text{ (poniamo)} \qquad (2')$$

dalla quale s rileva che qui a non è costante risp. x, poiché a è la reciproca del valore variabile x.

6. Se questa osservazione non persuade, si ragioni su un esempio numerico, così: siano x=5 un valore particolare di x e x=5+w (ove w è affatto arbitrario) un altro valore qualunque di x; ai quali corrispondono, pel rapporto $\frac{\varphi'(x)}{\varphi(x)}$, i seguenti valori: per x=5, $\frac{\varphi'(5)}{\varphi(5)} = \frac{1}{5}$; per x=5+w, $\frac{\varphi'(5+w)}{\varphi(5+w)} = \frac{1}{5+w}$. Domanda: sono uguali questi due valori? – Se w=0, sì, sono uguali; e allora si ha $\frac{\varphi'(5)}{\varphi(5)} = \frac{\varphi'(5)}{\varphi(5)}$; la formola (2) è valida anche per la funzione attuale; ma essa viene a dire il pane è pane, perché le due espressioni 5 e 5+w rappresentano (per w=0) un unico valore di k, cioè il 5. Ma se $w \gtreqless 0$, i due valori $\frac{1}{5}$ e $\frac{1}{5+w}$ non sono, o almeno non sembrano uguali; però, se la formola (2) è vera in tutti i casi, e quindi anche nel caso φ(x)=kx, è giocoforza concludere che qualunque sia w (diverso da zero) $\frac{1}{5}$ è uguale a $\frac{1}{5+w}$; or siccome due frazioni come queste, aventi uguali numeratori e denominatori differenti, si sono sempre ritenute differenti, la pretesa validità della (2) per la funzione φ(x)=kx, viene a dire il pane è vino. E questo si può ben dire, ma a condizione di cambiare il significato delle parole; come ben si può dire (e il sig. L.S. dice realmente) che $\frac{\varphi'(x)}{\varphi(x)}$=a è una cost. risp. ad x in tutti i casi, e cioè anche quando, come nella (2'), a non è indipendente da x, ma anzi è una funzione esplicita ($a = \frac{1}{x}$) della x variabile.

Se x assume i valori 1, 2, 3, ... 5, ..., la quantità a assume i valori $1, \frac{1}{2}, \frac{1}{3}, \dots, \frac{1}{5}, \dots$, cioè varia insieme con x; e poiché a sta a rappresentare il valore del rapporto $\frac{\varphi'(x)}{\varphi(x)}$, corrispondente a un valor dato di x, è questo rapporto che successivamente assume i valori $1, \frac{1}{2}, \frac{1}{3}, \dots, \frac{1}{5}, \dots$, ossia è questo rapporto che varia al variare di x: e se $\frac{\varphi'(x)}{\varphi(x)}$ varia al variare di x, come può esso in pari tempo essere cost. risp. x, soddisfacendo alla (2)? Non altrimenti che mutando il significato fin qui accettato delle parole. Senza una convenzione che chiarisca e spieghi una tal mutazione di significato, bisogna proprio persuadersi che la formola (2), vera per la funzione particolare φ(x)=c=cost (n.°4), non è vera per la funzione particolare φ(x)=kx e perciò non è vera per una funzione qualunque.

6'. La formola (2), non vera per φ(x)=kx, non sussiste nemmeno per φ(x)=logx, né per φ(x)=senx: perché il valor a del rapporto $\frac{\varphi'(x)}{\varphi(x)}$ in questi casi è dato rispettivamente da $a = \frac{1}{x\log x}$, e $a = \frac{1}{\tan x}$ [cfr. opuscolo L.S. formole (8) e (14)] e quindi, a variando con x, non è a cost. risp. x; invece la formola (2) sussiste ed è vera, qualunque sia x, per φ(x)=A*, perché il corrispondente valore a di $\frac{\varphi'(x)}{\varphi(x)}$ è dato [cfr. la formola (11) opuscolo L.S.] da a=logA, che è veramente una cost. risp. x.

- IV. -

7. Da quanto precede si raccoglie:
1.° che vi sono funzioni per le quali la formola (2) è valida; per es la φ(x)=C=cost e la φ(x)=A* (o più in generale la φ(x)=CA*), nelle quali C e A sono costanti risp. x (e quindi A potrebbe avere il valor particolare e, e C il valor particolare 1).
2.° che vi sono funzioni per le quali la formola (2) non è valida; per es. la φ(x)=kx (o più in generale la $\varphi(x)=kx^m$); la φ(x)=logx; la φ(x)=senx o in generale le funzioni trigonometriche.

Ora domando: fra tutte le altre possibili e imaginabili funzioni di una variabile non ve n'ha qualcuna, oltre quelle indicate nel comma 1°, per la quale sia valida la formola (2)? Rispondo: no non ve n'è alcun'altra. – Le sole funzioni, fra tutte le possibili e imaginabili, per le quali la formola (2) è valida sono: la φ(x)=cost, e la funzione esponenziale φ(x)=CA* (che include i casi particolari C=1, A=e=2,71…): cioè soltanto quelle indicate nel 1°comma.
8. A persuadersene basta porre il quesito: "quale è la forma di una funzione φ(x) per la quale si verifica la condizione (2)?" o in altri termini: "assegnare quali funzioni sono definite dalla relazione

$$\frac{\varphi'(x)}{\varphi(x)} = a = \text{cost. risp. x."} \qquad (2)$$

Ora se per ogni valore di x sta questa relazione, per due valori di x differenti di una quantità arbitraria w si avrà:

$$\frac{\varphi'(x)}{\varphi(x)} = \frac{\varphi'(x+w)}{\varphi(x+w)} = a; \qquad (3)$$

e viceversa se la (3) sta, per x e w qualunque, da (3) si deduce la (2), come giustamente dice il sig. L.S.
Basta dunque cercare quali funzioni φ sono definite da (3), come già si è determinato al n.°2 "quali funzioni φ sono definite da φ(x)=(x+w)" e al n°3 "quali funzioni φ sono definite da φ(x)+φ(w)=φ(x+w)."
Ora, dalla (3) per w e x, qualunque, si ha:

$$a=\frac{\varphi'(x)+\varphi''(x)+\frac{w^2}{1\cdot2}\varphi'''(x)+\frac{w^3}{1\cdot2\cdot3}\varphi^{IV}(x)+\dots}{\varphi(x)+w\varphi'(x)+\frac{w^2}{1\cdot2}\varphi''(x)+\frac{w^3}{1\cdot2\cdot3}\varphi'''(x)+\dots}$$

dalla quale si ottiene

$$a\left[\varphi(x)+w\varphi'(x)+\frac{w^2}{1\cdot2}\varphi''(x)+\frac{w^3}{1\cdot2\cdot3}\varphi'''(x)+\dots\right]=\varphi'(x)+w\varphi''(x)+\frac{w^2}{1\cdot2}\varphi'''(x)+\frac{w^3}{1\cdot2\cdot3}\varphi^{IV}(x)+\dots$$

Questa identità deve sussistere per qualsiasi valore di w; essa si risolve dunque necessariamente nelle seguenti relazioni di condizione, che si ottengono uguagliando fra loro i coefficienti delle uguali potenze di w e sono:

$$\varphi'(x)=a\varphi(x),\ \varphi''(x)=a\varphi'(x),\ \varphi'''(x)=a\varphi''(x),\ \varphi^{IV}(x)=a\varphi'''(x),\ \ldots \qquad (h)$$

Da queste relazioni si rileva

1° che se $\frac{\varphi'(x)}{\varphi(x)}=a=\text{cost}$, sarà anche $\frac{\varphi^{(i)}(x)}{\varphi^{(i-1)}(x)}=a$; ossia il rapporto fra una derivata qualunque (d'ordine i) e la derivata precedente (d'ordine i−1) è costante ed =a.

2° che eliminando in ciascuna delle relazioni (h) la funzione di x che essa ha in comune con la precedente, si ricava

$$\varphi'(x)=a\varphi(x);\ \varphi''(x)=a^2\varphi(x);\ \varphi'''(x)=a^3\varphi(x);\ \varphi^{IV}(x)=a^4\varphi(x),\ \ldots . \qquad (h')$$

D'altra parte pel Teorema di Taylor si ha

$$\varphi'(x+w)=\varphi(x)+w\varphi'(x)+\frac{w^2}{1\cdot 2}\varphi''(x)+\frac{w^3}{1\cdot 2\cdot 3}\varphi'''(x)+\ldots;$$

ossia, sostituendo a φ'(x), φ''(x), ecc. i valori trovati (h'), si ha

$$\varphi(x+w)=\varphi(x)\left[1+(aw)+\frac{(aw)^2}{1\cdot 2}+\frac{(aw)^3}{1\cdot 2\cdot 3}+\frac{(aw)^4}{1\cdot 2\cdot 3\cdot 4}+\ldots\right]$$

dalla quale (osservando che il coefficiente di φ(x), ossia la somma chiusa in parentesi, è notoriamente uguale ad $e^{aw}$) si ricava

$$\varphi(x+w)=\varphi(x)\ e^{aw}; \qquad (4)$$

relazione fra due valori, corrispondenti ad x ed x+w, di una funzione φ che soddisfa alla imposta condizione (2).

Per x=0, la (4) diviene $\varphi(x)=\varphi(0)e^{aw}$; questa formola non contiene più la x ed è valida per ogni valore di w; d'altronde φ(0) non contiene né x, né w, cioè è un numero che si può indicare con C; scrivendo dunque x invece di w, come di consueto, e ponendo φ(0)=C= cost si ottiene finalmente la formola:

$$\varphi(x)=Ce^{ak} \qquad (A)$$

in cui C ed a sono 2 costanti rispetto ad x (o due parametri) che è l'espressione o la forma generale della funzione cercata, della funzione cioè che risponde al quesito posto al principio di questo n.°8. In altre parole: la formola (2) è valida per le funzioni incluse nel tipo (A); non è valida per altre funzioni φ. Esaminando la (A) si riconosce che possono darsi due casi: cioè o il parametro a è nullo o è diverso da zero

1) Se a=0, per ogni valor finito di x è $e^{ax}=e^0=1$, e la (A) diviene

$$\varphi(x)=C=\text{cost. risp. } x \qquad (A_0)$$

2) Se a non è zero, si ponga $e^a=A$ ossia a=logA, e la (A) diviene

$$\left.\begin{array}{c}\varphi(x)=C\cdot A^x,\\ \text{che comprende la}\\ \varphi(x)=C\cdot e^x\end{array}\right\} \qquad (A')$$

pel caso particolare di A=e; e, come altro caso particolare (C=1), comprende le funzioni

$$\varphi(x)=A^x \text{ e } \varphi(x)=e^x \qquad (A'')$$

Dunque in conclusione, come dicevo, le sole funzioni per le quali la formula (2) è valida sono quelle della forma (A); ossia sono le ($A_0$), (A'), (A'') già indicate al n.°7 comma 1°.

- V. -

9. Ma innumerevoli funzioni si possono si possono immaginare, ed esistono, le quali non rientrano nel tipo (A), e sono perciò diverse dalle ($A_0$), (A'), (A''); ora, se soltanto per queste poche, contenute nel tipo (A), e per nessun'altra di quelle innumerevoli funzioni è valida la formola:

$$\frac{\varphi'(x)}{\varphi(x)} = a = \text{cost. risp. } x, \qquad (2)$$

non c'è a sorprendersi del fatto che "la Commissione (di cui del resto io non ho fatto parte e non conosco lo scritto), abbia potuto scrivere che il brano del Bordoni a pag. 46, da me (L.S.) citato, riguardi un capo particolare" (cfr. opuscolo L.S. p.14), sorprendente cosa avrebbe fatto la Commissione, se avesse scritto il contrario – se avesse cioè attribuito carattere di generalità, anzi di universalità, a una proprietà [espressa dalla (2)] che è particolare, appunto perché spetta soltanto ad alcune funzioni, a quella contenute nella forma (A) ossia alle sole funzioni ($A_0$), (A'), (A'').

Piola, Bordoni, Brioschi, nomi cari e due volte venerati, perché come matematici hanno contribuito al progresso della Scienza e come italiani hanno aggiunto onore e lustro a questa onorata e povera patria nostra, Brioschi, Bordoni, Piola, ne sono sicurissimo, avrebbero consentito, in quel che son venuto dicendo; certo con maggiore limpidità, certo con maggiore autorità, ma senza dubbio considerazioni non dissimili dalle mie avrebbero Essi espresso in argomento. Comunque sia, la Verità, per fortuna, dopo Galileo più non si piega all'"ipse dixit"; e pare a ogni modo assai ardito il pensare che – posto ad es. $\varphi(x)=\log x$, e quindi $\varphi'(x)=1/x$ – Quegli illustri Maestri, "per richiamare in onore lo studio delle basi filosofiche del Calcolo Superiore (opuscolo L.S. p.4)", fossero stati disposti a riguardare il rapporto

$$\frac{\varphi'(x)}{\varphi(x)} = \frac{1}{x\log x}$$

come una quantità costante risp. ad x, in conformità alla presunta generalità della formola (2).

Da Lagrange a oggi progressi enormi ha fatto l'Analisi Matematica, e da Euclide e Apollonio ad oggi ha fatto enormi progressi la Geometria (appena paragonabili ai progressi colossali delle scienze fisiche e chimiche dopo Volta e dopo Lavoisier fino ai giorni nostri); ma né oggi, né mai, per fermo, si troverà un Matematico serio che disconosca l'opera dei due Geometri giganteschi e del grandissimo Analista.

Per parte mia, nelle modeste considerazioni che sopra ho esposto, credo di non essermi allontanato di una linea dalla dottrina magistralmente insegnata da Lagrange.

- VI. -

Nell'"opuscolo L.S." sono riportate le parole del Vangelo "non veni legem solvere, sed adimplere", sta bene; e ivi si commentano con le parole : "la legge del progresso sta nell'aggiungere a quanto si è trovato e non già nel sopprimere basi indispensabili", d'accordo - ma aggiungere errori non è progresso e sopprimere errori non è regresso. A ogni modo, combattere lealmente un errore non è mancar di rispetto o di stima verso chi ha errato, ma è servire e onorare la Verità - un'alta idealità che non è da tutti, perché troppo spesso costa amarezze e richiede coraggio.

Il sig. L.S. ne dà un esempio egli stesso nella nota a pagina 7 del suo opuscolo; là ove egli per amor di verità, e pur trattandosi di un sì grande matematico, giustamente dichiara che Lagrange si è illuso nel ritenere di avere raggiunto nelle sue opere uno degli obiettivi che le aveva inspirate: quello cioè di stabilire il Calcolo Differenziale senza intervento di infinitamente piccoli.

Ora, se nel campo della Matematica, un matematico come Lagrange è potuto cadere in tale illusione, senza restarne punto diminuito davanti alla Scienza e alla Storia; non potrà un pensatore insigne, ma alle matematiche profano, essersi a sua volta illuso in una questione matematica? E che perciò? ne resterà forse diminuita la sua fama di illustre statista, di patriota intemerato, di filosofo coltissimo, di benemerito cittadino?

A ogni modo, fra un giovine che, sciente o incosciente, ha scritto su pei giornali articoli laudativi di un errore e destinati a mantenere il sig. L.S. in una fallace illusione - e un uomo maturo che, addossandosi una parte ingrata, ha contribuito a impedire che nei Rendiconti di un'Accademia scientifica o nell'organo degl'Ingegneri e dei Tecnici quella illusione fosse ufficialmente confermata e constatata, - fra quel giovine compiacente, dico, e quest'uomo rigido e severo, chi ha dato maggior prova di rispetto e di venerazione pel sig. L.S.?

E un'altra prova di venerazione e di rispetto intendo avergli dato nel vergar queste pagine. Avrei potuto non curarmi della Confutazione del sig. L.S. ad uno scritto che non mi riguarda e che non conosco; avrei potuto, pur volendomene curare, condensare in poche righe (*) le mie obiezioni a quella Confutazione; io, invece, pel desiderio vivissimo di convincere un Uomo illustre e venerando, che onora la mia città natale che altamente stimo, mi sono sforzato di rispondere in modo da rendermi possibilmente intelligibile anche ai non matematici di professione: soltanto a quel desiderio e a questo sforzo è da attribuirsi la prolissità del presente scritto (non destinato alla stampa), che in questa dichiarazione trova la sua spiegazione e la sua ragione di essere.

G. Jung

Milano 8 maggio 1901.

---

(*) Verbigrazia con le righe seguenti:

"Pongasi $\varphi(x)=x^m$, onde $\varphi'(x)=mx^{m-1}$ e $\frac{\varphi'(x)}{\varphi(x)}=\frac{m}{x}$; il rapporto $\frac{m}{x}$ varia con x, dunque il rapporto $\frac{\varphi'(x)}{\varphi(x)}$ (ad esso uguale) varia con x e perciò non è cost. risp. x: ciò vuol dire che per $\varphi(x)=x^m$ - ossia per almeno una funzione - la formola (2) non è vera, dunque la formola:

$$\frac{\varphi'(x)}{\varphi(x)}=a=\text{cost. risp. } x \qquad (2)$$

non è vera per tutte le funzioni, contrariamente a ciò che il sig. L.S. afferma; c.d.d."

**90**
Gian Antonio Maggi a **Giulio Krall**[124]

Milano 7 dicembre 1932

Caro Professor Krall,

Avendo assegnando come tesi orale di Laurea la Sua interessante Nota “Spiegazione energetica ecc.” del Bollettino di Bologna, trovai però conveniente di valermi, per la ricerca del minimo condizionato, del metodo dei moltiplicatori disponibili, e considerare così tutte e tre le equazioni scalari che traducono la costanza della quantità di moto areale, rispetto al comun centro di massa, delle quali non mi persuadevano le ragioni di metterne due in disparte. Ora, ultimamente ho dovuto inoltre riconoscere che il metodo, da Lei tenuto, della eliminazione della C, per ridursi al caso del minimo assoluto, conduce a questa conseguenza inammissibile, che, qualunque delle tre suddette equazioni si usi per tale eliminazione - né apparisce ragione per preferirne una, a tale scopo - i tre versori risultano paralleli all’asse a cui l’equazione impiegata si riferisce. Difatti, le cinque equazioni a pag. 7 dell’Estratto conservano la stessa forma, applicata agli angoli formati dai tre versori con uno qualunque dei tre assi. Credo opportuno di metterLa a parte di questa mia osservazione, astenendomi di farne oggetto di pubblicazione, perché Ella ne disponga come meglio crede.

Mi è grato intanto rinnovarLe i miei ringraziamenti per la buona memoria, che sempre mi dimostra serbare di me, e, coi migliori saluti, confermarmi

Aff[mo] Suo
Gian Antonio Maggi.

P.S. Il metodo dei moltiplicatori fornisce immediatamente l’equazione delle grandezze, nell’ipotesi dei versori a priori paralleli. Nel caso generale, conduce a stabilire dodici equazioni; ma, con spontanei accorgimenti, l’equazione delle grandezze e delle orientazioni si trova agevolmente e anche elegantemente.

[124] La corrispondenza con Giulio Krall è raccolta da Maggi in un fascicolo denominato *M.S. destinato alle Memorie della Società Astronomica Italiana, trattenuto, e mandato invece lettera al Krall, di cui unita la copia.- Si aggiungono il seguito della corrispondenza col Krall e appunti attinenti alla questione. Dicembre 1932.* In questo fascicolo si trova anche un *Teorema del 12 Dicembre*, con dimostrazione: *Valori delle variabili che rendono stazionaria una funzione, subordinatamente ad una o più equazioni di condizione, la rendono ancora stazionaria, subordinatamente all'aggiunta di altre equazioni di condizione, soddisfatte dai valori medesimi*; un *Corollario* al Teorema e l'*Applicazione al Problema del Krall*. Inoltre vi si trova uno scritto, probabilmente di Vivanti, datato da Maggi 18 aprile 1935, con una dimostrazione di un'affermazione sugli estremi condizionati (si veda la lettera #221).

**91**

Giulio Krall a Gian Antonio Maggi

[busta indirizzata: Illustre Signor Prof. G.A. Maggi - Milano - Corso Plebisciti n. 3; mittente: Ing. Krall - Roma - Via Gaeta n. 12]

Roma, 10-XII-1932-XI

Illustre Professore,

Ho ricevuto la Sua gentilissima con la cortese osservazione.

Nel ringraziarLa tengo però a rilevare che, mentre anch'io avevo inizialmente seguito la regola generale dei moltiplicatori, ò *[sic!]* poi trovato il criterio semplificativo che credevo dovesse risultare sufficientemente espresso a pag. 6 (dalla formula in poi) della nota del Bollettino, ma che è posto in piena evidenza nella I^ Nota Lincea, di cui Le trasmetto copia.

Esso criterio è ineccepibile. Ecco qual'è *[sic!]* il mio concetto:

Si tratta di minimizzare una certa funzione E condizionatamene con la relazione vettoriale

1) $\underline{n}c+\underline{k}p_{\psi}+\underline{k'}p_{\psi'}=\underline{K}$=vettore costante

Questa condizione equivale a 3 relazioni scalari, quelle ad esempio che si ottengono proiettando la 1) su $\underline{K}$ e su due altre direzioni normali $\underline{u_1}$ e $\underline{u_2}$

Si ottiene così:

1') $c\cos\beta+ p_{\psi}\cos\alpha+ p_{\psi'}\cos\alpha'=K$

1'') $c(\underline{n}x\underline{u_i})+ p_{\psi}(\underline{k}x\underline{u_i})+ p_{\psi'}(\underline{k}x\underline{u_i})=0$

Io procedo nel modo seguente:

a) Minimizzo la E tenendo conto della sola 1') con che risulta il parallelismo di $\underline{n}$, $\underline{k}$, $\underline{k'}$ a $\underline{K}$

b) Come conseguenza di quel parallelismo le 1'') risultano automaticamente soddisfatte.

Sicché è indifferente tenerne conto a priori (come si può fare col metodo dei moltiplicatori) o prescinderne in un primo tempo, salvo a constatarne, a calcoli fatti, il necessario sussistere.

Un procedimento analogo non riuscirebbe invece condizionando il minimo di E ad una delle 1''), ad es. la prima, e trascurando provvisoriamente le altre.

Risulterebbe infatti, lo riapprendo dalla Sua lettera,

$$\underline{n}=\underline{k}=\underline{k'}=\underline{u_i}$$

e la 1), per $K\neq 0$, come si suppone, non risulta soddisfatta.

(In una Nota in corso di stampa ò studiato il caso K=0. Si à, alla lunga, collisione.)

Mentre Le rinnovo i ringraziamenti per l'interesse ch'Ella accorda ai miei studi, mi lusingo che, quanto sopra Le sembrerà esauriente. – Troppo mi spiacerebbe che un cultore di Meccanica di così alta autorità quale Lei è, avesse ancora dubbi sulla coscienziosità del mio lavoro.

Le sarò, s'intende, infinitamente grato se, all'occasione Ella vorrà far sapere anche ai Colleghi, che fossero eventualmente a parte della sua riserva, che questa si è felicemente risolta. Ciò mi sta particolarmente a cuore, per ragioni morali e, indirettamente, anche materiali.

Credo d'averLa intrattenuta anche troppo, perciò concludo coi più distinti ossequi e saluti.

Suo dev$^{\underline{mo}}$ G. Krall

**92**
Gian Antonio Maggi a Giulio Krall[125]

Milano 14 dicembre 1932

Caro Professor Krall,

Sono ben lieto di affermarLe che, colle considerazioni della Sua lettera e della Nota lincea, di cui però non si trova cenno nella Sua Nota delle Memorie della Società Astronomica Italiana, ed altre prendono il posto nella Sua Nota del Bollettino, che, come Le scrivevo, non mi rendevano persuaso, col concorso anche di circostanze favorevoli, proprie del presente caso, riconosco la piena esattezza della Sua deduzione. Poiché ita il fatto, del quale però io ho sentito il bisogno di cercare la dimostrazione, che valori che rendono stazionaria una funzione, con parte delle imposte equazioni di condizione, e soddisfando le rimanenti, la rendono stazionaria con tutte quante. E la circostanza favorevole, a cui accenno, è che l'equazione relativa al modulo, scelta fra le tre, Le fornisce α=α'=β=0, coi quali valori sono soddisfatte le rimanenti due.

Non mancherò certo d'informare della Sua risposta un pajo di colleghi, che dovettero venire a parte della mia osservazione. Non occorre poi assicurarLa che, con così buon nome, non potevo toccarLe, in nessun caso, scapito di qualsiasi specie.

Le rinnovo con questo i più cordiali saluti, coi quali mi confermo

Aff.mo Suo
Gian Antonio Maggi.

**93**
Gian Antonio Maggi a **Charles-Ange Laisant**[126]

1496 (1899, 5) (M. Frolov) - Objection au §30 des "Études geometriques sur la théorie des parallèles" de Lobatschevsky (traduction de Hoüel). Les perpendiculaires DE, FG élevées des milieux des cotés AC, CB du triangle ABC étant supposées parallèles, la perpendiculaire HK élevée du milieu du troisième coté AB ne peut avoir avec aucune des deux un point d'intersection (Lobatschevsky entend toujours par là un point au fini), car autrement elles devraient se rencontrer aussi en ce même point (§29): ce qui est contraire à l'hypothèse qu'elles soient parallèles. Ce fait établi, deux cas pour chacune des deux droites DE, FG. Se présentent comme parallèles: la perpendiculaire HK est une parallèle à DE et (ou à FG), ou elle est un simple non sécante de cette droite (ou de l'autre). Le raisonnement de Lobatschevsky, exempt de tant pétition de principii, écarte le deuxième cas, et démontre que c'est le premier qui a lieu.

Pour plus de détails et un procédé quelque peu différent, on pourra voir le §111 des "Nouveaux principes de géométrie avec une théorie complète des parallèles" (Voir le Recueil complet des ouvrages de Lobatschevsky, publié par l'Université Impériale de Kasan, ouvrages écrits en langue russe, 1884).

D'ailleurs on peut reconnaître à coup d'œil l'exactitude de la proposition de Lobatschevsky à l'aide de la représentation conforme comme du plan lobatschevskyen

[125] In due copie di minuta.
[126] Redattore del periodico *Intermédiaire des mathématiciens*. Questa minuta si trovava nella cartella *Giudizii*.

(surface pseudosphérique), dans un cercle du plan euclidien, par laquelle les points à l'infini ont pour image le cercle limite, les droites (géodésiques) des cercles normaux au cercle limite, les cercles (y constant toute ligne à courbure géodésique constante) des cercles: et réciproquement. Suivant que le cercle circonscrit à l'image du triangle ABC se trouvera tout à l'intérieur du cercle limite, ou lui sera tangente, ou sécant, le cercle (lobatschevskyen) circonscrit au triangle objectif aura son centre effectif, au fini ou à l'infini, ($1^r$ et $2^m$ cas) ou idéal ($3^m$ cas): c'est à dire que les perpendiculaires élevées du milieu des cotés du triangle ABC se rencontreront en un point, au fini, ou à l'infini ($1^r$ et $2^m$ cas) ou ne se rencontreront pas. Cette belle remarque est de mon collègue M. Bianchi.

G.A. Maggi (Pise)

Mandata a M. C.A. Laisant
Rédacteur de l'Interm. des Mathém. Avenue Victor Hugo 169, Paris.
(il 30 Giugno 1899).

**94**

**Giulio Lazzeri** a Gian Antonio Maggi
[carta intestata: Periodico di matematica - Livorno - Via del Porticciolo 2]

Livorno, li 15 nov. 1900

Egregio Professore,

Non son riuscito a trovare l'opuscolo del Sig.$^r$ Bertrand di Ginevra di cui mi parlò; ma pensando alla strana formula da lei indicata, mi pare che si possa trovare con questo non meno strano ragionamento. Sia E la somma degli angoli esterni, I la somma degli angoli interni e T l'area di un triangolo. Poiché E e T ricoprono tutto il piano, si ha

$$E+T=4R;$$

inoltre $$6R=I+E$$

Sommando $$2R+T=I,$$

e quindi $$\frac{I}{2R}=1+\frac{T}{2R}$$

Questa, se non erro, è la formula del Bertrand.

Il ragionamento mi pare molto arrischiato, ma in ogni modo ho voluto comunicarglielo
Cordiali saluti dal

Suo Aff$^{mo}$
G. Lazzeri

**95**
**Tullio Levi-Civita** a Gian Antonio Maggi
[busta e carta intestate: Hôtel Riposo Martignoni
Prima cappella del Sacro Monte sopra Varese;
busta indirizzata a: Chiar.[mo] Sig. Prof. Gian Antonio Maggi
via Zanardelli, 57 - Viareggio]

3-VIII-16

Carissimo amico,

Per un complesso di circostanze di famiglia (derivanti principalmente dal fatto che mio padre, senza essere ora malato, ebbe lo scorso inverno qualche disturbo cardiaco) debbo quest'anno tornare presto a Padova. Vi sarò fra pochi giorni, lasciando mia moglie nel Veronese a compiere la villeggiatura, non molto fresca in verità, ma sempre preferibile allo scirocco patavino. Per questo motivo, mentre genericamente tanto io quanto mia moglie rammarichiamo molto di non incontrarti, non possiamo rimpiangere la tua determinazione di non capitare per ora in Lombardia.

Chi sa che, quando ci vedremo a Roma a fin d'anno, sia finalmente spuntata l'alba della pace. Così a impressione, lo riterrei probabile, ma sta di fatto che Cadorna la pensa in modo alquanto diverso. Ho infatti saputo dal collega Iorini che due settimane or sono Cadorna fu da lui udito affermare (durante la table d'hôtes all'albergo Excelsior di Varese) che la guerra potrà finire soltanto entro il 1917. Ma forse - e questo è il mio perfezionamento ottimistico - Cadorna considera soltanto il lato militare.

Apprendo con piacere che stai curando un nuovo trattato da far seguire alla Dinamica fisica. Anch'io <u>dovrei</u> occuparmi di rimaneggiare le mie litografie per la pubblicazione, ma il lavoro mi è così ingrato che sempre cerco dilazioni.

Circa la questione specifica su cui richiami la mia attenzione, premetto che, nei libri consultati [nel preparare alcune generalità sulla dinamica dei mezzi continui, per il corso di mecc. sup. dell'anno passato] non ho trovata quella dimostrazione comprensiva del teorema di univocità, a cui tu aspiri.

Il tuo postulato [$\frac{dk}{dt} \geq 0$ per $p>0$] è certo attendibile nelle condizioni supposte (forze di massa e sforzi superficiali nulli), ma non mi sembra facile darne una giustificazione fisica soddisfacente. Bisognerebbe in qualche modo far intervenire la irreversibilità di certi processi naturali; altrimenti sarebbe lecito cambiare t in $-t$, e allora il postulato richiederebbe $\frac{dk}{dt}=0$, il che non è vero pei gas (a meno di non ammettere <u>a priori</u> che si trovino in quiete, ciò che si deve evitare).

Ricordo che, nelle accennate lezioni, ho richiamato in generale l'equazione (globale) delle forze vive sotto la forma

(1) $$dT=dL-\alpha dt,$$

(dL lavoro complessivo di tutte le forze esterne)

riattaccandovi considerazioni energetiche sia pei sistemi conservativi (in particolare fluidi perfetti e solidi elastici) che pei sistemi dissipativi (in particolare fluidi viscosi).

Pei fluidi perfetti

$$\alpha=\int_\tau \frac{dk}{dt}\frac{p}{k}d\tau;$$

e, se si vuole che il fenomeno abbia carattere conservativo per qualunque porzione $\tau$ del sistema e per qualunque tipo di movimento, si è indotti a introdurre l'ipotesi matematica "**k** funzione della sola **p**" con che $p\frac{dk}{k}$ si può risguardare come differenziale di una funzione **e** della sola p, l'energia interna per unità di volume. L'analoga energia interna di tutto il sistema vale in conformità

$$E=\int_\tau e d\tau,$$

e la (1) si scrive espressivamente

$$dT=dL - dE.$$

Per dL=0, scende

$$T+E=\text{cost},$$

e di qua si può ricavare il teorema di univocità con ipotesi addizionali su E. Per es., nel caso di gas perfetti in regime isotermico, p=ck (con c costante). A meno di una inessenziale costante additiva, e=ck, E=cM (M massa totale del gas). In questo caso il teorema di univocità apparisce semplice conseguenza (della specifica equazione di stato p=cke) dell'invariabilità della massa. Ma si tratta sempre, come nell'esempio da te ricordato del moto irrotazionale, di casi particolari, in cui si sfruttano ipotesi addizionali. Resta da vedere se ciò è veramente imposto dalla natura delle cose, o se il tuo desiderato è suscettibile di una risposta spontanea e perspicua. Può essere che la difficoltà che ti ho segnalato provenga dal fatto che a me sfugge un'intuizione fisica che tu invece hai netta.

Anch'io ho scorso nel Nuovo Cimento la ricerca di Somigliana, ma non l'ho approfondita, pur trovandola notevole per l'impostazione generale del problema e per la grande semplicità e chiarezza della trattazione. Non sono però in grado di esprimere alcun parere circa il dubbio da te affacciato che possano non essere distorsioni di Volterra le deformazioni ai cui spostamenti si impongono discontinuità rappresentate da spostamenti rigidi.

Jung ricambia cordialmente i saluti; interpreto il pensiero di Cisotti (che vedrò domani) ringraziando e ricambiando anche da parte sua.

Mia moglie, che sta, se non bene, certo assai meglio dello scorso anno, desidera esserti ricordata, ricordando a sua volta simpaticamente la tua gentilezza cordialissima.

Io infine ti prego di presentare distinti ossequi in famiglia e ti invio una buona stretta di mano.

Tuo aff.[mo]

T. Levi-Civita

**96**
Tullio Levi-Civita a Gian Antonio Maggi
[busta indirizzata a: Chiar.$^{mo}$ Sig. Prof. Gian Antonio Maggi
via Zanardelli, 54 - Viareggio;
busta e carta intestate: Prof. T. Levi-Civita - Padova]

12-VIII-16

Carissimo amico,

Ti ringrazio della tua lettera cordiale e delle considerazioni interessanti sulle distorsioni elastiche. Questa mi hanno invogliato a rivedere la tua nota del 908[127] [in cui ho segnata la correzione materiale accennatami], la memoria di Volterra e gli articoli recenti di Somigliana.

Mi pare di aver visto chiaro e, per mio conto, sono rimasto soddisfatto di tutti tre.

Se non erro, le formule (5) di Somigliana bastano ad assicurare la continuità delle componenti di deformazione attraverso il diaframma, nella ipotesi che si assegnino agli spostamenti discontinuità rappresentate da spostamenti rigidi. Viceversa, se si ammette a priori che rimangano continue in un campo le caratteristiche di deformazione e loro derivate, risulta, come è mostrato con grande perspicuità nella tua nota, che gli spostamenti non possono essere affetti da discontinuità in un campo semplicemente connesso, mentre ecc.

Rimane, come tu giustissimamente rilevi, la questione di esistenza, trattandosi, nelle ricerche di Volterra, di dati a priori esuberanti. Ma la difficoltà è superata dalla costruzione effettiva delle 6 distorsioni elementari, con che il problema è riportato all'ordinaria teoria dell'elasticità.

Rispetto al $p\frac{dk}{dt}\geq 0$ [su cui - siamo d'accordo - non è il caso di faticare], io volevo dire semplicemente che, trattandosi di fenomeni reversibili, la stessa intuizione che vale per $dt>0$, dovrebbe valere per $dt<0$. E non vedevo, né vedo, quale potesse essere tale intuizione reversibile.

De Marchi,[128] che non si è mosso da Padova, mi incarica di trasmetterti saluti affettuosi e auguri pel figliolo.[129] Ben volentieri mi vi associo una volta ancora. Andrò presto a Casarino per un paio di giorni; intanto ti ringrazio sia a nome di mia moglie che di mio padre per l'interessamento cortese.

Con preghiera di ossequi in famiglia e di ricambio cordiale a Veronese, ti stringo la mano

Tuo aff.$^{mo}$
T. Levi-Civita

P.S. L'indirizzo espressamente comunicatomi da Signorini come preferibile è [o meglio era circa un mese fa, quando egli mi scrisse] Sottotenente A. S. Bosco Chiesanuova (Verona) Naturalmente mi sono attenuto anch'io a questo indirizzo per ringraziarlo della sua memoria. Aggiungo per notizia che Signorini appartiene al IX Regg. Art. da Fortezza.

---

[127] G.A. Maggi, "Sugli spostamenti elastici discontinui", *Rendiconti della R. Accademia dei Lincei*, s. V, v. XVII, pp. 571-576. Ristampata, con un'aggiunta, in *Selecta,* pp. 257-264.
[128] Potrebbe trattarsi di Luigi, geofisico, docente presso l'Università di Padova.
[129] Felice Maggi (Messina 1893-Desenzano 1929).

**97**
Tullio Levi-Civita a Gian Antonio Maggi
[foglietto, senza data, ma febbraio 1923 per la risposta di Maggi, #98]

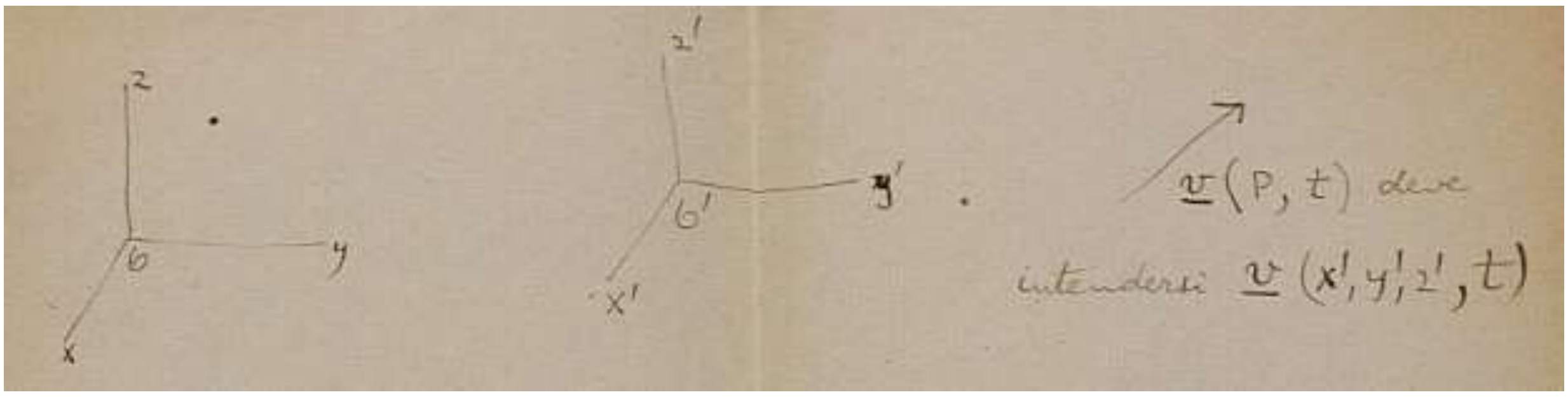

Ora, essendo gli assi Ox’y’z’ in moto traslatorio (anche non uniforme) rispetto ad Oxyz, la derivata di un vettore variabile è la stessa, sia che si faccia con referenza all’uno, oppure all’altro dei due triedri. Insomma $\frac{dv(P,t)}{dt}$ va concepito come la derivata di $\underline{v}$(x’,y’,z’,t) fatta con referenza ad Ox’y’z’, nell’ipotesi che x’,y’,z’ varino con legge assegnata. Se in particolare si suppone che P designi una particella, animata (rispetto ad Oxyz) dalla velocità $\underline{v}$ e quindi, rispetto ad Ox’y’z’, dalla velocità $\underline{v}-\underline{w}$, si ha la formula della sig.$^{na}$ Guglielmi[130]

$$\frac{d\underline{v}}{dt}=\frac{\delta\underline{v}}{\delta t}+\frac{d\underline{v}}{dP}(\underline{v}-\underline{w}).$$

**98**
Gian Antonio Maggi a Tullio Levi-Civita

Pisa 10 febbraio 1923

...
Ho veduto il tuo foglietto idrodinamico, che conta tra gli oggetti dei miei ringraziamenti. Certo, in $\underline{v}$(P,t) si può intendere <u>a priori</u> P riferito agli assi mobili; ma, con questo, $\frac{\delta\underline{v}}{\delta t}$ deve annullarsi, nell’ipotesi che il moto del fluido abbia carattere stazionario rispetto all’osservatore mobile, mentre, in principio a pag. (25), l’autrice scrive, in questa ipotesi:

$$\frac{\delta\underline{v}}{\delta t}=\underline{\Omega}\Lambda\underline{v}.$$

Mi sembra quindi che, nel dedurre

$$\frac{d\underline{v}}{dt}=\frac{\delta\underline{v}}{\delta t}+\frac{d\underline{v}}{dP}\frac{dP}{dt}$$

da $\underline{v}=\underline{v}$(P,t) si deve intendere P rispetto agli assi fissi.
Mi dirai... se prendo un abbaglio.
...
Al Prof. T. Levi-Civita - Roma.
Sulla Nota “Sul moto vorticoso dei liquidi” della Sig.$^{na}$ A. Guglielmi (Atti del R. Ist. Veneto 1921-1922).[131]

---
[130] Si veda la lettera successiva.
[131] Tale Nota non è pubblicata sul volume citato, dove appare solo una presentazione di Levi-Civita a p. 19.

**99**
Gian Antonio Maggi a Tullio Levi-Civita

Milano, 10 Novembre 1925.

Carissimo amico,

Grazie della sollecita risposta, colle indicazioni per la scoperta di Adams in “Nature” - da non confondersi colla scoperta di Nettuno - e grazie del ricambio dei saluti da Bellagio, che mi ha trasmesso il collega Cisotti. Sui moti rettilinei per gravitazione di un sistema continuo non ho in mente di aver mai veduto qualcosa, e aspetto di essere istruito dai tuoi risultati, di cui pregusto l’interesse.

Ti accludo il foglietto per Lorentz colla mia firma e quelle di Bianchi, Cisotti, Vivanti, Bruni, e Vaglia Cambiario della Banca d’Italia di centocinquanta lire, per raccogliere le quote, non ho potuto assolvere il mio compito prima d’oggi.[132] Inter nos, io c’entro per cento, Bruni per venti e gli altri per dieci. Te lo dico, a prova ch’io non ho mancato di procurare che la somma fosse discreta, tenendo il noto conto dei cambi, ahimè!, presenti, ma, a tal fine, avendo pure accennato al mio contributo, naturalmente, senza far insistenze, non ho ottenuto di più.

Ora permettimi di chiedere la tua autorevole opinione sopra un punto della Relatività - alle soglie della medesima - a proposito, particolarmente, della conferenza del Langevin, di cui ti parlerà l'amico Enriques, col quale, grazie a Mathesis, poi agli Intellettuali, ho avuto il piacere di trovarmi ripetutamente.[133]

Facciamo le ipotesi classiche dell'etere fisso occupante lo spazio S, e di uno spazio S' in moto assoluto, relativamente a S. Facciamo ancora l'ipotesi normale che la velocità della luce, emanante ad un istante, da un punto P, sia la stessa, non possieda o possieda il punto, una velocità. Sia P il punto di S quando emana la luce al supposto istante, e P' il punto, mobile con S', che porta rispetto a S' la posizione che P possiede in S. La velocità della luce emanante da P, in S, sarà la stessa in tutte le direzioni. Ma la velocità della luce relativa allo spazio S' si comporrà, in ogni direzione, colla regola del parallelogramma, colla eguale e contraria alla velocità dello spazio S' rispetto allo spazio S. Vale a dire P’ dello spazio S', considerato come fisso, si giudicherà emanare da P' la luce, con velocità di propagazione diversa nelle diverse direzioni. Questa differenza è ciò che doveva rilevare, pel movimento annuo del Globo Terrestre (S') nello spazio (S), l'interferometro di Michelson. E non essendo stato rilevato, il precedente ragionamento è risultato in difetto, e sono caduti i principii classici su cui è fondato. Non è così? Ora, secondo me, l'accennata composizione di velocità non ha a che vedere colla ipotesi balistica, la quale postula una composizione delle velocità assolute, della luce e della sorgente, nello spazio S. In conseguenza di che risulta eliminata la differenza, nella diverse direzioni, della velocità relativa ad S', considerato come fisso. Non vedo quindi come il Langevin ed altri (direi anche il Castelnuovo), allegano i fatti astronomici opposti alla ipotesi balistica a dimostrazione della eguaglianza in tutte le direzioni della velocità relativa, presentando questa eguaglianza come un fatto sperimentale, che non richiedeva la conferma del risultato negativo di Michelson. Ti sarò molto grato di una parola in proposito che mi rischiari ogni dubbio sulle mie interpretazioni di queste circostanze fondamentali.

---

[132] Si veda la lettera seguente.

[133] Probabilmente si riferisce all’associazione Fédération internationale des Unions intellectuelles (Kulturbund), che nel 1925 tenne il secondo congresso a Milano, alla quale aderivano, tra gli altri, Enriques e Langevin.

Il Cisotti e il Bianchi mi hanno informato della pubblicazione del Giorgi sulla esperienza del Miller,[134] che hanno ricevuto in questi giorni, e per mezzo loro ho potuto anche vederlo. Io non lo ho ricevuto. Non so se il Giorgi mi ha soppresso l'M., come l'Hadamard...

Chiudo con questo, presentandoti i migliori saluti di noi tutti, coi quali e colla preghiera dei miei più cordiali rispetti alla tua signora, ti abbraccio mi confermo

Tuo aff. amico
G.A. Maggi

**100**[135]
Tullio Levi-Civita a Gian Antonio Maggi
[busta indirizzata a: Chiar.mo Sig. Prof. Gian Antonio Maggi - Corso Plebisciti, 3 - Milano; busta e carta intestate: Prof. T. Levi-Civita - Via Sardegna, 50 - Roma (25)]

10-XI-1925

Carissimo amico,

Mi affretto ad accusare ricevuta della tua raccomandata con incluso vaglia di £ 150 e molto interessante accompagnatoria. Quanto ai biglietti e alla sottoscrizione lorentziana, dobbiamo trovarci domattina con Corbino per raccogliere le file e fare la spedizione in Olanda. Ti ringrazio per la tua risposta e per il tuo interessamento preventivo relativi alla mia notina.

Quanto alla questione relativista mi pare prima di tutto che tu la abbia impostata nel modo il più chiaro e preciso, sicché non potrei che parafrasare, sciupando.
Resta il punto essenziale, su cui tu richiami la mia attenzione, cioè per quale motivo si affermi che i fatti astronomici [ce ne sono di decisivi, specialmente spettroscopici, checché ne dica il La Rosa] che sono in contrasto coll'ipotesi balistica, renderebbero superflua l'esperienza di Michelson. Io proprio non lo so, e sono perfettamente del tuo avviso che, per abbandonare la teoria classica, bisogna accertare ad un tempo l'indipendenza della velocità della luce, sia da quella della sorgente che da quella dell'osservatore. Non so che cosa dica in proposito il Langevin; ma nel libro del Castelnuovo non trovo simile asserzione.

C'è qui (arrivato stamane da Torino) lo Stekloff il quale mi aveva anche portato i tuoi saluti.

Il Giorgi già da vari giorni mi chiese telefonicamente il tuo indirizzo per mandarti quelle tali pubblicazioni. Non capisco quindi come mai non ti siano ancora pervenute.[136] Ad ogni modo, vedendo il Giorgi, gli domanderò se ha spedito.

Ossequi in famiglia; saluti cordiali da mia moglie e un abbraccio affettuoso

Dal tuo T. Levi-Civita

[134] G. Giorgi, “Sulle esperienze di Miller”, *Atti della Pontificia Accademia delle scienze*, LXXVIII, 1925, pp. 170-181.
[135] Dopo questa si inserisce forse la #227, scritta da Maggi e mancante del destinatario, ma probabilmente diretta a Levi-Civita.
[136] Si veda la lettera precedente.

**101**

Tullio Levi-Civita a Gian Antonio Maggi
[busta indirizzata a: Chiar.mo Sig. Prof. Gian Antonio Maggi -
Corso Plebisciti, 3 - Milano (121);
busta e carta intestate: Prof. T. Levi-Civita - Via Sardegna, 50 - Roma (25)]

19-V-1930

Mio caro amico,

reduce da Padova, dove ho fatto una brevissima scappata, trovo la tua letterina del 16, concernente il "nostro" Cirenei.[137] Come ben supponevi, egli aveva scritto anche a me in termini analoghi, ed io gli risposi press'a poco così: "Ricordo perfettamente di aver riconosciuto, quando esaminai i Suoi apparecchi, che sono sensibili e costruiti con garbo. Ma non sono in grado di pronunciarmi circa il loro interesse didattico. Mi trovo perciò nella impossibilità di compiacere il desiderio espresso nella Sua lettera del... . Mi creda con distinti saluti...". Arrivando, ho trovato una ulteriore lettera, che mi affretto a comunicarti per notizia e consiglio. La mia idea sarebbe di prenderne atto e basta; e così farò, a meno che tu non trovi preferibile un diverso contegno, riscrivendomi in argomento.

Ancora non so se le sedute lincee avverranno, secondo la consuetudine, nell'ultima settimana del mese o saranno rinviate alla successiva. Comunque arrivederci presto e mille cordiali saluti

Tuo aff.
T. Levi-Civita

P.S. Quanto alla proposta per il s.c. *[socio corrispondente]*, essa era di spettanza dei soli soci della sezione mat., sicché a noi di mecc. la scheda non è stata inclusa di proposito.

**102**

Gian Antonio Maggi a Tullio Levi-Civita[138]
[busta indirizzata a: Ch.mo Sig. r Prof. Tullio Levi-Civita
della R. Università di Roma. (25). - Via Sardegna, 50;
mittente: Da G.A. Maggi - Corso Plebisciti, 3, Milano]

Milano 20 Maggio 1930

Carissimo Amico,

Grazie della sollecita tua risposta, e dell'accluso Cirenei, che ti rimando, accluso. Io ho ben presente di non essermi menomamente pronunciato col Cirenei, sui suoi strumenti ed esperimenti, e credo che quel mio mutismo sia confermato nella nostra Relazione. Ma avrei bisogno che me lo potessi confermare, per rispondere al Cirenei ch'egli s'inganna, e che io non ho mai né detto né scritto che i suoi strumenti siano idonei per l'insegnamento della Meccanica.

[137] Si confronti la lettera #85 scritta da Maggi a Guidi sull'argomento. Si veda anche la nota alla lettera #8.
[138] In due copie di minuta.

Le espressioni benevole da te dette a loro riguardo, che egli ora domanda di confermare, non sono che una delle molteplici dimostrazioni di quella tua bontà, che non trova paragone che nel tuo ingegno. Ed è tutto dire!

Se io m'inganno nel ricordare gli esperimenti del Cirenei, tu mi potrai correggere, ma io li ricordo così. Egli se la prendeva col principio della leva, da lui conciato nella forma (sotto cui l'ha più volte stampata), "ciò che si guadagna in forza, si perde in tempo"; e, sotto queste mentite spoglie, lo demoliva con esperimenti, che, di fatto, sarebbero stati diretti contro l'obbiezione al da mihi ubi consistam, terramque movebo, consistente nell'osservazione che, equilibrato il peso ideale della Terra con un peso accessibile all'uomo, e quindi con un braccio della potenza estremamente grande, si trova poi occorrente un tempo infinitamente grande, per girare la leva, pur di un angolo infinitamente piccolo. La ragione sta nell'enorme momento d'inerzia, rispetto all'asse. Ora, il Cirenei ci mostrava che la sua leva, dove la resistenza era rappresentata da un contrappeso, girava di un certo angolo, sempre in uno stesso tempo, comunque si accrescesse il braccio della potenza equilibrante, coll'applicarsi un peso più piccolo. E la ragione ne era chiaramente che, col forte momento d'inerzia, che la sua massiccia leva possedeva, per proprio conto, le piccole differenze di tempo, corrispondenti alla variazione del braccio e alla grandezza del pesetto variabile, diventavano insensibili alla abbastanza grossolana osservazione.

Così, questi esperimenti - come io li ricordo - che cominciavano col non aver nulla a che fare colla confutazione del principio della leva, sotto forma onesta, neppure concludevano qualcosa, per l'obbiezione al da mihi ubi consistam ecc., perché, eseguiti colle debite cure, dovevano, col fornire un tempo crescente, col diminuire del peso equilibrante, o il crescere del braccio, confermare l'obbiezione.

Non ho quindi ragione di dire che ci è voluta tutta la tua bontà, per dorare la pillola, che, così dorata, e da tanto doratore, egli si è ben guardato d'ingoiare, per farne bella mostra col colto pubblico?...

Finora io non ho risposto. Ho scritto invece in proposito all'Accademico Guidi,[139] che, tempo fa, mi domandò d'informarlo della nostra Relazione, avendo il Cirenei chiesto il giudizio dell'Accademia d'Italia!

Se, senza troppo disturbo, mi puoi confermare il mio mutismo... scritto, risponderò come ti dicevo: diversamente, il miglior partito mi sembrerebbe... di conservarlo intatto.

Di nuovo, le migliori cose a te e alla Signora, e, contando sempre sul piacere di rivederci presto, intanto un abbraccio a te dal

Tuo aff. amico
Gian Antonio Maggi.

[139] Si veda la lettera #85.

**103**
Tullio Levi-Civita a Gian Antonio Maggi[140]
[busta indirizzata a: Chiar.$^{mo}$ Sig. Prof. Gian Antonio Maggi
della R. Università di Milano - Valnegra (Bergamo)]

Padova 5-IX-1930

Carissimo amico,

Ricevo col più grande piacere la interessante tua lettera di ieri, che mi trova sbarcato qui da un paio di giorni appena. Prima di tutto esprimo, anche a nome di mia moglie, i più fervidi rallegramenti per l'imminente matrimonio della tua figliola.[141] Vorrai, te ne prego, presentarle le nostre felicitazioni coi più lieti auguri. Posso facilmente rendermi conto, essendo anch'io particolarmente amante della mobilità, della noia che ti deve aver causata la nevralgia ad una gamba, per di più in campagna; e affretto coll'augurio il momento in cui sarai libero anche dai postumi residui.

Anch'io mi sono portato qui, col proposito di studiarli un po', i più moderni libri sulla meccanica ondulatoria e quantistica (De Broglie, Frenkel, Dirac, quest'ultimo comperato a Londra "just out"). Già in viaggio mi sono addormentato parecchie volte sul Dirac, perché faticoso o almeno molto concettoso e richiedente pacata meditazione di tavolino. Non so se troverò qui l'energia necessaria. Ora sono un po' stanco dell'attività turistica e congressistica (essendo stato un mese in Inghilterra e poi a Stoccolma), e mi riposo. Il Frenkel (pur avendolo sfogliato solo in principio) mi pare che si legga bene. Sono poi arrivato circa a metà del De Broglie che presenta innegabili pregi di semplicità, chiarezza, precisione logica e riattacco, nei limiti del possibile, alle concezioni classiche. Però non dovunque è detta l'ultima parola. In particolare per es. l'associazione di una propagazione ondosa ad una famiglia di $\infty^{n-1}$ traiettorie, come la esponi tu, è, a mio gusto, molto preferibile per semplicità concettuale e formale. Anche nella subordinazione dell'ottica geometrica alla fisica non mi pare ancora raggiunta la forma definitiva. In attesa di andare più avanti, leggerò naturalmente col massimo interesse la tua recensione, quando uscirà nel bollettino del Pincherle.

Nei riguardi dell'attrito, già ero persuaso che la tua formulazione sia definitiva e che ben poco si possa aggiungervi dal punto di vista sintetico.

Speriamo che la polemica Gugino-Mineo cessi per esaurimento; sono comunque lieto che, per ogni eventualità, vi sia un giudice non sospetto di parzialità verso uno dei contendenti, che si è presa la briga di approfondire la questione.

A Stoccolma eravamo 5 italiani:
il Signorini (rappr. del Comitato di ricerche)
il Burzio (dell'Acc. Militare di Torino, rappr. del Min.$^{o}$ della Guerra)
il Panetti, che ha tenuto una conferenza generale sulle oscillazioni dei veicoli,
Finzi-Contini, assistente di idraulica a Milano; e io. Al precedente congresso dei geofisici gli italiani erano 16, e li abbiamo trovati quasi tutti ancora a Stoccolma.

Il nostro congresso si è svolto in una atmosfera molto simpatica di cordialità. A giudicare degli argomenti, molte comunicazioni devono essere state interessantissime, ma se ne sono udite assai poche, prima di tutto perché le 4 sezioni del congresso (idro e aerodinamica, elasticità, stabilità e vibrazioni, mecc. raz. e balistica) funzionavano contemporaneamente, e poi perché era diffusa la sensazione che fosse più proficuo parlare alla spicciolata coi colleghi, anziché ascoltare letture, specie in lingue non egualmente famigliari. C'è anche stata sovrapposizione fra la comunicazione di Signorini "Sulle deformazioni termoelastiche finite" cui ho assistito e la prima parte della conferenza generale del Panetti. Andai alla seconda, ma, naturalmente, non potei capire bene, mancandomi le premesse generali. Tutto ciò si potrà per altro leggere tra breve, perché il volume dei Rend. del Congresso è già in bozze.

[140] Questa lettera si trovava nelle *Rusticationes*.
[141] Si veda la lettera #145.

Il Signorini è andato ulteriormente al Nord; si tratterrà poi a Berlino e a Vienna, sicché rientrerà soltanto in Ottobre.

Ossequi in famiglia, coi migliori ricordi da parte di mia moglie; una cordialissima stretta di mano

Tuo aff.
T. Levi-Civita

**104**

Tullio Levi-Civita a Gian Antonio Maggi[142]

[busta indirizzata a: Chiar.mo Sig. Prof. Gian Antonio Maggi - Valnegra (Bergamo)]

Padova 19-IX-1930

Carissimo amico,

Il tuo bel volume,[143] dono in ogni caso graditissimo, non poteva giungere in un momento per me più opportuno e propizio.

I colleghi di qui mi fanno premura per includere in un fascicolo del Seminario patavino, che trovasi presentemente in corso di stampa, quella piccola osservazione sul comportamento asintotico dei potenziali newtoniani, dovuti a distribuzioni che si estendono all'∞, di cui ti accennai quando ci trovammo a Roma nel Giugno scorso. Proprio in questi giorni debbo pensare alla redazione, e mi è veramente prezioso l'aver sott'occhio in sintesi le nozioni classiche della teoria del potenziale, presentate da un pensatore che vi lascia cogliere il frutto maturo delle sue proficue meditazioni e della sua insigne esperienza.

Ben si intende che mi interessa anche tutto il resto: sia la teoria fenomenologica del campo elettromagnetico che è il nucleo centrale del volume, sia i succosi complementi sulla relatività ristretta. Mi riservo di gustarne pian piano il dettaglio fra pochi giorni, quando mi sarò sbrigato del piccolo contributo che debbo ultimare al più presto. Intanto ho rilevato da una rapida scorsa che hai con grande benevolenza fatto posto anche al mio criterio di raccordo fra l'elettrodinamica classica e la maxwelliana, e te ne sono schiettamente obbligato.

Al primo vederci avrò certo tratto dalla lettura più di uno spunto per avere da te istruttivi complementi su personali punti di vista o accenni, volutamente riservati e parchi.

Ora mi limito a porgerti i più vivi e cordiali ringraziamenti con rinnovati auguri e affettuosi saluti

Tuo obbl.mo T. Levi-Civita

---

[142] Questa lettera si trovava nelle *Rusticationes*.

[143] Probabilmente si tratta di: G.A. Maggi, *Teoria fenomenologica del Campo Elettromagnetico. Lezioni di Fisica Matematica*, Milano, Hoepli, 1931.

**105**
Tullio Levi-Civita a Gian Antonio Maggi[144]
[cartolina indirizzata a: Chiar.$^{mo}$ Sig. Prof. Gian Antonio Maggi - Valnegra (Bergamo)]

Padova 23-IX-1930

Mio caro amico,

Mille grazie per la commemorazione che ho letta con grande interesse ed anche, permettimi di aggiungere, con sincera ammirazione. Essa non soltanto è affettuosa, sentita, efficacissima pel sapiente e misurato concorso dell'elemento umano e dell'illustrazione scientifica, ma contiene altresì un saggio magnifico di pensiero e di stile nelle pag. 7 ed 8 in cui hai tratteggiato le caratteristiche della fisica teorica odierna.[145]

Coi più cordiali saluti, anche da parte di mia moglie, ti stringo la mano

Tuo aff. T. Levi-Civita

**106**
Tullio Levi-Civita a Gian Antonio Maggi
[busta indirizzata a: Chiar.$^{mo}$ Sig. Prof. Gian Antonio Maggi -
Corso Plebisciti, 3 - Milano (121);
busta e carta intestate: Prof. T. Levi-Civita - Via Sardegna, 50 - Roma (25)]

2-I-1931

Carissimo amico,

Mi sono lasciato prevenire dalla tua amichevole cortesia. Non ho giustificazioni; ma, pur secondo in tempo, ti porgo con eguale affetto i più sinceri e lieti auguri pel 1931, anche da parte di mia moglie.

Avevo una vaga speranza di vederti Domenica prossima alla seduta dei Lincei, ma non c'è cenno di ciò nella tua lettera, e temo quindi che la tua venuta non sia così prossima, tanto più forse che il soggiorno a Parigi della tua figliuola costituirà ora una ben forte attrazione verso Nord Ovest.

Nella seduta del 7 X$^{mbre}$ u.s. Castelnuovo lesse ai Lincei una bella rievocazione dell'opera del Cremona. So che ne ha già corretto le bozze e che la commemorazione apparirà nel fascicolo corrispondente alla detta seduta.[146]

Di commemorazioni tenute a Napoli non ho notizia. Me ne informerò (al più tardi da Signorini, che sarà di passaggio fra pochi giorni) e ti riferirò.

Debbo anch'io (assicuro all'Autore) scusarmi a proposito della nota di Gugino sulla determinazione delle reazioni. Il risultato è molto bello e notevole, e io lodai l'autore, presentandogli con piacere la Nota ai Lincei. Viceversa tu non solo avevi avuto l'idea, ma la avevi anche realizzata, in certo senso, con maggiore generalità, considerando come elemento primordiale il corpo rigido, anziché la consueta astrazione del punto materiale.[147] Il mio torto è di essermi a suo tempo lasciato sfuggire un così fondamentale § del tuo libro. Forse (ma in altro senso) un po' di torto ce l'hai anche tu, nel non aver richiamata con pubblicazioni speciali l'attenzione degli studiosi sopra un complemento che quanto più ci penso, tanto più mi pare essenziale del principio dei lavori virtuali.

[144] Questa lettera si trovava nelle *Rusticationes*.
[145] Potrebbe trattarsi di: G.A. Maggi, “In memoria del Prof. Aldo Pontremoli. Discorso”, *Annuario della R. Università di Milano per l’a.a. 1929-1930*, pp. 156-162.
[146] G. Castelnuovo, "Luigi Cremona nel centenario della nascita. Commemorazione", *Rendiconti della Reale Accademia dei Lincei*, s. VI, v. XII, 1930, pp. 613-618.

Ho preso atto con molto interesse di quanto mi comunichi circa la tua fervida attività al Seminario, all'Università, all'Istituto Lombardo.

Vado anch'io leggendo quà *[sic!]* e là sulla nuova meccanica, ma non conosco, cioè non ho neanche mai sfogliato, il libro del Bloch, di cui già questo autunno avevo sentito dir bene dal Palatini. In questi giorni ho letto le litografie del Persico, che sono sintetiche e molto perspicue: come forse già sai, il Persico passerà ora da Firenze a Torino. Qui (in seguito al concorso di Catania per Fisica sperimentale) abbiamo chiamato il giovane Rasetti per la spettroscopia.

Aspettavo a mandarti una comunicazione sulla teoria del Prandtl[148] (fatta ad Aachen l'anno scorso) assieme a quella tale osservazione sul comportamento asintotico del potenziale di masse estendentisi all'∞, di cui ti scrissi da Padova e di cui dovrei presto ricevere gli estratti. Giacché ne ho l'occasione, ti mando intanto la nota tedesca.

Ho presentato il mese scorso ai Lincei una nota (piuttosto elementare, come immagini facilmente) sui centri di gravità (!):[149] non dispero che possa interessarti, se è vero che la questioncina non era venuta in mente ad altri, ciò che risulterebbe dalle ricerche che ho fatto con la maggior diligenza possibile.

Saluti cordiali anche da mia moglie; rinnovati, ottimi auguri e una stretta di mano

Dal tuo aff.
T. Levi-Civita

[147] Si vedano le lettere #82 e 83.
[148] Gli Atti del Convegno *Vorträge aus dem Gebiete der Aerodynamik und verwandter Gebiete* tenuto ad Aachen nel 1929 sono editi da Springer Berlin e si trovano in rete (www.springer.com/de/book/9783662333945). Contengono tre contributi di L. Prandtl e uno di Levi-Civita.
[149] Potrebbe trattarsi di: T. Levi-Civita, "Sezioni piane di un corpo e direttrici ortobariche", *Rendiconti Accademia dei Lincei*, s. 6, v. XII, 1930, pp. 535-541.

**107**
Gian Antonio Maggi a Tullio Levi-Civita

Milano 28 Settembre 1934

Carissimo amico,

Ti mando insieme una copia della mia Nota col Finzi "Condizioni d'esistenza delle onde elettromagnetiche armoniche",[150] che vorrai aggradire colla consueta tua bontà. Di mio propriamente non ci sono che il motivato riconoscimento della inconsistenza della dimostrazione accennata nella citata Nota dell'Istituto "Soluzione del problema ecc.", e le macerie della precedente costruzione, fondata sulla dimostrazione medesima. Il resto è opera del Finzi. Per cui io rimetto a lui la distribuzione della Nota, non mandandola direttamente che a pochissimi - sei o sette, tra cui tu, che, ad ogni modo, avrei sempre considerato come privilegiato - ai quali mandai la suddetta Nota dell'Istituto Lombardo e le connesse lincee, prima di riconoscere l'errore. Fortunatamente avevo trattenute le Note lincee. Per combinazione, lo stesso giorno che il Finzi, postosi all'opera, su mio suggerimento, mi portava le sue conclusioni, io ricevevo dal Prof. Sona[151] notizia del suo risultato, oggetto poi della Nota ch'io presentai per lui ai Lincei, di cui quelle conclusioni sembrano la primitiva ragione. Il Sona, che, a suo tempo, avvalorò la mia costruzione col ricondurre le equazioni da me aggiunte alle condizioni del Love, nell'ultima Nota, aggiungendo $\frac{\delta\Omega}{\delta n} = 0$, sulla superficie d'onda, alle $F=\mathrm{grad}\Omega$, $\Delta_2\Omega=0$ da me trovate, dimostrava, con un calcolo che pensò essere assai *[azzeccato]*, che le mie onde non sono suscettibili, in generale, di propagazione. Sic transit gloria mundi! A questo risultato ho esplicitamente accennato in un'accompagnatoria, che lessi presentando la Nota col Finzi all'Istituto, e si trova stampata, a processo verbale, nei Rendiconti, oltre il richiamo che ne ò *[sic!]* fatto nella Nota medesima. I risultati del Finzi furono però trovati senza notizia da quello del Sona, e ne è manifesta la diversa impronta, in seguito semplicemente alla mia constatazione della inconsistente dimostrazione, della quale il Sona non fa cenno.

Avrai ricevuto a Pocol l'illustrata da Lanzo che ti mandai in ricambio della graditissima tua. Non ho poi notizia de' tuoi successivi spostamenti, per cui, per non sbagliare, indirizzo questa mia a Roma. Noi abbiamo proseguito la nostra villeggiatura di Lanzo sempre meglio, avendoci ivi raggiunto la Bianca[152] colle due bambine, reduci dalla spiaggia della Manche. In questo momento ho in casa tutto Parigi, e ben vorrei non fosse che per pochi giorni ancora.

Ti prego, terminando, di far aggradire alla gentile Signora i miei più cordiali rispetti, e tu ricevi un abbraccio dal

Tuo aff. amico
Gian Antonio Maggi

Al Ch.° Prof. Tullio Levi-Civita della R.Università di Roma
Via Sardegna 50.

---

[150] In *Rendiconti del R. Istituto Lombardo*, s. II, v. LXVII, 1934, pp. 331 e 363-370. Si tratta di una Nota in correzione della precedente "Soluzione del problema della riflessione e rifrazione delle onde elettromagnetiche armoniche di forma qualsivoglia ad una superficie piana", pubblicata sul medesimo volume alle pp. 199-204 e 331. Tale Nota ed altre di complemento ad essa si trovano pubblicate tra il 1933 e il 1936 sui *Rendiconti dei Lincei* e su quelli dell'Istituto Lombardo.

[151] Si veda la lettera di Maggi a Finzi.

[152] Cognata di Maggi.

**108**

Tullio Levi-Civita a Gian Antonio Maggi

[cartolina illustrata indirizzata a: Italia-Milano (21) - Corso Plebisciti, 3

Chiar.$^{mo}$ Sig. Prof. Gian Antonio Maggi]

Kiev 8-VI-1935

Coi più affettuosi saluti
T. Levi-Civita

Illustre e caro Maestro
*[Segue scritta in russo, trascritta e tradotta su un foglietto da Maggi:]*

Traduzione.

Le lezioni del celebre scienziato italiano Levi-Civita mi rammentano quel tempo felice della mia gioventù, quando, molti anni fa, ebbi la fortuna di ascoltare lezioni nella Università di Pisa, e vicino ai Suoi studenti, ascoltare le lezioni di Elettromagnetismo.

*[Maggi spiega inoltre:]*
Righe del prof. K[...] aggiunte a illustrata di Levi-Civita “Coi più affettuosi saluti” in data di Kiev, 8, VI, 1935, ricevuta il 17.
N.B. L’aggiunta in alto nella terza riga è nel testo.

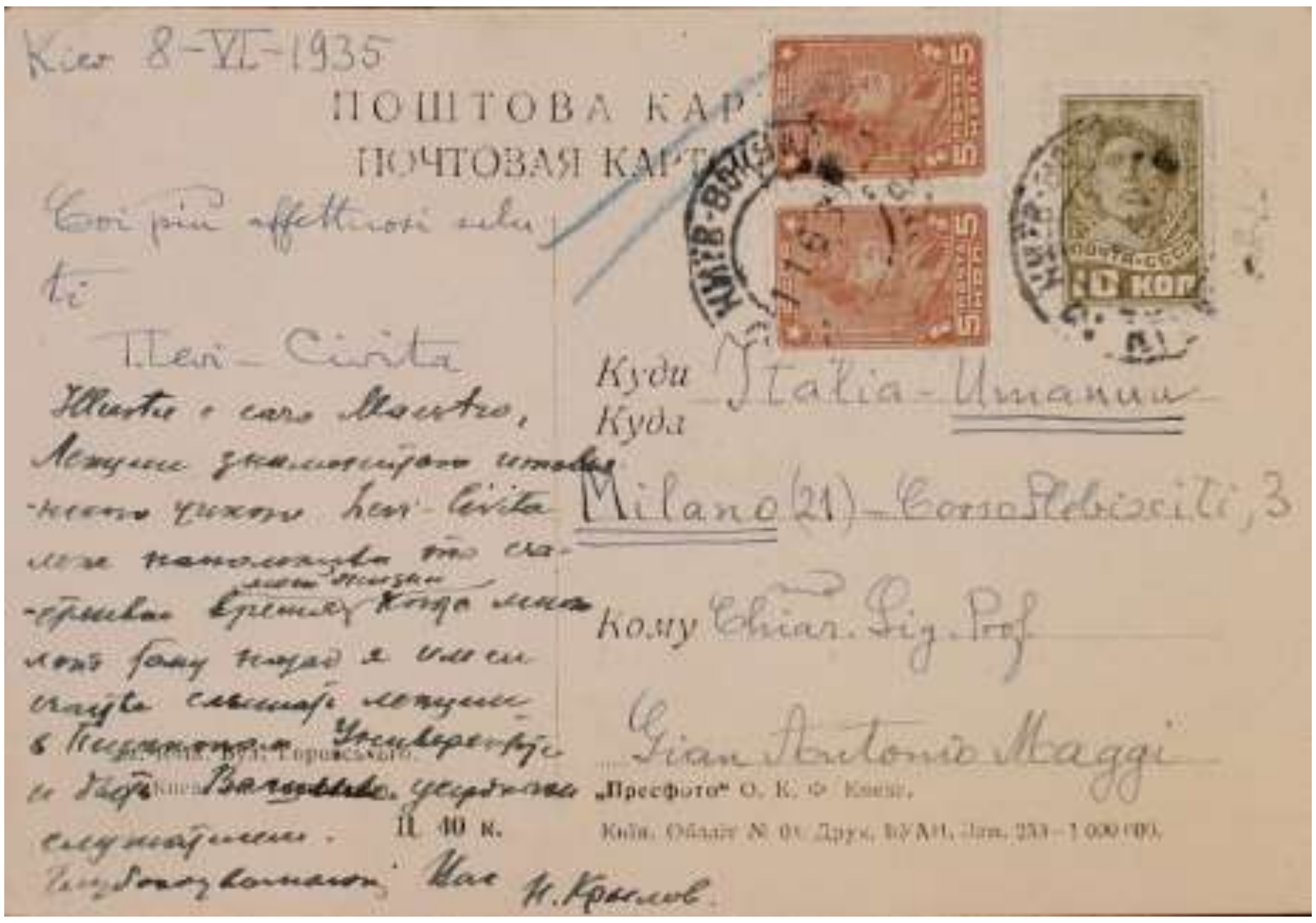

*Frammento probabilmente indirizzato al medesimo prof. K. della cartolina precedente:*

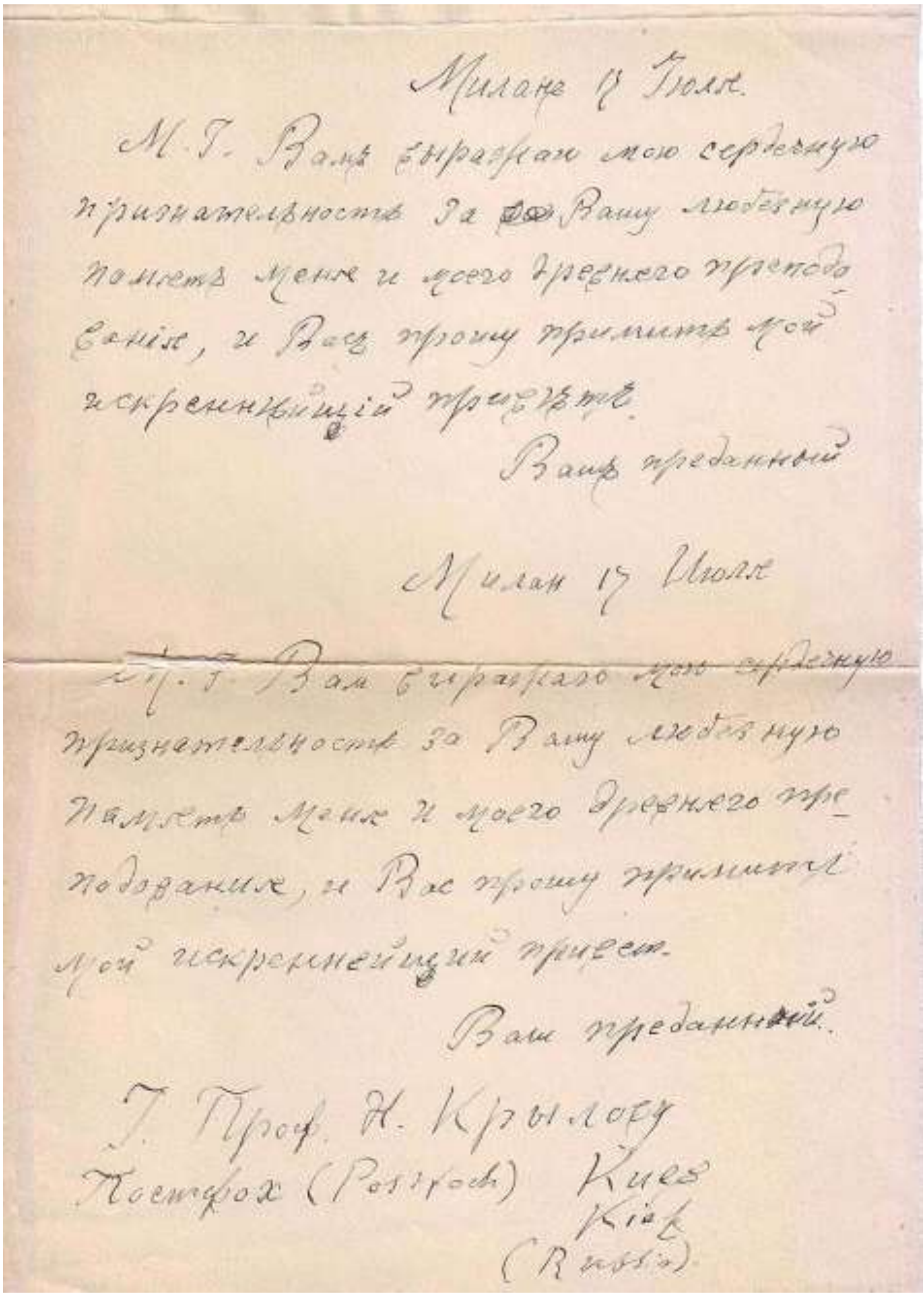

Милане 17 Июля.

М. Г. Вамъ выражаю мою сердечную признательность за Вашу любезную память меня и моего древнего преподаванія, и Васъ прошу принять мой искреннѣйшій привѣтъ.

Вашъ преданной

Милан 17 Июля

М. Г. Вам выражаю мою сердечную признательность за Вашу любезную память меня и моего древнего преподавания, и Вас прошу принять мой искреннейший привет.

Вам преданный.

Г. Проф. Н. Крылову
Поштфох (Postfach) Киев
Kiev
(Russia)

**109**
**Antonino Lo Voi** a Gian Antonio Maggi
[inviato a: Ch^mo^ Prof. Gian Antonio Maggi
R. Università Milano; data dal timbro postale: 29/11/1929]]

Omaggio con preghiera di darne giudizio.

## Sui fenomeni di capillarità

1. — In tutti i libri di fisica sperimentale che ho visti ho trovato spiegato i fenomeni di capillarità nel modo seguente.

Se la forza di adesione supera la forza di coesione, considerando una molecola che tocca le pareti, e si trova alla superficie terminale, si vede che essa è soggetta a due forze; una orizzontale dovuta alla adesione delle pareti; l'altra diretta verso la massa del liquido, dovuta alla forza di coesione.

Facendo la risultante delle due forze, si ha una forza diretta verso l'esterno, e si conclude che in vicinanza delle pareti, il liquido si deve sollevare.

Orbene, quasi tutti gli insegnanti si sono sempre accorti che questo ragionamento è sbagliato. Così; pensando l'adesione uguale alla coesione, la precedente costruzione della risultante delle due forze, porta all'assurdo che il liquido si dovrebbe sollevare.

Vogliamo qui dare una spiegazione elementare del fenomeno; e ricavarne risultati quantitativi.

2. — Consideriamo una molecola che tocchi la parete; e la supponiamo ad una distanza dalla superficie libera dal liquido, che superi $\frac{1}{10000}$ di mm. Siamo allora sicuri che le molecole del liquido che esercitano una azione di coesione sensibile sulla nostra molecola, riempiono *tutta* una emisfera avente il centro nella molecola. Dunque, la forza di coesione che si viene ad esercitare sulla molecola considerata *è orizzontale*, come lo è la forza di adesione dovuta alle pareti.

Ne consegue che se, per esempio, l'adesione supera la coesione, il liquido in vicinanza delle pareti verrà soggetto ad una pressione maggiore di quella che si ha in punti interni alla massa posti sullo stesso orizzonte. Poichè la pressione cresce per i peso della colonna del liquido sovrastante, tale colonna deve essere maggiore in vicinanza delle pareti.

3. — Consideriamo, per esempio, un tubo cilindrico, di raggio $r$ verticale, contenente un liquido la cui adesione alle pareti superi la coesione di una quantità $p$. Per ogni centimetro di altezza del liquido, si ha in totale un incremento della pressione misurato dal prodotto della lunghezza della circonferenza per $p$:

$$2\,\pi\, r\, p.$$

Se dunque il cilindro capillare ed il suo fondo è in comunicazione con una vaschetta, il livello del liquido dentro il cilindro supera il livello del liquido nella vaschetta; e precisamente, se chiamiamo con $h$ il dislivello e con $p^1$ il peso specifico del liquido, deve essere

$$p.^1\, r^2\, h = 2\, r\, p;$$

da cui

$$h = \frac{k}{r};$$

ove, si è posto $k = \frac{2\,p}{p^{\,1}}$.

Quindi, *il sollevamento del liquido è inversalmente proporzionale al raggio del cilindro capillare; direttamente proporzionale alla differenza fra la forza di adesione e la forza di coesione* (tale differenza può essere negativa); *ed è inversamente proporzionale al peso specifico del liquido.*

4. — Chiudiamo con una osservazione banale.

Alla superficie si forma un menisco e le forze di tensione superficiale avranno una risultante diretta verso l'alto. Tale risultante provocherà quindi un sollevamento $h^1$ (inversamente proporzionale al raggio) che viene a sommarsi al sollevamento $h$.

*Milano*, 16 - 11 - 29 - (A. VIII).

*Antonino Lo Voi*
Liceo Sc.

**110**
Gian Antonio Maggi a Antonino Lo Voi
[carta intestata: Seminario matematico e fisico - Milano]

Milano 5 dicembre 1929
Corso Plebisciti, 3.

Egregio Signore,

Accogliendo la cortese domanda del mio consiglio sul Suo articolo, del cui invio La ringrazio, eccoLe alcune mie riflessioni in proposito.

Nelle ordinarie esposizioni elementari dei fenomeni capillari, si sogliono distinguere le tre ipotesi che l'adesione del liquido alle pareti sia maggiore, minore o uguale alla coesione, nelle quali la risultante è rispettivamente volta all'infuori della parete, all'indentro, o nullo, e la superficie libera, come perpendicolare alla risultante di questa e della gravità, ne risulta, rispettivamente concava, convessa e piana. Non mi rendo quindi ragione della deduzione assurda, accennata nel Suo scritto, che il liquido si dovrebbe sollevare, nell'ipotesi dell'eguaglianza dell'adesione e della coesione. Perché l'altezza del livello del liquido nel tubo si subordina poi alle suddette tre forme della superficie libera, e invocando la tensione superficiale, la quale ha componente verticale volta in alto, nel primo caso, volta in basso, nel secondo, e nulla, nel terzo, si trova, nei primi due casi, rispettivamente, elevazione e depressione, e, nel terzo caso, il livello semplicemente dovuto alla gravità, che è il livello all'estremo dal tubo.

Sulla spiegazione proposta ai §§2,3 credo basti osservare che la risultante dell'adesione o della coesione relativa alla molecola interna, emisferica, aderente alla parete, riuscendo orizzontale, non può essere portata in conto, per render ragione della maggiore o minore altezza del livello del liquido nel tubo, che è subordinato alla gravità e alle componenti verticali, come la gravità, delle altre forze, tali perciò da sottrarli dalla gravità (elevazione) o cospirare con essa (depressione).

Difatti, diversamente da quanto è accennato al §4, l'elevazione, o la depressione, si spiegano cola sola tensione superficiale, che introduce la suddetta componente verticale, volta in alto, o in basso, secondo che la superficie libera è concava o convessa.

Aggradisca con questo i miei distinti saluti e mi creda, mi affermo

Dev$^{mo}$ Suo
G.A. Maggi.

**111**
Antonino Lo Voi a Gian Antonio Maggi
[busta indirizzata a: Ch$^{\underline{mo}}$ Prof. Antonio Maggi - Corso Plebiscito 3 - Milano]

Ch$^{\underline{mo}}$ Prof Maggi,

sentitamente La ringrazio della Sua lettera.

Anche il prof. La Rosa, pur riconoscendo vero quanto è contenuto nel n°1, dice che bisogna fare ricorso alla tensione superficiale. Fu allora che io aggiunsi l'osservazione del n°4.

Però, io penso, che la tensione superficiale è una forza interna che influisce grandemente sulla forma della superficie libera. Ma: "ad ogni forza interna corrisponde una forza interna uguale e contraria, se per un istante si suppone il liquido irrigidito"; pertanto, la tensione superficiale, io credo, che non possa spiegare il sollevamento; fatto eminentemente esterno.

Gradisca i più distinti ossequi e ringraziamenti

dev$^{\underline{mo}}$
Antonino Lo Voi

Milano 10-12-29-(VIII)

**112**
Gian Antonio Maggi a Antonino Lo Voi

Casa 21 Dicembre 1929

Egregio Signor Professore,

Non mi saprei veramente spiegare come il prof. La Rosa riconosca vero quanto è contenuto nel §1. Comunque sia, l'accennato assurdo non segue dall'apparizione dei fenomeni capillari che si trova in diversi trattati di mia conoscenza e che, nella mia lettera ho brevemente riassunto. Per quanto alla tensione superficiale, le tensioni superficiali relative ai varii elementi della superficie del liquido non si elidono mutuamente, come è facile riconoscere; e danno quindi luogo ad un risultante non nullo, il cui componente verticale, sottraendosi dal peso, col menisco concavo, e aggiungendovisi, col menisco convesso, rende ragione dell'elevazione e rispettivamente della depressione del livello del liquido nel tubo.

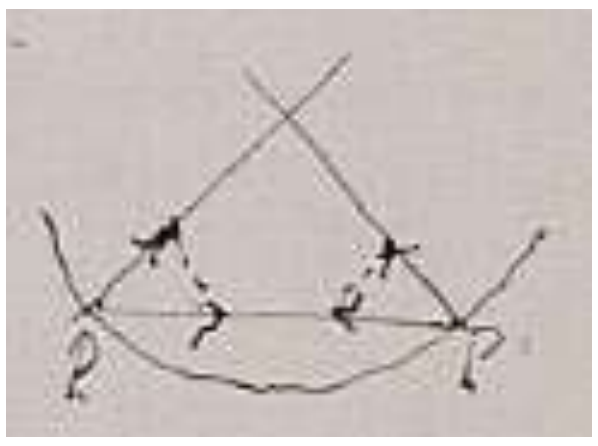

Con nuovi distinti saluti

Dev.$^{mo}$ Suo
G.A. Maggi.

**113**

**Gaetano Magnanini** a Gian Antonio Maggi
[busta indirizzata a: Illustr.mo Comm. Prof. O. G.A. Maggi - Stabile di Meccanica razionale Facoltà Matematica R. Università Milano; bollo postale: 3/12/1929; cartolina con ritratto in accompagnamento al dattiloscritto]

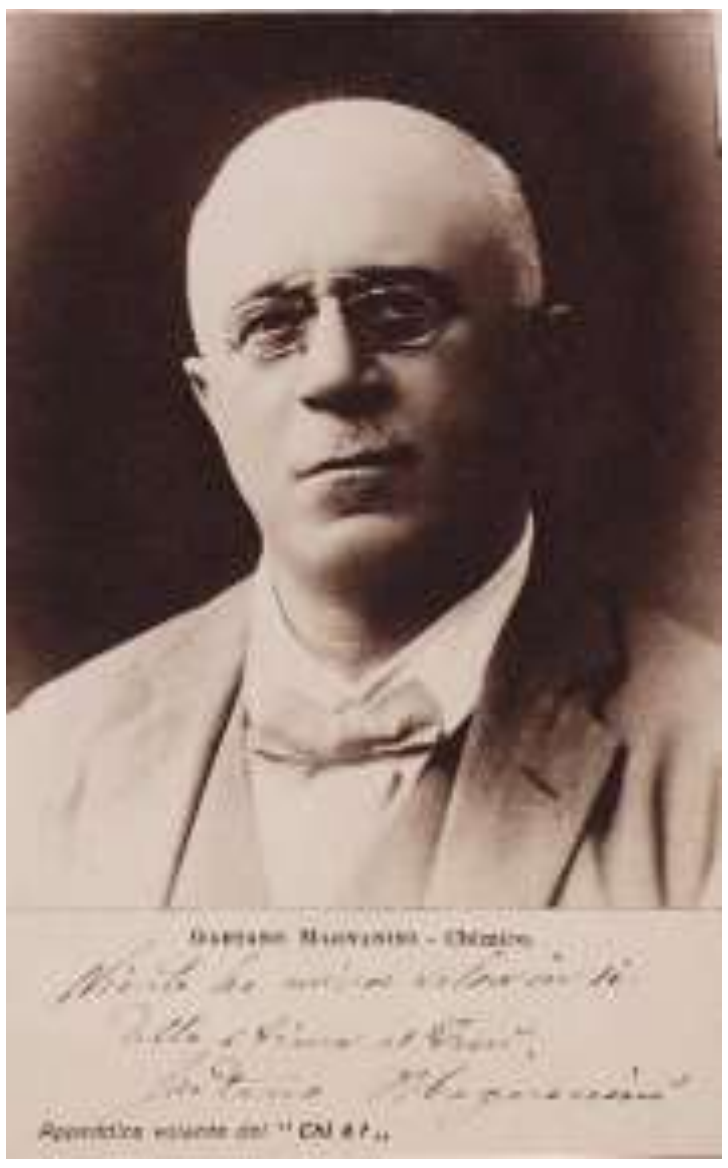

Niente ha minor valore in sé della stima altrui.
Gaetano Magnanini

Nessun rapporto ha il motto col comune apprezzamento soggettivo della stima altrui la quale è cara a tutti. Esso riguarda invece il valore oggettivo della stima altrui, considerato in sé, perché, come nella relatività dello spazio e del tempo, avviene della stima altrui che essa dipende dal sistema di riferimento.

3/11/29[153]

Cariss. ed Illustre
Magnanini desidera, non disturbarti, ne vuole lunghe cose sempre per non darti disturbo.
Ma desidera solamente sapere se ha torto e perché, ovvero se per caso ha ragione. Sono in lite con un relativista che non mi pare abbia capito la relatività. Saluti cari.

DICE EINSTEIN (Trad. Calisse) Edizione Zanichelli pag. 21 -

§9. RELATIVITÀ DELLA CONTEMPORANEITÀ

Finora abbiamo considerato un certo corpo come termine di riferimento, che abbiamo designato come sede ferroviaria. Supponiamo ora, che sul binario cammini con velocità costante v un lungo treno nella direzione indicata nella Fig. 1.

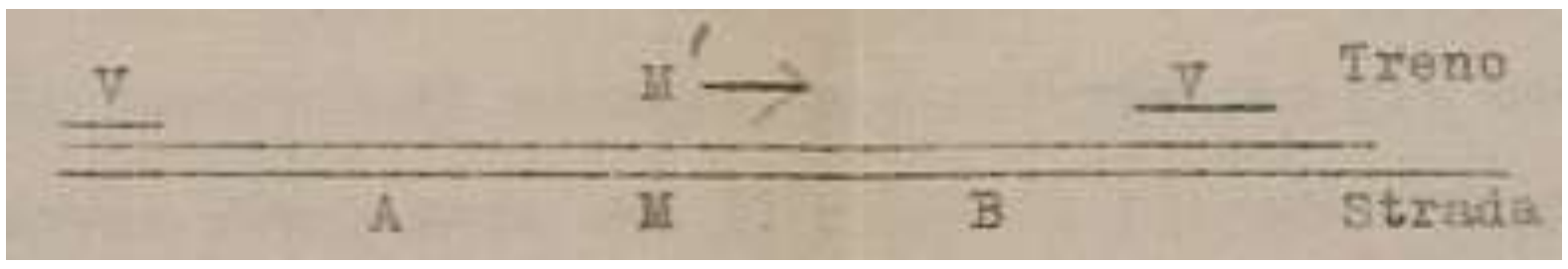

Le persone che viaggiano nel treno, si valgono del treno come corpo rigido (sistema di coordinate), al quale riferiscono tutti gli avvenimenti. Ogni avvenimento, che ha luogo lungo il binario, ha luogo altresì in un determinato punto del treno. Anche la definizione della contemporaneità si applica rispetto al treno precisamente come rispetto alla sede stradale.
Sorge allora spontanea la seguente questione:
"Due avvenimenti (p. es. i due colpi di folgore A e B), contemporanei rispetto alla strada, sono anche contemporanei rispetto al treno?" Possiamo subito mostrare, che la risposta deve essere negativa.

[153] Il mese è dicembre dal bollo postale.

Se diciamo, che i colpi di folgore A e B sono contemporanei, ciò significa che i raggi luminosi uscenti dai luoghi A e B s'incontrano nel punto medio M del tratto di strada AB. Ma agli avvenimenti A e B corrispondono anche due luoghi A e B sul treno. Sia M' il punto medio del tratto AB sul treno in movimento. Nello stesso istante del colpo (giudicato dalla strada) il punto M' coincide con M, ma esso si muove con la velocità del treno v verso destra. Se un osservatore situato in M' sul treno non possedesse questa velocità v, esso rimarrebbe immobile in M, e i raggi luminosi uscenti da A e B arriverebbero a lui contemporanei, ossia questi raggi s'incontrerebbero al suo posto. Ma effettivamente quell'osservatore (giudicato dalla strada) muove incontro al raggio uscente da B ed avanza rispetto a quello proveniente da A. Perciò egli deve vedere prima il raggio uscente da B e poi quello uscente da A. Gli osservatori adunque, che riferiscono gli avvenimenti al treno, devono concludere che il colpo di folgore B è avvenuto prima del colpo di folgore A. Arriviamo così al seguente importante risultato: GLI AVVENIMENTI CHE SONO CONTEMPORANEI RISPETTO ALLA STRADA, NON SONO CONTEMPORANEI RISPETTO AL TRENO, E VICEVERSA (relatività della contemporaneità). Ogni sistema di riferimento ha il suo proprio tempo. Un dato temporale non ha alcun senso, se non si assegna il sistema, al quale il tempo si riferisce.

%=%=%=%=%=%=%=%=%=%

MAGNANINI dice invece:
che, l'osservatore che si trova sul treno vede prima il raggio che sorte da B, <u>ma non deve concludere</u> che l'avvenimento in B ha preceduto quello in A, perché <u>se è intelligente, deve capire</u> che muovendosi verso B avrebbe dovuto vedere sempre prima l'avvenimento accaduto in B, anche se i due avvenimenti fossero stati, come sono stati, realmente contemporanei. Dunque la <u>relatività</u> della contemporaneità è una <u>illusione</u>!

**114**
Gian Antonio Maggi a Gaetano Magnanini
[carta intestata: R. Università di Milano - Corso Roma, 10]

Milano, 5 Dicembre 1929

Carissimo Magnanini

Al passo dell'Einstein da te citato non va cercato di più che mostrare come due avvenimenti che risultano contemporanei, giudicati da un riferimento fisso, possono non risultare più tali, giudicati da un riferimento mobile, rispetto al primo, considerato come fisso. Magnanini, secondo me, ha ragione di trar in scena l'uomo intelligente che scopre il trucco. Ciò non toglie che, per l'esatta teoria dello stesso Einstein, la relatività della contemporaneità non sia un'illusione, come tale che nessun uomo intelligente riuscirebbe ad accordare la contemporaneità di due eventi, giudicati da due riferimenti, in movimento, l'uno per rispetto all'altro. Essa risulta dalla trasformazione di Lorentz, che non è che una traduzione in formole dell'esperienza di Michelson; traduzione e deduzione alle quali non saprei che obbiezioni si potrebbero fare. E tanto è vero che il passo da te discusso non si identifica colla teoria esatta, che, per avere la differenza dei due tempi fornita da questa, bisogna dividere la differenza fornita da quella per $\sqrt{1-\left(\frac{v}{c}\right)^2}$, c veloc. della luce, v veloc. del treno; dove $\left(\frac{v}{c}\right)^2$ è inapprezzabile in confronto di uno, ma, per escludere l'identificazione non importa.
Credo con questo d'aver appagato anche il tuo desiderio di cose non lunghe; ma, se desideri dell'altro, resto, con piacere, a tua disposizione. Intanto aggradisci i migliori saluti, e credimi

Tuo aff. collega

*[scritto in penna rossa:]* Non mandata questa parte

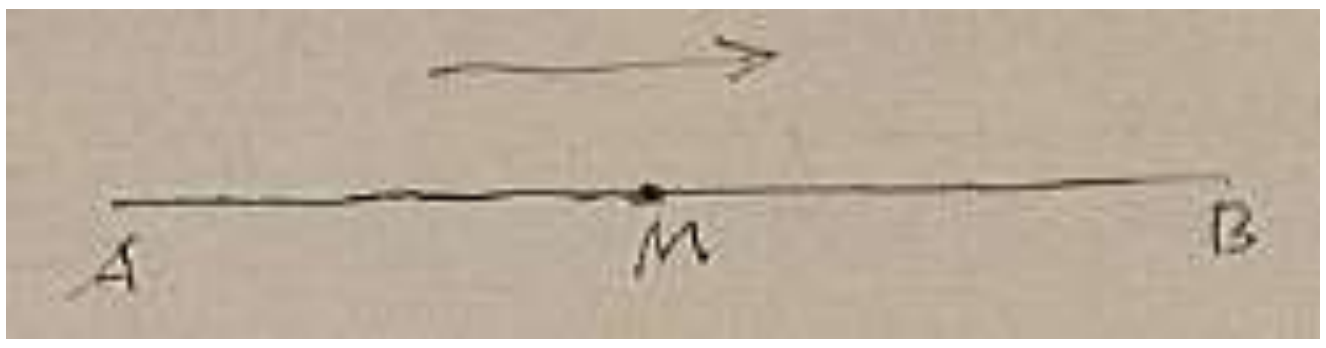

AB=s, AM=$\frac{s}{2}$
La luce impiega a percorrere lo spazio AM=BM=$\frac{s}{2}$ il tempo $\frac{s}{2}\frac{1}{c}$.
In questo tempo M', punto del treno coincidente con M al momento dello scoppio delle due folgori, avanza da A verso B dello spazio $\frac{s}{2}\frac{1}{c}$v.
La luce impiega a percorrere questo spazio il tempo $\frac{s}{2}\frac{v}{c}\frac{1}{c}=\frac{s}{2}\frac{v}{c^2}$.
La differenza dei tempi in cui sono percepite in M' le due folgori è quindi
$2\frac{s}{2}\frac{v}{c^2}=\frac{v}{c^2}s$.
La trasformazione di Lorentz fornisce invece per questa differenza

$$\frac{\frac{v}{c^2}s}{\sqrt{1-\left(\frac{v}{c}\right)^2}}.$$

Notisi che $\frac{s}{\sqrt{1-\left(\frac{v}{c}\right)^2}}$ è la distanza dei due punti A e B giudicata dal treno.

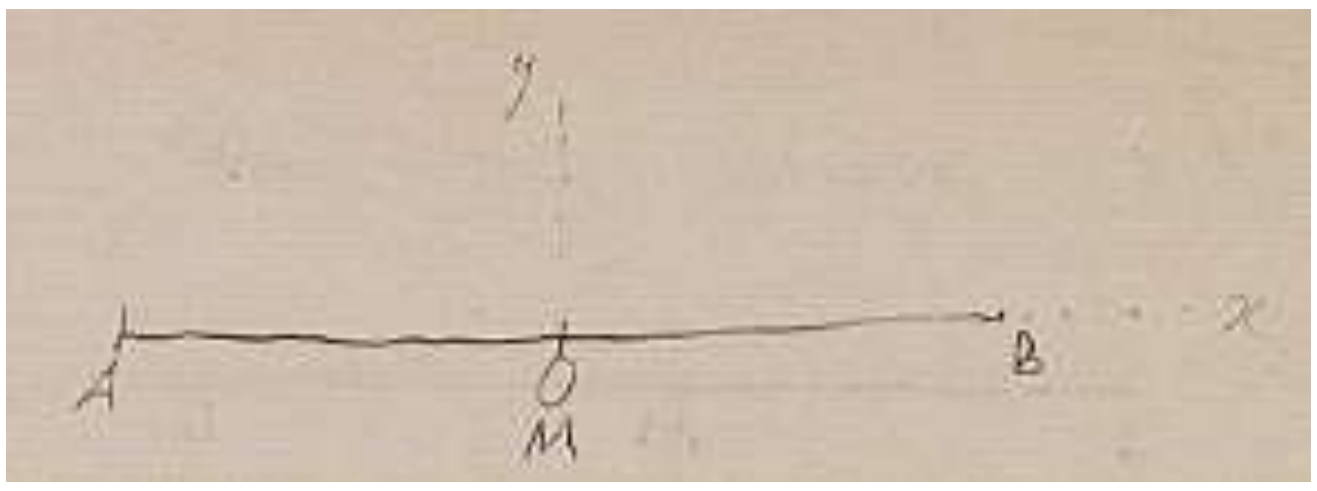

Il tempo t essendo assegnato all'istante dello scoppio delle due folgori in A e in B, nel riferimento alla linea, a quei due tempi, negli stessi punti, nel riferimento al treno, corrispondono

$$t'_1=\frac{t-\frac{s}{2}\frac{v}{c^2}}{\sqrt{1-\left(\frac{v}{c}\right)^2}}, \qquad t'_2=\frac{t+\frac{s}{2}\frac{v}{c^2}}{\sqrt{\phantom{x}}}$$

donde (invece di t−t=0)

$$t'_2 - t'_1 = \frac{s\frac{v}{c^2}}{\sqrt{1-\left(\frac{v}{c}\right)^2}}.$$

**115**
Gaetano Magnanini a Gian Antonio Maggi

Modena 7/12 29

Car.mo ed Illustre,

Ti ringrazio moltissimo degli schiarimenti. Per me che, pure avendo fatto il biennio di matematica, mi sento assai poco matematico; la cosa si presenta nella mia mente in questi termini:
La teoria è essenzialmente un sommo lavoro di matematica, cosiché le variazioni di spazio e tempo che accompagnano tale teoria devono essere giudicate come variazioni matematiche, non come variazioni fisiche. Tali variazioni matematiche sono indispensabili affinché il sistema di formole e di calcoli relativistici possano con successo essere applicati alla soluzione dei problemi del Cosmo; sopratutto *[sic!]* per i fenomeni celesti dove tutti i movimenti sono eminentemente relativi. -Non è dunque: che si contraggono realmente i regoli, e si rallentino gli orologi. Ma le cose avvengono quantitativamente come se tali contrazioni avvenissero. - Così si soddisfa al bisogno costante di quantitatizzare che è proprio dello spirito dell'uomo. Io ho capito bene o male così!-?!-Saluti cordialissimi aff. Magnanini.

**116**
Gian Antonio Maggi a Gaetano Magnanini
[carta intestata: Seminario matematico e fisico - Milano]

Milano, 9 Dicembre 1929

Carissimo Magnanini,

Trovo che quello che mi scrivi sta bene, e il tuo concetto della valutazione dello strumento matematico della Relatività per l'analisi del Cosmo dovrebbe soddisfare anche lo stesso Einstein, che non vuol sentire che la sua dottrina sia piuttosto matematica che fisica. Credo poi che ti possa interessare una breve aggiunta alla mia precedente lettera. Nel passo citato, come ti accennavo, l'Einstein adotta provasoriamente *[sic!]* un certo criterio, per decidere della contemporaneità, e mostra che, in base ad esso, due avvenimenti giudicati contemporanei con un riferimento non saranno giudicati tali con un riferimento diverso. Attribuendo una propagazione diversa alla luce nel riferimento alla linea ferroviaria e nel riferimento al treno, ne segue che l'osservatore intelligente del Magnanini riconosce il movimento del treno, e, come dicevo, scopre il trucco. La relatività della contemporaneità stabilita a questo modo si può ben chiamare un'illusione. Tutt'altra cosa è la relatività della contemporaneità che scaturisce dai principii della relatività che l'Einstein, come avrai veduto, nello stesso opuscolo, introduce più tardi. Con questi principii le stesse modalità appartengono alla propagazione della luce sia nel riferimento alla linea ferroviaria sia nel riferimento al treno. Però la caduta delle due folgori nei due punti indicati, se è contemporanea con un riferimento, non è più nell'altro, perché due avvenimenti contemporanei in due punti, con un riferimento, non sono più nell'altro, salvo che siano in uno stesso piano, perpendicolare alla veloc. di trascinamento.
Con sinceri cordiali saluti

Affez tuo
GAM.

**117**
Gaetano Magnanini a Gian Antonio Maggi
[busta indirizzata a: Illustre Comm. Professore G.A. Maggi
della R. Università di Milano - Corso Via Plebiscito;
data del timbro postale: 23/12/1929; dattiloscritto]

Estratto da LAMMEL FONDAMENTI DELLA RELATIVITÀ (ZANICHELLI) Pag.63
DICE LAMMEL:
"Immaginiamo due osservatori A e B in una medesima posizione, in uno stesso istante, e che in quella stessa posizione e stesso istante venga emesso un raggio di luce. I due osservatori si comportino diversamente:
A, rimane al suo posto; B, si muove verso destra dall'istante in cui venne emessa la luce, e così pure si muove con assai maggior velocità il fronte di luce L; cioè A, B e L che coincidevano in un dato istante iniziale; si vanno distanziando col tempo

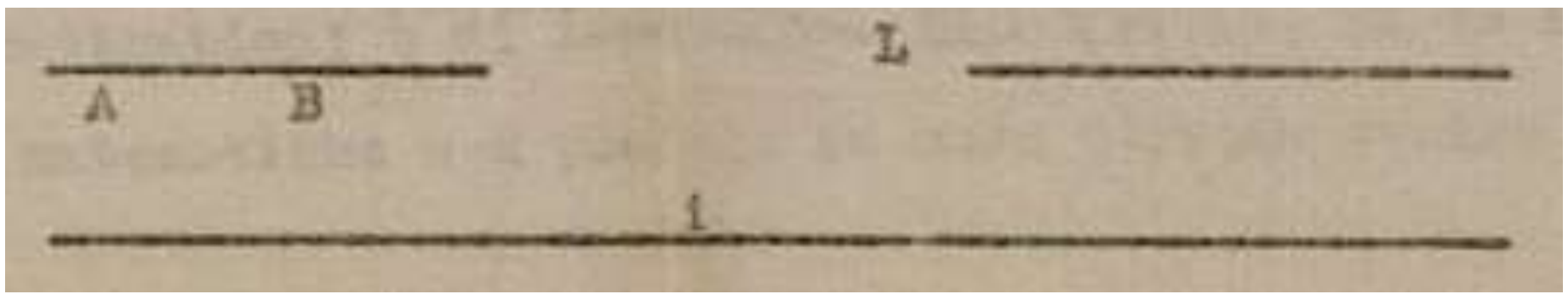

In un altro istante, trascorso il tempo t, i due osservatori determinano ciascuno la velocità della luce che trovano entrambi identica =C
Quale segreto nasconde quì la natura?
L'osservatore A dice: devo dividere il tratto AL trascorso dalla luce per il tempo t, ed ottengo la velocità AL/T=C; ma l'osservatore B ragionando analogamente dice: devo dividere il tratto BL che ha percorso la luce durante lo stesso tempo t e ho BL/T=C; grandezze diverse divise per lo stesso numero non possono dare quozienti eguali! Quì in Einstein si è condensato il pensiero di un epoca *[sic!]*. Se le due divisioni danno lo stesso numero sospettiamo che la grandezza t non sia la stessa cioè "Un determinato istante non è una cosa determinata. È pregiudizio ritenere che contemporaneo abbia un significato determinato."

---

DICE INVECE MAGNANINI:
La natura quì non nasconde nessun segreto, ma è chiara e sincera. La velocità della luce essendo sempre come tutte le velocità eguale allo spazio percorso nella unità di tempo, C ha lo stesso valore per A e per B; solamente che per A tale valore è AL/T=C invece l'osservatore in moto B deve voltarsi indietro, e, vedendosi distanziato, dovrà dubitare di essersi mosso verso L del tratto AB, nel qual caso sbaglia nel ritenere la velocità della luce BL/T=C perché nel tempo t, la luce ha percorso non solo il tratto BL ma anche il precedente AB. Che, se l'osservatore B non ha modo di assicurarsi del suo stato di moto, stia zitto per non ragionare da alchimista, nel quale caso invece di essersi condensato in Einstein il pensiero di un epoca, si sarebbe condensato nel Lambicco dell'osservatore, un semplice equivoco.
Concludendo: non è possibile rappresentare la relatività fisicamente come sarebbe il precedente tentativo e tanti altri; perché la relatività è teoria assolutamente matematica; è di importanza grandissima, ma le contrazioni di Lorentz essendo matematiche non possono in modo diretto venire dimostrate fisicamente.
Quando nella relatività si afferma esistere una dipendenza fra le distanze spaziali e temporali, e precisamente: che lo spazio e il tempo sono di tal natura da soddisfare alle equazioni di Lorentz; si afferma che nella relatività spazio e tempo dipendono l'uno dall'altro, e nello stesso tempo dal mondo reale. E sta bene: perché combinando il principio di relatività colla legge della propagazione della luce, si ottiene una legge matematica di trasformazione (trasformazione di Lorentz) delle coordinate rettangolari e del tempo, che soddisfa ai fenomeni naturali; cosicché: ogni legge del cosmo si trasforma in altra equivalente quando alle variabili x y z e t si sostituiscono nuove variabili x' y' z' e t', a condizione che la dipendenza matematica fra le prime e le seconde sia espressa dalla trasformazione di Lorentz. Ma rimane sempre aperto il problema alla ricerca matematica, poiché, non essendo le contrazioni dello spazio e del tempo fenomeni reali, si potrà tentare la non facile scoperta di un nuovo gruppo di trasformazioni, probabilmente assai più complicate, ciò che non lo renderebbe preferibile, ma egualmente soddisfacente alle leggi dei fenomeni naturali, e per il quale non si avrebbe come necessaria la condizione della relatività dello spazio e del tempo. Fino ad ora però ciò non si è trovato e solamente la relatività ha dato per soluzione una meccanica soddisfacente.

**118**
Gian Antonio Maggi a Gaetano Magnanini

Milano 31 Dicembre 1929

Carissimo Magnanini,

Effettivamente, secondo i principii della relatività, B trova la velocità di propagazione della luce, c, dividendo la lunghezza del tratto BL per il tempo che la luce impiega a descriverlo, lunghezza s' e tempo t', quali appartengono al riferimento con cui va unito, rispetto al riferimento con cui esso va unito, mobile, al modo indicato, con velocità di grandezza v, rispetto al riferimento con cui va unito A. Che se si indicano con s e t la lunghezza di tratto AL, e il tempo che la luce impiega a descriverlo, quali appartengono al riferimento con cui va unito A, si ha, per la trasformazione di Lorentz

$$(1) \qquad s' = \frac{s - vt}{\sqrt{\ }} \qquad t' = \frac{t - \frac{v}{c^2} s}{\sqrt{\ }}$$

Ne segue

$$s' - ct' = \sqrt{\frac{1 + \frac{v}{c}}{1 - \frac{v}{c}}} (s - ct).$$

Per cui si verificano, come reciproca conseguenza l'una dell'altra,

$$(2) \qquad s - ct = 0, \qquad s' - ct' = 0,$$

ossia si ha

$$(3) \qquad c = \frac{s}{t} = \frac{s'}{t'}.$$

Il Lammel sembra non accennare che a diverso apprezzamento del tempo. E, per verità, per avere (2) e (3) basta supporre

$s' = s - vt, \qquad t' = t - \frac{v}{c^2} s$

con che si ha

$$s' - ct' = \left(1 + \frac{v}{c}\right)(1 - ct).$$

Ma, per soddisfare tutte le occorrenti condizioni, sono necessarie le (1).
Per quanto al seguito accennato dal Lammel, penso ch'egli intenda che le (1) che conosce lui, conosce Magnanini, conosco io, sono un segreto per i non iniziati, fedeli a quell'ampia approssimazione che conta per esattezza nella quasi totalità dei problemi concreti.

Chiudo, con questo...

**119**
Gaetano Magnanini a Gian Antonio Maggi

*[Scritto in matita da Maggi:]* La risposta in Note.

Roma 5/1/30

Illustre e Car.$^{mo}$ Prof.

Essendo venuto da Napoli a quì, ho trovato la tua gent.$^{ma}$ replica che mi era stata respinta da Modena. Mi propongo di studiarla presto, riprendendo l'esame della relatività in cui sono entrato sia pure superficialmente da poco, perché la relatività è entrata nella Chimica. Ma quasi nulla ricordo del biennio di matematica superato quì in Roma sono già 45 anni passati! e sepolto prima da lustri di Laboratorio in operazioni scientifico-manuali di altra specie, poi da altri lustri di vita touristica *[sic!]* attraverso il mondo Europeo-Asiatico-Africano!.

Io mi gratto ora la testa (excuse my *[sic!]*) per vedere se posso capacitarmi di questa domanda che la mia testa fa a se medesima:

L'accorciamento delle lunghezze ed il ritardamento negli orologi sono reali come non credo, o sono solo apparenti nel senso che ti ho scritto. Se le cose stanno come mi scrivi, e devo accettare

$$c=\frac{s}{t}=\frac{s'}{t'}$$

non posso dare alla relatività una base reale cioè una serietà obbiettiva. Dove va allora il metodo induttivo Superba conquista delle nostre scienze? Si fa come il gambero.

Però in tentativi per ottenere equazioni relativistiche che hanno dimostrato di arrivare a spiegare fatti disparati, diversamente inspiegabili, si può permettere anche il passo del gambero.

In fondo anche le equazioni relativistiche hanno carattere empirico e come tali: tutto serve! Ma se si costruisce colla relatività una meccanica delle apparenze sacrificando l'induzione, io pretendo che si riconosca che quelle variazioni di s e di t sono bensì matematiche ma fisicamente apparenti.- Rinnovo gli auguri e saluti

Obbt. e Aff. Magnanini

**120**
Gian Antonio Maggi a Gaetano Magnanini

Milano 19 Gennajo 1930

Carissimo Magnanini,

Per rispondere alla gradita tua del 5 corrente, mi sembra che possa servire meglio di tutto un richiamo della questione che si può considerare come originaria.

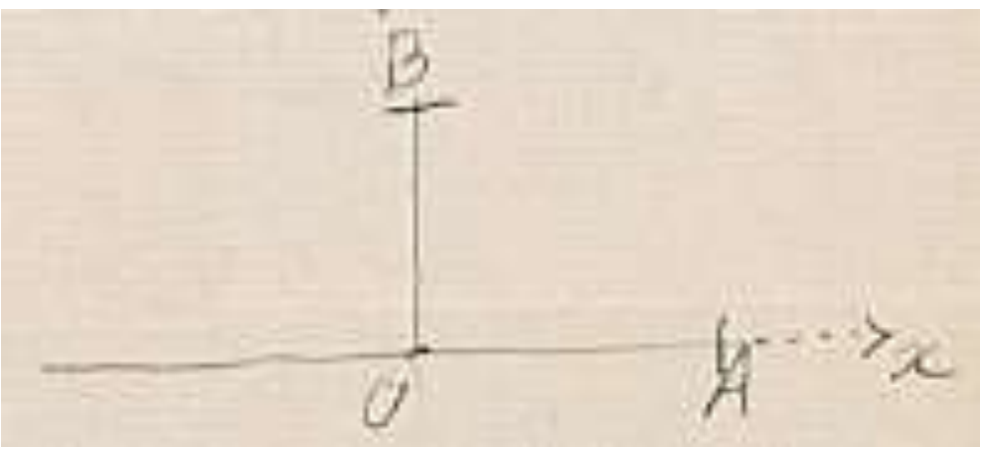

Un punto O, in cui si trova un osservatore, e due punti A e B, a cui sono applicati specchi piani, perpendicolari rispettivamente a OA e OB, appartengono ad una figura o spazio S', mobile di moto traslatorio, rettilineo, uniforme, con velocità di grandezza v, e nella direzione di OA, relativamente ad uno spazio S, nel quale si propaga la luce, con velocità di grandezza c, non perturbata dal suddetto movimento di S'. L'osservatore si propone di determinare la velocità di propagazione della luce nella direzione di OA e nella direzione di OB, giudicata dal riferimento S', dividendo il doppio della lunghezza del segmento OA o OB per il tempo impiegato da un segnale luminoso a partire da O, riflettersi in A o B, e tornare in O. Egli ha a disposizione l'esperienza e il calcolo, fondato sulla propagazione della luce, come avviene riferita a S.
Cominciamo col calcolo. Indicando con $O_1$,O,$O_2$ le posizioni occupate da O in S, ai momenti della partenza del segnale, del suo arrivo in A o in B e del suo ritorno, con t il tempo impiegato dal segnale a portarsi da $O_1$ in A o in B e da A o B in $O_2$ e con l la lunghezza di OA o OB, il calcolo fornisce agevolmente, nel primo caso,

(1) $$t=\frac{2l}{c}\frac{1}{1-\left(\frac{v}{c}\right)^2}, \qquad \frac{2l}{t}=c\left(1-\left(\frac{v}{c}\right)^2\right),$$

e, nel secondo caso,

(1) $$t=\frac{2l}{c}\frac{1}{1-\left(\frac{v}{c}\right)^2}, \qquad \frac{2l}{t}=c\sqrt{1-\left(\frac{v}{c}\right)^2}.$$

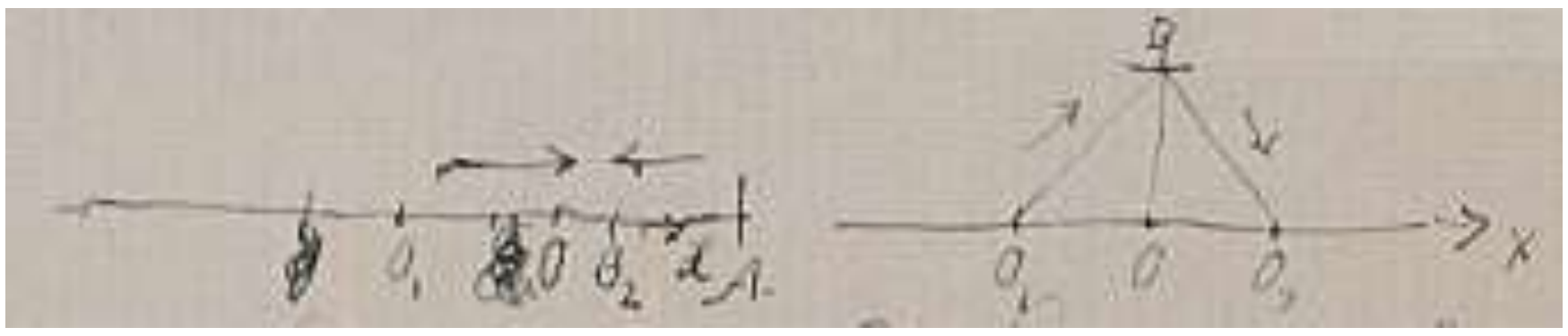

Ora, indicando con l', t' la lunghezza di OA o OB, e il tempo impiegato dal segnale a portarsi da O in A o B, e da A o B di nuovo in O, col riferimento a S' - cioè come risultano all'osservatore dall'esperienza eseguita nello spazio S', con riferimento a questo spazio - si dà per un risultato dell'esperienza che, in ogni caso, si trova

$$\frac{2l'}{t'}=c.$$

E questo risultato, da reputarsi fornito dall'esperienza, si riconduce al precedente risultato del calcolo, facendo, in (1) e (2),

$$t=\frac{t'}{\sqrt{1-\left(\frac{v}{c}\right)^2}},$$

in (1)

(α) $$l=l'\sqrt{1-\left(\frac{v}{c}\right)^2},$$

in (2)

$$l=l',$$

cioè postulando che il tempo impiegato dalla luce a compiere il [suo] tragitto abbia un valore diverso, secondo che si riferisce a S o a S', e la lunghezza di un segmento abbia lo stesso valore, se la direzione del segmento è perpendicolare a quella del movimento, e se la direzione è la stessa, un valore diverso a seconda dei due casi. Nessuno dei due valori appartiene al tempo o alla lunghezza a maggior ragione dell'altro. È postulato il fatto che il valore del tempo e della lunghezza non ha un significato intrinseco, ma è subordinato al riferimento.

Questi sono i concetti della teoria einsteiniana. Lorentz, come ben sai, aveva supposto, sul principio, un effetto contrazione, conformemente alla (α), che scaturisce dalla sua trasformazione, dando il nome di contrazione lorentziana, che si suol conservare anche nella Relatività.

Scrivimi, se altro desideri, ché mi farai sempre piacere e intanto ricevi una cordiale stretta di mano da

Tuo aff. collega G.A. Maggi.

**121**

**Manzoli** a Gian Antonio Maggi

[carta intestata: Casa veterani e invalidi Umberto 1°

Milano - 10 Via Giosuè Carducci - 10 - Il Presidente]

Milano 28 Febb

Carissimo amico,

Il mio Segretario e Cassiere di questa Casa Prof.r Arturo Finzi mi ha fatto dono del qui unito opuscolo nel quale risolve un problema anzi vari problemi che dal poco che io ancora ricordo mi sembra di grande importanza nella soluzione - mai finora cercata - di vari problemi di applicazioni di calcolo pratico.

A lei mi permetto di sottoporlo - con carico di restituzione - per la grande importanza che mi sembra abbia ed anche per conoscerne il di lei parere in merito.

Qualora poi ritenesse che il ripetuto studio come metodo fosse degno di un elogio suo e di codesta Università le sarò grato se volesse - per mio mezzo - esternarlo al Prof.r Finzi.
Qui nulla di nuovo che menomamente la possa interessare. Del resto io e mia moglie facciamo una vita tutta in casa e le mie occupazioni - ed anche il mio cuore è dedito tutto ai miei poveri Veterani nella ricerca di denari affinché almeno il pane non manchi loro.
Dirà lei: Ma il Ministero? Niente o poco ha nelle casse. Né io so più dove stendere la mano per trovarne.

E questi patriottici avanzi rappresentano i pochi avanzi rimasti delle nostre prime battaglie.
Così è carissimo amico. Si abbia coi miei i saluti di mia moglie.

Affmo [T.G.G.] Manzoli
Via Circo 14

**122**
Gian Antonio Maggi a **Roberto Marcolongo**[154]

A Marcolongo - 30 Novembre 1917
E ora lasciami discorrere d'altro. Ricorderai una breve Nota del Burgatti, sul principio d'Archimede nei solidi, negli Atti dell'Accademia di Bologna.[155] Io non la posseggo per cui credo di non averla ricevuta. Ma veduta, per proprio conto, dai miei colleghi tecnici, questi rilevarono la discordanza tra la conclusione che i corpi estranei devono essere espulsi dai mezzi elastici, e i fatti più accertati; e, d'accordo, crediamo di aver riconosciuto il difetto del precedente ragionamento. Il qual difetto mi sembra consistere in ciò che il Burgatti, immaginato in un mezzo elastico in equilibrio un contorno, che ne limita una regione, afferma che, concependo sostenuta questa regione del mezzo con un corpo estraneo, si manterranno *[...]*[156] riluttante e il *[...]* applicate ai *[...]* del contorno. Questo è vero per un fluido, e ne segue che il corpo estraneo, ove sia libero di muoversi, non potrà mantenersi in equilibrio. Ma, se non erro, non è egualmente vero per un mezzo elastico, nel qual caso il risultante e il momento delle pressioni esercitate sopra un contorno, che s'immagini tracciato nell'interno, per effetto della parte circostante, non sono più inerenti al contorno; e, concepito il contorno riempito da un corpo estraneo, si potrà sempre immaginare stabilito un sistema di pressioni che mantenga il corpo, comunque mobile, in equilibrio. Siccome le conclusioni del Burgatti, accettate senza discussione, hanno fatto capolino in scritti di materia d'applicazione, avevo preso un certo impegno di far presente al Burgatti stesso le nostre obbiezioni. Con questa occasione di scriverti, le sottopongo al tuo giudizio. Che se trovi che abbiamo torto *[...]* mostrare il nostro *[...]* dividi la nostra *[...]* nessuno meglio di te potrà scrivere al Burgatti per richiamare sull'argomento la sua attenzione. Ad ogni modo ti prego di volermi rispondere in proposito una parola.

**123**
Gian Antonio Maggi a Roberto Marcolongo

Pisa 15 Febbrajo 1918.
Ti ringrazio della trascritta risposta del Burgatti, la quale, d'accordo, non è risposta che per modo di dire. E che cosa possiamo intender noi che tu non intenda? Il fatto sta che difficilmente si attribuirebbe a quelle parole un significato qual che si sia. Il collega sembra essersi scordato della differenza dei risultati, nel posare una pietra sopra uno specchio d'acqua o sopra un suolo elastico. Ché, nel primo caso, la pietra fa un buco nell'acqua, e, nel secondo caso, il suolo cede, e acquista una nuova posizione d'equilibrio, in relazione col peso della pietra. "Muterà poi..." egli dice, con quel che segue. Ma io faccio seguire invece: muterà la distribuzione dello spostamento, e precisamente, muteranno gli spostamenti alla superficie $\sigma$, come occorre perché, col nuovo assetto, il corpo resti in equilibrio nel mezzo considerato.
Al prof. Marcolongo V. Sua lettera in data di Napoli, 12 II.[157]

---

[154] La grafia, tranne che nell'intestazione, non sembra quella di Maggi e la lettera è in pessime condizioni di conservazione.
[155] Potrebbe trattarsi di: P. Burgatti, "Il principio d'Archimede nei mezzi solidi", Nota letta alla R. Accademia delle Scienze dell'Istituto di Bologna nella sessione del 21 maggio 1916.
[156] Dove manca il testo, la carta è strappata.
[157] Non è presente nel Fondo.

**124**

Roberto Marcolongo a Gian Antonio Maggi
[busta indirizzata a: Illustre Prof. G.A. Maggi della R[a] Università -
Via Plebisciti 3 - Milano (121);
scritto da Maggi sulla busta: *Contiene la risposta in data 19, XI, 1930.*;
busta e carta intestate: R. Università di Napoli - Istituto di Meccanica Razionale]

12.XI.930.

Caro Maggi.

La corrente continua ed ad alto potenziale di bestialità degli esami che si scarica su di me da che sono arrivato, cioè da sedici giorni, non mi ha consentito di rispondere alla tua letterina, né di ringraziarti del dono delle tue lezioni. Dono squisito e graditissimo come sempre e come mi sono sempre state le tue cose da cui ho tanto appreso e avrò sempre da apprendere.

Ho cominciato a leggerle nei brevi, troppo brevi momenti di libertà e spero finirle di studiare appena ridiventerò padrone del mio tempo. Fermandomi a leggere il capitolo a pag. 101 "Deriv. della funz. pot. di un doppio strato" mi sono ricordato delle mie lezioni sullo stesso argomento fatte qui vari anni or sono e delle memorie famose di Liapunoff e di Ernest Neumann. Se non ho mal compreso, tu non ti sei preoccupato (e in un corso che non ha solamente per oggetto la funz. pot. ciò è perfettamente legittimo) della grossa questione della esistenza delle derivate normali e delle altre, per la quale non basta la continuità del momento e la sua derivabilità. Non so se tu conosca (ma lo conoscerai certamente) questo esempio.

Su di un cerchio di raggio a è deposto uno strato doppio la cui densità vari per cerchi concentrici proporzionalmente al raggio ρ; talché sia h=kρ. Elevo in O la normale z al piano del cerchio e su questa assumo un punto P distante da O di ζ>0.

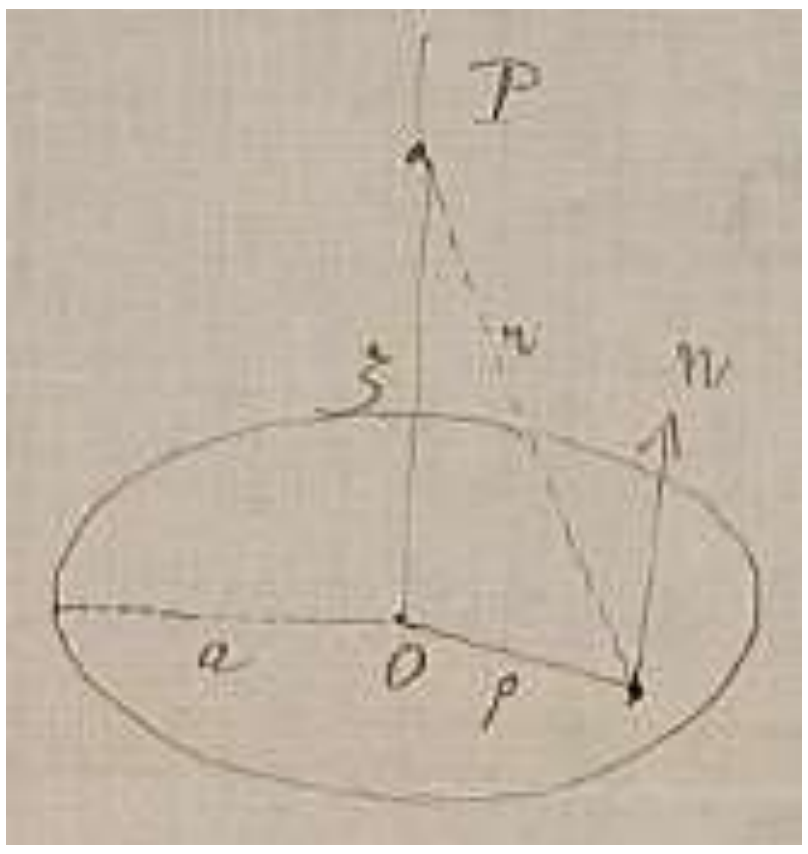

La funzione potenziale in P del doppio strato è data da

$$V_p = k\int_\sigma \frac{\zeta}{r^3}\,d\sigma \qquad (d\sigma = \text{den. d'area})$$

$$= 3k\pi\zeta\int_0^a \frac{\rho^3 d\rho}{(\rho^3+\zeta^2)^{3/2}}$$

ossia

$$V_p = 3\pi k\zeta\left\{\log\left(a+\sqrt{a^3+\zeta^2}\right) - \log\zeta - \frac{a\zeta^2}{\sqrt{a^3+\zeta^2}}\right\}$$

$V_p$ è finita anche per ζ=0. Ma la derivata rispetto a ζ (deriv. normale) per ζ=0 è infinita sia per ζ>0 che per ζ<0. Dunque in O non esistono le derivate normali né da una parte, né dall'altra, pur essendo k la plus honnête fonction du monde!

Questo anno io torno ad esporre la meccanica delle matrici oramai messa in pieno accordo con le omografie. Ho al proposito (come ti potrà dire il Cisotti) estesa una grossa memoria di cui presero visione Burgatti, Vataghin, Burali: E ciò costituì un'altra delle mie molte dispute con Burali. Siccome non ho più tempo di litigare con nessuno (mi bastano Picone, il gran Mauro, ed il Signorini) così ho messo a dormire ogni cosa.

Avrai già veduto Vittorio;[158] ed io ti avrò superlativamente annoiato! Saluti agli amici ed un abbraccio dal tuo centenario R. Marcolongo.

**125**

Gian Antonio Maggi a Roberto Marcolongo

Milano 19 Novembre 1930

Carissimo Marcolongo,

Ti ringrazio molto della tua cara lettera, che ho trovato jeri sera sul mio banco, tornando da casa Vivanti, dove mi ero congedato dal nostro Martinetti, coll'incarico di recarti i miei migliori saluti. Ti ringrazio poi in particolare dell'amichevole gradimento del mio volumetto. Per quanto alle tue osservazioni a proposito delle derivate della funzione potenziale di doppio strato, ebbi per le mani la famosa memoria del Ljapunoff, pesante, a mio modo di sentire... il peso, come, in generale, la produzione di quel valentuomo, ben meritevole di miglior sorte. Ma non vedo bene, dal tuo scritto, in quanto la grossa questione dell'esistenza ecc. concerne la mia deduzione della formola (3) a pag. 104, e delle sue conseguenze. Né mi sembra che vi rechi pregiudizio l'esempio che tu mi citi. Perché, io dichiaro esplicitamente di intendere che la h ammette derivate prime limitate e continue, per stabilire la suddetta (3), e derivate seconde idem, per le dedotte conseguenze. Ora la $h=k\rho$ dell'esempio, essendo $\rho=\sqrt{a^2+b^2}$, $\frac{\delta\rho}{\delta a}=\frac{a}{\sqrt{a^2+b^2}}$, $\frac{\delta\rho}{\delta b}=\frac{b}{\sqrt{a^2+b^2}}$, nel punto (a=0, b=0), cioè $\zeta=0$, dove si verifica l'accennata singolarità, non ammette più derivate di valore determinato. Per modo che io la chiamerei la plus honnête fonction du monde, a patto di chiamare parimente honnête, con relativo contorno, una ragazza, che, pur una sola volta, abbia perso un *[*ferro*]*.

Rinnovo con questo a te e ai Tuoi i più cordiali saluti, e con abbraccio mi confermo

Tuo aff^mo^

Nino Maggi.

[158] Martinetti; si veda la lettera successiva.

**126**
Gian Antonio Maggi a **Luigi Medici**

Milano 3 Gennaio 1937

Carissimo Gino,

Ho letto, da cima a fondo, il tuo volumetto, donatomi con così care parole di dedica, e letta per inteso con sincero interesse e piacere.[159] E con questo intendo anche di dire che, per essere libro di memorie personali, ma attinenti a materia originale, punto standardizzato alla '900 - cose e persone, e, tra queste, in prima linea, l'autore - lo reputo destinato ad essere letto con interesse e piacere, anche da chi non sia così legato, per proprio conto, con quella memoria da vedermi da te il dolce nome di zio Antonio. Rimando al discorso voce alcune osservazioni, che ti proveranno, se non altro, l'attenzione che ho posto nella lettura.

In attesa quindi di trovarmi con te, per riprendere a voce il discorso, mando a te e alla Felice,[160] desiderandovi noi tutto, le più affettuose cose, e, con un abbraccio, ti raccomando di voler sempre lo stesso bene al

Tuo Zio Antonio.

All'Ill. Sig. Avv. Prof. Luigi Medici
Via Moscova Città

**127**
Gian Antonio Maggi a Luigi Medici

P.M.
A Gino Medici
Leggendo il suo M.S della *[*selezione*]* del Porta della Meneghina

Parrebbe dal M.S. opinione del Cherubini che itt sia plurale di in. Non so che cosa ne pensano Salvioni e C. Per me è palese errore, di cui non è difficile risalire alla causa.

itt è plurale di ett. Pitt plur. di pett nel Porta, e nella frase andà coi pitt per ari: sonitt plur di sonett, nel Porta: boffitt plur. di boffett, nella frase cribi e boffitt! Però, nei nomi che non contengono un diminutivo, questa forma è ormai caduta in disuso. Si conserva invece non come i suaccennati: fiolitt plur. di fiolett: angiolitt plur. di angiolett: diavolitt plur. di diavolett: poveritt plur. di poverett. Ma, in altri diminutivi, e qua sta il punto!, pel plurale in itt non si ha il singolare in ett: passaritt ma passarin: barabitt ma barabin: Martinitt ma Martinin: cagnitt ma cagnin ecc.

---

[159] Potrebbe trattarsi di: L. Medici, *Una famiglia dell'Ottocento lombardo: anime, memorie e cose quasi romantiche*, Officina grafica L. Sala & Figli, Milano, 1936-XV. A p. 134 si legge: "Se Francesco Casorati fu celebre nell'arte medica - di lui si conservano studi importanti sulle febbri intermittenti, ecc. - il figlio, Felice, toccò invece il vertice della gloria nel campo delle matematiche superiori. Pavia ne serba le sembianze in un marmo che pure oggi, nel mio pellegrinaggio di amore, rividi sotto i portici dell'Ateneo; e l'Istituto Tecnico, a lui intitolato, ne tramanda alle scolaresche il nome. Opere celebri rimangono a perpetuarne la fama di chiaro docente di calcolo sublime, come allora si diceva, e di geodesia; e i suoi discepoli - pochi in vero superstiti - dalle cattedre ideali più illustri delle università italiane, ne ravvivano il ricordo. E basta citare qui, per tutti, G. Antonio Maggi, anima generosa di ambrosiano autentico, poliglotta, scrittore d'alto volo, scienziato di fama mondiale."

[160] Felicita Maria d'Incisa di Camerana.

Visin, cassin, spin è giusto che facciano al plurale in in, perché, mentre si può pensare a un passerott, barabatt, cagnett, soltanto non usati, come sono usati invece, fiolett, angiolett, diavolett, niente di simile si può immaginare per tali nomi, di tutt'altra fattura. In questi casi il singolare si fa in in. E sebbene si possa fare in in anche il plurale (duu angiolitt o duu angiolin), ma essendo pur vero che, nei precedenti esempi, si può sostituire il singolare in in al singolare in ett (on fiolin, on fiolett) è derivata di qui l'opinione che itt sia il plurale di in.

**128**

Gian Antonio Maggi a **Giacinto Morera**

Pisa 16 Giugno 1908
(A Morera)

Carissimo amico,

Ti mando la traduzione del Saggio del sig. Saltikov. E non ti prenda, per caso, scrupolo d'avermi procurato troppa fatica, perché questa si è ridotta allo scrivere, e, se dovrai qualche volta lottare con la mia zampa, ti convincerai facilmente che ho scritto currenti calamo. Ho preferito questo a far un sunto, perché non avendo in questa materia una particolare competenza, difficilmente avrei potuto farlo in modo da riuscire a te competentissimo di qualche servizio. Intanto ho seguito, per mio proprio conto, coll'interesse di cui mi sembra meritevole, l'opera dell'egregio Saltikov. Il quale, in questi giorni, mi mandò da Kiev, la collezione delle sue memorie, per cui, pochi giorni dopo, potevo dispensarti di mandarmi l'originale. Ti rimando questo, insieme coll'altro opuscolo, che mi riserbo di leggere nella copia mandatami dall'autore.

Non abbiamo ancora perfettamente fissato il posto della nostra villeggiatura, ma spero che non ci mancherà l'occasione di ritrovarsi nel corso delle vacanze. Chiudo per ora, mandandoti i nostri migliori saluti, e pregandoti di ricordarmi ai colleghi di costì, specialmente al Somigliana. In particolare ricevi una cordiale stretta di mano del

Tuo aff mo
G.A. Maggi.

## 129
Gian Antonio Maggi a **Raffaello Nasini**

Casa 29 Giugno 1907

Egregio collega,

Le rendo l'opuscolo, con un po'di traduzione, un po' di resoconto, e, le altro desidera disponga di me, sicuro di farmi sempre piacere.

Non mi pare che regga l'applicazione che l'autore fa del teorema della minima azione ai processi chimici. Questo teorema confronta il movimento effettivo coi movimenti virtuali (possibili cogli stessi vincoli, nel caso dei sistemi olonomi) divisi, che hanno comune con esso le posizioni estreme del sistema, e la forza viva su di esse, e soddisfanno alla stessa equazione della conservazione dell'energia (con che la forza viva sarà la stessa anche nell'altro estremo). Il teorema afferma che al movimento effettivo, almeno per un intervallo di tempo abbastanza ristretto, compete il minimo valore di $\int_{t'}^{t''} Tdt$ (T= forzaviva, t=tempo). Ora questo integrale non rappresenta l'importo di forza viva, o di energia cinetica, corrispondente al movimento considerato, né si potrebbe surrogarvi l'importo medio $\frac{\int_{t'}^{t''} Tdt}{t''-t'}$, perché l'intervallo (t't") non riesce generalmente lo stesso pei varii movimenti considerati - circostanza essenziale, di cui l'autore non tiene conto.

Non saprei poi (Ella sarà ben miglior giudice di me) come dal minimo importo dell'energia cinetica l'autore inferisce il massimo svolgimento di calore. Io ho un forte sospetto ch'egli v'inserisca come anello di congiunzione il massimo importo di energia potenziale - cioè di energia di affinità chimica - e misuri questa col calore svolto dalla combinazione. Ciò sta bene, ed è quello che fanno Thomsen e Berthelot, nei passi riportati dall'autore, ma corrisponde ad un punto di vista diverso dal suo, simile a quello di chi inferisce l'altezza della caduta di un grave del calore svolto dall'urto del medesimo contro il suolo. Dal punto di vista dell'autore, bisognerebbe, in questo caso, considerare come importo dell'energia cinetica la produzione di forza viva, proporzionale al calore svolto dall'urto. Ma l'esempio non si presta, con egual efficacia, all'inversione, perché qui si tratta di un sistema libero, e, nel caso dell'autore, ci stanno i vincoli della temperatura a delle pressioni assegnate.

Mille ossequi alla gentile Signora. A Lei i migliori saluti dall'

Aff$^{mo}$ suo
G.A. Maggi.

Al Ch$^{mo}$ Prof R. Nasini
S.M.

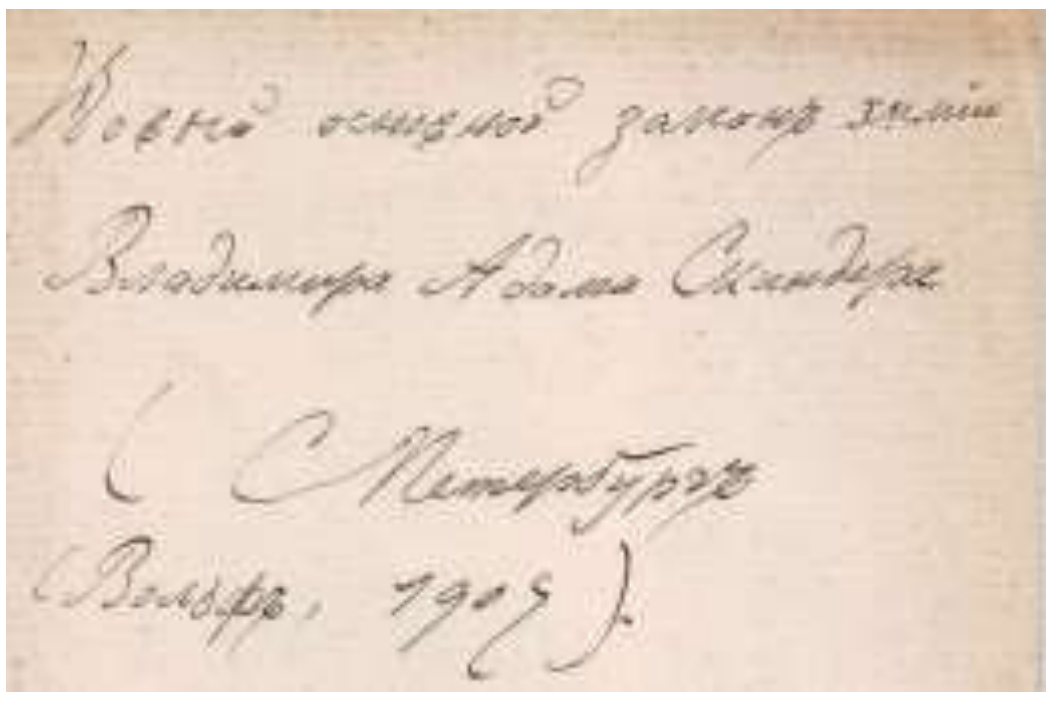

Opuscolo mandatomi in esame da Nasini da parte di Cannizzaro.[161]

[161] Potrebbe trattarsi di Stanislao, maestro di Nasini.

**130**
**Alpinolo Natucci** a Gian Antonio Maggi
[cartolina postale indirizzata: Al Chiarissimo Signor Prof Gian Antonio Maggi
Corso Plebiscito n° 3 - Milano]

Chiavari 28 dicembre 1932

Chiarissimo Sig. Professore,

Le scrivo, innanzi tutto per presentarle fervidi auguri per il nuovo anno. In secondo luogo vorrei comunicarle una mia idea riguardo alla forza centrifuga. Essa si riconnette al Suo scritto in: "Selecta" pag 93 che ho letto con molto interesse. In esso ella conclude che "la forza centrifuga di un punto mobile è una forza che il punto esercita sul materiale, e non una forza applicata allo stesso punto", e di questo io sono pienamente convinto. Poi Ella aggiunge: vi sono casi in cui si chiama… f. centrifuga… una forza… applicata allo stesso punto. E cita la nota esperienza della sferetta che si allontana sull'asse di rotazione. Ora a me, riflettendo su questo caso, è sembrato che anche ora, come in altri casi; la forza centrifuga possa farsi dipendere colla inerzia del corpo. Il corpo tenderebbe a muoversi in linea retta per inerzia, ma è costretto a muoversi in circolo, dunque esso tende ad accrescere il raggio di questo circolo, fino all'infinito fino, cioè, a tramutarlo in una retta. La sferetta si allontana dunque più che può.

Non so se Ella sa che sto preparando un libro di fisica per le scuole medie, che deve essere edito dalla Casa editrice E.S.T.[162] Ne ho già consegnato i manoscritti di 2 volumi (Meccanica-Acustica-Termologia- Ottica), e nei primi mesi del 1933 preparerò il 3°.

Avrebbe difficoltà a premettervi una breve prefazione? Mi permetto fare questa domanda, perché la Casa ne sarebbe molto contenta ed io pure.

Con ossequi, ed affettuosi saluti

Dev$^{mo}$ Suo A. Natucci

[162] *Elementi di Fisica moderna*.

**131**
Gian Antonio Maggi a Alpinolo Natucci
[busta indirizzata a: Prof. Alpinolo Natucci - Chiavari.]

Milano 29 Dicembre 1932

Caro Professor Natucci,

La ringrazio, prima di tutto, de' suoi buoni augurii, e cordialmente li ricambio.

Su quanto mi scrive della forza centrifuga, il secondo significato, di forza applicata al punto, compare coll'equilibrio e col movimento del punto, relativo ad assi in movimento rotatorio uniforme. Esempio più alla mano è quello dell'equilibrio, col quale si ha la parte della forza centrifuga nella gravità. Io ho voluto recare un esempio di movimento. Ora pur ammesso, in questo movimento relativo, una tendenza della sferetta a scostarsi dal centro, sul raggio, per virtù d'inerzia, essendo però puramente geometrico, e non materiale, il circolo che il centro della sferetta è vincolato a descrivere, la sferetta non preme, in direzione del raggio, sopra alcun corpo, e non si può, quindi parlare del primo significato. Le pressioni della sferetta si esercitano perpendicolarmente alla asticina (supposto trascurabile l'attrito), e quindi perpendicolarmente al raggio.

È poi da attribuirsi la così detta disaggregazione delle parti di un corpo rapidamente ruotante, per forza centrifuga, all'inerzia, in virtù della quale le particelle sfuggono, <u>secondo la tangente</u> alla orbita circolare fino a quel momento descritta, quando la coesione del corpo non è più capace di fornire la forza <u>centripeta</u>, necessaria al mantenimento della compagine. Ma è questa una questione diversa dalla precedente.

Per la proposta che mi fa di una prefazione al Trattato di Fisica, di cui sento con piacere che attende alla pubblicazione, Le sono grato dei sentimenti di cui la riconosco ispirato, ma sono spiacente di non poter soddisfare il suo desiderio. Perché la mia competenza in materia non potrebbe oltrepassare il capitolo sulla Meccanica, e anche per questo non mi sentirei di mettermi in concorrenza con altri che abbiano assai maggior esperienza dell'insegnamento medio.

Di nuovo aggradisca i miei augurii, coi quali e coi migliori saluti, mi confermo

Aff suo

**132**
**Pietro Pagnini** a Gian Antonio Maggi
[stampato con correzioni]

Con corr[illegible]

# Nota sull' esperienza del Michelson

Nella nota disposizione interferenziale adoperata per la verifica del trascinamento terrestre di velocità v e qui rappresentata schematicamente, il Michelson per i due percorsi $OR_1O'''$ e $OR_2O''$ prende i tempi

$$t_1 = t_1' + t_1'' = \frac{1}{c+v} + \frac{1}{c-v} = \frac{21}{c}\left(1 + \frac{v^2}{c^2}\right)$$

e

$$t_2 = \frac{2s}{c} = \frac{21}{c}\left(1 + \frac{1}{2}\frac{v^2}{c^2}\right) \qquad (\alpha)$$

e quindi la differenza

$$\Delta t = t_1 - t_2 = 1\frac{v^2}{c^3} \ (1).$$

Se calcoliamo direttamente ~~i cammini ottici~~ lo spostamento $O''O'''$, ponendo $d$ = spostamento di $R_1$ quando la perturbazione va incontro da O a $R_1$, $d'$ = spostamento di $O'$ al ritorno da $R_1$ in $O'''$, e $p$ lo spostamento uguale all'andata ed al ritorno per $OR_2$ e $R_2O''$, avremo

$$d = \frac{vl}{c+v}, \quad d' = \frac{vl}{c-v}, \quad p = \frac{vl}{\sqrt{c^2 - v^2}}$$

e quindi

$$\Delta s = O''O''' = d + d' - 2p = 1\frac{v^3}{c^3} \ (2).$$

$$= \frac{l v^2}{c^2}\cdot V$$

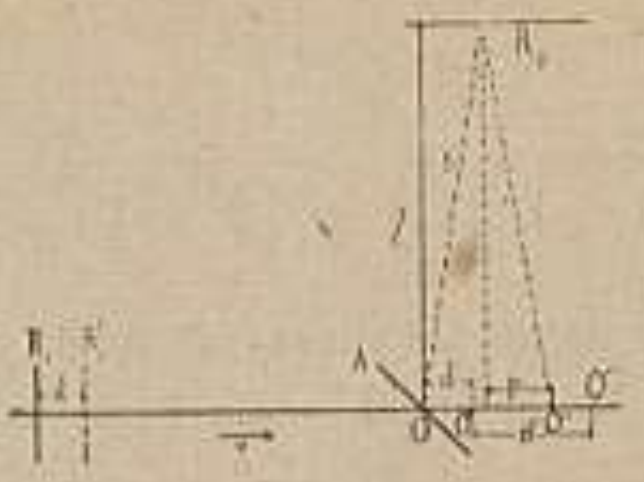

In conclusione da punti di partenza differenti e con considerazioni di cinematica elementare si giunge alla relazione

$$\Delta' s = v \Delta t \quad (3)$$

Il Michelson (*) invece ottenuto lo stesso $\Delta t$ della (1) lo moltiplica per $c$ velocità della luce per avere

$$\Delta s = c \Delta t \quad (4)$$

~~in palese contradizione colla (3) ottenuta direttamente~~, e per il numero $\Delta f$ delle frange spostate durante la rotazione di 90° dell' apparecchio dà in conseguenza

$$\Delta f = \frac{2c\Delta t}{\lambda} = \frac{2l}{\lambda} \frac{v^2}{c^2} \quad (5),$$

che corrisponde a $\Delta f = 0,4$ coi suoi dati sperimentali, mentre dalla (3) divisa per $\lambda$ abbiamo

$$\Delta' f = \frac{2\Delta s}{\lambda} = \frac{2v\Delta t}{\lambda} = \frac{2l}{\lambda} \frac{v^3}{c^3} \quad (6),$$

e quindi $\Delta f = 0,00004$, che rende pienamente ragione del resultato negativo dell' esperienza.

E siccome la (5) equivale all' altra formula spesso usata

$$\Delta f = \frac{2\Delta t}{T} = \frac{2l}{\lambda} \frac{v^2}{c^2} \quad (7)$$

ed all' altra

$$\Delta f = \Delta n = \frac{l}{\lambda\left(1 + \frac{v}{c}\right)} + \frac{l}{\lambda\left(1 - \frac{v}{c}\right)} - \frac{2l}{\lambda\sqrt{1 - \frac{v^2}{c^2}}} = \frac{2l}{\lambda} \frac{v^2}{c^2} \quad (8)$$

dove facciamo uso della lunghezza di onda nei due percorsi $l_1 = l_2 = l$, *questa coincidenza ci rivela l' errore* del procedimento del Michelson.

Infatti in tutti gli usi correnti della disposizione del Michelson *si tratta sempre di osservare le frange dovute ad uno spessore* $\Delta' s$ *di lamina d' aria compreso fra lo specchio* $R_2$ *e l' immagine coniugata* $R_1^{\circ}$ *di* $R_1$ *rispetto ad* O. Perciò potremo concludere che anche nel nostro caso lo spessore di questa lamina

(*) Phil. Mag. Vol. 24, 1887, pag. 453.

$\Delta' s$ ottenuto nel tempo $\Delta t$ colla velocità $v$ di trascinamento terrestre, essendo $\Delta' s = v \Delta t$ e non $\Delta s = c \Delta t$, *il numero delle frange corrispondenti è quello dato dalla* (6) *e* non quello corrispondente alle (5) (7) (8).

In conclusione le *formole* (4) (5) (7) (8) *dedotte dal calcolo del Michelson* ***non sono valide*** *per le frange localizzate in una lamina di spessore* $\Delta s$ *dovuta al trascinamento terrestre quando* $O R_1 = O R_2$ *ossia nelle condizioni che gli sperimentatori hanno comunemente adottato, ma invece* ***lo sono*** *per quelle che si osservano quando si faccia subire l'effetto di questo trascinamento ad una lamina di spessore* $O R_1 - O R_2 = D$ *ossia proprio nel caso che* $O R_1$ *non corrisponda alla condizione imposta dal Michelson di essere uguale a* $O R_2$.

Se dunque nell'ipotesi $O R_1 = O R_2$ vogliamo calcolare quale sia il tempo impiegato dalla perturbazione a percorrere $\Delta' s$ con la velocità $c$ dovremo porre la differenza invece della somma ($\alpha$) dei tempi di andata e ritorno per ciascun tratto e quindi per $O R_1 O'''$

$$\Delta t_1 = t_1'' - t_1' = \frac{2 l v}{c^2} \left(1 + \frac{v^2}{c^2}\right) \quad (9)$$

e per l'altro percorso $O R_2 O''$

$$\Delta t_2 = \frac{2 p}{c} = \frac{2 l v}{c^2} \left(1 + \frac{1}{2} \frac{v^2}{c^2}\right) \quad (10)$$

le quali daranno per il tempo cercato

$$\Delta' t = \frac{\Delta_s}{c} = \Delta t_1 - \Delta t_2 = l \frac{v^3}{c^4} \quad (11)$$

ossia la relazione $\qquad \Delta' s = c \Delta' t \quad (12)$

la quale confrontata colla (4) *ne mostra l'incompatibilità* perchè

$$\Delta' t = \Delta t \frac{v}{c}$$

Giunti a questo punto ci potremmo domandare se modificando convenientemente la disposizione sperimentale del Michelson si possa arrivare alla verifica del moto terrestre.

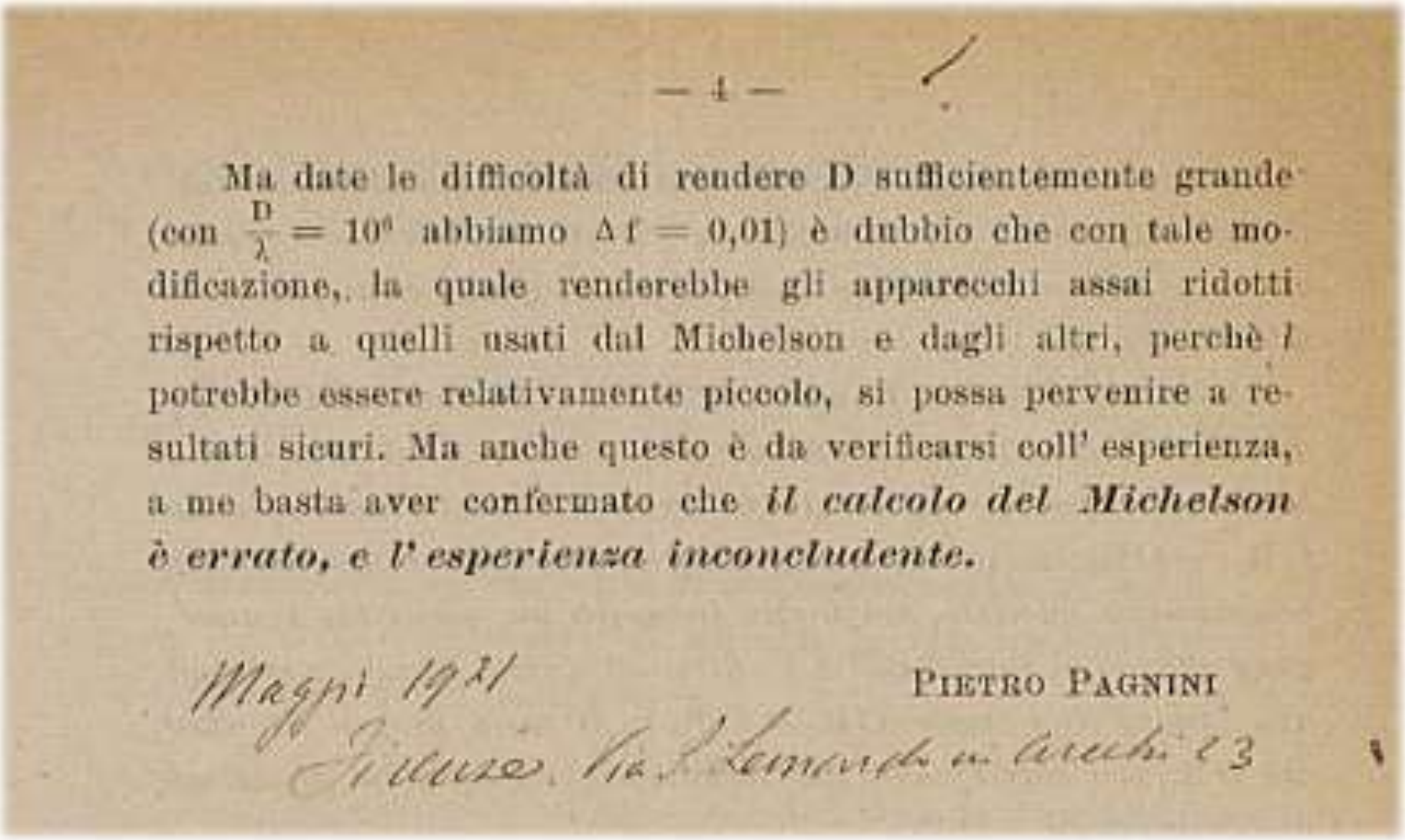

— 4 —

Ma date le difficoltà di rendere D sufficientemente grande (con $\frac{D}{\lambda} = 10^6$ abbiamo $\Delta f = 0{,}01$) è dubbio che con tale modificazione, la quale renderebbe gli apparecchi assai ridotti rispetto a quelli usati dal Michelson e dagli altri, perchè *l* potrebbe essere relativamente piccolo, si possa pervenire a resultati sicuri. Ma anche questo è da verificarsi coll' esperienza, a me basta aver confermato che ***il calcolo del Michelson è errato, e l' esperienza inconcludente.***

Maggio 1921 PIETRO PAGNINI

Firenze. Via S. Leonardo in Arcetri 23

**133**

Gian Antonio Maggi A Pietro Pagnini

Pisa 21 Giugno 1921

Egregio Signor Pagnini,

La ringrazio delle due copie della Sua Nota. Mi mancò il tempo, colpa varie occupazioni, non la buona volontà, di esaminare il manoscritto, ne convengo colle osservazioni indicate nella seconda copia. Ma, riconosciuto che il Michelson moltiplica $\Delta t$ per c, per avere il cammino percorso dalla luce nel tempo $\Delta t$, mi sembra poi fuor di dubbio che nei ragionamenti sulle interferenze, va introdotto, non $\Delta s' = v\Delta t$, che è il cammino percorso dal globo terrestre, ma $\Delta s = c\Delta t$; e questo è appunto ciò che fa il Michelson.

Mi spiace di non poter dir altro, nel provarLe che ho preso in considerazione i suoi gentili invii.

Intanto Le resto grato della memoria, e, con distinti saluti, mi confermo

Dev.$^{mo}$

GAM.

Al Sig Pietro Pagnini - Firenze.

**134**
Gian Antonio Maggi A Pietro Pagnini

Pisa 20 Ottobre 1922
(Sig. Pietro Pagnini
Via S. Leonardo in Arcetri 23 Firenze.)

Egregio Signor Pagnini,

Soltanto ora posso dirLe qualcosa sul Suo M.S., perché giunto durante una mia assenza, per una decina di giorni, da Pisa.

Il risultato a cui Ella giunge non manca intrinsecamente di interesse, ma non mi sembra facile attribuirgli una portata all'infuori del problema della barca o del nuotatore trascinato dalla corrente, così da valersene per la spiegazione dell'esperimento di Michelson. Occorre, a tal fine, precisare in quanto il movimento del Suo punto rispecchi la propagazione delle radiazioni luminose, in quell'esperimento.

A me non pare che l'analogia colla barca o col notatore che traversa un fiume fornisca la più preferibile dimostrazione dell'esperimento di Michelson. Un discorso ben altrimenti persuasivo mi sembra quello che si può leggere nel Weyl Raum-Zeit-Materie,[163] a pag. 142.

Mi è grato intanto presentarLe i miei migliori saluti, e confermarmi

Dev.mo Suo
G.A. Maggi.

**135**
**Alessandro Paoli** [a Gian Antonio Maggi]
[biglietto da visita: Prof. Alessandro Paoli -
s.l., ma probabilmente Pisa; s.d., ma antecedente al 1904]

Chiarissimo e illustre Collega,

Le mando la Memoria[164] del Sig D.r Quartaroli - con la dichiarazione del nostro collega Prof. Sestini. Il D.r Quartaroli ha introdotto le modificazioni suggerite dal Prof Maggi.

Pieno di ossequio sono Suo Devotissimo

A.P.

---

[163] H. Weyl, *Raum. Zeit. Materie. Vorlesung über allgemaine Relativitätstheorie*, Springer, Berlin, 1919.
[164] Non è presente nel Fondo.

**136**
**Roberto Parodi** a Gian Antonio Maggi
[busta indirizzata: All'Ill.$^{mo}$ Professore G.A. Maggi -
Docente di Meccanica Razionale - Milano;
busta e carta intestate: R. Istituto Nautico di Genova]

Genova 19 maggio 1926

Chiarissimo Professore,

Avevo inviato alla Redazione del Nuovo Cimento una mia breve nota, nella quale facevo seguito al Suo articolo sulla forza centrifuga,[165] dimostrando come la sua natura possa essere anche benissimo presentata elementarmente, senza le artificiosità del Prof. Bernini.

La Redazione del N.C., ignoro la ragione, non credé più opportuno continuare nell'argomento e mi restituì il manoscritto.

Non voglio discutere il provvedimento che a me tornò poco gradito, mi permetto però colla presente inviarLe in esame uno opuscolo nel quale sostengo le mie particolari vedute.

Voglia compiacermi di esaminarlo.

Con Ossequi

Prof. Parodi Roberto.

**137**
Gian Antonio Maggi a Roberto Parodi

Milano 29 Maggio 1929 - (Corso Plebisciti, 3).

Pregiatissimo Signor Professore,

Soltanto con parecchi giorni di ritardo ho ricevuto, attraverso Pisa e gli uffici dell'Università, la Sua lettera in data del 19 corr, poi l'opuscolo annunciatovi, e del cortese invio cordialmente La ringrazio.

I miei concetti di forza centrifuga sono quelli che espongo nell'articolo da Lei accennato, conformi a quelli del mio insegnamento di Meccanica Razionale, già tradotti nei miei libri sulla stessa materia. Affermando l'esistenza di una specie di forza centrifuga applicata dal mobile considerato ad un altro corpo, secondo il principio dell'eguaglianza dell'azione e della reazione, io non posso certo convenire che questa sia, senz'altro, asserzione manifestamente inesatta. E affermando poi l'esistenza di un'altra specie di forza centrifuga, applicata allo stesso corpo, la quale è forza apparente, collo stesso significato con cui è apparente il movimento o l'equilibrio relativo ad assi mobili, considerati come fissi, tanto che mostra come si possa eliminare colla sostituzione del riferimento fisso al mobile, non posso neppure convenire col ricondurre la forza centrifuga alla forza d'inerzia asserita come una forza vera e propria.

Io non vedo come dall'applicazione di due forze motrici eguali e contrarie ad uno stesso punto materiale possa seguire altro che un'accelerazione nulla, come non è l'accelerazione (centripeta) del punto in movimento circolare uniforme, nel riferimento da cui risulta questo movimento. Nulla diventa l'accelerazione del punto, riferendone la posizione ad una terna d'assi che possiede movimento rotatorio conforme al suddetto movimento circolare uniforme, ma con questo riferimento, scompare anche il movimento del punto. Considerando come reale il movimento circolare uniforme, sembra naturale di considerare come apparente l'equilibrio relativo agli assi mobili, considerati come fissi, e come apparente la forza (centrifuga, applicata al punto) colla quale si consegue l'equilibrio medesimo.

[165] G.A. Maggi, "Che cos'è la forza centrifuga?", *Il Nuovo Cimento*, 1926, pp. 21-30.

Così mi spiace di dover cominciare a dissentire sul concepire la forza d'inerzia - definita come eguale e contraria alla forza motrice applicata allo stesso mobile - come forza vera e propria. Una vis inertiæ, proporzionale alla massa, intesa come opposizione del corpo, quantum in se est, al mutamento di stato, è, per verità, introdotta da Newton (principia, Definitio III). Ma del discorso di Newton a me pare che non si possa inferirne altro che una preparazione alle vires impressæ, dalle quali il corpo cogitur statum... mutare (Lex Motus I). E perché la coazione abbia successo, mi sembra che bisogna pur intenda l'accennata opposizione altrimenti che traducibile con una forza, applicata allo stesso corpo, eguale e contraria, della quale in Newton non apparisce, ch'io sappia, cenno.

Le rinnovo intanto i miei ringraziamenti, coi quali e coi migliori saluti La prego di credermi

Dev.$^{mo}$ Suo G.A. Maggi.

## 138
## Gian Antonio Maggi a **Ernesto Pascal**

Pisa 25 Maggio 1910

Caro Collega,

Conformemente all'invito trasmessomi dal collega Nicoletti, ti mando, coll'accluso foglietto, il sommario del mio corso di Fisica Matematica: poche parole, perché di più non ne comporterebbe la qualità dell'argomento.

Siccome poi ho codesta buona occasione di scriverti, ti faccio anche il seguente discorso, del quale, per la ragione che pur dico, avevo oramai deposto il pensiero. Sollecitato dal prof. Saltykoff dell'Università di Kharkoff, mandai a quella Società Matematica, per essere inserito nei Rendiconti fin dall'estate dello scorso anno, una breve necrologia del nostro compianto Morera. La scrissi in lutto, e perciò avevo pensato, dopo che fosse pubblicata, di farne una traduzione, e offrirla al Giornale di Matematica, dove il Morera pubblicò il suo primo lavoro. A ciò mi moveva il desiderio di rendere leggibile quel mio piccolo omaggio, reso alla cara memoria di lui. Se non che il mio articolo fu presentato alla Società, al ricominciare delle sue adunanze, alla fine d'Ottobre (vecchio stile), ma soltanto ora mi viene annunciato che ne riceverò a giorni gli estratti. Questo lungo tempo trascorso era la ragione che mi aveva fatto abbandonare il proposito dell'offerta: e tu farai del mio discorso, affatto liberamente il conto che creda.[166]

...

Al prof. E. Pascal - Napoli.

[166] Il necrologio in questione venne pubblicato: "Giacinto Morera. Commemorazione (traduzione dal russo dell'Autore)", *Giornale di Matematiche di Battaglini,* v. XLVIII, pp. 317-324.

**139**
Gian Antonio Maggi a **Giuseppina Pastori**

Casa 4 Gennaio 1937

Gentilissima Signorina,

La ringrazio vivamente delle Sue gentili righe, e Le ricambio i migliori e più cordiali augurii. E poiché ne ho questa buona occasione, Le dirò anche che mi spiacque sinceramente di non sentir ricordato il suo nome da P Gemelli, nell'interessante sua conferenza al Seminario, pur di materia attinente ai bei volumi, che uniscono col nome di P. Gemelli il Suo, e di cui mi fece, a suo tempo, gratissimo dono. Come Ella vede, quel silenzio non fu ragione di dimenticare la Sua bella operosità.

Di nuovo, e da parte anche dai miei, mille auguri e cordiali saluti, coi quali La prego di credermi sempre

Aff Suo

Gent.[ma] Sig.[na] Giuseppina Pastori
Prof. della Università Cattolica del Sacro Cuore
Via Monte di Pietà 3 - Milano.

**140**
Giuseppina Pastori a Gian Antonio Maggi
[busta indirizzata a: Illustre Professore
Nob. Gian Antonio Maggi - Prof. Emerito della R[a] Università -
C. Plebisciti 3 - Milano; mittente: Pastori - Milano - V. M[te] Pietà 3]

Milano 6 genn. 1937 XV

Illustre Professore,

Grazie infinite del Suo biglietto. In un momento così arduo e triste, la Sua autorevole parola è per me preziosa.

Veramente, assai più che l'immeritato silenzio altrui sul mio lavoro passato, mi duole la forzata impossibilità di lavorare ancora. Ella sa ben valutare la passione del lavoro scientifico, e che cosa significhi la privazione di ogni mezzo di ricerca. L'infortunio inatteso non fu causato da me; spero dunque che la verità si faccia strada e che gli errori si riparino.

Con animo vivamente grato La prego accogliere l'espressione del più cordiale ossequio

Sua dev.[ma]
Giuseppina Pastori.

**141**
Gian Antonio Maggi a **Giuseppe Peano**

(Mandata al prof. G. Peano con preghiera di pubblicarla sulla Rivista di Matematica il giorno 8 Marzo 93)[167]

Chiarissimo collega,
...

Al quesito proposto dal mio egregio amico Sig F. Crotti (Rivista Matematica Vol II pag. 176) io rispondo che, se la formola di Simpson applicata al trapezio curvilineo compreso fra le ordinate $y_1$ e $y_{21}$ dà 9600 e 9720 applicata al trapezio compreso tra le ordinate $y_2$ e $y_{20}$, questo non è alcun assurdo che possa intaccare la bontà di quella formula. A seconda del caso, si sostituisce alla curva del trapezio una successione di 10 o di 9 archi di parabola coll'asse perpendicolare all'asse delle x, convessi i primi e concavi i secondi verso questo asse. E quel risultato significa puramente che la somma delle aree delle figure limitate dalla seconda successione e dalla parte della prima avente gli estremi comuni con essa, cioè 18 volte l'area della figura chiusa fra mezzo arco parabolico concavo e mezzo arco parabolico convesso, supera di 120 la somma delle aree dei due trapezii le cui curve sono le parti rimanenti della prima, pari all'area del trapezio la cui curva è un arco parabolico convesso. Che se qualcuno coi dati del quesito applica la formula di Simpson ad un determinato trapezio chiuso tra le ordinate $y_1$ e $y_{21}$ e alla sua parte chiusa fra le ordinate $y_2$ e $y_{20}$, vuol dire che nella valutazione del tutto e della parte s'accontenta d'una approssimazione che implica errori la cui differenza supera 120. Non vi è certamente limite all'errore che si può commettere applicando la formola di Simpson senza preoccuparsi che la successione degli archi parabolici non si scosti fuori di certi termini dalla curva del trapezio considerato. Piuttosto si può osservare che, coi dati del quesito, la curva e la parte di essa compresa fra le ordinate $y_2$ e $y_{20}$ non potranno adattarsi egualmente bene alla relativa successione d'archi parabolici.

Ma l'Ing. Crotti, come apparisce dalla sua Memoria "Sopra alcune formule di planimetria e stereometria" (Atti del collegio degli ingegneri ed Architetti di Milano Anno XVIII) muove dall'idea che la curva del trapezio considerato sia incognita o non si voglia conoscere: e discorre come se il valore fornito dalla formula di Simpson fosse un valore più plausibile dell'area cercata. Io non vedo veramente come un tal problema si possa porre; ad ogni modo non è quello che si propone di risolvere la formola di Simpson.

Questa essendo la mia opinione, conforme a quella esposta nella Rivista *[qui c'è un rimando, ma ne manca il testo]* dai prof. Jadanza e Bardelli, e suppongo anche alla Sua,[168] quantunque Ella si limiti a rispondere ad alcune osservazioni del prof. Bardelli; ne' mai avendo potuto *[*discorrere in tema diverso*]*, ho ben ragione di meravigliarmi nel trovarmi citato nella Memoria ora uscita dell'Ing. Crotti "Alcune considerazioni sulla recente edizione (12[a]) del Manuale del prof. Colombo" (Atti del Collegio degli Ingegneri ecc. Vol. XXV) come venuto alla conclusione che alla formola di Simpson manca ogni fondamento teorico e logico, o qualcosa d'analogo. Son posto, è vero, in ottima compagnia, tale da farmi vivamente desiderare, se quelli sono realmente di diverso parere (è nominato anche Lei), di

[167] Si confrontino le note relative alla corrispondenza con Crotti (in particolare alla lettera #50).
[168] Peano scrive con Jadanza l'articolo "Una difesa delle formule di Simpson, ed alcune formule di quadratura poco note", *Il Politecnico*, 1893; si dà notizia di questo articolo anche sulla *Rivista di Matematica edita da G. Peano*, 1993, v. III, p. 137.

conoscere le loro ragioni. Ma amicus Plato sed magis amica veritas; e perciò non mi si farà appunto di tenerci che chi notasse il mio nome non mi attribuisca un giudizio al mio assolutamente opposto.
...
Messina 6 Marzo 1893

Dev.mo Suo
G.A Maggi

**142**
Giuseppe Peano a Gian Antonio Maggi
[cartolina postale indirizzata: Al Ch.mo prof. Gian Antonio Maggi - Università di Messina]

Torino 14 marzo 93.

Egregio collega,

La sua lettera sarà pubblicata nel prossimo fascicolo della rivista.[169] Anch'io mi vedo menzionato dal sig. Crotti, che mi attribuisce un'opinione molto diversa dalla vera. La mia opinione sulla formula di Simpson trovasi nel mio libro Applicazioni geometriche del calcolo infinitesimale, e nelle mie lezioni di Analisi infinitesimale, che sto pubblicando, e di cui a giorni Ella riceverà in omaggio copia del primo volume.

Mi creda Suo devotissimo.

G. Peano.

**143**
**Luisa Pelosi** a Gian Antonio Maggi
[busta indirizzata a: Chiar.mo Prof. Gian Antonio Maggi
della R. Università - Corso Plebiscito 3 - Milano (21)]

Torino 12 12 1925

Chiarissimo Professore,

La ringrazio vivamente per il gentile invio della Sua memoria che Le rimanderò fra un paio di settimane, cioè appena fatto l'esame di laurea.

Ho fatto passare i calcoli contenuti nel Suo lavoro e mi pare che a pag. 376, nelle formule (6) vi sia un errore di stampa perché a me parrebbe che invece di scrivere:

$$f_2 = \frac{\alpha}{a}\frac{y}{r^2} \text{ ed } f_2 = -\frac{\alpha}{a}\frac{x}{r^2}$$

occorra scrivere

$$f_2 = -\frac{\alpha}{a}\frac{y}{r^2} \text{ ed } f_2 = \frac{\alpha}{a}\frac{x}{r^2}$$

inoltre nella linea 3 della stessa pagina, invece di §2, va scritto evidentemente §1.

Così pure nella espressione del vettore $F_2(P)$ data a pag. 378 mi pare che devono essere cambiati di segno i termini del secondo membro.

---

[169] "(Sulla formola di Simpson). Estratto da una lettera al Prof. Peano", *Rivista di Matematica edita da G. Peano*, v. III, pp. 60-61.

Profittando della Sua cortesia mi permetto di chiederle uno schiarimento a proposito della relazione

$$\frac{\delta X}{\delta x}+\frac{\delta Y}{\delta y}+\frac{\delta Z}{\delta z}=\frac{1}{R_r}+\frac{1}{S_r},$$

scritta nella penultima riga di pag. 375. Di questa relazione ho trovato una dimostrazione vettoriale semplicissima, che è la seguente.

Se indico con N un vettore unitario normale alla superficie che si considera (il quale ha per componenti X, Y, Z) allora il primo membro della precedente relazione vale divN, cioè I, $\frac{dN}{dP}$ e questo invariante primo vale precisamente la curvatura media della superficie in P.

Ho cercato nella Geometria differenziale del Bianchi,[170] la dimostrazione cartesiana, della citata relazione, ma non vi è.

Ho cercato per mio conto una dimostrazione cartesiana, ma quella che ho trovato non mi pare molto semplice.

Sono partita dalle note formule di Rodrigues

$$\frac{\delta X}{\delta u}=\frac{1}{r}\frac{\delta x}{\delta u},\quad \frac{\delta Y}{\delta u}=\frac{1}{r}\frac{\delta y}{\delta u},\quad \frac{\delta Z}{\delta u}=\frac{1}{r}\frac{\delta z}{\delta u},$$

che ho trasformato facendo figurare le derivate di X, Y, Z rispetto ad x, y, z, e allora vengo ad ottenere un determinante, che è

$$\begin{vmatrix} \frac{\delta X}{\delta x}-\frac{1}{r} & \frac{\delta X}{\delta y} & \frac{\delta X}{\delta z} \\ \frac{\delta Y}{\delta x} & \frac{\delta Y}{\delta y}-\frac{1}{r} & \frac{\delta Y}{\delta z} \\ \frac{\delta Z}{\delta x} & \frac{\delta Z}{\delta y} & \frac{\delta Z}{\delta z}-\frac{1}{r} \end{vmatrix}=0$$

che sviluppato da un’equazione di II° grado rispetto ad 1/r, e la somma delle radici di essa vale $\frac{\delta X}{\delta x}+\frac{\delta Y}{\delta y}+\frac{\delta Z}{\delta z}$ e così risulta dimostrata la relazione citata.

Parrebbe però che debba esserci qualche dimostrazione più semplice; se Lei sapesse indicarmela le sarei più che infinitamente grata.

Perdoni nuovamente la mia grande libertà, e con molti ringraziamenti e rispettosi saluti mi creda

Sua Dev
Luisa Pelosi

Torino. Via Garibaldi 3.

---

[170] L. Bianchi, *Lezioni di geometria differenziale*, Enrico Spoerri, Pisa, 1894.

**144**
Gian Antonio Maggi a Luisa Pelosi

Milano 14 Dicembre 1925.

Gentilissima Signorina,

La ringrazio dell'avvertito scambio di segni e dell'errato §, dei quali ho preso nota. La dimostrazione vettoriale, ch'Ella mi accenna, della nota formola è certamente semplicissima, salvo le riserve che si possono fare sulla semplicità intrinseca dell'algoritmo delle omografie. Una dimostrazione cartesiana conforme alla Sua, ma abbreviata, si trova partendo da

$$dX=\frac{1}{r}dx,\quad dY=\frac{1}{r}dy,\quad dZ=\frac{1}{r}dz,$$

dove x, y e z s'intendono appartenere ad una linea di curvatura, e introducendovi

$$dX=\frac{\delta X}{\delta x}dx+\frac{\delta X}{\delta y}dy+\frac{\delta X}{\delta z}dz,\quad dY=\frac{\delta Y}{\delta x}dx+\ldots,\quad dZ=\frac{\delta Z}{\delta x}dx+\ldots$$

Io ho sostituito con questa, a suo tempo, per conto mio, la più complicata dimostrazione che si trova in Knoblauch, Einleitung in die allg. Theorie der Krummen Flächen (Liepzig, Teubner, 1888) pag. 45, [171] conducente, per una via più lunga, allo stesso determinante eguagliato a 0. Nel trattato del Bianchi[172] effettivamente la formola non si trova, e mi ricordo che, avendogliene parlato, egli non mostrò di sentire il bisogno di farvela comparire.

Coi miei cordiali saluti e buoni augurii, mi creda sempre

Dev.[mo] Suo

Gent. Sig.[a] Luisa Pelosi
Via Garibaldi, 3, Torino.

---

[171] J. Knoblauch, *Einleitung in die allgemeine Theorie der Krummen Flächen*, Liepzig, Teubner, 1888.
[172] Si veda la nota alla lettera precedente.

**145**
Gian Antonio Maggi a **Salvatore Pincherle**[173]

Carissimo Pincherle,

Ricordando l'accenno di un tuo desiderio, nelle conversazioni romane di gradita memoria, ho dato un'occhiata alla Nota del Mineo e del Gugino nel Bollettino, e riveduto la Nota del secondo, oggetto degli appunti del primo. Il Mineo mi fa un po' l'impressione di cercar il pelo nell'uovo, e dai suoi appunti mi sembra che il Gugino si difenda abbastanza bene, finendo anche per rilevare, a carico del Mineo, una eguaglianza inesatta, che però apparisce corretta nella successiva Nota.

È poi singolare che il Mineo, cercando la festuca, non abbia avvertito la trave, scovata dal Cartan, grazie alla quale riesce infirmato il teorema della prima Nota lincea del Gugino. Il quale, nella Nota lincea "Sulla validità ecc." oggetto degli appunti della seconda Nota del Mineo, corre, non senza dar prova d'ingegno, al riparo. Ma il primitivo teorema, limitato all'ipotesi della velocità nulla, cioè ad un istante particolare del movimento, perde la maggior parte del suo significato. E il teorema rinnovato, colla nuova definizione dell'effetto cinetico-dinamico, non mi sembra poter conservare il significato del primitivo, per artificiosità e complicazione della nuova definizione, in confronto della semplicissima antica.

Dirò ancora che il primitivo teorema, demolito poi dal rilievo del Cartan, non mi persuadeva, per un certo contrasto che mi pareva intravvedere tra quel massimo ad ogni istante, e il minimo della forza viva totale in un intervallo di tempo, che forma il contenuto del principio della minima azione. Tenni però il dubbio per me, tanto più che l'autore non mi aveva ancora mandato il suo lavoro e non avevo quindi l'occasione di farmi vivo con lui. Io mi permetto di credere che il Cartan fu indotto da simile impressione a rifare i calcoli fino all'origine, e scoprire a quel modo termini dal 4° ordine in $\tau$, che avrebbero dovuto comparire nella formola finale, e che, a onor del vero, restano insidiosamente nascosti. Contiamo trattenerci a Valnegra fino alla metà del prossimo Settembre, richiamati allora a Milano dalle nozze della mia figliola minore coll'Ing Alberto Ferratini, cognato del collega e amico Chisini.[174] Egli è pure qua, colla famiglia, e sono la nostra sola, ma anche ottima compagnia.

Sapendo che ti scrivo, m'incarica de' suoi saluti. A te e alla tua gentile Signora io invio le migliori cose, e, con una cordiale stretta di mano, mi confermo

Tuo aff$^{mo}$ amico
G.A. Maggi

Valnegra ( Bergamo)
9 agosto 1930.
(Spedita il 14)

P.S Ho per le mani l'Introduction ecc. del De Broglie,[175] e penso che ne potrei fare una recensione pel Bollettino, punto, avendo visto che è tra i libri ricevuti, non abbiate già procurato con qualcuno che potrà fare anche meglio di me.
Di nuovo le più cordiali cose.

---

[173] Questa lettera si trovava nelle *Rusticationes*.
[174] La figlia minore di Maggi si chiamava Bianca e l'ing. Ferratini era fratello di Ada, moglie di Chisini.
[175] Potrebbe trattarsi del volume: *Introduction à l'étude de la mécanique ondulatoire*, Hermann, Paris, 1930.

**146**

Salvatore Pincherle a Gian Antonio Maggi[176]

[cartolina postale indirizzata a: Ch$^{mo}$ Signor Prof. Comm. G.A. Maggi
Valnegra - (Bergamo)]

Montese (Modena) 19.8.30

Carissimo amico,

Grazie per le esaurienti informazioni sulla polemica M.-G. Ben volentieri pubblicherò sul Bollettino la tua recensione sul De Broglie, e ringraziamenti cordiali anticipati.

Anche a nome di mia moglie, che m'incarica dei tuoi saluti distinti; vivi rallegramenti per il fidanzamento della tua gentile figliuola; rallegramenti e saluti a Chisini e Signora, entrambi miei antichi scolari.[177]

A te la più affettuosa stretta di mano dal tuo

S. Pincherle

**147**

Gian Antonio Maggi a **Vittorio Pizzini**

*[Scritto in rosso:]* Traduzione di due biglietti a Vittorio Pizzini indovinando più che decifrando il corsivo tedesco della scrittrice Frau Line Hernssel
Lipsia 05
Hofstr. 12 II.

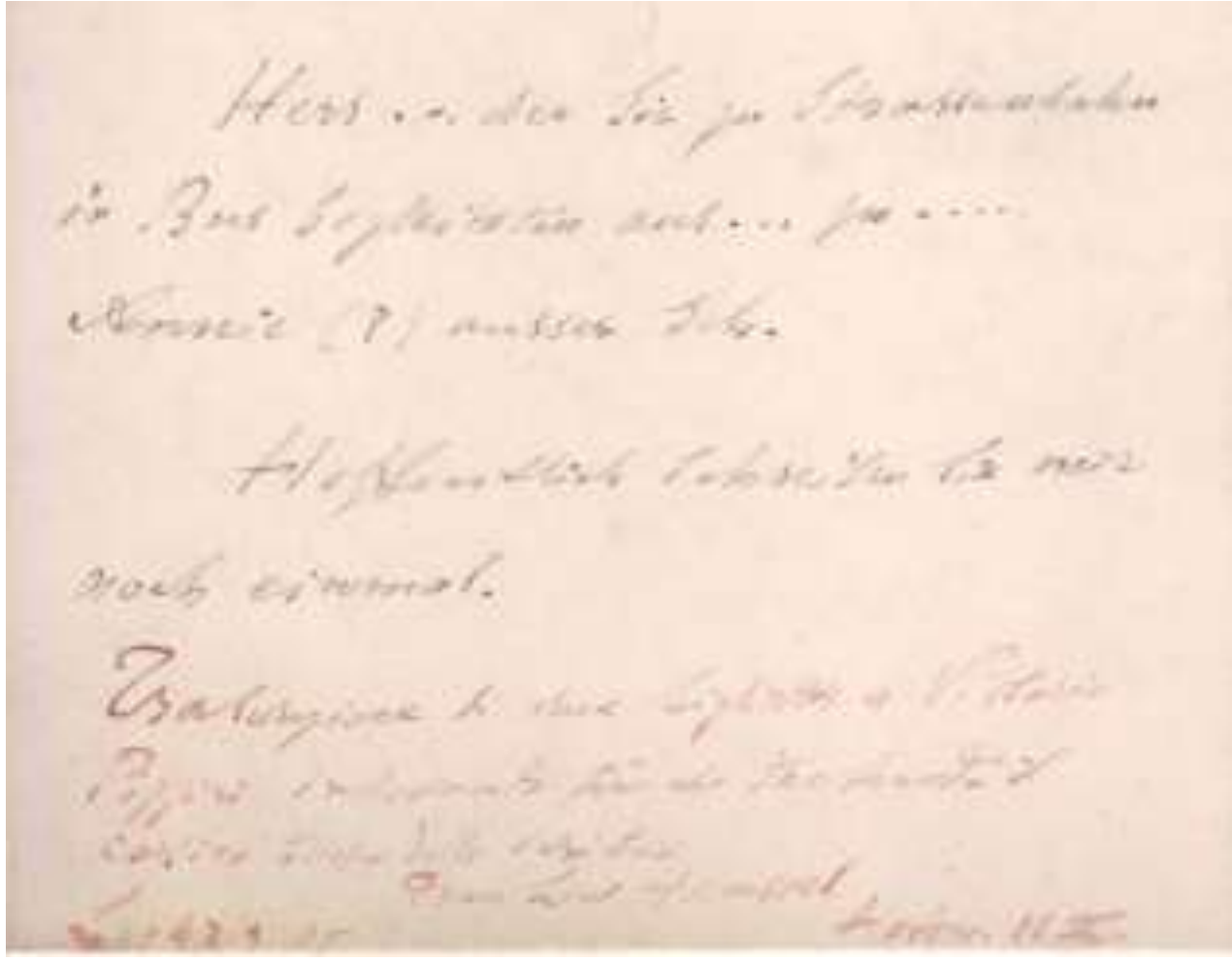

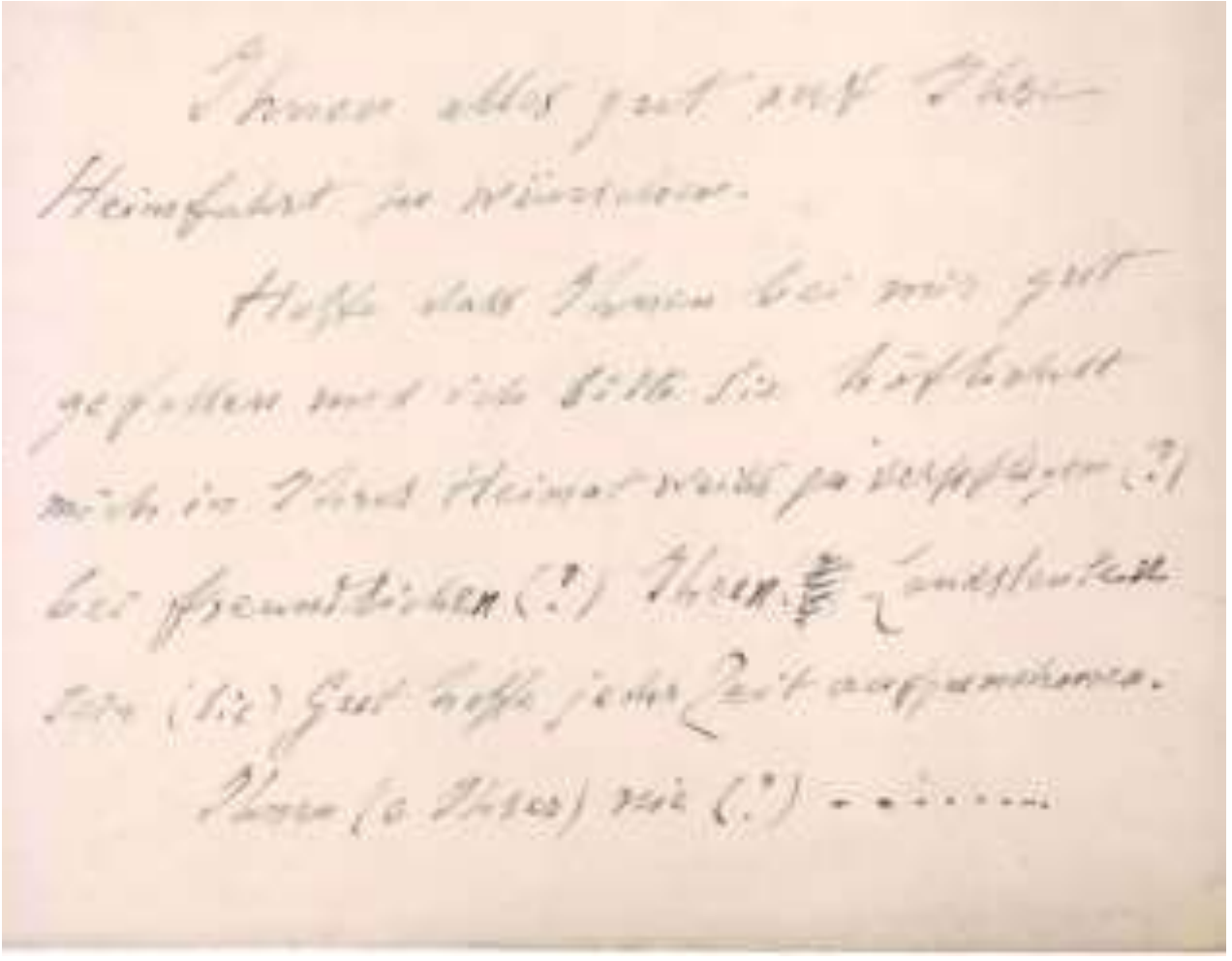

[176] Questa lettera si trovava nelle *Rusticationes*.

[177] Pincherle si riferisce alla Facoltà di Matematica dell'Università di Bologna, dove Chisini si laureò nel 1912 e la sua futura moglie, Ada Ferratini, si iscrisse nell'a.a. 1913/14 (nel 1916/17 era ancora iscritta al terzo anno, poi probabilmente si ritirò, perché non risulta essersi laureata).

**148**
Gian Antonio Maggi a **Virgilio Polara**[178]

Casa Lunedì (7 Febbraio 1910)

Caro Dr Polara,

Tornando sul Suo saggio, trovo da proporle questa difficoltà.

Nell'esperimento della macchina d'Atwood, ora si trascuri la massa della carrucola, le differenze fra la accelerazione dei due pesi nel relativo movimento e nel movimento libero per gravità hanno grandezze a norma del principio dell'eguaglianza dell'azione e della reazione. Invece le loro orientazioni sono eguali, e non eguali e contrarie. Sta bene che la tensione applicata a ciascun corpo apparisce tendente, diciamo così, ad accorciare il filo: ma, per la particolare disposizione del sistema, direzione e senso riescono in tal modo i medesimi.

Naturalmente, si rientra nella meccanica ordinaria, attribuendo quelle differenze di accelerazione alle resistenza del sostegno, che si esercita pel tramite del filo.

Questa osservazione mi riconduce al mio primo concetto che la mutua azione dei due pesi non si adatti troppo bene alla teoria della macchina d'Atwood.

A simili osservazioni si presta la macchina d'Atwood con piani inclinati, dove intervengono le resistenze di questi altri sostegni. Difatti, ch'io sappia, il Mach non invoca questo genere di esperimenti, ai quali invece accennò il Vailati.

Ci pensi. Io non posso che incoraggiarLa vivamente a proseguire in codesto lavoro.
Con tanti saluti

Aff$^{mo}$ Suo

(Cfr. Cartella "Note Varie")

[178] Si vedano anche gli scritti #230 e #231, che potrebbero far parte di questo carteggio.

**149**

Virgilio Polara a Gian Antonio Maggi

[busta indirizzata a: Ill.mo Prof. G.A. Maggi - Via del Risorgimento 1 - Pisa; consegnata a mano]

7 febbraio 1910

Ill.mo Professore,

La ringrazio sentitamente del vivo interessamento mostrato a mio riguardo e della gentile premura avuta nel consigliarmi.

Quanto alla Sua osservazione, che cioè le differenze di accelerazioni dei due pesi nel relativo movimento e nel movimento libero per gravità sono dirette nello stesso senso - e non in senso contrario - a me pare che, considerando il sistema complessivo, ciò non sia vero. Noi possiamo infatti riferirci ad un punto qualsiasi del filo, per es. al punto culminante nella gola della carrucola. Allora è manifesto che le differenze di accelerazioni appaiono dirette in senso contrario. Del resto mi pare evidente che queste due accelerazioni tenderebbero da sole a muovere la carrucola in senso opposto, tanto che sono l'effetto di due tensioni in senso contrario. Le tensioni, naturalmente, possiamo riferirle ad uno stesso punto del filo, con che il loro senso opposto risulta evidente.

Quanto poi al fatto che, introducendo la tensione del filo, si rientra nella Meccanica ordinaria, Le faccio osservare che da questa tensione io ho a Lei fatto cenno solo per esprimere il meccanismo vero del fenomeno, ciò per farle vedere che in realtà la mutua variazione di accelerazione dipende dal mutuo rapporto delle masse dei due corpi; e ciò in relazione a quanto Ella mi obbiettò in principio.

Si intende che nella trattazione io piglierò il fatto in sé, definirò quel rapporto costante fra le variazioni di accelerazioni come il rapporto inverso delle masse e solo per dare una giustificazione immediata di questa definizione - cioè per far vedere come la massa d'un corpo sia qualcosa come l'inerzia al moto nel corpo stesso - invocherò il fatto, quasi evidente di per sé stesso per la disposizione che si vede assumere al filo, che la tensione agisce egualmente sui due corpi. Mentre nel metodo ordinario per avere la misura della massa occorre la misura statica della forza, qui si accenna semplicemente alla nozione di tensione come chiarimento.

Solo stamane, in seguito alla Sua citazione, ho letto le brevi osservazioni del Vailati. Con vivo piacere ho riscontrato che egli era del mio ordine di idee. Egli però non fa alcun cenno sul modo di introdurre la massa servendosi della Macchina d'Atwood o del piano inclinato, limitandosi solo a dire che questi apparecchi potrebbero servire. Egli inoltre non ha rilevato il fatto, che a me pare fondamentale, della maggiore generalità e comodità che presenta questo metodo in confronto all'altro che per definire la massa parte da postulati sui corpi isolati.

Molto volentieri sarei venuto di persona a ringraziarla e a darle questi brevi appunti: ma non è stato possibile per oggi.

Scusi, egregio Professore, se profitto un po' troppo della Sua gentilezza.

Mi creda distintamente.

Dev.

Virgilio Polara

**150**
Gian Antonio Maggi a **Cino Poli**

Milano 30 Giugno 1927

Egregio Professor Poli,

Leggendo la Sua Nota nell'ultimo fascicolo dei Rendic. dei Lincei,[179] vi ritrovo, con differenze che conferiscono manifesta impronta di originalità al Suo lavoro, un'analogia col mio svolgimento dal conforme soggetto, quale si trova nel Capitolo III della mia Dinamica dei Sistemi, che non saprei rinunciare a far presente alla Sua benevola attenzione, fosse solo per compiacermi dell'accordo, inteso a vedute che non sono le più comuni. È pur vero che chiamo postulati le ipotesi con cui introduco, nello stesso Capitolo, i vincoli, (pag. 129) e le reazioni vincolari (pag. 260). Ma, alla stregua del concetto generale così dei vincoli come delle pressioni, quale è stabilita nel Capitolo I (§2 e §12), tali ipotesi o postulati ricevono, in sostanza, il significato di definizione dei più precisi vincoli, che occorrono nel considerato problema, e delle pressioni supposte capaci di sostituire i vincoli medesimi, nel sistema, concepito come libero.

Io poi non faccio distinzioni di vincoli lisci, perché comprendo le pressioni di attrito (Op. cit. §85 ed Elementi di Statica §53) tra le forze impresse.

Proposta, che trova riscontro nella teoria dell'attrito dal Painlevé, e non so perché, o malgrado di un'innegabile semplificazione di importanti enunciati, non che non accolta, non sia mai stata neppure presa in considerazione.

Ben minor conto ho sempre fatto dell'adozione della mia sostituzione del sistema dei corpi rigidi al tradizionale sistema di punti materiali, quantunque la minor semplicità delle formole mi sembri largamente compensata dall'immediata applicazione a problemi concreti, come quelli dei sei Esempii, che nel ricordato mio libro, accompagnano l'esposizione generale.

Mi è grato intanto presentarLe i miei migliori saluti e pregarLa di credermi

Aff Suo
G.A.M.

**151**
Cino Poli a Gian Antonio Maggi
[busta indirizzata a: Ill$^{mo}$ Signor Prof. G.A. Maggi - Lanzo d'Intelvi (Como); mittente: D.$^{r}$ Cino Poli - Torino - Via Vico, 2 - Tel. 47.960]

Torino 24 settembre 1927

Ill.$^{mo}$ Signor Professore,

sono profondamente mortificato di non aver finora risposto alla gentilissima Sua lettera. La prego di perdonarmi l'involontaria mancanza. Quando la ricevetti ero assorbito dagli esami, poi lasciai Torino e mi mancò il tempo di consultare il Suo trattato, per poterLe rispondere esaurientemente.

---

[179] C. Poli, "Sui principî della meccanica analitica", v. V, 1927, pp. 656-661.

Per ciò che riguarda le Sue osservazioni sulla mia Nota,[180] riconosco di buon grado che la mia trattazione è formalmente identica a quella già da Lei data. Temo però di non aver chiaramente espresso nella mia Nota che non intendevo far cosa originale con la trattazione analitica del principio dei lavori virtuali o di quello di d'Alembert: ma sibbene modificare l'interpretazione che se ne dà generalmente.

Siamo d'accordo che quei principi si riducono in sostanza ad enunciare le due equazioni seguenti:

(I) $m\frac{d^2x}{dt^2} = X + X'$

(II) $\sum X'\delta x = 0$

Ma quale è il loro contenuto fisico? Secondo me non si può dire che esprimono delle proprietà dei vincoli che dobbiamo ammettere come postulato o dato di esperienza. Mi spiego. Dobbiamo distinguere fra vincoli fisici, concreti, e vincoli astratti, matematici. Fisicamente un punto vincolato è un punto soggetto a certe forze provenienti da azioni di altri corpi; cioè non differisce dal punto libero. Che differenza c'è fra il moto circolare di una sferetta attratta da un punto fisso con azione newtoniana e il moto della medesima sferetta legata a quel punto con un filo? La diversa origine o natura fisica della forza, che analiticamente porta a problemi di natura diversa perché nel primo caso la forza non dipende dalla velocità e si suppone nota; nel secondo è funzione della velocità e generalmente non si conosce, ma si sa che contiene una componente che equilibra sempre la reazione centrifuga. Se il filo è elastico, quel fatto non è sufficientemente per permettere l'integrazione della equazione del moto; bisogna conoscere l'espressione della forza, in funzione della distanza dal centro fisso, come nel 1° caso. Se il filo è inestendibile (meglio il suo allungamento relativo è minimo, trascurabile), ne deduciamo la conoscenza della traiettoria e quindi la possibilità di ignorare la componente normale delle forze applicate (sia che provengano dal filo, come qualsiasi altra). Ma non possiamo dire che si possono ignorare le reazioni del vincolo, perché in generale vi sono anche delle componenti tangenziali (p. es. se invece del filo si ha una asticciuola girevole con attrito).

Dunque io credo che non si debba postulare che le reazioni dei vincoli soddisfacciano alla (II), ché anzi la distinzione tra forze applicate e reazioni vincolari è artificiosa. Ma che solo si può dire che non ponendo il secondo membro delle equazioni del moto come è indicato dalla (I), ove le X' soddisfacciano alle (II), per la integrazione delle equazioni del moto è sufficiente conoscere le X e le equazioni che definiscono le δx che compaiono nella (II). Questo risultato è utile tutte le volte che per una ragione qualsiasi si è in grado di sapere che il sistema non ammette certamente certi spostamenti.

Non so se sono riuscito più chiaro in queste poche righe che non nella mia Nota: temo di no. Forse ho anche impostato poco bene la quistione che non è tanto sul significato o modo di presentare il principio dei lavori virtuali, quanto piuttosto: che cosa dobbiamo intendere col nome di vincoli?

Le rinnovo Le mie scuse vivissime, insieme ai ringraziamenti per le Sue osservazioni e Le sarò gratissimo sempre se vorrà esprimermi il Suo autorevole pensiero.

Gradisca i miei rispettosi ossequi.

Suo dev$^{mo}$
Cino Poli

[180] Si veda la lettera precedente.

**152**
Gian Antonio Maggi a Cino Poli

Risposta al Prof. Cino Poli.

La ringrazio della Sua lettera, che ho letto con interesse. Che se, col prof. Somigliana, ho accennato a non aver ricevuto da Lei risposta alla mia lettera della scorsa estate, è perché mi sembra desiderabile che gli studiosi, almeno finché non si esce troppo dal campo dei proprii studii più speciali, s'intendano anche colla corrispondenza personale, piuttosto che limitarsi ad un mero puro scambio di stampati. Vedo quanto Ella scrive. Tutto sta sul significato che attribuisce al termine postulato. Nella Matematica Pura è una proposizione intorno alla quale non può più esserci da discorrere. Nella Matematica Applicata difficilmente può scompagnarsi il concetto del postulato da quello di un'ipotesi, che, se non da ragioni a priori, vuol essere giustificata dall'opportunità delle conseguenze.

Ella vorrebbe che dell'equazione

(II) $\sum X' \delta x = 0$

non si dicesse più di questo, che serva col concorso delle

$$m\frac{d^2x}{dt^2} = X + X' \quad \text{e} \quad \sum X' \delta x = 0$$

a determinare le X', formandosi tante equazioni indipendenti quante sono le incognite, le x e le x'. A pag. 163 del ricordato mio libro io dico appunto "Questa considerazione rende ragione, fino ad un certo punto, del secondo postulato". Ma poi accenno alle espressioni che ne risultano per le X', come quello che, caso per caso, appariscono come le più indicate per forze atte a sostituire il vincolo. Questo non mi sembra dare alla proposizione in discorso base sperimentale, ma puramente corroborare l'opportunità analitica con una opportunità meccanica, tanto più significativa che le stesse espressioni particolari, relative ai singoli casi, sono la ragione storica del principio dei lavori virtuali, di cui non occorre richiamare la connessione colle (II).

Per quanto poi a distinguere tra forze imposte (o applicate) e forze vincolari, io ne elimino ogni sostanziale differenza, col concetto della pressione, dipendente, in generale anche da tempo e velocità, ch'io introduco §12 *[sic!]*, come una delle due convenzioni destinate alla effettuazione del Calcolo del Movimento, l'altra convenzione consistendo nei vincoli imposti al sistema (§2): intendendo che, a comporre le forze impresse, concorrono, in generale, pressioni preventivamente date, così chiamate in confronto delle pressioni vincolari; da calcolare nell'indicato modo.

Ella mi farà sinceramente piacere, quando abbia occasione di venire a Milano, a cercare di me, per discorrerne insieme anche a voce. Intanto La prego di aggradire i miei migliori saluti e di credermi

Aff^mo Suo
G.A. Maggi

**153**
Gian Antonio Maggi a Cino Poli

Milano 30 Dicembre 1927

Egregio Professor Poli,

La ringrazio vivamente del cortese invio del suo bel volume “Meccanica Generale e Applicata”,[181] di cui ho preso visione, rilevandone pregi di composizione e d’esposizione, e, con piacere, senza pregiudizio, della Sua originalità, un notevole accordo, sotto parecchi aspetti, coi miei concetti.

Così, convengo con Lei nello sviluppo, maggiore del solito, della teoria del movimento relativo. E quello che Ella introduce come velocità di un sistema corrisponde al mio atto di movimento. Mi permetta poi un’osservazione, che al postutto Le prova la mia attenzione? Perché limitare la definizione del moto traslatorio colla condizione che le traiettorie di tutti i punti siano rette parallele, anziché linee qualsivogliano riducibili l’una all’altra con uno spostamento traslatorio? Io definisco il movimento traslatorio come quello in cui è traslatorio lo spostamento dalla posizione inziale alla posizione al tempo generico, coll’uso di quello spostamento, che non ho trovato troppa adesione, e pur tuttavia a me sembra concetto efficace.

Nuovi ringraziamenti, e cordiali saluti e augurii coi quali La prego di credermi

Aff Suo
G.A. Maggi.

**154**
Gian Antonio Maggi a **Giovanni Polvani**

*[in penna rossa:* Non mandata*]*

Casa 22 Febbrajo 1932

Carissimo Polvani,

Confermo che, col maggior interesse, ho letto e meditato il suo bel Discorso, bello per l’argomento e per lo svolgimento, né ho mancato di gustare il toscano sapore della forma. E, per cominciare dalla conclusione, convengo pienamente con Lei e con Lattanzio, ed è mia antica opinione, tanto che nella prefazione dei miei Principii ecc. del 1896, dico “Per quanto riguarda la sostanza dei principii ... mi sono attenuto a quell’idea della Filosofia Naturale ... che riduce la spiegazione di un fenomeno alla determinazione della logica connessione di un fatto ideale, parziale immagine del fenomeno, in quanto che ne riproduce alcuni aspetti, con certi fatti che si assumono come postulati cardinali, e hanno la loro ragion d’essere in una simile riproduzione di leggi generali della Natura. Idea, che, mentre lascia sempre adito alla ricerca di una spiegazione più soddisfacente ... mette da parte quella spiegazione vera, che fa respingere tutte le possibili, per distruggere in alcuni ogni fede nella Scienza, e trarre altri a cercare ciò che presumibilmente non è conseguibile”.

Questo concetto, che, a quei tempi antichi, cominciava a farsi strada, colle riforme del Kirchhoff e del Hertz dell’esposizione della Meccanica, al tempo stesso che certamente riconosce la “scientia cum ignoratione conjuncta et temperata”,[182] mi sembra però fare, a codtsta conjunctio e temperatio una parte alquanto diversa da quella che si è fatta dal concetto impostosi in questi ultimi anni, direi dopo la guerra, contemplato nel suo discorso. L’antico concetto mi sembra tracciare più distintamente la separazione tra la scienza e l’ignoranza, e consolarci con una maggior ragione di vero, senza rimangiarmi con questo la rinuncia alla spiegazione vera.

[181] C. Poli, *Meccanica generale e applicata*, Torino, UTET, 1927.

Prendiamo, per esempio, la teoria idrodinamica dei vortici, che reca i nomi di Helmholtz e di Thomson. Ammettiamo tutto quello che ci occorre per scrivere le equazioni idrodinamiche di Lagrange: in particolare la continuità del fluido, che è una semplice apparenza, e l’ipotesi del fluido perfetto, ch’è una semplice approssimazione. Un calcolo impeccabile ci serve i vortici, e le loro proprietà, che le esperienze di Tait (Recent advances in Science, Lecture XII)[183] verificano appuntino cogli anelli di fumo. E sotto lo stesso aspetto, si può, richiamare la teoria della doppia rifrazione, nell’Ottica Fisica, la teoria delle onde elettromagnetiche ecc.

I concetti di probabilità e di statistica s’introducono colla teoria molecolare dei gas, non più giovane della teoria dei vortici. Ma ciò che mi sembra caratteristico nell’uso che di questi concetti fa il nuovo indirizzo è, nella strepitosa folla, ch’Ella giustamente ci ricorda ascendere a miliardi di miliardi di miliardi, la ricerca dell’individuo. L’atomo; e, per l’atomo, l’elettrone e il protone. Se non che, tra l’altre cose, le due condizioni iniziali che determinano il movimento di un corpo, soggetto a forze date, coll’elettrone non sono atte a ricevere ambedue ad un tempo, un significato concreto. Il momento dell’elettrone ne risulta rimesso al caso. Ma il concetto statistico provvede ad aggiustare le cose. Ditemi, dal mucchio dei caratteri di stampa, mescolati ben bene a tondo, la lettera che estrarrò. È come domandare dove andrà l’elettrone. Ma estraggiamo dal mucchio duemilaottocentottantuno lettere, e vi dirò che vi si troveranno tante a, e tante b, tante c ecc.

Sta bene. Ma, mentre l’indagine dell’elettrone è un passo avanzato nella ricerca della spiegazione vera, a cui accennavo in principio, io mi domando, se, col rimettere il suo movimento al caso, per poi rifugiarsi nei grandi numeri, si fa un conforme passo avanzato nel proposto possesso del vero.

Per cui, senza menomare l’interesse dei nuovi concetti, io semplicemente domando che il momentaneo abbandono della Fisica Matematica così detta classica sia considerato un cambiamento di rotta, piuttosto che la conquista di una posizione che la metta fuori di combattimento. Domanda, che il Suo Discorso mi offre l’occasione, ma, distinguiamo bene, non la ragione di fare.

Con questo Le rinnovo ringraziamenti e rallegramenti, e con una cordiale stretta di mano mi confermo Aff Suo

[182] L.C.F. Lattanzio, *Divinæ institutiones*, libro III, cap. VI.

[183] P.G. Tait, *Lectures on some recent advances in physical science with a special lecture on force*, Londra, MacMillan, 1876; la Lecture XII si intitola: “Structure of matter”.

**155**

**Georg Sabinine** a Gian Antonio Maggi[184]

[busta indirizzata: A Meur Gian Antonio Maggi professeur à l'université de Pise
Via Nuova, 16 - Pisa;
mittente: Russie - Moscou de la part de Sabinin e scritte in cirillico]

Monsieur et tres honoré Collegue!

Au congrès de Paris, le professeur Vassilief (de Kasan) a eu l'occasion de vous présenter ma fille, qui s'occupe aussi de mathématiques. C'est à cause de cela, que je me permets de vous adresser une prière à son sujet.

Voici ce dont il s'agit: elle a composé deux articles qui concernent la Géometrie elementaire traitant la théorie des parallèles. L'un porte sur la détérmination de la somme des angles du triangle rectiligne, l'autre a pour sujet le postulatum d'Euclide.

Ma fille desirerait que ces articles paraissent dans un journal imprimant les articles écrits en français.

Ayez la bonté de prendre connaissance de ces articles que je viens de lire et que je trouve exposés exactement. Or mon opinion ne saurait pas suffire, ma specialité n'étant pas la théorie que l'on trouve publiée pour la première fois dans la brochure connue de Lobatschevsky et qu'ont poussé plus loin en la développant avec details les Géometres Beltrami, Rieman, Helmholmtz, Klein, Houel, Pointcarret, et les outres. C'est pour cela que je vous demande d'en être le juge.

Si vous etes après lecteur de mon avis, soyez assez bon pour user de votre influence aupres d'un rédacteur du journal que vous connaissez, en les très faisant accepter.

Nous vous en serons , ma fille et moi, très reconnaissants.

Si non, veuillez, je vous prie, nous renvoyer les manuscrits.

Veuillez agréer Monsieur très honoré Collegue, l'assurance de ma considération la plus distinguée

G. Sabinine

Mon adresse: Russie, Moscou, Arbatte, ruelle Nicoleky, maison Sabolotsky (28)

$19\frac{X}{14}00$

[184] Alcuni nomi e alcuni termini hanno una grafia errata e alcuni accenti non sono coerenti con una corretta grammatica, ma vengono lasciati come scritti dall'autore; non si inseriscono i *[sic!]* per non appesantire il testo. Così per tutte le lettere di questa corrispondenza.

**156**

Gian Antonio Maggi a Georg Sabinine

Monsieur et très honoré Collègue

J'espère que vous aurez reçu ma lettre du 28 Octobre.[185] Je regrette de n'avoir pu prendre connaissance plus tôt du manuscrit que vous avez bien voulu m'envoyer, d'autant plus que je ne pourrais pas m'empêcher d'y faire les quelques observations qui suivent.

Je ne vois pas de quelle façon la démonstration dont la première Note est l'objet serait à l'abri des objections sur lesquelles l'auteur, dans son introduction, montre d'être parfaitement renseigné. Dans les relations comme
$B=A_1+\Delta+H$,
$B+E_1=A_1+A_3+\Delta$ etc.
(pag. 8) il n'est pas question de comparaison d'espaces angulaires entre eux, mais de la composition d'une aire à l'aide de ces espaces et de l'aire du triangle; je ne saurais pas comment cela seront possible sans attribuer à l'espace angulaire la signification d'une aire infinie. On ne poussait alors admettre raisonnement à l'aide duquel on parvient aux deux formules (4) et (6) sans méconnaitre la valeur d'objections à l'usage d'un aire infinie telles que celles que nous lisons dans le 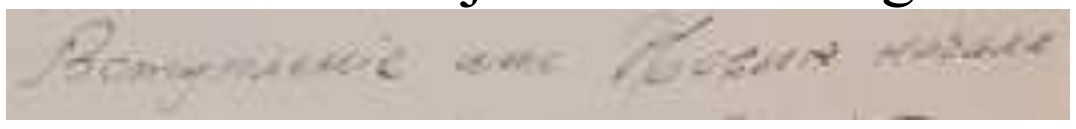 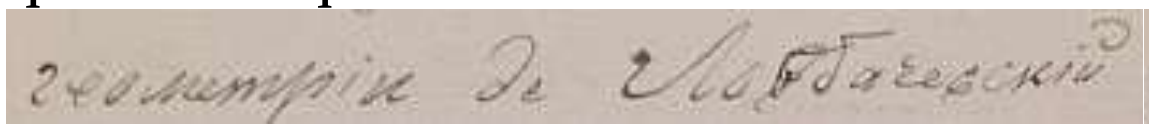 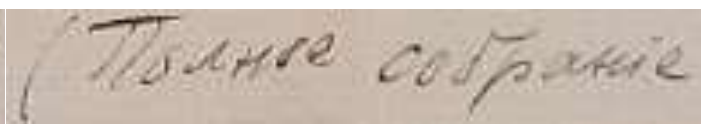 pag. 222 et suivantes).
An demeurant je n'ai aucun besoin de rappeler que la question de la possibilité de démontrer le postulat des parallèles à l'aide des autres postulats de la géométrie euclidienne est considérée comme résolue negativement depuis que Beltrami a démontré qu'une espace dans lequel la géométrie de Лобачевский[186] est verifiée est materiellement possible, celle géométrie coincidant avec celle de la surface pseudosphérique. Je ne voudrais jamais faire appel à "ipse dixit". Mais il n'est pas moins vrai que ce serait tout un edifice solidement bâti qu'il faudrait désormais demolir pour arriver à une conclusion contraire.

Quant à la deuxième Nota, la démonstration est fondée sur ce que la limite d'une sphère lorsque son rayon croit infiniment est un plan. Cela revient à admettre d'avance le postulat des parallèles, Лобачевский ayant démontré que, avec un angle de parallelisme différent du droit, la surface à laquelle une sphère s'approche infiniment en augmentant son rayon est une horisphère 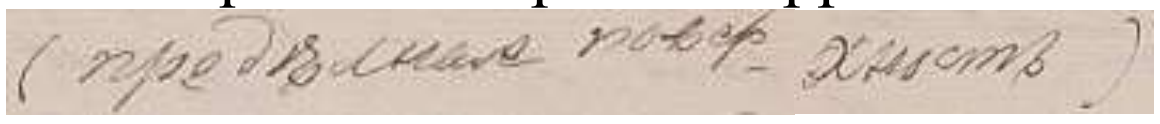 c'est-à-dire la surface engendrée par la rotation d'un horicycle 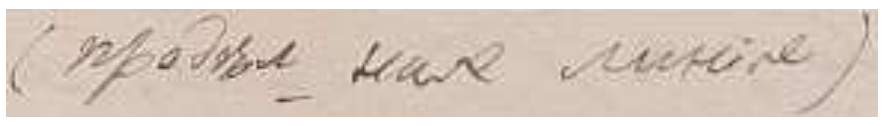 autour d'un axe 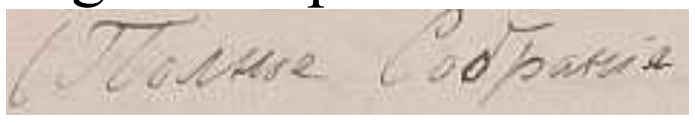 pag. 317 et suivantes).

Je garde à vos ordres le manuscrit afin de l'avoir sous les jeux, si vous croyerez bien de me répondre par quelques observations.

Je vous prie, Monsieur et très honoré Collègue, de présenter de ma part mes salutations bien respectueuses à Mademoiselle Sabinine et d'agréer l'assurance de ma considération la plus distinguée.

Votre tout dévoué.

Pisa le 20 Décembre 1900.

---

[185] Non è presente nel Fondo.
[186] Lobatschevsky.

**157**
Georg Sabinine a Gian Antonio Maggi

Moscou le 5 Mars 1901.

Monsieur et très honoré confrère!

D'après votre dernière lettre, je crois mon devoir de vous informer que jusqu'à présent je n'ai reçu aucun mot à l'égard du manuscrit de ma fille. J'ai reçu de vous le 11 Janvier 1901 une votre petit mémoir "Sulla theorie del pendulo"[187] et votre programme de mecanique analytique.

Ainsi je n'eu pas reçu la lettre que vous m'avez envoyée le 20 Decembre 1900.

Permettez moi de vous prier de m'envoyer une autre lettre contenante tout ce que vous avez écrit à l'égard du manuscrit de ma fille dans votre lettre (le 20 Decembre 1900) qui renfermait votre réponse sur ma lettre (17 Octobre 1900) et qui, comme l'on supposez, est perdue.

En attendant j'ai l'honneur d'envoyer à vous mes salutations les plus distinguées
Votre tout devoué

G. Sabinine

Mon adresse: Russie, Moscou, Arbatte, Nicolsky road, maison Sabolotsky 28

**158**
Georg Sabinine a Gian Antonio Maggi

3 Avril 1901.

Monsieur et tres honoré Collegue!

Je vous suis très reconnaissant de votre lettre (du 13 Mars 1901)[188]. Il m'est très agreable d'y repondre enfin de vous exposer mes considérations sur le travail da ma fille; ces considérationssont les suivantes.

1) Pour considerer la question de la theorie des parallèles comme résolue définitivement, il est necessaire de constater le postulatum d'Euclide par une de deux démonstrations parmi lesquelles la démonstration directe est renfermée dans la déduction du postulatum d'Euclide des axiomes precédants à lui et acceptées par tous les Géomètres; l'autre démonstration consiste dans la deduction du postulatum d'Euclide d'autre postulatum precedant à lui et accepté pour l'axiome ou dans ce que le même postulatum d'Euclide doit être accepté pour l'axiome; mais cette l'autre démonstration ne peut avoir lieu qu'en cas où, de la manière exacte, on parviendra à la conclusion que les axiomes precedants au postulatum d'Euclide et acceptés par tous les Géomètres sont tant insuffisants qu'il est impossible d'en déduire le postulatum d'Euclide; chacune de ces deux démonstrations est le sujet à considération des Géomètres. Beltrami, comme vous ecrivez vraiment, dans son ouvrage bien connu à démontré qu'un espace, dans lequel la Géometrie de Lobatschevsky est verifiée, est materiellement possible, mais de cette démonstration de Beltrami ne descende pas la conclusion que les axiomes precedants au postulatum d'Euclide et acceptés par tous les Géomètres sont tant insuffisants qu'il est impossible d'en déduire le postulatum d'Euclide.

[187] "Sulla teoria del pendolo", *Giornale di Matematiche di Battaglini*, v. XXXVIII, pp. 1-6.
[188] Non è presente nel Fondo.

Cette conclusion ne descende pas pareillement des ouvrages de Lobatschevsky et des Géomètres dont les recherches relatives à la théorie de Lobatschevsky ont poussé cette théorie plus loin.
2) Lobatschevsky, qui adopte l'analyse infinitesimale, il me semble, sur ces pages de son ouvrages qui sont indiquées par vous, fait les objections sur la comparaison des espaces infinis, comme des espaces finis, mais non sur l'usage des espaces infinis en général; il me semble encore que Lobatschevsky fait les objections sur ce que l'on determine l'angle comme la partie indefinie du plan, comprise entre deux droites qui se coupent, mais non sur ce que l'on attribue la signification d'une aire infinie à l'espace d'un angle qui est entendu comme l'écartement de deux droits ou l'inclinaison d'une sur l'autre de deux droites qui se coupent. Cela reçu, il me semble que sans méconnaitre les valeurs des objections sur l'espace d'une aire infinie, telles que l'on trouve dans l'ouvrage de Lobatschevsky, 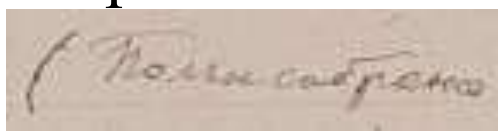 pages 222 et suivantes), on peut admettre le raisonnement, à l'aide duquel ma fille déduit les deux formules (4) et (6), à plus forte raison, que cette déduction n'est fondée, que sur l'égalité

$$(1)\qquad \frac{A}{D}=\frac{a}{d}$$

qui, comme-on le sais, est démontré avec toute la rigueur par réduction à l'absurde ou par la *[*mesure*]* des limites. An demeurant sans attribuer à l'espace angulaire la signification d'une aire infinie, on ne peut pas démontrer que tous les angles droits sont égaux, si ce n'est dans le cas où l'on adapte une de ces deux admissions, des quelles une consiste en ce que, si on fait tourner l'oblique quelconque AD autour du point A de manière que le plus petits des deux angles adjacents grandisse, l'autre diminuera d'autant, l'autre admission est celle-ci: il y a la courbe *[...]*, que la circonférence, mais chacune de deux admissions n'est *[*assume des axiomes*]* acceptés par tous les Géomètres…
3) Quand à la seconde note de ma fille, sans doute, il sera mieux qu'elle ne soit pas imprimée, car dans element de la Géométrie il est impossible de démontrer, que la limite d'une sphère lorsque son rayons crois infiniment est le plan.

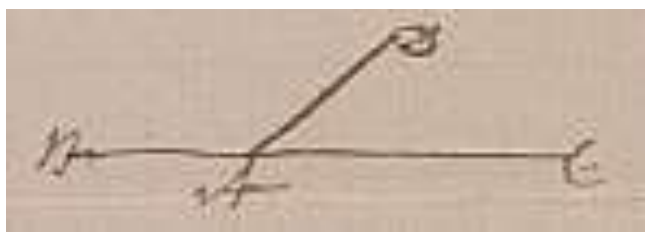

Si tout de même vous ne trouvez pas possible que la première note de ma fille soit imprimée dans quelque journal mathématique, permettez moi d'esperer, Monsieur et très honoré Collegue que vous aurez la bonté de me répondre par quelque objection.

Ma fille très sensible à votre bon souvenir charge de transmettre les salutations les plus empressées.

En attendant j'ai l'honneur d'envoyer à vous mes salutations les plus distinguées.

Votre tout devoué

G. Sabinine

Mon adresse: Russie, Mosca, Arbate, Nicolsky *[...]* (ruell), maison Sabolotsky, 28

**159**
Gian Antonio Maggi a Georg Sabinine

Pise le 27 Mai 1901

Monsieur et très-honoré confrère,

J'ai attendu quelque temps à répondre à votre lettre, en essayant de trouver une façon de mieux exprimer ma pensée, plutôt qu'en cherchant de nouvelles raisons à son appui.

Je ne puis qu'insister sur cela que j'envisage la question des parallèles comme définitivement résolue par une démonstration de cette espèce que vous appelez la deuxième. C'est bien le résultat atteint par M. Beltrami, en constant l'identité de la géométrie de Lobatschevsky et de la géométrie de la pseudosphère. En effet, après cette constatation, nous ne savons pas seulement quelles conséquences découlent de l'omission du cinquième postulat d'Euclide - c'est le point auquel Lobatschevsky s'est arrêté - mais, bien plus nous possédons une surface effective, dont nous pouvons construire un modèle (comme a fait M. Beltrami) sur laquelle, à l'aide de figures sensibles, que nous pouvons tracer et mesurer sur le modèle, nous reconnaissons que tous les postulats d'Euclide jusqu'au cinquième, non compris, y gardent leur validité, tandis qu'il s'y vérifie les propriétés lesquelles excluent la validité du cinquième postulat. Par exemple, la somme des angles internes d'un triangle rectiligne (géodesique) quelconque de la pseudosphère est moins que deux droits: ce que nous pouvons constater sur le modèle par des mesures goniométriques. D'après Lobatschevsky cela devait arriver pour les triangles rectilignes du plan, en supprimant le cinquième postulat, et réciproquement ce postulat ne peut pas se vérifier si cela a lieu; mais il restait à prouver que cela <u>pouvait</u> arriver, c'est-à-dire qu'une surface sur laquelle les premiers quatre postulats d'Euclide gardent leur validité, et cependant la somme des angles internes d'un triangle rectiligne est moindre que deux droits n'est pas absurde. C'est la possibilité d'une telle surface que l'étude de la pseudosphère a démontrée: et alors comment démontrerait on encore que la validité des quatre premiers postulats d'Euclide entraîne nécessairement que la somme des angles internes d'un triangle rectiligne est égale à deux droits?

Quant à l'usage des aires infinies, ces aires, dans les raisonnements en question, se trouvent ajoutées entr'elles et à des aires finies, et il ne s'agit pas de l'*[*accorder*]* ou le défendre <u>à priori</u>, mais de demander une définition de ce qu'il faut entendre par une aire infinie en comparaison avec une aire finie, ce qui est bien nécessaire pour en tirer des conclusions ayant une signification determinée. À l'aide de

$$\frac{\mathrm{A}}{\mathrm{D}}=\frac{\mathrm{a}}{\mathrm{d}}$$

nous pouvons définir le rapport de deux espaces angulaires par le rapports des angles correspondants. Mais qu'est ce qu'il faut entendre par le rapport $\frac{\Delta}{\mathrm{D}}$? Ce rapport doit representer un nombre, et je ne vois pas comment ce nombre pourrait être autrement fixé qu'en établissant l'une ou l'autre de ces deux choses: ou combien de mètres carrés vaut D ou à combien de degrés correspond Δ. Je ne sais pas qu'on ait jamais fait autre que considérer la mesure de D en mètres carrés comme infinie. Cela revient tout court à $\frac{\Delta}{\mathrm{D}}=0$. Mais cette supposition implique un passage à la limite dans des conditions que la question ne pourrait pas préciser: c'est precisément l'objection de Lobatschevsky:

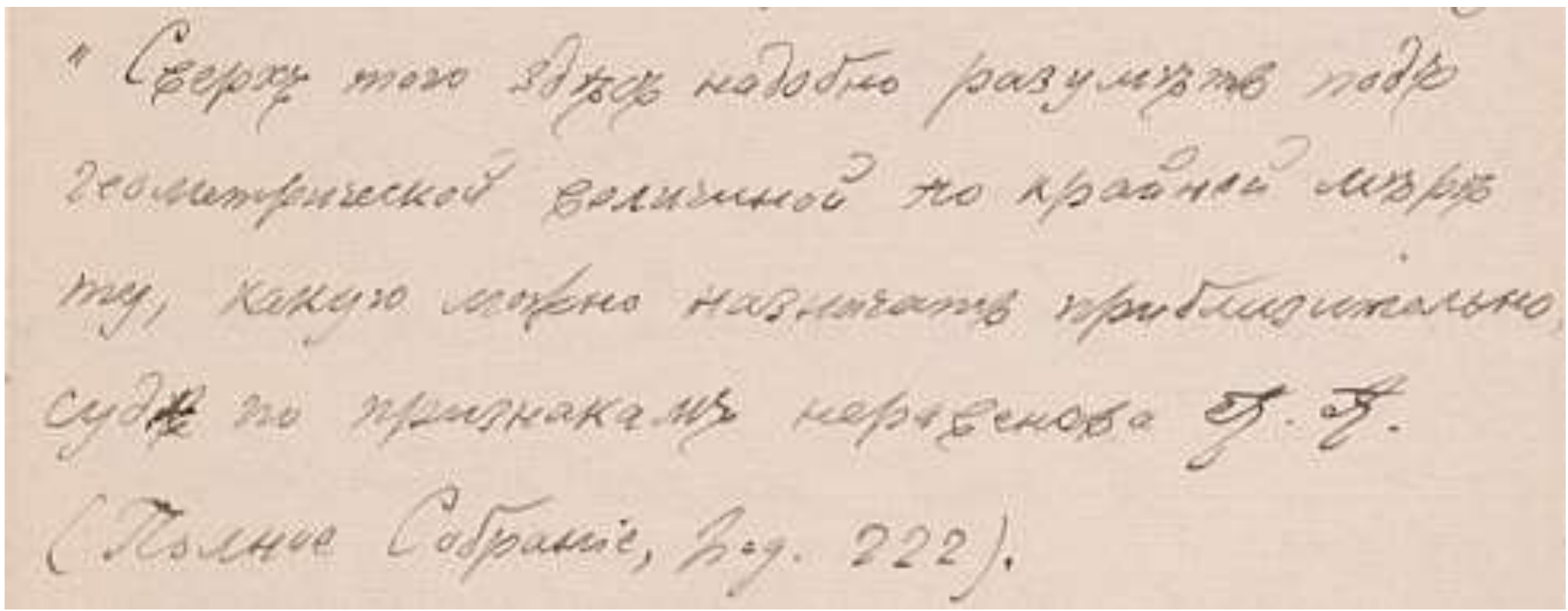

« Сверхъ того здѣсь надобно разумѣть подъ геометрической величиной по крайней мѣрѣ ту, какую можно назначать приблизительно, судя по признакамъ неравенства Н. Л.
(Полное Собраніе, стр. 222).

À l'égard enfin de la démonstration de l'égalité des angles droits, je ne pourrais que faire les mêmes objections qu'auparavant à un raisonnement conduit avec le même procédé.

Je regrette extrêmement de ne pouvoir pas proposer la publication du manuscrit que vous m'avez l'honneur de m'envoyer. J'ajouterai, à cet égard, que mon avis est celui de mes confrères, et que je ne pourrais pas attendre un avis contraire des directeurs de journaux de mathématiques avec lesquels je suis en correspondance. J'espère que vous ne me ferez pas attendre longtemps l'occasion de vous prouver tout mon désir de vous servir et de répondre à votre obligeance.

Je vous prie, en même temps, d'agréer l'assurance de ma considération bien distinguée, et de présenter à Mademoiselle votre fille mes salutations les plus respectueuses.

Votre tout dévoué
Gian Antonio Maggi.

# 160
Georg Sabinine a Gian Antonio Maggi[189]

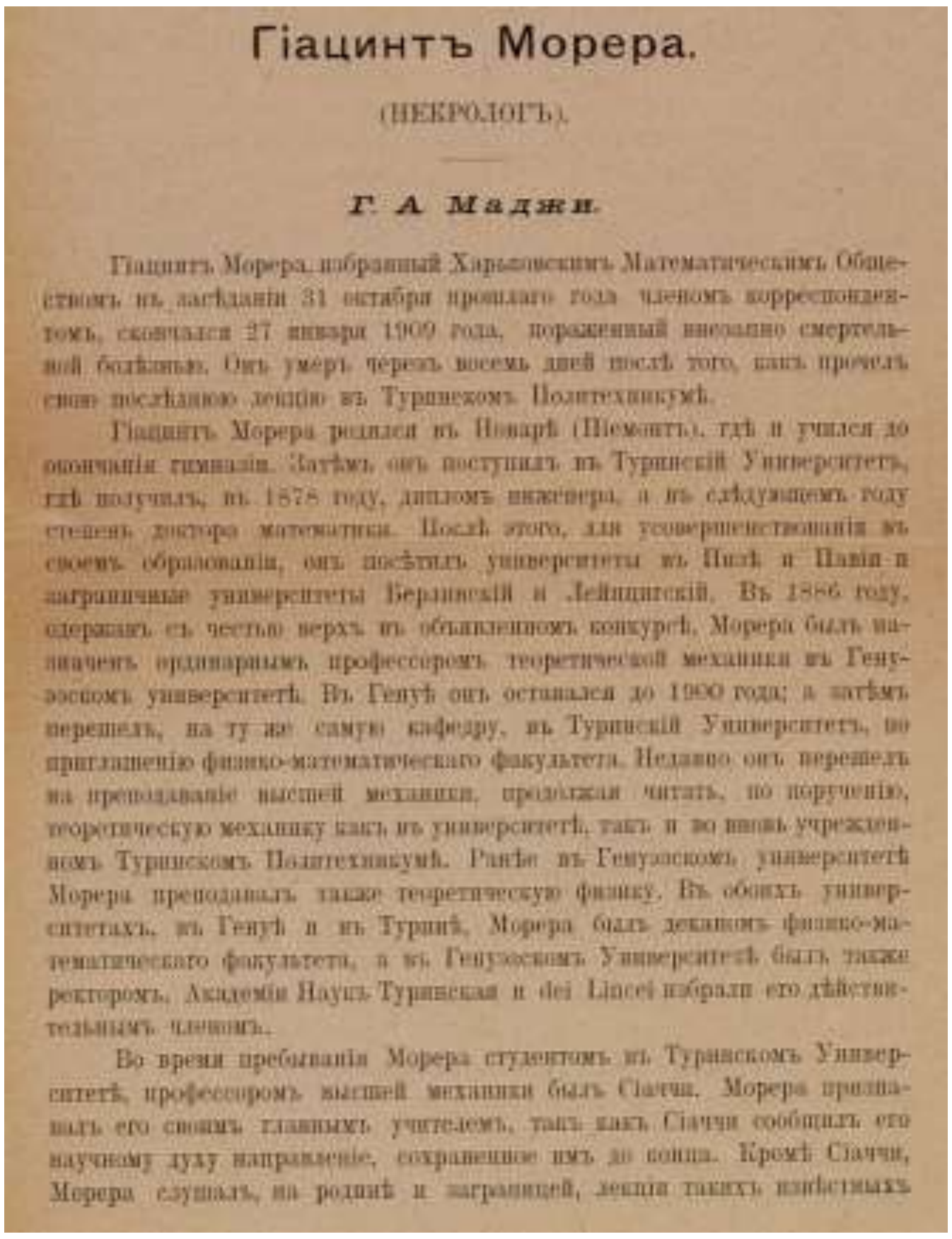

## Гіацинтъ Морера.

(НЕКРОЛОГЪ).

**Г. А. Маджи.**

Гіацинтъ Морера, избранный Харьковскимъ Математическимъ Обществомъ въ засѣданіи 31 октября прошлаго года членомъ корреспондентомъ, скончался 27 января 1909 года, пораженный внезапно смертельной болѣзнью. Онъ умеръ черезъ восемь дней послѣ того, какъ прочелъ свою послѣднюю лекцію въ Туринскомъ Политехникумѣ.

Гіацинтъ Морера родился въ Новарѣ (Піемонтъ), гдѣ и учился до окончанія гимназіи. Затѣмъ онъ поступилъ въ Туринскій Университетъ, гдѣ получилъ, въ 1878 году, дипломъ инженера, а въ слѣдующемъ году степень доктора математики. Послѣ этого, для усовершенствованія въ своемъ образованіи, онъ посѣтилъ университеты въ Пизѣ и Павіи и заграничные университеты Берлинскій и Лейпцигскій. Въ 1886 году, одержавъ съ честью верхъ въ объявленномъ конкурсѣ, Морера былъ назначенъ ординарнымъ профессоромъ теоретической механики въ Генуэзскомъ университетѣ. Въ Генуѣ онъ оставался до 1900 года; а затѣмъ перешелъ, на ту же самую кафедру, въ Туринскій Университетъ, по приглашенію физико-математическаго факультета. Недавно онъ перешелъ на преподаваніе высшей механики, продолжая читать, по порученію, теоретическую механику какъ въ университетѣ, такъ и во вновь учрежденномъ Туринскомъ Политехникумѣ. Ранѣе въ Генуэзскомъ университетѣ Морера преподавалъ также теоретическую физику. Въ обоихъ университетахъ, въ Генуѣ и въ Туринѣ, Морера былъ деканомъ физико-математическаго факультета, а въ Генуэзскомъ Университетѣ былъ также ректоромъ. Академія Наукъ Туринская и dei Lincei избрали его дѣйствительнымъ членомъ.

Во время пребыванія Морера студентомъ въ Туринскомъ Университетѣ, профессоромъ высшей механики былъ Сіаччи. Морера признавалъ его своимъ главнымъ учителемъ, такъ какъ Сіаччи сообщилъ его научному духу направленіе, сохраненное имъ до конца. Кромѣ Сіаччи, Морера слушалъ, на родинѣ и заграницей, лекціи такихъ извѣстныхъ

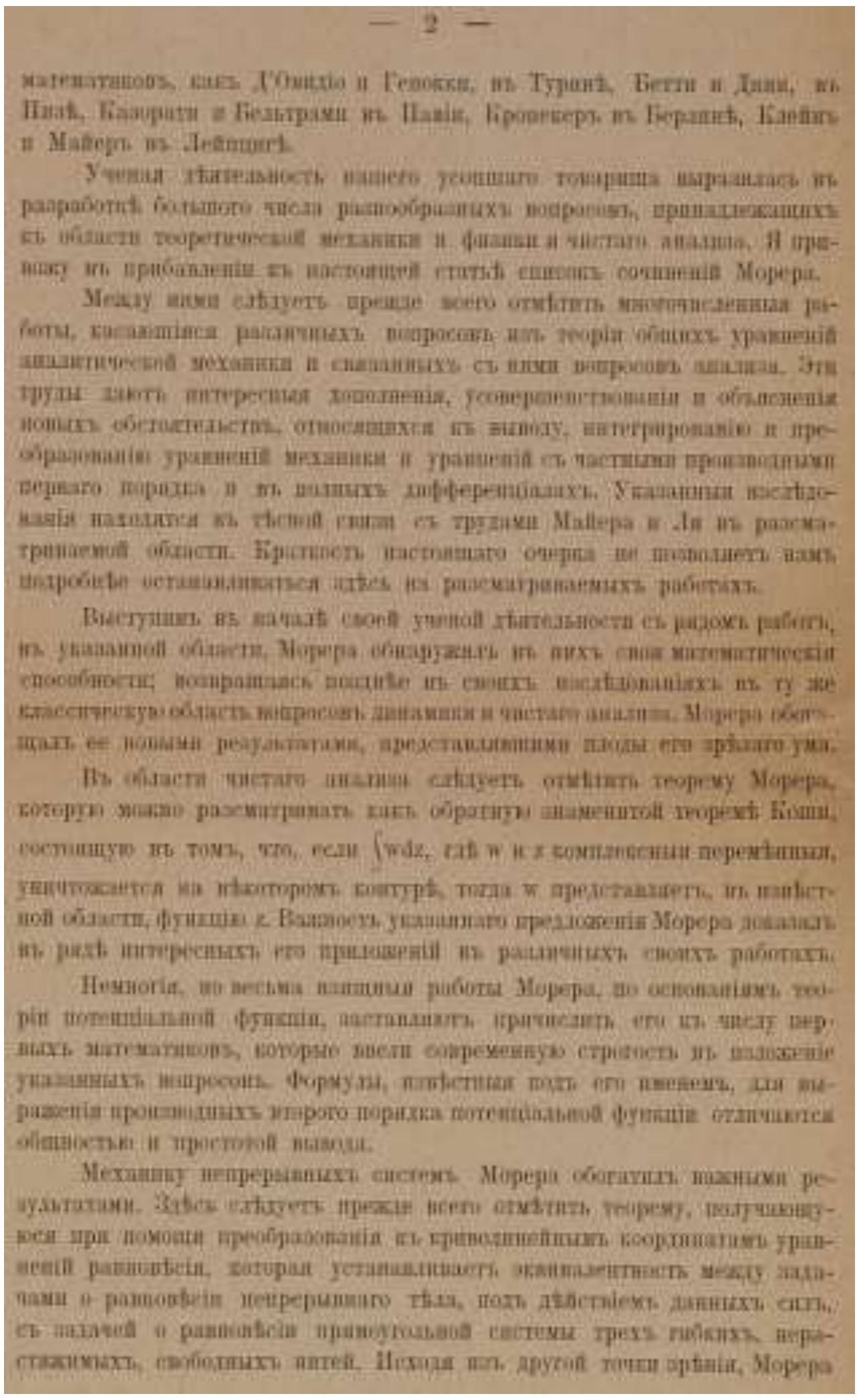

— 2 —

математиковъ, какъ Д'Овидіо и Генокки, въ Туринѣ, Бетти и Дини, въ Пизѣ, Казорати и Бельтрами въ Павіи, Кронекеръ въ Берлинѣ, Клейнъ и Майеръ въ Лейпцигѣ.

Ученая дѣятельность нашего усопшаго товарища выразилась въ разработкѣ большого числа разнообразныхъ вопросовъ, принадлежащихъ къ области теоретической механики и физики и чистаго анализа. Я прилагаю въ прибавленіи къ настоящей статьѣ списокъ сочиненій Морера.

Между ними слѣдуетъ прежде всего отмѣтить многочисленныя работы, касающіяся различныхъ вопросовъ изъ теоріи общихъ уравненій аналитической механики и связанныхъ съ ними вопросовъ анализа. Эти труды даютъ интересныя дополненія, усовершенствованія и объясненія новыхъ обстоятельствъ, относящихся къ выводу, интегрированію и преобразованію уравненій механики и уравненій съ частными производными перваго порядка и въ полныхъ дифференціалахъ. Указанныя изслѣдованія находятся въ тѣсной связи съ трудами Майера и Ли въ разсматриваемой области. Краткость настоящаго очерка не позволяетъ намъ подробнѣе останавливаться здѣсь на разсматриваемыхъ работахъ.

Выступивъ въ началѣ своей ученой дѣятельности съ рядомъ работъ, въ указанной области, Морера обнаружилъ въ нихъ свои математическія способности; возвращаясь позднѣе въ своихъ изслѣдованіяхъ въ ту же классическую область вопросовъ динамики и чистаго анализа, Морера обогащалъ ее новыми результатами, представлявшими плоды его зрѣлаго ума.

Въ области чистаго анализа слѣдуетъ отмѣтить теорему Морера, которую можно разсматривать какъ обратную знаменитой теоремѣ Коши, состоящую въ томъ, что, если $\int w dz$, гдѣ $w$ и $z$ комплексныя перемѣнныя, уничтожается на нѣкоторомъ контурѣ, тогда $w$ представляетъ, въ извѣстной области, функцію $z$. Важность указаннаго предложенія Морера доказалъ въ рядѣ интересныхъ его приложеній въ различныхъ своихъ работахъ.

Немногія, но весьма изящныя работы Морера, по основаніямъ теоріи потенціальной функціи, заставляютъ причислить его къ числу первыхъ математиковъ, которые ввели современную строгость въ изложеніе указанныхъ вопросовъ. Формулы, извѣстныя подъ его именемъ, для выраженія производныхъ втораго порядка потенціальной функціи отличаются общностью и простотой вывода.

Механику непрерывныхъ системъ Морера обогатилъ важными результатами. Здѣсь слѣдуетъ прежде всего отмѣтить теорему, получающуюся при помощи преобразованія къ криволинейнымъ координатамъ уравненій равновѣсія, которая устанавливаетъ эквивалентность между задачами о равновѣсіи непрерывнаго тѣла, подъ дѣйствіемъ данныхъ силъ, съ задачей о равновѣсіи прямоугольной системы трехъ гибкихъ, нерастяжимыхъ, свободныхъ нитей. Исходя изъ другой точки зрѣнія, Морера рѣшилъ задачу интегрированія разсматриваемыхъ уравненій равновѣсія при помощи произвольныхъ функцій. Послѣднимъ результатомъ онъ воспользовался въ одной изъ своихъ послѣднихъ статей, относительно равновѣсія упругихъ тѣлъ и получилъ интересные выводы.

[189] Sabinine invia copia dell'articolo di Maggi: "Giazint Morera. Nekrolog" (in russo), apparso sul periodico *Comm. de la Societé Mathémathique de Kharkoff*, s. II, v. XI, pp. 243-248.

Среди позднѣйшихъ работъ Мореры особенно выдѣляется его изслѣдованія о притяженіи неоднороднаго эллипсоида и эллипсоидальныхъ слоевъ, о задачѣ Дирихле въ эллипсоидальной области и объ «эллипсоидальныхъ гармоническихъ функціяхъ». Морера вводитъ послѣднія функціи, распространяя прежній полученный имъ результатъ, относящійся къ изслѣдованіямъ Пицетти о геоидѣ. Названныя функціи получаются при разсмотрѣніи потенціальной функціи нѣкоторыхъ неоднородныхъ эллипсоидовъ подобно тому, какъ, по способу Томсона и Тета, получаются сферическія гармоническія функціи, при разсмотрѣніи потенціальной функціи однородной сферы. Аналогично послѣднимъ, эллипсоидальныя функціи прилагаются къ изученію задачъ, касающихся эллипсоида вмѣсто сферы. Введенными функціями Морера предполагалъ также воспользоваться для рѣшенія задачи о равновѣсіи эллипсоидальнаго упругаго тѣла [1].

Названія остальныхъ работъ Мореры находятся въ спискѣ, приложенномъ въ концѣ настоящаго краткаго очерка.

Поскольку въ твореніи Мореры важность предмета соединяется съ талантливостью и сжатостью изложенія, часто на весьма небольшомъ числѣ страницъ. Эти качества служили отличительной чертой его ученой дѣятельности.

Съ этими проявленіями научнаго творчества вполнѣ гармонировалъ весь духовный обликъ Мореры какъ человѣка. Онъ обладалъ краткой остроумной рѣчью; не обращая вниманія на внѣшнія черты, онъ особенно цѣнилъ нравственныя качества человѣка; и самъ отличался честностью, строгостью, проницательностью, какъ въ своей научной дѣятельности, такъ и въ общественныхъ и служебныхъ отношеніяхъ. Морера былъ притомъ истиннымъ другомъ для своихъ товарищей, и любимое имъ семейство могло гордиться тѣмъ, что раздѣляло съ наукою его время и заботы. Поэтому мы оплакиваемъ, съ кончиною Мореры, потерю не только выдающагося ученаго, но также и хорошаго человѣка.

---

[1] Сомильяна, въ некрологѣ Мореры въ «Nuovo Cimento» (апрѣль), сообщаетъ, что послѣ Мореры осталась почти законченная рукопись, посвященная послѣднему вопросу.

## Списокъ сочиненій Гіацинта Мореры.

### Уравненія динамики и уравненія съ частными производными перваго порядка и въ полныхъ дифференціалахъ.

1. Sopra una formola di meccanica analitica (Rendiconti del R. Istituto Lombardo di Scienze e Lettere, 1882).

2. Teorema fondamentale nella teoria delle equazioni canoniche del moto (Ibid.).

3. Il metodo di Pfaff per l'integrazione delle equazioni a derivate parziali del 1° ordine (Ibid., 1883).

4. Sul problema di Pfaff (Atti della R. Accademia delle Scienze di Torino, 1883).

5. Ueber die Integration der vollständigen Differentiale (Mathematische Annalen, 1886).

6. Sulla integrazione delle equazioni a derivate parziali del 1° ordine (Giornale della Società di Lettere e di Conversazioni Scientifiche (Genova), 1887).

7. Sull' integrazione delle equazioni ai differenziali totali del 2° ordine (Memorie della R. Accademia delle Scienze di Torino, 1902).

8. Intorno ai sistemi di equazioni a derivate parziali del 1° ordine in involuzione (Rendic. del R. Ist. Lomb., 1903).

9. I sistemi canonici d'equazioni ai differenziali totali nella teoria dei gruppi di trasformazioni (Atti della R. Acc. delle Scienze di Torino, 1903).

10. Sulle equazioni dinamiche di Lagrange (Ibid.).

11. Sulle equazioni dinamiche di Hamilton (Ibid., 1904).

12. Sulla trasformazione delle equazioni dinamiche di Hamilton (Rendiconti della R. Accademia dei Lincei, 1905).

### Алгебраическій анализъ.

13. Sulle proprietà invariantive del sistema di una forma lineare e di una forma bilineare alternata (Atti della R. Acc. delle Scienze di Torino, 1883).

14. Un piccolo contributo alla teoria delle forme quadratiche (Rendiconti del R. Ist. Lomb., 1886).

### Теорія аналитическихъ функцій.

15. Sopra una nuova costruzione geometrica del teorema dell' addizione degli integrali ellittici (Atti della R. Acc. delle Scienze di Torino, 1880).

16. Intorno alla risoluzione di certe equazioni modulari (Rendic. del R. Ist. Lomb., 1885).

17. Ueber einige Bildungsgesetze in der Theorie der Theilung und der Transformation elliptischer Functionen (Mathem. Annalen, 1885).

18. Un teorema fondamentale nella teoria delle funzioni di una variabile complessa (Rendic. del R. Ist. Lomb., 1886).

19. Sulla rappresentazione delle funzioni di una variabile complessa per mezzo di espressioni analitiche infinite (Atti della R. Acc. delle Scienze di Torino, 1886).

20. Intorno all' integrale di Cauchy (Rendic. del R. Ist. Lomb., 1889).

21. Sulla definizione di una funzione di variabile complessa (Ibid., 1901).

### Основанія теоріи потенціальной функціи.

22. Sulle derivate seconde della funzione potenziale di spazio (Rendic. del R. Ist. Lomb., 1885).

23. Intorno alle derivate normali della funzione potenziale di superficie (Ibid.).

24. Sopra una formola di calcolo integrale (Rendiconti del Circolo Matematico di Palermo, 1897).

### Равновѣсіе непрерывныхъ системъ.

25. Sull' equilibrio delle superficie flessibili e inestendibili (Rendic. della R. Acc. dei Lincei, 1883).

26. Sulle equazioni generali per l'equilibrio dei sistemi continui a tre dimensioni (Atti della R. Acc. delle Scienze di Torino, 1884).

27. Soluzione generale delle equazioni indefinite dell' equilibrio di un corpo continuo (Rendic. della R. Acc. dei Lincei, 1892).

28. Intorno all' equilibrio dei corpi elastici isotropi (Atti della R. Acc. delle Scienze di Torino, 1907).

### Притяженіе эллипсоида и эллипсоидальныя гармоническія функціи.

29. Alcune considerazioni relative alla Nota del Prof. Pizzetti „Sull' espressione della gravità alla superficie del geoide supposto ellissoidico" (Rendic. della R. Acc. dei Lincei, 1894).

30. Sull' attrazione dell' ellissoide eterogeneo (Atti della R. Acc. delle Scienze di Torino, 1904).

31. Sull' attrazione degli strati ellissoidali e sulle funzioni armoniche ellissoidali (Ibid., 1905).

32. Sull' attrazione degli ellissoidi e sulle funzioni armoniche ellissoidali di 2ª specie (Memorie della R. Acc. delle Scienze di Torino, 1905).

33. Alcune considerazioni sulle funzioni armoniche ellissoidali (Rendic. della R. Acc. dei Lincei, 1906).

34. Sulla funzione potenziale di un doppio strato ellissoidico (Ibid., 1908).

### Вопросы изъ области теоретической механики и физики.

35. Sul moto di un punto attratto da due centri fissi colla legge di Newton (Giornale di Matematiche, 1880).

36. Sulla separazione delle variabili nelle equazioni del moto di un punto sopra una superficie (Atti della R. Acc. delle Scienze di Torino, 1881).

37. Sul problema della corda vibrante (Ibid., 1888).

38. Sui moti elicoidali dei fluidi (Rendic. della R. Acc. dei Lincei, 1889).

39. Studii di termodinamica (Ibid., 1891).

40. Sui sistemi di forze che ammettono la funzione delle forze (Rendic. del Circolo Matem. di Palermo, 1891).

41. Stabilità della configurazioni di equilibrio di un liquido in un tubo capillare di rotazione (Ibid., 1902).

42. Intorno alle oscillazioni elettriche (Nuovo Cimento, 1902).

43. Sull' espressione analitica del principio di Huyghens (Ibid., 1905).

44. Sulla teoria dell' ellissoide fluido in equilibrio di Jacobi—Lettera al prof. Volterra (Ibid., 1908).

### Вопросы изъ области анализа.

45. Osservazione relativa al resto della serie di Taylor (Rivista di Matematica (Peano), 1892).

46. Alcune considerazioni relative all' equazione differenziale $\frac{du}{dx} + \frac{dv}{dy} + \frac{dw}{dz} = 0$ (Rendic. del Circolo Matem. di Palermo, 1893).

47. Dimostrazione di una formola di calcolo integrale (Rivista di Matem., 1896).

48. Sui polinomi di Legendre (Rendic. del Circolo Matem. di Palermo, 1897).

### Прочее.

49. L'insegnamento delle Scienze Matematiche nelle Università Italiane. Discorso inaugurale per l'anno accademico 1888—89 nella Università di Genova.

50. Francesco Siacci—Commemorazione (Atti della R. Acc. delle Scienze di Torino, 1908).

---

**161**
Gian Antonio Maggi a Georg Sabinine

Pise le 23 Juin 1910

Cher Monsieur et très-honoré Collègue,

Je vous remercie et vous prie de vouloir bien présenter mes remerciements à la Société Mathématique pour les extraits de ma communication, que je viens de recevoir. Je vous suis extrêmement obligé pour les soins que vous avez pris d'améliorer mon écrit, en le conservant presque intégralement. Cependant il n'aurait pas été inutile de m'envoyer les épreuves que vous m'aviez annoncées. À'part quelque détail sur la carrière de Morera, qui n'a aucune importance, je trouve à pag. 2 ligne 26 нѣкоторыи, où j'avais écrit какой нибудь (На нѣ которомъ контурѣ substitué à по какой нибудь окрестности). Je vois bien que je me suis trompé sur l'usage de какой нибудь, en traduisant à la lettre quel qu'il soit (en italien qualsivoglia, le mot employé par Morera). Je croyais d'employer un mot plus expressif, mais ayant le même sens du каждый. J'apprends qu'il ne signifie autrement que quelque. Je me permets donc de vous prier instamment (pourvu que je ne me trompe une deuxième fois) d'avoir la complaisance d'insérer dans le prochain Numéro du Bulletin la petite поправка.

| Напечатано | Слѣдуетъ читать |
|---|---|
| 2 стр. 26 стр. на нѣкоторомъ | на каждомъ. |

Je vous demande bien pardon de la peine que je me prend la liberté de vous donner.

Agréez, Cher Monsieur, avec mes excuses et, de nouveau, avec mes remerciements, mes salutations bien empressées, et veuillez présenter à M.me Saltykoff mes hommages les plus distingués.

Votre tout dévoué
Gian Antonio Maggi.

**162**
**Luigi Sala** a Gian Antonio Maggi[190]

Egregio Signore,

Milano 15 settembre 1885

Sono anzitutto gratissimo a Lei per la cortesia che ha avuto di assumersi l'incarico di pronunciare un giudizio sulla questione di cui trattano le qui accluse carte, la quale mi pare di una certa importanza.[191]

Vedrà dalla mia lettera 2 novembre 1884[192] diretta dal Prof. Ferrini[193] che desidero altro che la verità, ma la verità dimostrata. Non si dubiti che farò sempre buon viso anche ad un giudizio contrario alla mia tesi, quando sia dimostrato conforme alla verità, la quale per me sta in cima ad ogni cosa.

Credo poi che gioverà molto alla dimostrazione della verità, in modo da poter essere agevolmente anche da me compresa, il dare una diretta risposta alle due domande, colle quali chiudo l'ultimo mio scritto del 4 agosto 1885.[194]

Coi sensi della più profonda e, ad un tempo, affettuosa considerazione ho il pregio di dirmi

Devotissimo Suo
Luigi Sala

P.S. Unisco il trattato di Termodinamica del Saint-Robert (edizione seconda)[195]

**163**
Gian Antonio Maggi a Luigi Sala

Mezzago $\frac{1}{10}$85.[196]

Pregiatissimo Signore.

Ho esaminato la questione, veramente assai interessante, su cui Ella ha voluto farmi l'onore di domandare il mio parere; e mi pregio comunicarLe le seguenti considerazioni.

1. Innanzi tutto gioverà ben tener presente che, in termodinamica, una macchina calorica si concepisce come essenzialmente costituita da un corpo, il quale compie un ciclo di trasformazioni. Ogni trasformazione implica generalmente una trasmissione di calore fra il corpo considerato ed altri corpi capaci di cedergli calore (sorgenti) o di assorbirne

---

[190] La corrispondenza con il Commendator Luigi Sala è racchiusa in un fascicolo intitolato: *"Sul Coefficiente Economico delle Macchine termiche." Giudizio sopra una questione proposta dal Comm. Luigi Sala di Milano. Mezzago - Settembre - 1885.*
[191] Non sono qui conservate.
[192] Non è presente nel Fondo.
[193] Potrebbe trattarsi di Rinaldo Ferrini.
[194] Non è presente nel Fondo.
[195] P. de Saint-Robert, *Principes de Thermodynamique*, (II ed.) Loecher, Torino e Firenze, 1870.
[196] Probabilmente a questa lettera è seguita una risposta di Sala in data 31/10/1885, che non è presente nel Fondo. Infatti, la lettera #229 (priva di destinatario), per l'argomento trattato, potrebbe far parte di questa corrispondenza.

(refrigeranti), per effetto del quale varia il calore interno del corpo, e un lavoro dinamico viene sviluppato o consumato. Compiuto il ciclo, il calore interno del corpo riprende il valore iniziale; la somma degli aumenti di calore interno dev'essere compensata da una somma eguale di diminuzioni. Invece il lavoro dinamico sviluppato potrà superare il lavoro assorbito. In questo caso, l'eccesso del primo lavoro sul secondo, rappresenta un lavoro prodotto per effetto del ciclo; e, pel principio dell'equivalenza del lavoro e del calore, questo lavoro sarà prodotto a spese di un'equivalente quantità di calore, la quale non potrà essere che la differenza fra il calore attinto dalle sorgenti, e il calore ceduto ai refrigeranti, durante il ciclo. Perciò quando il corpo, che compie il ciclo, si considera come una macchina, che produce lavoro, a spese di calore, questa quantità di calore sarà quella che la macchina utilizza, in ogni ciclo.

2. Ciò posto, il Saint Robert, quando, a pag. 91 del suo libro, afferma che una macchina termica, che funziona fra i limiti di temperatura (assoluta)

$$T=150°+274°=424$$
$$t=50°+274°=324,$$

per quanto sia perfetta, non potrà utilizzare più di $\frac{150-50}{424}=\frac{100}{424}$ del calore fornito dal combustibile, intende che, supposta 150° Centigr. la temperatura della sorgente, o la più alta temperatura che presenta una sorgente, e 50° la temperatura del refrigerante, o la più bassa temperatura che presenta un refrigerante, per effetto di un ciclo, non si potrà convertire in lavoro più di $\frac{100}{424}$ del calore fornito, durante il ciclo, dalle sorgenti.

Questa proposizione non implica necessariamente che il calore fornito dalle sorgenti elevi la temperatura del corpo; quindi non potrebbe essere subordinata alla condizione che le sorgenti riscaldassero il corpo, a partire dalla temperatura di −274°.

3. I principii di termodinamica sui quali è fondata l'affermazione in discorso si possono riassumere nelle tre proposizioni seguenti:

Proposiz. 1ª. Il rapporto della quantità di calore convertita in lavoro, alla quantità di calore attinta dalle sorgenti, in un ciclo qualsiasi, non potrà ricevere un valore maggiore di quello che riceve in un ciclo di Carnot (Saint Robert pag. 81), dove la sorgente *sia quella fra le sorgenti del ciclo considerato, che ha la più alta temperatura, e il refrigerante sia quello, che ha la temperatura più bassa* *[Maggi scrive in alternativa, da * a *:* sia mantenuta alla più alta temperatura che presenta una sorgente nel ciclo considerato, e il refrigerante sia mantenuto alla più bassa temp., che presenta un refrigerante.*]*

Proposiz. 2ª. Il rapporto della quantità di calore convertita in lavoro, alla quantità di calore attinta dalla sorgente, in un ciclo di Carnot, dipende unicamente dalla temperatura della sorgente e da quella del refrigerante; fissate queste temperature, riceve lo stesso valore, qualunque sia la natura del corpo, che compie il ciclo.

Proposiz. 3ª. Supposto che un gas perfetto compia un ciclo di Carnot, dalle leggi di Mariotte e di Gay Lussac, e dal principio dell'equivalenza del calore e del lavoro, si deduce che il rapporto della quantità di calore convertita in lavoro, per effetto del ciclo, alla quantità di calore fornita al gas, durante il ciclo, dalla sorgente, è eguale al rapporto della differenza fra la temperatura della sorgente e quella del refrigerante alla temperatura della sorgente, aumentata di una temperatura espressa dall'inversa del coefficiente di dilatazione termica dei gas perfetti.

Misurando la temperatura, nel modo ordinario, in gradi centigradi, contati a partire dalla temperatura del ghiaccio fondente, il coefficiente di dilatazione termica dei gas perfetti risulta prossimamente $\frac{1}{274}$.

Perciò, adottando questo valore come esatto, e rappresentando $T_c$ e $t_c$, le temperature della sorgente e del refrigerante, misurate nel modo ordinario, e con Q, q le quantità di calore attinte dalla sorgente, e versate nel refrigerante, durante un ciclo, in virtù della precedente proposizione, sarà

$$\frac{Q-q}{Q}=\frac{T_c-t_c}{T_c+274} \qquad (1).$$

4. La questione principale concerne appunto questa relazione. Ma, prima di venire ad essa, trovo opportuno di fare la seguente osservazione.
V.S. dichiara nell'Avvertenza finale alla Sua Nota in data 4 Agosto p.p. che, per l'applicazione rigorosa della relazione in discorso, si dovrà intendere che, nella macchina calorica agisca, come corpo intermediario, un gas perfetto; e aggiunge che il Saint Robert a pag. 89, 90, 91 applica all'ingrosso alle macchine a vapore le formole desunte dall'equazione dei gas perfetti (che riunisce le leggi di Mariotte e Gay Lussac). Ora sta bene che la relazione ( ) *[sic!]* si deduca direttamente, valendosi dell'ipotesi particolare che il ciclo di Carnot sia compiuto da un gas perfetto; ma Ella vede che, una volta che noi l'ammettiamo per questo caso particolare, per la proposizione 2$^a$, la relazione stessa si applicherà rigorosamente al caso che il ciclo di Carnot sia compiuto da un corpo di qualsiasi altra specie; mentre, per la proposizione 1$^a$, se ne deduce che $\frac{T_c-t_c}{T_c+274}$ è rigorosamente il limite massimo, che può raggiungere il rapporto del calore convertito in lavoro, per effetto di un ciclo, al calore attinto dalla sorgente, durante il ciclo, quando $T_c$ rappresenta la temperatura della sorgente, o della sorgente la cui temperatura è la più alta, se ve ne fossero parecchie, e $t_c$, la temperatura del refrigerante, o del refrigerante a temperatura più bassa.

Perciò, una volta ammesso che quando un gas perfetto compie un ciclo di Carnot, attingendo calore alla temperatura di 150° Cent. e cedendone a 50° Cent., il rapporto del calore convertito in lavoro al calore attinto dalla sorgente sia $\frac{150-50}{424}=\frac{100}{424}$, in virtù di quella proposizione, si potrà rigorosamente affermare che una macchina a vapore, dove la caldaja sia mantenuta a 150°, e il condensatore a 50°, non potrà utilizzare, e cioè convertire in lavoro, per effetto d'ogni ciclo, più di $\frac{100}{424}$ del calore fornito dal combustibile, durante il ciclo; né farà alcun ostacolo che, nella macchina a vapore, il corpo intermediario sia una massa d'acqua, che non si trova mai allo stato di gas perfetto, e che il ciclo da essa compiuto non sia un ciclo di Carnot.

5. Se non che quella è la proposizione da Lei impugnata; poiché pare a V.S. che, facendo il rapporto di 150−50 a 150+274, si confronti il calore convertito in lavoro, anziché col calore fornito dalla sorgente, con questa quantità di calore aumentata del calore preesistente nel corpo intermediario; per modo che $\frac{150-50}{150+274}$ non possa rappresentare il calore utilizzato, in parti del calore attinto dalla sorgente, se non nel caso che il calore preesistente nel corpo sia nullo, ossia che il corpo cominci ad attingere calore allo zero assoluto, ipotesi lecita in teoria, ma inammissibile in pratica. Questa è la questione principale.

Ora, supponiamo che un gas perfetto compia un ciclo di Carnot, attingendo calore da una sorgente mantenuta a 150° Cent. e cedendone ad un refrigerante a 50°; e proponiamoci di determinare la quantità di calore Q, che attingerà dalla sorgente, e la quantità q che verserà nel refrigerante, ad ogni ciclo, donde dedurremo la quantità Q−q convertita in lavoro, per confrontarla con Q.

Per la definizione del ciclo di Carnot, il gas, preso alla temperatura di 50° del refrigerante, si porterà alla temperatura di 150° a cui è mantenuta la sorgente, mediante una semplice compressione, applicandovi ad ogni istante una pressione infinitamente poco superiore alla sua, senza comunicarvi, né sottrarne calore; quindi, si lascerà espandere, opponendogli ad ogni istante una pressione infinitamente poco inferiore alla sua, a contatto della sorgente, la quale gli cederà calore, in tal misura che, durante l’espansione, la sua temperatura non si abbasserà, e neppure si eleverà, ma si manterrà invariabilmente eguale a 150°; poi, separato dalla sorgente, lasciandolo espandere ulteriormente, allo stesso modo, senza comunicarvi né sottrarne calore, si farà raffreddare fino alla temperatura iniziale di 50°; finalmente, si comprimerà, finché il suo volume riprenda il valore iniziale, applicandovi una pressione, ad ogni istante, infinitamente poco superiore alla sua, e mantenendolo a contatto del refrigerante, il quale ne sottrarrà il calore sviluppato dalla compressione, in tal maniera che la sua temperatura si conserverà invariabilmente di 50°. Allora, la temperatura, e il volume del gas, avendo ripreso il valor iniziale, l’avrà ripreso anche la pressione; quindi il gas sarà ritornato nello stato iniziale, e il ciclo sarà chiuso.

Il gas attingerà la quantità di calore Q durante la seconda operazione, e cederà la quantità q durante la quarta.

Osserviamo che mentre il gas attinge dalla sorg la quantità Q di calore, la sua temperatura, previamente elevata a 150°, si mantiene invariabilmente a 150°.

Ora il calore interno di un gas perfetto non dipende che dalla temperatura. Quindi mentre il gas assorbe la quantità di calore Q il suo calore interno non varia; e perciò da una parte la quantità di calore di Q, deve spendersi tutta quanta nel produrre il lavoro sviluppato dall’espansione corrispondente; dall’altra questo lavoro deve prodursi a sole spese di questa quantità di calore. Segue da ciò che, rappresentando con $v_2$, $v_3$ il volume del gas al principio e alla fine dell’espansione, con v, p il volume e la pressione ad un istante qualunque dell’espansione medesima, con dv l’aumento di volume corrispondente a un decremento infinitesimo di p, e finalmente con J l’equivalente meccanico del calore, sarà

$$Q = \frac{1}{J}\int_{v_2}^{v_3} p\,dv;$$

relazione donde ci vagliamo, per calcolare Q.

Per effettuare questo calcolo; rammentiamo che, rappresentando con $t_c$ una temperatura qualsiasi espressa in gradi centigradi nel modo ordinario, la pressione p che possiede un gas, alla temperatura $t_c$ e sotto il volume v è data in termini di $t_c$ e di v dalla relazione

$$p = R\frac{t_c + 274}{v},$$

dove R è una costante.

Supposto che p rappresenti la pressione del gas ad un istante qualunque, durante l'espansione considerata, siccome, durante l'espansione medesima, la temperatura del gas si mantiene costantemente a 150°, sarà $t_c$=150°; e perciò la pressione in discorso sarà data in termini del volume corrispondente della relazione

$$p=\frac{R(150+274)}{v}.$$

Introducendo questa espressione di p nella ( ), ne ricaviamo

$$Q=\frac{1}{J}\int_{v_2}^{v_3}\frac{R(150+274)}{v}dv=\frac{R(150+274)}{J}\int_{v_2}^{v_3}\frac{dv}{v},$$

donde perché

$$\int_{v_2}^{v_3}\frac{dv}{v}=\log v_3-\log v_2=\log\frac{v_3}{v_2},$$

Si conclude

$$Q=R(150°+274°)\log\frac{v_2}{v_3}.$$

Lo stesso ragionamento ci serve per trovare q. Basta che osserviamo che q è la quantità di calore che il gas passando dallo stato iniziale a quello in cui si trova al principio dalla quarta operazione, o alla fine della terza, a contatto di una sorgente mantenuta a 50°, dovrebbe attingere dalla sorgente medesima, perché la sua temperatura si conservasse costantemente eguale a 50°.

Ciò posto, rappresentando con $v_1$, $v_4$ il volume del gas nello stato iniziale, e alla fine della terza operazione, dal risultato precedentemente ottenuto concludiamo immediatamente

$$q=R(50°+274°)\log\frac{v_4}{v_1}.$$

Ora si dimostra facilmente che sarà

$$\frac{v_3}{v_2}=\frac{v_4}{v_1}.$$

Difatti, la pressione p e il volume v che presenta ad ogni istante un gas, mentre si espande o si comprime senza comunicazione, né dispersione di calore, sono collegati dall'equazione

$$pv^{\alpha}=\text{Cost.}$$

dove α è un numero.

Perciò si ha

$$p_1v_1^{\alpha}=p_2v_2^{\alpha},\qquad p_3v_3^{\alpha}=p_4v_4^{\alpha}$$

$$\frac{p_1}{p_2}=\left(\frac{v_2}{v_1}\right)^{\alpha},\qquad \frac{p_4}{p_3}=\left(\frac{v_3}{v_4}\right)^{\alpha}.$$

D'altra parte si ha

$$\frac{p_1}{p_2}\frac{50°+274°}{150°+274°}=\frac{v_2}{v_1},\qquad \frac{p_4}{p_3}\frac{50°+274°}{150°+274°}=\frac{v_3}{v_4};$$

e dividendo queste equazioni per le precedenti, se ne ricava

$$\left(\frac{v_2}{v_1}\right)^{1-\alpha}=\left(\frac{v_3}{v_4}\right)^{1-\alpha}=\frac{50+274}{150+274},$$

donde

$$\frac{v_2}{v_1}=\frac{v_3}{v_4}$$

e

$$\frac{v_3}{v_2}=\frac{v_4}{v_{14}}$$

c.v.d.

Posto allora $\frac{R}{J}\log\frac{v_3}{v_2}=\frac{R}{J}\log\frac{v_4}{v_1}=K$; per modo che K è una quant. che dipende dall'esposizione che subisce il gas, mentre attinge la qu. di cal. Q per le ( ), ( ) si ha

$$Q=K(150°+274°) \qquad q=K(50°+274°)$$

donde

$$Q-q=K(150°-50°)$$

e finalmente

$$\frac{Q-q}{Q}=\frac{150-50}{150+274}.$$

Così, cercando il rapporto del calore convertito in lavoro, per effetto d'un ciclo, al calore fornito dalla sorgente, nell'ipotesi che un gas perfetto compia un ciclo di Carnot, attingendo calore a 150°, e cedendone a 50°, noi arriviamo a quella conclusione, donde, per quanto abbiamo precedentemente veduto, valendosi delle proposizioni 1[a] e 2[a] (§), si deduce direttamente la proposizione affermata dal S−R.

Questo è quello che Ella chiama il problema del S−R. Ella vede che questo problema acclude necessariamente il limite che la quantità di calore Q, colla quale si confronta la quantità di calore convertita in lavoro, sia fornita al corpo intermediario dalla sorgente, e che di questo limite si tiene perfettamente calcolo nella soluzione. Difatti si determina Q, dividendo per J il lavoro compiuto dall'espansione del gas nella seconda operazione del ciclo; il qual lavoro, siccome durante quell'espansione il calore interno del gas si mantiene costante, dev'essere fatto a sole spese del calore fornito dalla sorgente, durante l'espansione medesima, e, assorbito tutto quanto, risultando così esattamente equivalente ad esso.

Notiamo bene che la quantità di calore che la sorgente fornisce al gas si spende tutta quanta nel lavoro d'espansione, e non eleva la temperatura del gas.

6. La circostanza che, nel nostro problema, il corpo si riscalda, senza comunicazione di calore, finché raggiunge la temperatura della sorgente, e, posto soltanto allora in comunicazione colla sorgente, ne attinge calore, mentre la sua temperatura si mantiene costante, è di capitale importanza, nella presente questione. Ella vede che, una volta stabilito questo punto, non si può parlare di un caso, in cui il corpo attinga calore, partendo da –274°, nel quale sarebbe lecito l'uso delle temperature computate da –274°, e di un caso, in cui il corpo attinga calore partendo da 0°, o da un'altra temperatura prefissata, nel quale alle temperature assolute, si dovrebbero sostituire le temperature contate a partire da 0°, o da questa temperatura: quella temperatura che Ella chiama relativa, e che distingue coll'indice r.

Sta bene che in una macchina dove la sorgente sia mantenuta a 150°, e il refrigerante a 50°, si possa immaginare che il corpo intermediario sia riscaldato dalla sorgente, a partire da una temperatura qualsiasi. Ma il nostro problema ha per oggetto di trovare il limite superiore del rapporto fra il calore convertito in lavoro e il calore attinto dalla sorgente, durante un ciclo; e a tal fine dobbiamo trovare il valore che riceve questo rapporto nell'ipotesi che il corpo intermediario compia un ciclo di Carnot. Allora, poiché la sorgente è mantenuta a 150°, bisognerà supporre che il corpo, previamente elevato a questa temperatura, senza comunicazione di calore, attinge calore, espandendosi, mentre la sua temperatura si mantiene invariabilmente a 150°. Facendo un'ipotesi diversa, e cioè supponendo che il corpo intermediario sia posto in comunicazione colla sorgente ad una temperatura inferiore, il ciclo non sarà più invertibile, e il rapporto fra il calore utilizzato, e il calore attinto dalla sorgente, in un ciclo, non potrà raggiungere il limite massimo. V.S., per determinare il rapporto in discorso, suppone che una massa di gas, presa a 0°, attingendo calore dalla sorgente, si riscaldi fino a 150°, e trova per valore di quel rapporto $\frac{2}{3} > \frac{100}{424}$, valore corrispondente al ciclo di Carnot. Ma, come mi permetterò di farLe osservare fra poco, le operazioni da Lei considerate non riconducono il gas alle condizioni iniziali, e perciò non costituiscono un ciclo.

Poiché Q è una quantità di calore che il corpo intermediario attinge alla temperatura della sorgente calorifica, la relazione

$$\frac{Q-q}{Q} = \frac{T_c - t_c}{T_c + 274}$$

non potrebbe implicare la condizione che il corpo l'attinga, a partire dallo zero assoluto, senza essere assurda. E poiché, come abbiamo veduto, essa è una legittima conseguenza del principio dell'eguaglianza del calore e del lavoro, ciò non si potrebbe ammettere, senza impugnare questo principio, che è la pietra angolare della termodinamica?

7. Ma perché la relazione in discorso dovrebbe implicare quell'ipotesi? Ella afferma che la relazione

$$\frac{Q}{q} = \frac{H}{h} = \frac{T_c + 274}{t_c + 274},$$

dove H, h rappresentano il calore interno che contiene il gas, rispettivamente a $T_c$ e a $t_c$, si può considerare come l'equaz. fondam. della Sua dimostrazione.

Esaminiamo separatamente le tre proposizioni

$$\frac{Q}{q} = \frac{T_c + 274}{t_c + 274}, \qquad \frac{H}{h} = \frac{T_c + 274}{t_c + 274}, \qquad \frac{Q}{q} = \frac{H}{h}.$$

Dalla prima proposizione, che in sostanza non è altro che la ( ), segue

$$\frac{Q}{T_c + 274} = \frac{q}{t_c + 274},$$

donde, rappresentando con K il valor comune dei due rapporti

$$Q=K(T_c+274). \qquad (3)$$

Q rappresenta la quantità di calore che il corpo considerato attinge dalla sorgente. Ora dalla (3), ammesso che Q rappresenti una quantità di calore, che, attinta dal gas, eleva la temperatura da un certo valore iniziale fino a $T_c$, segue senz'alcun dubbio che questa temperatura iniziale non potrebbe essere che –274°, perché per $T_c$=–274°, se ne deduce Q=0. Ma, nulla ci obbliga a fare questa ipotesi, per la sola ragione che la quantità di calore Q è collegata colla temperatura $T_c$ da quella relazione. Nel caso nostro, ben diversamente, Q è una quantità di calore che il gas assorbe, mentre la sua temperatura si mantiene costantemente eguale a $T_c$.

Non vedrei quali prove si potrebbero ricavare dalla seconda preposizione. Noi immaginiamo di prendere il gas alla temperatura iniziale $t_c$, cioè con una certa quantità di calore interno h, e non concerne il nostro problema il modo in cui avrà assorbito questa quantità di calore.

Finalmente, consideriamo la terza proporzione. Essa esprime il vincolo che collega le quantità di calore Q, q che il gas attinge dalla sorgente e versa nel condensatore, colla quantità iniziale di calore interna h e col calore interno H, che contiene il gas, in seguito alla compr. adiabatica, a cui si assoggetta nella prima operazione del ciclo. V.S. insiste sulle conseguenze di questo vincolo; ma per dedurne che il gas dovrà attingere calore a partire da –274°, occorrerebbe ch'esso ci obbligasse ad attribuire a Q il significato di una quantità di calore che il gas attinge mentre si eleva la sua temperatura, o ad h quello di una quantità di calore che il gas deve attingere dalla sorgente, durante il ciclo.

Se non che da questa proporzione si ricava la relazione

$$\frac{Q-q}{Q}=\frac{H-h}{H}$$

che in sostanza è l'origine della controversia. $\frac{H-h}{H}$ è il rapporto della quantità di calore consumato nell'espansione adiabatica, alla quantità di calore interno, che il gas possiede prima di questa espansione. Facendo questo rapporto, non v'ha dubbio, noi confrontiamo una quantità di calore convertito in lavoro, con una quantità di calore, che si compone di una parte generata a nostre spese, e di una parte preesistente. Ma, badiamo bene, la quantità di calore consumata nell'espansione adiabatica non è la quantità di calore convertita in lavoro, per effetto dell'intero ciclo; e perciò nulla toglie che dividendola per la quantità totale di calore interno che il gas possiede prima dell'espansione, si possa ottenere lo stesso risultato come dividendo la quantità di calore convertita in lavoro, per effetto del ciclo - vale a dire la differenza fra la quantità di calore attinta dalla sorgente, e la quantità versata nel refrigerante - per la quantità di calore attinta dalla sorgente.

Osserviamo altresì che la quantità di calore interno, che il gas possiede prima dell'espansione adiabatica, in un ciclo di Carnot, è la somma del calore interno preesistente h, e del calore sviluppato dalla compressione adiabatica. La sorgente non vi presta alcun contributo. Inoltre, essendo

$$Q–q=\frac{R}{c}(H-h)\log\frac{v_3}{v_2},$$

dove $\frac{v_3}{v_2}$ si può prendere ad arbitrio, la quantità di calore Q–q convertita in lavoro, per effetto del ciclo, fissata la quantità di calore H–h consumata dall'espansione adiabatica, potrà ricevere ancora un valore qualsivoglia.

8. Ritornando all'esempio numerico di una macchina che funzione fra i limiti di temperatura 50° e 150°, ho già osservato che, quando V.S. suppone che il gas sia riscaldato per effetto della sorgente da 0° a 150°, esclude che la macchina funzioni secondo un ciclo di Carnot. Ella poi immagina che il gas si espanda adiabaticamente, finché la sua temperatura si riduce a 50°, e ammesso che la quantità di calore fornito dalla sorgente stia alla quantità consumata da questa espansione come 150 a 150–50=100, ciò che torna a supporre che il gas si riscaldi a volume costante, in modo che tutto il calore, che attinge, s'impieghi ad aumentare il calore interno, conclude che il rapporto del calore utilizzato al calore fornito dalla sorgente sarà $\frac{100}{150}=\frac{2}{3}$. Questo sarà senza dubbio il rapporto del calore convertito in lavoro dalla macchina al calore fornito dalla sorgente, supposto che il processo si arresti all'espansione adiabatica; ma, in tal caso, il gas non compie un ciclo; poiché, in seguito a quell'espansione, il gas si troverà ad un volume maggiore dell'iniziale, e alla temperatura di 50° maggiore della temperatura iniziale di 0°. Per compire il ciclo dovremo aggiungere altre operazioni, che riconducano il volume e la temperatura al valore iniziale, e a tal fine bisognerà spendere lavoro, e disperdere calore. Per esempio, si potrà comprimere il gas, a contatto di un refrigerante mantenuto a 50°, capace di assorbire il calore sviluppato dalla compressione, in modo che la sua temperatura si mantenga invariabilmente a 50°, finché si riduce a quel volume minore del volume iniziale, al quale si riduce, prendendolo nello stato iniziale, e comprimendolo adiabaticamente, finché la sua temperatura da 0° sale a 50°; raggiunto che abbia quel volume, lasciandolo espandere adiabaticamente finché la sua temperatura si riduca a 0°, è chiaro che riprenderà lo stato iniziale, per modo che il ciclo sarà compiuto.

Allora rappresentando con q la quantità di calore versata nel condensatore, e con c il calore specifico dei gas volume costante (supposta la massa del gas uguale all'unità) la quantità di calore convertita in lavoro per effetto del ciclo sarà C·150–q, e secondo i principii della termodinamica precedentemente esposti, dovrà essere $\frac{C\cdot 150-q}{C\cdot 150}<\frac{150-50}{150+274}=\frac{100}{424}$. Io, facendo il calcolo, ho trovato $\frac{C\cdot 150-q}{C\cdot 150}<0{,}136$.

9. La risposta alle domande colle quali Ella chiude la Sua Nota in data 4 Agosto p.p scaturisce spontaneamente delle precedenti conclusioni.

Ella domanda, in primo luogo "Come può il calor naturale ancor sussistente dopo il salto aver generato lavoro, quando si sa che il lavoro non può essere prodotto che dalle calorie scomparse durante il salto di temperatura, appunto perché trasformate in lavoro?"

A questo proposito, osserviamo, che quando un gas percorre un ciclo di Carnot, il calore interno ch'essa contiene al principio del ciclo non genera lavoro; in questo senso che, siccome il gas non contiene mai una quantità di calore interno minore della quantità iniziale, nessuna parte di essa scompare per convertirsi in lavoro. Ciò non toglie però che la grandezza del calore interno iniziale debba concorrere a determinare quella del lavoro che corrisponde alle singole trasformazioni del ciclo, come pure quella del lavoro risultante prodotto da un intero ciclo. Questa non sarebbe certamente una ragione perché, trattandosi di misurare il coefficiente economico di una macchina calorica, si confrontasse il calore utilizzato col calore fornito dalla sorgente aumentato del calore preesistente nel corpo intermediario, Ma, ripetiamo, siffatto errore non si commette, quando, per trovare il rapporto del calore utilizzato al calore fornito dalla sorgente, in un ciclo di Carnot, dividiamo H–h per H, poiché H è il calore interno originario aumentato del calore sviluppato dalla compressione adiabatica, e, ciò che importa ancor più di ben osservare, H–h, quantità di calore che scompare durante il salto di temperatura corrispondente all'espansione adiabatica, non è il calore utilizzato dal ciclo. H–h è il calore che si converte nel lavoro prodotto dall'espansione adiabatica; invece il calore utilizzato del ciclo è la quantità di calore che si converte nel lavoro risultante prodotto dall'intero ciclo.

Finalmente alla domanda 2ª. "Perché non si potranno valutare distintamente gli effetti dinamici di una quantità di calore artificialmente aggiunte al preesistente calor naturale dei corpi intermediarii operanti nelle macchine termodinamiche? o, in altri termini, ... perché non si potrà calcolare il lavoro entro i limiti del calore stato aggiunto al calor naturale preesistente nel corpo intermediario di una macchina calorica?" risponderò che, secondo i principii della termodinamica, supposto che la macchina considerata funzioni secondo un ciclo di Carnot, rappresentando con $T_c$, $t_c$ la temperatura centigrada, computata nel modo ordinario, della sorgente calorifica e del refrigerante, e con J l'equivalente meccanico del calore, la macchina convertirà in lavoro, per effetto d'ogni ciclo, la frazione $\frac{T_c - t_c}{T_c + 274}$ del calore che il corpo intermediario, durante il ciclo, attinge dalla sorgente, per modo che $J\frac{T_c - t_c}{T_c + 274}$ sarà il lavoro prodotto, per effetto d'un ciclo, da ogni caloria attinta dalla sorgente; donde segue, per medesimi principii, che, qualunque sia il ciclo secondo cui funziona una macchina, rappresentando con $T_c$, $t_c$ la massima temperatura d'una sorgente, e la minima d'un refrigerante, la frazione del calore <u>ceduto</u> <u>dalla</u> <u>sorgente</u> durante un ciclo, convertita in lavoro, per effetto del ciclo, non potrà superare $\frac{T_c - t_c}{T_c + 274}$, e il lavoro prodotto, per effetto d'ogni ciclo, da ogni caloria fornita dalla sorgente, durante un ciclo, non potrà superare $J\frac{T_c - t_c}{T_c + 274}$.

Qui finalmente faccio punto; ed è tempo, poiché Ella mi potrà forse fare l'appunto di non aver saputo essere più breve. Ella veda che non ho potuto dare un giudizio conforme al Suo; ma non ho esitato per questo ad esporlo francamente, poiché sono troppo persuaso che la verità sta per Lei in cima ad ogni cosa (la profonda stima che nutro per Lei, La dispensava senza dubbio dal farmi la Sua professione di fede), per temere ch'Ella se ne avrebbe a male. Io non so se avrò dimostrato abbastanza chiaramente il mio assunto; ma se Le occorresse qualche spiegazione, Ella sa che sarà sempre per me un vero piacere il poterLa servire. Ad ogni modo Le sarò obbligato se vorrà comunicarmi le Sue conchiusioni. Di nuovo La ringrazio di avermi chiamato a partecipare a questa dotta controversia, e colgo l'occasione per riverirLa distintamente e per affermarmi

Dev[mo] Suo
G.A.M.

Determiniamo il coefficiente economico di una macchina che funziona secondo un ciclo di questa specie. Il gas non assorbe calore che durante la prima operazione, e poiché passa dalla temp. $t_1$ alla temp. $t_2$, mentre la sua temp. si mantiene costante assorbe una quantità di calore $Q=C(t_2-t_1)$ dove C è il cal. spec. a volume costante. Il calore utilizzato, ossia trasformato in lavoro, durante l'intero ciclo, sarà, come in ogni ciclo, la differenza tra il calore Q attinto dalla sorgente e il calore q disperso nel refrigerante, perché compiuto il ciclo il calore interno riprende il valore iniziale, e quindi né parte del calore attinto potrà essere impiegato ad aumentare il calore interno, né parte del lavoro sviluppato (s'intende sempre definitivamente, durante l'intero ciclo) potrà essere prodotto a spese di calore interno. Ora per una formula precedentemente dimostrata, siccome il gas passa dallo stato C allo stato D mediante un'operazione isotermica, sarà

$$q=\frac{R}{J}(t_3+274)\log\frac{v_2}{v_3}$$

Ma si ha

$$\left(\frac{v_3}{v_1}\right)^{\frac{C}{c}-1}=\frac{t_1+274}{t_3+274}$$

$$\left(\frac{v_2}{v_1}\right)^{\frac{C}{c}-1}=\frac{t_2+274}{t_1+274}$$

Quindi

$$\left(\frac{v_2}{v3}\right)^{\frac{C}{c}-1}=\frac{t_2+274}{t_1+274}$$

$$\left(\frac{C}{c}-1\right)\log\frac{v_2}{v_3}=\log\frac{t_2+274}{t_1+274}$$

$$\log\frac{v_2}{v_3}=\left(\frac{c}{C-c}\right)\log\frac{t_2+274}{t_1+274}$$

e sostituendo in q

$$q=c\frac{1}{J}\frac{R}{C-c}(t_3+274)\log\frac{t_2+274}{t_1+274}$$

Ma si ha per una nota formula

$$\frac{R}{C-c}=J$$

quindi

$$q=c(t_3+274)\log\frac{t_2+274}{t_1+274}$$

donde si conclude

$$Q-q=c\left\{t_2-t_1-(t_3+274)\log\frac{t_2+274}{t_1+274}\right\}$$

$$\frac{Q-q}{Q}=\frac{t_2-t_1-(t_3+274)\log\frac{t_2+274}{t_1+274}}{t_2-t_1}=1-\frac{t_3+274}{t_2-t_1}\log\frac{t_2+274}{t_1+274}.$$

Ora facciamo $t_1$=0, $t_2$=150, $t_3$=50, come vuol il nostro esempio numerico, e troviamo

$\frac{Q-q}{Q}=1-\frac{324}{150}\log\frac{424}{274}$.

Eseguendo il calcolo si trova $\frac{324}{150}$=2,16. Si ha poi $\log\frac{424}{274}$ =log1,547...>log1,5>0,40.

Quindi $\frac{324}{150}\log\frac{424}{274}$ >2,16x=0,40=0,864 e finalmente $\frac{Q-q}{Q}$<1–0,864=0,136.

Concludiamo così che nel caso considerato si trasforma in lavoro, durante l'intero ciclo, meno di 0,136 del calore ceduto dalla sorgente, mentre, nel caso di un ciclo di Carnot se ne trasforma in lavoro $\frac{100}{424}$=0,235.

**164**
Gian Antonio Maggi a Luigi Sala

*[a lato:]* Sulla Teoria Dinamica dell'Aratro Mem dell'Ing Luigi Ferrario (Politecnico 1889) Lettera al Comm. Luigi Sala.

Ammettiamo, senza discuterla, l'ipotesi dell'ing. F che la falda di terra smossa dall'aratro si possa considerare come un corpo invariabile, che possiede un movimento di traslazione uniforme con velocità data, e un movimento di rotazione intorno ad una retta (AB) parallela alla direzione della traslazione.

Con ciò resta ancora e solamente indeterminata la legge della rotazione, e possiamo procurare di determinarla in modo che sia

$$\delta\int dt\, S\, mv^2 = 0$$

dove la somma S si estende a tutti i punti invariab. uniti da massa m componenti la falda, e l'integr. $\int$ va preso fra due valori qualunque del tempo, corrispondenti al passaggio della falda da una sua porzione ad un altra *[sic!]* qualsivoglia.
Perciò osserviamo che indicando con M la massa della falda, con K il suo momento d'inerzia per rispetto all'asse di rotaz. con C la velocità di traslaz. e con ω la velocità angolare, funzione incognita del tempo t, *[*si avrà*]*, nelle nostre ipotesi,

$$Smv^2 = Mc^2 + K\omega^2.$$

*[*Avremo*]*, per la funz. ω di t che soddisfa a (), conformemente ai principii del calcolo delle variaz.:

$$\omega=0.$$

E per conseguenza non sarà soddisfatta () se non supponendo nulla la rotazione, e la falda dotata del solo moto di traslazione.

L'Ing. F. fatte le ipotesi suddette, trova la linea che deve percorrere un punto grave libero per passare da un posto ad un altro in modo che sia $\delta\int mv^2 dt = 0$, e suppone che questa condiz. sia egualmente verificata dall'intersezione di un cilindro avente quella linea per direttrice e le generatrici orizzontali, con un cilindro circolare coll'asse orizzontale e perpendicolare alle generatrici della superf. precedente: perché (a parità di valori iniziali) $mv^2$ avrà uno stesso valore sulle due linee allo stesso livello. Ma questa eguaglianza non si verificherà per lo stesso valore di t (o di S). O quindi della *[*circostanza*]* che è $\delta\int mv^2 dt = 0$ per una linea ne segua che debba esserlo anche per l'altra.

Messina 1890.

**165**
Luigi Sala a Gian Antonio Maggi

Egregio amico,

Milano il 2 agosto 1896

Ricevo soltanto in questo momento il quesito sulla nota questione quale venne formolata da mio nipote Luigi Ferrario, e mi affretto a mandarglielo.

Nello stesso incontro, alfin di meglio interessarLa a leggere, con suo commodo, la mia memoria manoscritta del 15 pp. dicembre su alcune modificazioni teoriche e pratiche del ciclo di Carnot, Le unisco una mia commemorazione del Piatti (ottobre 1867)[197] ed un opuscolo 17 e 18 novembre 1894 intitolato Onoranze a Gio Batta Piatti. Mi farà cosa gradita a voler conservare per sé questi due stampati riguardanti persona citata nel mio manoscritto con lode speciale.

Colla più distinta considerazione mi è caro dirmi.

Devotissimo suo amico
Luigi Sala

Copia da consegnarsi all'Egregio Prof. Gian Antonio Maggi pel suo giudizio scritto.

Quesito 1^mo^

Abbiansi le espressioni:

(1) $\frac{Q}{T} = \frac{mT}{T} = \frac{q}{t} = \frac{mt}{t} = m$

(2) $\frac{V}{T_I} = \frac{gT_I}{T_I} = \frac{v}{t_I} = \frac{gt_I}{t_I} = g$

La (1) è desunta dalla ben nota formula $\frac{Q}{T} = \frac{q}{t}$ del ciclo di Carnot, e la (2) dalle formule del moto uniformemente accelerato; avvertendosi, ad ogni buon fine, che in questa seconda espressione i simboli V, v, $T_1$, $t_1$ significano rispettivamente due diverse velocità finali e costanti, un'accelerazione costante e due tempi diversi finali e costanti.

È fuor di dubbio che, nella seconda espressione, la derivata dell'accelerazione costante g è eguale a zero. Ora, in vista della grande analogia presentata dalle due espressioni, si domanda se anche nella prima espressione, la derivata della m, la quale m è evidentemente costante rispetto alle due temperature assolute T e t, sia o no eguale zero.

Quesito 2^do^

Nel moto uniformemente accelerato le velocità sono proporzionali ai tempi perché, in questo moto, tra velocità e tempo vi ha un rapporto costante.

Nel ciclo di Carnot, essendo la quantità di calore proporzionale alle temperature assolute, si domanda se anche fra le quantità di calore e le temperature assolute vi debba essere un rapporto costante.

---

[197] L. Sala, "L'ingegnere G.B. Piatti. Ferrovia ad aria compressa. Traforo del Cenisio. Commemorazione - Milano: Presso la Società per la pubblicazione degli Annali Universali delle Scienze e dell'Industria, 1867", *Annali Universali di Statistica*, Ottobre 1867.

Occorre appena di avvertire che qui intendesi parlare di costanza relativa e non assoluta.

Appendice 1$^{ma}$ ai quesiti 1$^{mo}$ e 2$^{do}$

Data la formola del ciclo di Carnot $\frac{Q}{T}=\frac{q}{t}=m$, la quale, decomposta ne' suoi fattori diventa $\frac{mT}{T}=\frac{mt}{t}=m$, Ferrario dice: - il rapporto m è una variabile continua perché possiamo variare simultaneamente le quantità di calore Q e q tenendo costante le temperature assolute T e t.

Sala risponde: - ammettesi che si possano variare le Q e q tenendo costanti le T e t, ma si possono variare non già in modo continuo e durante una stessa percorrenza di ciclo, ma soltanto col sostituire costanti arbitrarie ad altre costanti arbitrarie passando da una percorrenza all'altra del ciclo; cosicché il rapporto m invece di essere una variabile continua è una costante arbitraria o relativa.

Avviene qui quello che si verifica nella formola del moto uniformemente accelerato. $\frac{V}{T_I}=\frac{v}{t_I}=g$ da cui $\frac{gT_I}{T_I}=\frac{gt_I}{t_I}=g$.

Si possono certo variare le velocità V e v tenendo costanti i tempi $T_1$, $t_1$, ma variare non già in modo continuo, bensì colla sostituzione di costanti arbitrarie ad altre costanti arbitrarie passando da un moto uniformemente accelerato all'altro.

Nelle velocità V e v si può sostituire una serie di g, $g_I$, $g_{II}$, ... $g_n$ tutte costanti, come si può sostituire nelle qualità di calore Q, q una serie di m, $m_I$, $m_{II}$, ... $m_n$ tutte costanti: ma queste sostituzioni possono avvenire, ripetesi, soltanto nel passaggio dall'uno all'altro moto uniformemente accelerato, dall'una all'altra percorrenza di ciclo.

Appendice 2$^{da}$ ai quesiti 1$^{mo}$ e 2$^{do}$

Il Bordoni nella prima delle sue - Lezioni di calcolo sublime[198] a pag. 8 n°5 del 1$^{mo}$ vol. domanda: - Qual dev'essere il significato della φ perché la funzione φ(x+ ω) eguagli la somma delle due φ(x), φ (ω) qualunque siano x, ω: - od, in altri termini, quale dev'essere il significato della φ, perché l'equazione

(1) φ(x)+ φ(ω)= φ(x+ ω)

venga soddisfatta, qualunque siano i valori delle indeterminate x, ω.

Risponde il Bordoni col dimostrare luminosamente, mediante il teorema di Taylor, che il significato della φ consiste nel dover essere la stessa φ eguale ad una costante arbitraria K, che moltiplica tanto la somma delle due indeterminate x, ω quanto ciascuna di esse separatamente.

Abbiasi ora l'equazione del ciclo di Carnot

(2) $\frac{Q}{T}=\frac{q}{t}$, la quale, decomposta ne' suoi fattori, assume la forma,

(3) $\frac{mT}{T}=\frac{mt}{t}$, da cui

(4) $\frac{mT}{T}-\frac{mt}{t}=0$.

Dal primo membro dell'equazione (4) si può ricavare

[198] A. Bordoni, *Lezioni di calcolo sublime*, in 2 volumi, Giusti Fonditore, Milano, 1831.

(5) $m\left(\frac{T}{T}\right)+m\left(-\frac{t}{t}\right)=m\left(\frac{T}{T}+\left(-\frac{t}{t}\right)\right)$

Paragonando questa equazione colla (1), si vede che facendo $\frac{T}{T}=x$, $\left(-\frac{t}{t}\right)=\omega$, anche la m dev'essere eguale a φ.

Ora, essendo la φ eguale alla costante arbitraria K, dovrà anche la m corrispondere ad una costante arbitraria e relativa.

Dalla dimostrazione di Bordoni emerge pure che le derivate della costante arbitraria K <u>sono</u> <u>nulle</u>, <u>com'era</u> <u>d'altronde</u>, egli dice, <u>facile</u> <u>a</u> <u>prevedere</u>.

Appendice $3^{za}$ ai quesiti $1^{mo}$ e $2^{do}$

Dal ciclo di Carnot si deduce, pel calore esterno, la formola:

(1) $\frac{Q}{T}=\frac{q}{t}=m$ da cui

(2) $\frac{Q}{q}=\frac{T}{t}=K$

e pel calore interno l'analoga formula,

(3) $\frac{H}{T}=\frac{h}{t}=m'$, da cui

(4) $\frac{H}{h}=\frac{T}{t}=K$

Riunendo le espressioni (2) e (4) come fa il Saint Robert a pag 137 e 138 ed aggiungendovi la decomposizione delle quantità Q, q, H, h nei rispettivi fattori, facendo, cioè, Q=mT, q=mt, H=m'T, h= m't si ha l'espressione:

(5) $\frac{Q}{q}=\frac{mT}{mt}=\frac{H}{h}=\frac{m'T}{m't}=\frac{T}{t}=K$

da cui risulta che le quantità di calore esterno sono proporzionali alle quantità di calore interno.

Ora, essendo evidente che se nell'equazione $\frac{mT}{T}=\frac{mt}{t}$ si cangiasse la m di mt, a numeratore nel secondo membro, in $m_1$, quantità diversa dalla m, non può più essere soddisfatta l'equazione stessa, e quindi non può più sussistere la proporzionalità tra le quantità di calore e le temperature, ne consegue che anche le quantità $\frac{mT}{mt}=\frac{m'T}{m't}$, significanti, nell'espressione (5), la proporzionalità tra calore esterno e calore interno, cesseranno d'essere proporzionali quando, nel primo membro di questa equazione, la m del mt a denominatore avesse variato diventando $m_1$, quantità diversa da m.

Osservazione

Qualunque sia stata la precedente variabilità di m si vede che quando questa viene assunta come coefficiente di due temperature diverse T e t restando eguale a sé stessa, deve considerarsi come costante rispetto a T e t nel limite di questa variazione di temperatura.

Appendice $4^{ta}$ nei quesiti $1^{mo}$ e $2^{do}$

Data, pel calore esterno, la formola del ciclo di Carnot $\frac{Q}{T}=\frac{q}{t}=m$, la quale, decomposta ne' suoi fattori, diventa $\frac{mT}{T}=\frac{mt}{t}=m$, osservo che soltanto un'accurata analisi del periodo di formazione delle Q=mT e q=mt e del loro stato finale può far conoscere la causa dell'errore di Ferrario.

Durante il periodo di formazione della Q=mT tanto la Q come la m, di cui la stessa Q è funzione, sono variabili continue in funzione dello sforzo decrescente di espansione del gas, che assorbe il calore della fonte esterna; ma appena è cessato questo sforzo, cessano di essere variabili continue per diventare la Q una complessiva quantità finale di calore determinata e fissa, ed anche il fattore m, a cui ho dato il nome di massa virtuale, diventa esso pure una quantità finale, determinata e fissa.

Quanto a q=mt, la q rappresenta, quel calore che, col consumo di un lavoro esterno, si va formando per essere mano mano gettata sul refrigerante. Anche qui tanto la q come la m sono variabili continue in funzione dello sforzo crescente di compressione che riconduce il gas allo stato iniziale di pressione.

Cessato questo sforzo, la complessiva quantità di calore q = mt, stata gettata sul refrigerante, diventa una quantità di calore finale determinata e fissa, e per conseguenza diventa una quantità finale determinata e fissa anche il fattore m, da me chiamato massa virtuale ed eguale al coefficiente m di T.

Osservazione

Nel confronto altrove fatto tra la m del calore esterno, la m' del calore interno e la g accelerazione costante nel moto uniformemente accelerato, vi ha questa sola differenza, che m' e g sono sempre state costanti sino dall'origine mentre la m, fattore del calore esterno, fu variabile continua durante il periodo di formazione delle Q e q e diventa costante soltanto dopo che, cessato ogni sforzo decrescente o crescente, le quantità di calore Q e q sono diventate quantità di calore finali, determinate e fisse.

La causa dell'errore, in cui persiste in Ferrario, sta dunque nel confondere il periodo di formazione delle Q e q col loro stato finale. Il rapporto fornito dalla formola del ciclo di Carnot riguarda evidentemente, da una parte, temperature assolute e costanti e dall'altra, quantità di calore finali, determinate e fisse, non essendo esse più soggette ad alcun sforzo od altra causa qualsiasi di variazione.

Le variazioni, lo ripeteremo per l'ultima volta, non possono avvenire che passando da una percorrenza all'altra del ciclo; e avvengono non già in modo continuo, ma mediante sostituzione di costanti arbitrarie ad altre costanti arbitrarie come si verifica pel moto uniformemente accelerato ed anche pel calore interno se si cangia la massa del corpo intermediario.

**166**
Luigi Sala a Gian Antonio Maggi

Egregio Professore

Milano il 10 agosto 1896

Questa mattina parto per le acque del Masino e non sarò di ritorno che alla fine di agosto. Ciò mi ritarda il piacere di conferire con Lei sui noti argomenti. Ho però in compenso la soddisfazione di vedere offerta a Lei l'opportunità di rivolgere più lungamente l'attenzione sopra questioni e materie che non pajono prive d'importanza.

Ella può bene immaginare quanto sarei contento e quanto Le sarei grato se, al mio ritorno, un più lungo lasso di tempo mi avesse fatto conseguire in maggior copia le preziose sue osservazioni.

Rinnovando i sensi della più distinta considerazione mi confermo

Devotissimo suo Amico
Luigi Sala

**167**
Gian Antonio Maggi a Luigi Sala

Preg$^{mo}$ Signore

M'ero riservato di venire a *[*rivederLa*]* prima di partire per la campagna, quando ricevetti il Suo annuncio che alla sua volta lasciava Milano per tutto il mese. Mi spiacerebbe che avesse desiderato una più pronta risposta ai proposti quesiti,[199] dei quali, per precedenti impegni, non ho potuto occuparmi prima. Ed ora ecco ciò che pare a me di rispondervi.
Le due formule

$$\frac{Q}{T} = \frac{q}{t}, \quad (1)$$

$$\frac{V}{T} = \frac{v}{t} \quad (2),$$

la prima riferita ad un determinato ciclo di Carnot, e la seconda ad un determinato moto unif. acc., quantunque si riducano ambedue ad una semplice proporzione, risultano, sotto il punto di vista delle questioni in discorso, se ben si considera la natura degli elementi che vi figurano, d'indole sostanzialmente diversa. Dalla (2), i termini di ciascun rapporto sono particolari valori di due variabili - indichiamole con t' e v' - continue: mentre in (1), tenuto fermo che ci riferiamo ad *[sic!]* ciclo di Carnot determinato, la stessa circostanza non si verifica: ma, se indichiamo con t' e q' la temperatura d'una sorgente e la quantità (assoluta) di calore corrispondente, (ceduta o assorbita, a seconda del caso) [da essa ceduta] [(così possiamo chiamare anche la quantità di calore assorbita dal refrigerante, purché, a suo luogo, teniamo conto del segno]),[200] le variabili t' e q' sono suscettibili - [badiamo bene, pel ciclo cui si riferisce la (1)] - puramente dei valori T e t, Q e q.

[199] Si veda la lettera #165.
[200] Le parentesi sono così nel testo.

Perciò g si può considerare come il valor costante di una funzione d'una variabile continua, il tempo, rispetto alla quale si può derivare e la derivata (totale) risulta zero. Invece [finché il ciclo resta quello che si è inteso][qualunque del resto esso sia,] non si può egualmente considerare la derivata di m rispetto alla temperatura d'una sorgente: la quale non è suscettibile che di due valori, e non è quindi una variabile continua, com'è necessario supporre per discorrere di derivata.

Quindi, in risposta al I Quesito, intesa la (1) riferita ad un determinato ciclo di Carnot, non si può dire che la derivata della m rispetto alla temperatura d'una sorgente sia zero, perché non a priori si può discorrere, in tal caso della derivata medesima.

Per lo stesso ordine di considerazioni in risposta al II Ques., non ò opportunamente chiamata costante la m: come non sarebbe ben detto che il rapporto di 4 a 2 e di 8 a 4 è costante. Abbiamo in questo come in quel caso due valori eguali, e non ci può essere termine più proprio di eguale per esprimere la loro relazione. Il termine costante lascia supporre la possibilità a priori d'una molteplicità di valori del rapporto in discorso, quale appunto si avrebbe, quando i termini fossero variabili continue. Ciò che ha luogo pel rapporto dell'increm. di veloc. al corrispondente increm. del tempo, in un determinato movim. unif. acc. - questi increm. essendo variabili continue - e non ha luogo sul rapporto della quantità di calore ceduta da una sorgente alla temperatura di questa, in un determinato ciclo di Carnot, le quali non sono suscettibili che di due valori - quello che spetta alla riscaldante e quello che spetta alla refrigerante.

Venendo infine alla domanda contenuta nel foglio da Lei inviatomi da più tardi, trovo da osservare che quando si afferma che nella (1) Q può avere un valore qualunque corrispondentemente ad un determinato valore della temperatura T, mentre sulla (2) V ha un unico valore corrispondentemente ad un determinato valore del tempo T, le due relazioni (1) e (2) sono considerate da un punto di vista diverso. E cioè la (2) si riferisce ad un determinato moto unif. accelerato - per esempio, se si tratta della caduta libera di un punto grave, si riferisce ad un determinato luogo - mentre si applica la (1) a cicli di Carnot necessariamente diversi, poiché, fissato T, cambiar Q significa passare da un ciclo di Carnot ad un altro.

Togliamo la restrizione che si riferisce ad un moto unif. acc. determinato, e allora in (2) V potrà avere un valore qualunque corrispondentemente ad un determinato valore del tempo T. Non è il notissimo fatto che la veloc. acquisita da un punto grave liberamente cadente dopo un minuto secondo è diversa secondo la latitudine e la latitudine[201] del luogo? Ciò che collima con quanto è da Lei notato nell'App. I.

*[In un riquadro:* D'altra parte, se togliamo sulla (1) le restriz. che si riferisce ad un ciclo di Carnot determinato, noi vi possiamo considerare la temperatura d'una sorgente e la corrispondente quantità di calore come variabili continue, e nulla impedisce di derivare m per rispetto all'una e all'altra separatamente: che si suppone il corpo un gas perfetto, e stabilito puramente il rapporto $V_2$:$V_1$ dei volumi al principio e alla fine della comunicazione di calore Q alla temp. T, si avrà una successione continua di cicli di Carnot, nei quali varia in modo continuo Q e T giusta la relazione

$$\frac{Q}{T} = J \log \frac{v_2}{v_1} = \text{Cost.}$$

E in questo caso la derivata totale di m rispetto a T sarà zero, e m sarà opportunamente chiamata costante.*]*

[201] Forse Maggi intendeva scrivere latitudine e altitudine.

Realmente, se la (1) si riferisce ad un determinato ciclo di Carnot, come la (2) ad un determ. moto unif. acc. (alla caduta libera d'un punto grave in un determ. luogo) data T in (1) è determinato Q, non altrimenti che dato T in (2) è determinato V. Ma, per le ragioni sopra esposte, non per questo dico che le due formule siano paragonabili fra loro, così da poter valersi delle conseguenze dell'una per stabilire quelle dell'altra, o doversi meravigliare che conducano conclusioni diverse.

Con questo resto a' Suoi comandi, prefiggendomi in ogni caso di salutarLa prima che finiscano le vacanze. D'altro non Le parlo per non farLa più a lungo aspettare: e solo mi limito a dirLe che ho letto con molto interesse la biografia del Piatti, onde ne La ringrazio una volta di più
Gradisca ecc.

Mezzago 25 Agosto 96
(Lettera al Comm. L. Sala).

**168**
Gian Antonio Maggi a Luigi Sala

Al Comm. Luigi Sala - Milano.
Pisa 6 Marzo 97

Pregiatissimo Signore

Per mezzo della Mamma,[202] Le rinvio il manoscritto, che, in attesa di questa occasione raccomandata dalla maggior sicurezza, non vorrei aver trattenuto troppo a lungo.

Il problema ch'Ella si è proposto è di surrogare al ciclo, ch'io seguiterò a chiamare di Carnot, in omaggio alla consuetudine, un altro più vantaggioso. Ora, per quanto riguarda il coefficiente economico, se istanno i principii della Termodinamica, quello del ciclo di Carnot, a parità di temperature esterne, non può essere superato. Ella ricorda (pag. 7) questa proporzione, introducendovi la riserva che il ciclo comprenda operazioni isotermiche. Ma realmente essa sta senza alcuna restrizione; e infirmarla, sia pure col restringerne la generalità, significa sacrificare per lo meno il 2° principio ( o principio di Clausius).

In pratica potrà tornar più proficuo, secondo le esigenze del caso, ricorrere ad altri cicli: e se ne possono citare esempii diversi. La macchina a vapore, quelle ad aria calda di Ericsson e di Stirling. Ma questo per circostanze d'altra natura di quelle che considera la teoria per stabilire quel teorema.

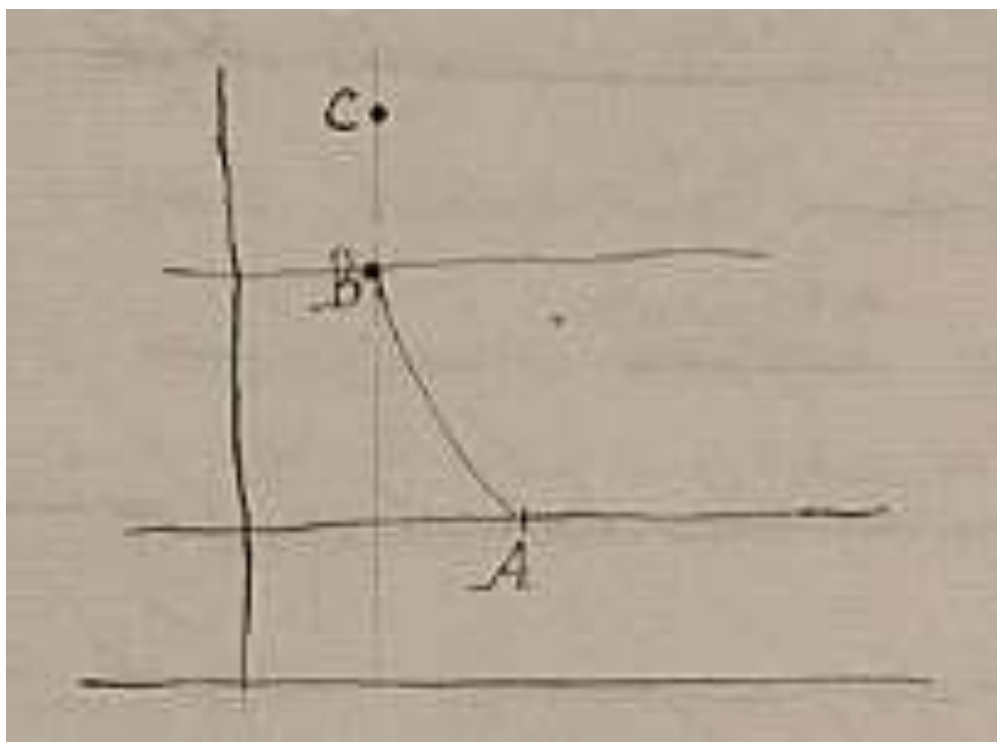

Ella propone, per conseguire un coefficiente economico maggiore, il procedimento da Lei chiamato ciclo M; ma un tal ciclo non vedo come sia possibile. Trasportato il corpo intermediano dallo stato A allo stato B con una trasformazione adiabatica, e da B in C con un riscaldamento a volume costante, non si può ricondurlo da C con una semplice trasformazione adiabatica in A: poiché in tal caso due distinti linee adiabatiche avrebbero comune il punto A, ciò che è impossibile, come emerge senz’altro dell’essere la loro equazione in coordinate p, v, e designando con S (p,q) l’entropia nello stato definito da p, v:

$$S(p,v)=\text{Cost.}$$

[202] La madre di Maggi si chiamava Clara Anelli.

Poi quanto all’impianto di Parigi, qui la questione sostanziale consiste nello spendere un lavoro in un posto per utilizzarlo in un posto lontano: che se si riuscisse ad utilizzare tutto il lavoro speso, ciò vorrebbe dire che il lavoro totale, al termine del procedimento sarebbe nullo Questione per tanto assai diversa, malgrado i punti di contatto, da quella del ciclo di Carnot e analoghi, dove si tratta, alla fine di un procedimento risultante di produrre un lavoro spendendo una corrispondente quantità di calore. Perciò quell’esempio non mi pare il più adatto per appoggiare la proposta del nuovo ciclo, il quale dovrebbe essere un ciclo termico, cioè rispondente allo scopo del ciclo di Carnot. Al quale proposito, occorre appena osservare che, nella questione della “trasmissione della pressione”, ciò ch’è chiamato coefficiente di rendimento per un significato affatto diverso dal coefficiente economico di ciclo termico.

Mi spiace che il poco tempo, che mi lasciano le mie occupazioni nel corso dell’anno scolastico, non mi abbia permesso di rendermi sufficiente conto dell’ottima parte, dedicata a considerazioni di meccanica molecolare. Ma non vedrei come ne potrebbero venir modificate le precedenti mie conclusioni: ben inteso, se ammettiamo che reggano gli ordinarii principii.

Sempre a’ Suoi comandi, mi è grato ecc.

## 169
Gian Antonio Maggi a **Nikolai Nikolájewitsch Saltykow**

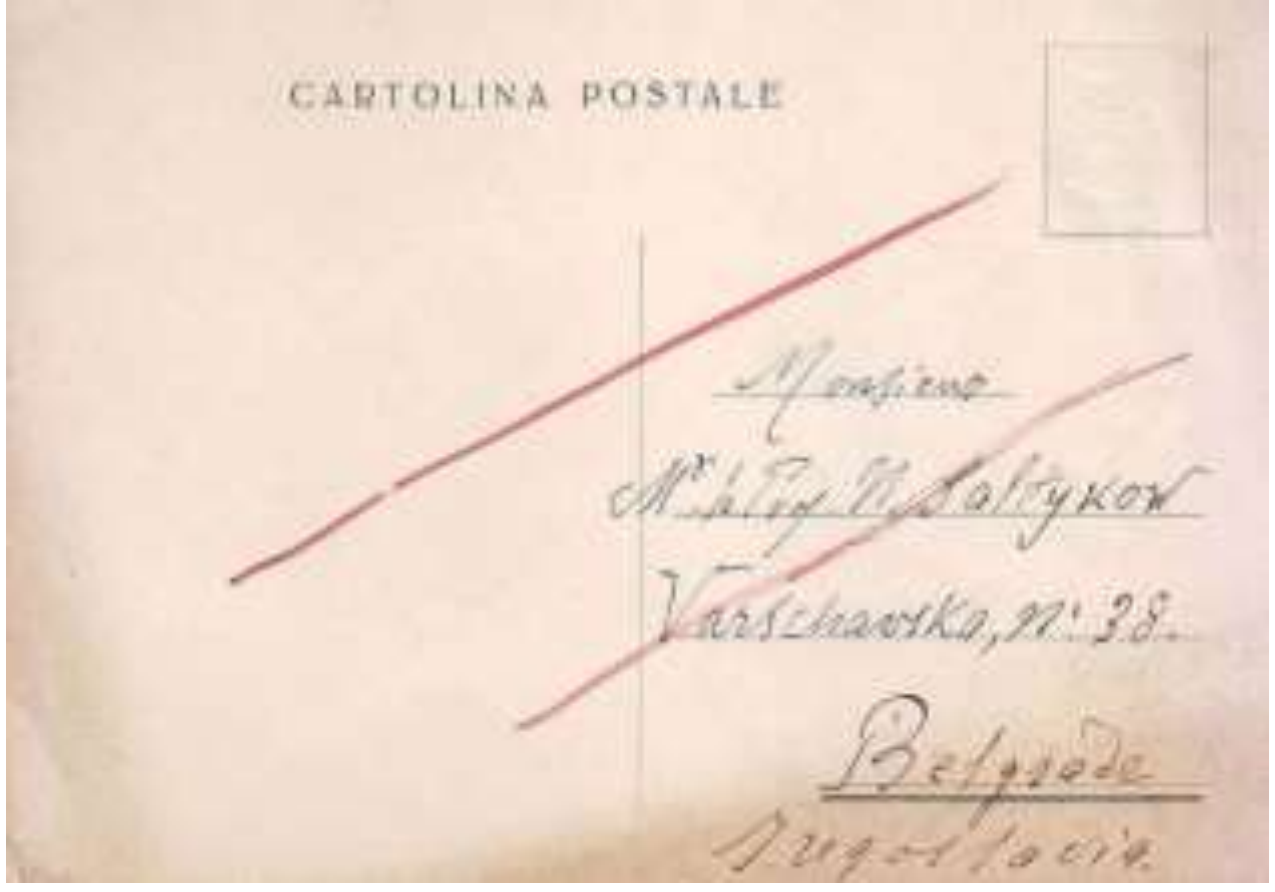
CARTOLINA POSTALE
Monsieur
Varschavska, N° 38.
Belgrade
Jugoslavia

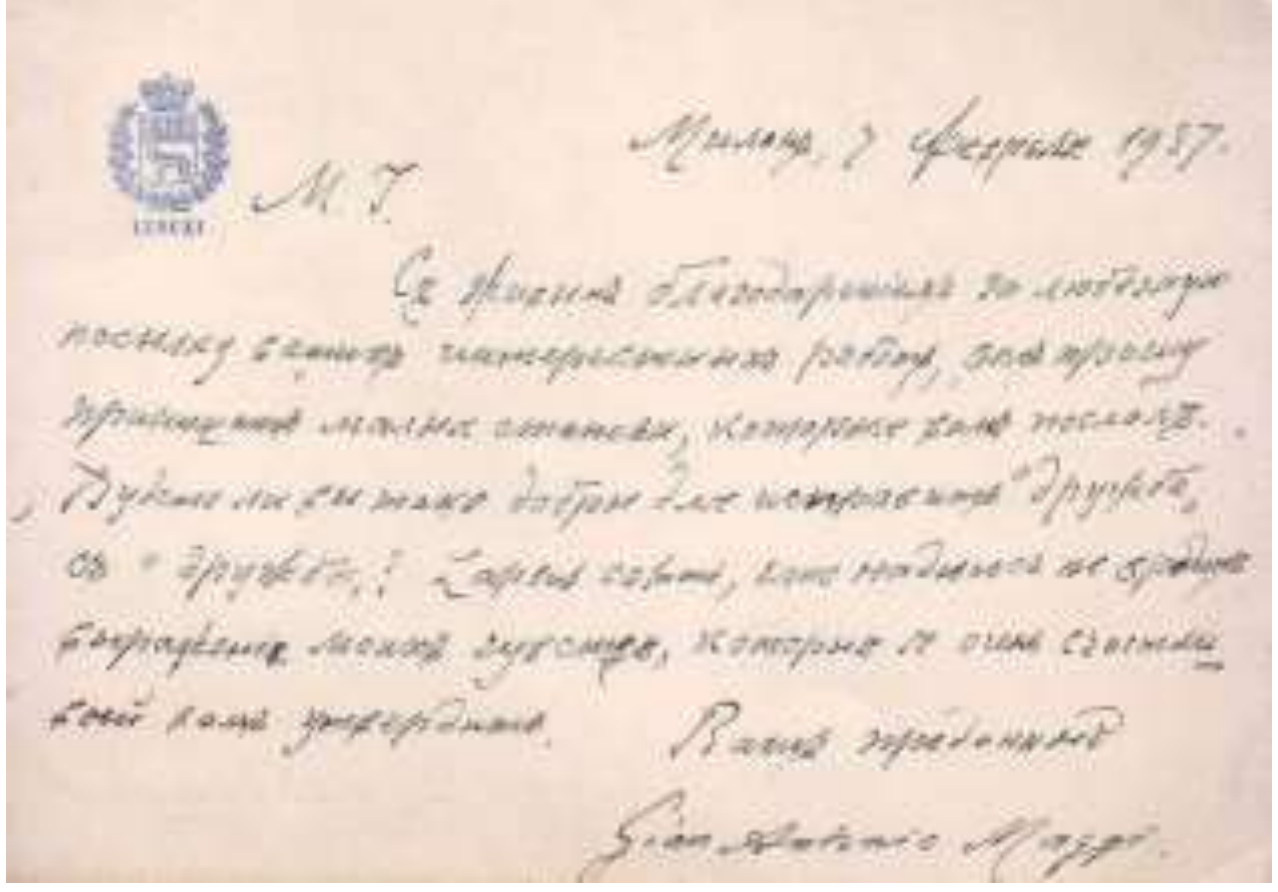
Gian Antonio Maggi

**170**
**Angelo Scala** a Gian Antonio Maggi

Napoli, 14 Marzo 1909

Egregio Sig. Professore.

Fin dalla fine dello scorso Gennaio mi ero proposto di scriverle a proposito del mio lavoro, ma disgraziatamente nei primi giorni di Febbraio son caduto ammalato prima di influenza, poi di catarro bronchiale e ho dovuto restare a letto più di mese. Le posso però annunziare che fin dalla fine di Gennaio ero giunto a qualche conclusione relativamente al punto più importante del mio lavoro, cioè alla possibilità o meno della deformazione tipica del 3° tipo, quale è definita dal Gebbia nella sua nota Memoria.[203] Questa possibilità dipende direttamente da quella della deformazione per interposizione e soppressione di materia dal Gebbia esposta ai n[i] 15-30 della Memoria stessa e contrastante, almeno apparentemente, col noto teorema del prof. Volterra sugli spostamenti polidromi, da Lei esposto in modo più intuitivo e interpretato in maniera più consentanea al significato fisico dello spostamento come teorema sugli spostamenti discontinui. Secondo questo teorema, in un corpo semplicemente connesso non sono possibili spostamenti discontinui, corrispondenti a una distribuzione regolare di parametri di dilatazione, mentre per la deformazione per interposizione e soppressione di materia quale è definita dal Gebbia anche per un corpo semplicemente connesso gli spostamenti hanno una superficie di discontinuità. A me sembra però che la contraddizione sia soltanto apparente. Infatti il teorema del prof. Volterra richiede esplicitamente per la sua dimostrazione la continuità dei parametri di dilatazione e delle loro derivate prime e seconde, mentre la teoria del Gebbia, esposta nel modo più rigoroso compatibile con il grado d'approssimazione della teoria matematica dell'equilibrio elastico, non richiede che la continuità delle derivate prime e seconde dello spostamento, anzi, come mi pare d'esser riuscito a dimostrare, per le derivate seconde non è neppur necessaria la continuità, ma solo l'esistenza e l'integrabilità nel campo occupato dal corpo nella configurazione anteriore. Queste ipotesi meno restrittive possono benissimo render possibile la deformazione per interposizione e soppressione di materia quale è definita dal Gebbia anche per un corpo semplicemente connesso, e, del resto, la costruzione tassativa di tale deformazione è, a parer mio, rigorosissima.

Avevo intenzione di mandarle fra qualche giorno la prima parte del mio lavoro contenente l'esposizione della teoria del Gebbia fino al punto controverso in discorso. In questa esposizione ho cercato di dare ai ragionamenti del Gebbia il necessario rigore, dove questo mancava, e di ridurre al minimo possibile le condizioni e ipotesi restrittive.

Ma avrei prima bisogno del Suo aiuto per chiarire un concetto che, mentre non infirma le conclusioni e non altera il rigore dei ragionamenti, mi procura tuttavia qualche imbarazzo nell'esposizione. E il concetto di cui intendo parlare è quello dei campi di forza introdotto dal Gebbia al n° 3 della sua Memoria. Egli avverte esplicitamente che non assegna le forze esterne né come funzioni delle coordinate dei punti della configurazione anteriore, né come funzioni delle coordinate dei punti della configurazione posteriore. Ma nelle equazioni d'equilibrio da lui trovate, egli fa comparire le forze esterne come funzioni delle coordinate dei punti della configurazione posteriore e cerca di ridursi alla configurazione anteriore, mediante uno sviluppo in serie di Taylor fatto nell'ipotesi che le funzioni esprimenti le componenti della forza esterna unitaria sian funzioni regolari. Ma funzioni regolari di che cosa? Delle coordinate dei punti della configurazione anteriore o posteriore no di certo, anche perché ciò sarebbe inutile, essendo il punto di partenza dello sviluppo un punto della configurazione anteriore e il punto d'arrivo (diciamo così per brevità) il punto corrispondente della configurazione posteriore. Tutto farebbe dunque supporre che il Gebbia col concetto dei campi di forza intenda considerar come variabili le forze esterne durante la deformazione, ma allora i valori delle funzioni che le rappresentano dovrebbero esser nulli nella configurazione anteriore, il che egli non suppone, né può supporre. Come Ella ben vede, è un concetto che per me non è affatto chiaro ed è perciò che io mi rivolgo alla Sua cortesia per qualche schiarimento.

[203] Probabilmente si tratta di: M. Gebbia, "Le deformazioni tipiche dei corpi solidi elastici", *Annali di Matematica pura ed applicata*, s. III, 1904, t. X, pp. 141-230.

Le chiedo scusa del disturbo che le arreco e La ringrazio anticipatamente della Sua gentilezza.

Accolga i miei più distinti saluti e mi creda.

Suo Dev$^{mo}$
Angelo Scala
Via Parma al Vasto n° 33 - Napoli.

P.S.- Le sarei gratissimo se Ella volesse ancora dirmi se mi conviene ammettere con Lamé le forze interne nulle quando il corpo è allo stato naturale ovvero limitarmi all'ipotesi di forze interne centrali (non nulle). E in questo ultimo caso come potrei definire con esattezza lo stato naturale?

**171**
Gian Antonio Maggi a Angelo Scala

Mandata sotto forma alquanto modificata dal terzo capoverso in poi.

Pisa 28 Marzo 1909

Caro Scala,

Mi ha fatto molto piacere di ricevere notizie Sue e del Suo lavoro, e di sentire che sia giunto a qualche risultato positivo; mi spiace bensì che ragione del Suo ritardo a fornirmene sia stata una abbastanza lunga indisposizione: procuri di mantenersi sano!

Ho sempre inteso che i risultati del Gebbia potessero giustificarsi di fronte al teorema di Volterra coll'accertare che non implicano le condizioni restrittive di questo teorema: e perciò crederei che Elle abbia oramai risoluto quella sostanziale difficoltà. Come Le accennai, il prof. Somigliana comunicò al Congresso dei Matematici, a Roma, un interessante esempio di spostamento elastico discontinuo di un corpo semplicemente connesso. Gli Atti del Congresso si stanno pubblicando, ed Ella potrà utilmente consultare quella comunicazione, se esce in tempo il secondo volume, che dovrebbe contenerla. Le consiglierei pure di leggere la Nota del compianto prof. Morera (Atti dell'Acc. delle Scienze di Torino, 1906-07) "Intorno all'equilibrio dei corpi elastici isotropi", che, specialmente nelle conclusioni, richiama la questione in discorso.

Per quanto ai "campi di forza" bisogna che Ella tenga ben presente che le forze di massa (forze limite) e le pressioni considerate sono quelle che corrispondono all'equilibrio del corpo nella configurazione posteriore, per cui il campo di forza è effettivamente sempre rappresentato dalla configurazione posteriore. Nella forma primitiva delle equazioni

$$\frac{\delta X_x}{\delta x} + \frac{\delta X_y}{\delta y} + \frac{\delta X_z}{\delta z} + X = 0$$

...

$$X_x \cos(nx) + X_y \cos(ny) + X_z \cos(nz) + X_n = 0$$

...

forze di massa (forze limite) e pressioni compariscono appunto come funzioni del punto (x,y,z) della configurazione finale. Se non che, essendo questa configurazione incognita, e supponendosi data invece la configurazione anteriore (lascio da parte le altre configurazioni, dal Gebbia indicate come fondamentali), si procura di rappresentare quelle forze e pressioni come funzioni del punto ($x_0$, $y_0$, $z_0$) della configurazione iniziale, o anteriore, corrispondente a (x,y,z) secondo le relazioni

$$x=x_0+\xi, \qquad y=y_0+\eta \qquad z=z_0+\zeta$$

dove $\xi$, $\eta$, $\zeta$ sono supposti infinitesimali. Gli sviluppi colla formola di Taylor che servono a questo scopo suppongono quindi la regolarità rispetto alle coordinate. Questi sviluppi rappresentano sempre le forze competenti all'equilibrio nella configurazione posteriore, e quindi elementi attinenti ai punti di questa configurazione - espressi però come funzioni dei punti corrispondenti, in base alla (2), della configurazione anteriore. Ella potrebbe vedere in proposito i §§512 e segg. della mia Meccanica. La Sua objezione, per verità, non mi riesce perfettamente chiara: perché esclude come inutile la considerazione della regolarità delle funzioni rispetto alle coordinate? Veda se Le giova tener ben presente la suddetta circostanza che il campo di forza, riferito sia alla configurazione anteriore che alla configurazione posteriore, appartiene sempre di fatto alla configurazione posteriore.

Con questo riesce chiarito anche il concetto dello "stato naturale", che va semplicemente inteso come lo stato del corpo in equilibrio nella configurazione anteriore con forze di massa (forze limite) a pressioni tutte nulle. Chiaramente non ha a che vedere con esso il campo di forze corrispondente all'equilibrio nella configurazione posteriore, riferito alla configurazione anteriore. Io ho procurato di meglio chiarire questi concetti coll'introdurre le forze di massa (forze limite) e pressioni corrispondenti all'equilibrio nella configurazione posteriore come aumenti di forze di massa e di pressioni corrispondenti all'equilibrio nella configurazione anteriore. (Meccanica, s.c.).

Le consiglierei infine di non fare la limitazione, per quanto non sostanziale, delle forze di massa nulle. Può fare invece, se Le serve, l'ipotesi delle forze centrali, come quella che, conformemente al concetto della delle forze limite, giova principalmente per le applicazioni.

Mi scriva sempre che creda, sicuro di farmi piacere (sempre accluso francobollo), e riceva intanto i miei cordiali sluti.

Aff^mo^ Suo
G.A. Maggi.

**172**
Angelo Scala a Gian Antonio Maggi
[cartolina postale con risposta pagata indirizzata: All'Egregio Sig. Prof.
Prof. Gian Antonio Maggi - R. Università di Pisa]

Napoli, 14 giugno 1909

Egregio Sig. Professore,

Ho presentato fin dalla fine del mese scorso la mia tesi di laurea e credo che Ella a quest'ora vi avrà già dato una occhiata. Avrei voluto far qualcosa di meglio, ma il tempo stringeva e il bisogno urgente di laurearmi in Luglio (già tardi, perché con un anno di ritardo) mi ha costretto a presentare senza altro i pochi risultati ottenuti fino a tutto lo scorso Aprile.

Confido che chi esamina il mio lavoro sia, a causa di queste considerazioni, se non contento, almeno indulgente. Ma, ove mai il lavoro stesso richiedesse tali e tante modificazioni da obbligarmi assolutamente a ritirarlo (cosa per me dolorosissima e dannosa), Le sarei gratissimo se Ella volesse avvertirmene non più tardi della fine del mese in corso (avendo io stabilito di partire da Napoli il 30 Giugno, per presentarmi all'esame di laurea che mi è stato scritto avrà luogo il 5 Luglio) per evitarmi un inutile e dispendioso viaggio.

La ringrazio vivamente del favore e Le chiedo scusa del disturbo che Le arreco. Accolga i miei più distinti e cordiali saluti e mi creda

Suo Dev$^{mo}$
Angelo Scala
Via Parma al Vasto n° 33 - Napoli.

**173**
Gian Antonio Maggi a Angelo Scala

Pisa 26 Giugno1909
Al Sig.$^{r}$ Angelo Scala - Via Parma al Vasto 33 - Napoli.

Caro Scala,

Perché non mandarmi in tempo utile un abbozzo della sua tesi? Io non potei vederla prima di jeri l'altro, di ritorno da un viaggio in Sicilia, e sono davvero dolente di doverLe dire, d'accordo coi professori Bianchi, Pizzetti, Bertini, che, per quanto desideriamo giovare ad un nostro studente che si è sempre fatto buona prova, e ci dolga del pregiudizio che può arrecargli il ritardo della laurea, non possiamo in alcun modo giudicarlo sufficiente. Ella non potrebbe altrimenti che convenire che essa è poco altro che una traduzione delle Memorie del prof. Gebbia, dalle quali, secondo il disegno ch'io Le avevo indicato, avrebbe dovuto ricavarne un diverso lavoro. Delle questioni ch'io Le proposi ha toccato la possibilità delle deformazioni del terzo tipo, e il comportamento delle pressioni elastiche all'infinito, Ma la prima questione non è chiarita, perché resta da verificare che nella fattispecie *[*compariscono*]* quelle singolarità per cui cade in difetto il teorema di Volterra. Le formole del Somigliana potrebbero nascere da una compensazione come per effetto di compensazione si annulla la deformazione forte del dato campo finito. E per quanto alla seconda, Ella, sviluppando il ragionamento del prof. Gebbia, ne riproduce la svista che dal comportamento di $\int_\Sigma L\, d\Sigma$ all'infinito si possa dedurre senz'altro quello di L. Supponga L=pcos (nx), con p costante che è il caso di una pressione costante normale alla superficie Σ, e troverà $\int_\Sigma L\, d\Sigma = p\int_\Sigma \cos(nx)\, d\Sigma = 0$, con qualunque valore di p.
Io non dubito che, ritornando sull'argomento, Ella potrà ricomporre una tesi meritando l'approvazione. Procuri di commentare la prima questione, e di correggere la seconda: io credo che basti assumere la proprietà asintotica in discorso come postulato, i punti all'infinito non comparendo fra i dati assegnati dal problema principale. In ogni caso, occorrerà che Ella almeno riprenda l'esposizione del prof Gebbia, e vi connetta quello studio comparativo ricordando principalmente le recenti ricerche del Fredholm, qualche considerazione critica, per presentare se non altro, anche come lavoro di compilazione, una relazione avente carattere d'argomento.

Si faccia animo, e creda che noi tutti, rammentando i suoi buoni studii, molto spiacenti di non poter dare per ora altro giudizio, ci teniamo, sopra tutto nel Suo interesse, ch'Ella li coroni onorevolmente.

Sarei lieto di fornirle altri consigli. Intanto creda al mio dispiacere di non poter altrimenti risponderLe e riceva i miei cordiali saluti

Aff$^{mo}$ Suo

**174**

Angelo Scala a Gian Antonio Maggi

Napoli, 30 Giugno 1909

Egregio Sig. Professore,

Aderendo al desiderio da Lei espresso nella Sua lettera, io ritiro la mia tesi di laurea per ripresentarla modificata e corretta nel prossimo Ottobre. Non Le nascondo però che sono profondamente scoraggiato, e non senza ragione. Aver lavorato quattro anni in modo non indegno di qualche considerazione, aver atteso un altro anno il conseguimento della laurea e vedersi ora in procinto di dovere attendere altri quattro mesi almeno, con immenso danno proprio e in maggior parte non per colpa propria, ma per avere avuta la disgrazia di incappare in un argomento poverissimo di risorse, è cosa che addolora ed avvilisce. E spero che in Ottobre non si ripeta il medesimo fatto; ma se ciò accadesse, sarò costretto a cambiare argomento, se pure ragioni di famiglia e la necessità di procacciarmi al più presto una stabile occupazione che mi metta in grado di provvedere al mio avvenire non mi obbligheranno a dire per sempre addio agli studî di Matematica e a rivolgermi a più proficuo lavoro (quale sarebbe p. es. la preparazione a concorsi per impieghi nelle pubbliche amministrazioni), considerando come perduti cinque anni di fatiche. Ciò mi sarà senza dubbio doloroso, ma che fare?

Mi perdoni questo piccolo sfogo, e veniamo al mio lavoro.
Quanto alle due questioni a cui Ella accenna nella Sua lettera, dicendo insufficientemente trattata la prima, e assolutamente sbagliata la seconda, farò quanto potrò per accontentarla. Ma quanto alla prima questione, io son d'avviso che l'indole troppo generale del lavoro non consenta d'andare più in là di quel che io ho fatto: almeno le confesso che non ne vedo il modo. Non vedo, in verità, come si possa verificare che nella fattispecie compariscono quelle singolarità per le quali cade in difetto il teorema di Volterra altro che nei singoli casi speciali; nel caso generale credo che non possa stabilirsi altro se non che l'inapplicabilità del teorema di Volterra può benissimo aver luogo nelle ipotesi su cui è fondata la teoria del prof. Gebbia. E preferirei quindi per ora non occuparmi più oltre di questa questione, sciupandovi intorno un tempo prezioso e col rischio di imbrogliare ogni cosa.

Quanto all'altra questione, convengo che Ella ha pienamente ragione e correggerò il punto sbagliato; ma finora non vedo il modo e Le sarei ben grato di nuovi e più ampî consigli in proposito. Come pure Le sarei riconoscentissimo se, trovando nel lavoro che ho presentato qualche altra inesattezza, volesse accennarmi il modo correggerla.

Non Le nascondo però che l'avere io già studiato l'argomento, e trovatolo così povero di risorse, mi toglie ogni illusione circa al partito che potrò ancora trarne. Vedrà che, al più, non verrà fuor che un lavoro di compilazione, ottenuto rifondendo l'esposizione del prof. Gebbia e mettendone il metodo e i risultati in relazione col metodo e coi risultati del Fredholm. Basterà questo? Io spero di sì. Se no, sarà per me una vera rovina, perché o sarò costretto a sostenere senz'altro l'esame di laurea in Ottobre, contentandomi di conseguir la laurea con una scarsa votazione, o, se neanche questa volesse accordarmi la Commissione Esaminatrice, mi vedrò obbligato a cambiar l'argomento, perdendo qualche altro anno intorno a un nuovo tema, o forse anche a dire addio agli studî di Matematica. Vede bene che per me è questione di vita o di morte.

Confido che non mi mancheranno i Suoi preziosi consigli ed Ella mi perdonerà se abuserò della Sua gentilezza e della Sua pazienza importunandola un po' troppo spesso. Le sarò gratissimo se, andando via da Pisa in queste vacanze, Ella vorrà comunicarmi il Suo indirizzo.

Perdoni il disturbo che Le arreco e accolga i miei più distinti saluti.

Suo Dev$^{\underline{mo}}$
Angelo Scala

**175**
Gian Antonio Maggi a Angelo Scala

Pisa 3 Luglio 1909

Caro Scala,

Prima di tutto, non si trattava tanto di aderire ad un mio desiderio quanto di approfittare di un avvertimento ch'io Le davo per Suo bene, a nome anche degli altri suoi professori, che Le nominai.

Piuttosto c'è da meravigliarsi che Ella abbia pensato per un momento che la Commissione potesse tener per buona una tesi, la quale si riduceva, ripeto, a poco altro che ad una trascrizione delle Memorie del Prof Gebbia, suggeriteLe per studio.

Non è neppure a posto la Sua esposizione che Sua gran disgrazia sia stata d'incappare in un argomento poverissimo di risorse. Io, sulla Sua domanda che Le dessi un tema attinente ai miei insegnamenti, Le avevo proposto una specie di saggio sulla formole risolutive delle equazioni dell'equilibrio elastico, in relazione colle recenti pubblicazioni del Gebbia, del Fredholm, del Volterra. Una lucida sposizione della questione, nei più ampi limiti che le assegnano gli studii del Volterra, da richiamarsi a suo luogo: un resoconto, con osservazioni critiche, dei risultati anteriori alle suddette pubblicazioni: un'applicazione del concetto a cui è informato il procedimento del Gebbia: semplificazioni ed emendamenti alla sua trattazione: coordinamento della medesima alle ricerche del Fredholm con qualche studio di codeste. Tutto ciò offriva elementi per una tesi soddisfacente. Certo una tesi più atta a provare l'acume critico e lo spirito metodico del candidato che non il suo talento inventivo. Ma se Ella voleva questo e non quello, doveva pensarci prima! Accettato il tema, come può chiamarsi vittima della sua sterilità, quando non prova che d'aver letto il solo Gebbia, e senza la meditazione sufficiente per riconoscere che il comportamento della questione della possibilità della deformazione del terzo tipo non può essere altrimenti che la continuazione dell'esame della relativa espressione analitica che il Gebbia, da Lei riprodotto, deve pure spiegare fino al punto che occorre per la sue conclusioni - per rilevare l'errore dalla deduzione sul comportamento all'infinito (circostanze, l'una e l'altra a cui io Le ho pure accennato) - per fare qualcosa di più che riprodurre, per la massima parte la trattazione del Gebbia, nella stessa sua forma? S'Ella abbandonerà la prima questione, rinuncerà *[*ad una risorsa*]* del tema. Per quanto alla seconda, Le ho già scritto ch'io credo che le proprietà asintotiche in discorso possono essere assunte come ipotesi, ciò che non restringe la generalità del problema fondamentale, perché le deformazioni finali del dato corpo finito si elidono l'una coll'altra.

Quello che dico di *[*occuparsi*]* a fare è il meno che poteva ricavare dal tema. Perché non presentare almeno codesto? Ella teme che sia poco. O come ha potuto non temere di presentare quel manoscritto che non contiene nulla? S'Ella ha fatto buona prova ne' Suoi studii, è una ragione perché reclamiamo che impieghi le Sue attitudini o terminarli anche bene. Del danno che può arrecarle il ritardo della Laurea ci spiace molto, ma Ella intenderà che simili considerazioni non devono avere che assai scarsa parte fra i criterii *[*di*]* qualunque Commissione Esaminatrice.

S'Ella vuol mutare tema, non se ne trattenga, badi bene!, per suoi riguardi a me. Diversamente, sono ben disposto a fornirle i consigli ch'Ella mi chiederà: il mio indirizzo, dopo la metà del mese, sarà Chiomonte (Provincia di Torino). Debbo però lamentare ch'Ella finora non abbia seguito, e neppure chiesto, maggiormente i miei consigli, e il primo che Le debba dare è di procurare di presentare un lavoro, sia pure modesto, ma serio.

Mi creda, con tanti saluti,

Aff.mo Suo

*[frammento che probabilmente fa parte di una precedente stesura di questa lettera]*

Infine non è ancora a posto ch'Ella addossi tutta all'argomento la colpa del quinto anno consumato, dopo che, quando venne a vedermi dopo le vacanze, nelle quali non s'era mai fatto vivo, mi disse di non aver potuto far altro che rendersi conto delle Memorie del Gebbia (in quattro mesi): e non mi scrisse più che il 14 Marzo per dirmi che una non breve malattia l'aveva obbligato a rimandare la sua lettera, per chiedermi alcune cosette alle quali ho subito risposto e annunciarmi l'invio della prima parte del lavoro con alcuni risultati - quali? poiché il poco che sia sulla possibilità etc. è precisamente ciò che le avevo detto io: né mi mandò questo altro, tanto che io non non *[sic!]* vidi la sua tesi, colla quale Ella resta al punto a cui diceva di esser giunto dopo le vacanze, che al mio ritorno di Sicilia, già passata per le mani dei miei Colleghi. Ella ha per lo meno lavorato con grande lentezza e omettendo di render conto di ciò che tentava di fare, ha rinunciato all'ajuto che altri gli poteva dare.

**176**

Gian Antonio Maggi a **Elena Schtschukin**

[carta intestata: R. Università di Milano - Istituto Matematico - Via C. Saldini, 50 (Città degli Studi) - Telef. 292-393 - Il Direttore]

Milano 23 Marzo 1932
Alla Signora Elena Schtschukin *[*Bodassen*]*
Leningrado *[...]*

Gentilissima Signora,

Di Lei, a Pisa, io posso attestare che vi fu iscritta come aspirante alla Laurea in Matematica, che frequentò i miei corsi di Meccanica e di Fisica Matematica. (Ella poi sa che non se ne presentò mai agli esami) e che ha partecipato alla Commissione dell'esame di Calcolo, Presidente il Dino, nel quale Ella ottenne la promozione. Questo sono pronto a dichiararlo al Sig.$^{r}$ ,[204] o a chi per lui, su conforme domanda. Perché io non saprei prendere l'iniziativa di scrivere a persona che non ho il bene di conoscere, e che non mi conosce, per mandare queste notizie di un ventennio fa che non mi persuado possa servire a qualche cosa. Né mi persuade che potrebbero servire le più esplicite dichiarazioni, a cui Ella accenna, per lo stesso lungo tempo trascorso, quando Ella si trova così a contatto di persone competenti, in grado di poter giudicare dalla sua attitudine presente, e dalla sua presente preparazione, che son quelle che veramente hanno valore. Ad ogni modo, tali dichiarazioni io, con dispiacere, non gliele posso rilasciare, per almeno tante ragioni quante ne ebbe, per non rilasciarle, il prof. Levi-Civita, che pure non mancò, dove poté, per quanto mi scrisse, d'interessarsi a Suo vantaggio.

Mi rallegra piuttosto di sentire che Ella abbia potuto acquistarsi la Laurea, e Le faccio i migliori augurii per l'avvenire.

Unisco con questo i miei saluti e mi confermo

Aff$^{mo}$ Suo
G.A.M.

[204] Lo spazio bianco è lasciato da Maggi.

**177**
Gian Antonio Maggi a **Corrado Segre**

Pisa 12 Gennajo 1919

Caro Segre,

Ti ringrazio vivamente del gentile invio del tuo Discorso Inaugurale,[205] una delle maggiori soluzioni che abbia trovato un matematico, per adempiere quel mandato: πρόβλημα Ἀρχιμήδειου, diceva, in simili casi, Cicerone. Cominciato a leggerlo, non lo lasciai più fin che non ne fui arrivato in fondo, e questo vuol dire l'interesse che mi ha suscitato, e il piacere che mi ha procurato. L'ho però anche letto colla dovuta attenzione, tanto che ho trovato di farvi una piccola osservazione, che mi permetto comunicarti. Il teorema di G. *[sic!]* Bernoulli afferma propriamente che, col crescere del numero delle prove, tende alla certezza la probabilità che il rapporto tra il numero delle volte in cui si presenta l'evento al numero totale della prova differisca dalla probabilità dell'evento meno di un numero prefissato, piccolo a piacere. (Il nostro compianto Pizzetti ne diede recentemente, nei Rendic. dei Lincei, una semplice quanto rigorosa dimostrazione). È pur vero poi che, alla stregua dei fatti, la differenza tra quel rapporto e la probabilità tende ad annullarsi, col crescere del numero delle prove. E se, per questo o per altro, ti sembra che debba intendersi il significato da me ricordato, ti prego di non accusarmi di eccessiva sottigliezza. Grazie di nuovo, e augurii di buona continuazione dell'anno, e una cordiale stretta di mano dal

Tuo aff. collega
G.A. Maggi.

**178**
Gian Antonio Maggi a **Antonio Signorini**[206]

Lanzo d'Intelvi (Como) 8 Sett. 1936

Carissimo Signorini,

Mi ha fatto un gran piacere ch'Ella sia stato chiamato a far parte della Commissione per la Meccanica Razionale, oltre il resto, che non occorre dire, perché trovi un così perfetto giudice la concorrente Sig.$^{na}$ Maria Pastori.[207] Mia concittadina, che trovai a Pisa, vincitrice del Concorso alla Scuola Normale, e compiuti egregiamente gli studii a Pisa, ritrovai a Milano, per diventare col tempo mia Assistente. Per le sue qualità personali e didattiche mi rimetto al documento che le rilasciai, come già mia Assistente, il quale traduce la mia più schietta opinione, a suo gran favore. Per la parte scientifica, vorrei che, all'infuori dell'interesse intrinseco, la Commissione facesse adeguato conto della coscienza, della serietà di proposito, dell'estesa e soda cultura che ne risulta dimostrata.

Prendesi la libertà di raccomandarla alla Sua autorevole attenzione a chi tanto si compiace ch'Ella gli ricordi i cari vincoli antichi, e accolga i più cordiali saluti, coi quali La prego di credermi sempre l'aff$^{mo}$ Suo

[205] C. Segre, "Le previsioni - Discorso Inaugurale" nov. 1918, *Annuario della R. Univ. di Torino*, [1918/19], pp. 11-25.
[206] Questa lettera si trovava nelle *Rusticationes*.
[207] Nella trascrizione dei *Giudizii* si veda l'attestazione a favore di M. Pastori scritta da Maggi (**XXIV**).

**179**
Gian Antonio Maggi a Antonio Signorini
[carta intestata: Lyncei]

Milano 25 Novembre 1936

Carissimo Signorini,

Le sono sinceramente grato dell'onore che volle fare alla mia raccomandazione, che, per quanto mi raccontano, fu il salvamento della mia raccomandata. A malgrado degli avversi giudizii, io sono persuaso che non hanno errato i giudici come Lei, che le hanno accordato la maturità: se, nei Concorsi a cattedra, va conferita adeguata importanza all'efficacia dell'insegnamento.

La ringrazio dei cari saluti recatemi da parte Sua dal collega Masotti, e cordialmente li ricambio, confermandomi

Aff$^{mo}$ Suo
Gian Antonio Maggi.

**180**
Gian Antonio Maggi a **Anna Maria Simonatti**[208]

Gentilissima Signorina,

Dal Suo manoscritto debbo desumere ch'Ella ha utilmente lavorato, ciò che, è ragione, oltre che per lodarLa, per contare che lavorerà anche quanto decorre ancora, per trarne una discreta tesi.

Intanto occorre meglio ordinare l'intera esposizione. Mettere in evidenza i diversi problemi che vi sono trattati, e, al tempo stesso, la loro mutua connessione. Ho indicato la divisione in Introduzione, Parte I, Parte II, Parte III.
L'Introduzione, in cui Ella si vale abbondantemente delle mie lezioni, potrebbe anche essere alquanto ridotta.

La Parte I fa riferimento, sulla traccia del foglio inserito. Mi sembra preferibile valersi della forma vettoriale, per analogia col seguente, in cui Ella fa uso di questa forma. Manterrà il metodo scalare, dove è di per sé stesso indicato.

Dalla Parte II occorrerebbe che, sul problema in discorso della riflessione e rifrazione delle onde sferiche ad una superficie piana, aggiungesse qualcosa a quanto io ne accenno nella mia Nota,[209] mentre dice anche di meno, e, tirando via, cade in qualche inesattezza. Questa Parte, così come sta, potrebbe essere levata; ben inteso, togliendo dall'Introduzione l'accennato proposito di trattare il suddetto problema, e indicandovi invece quello di trattare il problema che forma l'oggetto della Parte III. Come pure, il discorso concernente, in generale, la possibilità della propagazione di onde di forma qualsivoglia, in un mezzo isotropo, dovrebbe essere trasportato, o nell'Introduzione, o in principio alla Parte III, diventata Parte II. In ogni caso poi il suddetto problema della riflessione e rifrazione ecc., debitamente sviluppato (se non interamente risoluto), dovrebbe sempre formare la Parte finale allo stesso modo che è accennato alla mia Nota.

[208] Le formule presenti in questa corrispondenza non sono sempre leggibili con chiarezza; i vettori a volte sono indicati con la sottolineatura e volte no: si è mantenuta la forma originale.
[209] La tesi riguarda un approfondimento della Nota di Maggi: "Sulla propagazione delle onde luminose di forma qualsivoglia nei mezzi isotropi", *Rendic. dei Lincei*, v. XXIX, 1920, pp. 371-378; si veda anche la corrispondenza di Maggi con Daniele.

Nella Parte III converrà dedurre le equazioni a cui debbono soddisfare $\underline{E}_1$, $\underline{E}_2$, $\underline{H}_1$, $\underline{H}_2$, senza preoccuparsi, sul principio, delle condizioni div$\underline{E}$=0, divH=0.

Supposto però div$\underline{E}$, Ella troverà, come condizioni equivalenti

$$(1)\ \mathrm{div}\underline{E}_1+\frac{\alpha}{a}\underline{r}_I\times\underline{E}_2=0, \qquad (2)\ \mathrm{div}\underline{E}_2-\frac{\alpha}{a}\underline{r}_I\times E1=0,$$

(indicando con $\underline{r}_I$ il vettore avente l'orientaz. della normale alle onde e grandezza 1).
Dalle quali, supponendo

$$\underline{r}_I\times\underline{E}_2=0, \qquad \underline{r}_I\times\underline{E}_1=0,$$

cioè le onde elettriche trasversali, segue

$$(3)\ \mathrm{div}\underline{E}_1=0, \qquad (4)\ \mathrm{div}\underline{E}_2=0,$$

e reciprocamente.

Egualmente per le $\underline{H}_1$, $\underline{H}_2$, o le onde magnetiche.
Introducendo infine questa ipotesi nelle equazioni preesistenti, si arriverà alle equazioni da Lei trovate, che, non hanno la generalità dalle equazioni di Hertz, perché da div$E_1$=0 non seguono necessariamente le (3) (4), ma soltanto le (1), (2). Però l'ipotesi della trasversalità è quella che riporta alle note condizioni sperimentati, ed ha quindi una particolare importanza: ciò che Ella deve rilevare.

È abbastanza singolare che Ella si vale più volte dell'ipotesi $\underline{r}_I\times\underline{E}_1=0$, $\underline{r}_I\times\underline{E}_2=0$, senza indicarlo esplicitamente. Come troverebbe diversamente $E_2\wedge\underline{r}_I\wedge\underline{r}_I=-E_2$ (pag. XXXII)?

Per maggiori indicazioni Le serviranno il foglio inserito, e le postille, abbastanza abbondanti.

Per mandarmi la nuova redazione, può contare che mi tratterò a Gavinana le prime tre settimane del prossimo Settembre. Poi, dopo un passaggio per Pisa, ripartirò per Ancona. Con questo, Le presento i migliori saluti, e mi confermo

Aff[mo] Suo
G.A. Maggi.

Gavinana, Agosto 1924.

Alla Sig[na] Anna Maria Simonatti, Via Roma, 6, Livorno (Mandata la lettera il 3 Sett.).

**181**

Gian Antonio Maggi a Anna Maria Simonatti

Milano 28 Maggio 1926.

Gentilissima Signorina,

Proseguendo il calcolo della formola a riga 4 pag. 45 della Sua Tesi, introducendovi

$$\mathrm{div}E_2 = +\frac{\alpha}{a} E_1 \times r_I,$$

conformemente all'ipotesi divE=0 (pag. 46, formola (43)), valendomi di

$$E_2 \wedge r_I \wedge r_I = r_I(E_2 \times r_I) - E_2$$

(com'Ella può subito verificare), e del Suo risultato a pag. 49-53

$$\mathrm{rot}(E_1 \wedge r_I) + \mathrm{rot}E_1 \wedge r_I = E_1\left(\frac{1}{R+r} + \frac{1}{S+r}\right) + 2\frac{\delta E_1}{\delta r} - r_I \mathrm{div}E_1 - \mathrm{grad}(E_1 \times r_I).$$

dove si porrà $\mathrm{div}E_1 = -\frac{\alpha}{a} E_2 \times r_I$ (formola (43))

io trovo in fine

$$(1)\quad \frac{\alpha^2}{a} 2r_I(E_2 \times r_I) = a\Delta_2 E_2 - 2\alpha\frac{dE_2}{dr} - \alpha E_1\left(\frac{1}{R+r} + \frac{1}{S+r}\right)$$

e allo stesso modo, in base alle (41), (42)

$$(2)\quad \begin{cases} \frac{\alpha^2}{a} 2r_I(E_1 \times r_I) = a\Delta_2 E_1 + 2\alpha\frac{\delta E_2}{\delta r} + E_2\left(\frac{1}{R+r} + \frac{1}{S+r}\right) \\ \frac{\alpha^2}{a} 2r_I(H_2 \times r_I) = \text{etc.} \\ \frac{\alpha^2}{a} 2r_I(H_1 \times r_I) = \text{etc.} \end{cases}$$

Ora i secondi membri devono annullarsi in conseguenza delle equazioni della mia Nota dei Lincei (22 Dicembre 1920)[210] a cui debbono soddisfare $E_1$, $E_2$, $H_1$, $H_2$ in conseguenza delle

$$\frac{\delta^2 E}{\delta t^2} = a^2\Delta_2 E \quad \text{colla condizione divE=0}$$

$$\frac{\delta^2 H}{\delta t^2} = a^2\Delta_2 H \quad \text{“} \qquad \text{“} \qquad \text{dir H=0}$$

(non dimenticare questa condizione). Se ne conclude

$$(3)\quad \begin{cases} E_1 \times r_I = 0 \qquad E_2 \times r_I = 0 \;\text{ cioè }\; E \times r_I = 0 \\ H_1 \times r_I = 0 \qquad H_2 \times r_I = 0 \;\;\text{"}\;\; H \times r_I = 0. \end{cases}$$

Vale a dire che E e H, stabilito che debbano essere incompressionali, risultano in conseguenza delle equazioni di Hertz trasversali, qualunque sia la forma delle onde armoniche definite dalle equazioni alle ultime due righe di pag. 41.

Questo è un risultato certamente importante, e che io credo nuovo. Per cui, riserbandomi nel frattempo di fare qualche ricerca in proposito, io Le proporrei di verificare prima di tutto l'indicata analisi, e poi redigere una breve esposizione, limitandosi ai passaggi principali (evitando la comparsa delle componenti), ch'io potrei comunicare, come Suo lavoro.

[210] Si veda la nota precedente.

La traccia dovrebbe essere la seguente
1) Ricordare brevemente la mia Nota, il relativo problema, e il risultato delle due equazioni differenziali.
2) Introdurre le E, H date dalle formule armoniche di pag. 41 nelle equazioni di Hertz, e accennando ai passaggi principali arrivando alle (1) e (2) del presente.
3) Dedurre coll'indicata considerazione le (3), e la seguente conclusione.
Le rimando con questo il Suo M.S., e mi farà piacere ad accennarmene sollecitamente la ricevuta e ad assicurarmi che si è messa all'opera.
Intanto La saluto cordialmente e mi confermo

Aff$^{mo}$ Suo

**182**
Gian Antonio Maggi a Anna Maria Simonatti
[carta intestata: Lyncei]

Roma 4 Giugno 1926[211]

(Non mandata)

Gentilissima Signorina,

Prima di tutto Le rinnovo rallegramenti e auguri cordiali pel Suo fidanzamento. Lieto e importante avvenimento che intendo bene come, in questo momento, debba, in prevalenza, attirare i Suoi pensieri. Mi permetto tuttavia di insistere amichevolmente perché voglia anche pensare alla redazione del breve scritto, che Le proposi, e mandarmelo, possibilmente, in tempo ch'io possa farne la comunicazione ai Lincei per l'adunanza del 20 corrente.

L'introduzione, che Le accennavo, potrebbe ridursi a coteste parole:

"Il prof. G.A. Maggi nella Sua Nota "Sulla propagazione delle onde luminose di forma qualsivoglia in un mezzo isotropo" (*) *[rimando a piè di pagina:* (*) Questi Rendiconti, Vol. Serie...*]*, propostosi di trovare se e sotto quali condizioni l'equazione di d'Alembert

$$\frac{\delta^2\varphi}{\delta t^2} = a^2\Delta_2\varphi$$

è soddisfatta da un'espressione della forma

(1) $\varphi = f_1(p)\cos\alpha\left(\frac{r}{a} - t\right) + f_2(p)\sin\alpha\left(\frac{r}{a} - t\right)$

dove α e a sono costanti, e, immaginata una famiglia di superficie parallele di forma qualsivoglia, r indica la misura della normale compresa tra una superficie base e la superficie passante pel punto P, ha trovato per condizione necessaria e sufficiente che le funzioni $f_1(p)$, $f_2(p)$ soddisfacciano il sistema di equazioni differenziali

[211] È in due copie di minuta, entrambe con la dicitura *Non mandata*.

$$(2)\begin{cases} a\Delta_2 f_1 - 2\alpha \dfrac{\delta f_2}{\delta r} - \alpha\left(\dfrac{1}{R+r}+\dfrac{1}{S+r}\right)f_2 = 0 \\ a\Delta_2 f_2 + 2\alpha \dfrac{\delta f_1}{\delta r} + \alpha\left(\dfrac{1}{R+r}+\dfrac{1}{S+r}\right)f_1 = 0, \end{cases}$$

dove R e S indicano i raggi di curvatura principali alla base della normale r relativa al punto P.

Io mi propongo una simile questione per le equazioni di Hertz del campo fisso, occupato da coibente isotropo neutro

$$(3)\ \varepsilon\frac{\delta E}{\delta t} = \mathrm{crot}H, \qquad \mu\frac{\delta H}{\delta t} = -\mathrm{crot}E, \qquad \mathrm{div}E = 0, \qquad \mathrm{div}H = 0,$$

introducendo in queste equazioni per E e per H le espressioni conformi a (1)

$$(4)\ \begin{aligned} E &= E_1(p)\cos\alpha\left(\frac{r}{a}-t\right) + E_2(p)\sin\alpha\left(\frac{r}{a}-t\right) \\ H &= H_1(p)\cos\alpha\left(\frac{r}{a}-t\right) + H_2(p)\sin\alpha\left(\frac{r}{a}-t\right). \end{aligned}$$

Si trovano con queste un sistema di equazioni differenziali del 1° ordine per le quattro funzioni $E_1(p)$, $E_2(p)$, $H_1(p)$, $H_2(p)$ da cui eliminando successivamente la seconda e la prima coppia, si deducono due coppie di soluzioni, le quali, in base alle (2) che possono essere soddisfatte dalle stesse funzioni, permettono di enunciare il risultato notevole che i quattro vettori devono essere tangenti alla superficie in P e, per conseguenza, la propagazione di E e H per onde trasversali.

Dopo di questo Ella potrà passare al calcolo, che, come Le dicevo, basterà ridurre ai passaggi principali, senza entrare nei particolari. E per ciò non avrà in questa parte, che da trascrivere abbreviata l'analisi della Sua tesi. Levando la condizione $E_1 \times r_I t=0$ etc., come Le indicavo, arriverà a $a\Delta_2 E_1 + 2\alpha\frac{\delta E_2}{\delta r} + \alpha\left(\frac{1}{R+r}+\frac{1}{S+r}\right)E_2 = \alpha 2\alpha r_I (E_1 \times r_I)$ e potrà concludere con codesti termini:

"Ora, poiché dalle (3) conseguono

$$\frac{\delta^2 E}{dt^2} = a^2\Delta_2 E, \qquad \frac{\delta^2 H}{\delta t^2} = a^2\Delta_2 H,$$

le $E_1$, $E_2$, $H_1$, $H_2$ debbono soddisfare le (3). Quindi

$$E_1\times r_I=0, \quad E_2\times r_I=0, \quad H_1\times r_I=0, \quad H_2\times r_I=0$$

ossia, per (4)

$$E\times r_I=0, \quad H\times r_I=0.$$

*[seconda copia di minuta]*

(Titolo)

<u>Proprietà delle onde hertziane armoniche di forma qualsivoglia, in un mezzo isotropo</u>

(Introduzione)

(Non mandata) 1

Il prof. G.A. Maggi nella Sua Nota "Sulla propagazione delle onde luminose di forma qualsivoglia in un mezzo isotropo" (*), propostosi di trovare, se, e sotto quali condizioni, l'equazione di d'Alembert

(1) $$\frac{\delta^2\varphi}{\delta t^2} = a^2\Delta_2$$

è soddisfatta da un'espressione della forma

(2) $\varphi=f_1(p)\cos\alpha\left(\frac{r}{a}-t\right)+f_2(p)\sin\alpha\left(\frac{r}{a}-t\right)$,

dove α e a sono costanti, e, immaginata una famiglia di superficie parallele di forma qualsivoglia, r indica la misura della normale compresa tra una superficie base e la superficie passante pel punto P, ha trovato, per condizioni necessarie e sufficiente che le funzioni $f_1(p)$, $f_2(p)$ soddisfacciano il sistema di equazioni differenziali

(3) .....

dove R e S indicano i raggi di curvatura principali alla base della normale r relativa al punto P.

Io mi propongo una simile questione per le equazioni di Hertz del campo fisso, occupato da coibente isotropo, neutro,

(4) $\varepsilon\frac{\delta E}{\delta t}=c\,\mathrm{rot}H,\qquad \mu\frac{\delta H}{\delta t}=-c\,\mathrm{rot}E,\qquad \mathrm{div}E=0,\qquad \mathrm{div}H=0,$

introducendo in queste equazioni per E e per H le espressioni conformi a (1), e cioè rappresentanti una propagazione degli stessi vettori per onde parallele della supposta forma, armoniche, con velocità di propagazione a,

(5) $$\begin{aligned}E&=E_1(p)\cos\alpha\left(\frac{r}{a}-t\right)+E_2(p)\sin\alpha\left(\frac{r}{a}-t\right)\\ H&=H_1(p)\cos\alpha\left(\frac{r}{a}-t\right)+H_2(p)\sin\alpha\left(\frac{r}{a}-t\right).\end{aligned}$$

Si trova con questo un sistema di equazioni differenziali del 1° ordine per le quattro funzioni $E_1(p)$, $E_2(p)$, $H_1(p)$, $H_2(p)$ abbastanza complicato, da cui però eliminando successivamente la seconda e la prima coppia, si deducono due coppie di equazioni di forma notevole, le quali, in base alle (3), che debbono essere soddisfatte dalle stesse funzioni, permettono di enunciare il pur notevole risultato che le accennate onde, relative alla propagazione di E e di H, sono trasversali.

...

(Non mandata) 2

(Conclusione)

Ora poiché dalla (4) conseguono notoriamente

$$\frac{\delta^2 E}{dt^2}=a^2\Delta_2 E,\qquad \frac{\delta^2 H}{\delta t^2}=a^2\Delta_2 H,$$

le $E_1$, $E_2$, $H_1$, $H_2$ debbono soddisfare le (3). Quindi

$$E_1\times r_I=0,\quad E_2\times r_I=0,\quad H_1\times r_I=0,\quad H_2\times r_I=0$$

ossia, per (5)

$$E\times r_I=0,\quad H\times r_I=0.$$

**183**

Gian Antonio Maggi a Anna Maria Simonatti

Milano 16 Giugno 1926

Gent. Sig. Dr. Anna Maria Simonatti
Via Roma 6 Livorno.

Gentilissima Signorina,

Sono davvero spiacente di averLa inutilmente incomodata, tanto più in tempo in cui era inopportuno più che mai. Perché ripreso, dopo la Sua lettera, che ricevetti jeri mattino, il calcolo, come Le scrissi che mi ero proposto di fare, cogli elementi che mi mancavano a Roma, ho constatato che, coll'accennata correzione, le note formule si riducono alle equazioni della mia Nota, e ne sparisce il risultato ch'io destinavo alla proposta comunicazione. Che se ben si guarda, il processo con cui queste equazioni si deducono dalle equazioni di Hertz è, in sostanza, quello con cui se ne deducono per E e per H l'equazione di d'Alambert. Mi fornì le equazioni della mia prima lettera un singolare errore di trascrizione, non di calcolo, che ricevetti a Roma, ritrovando sulle formole, per scriverLe, come poi feci da Pisa, e mi era parso allora che poteva stare la correzione che recai sulla mia lettera a Livorno, colla quale non si mutava che un fattore innocuo. Ma verificato il segno, da cui dipendeva il risultato, in luogo di $2\underline{r}_I$ risulta 0, e scompare il prodotto $E_2 \times \underline{r}_I$.

Ella si ricorderà che, nella mia prima lettera, Le proponevo, innanzi tutto, di verificare i calcoli. Nulla di male che le Sue presenti occupazioni Le abbiano risparmiato un a lavoro negativo. Ma, per trarre pure qualche frutto da tutto questo, Ella non mancherà di rilevare come i calcoli espongono qualche volta a insospettati errori, e lo scrupolo che occorre imporsi prima di comunicare al pubblico un proprio risultato.

Di nuovo Le esprimo il mio dispiacere, e coi migliori saluti, dei quali La prego far parte ai Suoi genitori, mi confermo

Aff Suo

P.S. Gradirei un cenno di ricevuta di questa lettera, colla quale La prego pure di ricordarsi di ritirare dalla Tesi la mia lettera, che casomai sarà meglio distruggere che conservare.

**184**[212]
Gian Antonio Maggi a **Carlo Somigliana**

Pisa 29 Giugno 1916

Carissimo Somigliana,

Se fosti a Milano, avrai trovato, alla tua posta, il mio biglietto. Mi spiacque assai di non trovarti, e conto di essere più fortunato in qualche ventura gita a Milano o a Torino. Noi passeremo il Luglio e l'Agosto a Viareggio. Non dimenticarci, se passi da queste parti, che tutti ti vedranno con gran piacere.
Ho riletto, in questi giorni, riprodotte dal "Nuovo Cimento", le tue veramente belle ed esaurienti ricerche sulle distorsioni elastiche: e conto farne soggetto di particolar considerazione, pel corso di Fisica Matematica dell'anno venturo. Per tuo merito, quella teoria diventa un sodo e completo capitolo della Statica Elastica.

Ho poi bisogno che mi levi un dubbio. Le tue formole del Cap. II della 2° Nota mettono completamente a posto le distorsioni di Weingarten. Le relative condizioni ne risultano, senz'altro, verificate, nell'ipotesi che il taglio sia rappresentato da uno spostamento rigido infinitesimale. Ma, per le distorsioni di Volterra, occorre di più che siano continue, attraverso la superficie di discontinuità, le derivate prima e seconde dei coefficienti di deformazione. Che cosa si può affermare a questo proposito? Intanto, una volta che tu hai dimostrato che, colla discontinuità dello spostamento, è determinata la discontinuità delle derivate prima e seconda, mi sembra che non si possano più fare, senza riserva, ipotesi arbitrarie sullo spostamento e sui coefficienti di deformazione e loro derivate.

Certo, quella proprietà si deduce dalle

$$Du=l+qz-ry, \qquad Dv=\ldots \qquad Dw=\ldots$$

(riferite ad un terna d'assi qualsivoglia), identificando

$$D\frac{\delta u}{\delta x}, \qquad D\frac{\delta^2 u}{\delta x^2} \text{ etc.} \qquad \text{con } \frac{\delta Du}{\delta x}, \qquad \frac{\delta^2 Du}{\delta x^2} \text{etc.}$$

Ma questo scambio si può fare nell'ipotesi che la superficie di discontinuità si possa spostare, senza modificare la deformazione, come un diaframma: la qual ipotesi, a parte le coordinate tangenziali, risulta lecita, ammesse le condizioni restrittive di Volterra. È ciò che io intendo di fare nella mia Nota dei Lincei del 1908[213] (troppo sommaria, e dove uso costanti va letto funzioni, come poi subito si vede). È un dubbio che presumibilmente tu mi puoi levare con due parole: e perciò scuserai la mia pigrizia a cercare da me lo schiarimento.

Non mi era sfuggita la Nota della Sig.na Armanni,[214] che, per trovare spostamenti elastici derivati da quelli di Volterra ha avuto bisogno di abbandonare l'infinitesimale. La Sig.na Armanni, fedele allo spostamento polidromo, probabilmente non sospetta, nell'enunciato del teorema fondamentale da cui prende le mosse, l'esclusione dello spostamento discontinuo, rappresentabile (non dico interpretabile) o no, come spostamento polidromo. Per cui, forse, a metterla in guardia, poteva bastare la mia Noticina del 1905 nei

[212] La corrispondenza con Somigliana fino al 1917 si trova racchiusa in un fascicolo denominato: *Sugli spostamenti elastici discontinui e le distorsioni di Volterra. Corrispondenza con Somigliana e Levi-Civita e appunti*, ma non vi ho trovato alcuna lettera di Levi-Civita.
[213] G.A. Maggi, "Sugli spostamenti elastici discontinui", *Rendiconti della R. Accademia dei Lincei*, 1908, pp. 257-264.
[214] G. Armanni, "Sulle deformazioni finite del solidi elastici isotropi", *Nuovo Cimento*, 1915, pp. 427-447.

Rendiconti dei Lincei,[215] che vi dorme il sonno dei giusti; ma, per questo, l'ulteriore sviluppo della teoria ha mancato di renderle piena ragione.

Colla speranza di non tardar molto a vederti, per discorrere di queste e d'altre cose a voce, ti mando intanto, da parte anche dei miei, i più cordiali saluti. Un'affettuosa stretta di mano ricevi in particolare, dal

Tuo Aff$^{mo}$
G.A. Maggi.

**185**
Carlo Somigliana a Gian Antonio Maggi
[busta indirizzata a: Chiar.$^{mo}$ Signor Prof. G. Antonio Maggi
dell'Università - Pisa]

Casanova Lanza (Como) 6.7.916

Carissimo Maggi,

Rincrescermi assai di non averti potuto vedere a Milano e ti ringrazio della visita. Io andrò ancora a Torino nella settimana prossima, ma non so se potrò vederti. Risponderei a voce alla domanda che mi fai circa la teoria della distorsione e probabilmente meglio di quello che posso fare per lettera.

Se ho ben compreso (ma non ne sono sicuro) la difficoltà su cui richiami la mia attenzione riguarda la legittimità degli scambi

$$D\left[\frac{\delta^2 u}{\delta x^2}\right], \frac{\delta^2}{\delta x^2} D[u], \text{ecc.}$$

Ora tali scambi <u>non</u> sono <u>in generale</u> leciti; quantunque lo siano in moltissimi casi. Io ho ammessa la legittimità loro da principio, poi l'ho giustificata in fine con altre considerazioni. Perciò la mia trattazione è ben lontana dalla perfezione; se ora rifacessi il lavoro, seguirei altra via. Tuttavia le conclusioni vanno bene, e il fine giustifica i mezzi… anche in matematica! Spero di vederti presto e ne riparleremo.
Ho terminato negli scorsi giorni a Torino lo studio delle discontinuità di quegli integrali che io chiamo potenziali elastici e ho dato la dimostrazione di quelle formole di cui mi sono servito nelle due Note sulle distorsioni. La stampa si farà presto e te ne manderò un estratto.[216]

Con ciò il mio programma di teoria dell'elasticità è esaurito; e non me ne occuperò più, almeno per un gran pezzo. Ne sono un po' stufo e poi a furia di ripensare alle stesse quistioni, mi pare che ci si avvicini all'imbecillimento.

Per passare ad argomento ben diverso nel prossimo mese di Agosto andrò a Macugnaga, per compiere, insieme ad un ingegnere del Politecnico,[217] il rilievo topografico di quel ghiacciajo. Se vuoi fare una bella e comoda gita di montagna, vieni a trovarmi. Il panorama del Monte Rosa da questo versante orientale è veramente grandioso.

---

[215] G.A. Maggi, "Sull'interpretazione del nuovo teorema di Volterra sulla teoria dell'elasticità", *Rendiconti della R. Accademia dei Lincei*, 1905, pp. 255-256.
[216] C. Somigliana, "Sulle discontinuità dei potenziali elastici," *Atti della R. Accademia delle Scienze di Torino*, 51, 1915-16, pp. 1330-1352.
[217] Potrebbe trattarsi dell'ing. Carlo Jorio.

Se vieni a Torino nella settimana prossima, ricordati di avvertirmi. Intanto ti saluto caramente.

Tuo Aff$^{mo}$
C. Somigliana

**186**
Gian Antonio Maggi a Carlo Somigliana

Viareggio (Lucca) 10 Luglio 1916

Carissimo Somigliana,

Pur dispiacendomi di trattenerti in un campo, dei cui frutti mi dici di essere oramai sazio, non posso esimermi di tornare sull'argomento della mia ultima lettera. Poiché vedo che non mi sono spiegato affatto; ed è un risultato che un professore non potrebbe facilmente sopportare in santa pace. Non mi sono spiegato, perché l'accennato scambio non contava, nel mio discorso, che come materia affatto secondaria, e io non ho punto inteso di muovere alcuna difficoltà al tuo lavoro. Che se si presta alle riserve, a cui accenni, io non mi sono abbastanza inoltrato nei particolari per rilevarle.

Il dubbio che ti pregavo di chiarirmi è questo: se il tuo risultato, che, colla discontinuità dello spostamento, risulta determinata la discontinuità delle derivate prime e seconde dello spostamento medesimo, non contrasta, o almeno non limita, la possibilità, degli spostamenti di Volterra. Per questi spostamenti, infatti, la discontinuità non può essere altrimenti che rappresentata da uno spostamento rigido, con che risulta una discontinuità di forma data; e di più si domanda che siano continui i coefficienti di deformazione, e le loro derivate prime e seconde.

Non segue quindi dall'accennato tuo risultato che le condizioni con cui si definisce lo spostamento di Volterra riescono sovrabbondanti, salvo che la suddetta triplice continuità non sia implicita nella supposta forma della discontinuità? E vi sta implicita? O io m'inganno?

Ho ricevuto la cara tua stamane al momento di lasciare Pisa per Viareggio. Non è facile che ti raggiunga sulle Alpi, ma spero di poter combinare di vederci a più modesti livelli.

Intanto ti rinnovo le più cordiali cose e mi confermo

Tuo aff$^{mo}$
G.A. Maggi.

**187**
Carlo Somigliana a Gian Antonio Maggi
{busta intestata: Chiar^mo Signor Prof. G. Antonio Maggi
(dell'Università di Pisa) - Viareggio - Via Zanardelli, 54 (Lucca)]

Casanova Lanza (Como) 20.VIII.916

Carissimo Maggi,

Credo di avere ora completamente afferrato la tua osservazione; e mi sembra perfettamente giusta.

Dato un corpo, molteplicemente connesso come lo considera il Volterra, e prescritto un certo spostamento relativo rigido lungo un certo taglio, la distorsione (come io la definisco) è completamente determinata. Vale a dire esiste una sola deformazione per cui si verifica la discontinuità data, ed è per di più tale che lungo le due faccie *[sic!]* del taglio, le funzioni classiche si facciano equilibrio.

Queste potranno quindi, secondo i casi, avere derivate prime e seconde continue o discontinue; non è in nostro arbitrio di assegnare tali continuità, o discontinuità.

Ne viene di conseguenza che potrebbero esistere distorsioni, determinate da spostamenti rigidi, aventi funzioni con derivate discontinue. In questi casi le distorsioni del Volterra non esisterebbero. E si dovrebbe quindi concludere - mettendosi nel suo ordine di idee - che per tali corpi, non semplicemente connessi, le distorsioni non esistono.

Sarebbe interessante verificare con qualche esempio la possibilità di queste nuove discordanze con ciò che la più semplice intuizione ci suggerisce.

Ti ringrazio di aver richiamato la mia attenzione su questi fatti, che mi confermano sempre più nel mio modo di vedere in tale argomento.

Mille auguri di tranquille vacanze sulla bella spiaggia di Viareggio. Io andrò in montagna la settimana prossima.

Tuo aff.^mo
C. Somigliana

**188**
Gian Antonio Maggi a Carlo Somigliana

Viareggio 27 Luglio 1916

Carissimo Somigliana,

Ti ringrazio della tua risposta. Se la mia osservazione è giusta, mi pare ne segua che le distorsioni di Volterra non possono esistere che eccezionalmente, salvo che si levi la sovrabbondanza dalle condizioni, col dimostrare che dalla discontinuità dello spostamento, determinata con spostamento rigido della faccia del taglio, segue, oltre la continuità dei coefficienti di deformazione (o di tensione), anche quella delle loro derivate prime e seconde.

Per conto mio, tu hai messo le cose a posto col trasportare la condizione di continuità alle
$X_V, Y_V, Z_V$ mediante le

$$X_V+X_{V'}=0 \qquad Y_V+Y_{V'}=0 \qquad Z_V+Z_{V'}=0.$$

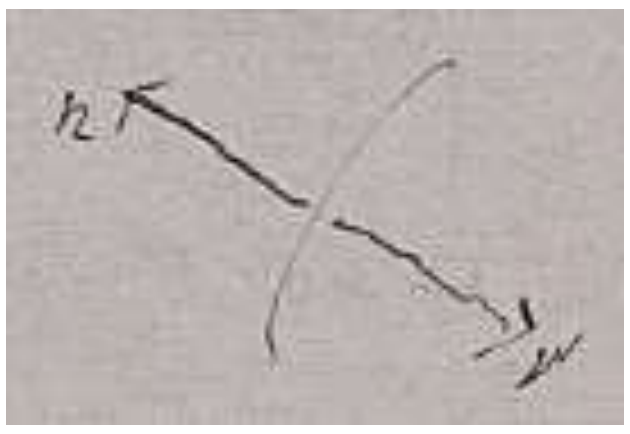

Mentre il Weingarten domanda la continuità della

$$X_x, Y_y, Z_z, Y_z, Z_x, X_y,$$

come se la discontinuità di codeste conducesse alla incongruenza, che recherebbe invece la discontinuità di quelle. Non mi sentirei d'affermare che il Weingarten sia caduto in un equivoco: ma la mia impressione resta che il problema dell'equilibrio elastico per spostamento discontinuo, da quella condizione non giustificata da necessità (come tu benissimo rilevi) sia stato mandato fuori di carreggiata.

**189**

Carlo Somigliana a Gian Antonio Maggi
[cartolina illustrata da Macugnaga indirizzata a: Chiar.$^{\text{mo}}$ Signor
Prof. G.A. Maggi (dell'Un.$^{\text{tà}}$ di Pisa) - Viareggio]

Macugnaga 9.VIII.16

Ricevo qui la tua carm̃a[218]. Siamo perfettamente d'accordo; il tuo concetto merita di essere sviluppato
Cordialmente tuo C. Somigliana

**190**

Gian Antonio Maggi [a Carlo Somigliana]
[appunti delle varie Note sulle distorsioni elastiche]

N.B.

Nella sua Nota "Sull'equilibrio dei corpi elastici più volte connessi" nei Rendic. dei Lincei, Vol. XIV Fasc 4 (1905), a pag. 193, il Volterra (pag. 199) dimostra un teorema, nel quale, tracciata in un corpo elastico isotropo una superficie a piacere, si può assegnare uno spostamento d'equilibrio elastico, che, sulla superficie, presenta la discontinuità, rappresentata da un spostamento rigido della faccia del taglio relativa alla stessa superficie, e possiede "deformazione regolare" - cioè coeff. di deformazioni continui, unitamente colle derivate prime e seconde. L'eccezione pel contorno, se dopo estendendo la superficie fuori del corpo, come lo stesso Volterra dice di fare in principio al seguente §4.

Questo risultato sembra non conciliarsi colla prima Nota dello stesso Volterra "Un teorema sulla teoria dell'Elasticità" (Ibid. pag. 127), dove è dedotto il risultato che, se lo spazio è semplicemente commesso, è sufficiente supporre i coefficienti di deformazione continuo, unitamente colle derivate prime e seconde, perché lo spostamento elastico sia nullo quando sono nulle le forze di massa e le pressioni (pag. 128). Ciò che implica che dall'ipotesi stessa sulle componenti di deformazione segue, nel caso del corpo semplicemente connesso, la continuità dello spostamento. Ciò che si trova nella mia Nota "Sugli spostamenti elastici discontinui" (Rendic. dei Lincei 1908).

Segue retro

D'altra parte, non saprei come giustificare l'asserzione dello stesso Volterra che le u, v, w, definite cogli integrali $\int_{\sigma}$, o $\int_{\sigma'}$, dove σ e σ' rappresentano due calotte componenti una superficie chiusa, soddisfacciano, salvo il contorno, le equazioni della elasticità, come gl'integrali estesi all'intera superficie chiusa (V. pag. 198, in fondo). La deduzione delle formole di Somigliana, che ho riportato, non estende ai pezzi di contorno chiuso la proprietà di soddisfare delle equazioni, propria dell'interno contorno: né, come sarebbe il caso della formola di Green, questa proprietà emerge dalla funzione integranda, in modo da poter estenderla ad un pezzo qualsivoglia del contorno.

[218] Carissima.

Pisa 27 Settembre 1916.

**191**[219]

Carlo Somigliana a Gian Antonio Maggi

Torino 17.XI.916

Carissimo Maggi,

Ti ho spedita la Memoria del Prof. Grdina sulla relatività,[220] e per compensarti del ritardo, te ne ho aggiunto un'altra dello stesso autore e nella stessa lingua. Tu mi dirai poi a che cosa si riferisce.

Domenica sera, a Roma, mi sono incontrato con Almansi. È stata una fortuna perché ho potuto constatare che del deplorato crimine egli ha una responsabilità limitata, per quanto egli sostenga che la introduzione della Memoria della Sig.[na] Armanni non contiene che una inesattezza di espressione. Spero che sarà possibile intenderci per finire onorevolmente la vertenza, senza spargimento di sangue. Quando io, ad un certo punto, ho esclamato: quel Volterra, in certe cose, è un gran testone! l'Almansi ha approvato con entusiasmo. Non siamo quindi lontani dall'intenderci.[221]

Circa quella possibilità di derivare le relazioni che esprimono le discontinuità per ottenere le discontinuità delle derivate, ho riguardato i miei appunti in questi giorni e ritengo che la operazione sia lecita per la prima derivazione, in concordanza colla dimostrazione, di cui mi hai parlato. Ma non lo sia più per le derivate successive. Ma converrà studiare ancora la quistione.

Ti saluto caramente

Tuo aff[mo]

Carlo Somigliana

---

[219] Le prossime tre lettere sono racchiuse, unitamente a quattro fogli di appunti del Maggi riguardo alle "formule di Somigliana", in un fascicolo intitolato: "*Sulla derivabilità di* $f^{+} - f^{-} = \varphi(x)$. *Corrispondenza con Somigliana. Dicembre 1916. (Calcolo di Pizzetti).*"

[220] Docente all'Ekaterinoslavl' Higher Mining School, a partire dal 1912 pubblicò una serie di articoli contrari alla Teoria della Relatività; ad es. Ia.I. Grdina, "K voprosu o masse elektrona" [A proposito della questione della massa dell'elettrone], Ekaterinoslav, 1912.

[221] A dicembre, Somigliana pubblica la Nota "Sulla teoria delle distorsioni. Al prof. E. Almansi", *Rendiconti della R. Accademia dei Lincei,* v. XXV, 1916, pp. 475-477 in risposta a quella di E. Almansi, "La teoria delle distorsioni e le deformazioni finite dei solidi elastici", *ib.,* pp. 190-191 che era stata presentata il 26 settembre. Probabilmente non si sono intesi.

**192**

Carlo Somigliana a Gian Antonio Maggi

Torino 5 dic. 916

Caro Maggi,

Mille grazie delle numerose tue comunicazioni. Quanto mi dici della Memoria del Grdina[222] sulla relatività mi pare interessante, e mi parrebbe assai utile che tu ne facessi argomento di una recensione da pubblicarsi nel Nuovo Cimento. Il russo è così poco conosciuto da noi, e tanto meno diffuse sono le sue pubblicazioni scientifiche in quella lingua.

Le tue dimostrazioni mi paiono soddisfacenti; né saprei in qual modo fare obiezioni, ma ecco il fatto che contraddice alla possibilità di scambio delle due operazioni di limite, derivazione.

Da una funz. pot. di superficie

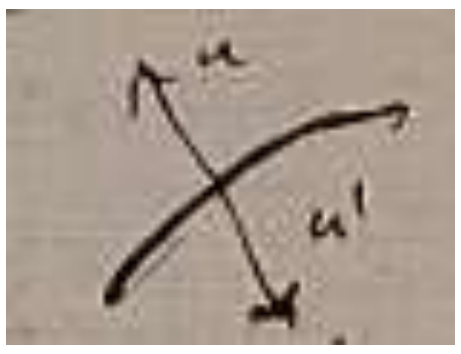

$$V=\int h\frac{ds}{r}$$

abbiamo per la continuità della funzione

$$V_u=V_{u'}$$

e derivando tangenzialmente

$$\frac{\delta V_u}{\delta x}=\frac{\delta V_{u'}}{\delta x} \qquad \frac{\delta V_u}{\delta y}=\frac{\delta V_u}{\delta y}$$

equazioni che sono giuste, cogli assi che io chiamo canonici. Ma derivando un'altra volta si dovrebbe avere

$$\frac{\delta^2 V_u}{\delta x^2}=\frac{\delta^2 V_{u'}}{\delta x^2} \qquad \frac{\delta^2 V_u}{\delta y^2}=\frac{\delta^2 V_{u'}}{\delta y^2} \qquad \frac{\delta^2 V_u}{\delta x\delta y}=\frac{\delta^2 V_{u'}}{\delta x\delta y}$$

Ora la terza di queste equazioni è giusta; ma non lo sono le prime due, poiché le discontinuità di quelle derivate seconde sono date invece dalle relazioni

$$D\left[\frac{\delta^2 V}{\delta x^2}\right]=-\frac{4\pi}{R_1}h$$

$$D\left[\frac{\delta^2 V}{\delta y^2}\right]=-\frac{4\pi}{R_2}h$$

$R_1$, $R_2$, essendo i raggi di curvatura.

Queste proprietà si possono verificare sulla funz. pot. della sup. sferica omogenea. Era questo l'esempio che io avevo considerato dapprima.

A te la soluzione di questa non facile quistione.

Empiricamente si trova che la prima derivazione è lecita; così pel caso precedente della funz. pot. di superficie ed anche per quella di doppio strato, poiché da

$$W_u-W_{u'}=4\pi g$$

si ricavano le due relazioni vere

$$\frac{\delta W_u}{\delta x}-\frac{\delta W_{u'}}{\delta x}=4\pi\frac{\delta q}{\delta x} \qquad \frac{\delta W_u}{\delta y}-\frac{\delta W_{u'}}{\delta y}=4\pi\frac{\delta q}{\delta y}$$

Ma, come vedi, procedendo oltre nascono i guai.

---

[222] Si veda la lettera precedente.

Nelle mie Note le cose alla fine vanno bene, poiché non si fa uso dei valori speciali delle discontinuità, ma solo del fatto che quelle discontinuità risultano determinate. Così almeno io credo.
Saluti cordiali

Tuo aff$^{mo}$
C. Somigliana

**193**
Gian Antonio Maggi a Carlo Somigliana

Pisa 23 Dicembre 1916

Carissimo Somigliana,

Prima di tutto, augurii per le Feste, e ringraziamenti per la casa tua del 5 corr.
A proposito della quale, credo di essere venuto a capo della questione relativa alle derivabilità di

$$f^{+}-f^{-}=\varphi,$$

rispetto ad una coordinata superficiale.
Sta, per qualunque valore di n, la relazione

$$\left(\frac{\delta^n f}{\delta u^n}\right)^{+}-\left(\frac{\delta^n f}{\delta u^n}\right)^{-}=\left(\frac{\delta^n \varphi}{\delta u^n}\right),$$

per una coordinata curvilinea u, e, introducendo il relativo arco s.

(1)
$$\left(\frac{\delta^n f}{\delta s^n}\right)^{+}-\left(\frac{\delta^n f}{\delta s^n}\right)^{-}=\left(\frac{\delta^n \varphi}{\delta s^n}\right),$$

quando si verifichino le condizioni d'esistenza della derivata, e di continuità, dalle due parti della superficie, indicate nell'accennato teorema.

Ma la coordinata cartesiana x, colla condizione che l'asse delle x sia tangente alla superficie, nel punto considerato, non è equiparabile alla suddetta s. Ed è la ragione per cui la (1), dove si faccia s=x, sta per n=1, ma non sta più, cominciando da n=2.
Difatti si ha

$$\frac{\delta f}{\delta s}=\frac{\delta f}{\delta x}\frac{\delta x}{\delta s}+\frac{\delta f}{\delta y}\frac{\delta y}{\delta s}+\frac{\delta f}{\delta z}\frac{\delta z}{\delta s},$$

per cui se l'asse delle x è tangente alla linea u, o s,

$$\frac{\delta f}{\delta s}=\frac{\delta f}{\delta x}.$$

Ma

$$\frac{\delta^2 f}{\delta s^2}=\frac{\delta^2 f}{\delta x^2}\left(\frac{\delta x}{\delta s}\right)^2+\ldots+2\frac{\delta^2 f}{\delta y\delta z}\frac{\delta y}{\delta s}\frac{\delta z}{\delta s}+\ldots\ldots+\frac{\delta f}{\delta x}\frac{\delta^2 x}{\delta s^2}+\frac{\delta f}{\delta y}\frac{\delta^2 y}{\delta s^2}+\frac{\delta f}{\delta z}\frac{\delta^2 z}{\delta s^2}$$

per cui se l'asse delle x è tangente alla linea u, o s,

$$\frac{\delta^2 f}{\delta s^2}=\frac{\delta^2 f}{\delta x^2}+\frac{\delta f}{\delta x}\frac{\delta^2 x}{\delta s^2}+\frac{\delta f}{\delta y}\frac{\delta^2 y}{\delta s^2}+\frac{\delta f}{\delta z}\frac{\delta^2 z}{\delta s^2}.$$

Ossia, indicando con $\alpha_2$, $\beta_2$, $\nu_2$, i coseni di direzione della normale principale alla linea, volta verso il centro di curvatura, e con R>0 il raggio di curvatura.

$$\frac{\delta^2 f}{\delta s^2}=\frac{\delta^2 f}{\delta x^2}+\frac{1}{R}\left(\frac{\delta f}{\delta x}\alpha_2+\frac{\delta f}{\delta y}\beta_2+\frac{\delta f}{\delta z}\nu_2\right).$$

Se la linea è linea di curvatura, relativa al punto considerato e l'asse delle z normale alla superficie e volto verso il centro di curvatura principale, corrispondente alla suddetta linea di curvatura, indicando con $R_1$ il relativo raggio di curvatura principale, si ha di qua, pel Teorema di Meunier,[223]

(2) $$\frac{\delta^2 f}{\delta s^2}=\frac{\delta^2 f}{\delta x^2}+\frac{1}{R}\left(\frac{\delta f}{\delta x}\alpha_2+\frac{\delta f}{\delta y}\beta_2\right)+\frac{1}{R_1}\frac{\delta f}{\delta z}.$$

Sia ora

$$f=\int_\sigma \frac{h d\sigma}{r}.$$

Allora

$$D\frac{\delta f}{\delta x}=0 \qquad D\frac{\delta f}{\delta y}=0 \qquad D\frac{\delta f}{\delta z}=-4\pi h \qquad \frac{\delta^2 f}{\delta s^2}=0,$$

Con che la (2) fornisce

$$D\frac{\delta^2 f}{\delta x^2}=\frac{4\pi h}{R_1},$$

che è il risultato da te indicato, tenendo conto della diversa ipotesi relativa al segno di $R_1$. Nello stesso caso da

$$f^+-f^-=0,$$

si ricava senz'altro, nell'ipotesi dell'asse delle x tangente

$$\left(\frac{\delta f}{\delta x}\right)^+-\left(\frac{\delta f}{\delta x}\right)^-=0.$$

Ma, per proceder oltre, deriveremo

(3) $$\left(\frac{\delta f}{\delta x}\right)^+-\left(\frac{\delta f}{\delta x}\right)^-=-4\pi h\cos\widehat{nx},$$

con che si ha

(4) $$\left(\frac{\delta^2 f}{\delta x^2}\right)^+-\left(\frac{\delta^2 f}{\delta x^2}\right)^-=-4\pi\frac{\delta h}{\delta x}\cos\widehat{nx}-4\pi h\frac{\delta\cos\widehat{nx}}{\delta x}.$$

Introducendo in questa relazione l'ipotesi che l'asse delle x sia tangente ad una linea di curvatura, si ricava poi subito sul precedente risultato.

Dalla (4) dell'asse che si azzera derivando (3) rispetto a y, si ricavano la formole (8) della tua Nota "Sulle Derivate seconde della fz. Pot. di Superficie". Mi son rivolto per questo all'amico Pizzetti, che ha maggior famigliarità di me coi lavori in superficie, e jeri sera mi ha portato il non breve calcolo, condotto con molta maestria ed eleganza. È quanto aspettavo per scriverti, perché non mancasse una verificazione delle formole generali.

Conto sul piacere di vederti dopo Natale, a Milano o a Torino. Intanto ti rinnovo gli augurii, e coi migliori saluti di noi tutti e una cordiale stretta di mano ti prego di credermi sempre

Tuo Aff$^{mo}$ G.A. Maggi

[223] Meusnier de la Place.

**194**
Gian Antonio Maggi a Carlo Somigliana

Pisa 1 Ottobre 1917

Carissimo Somigliana,

Ti ringrazio della tua cartolina, in data del 27.[224] Anche il mio ricordo, se non è un sogno matematico, non è più che di qualche cenno sull'indicata questione. Questa mattina stessa ho consegnato il M.S della relativa Nota al "Nuovo Cimento",[225] previo passaggio per le mani del nostro Levi Civita, in osservanza del precetto: e, per quanto mi scrive, pare che non gli sia dispiaciuto. Non vale però la pena di discorrere, prima di mandartela stampata.

La tua cartolina deve essersi incrociata con una mia, in data del 29, diretta a Casanova Lanza. Ti comunico, senz'altro, i risultati che io ho trovato, per la discontinuità delle derivate seconde della funzione potenziale di doppio strato $W = \int_\sigma g \frac{\delta \frac{1}{s}}{\delta n} d\sigma$.

Mi valgo delle formule

$$D\frac{\delta^2 W}{\delta x^2} = \frac{\delta^2 DW}{\delta s^2} - \frac{\delta DW}{\delta y}\frac{\delta^2 y}{\delta s^2} - D\frac{\delta W}{\delta z}\frac{\delta^2 z}{\delta s^2}$$

$$D\frac{\delta^2 W}{\delta x \delta y} = \frac{\delta}{\delta s}\frac{\delta}{\delta s'} DW - \frac{\delta DW}{\delta x}\frac{\delta}{\delta s}\frac{\delta x}{\delta s'} - D\frac{\delta W}{\delta z}\frac{\delta}{\delta s}\frac{\delta z}{\delta s'}$$

$$D\frac{\delta^2 W}{\delta x \delta z} = \frac{\delta}{\delta x} D\frac{\delta W}{\delta z},$$

nelle quale gli assi delle x e delle y sono supposti avere la direzione delle tangenti alle linee s e s'. Vi introduco

$$DW=4\pi g, \qquad D\frac{\delta W}{\delta z} = D\frac{\delta W}{\delta n} = 0,$$

e nell'ipotesi che le due linee siano le linee di curvatura trovo

$$\frac{1}{4\pi} D\frac{\delta^2 W}{\delta x^2} = \frac{\delta^2 g}{\delta s_1^2} - \frac{R_1 R_2}{R_2 - R_1}\frac{\delta \frac{1}{R_1}}{\delta s_2}\frac{\delta g}{\delta s_2}$$

$$\frac{1}{4\pi} D\frac{\delta^2 W}{\delta x \delta y} = \frac{1}{2}\left(\frac{\delta}{\delta s_1}\frac{\delta g}{\delta s_2} + \frac{\delta}{\delta s_2}\frac{\delta}{\delta s_1}\right) - \frac{1}{2}\frac{R_1 R_2}{R_2 - R_1}\left(\frac{\delta \frac{1}{R_1}}{\delta s_2}\frac{\delta g}{\delta s_1} + \frac{\delta \frac{1}{R_2}}{\delta s_1}\frac{\delta g}{\delta s_2}\right)$$

$$D\frac{\delta^2 W}{\delta x \delta z} = 0.$$

[224] Non è presente nel Fondo, così come quella del 29 nominata in seguito.
[225] Potrebbe trattarsi della Nota: "Sul collegamento di una funzione data nei punti di una superficie con una funzione dei punti dello spazio, e sua applicazione alla teoria della funzione potenziale di doppio strato", *Il Nuovo Cimento*, 1917, s. VI, t. XIV, pp. 5-11.

E non penso che basterebbe dedurre, in modo analogo, l'espressione corrispondente per $D\frac{\delta^2 W}{\delta x \delta y}$ perché non occorresse altro. Precisamente quello che fai per le fz pot. di superf. nel qual caso ha verificato un accordo completo.

La prima si riduce alla prima delle tue formule a pag. 8 della Nota Sulle Derivate Sec. della Fz. Pot. di Doppio Strato, facendovi pel primo membro $\frac{\delta^2 W}{\delta x^2}$ in luogo di $\frac{\delta^2 W}{\delta s_1^2}$: la seconda alla *[*quarta*]*, facendovi, pel primo membro, $\frac{\delta^2 W}{\delta x \delta y}$ in luogo di $\frac{\delta^2 W}{\delta s_1 \delta s_2}$, e nel secondo membro $\frac{1}{2}\left(\frac{\delta}{\delta s_1}\frac{\delta g}{\delta s_2}+\frac{\delta}{\delta s_2}\frac{\delta g}{\delta s_1}\right)$ in luogo di $\frac{\delta^2 g}{\delta s_1 \delta s_2}$. La terza, che si può scrivere

$$D\frac{\delta^2 W}{\delta s_1 \delta n}=0$$

non si concilia colla quinta delle formole suddette.

Io non saprei trovar un errore nelle mie formole che ho ripetutamente controllato e variamente comprovato. Mi pare invece che nelle formole della tua Nota si deve scrivere $\frac{1}{2}\left(\frac{\delta}{\delta s_1}\frac{\delta g}{\delta s_2}+\frac{\delta}{\delta s_1}\frac{\delta g}{\delta s_2}\right)$ invece di $\frac{\delta^2 g}{\delta s_1 \delta s_2}$, perché definito $\frac{\delta}{\delta s_1}\frac{\delta}{\delta s_2}$ con $\frac{1}{\sqrt{E}}\frac{\delta}{\delta u}\frac{1}{\sqrt{G}}\frac{\delta}{\delta v}$ e $\frac{\delta}{\delta s_2}\frac{\delta}{\delta s_1}$ con $\frac{1}{\sqrt{G}}\frac{\delta}{\delta v}\frac{1}{\sqrt{E}}\frac{\delta}{\delta u}$ non è generalmente $\frac{\delta}{\delta s_1}\frac{\delta}{\delta s_2}=\frac{\delta}{\delta s_2}\frac{\delta}{\delta s_1}$.

Con questa modificazione credo che la formola in fine a pag. 7 regga perfettamente. Ma, per quanto all'uso di

$$\frac{\delta^2 W}{\delta x^2}=(\alpha_1 Ds_1+\alpha_2 Ds_2+\alpha D_n)^2 W$$

per dedurre le rimanenti formule, oltre la precedente circostanza, e quella che il secondo membro mi sembra *[...]* un termine, relativo alle variazione di $\alpha_1$, $\alpha_2$, $\alpha$ con x, noto che non si tratta di $D\frac{\delta^2 W}{\delta s_1^2}, D\frac{\delta^2 W}{\delta s_1 \delta s_2}$ che sono $\frac{\delta^2 DW}{\delta s_1^2}, \frac{\delta^2 DW}{\delta s_1 \delta s_2}$, ma di $D\frac{\delta^2 W}{\delta x^2}, D\frac{\delta^2 W}{\delta x \delta y}$, che ne differiscono come risulta dalle formule ( ).

Io destino la nuova Nota,[226] come si diceva, a completare la prima, che, per quanto esatta, non mi è riuscita troppo ben redatta, e penso che la forma possa compromettere la sostanza, la quale mi pare invece meritevole di qualche attenzione, come dimostrano le indicate applicazioni. Potrei terminarla così: "Formole per la discontinuità etc, per quanto ci consta, furono date per la prima volta da Somigliana (*)[(*) Nota etc.]. Conformemente ai precedenti risultati, vanno ad esse recate le modificazioni, che risultano dal confronto colle nostre ( )".[227]

---

[226] "Nuove applicazioni di una formula commutativa", *Rendiconti della R. Accademia dei Lincei*, v. XXVI, 1917, pp. 201-206; la prima Nota alla quale Maggi si riferisce poco dopo è "Sopra un formula commutativa e alcune sue applicazioni", *ib.*, pp. 189-194.

[227] A. p. 206 si legge: "Formole nuove, per quanto consta all'autore, per la discontinuità delle derivate seconde della funzione potenziale di doppio strato furono date da Somigliana, in una sua recente Nota. Stando ai precedenti risultati, vanno recate a tali formole le modificazioni che risultano dal confronto colle nostre (22)-(26)."

Ma se le mie condizioni ti persuadono, e credi di servirtene, per tornare sull'argomento con un nuovo scritto, io sono ben disposto ad aspettarlo (purché non sia alla fine della guerra europea!) per constatare l'accordo.
*[a matita:* Seguono i miei progetti di viaggio etc.*]*

**195**[228]
Gian Antonio Maggi a Carlo Somigliana

*[in matita turchina:* Per la cartella "Giudizii"*]*

Cireglio (Pistoia) 22 Agosto 1922

Carissimo,

Grazie del gentile invio della tua Nota[229]. Essa non aveva mancato di richiamare la mia attenzione, quando comparve nei Rendiconti dei Lincei, e l'avevo attentamente letta, riservandomi di farti presenti alcune mie riflessioni, che non trovo da mutare alla rinnovata lettura.

La principale di queste è che la trasformazione diciamo pure di Voigt implica, col mutamento delle coordinate, anche quello del tempo. Non ho potuto procurarmi la Nota del Voigt, per cercare quale interpretazione egli presta a questo mutamento, colpa la chiusura estiva della Scuola Normale, che rende malagevole - non dirò impossibile - ricorrere alle raccolte di quella Biblioteca. Nella tua Nota non ne trovo cenno. Ora, a me sembra essenziale, per attribuire un significato alla trasformazione di Lorentz, o di Voigt, nella ordinaria meccanica newtoniana, conciliare colle definizioni newtoniane, il diverso valore del tempo in riferimenti diversi.

La priorità del Voigt è incontestabile pur di discendere alla trasformazione dell'equazione di d'Alembert. Giova però tener presente che il Lorentz ha indicato la trasformazione, che si nomina da lui, per la trasformazione in sé stessa delle equazioni di Maxwell-Hertz, e che questa è la prima applicazione che ha conferito al risultato del Lorentz il suo cospicuo valore. Quindi, se sta bene ricordare il Voigt, trovo pur giusto mantenere al Lorentz un posto conforme al suo ben più ampio contributo.

È poi certamente interessante il risultato che la trasformazione di Lorentz si presenta come un caso specialissimo delle trasformazioni rappresentate dalle formule (3),[230] *[richiamo a margine destro di pagina:* * Questo caso però è il solo che serve per la teoria delle relatività, la quale pone la restrizione che, colla trasformazione delle coord. e del tempo, si trasformi in sé stessa l'espressione $x^2 - y^2$. Stabilito che sia $x^2 - y^2 = x'^2 - y'^2$, le (3) forniscono $\varphi(x'+y')\psi(x'-y') = x'^2 - y'^2$ ossia $\varphi(u)\psi(v) = uv$, donde necessariamente $\varphi(u) = au, \psi(v) = \frac{1}{a}v$, che sono le (4).[231]*]*

[228] Da qui in poi le lettere sono state dallo stesso Maggi racchiuse in un fascicolo denominato: *Corrispondenza con Somigliana sulla sua Nota "Sulla trasformazione di Lorentz" (Rendic. dei Lincei, 21 Maggio 1922).*
[229] Si riferisce alla Nota: C. Somigliana, "Sulla trasformazione di Lorentz", *Rendiconti della R. Accademia Nazionale dei Lincei*, v. XXXI, 1922.
[230] *Ibid.*, p. 409: $x + y = \varphi(x'+y')$ $x - y = \psi(x'-y')$.
[231] *Ibid.*, p. 412: $x + y = a(x'+y')$ $x - y = b(x'-y')$, dove a e b sono costanti non nulle.

*[parte cancellata da Maggi:*
Ma la trasformaz. di Lorentz, come apparisce dalla relazione a riga 6 di pag. 416, è la sola , tra quelle, che trasforma in sé stesso il primo membro della equaz. di d'Alembert, ciò che torna trasformare in sé stesso $x^2 - y^2$ (o $x^2 - ct^2$): e questa è circostanza sostanziale per la teoria della relatività. Vedasi, ad esempio, Becquerel, Le Principe de Relativité, pag. 30.*]*

La trasformazione considerata in II non soddisfa a questa condizione, e non potrebbe trovar posto nella teoria della relatività.

Sarà con grande interesse e piacere che sentirò le tue osservazioni su queste mie riflessioni, che, se non altro, ti faranno fede della mia attenta lettura.

Ci siamo lasciati a Roma senza più ritornare sulle onde di forma qualsivoglia. Ci siamo messi d'accordo? È pur singolare che, all'infuori del Levi-Civita e del Laura, che s'è alquanto lagnato perché non lo ricordai nella mia Nota[232] (l'ho poi nominato esplicitamente in un breve riassunto, che il prof. Vallauri mi chiese pel Giornale di Elettrotecnica)[233] nessuno sembri aver preso in adeguata considerazione una questione importante come la conciliazione dell'ottica geometrica coll'ottica fisica. La Relatività non lascia più tempo per occuparsi d'altro!

Ti auguro buona villeggiatura, dovunque tu ti trovi. Noi contiamo restare quassù tutto il Settembre. I migliori saluti dai miei, e una cordiale stretta di mano dal

Tuo aff.mo
G.A. Maggi.

**196**
Carlo Somigliana a Gian Antonio Maggi

Casanova Lanza (Como) 20 Sett. 22

Carissimo Maggi,

ti sono grato della buona volontà con cui hai preso in esame la mia Nota. In un argomento come questo è necessario discutere con un competente, che non sia, nello stesso tempo, un infatuato relativista. Con quelli di questa categoria si ragiona male. È necessario invece porsi da un punto di vista obiettivo e prendere con risolutezza per quel che valgono i postulati della relatività.

Con le tue osservazioni io concordo in via generale. In quanto però riguarda il mutamento del tempo io non sento la necessità di una interpretazione fisica. La trasformazione ha sopratutto *[sic!]* valore analitico, nel senso che lega fra loro gli integrali della eq.e di propagazione. Quello che avviene in un dato tempo e luogo, quando la sorgente è fissa, si riproduce in un altro tempo e luogo quando la sorgente si muove. Il tempo concettualmente non cambia; resta sempre immutato, invariabile, eterno.

---

[232] "Sulla propagazione delle onde di forma qualsivoglia nei mezzi isotropi", *Rendiconti della R. Accademia Nazionale dei Lincei*, v. XXIX, 1920, pp. 371-378; si veda la lettera #56 nella quale Maggi spiega a Daniele come fosse nato lo scritto in questione.

[233] A p. 296 del n. 13 (1922) di tale *Giornale* si legge: "Il prof. Somigliana, estendendo ad una serie di onde comunque stabilita in un mezzo isotropo, un risultato precedentemente ottenuto dal prof. Laura, per una serie di onde che invadono progressivamente lo stesso mezzo, inteso originariamente in quiete, perviene alla conclusione che non sono possibili che onde sferiche o piane o cilindriche."

Il problema è di trovare degli integrali della eq. di propagazione per cui la sorgente cammini. La trasformazione da il modo di costruirli; ma non è escluso che vi si possa arrivare anche per altra via, e quindi venga quel mutamento apparente del tempo.
Ciò che mi sembra discendere da quelle considerazioni, e che non ho detto ancora in modo esplicito, è che la relatività, inventata per giustificare l'esito negativo della esperienza di Michelson, in realtà non lo giustifica affatto.

La teoria newtoniana della propagazione delle onde infatti non porta per nulla alla conclusione che la velocità di propagazione della luce debba essere diversa per l'osservatore in moto, rispetto a quello fermo, come le fanno dire arbitrariamente i relativisti. La teoria newtoniana, basata sulle equazioni, dice soltanto che il fenomeno nei due casi è diverso; e le modalità di queste differenze sono ben più complicate che non quelle risultanti da un semplice cambiamento di velocità. Del resto anche per l'osservatore newtoniano in moto (che è un essere intelligente) la velocità di propagazione non varia colla sua velocità propria di movimento, come per i relativisti. Analiticamente ciò risulta dalla equazione

$$x^2 - c^2t^2 = x'^2 - c^2t'^2$$

valida pei relativisti ed anche per noi che non lo siamo.

Osservazioni di questo genere ho già fatto a Levi-Civita un po' sommariamente, e non mi è sembrato che avesse obiezioni. Puoi tu indicarmene qualcuna?

Vedo che non ho più spazio per parlare dei tuoi integrali ad onde parallele. Lo farò un'altra volta.

Mille cordiali saluti.

Tuo aff.$^{mo}$
Carlo Somigliana

**197**
Gian Antonio Maggi a Carlo Somigliana

Cireglio (Pistoia) 26 Settembre 1922

Carissimo Somigliana,

Ti sono grato, alla mia volta, dell'attenzione che hai posto alle mie riflessioni, e sento con piacere che oggetto di discussione non resti che il significato del diverso valore del tempo. Ma, come accennavo, è la questione principale. Nessun dubbio intanto che io sia disposto a esaminare le cose da un punto di vista perfettamente obbiettivo. Io non presto alla teoria della relatività che un'adesione relativa. Poiché sta il fatto che quella teoria, attraente per ampiezza di quadro e acutezza di vedute, complica estremamente la trattazione razionale dei problemi fisici, per arrivare a perturbazioni dei risultati della teoria classica che, nella quasi totalità dei casi, sono inaccessibili alla verificazione sperimentale. Gli stessi risultati verificabili si fondano sopra esperienze che l'estrema delicatezza rende discutibili, e la celebre coincidenza dei 43" dello spostamento del perielio di Mercurio è pur da reputarsi non altrimenti che casuale, se si considera che lo spostamento effettivo è di circa 550", dei quali rende ragione fino a 43" l'ordinaria teoria newtoniana. Per modo che non si saprebbe il numero che fornirebbe la teoria della relatività, applicata alla spiegazione del fenomeno completo, tenendo conto cioè del contributo dei pianeti alla modificazione della geometria del $S_4$ cosmico, quando fosse possibile cavar qualcosa da quelle complicatissime equazioni.

Io aveva preparato un'osservazione di questo genere per il referendum promosso dal Hoepli; ma vi rinunciai, constandomi che simile obbiezione fu mossa da altri.

Porrei ora la questione così. Si tratta di trovare una sostituzione delle x,y,z,t nelle x',y',z',t', che trasforma in sé stesso il dalambertiano □. Occorre questa limitazione, se si vuol servirsi della stessa sostituzione per la trasformazione in sé stessa dell'equazione completa □φ=ψ, dove ψ è funzione prestabilita. Si trova, sufficiente e necessaria, la trasformazione di Lorentz, come, ridotto il caso alla sostituzione di x,t in x',t', emerge dalla formola generale a linea 6, pag. 410 della tua Nota. Per le applicazioni, occorre attribuire un significato concreto alle x',t', in confronto delle x,t. La teoria della relatività parte dell'ipotesi che x' rappresenta lo stesso punto come x, e t' lo stesso istante come t. Tu salvi il tempo newtoniano “immutato, invariabile, eterno”, col intendere che t' rappresenti un istante, in massima, diverso dall'istante rappresentato da t.

A questo proposito, io trovo da domandare che cosa significano più le formole che formalmente coincidono con quelle della teoria della relatività, e sono da questa teoria applicate, al noto moto, per la spiegazione dei fenomeni? Che cosa significa, ad esempio, la formola

$$v = \frac{v'+\lambda c}{1+\frac{\lambda}{c}v'},$$

in un moto vario, quando v e v' debbano assegnarsi ad istanti diversi (corrispondenti a t e a t'). Sta bene che, per l'applicazione al fenomeno di Fizeau,[234] serve l'ipotesi di v e v' indipendenti dal tempo, inteso il moto dell'acqua uniforme. Ma dobbiamo limitare la formola a questo caso? E per la stessa spiegazione dell'esperienza di Michelson come terrà più l'equazione

$$x^2 - c^2t^2 = x'^2 - c^2t'^2,$$

quando t e t' non s'intendono più rappresentare lo stesso tempo? Come tuttavia c rappresenta la velocità della luce per ambedue i riferimenti?

Più o meno facile che sia rispondere, è certo che queste e simili domande si presentano naturalmente. E a te si offre il compito di completare con una precisa risposta il nuovo significato delle formole relativistiche, che proponi colla tua Nota. Non è piccola posta della partita rimettere interamente sul piedestallo la dottrina newtoniana!

Resteremo qua fino ai primi d'Ottobre. Conto trovarmi a Milano il 7, per le nozze del figlio maggiore di mio fratello Carlo. Se una visita mezzo promessa, della mia figliola maritata[235] non mi richiamerà immediatamente a Pisa, spero vederti a Milano o a Torino. Intanto ti presento i migliori saluti, da parte anche dei miei. Ricordami, con rispetti, alla tua famiglia. Tu ricevi una cordiale stretta di mano dal

Tuo aff.mo
Gian Antonio Maggi.

*[nota cancellata da Maggi:* P.S. Il Levi-Civita, dal quale ho ricevuto lettera jeri l'altro,[236] per quanto mi scrive, divide il mio modo di vedere.*]*

---

234 Si riferisce all'esperienza di Fizeau (1851), mediante la quale si rileva l'influenza che il moto di una corrente d'acqua ha sulla velocità della luce e si misura il cosiddetto “coefficiente di trascinamento”, dando così prova sperimentale alla formula teorica di Fresnel.

235 Potrebbe trattarsi di Clara.

236 Nonè presente nel Fondo.

**198**
Gian Antonio Maggi a Carlo Somigliana
[cartolina]

Cartolina.

Cireglio (Pistoia) 27 Sett. 1922

Carissimo,

A scanso di equivoci, quando, nella mia lettera di jeri, aggiungo “necessario” a “sufficiente”, a proposito della trasformazione di Lorentz, come atto a trasformare in sé stesso il □, intendo subordinatamente a sottintese condizioni, come quella che le nuove variabili siano funzioni lineari intere delle antiche (usata nella relatività) o quella che la sostituzione debba essere compresa nella citata formula generale.

Con una sonora cordiale stretta di mano

Aff$^{mo}$ tuo G.A. Maggi.

Ch.$^{mo}$ Prof. C. Somigliana ecc. Casanova Lanza (Como).

**199**
Gian Antonio Maggi a Carlo Somigliana

Pisa 22 Aprile 1923

Carissimo Somigliana,

Il fascicolo uscito ora dei Rendiconti dei Lincei contiene, come avrai veduto, le mie riflessioni intorno alla tua Nota Sulla trasformazione di Lorentz, alle quali ti accennai a Catania. Ti manderò l'estratto, appena li avrò ricevuti dai Lincei. Ho aspettato, in questi giorni, una tua riga di risposta al biglietto che lasciai per te al Bureau dell'Albergo, al momento di salire in omnibus per la Stazione di Catania. Penso oramai che si siano dimenticati di consegnartelo. Era cosa riguardante l'adunanza del Congresso ch'io fui chiamato a presiedere, ma di così piccola importanza che non vale la pena di ricordarla altrimenti.

Piuttosto, col ricordo delle suddette riflessioni, riprendo, se non ti spiace, il discorso interrotto alla posta dell'Albergo della Grande Brettagna, per chiarire un presumibile equivoco sul significato dell'indipendenza della velocità di propagazione della luce dalla velocità della sorgente. Tu dicevi che dalla teoria del Voigt scaturisce questa indipendenza che la teoria della relatività assume come suo canone particolare. Ora, è pacifico, colla teoria ordinaria, che indipendente dalla velocità della sorgente è la velocità assoluta della propagazione della luce, cioè la velocità di propagazione rispetto al mezzo di propagazione, supposto fisso. Ma, o in conseguenza, colla stessa teoria ordinaria, risulta diversa a seconda del senso e della direzione la velocità di propagazione della luce emanata da un punto mobile relativa ad uno spazio che partecipa alla velocità del mobile, per modo che rispetto ad esso il punto, che funge da sorgente luminosa, è fisso. La conseguenza è immediata se vogliamo ragionare coi criterii comuni, e non inchinandoci a priori ai misteri della relatività. Il segmento AB, possieda una velocità di grandezza v nella direzione e nel verso da A a B. La veloc. di propagazione della luce, nella direzione del segmento e nel senso da A a B, vale a dire nel senso del movim. del segmento, risulta c−v, perché dallo spazio mobile rispetto al quale il segmento è fisso, è assunto come AB il segmento $AB = ct_1 = AB + BB_1 = AB + vt_1$; invece la veloc. di propagaz. della luce, nella stessa direzione e nel senso contrario al precedente, cioè nel senso contrario al movim. del segmento, risulta c+v, perché dal suddetto spazio mobile è assunto come AB il segmento $A_2B = ct_2 = AB - AA_2 = AB - vt_2$. S'intende che c rappresenta la velocità della luce rispetto al mezzo fisso, la quale, nel supposto caso di A mobile si assume la stessa come nel caso che A fosse fisso.

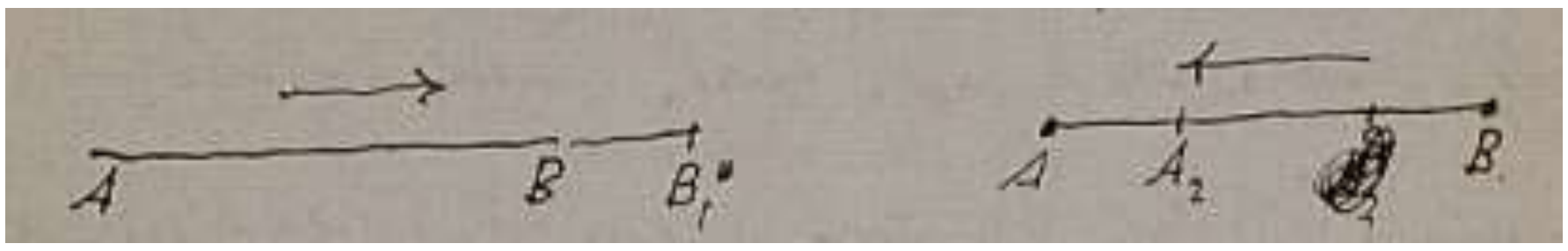

La teoria della relatività postula che la veloc. di propagazione relativa nei due suddetti sensi, come in ogni direzione e senso, abbia il valore comune c, che la teoria ordinaria assegna alla velocità di propagazione, nel mezzo fisso, e rinuncia con questo al meccanismo di un mezzo di propagazione, perché risultano perfettamente equivalenti rispetto alla propagazione della luce lo spazio rispetto al quale il punto A è fisso e quello rispetto al quale è mobile. Questo postulato poi assorbe l'ipotesi della teoria ordinaria che la veloc. di propagaz. sia la stessa per una sorgente fissa come per una sorgente mobile, ma questa ipotesi cambia di significato, riuscendo una conseguenza della equivalenza di due spazii mobili l'uno rispetto all'altro, che è canone fondamentale della teoria della relatività.

Per venire alla conclusione, nella teoria del Voigt, codesto moto di uno spazio rispetto al quale la sorgente mobile è fissa non ha con essa a che vedere, e non se ne fa cenno. Non mi sembra quindi che vi si possa trarre l'indipendenza della veloc. di propagaz. della luce dalla veloc. della sorgente nel significato di veloc. relativa, che appartiene alla teoria della relatività.

Con cordiali saluti ecc.

G.A.M.

**200**

Carlo Somigliana a Gian Antonio Maggi

Torino 28.IV.23

Carissimo Maggi,

Effettivamente a Catania non mi hanno consegnato alcun tuo biglietto, e sono dolente quindi di non aver potuto risponderti. Meno male se si trattava di cosa di poco momento. Se riguarda i verbali delle sedute della Sezione matematica, siamo sempre in tempo a scriverne all'Almagià, che dovrà pubblicare gli Atti del Congresso.
Ti ho mandato ieri le bozze di un mio articolo sulla Relatività destinato alla "Scientia" di Rignano.[237] Da esse vedrai come ho impostato la discussione sulla velocità, costante o variabile, della propagazione della luce. In sostanza, quanto mi scrivi credo che concorda colle mie idee. Ciò che realmente asseriscono i relativisti non ha mai un valore intrinseco, reale; ne ha generalmente soltanto uno formale o meglio verbale, perciò la discussione è difficile con loro. Ma io credo di essermi messo abbastanza al sicuro da difficoltà che derivino da questo stato di cose.

[237] "I fondamenti della relatività", *Scientia*, v. 34, 1923, pp. 1-10.

Ad ogni modo nessuno finora aveva messo in chiaro che, sia pure in modo non identico, la costanza della velocità della luce c'era già nella teoria ondulatoria. E l'affermarlo chiaramente mi pare molto utile, poiché, anche volendo distruggere e rinnovare, è almeno necessario sapere cosa si distrugge; e non è onesto passare sopra allegramente (direi quasi, relativisticamente) a quanto è già stato trovato dalla scienza classica. Credimi coi più cordiali saluti

tuo aff.[mo] C. Somigliana

P.S. Ho letto con piacere la tua Nota nei Rend. dei Lincei,[238] e la tua conclusione circa l'opportunità di studiare newtonianamente i fenomeni di propagazione con sorgente mobile, è pure la mia. Né ritengo assoluta la mia asserzione sulla possibilità di una interpretazione newtoniana dei risultati della relatività. Essa aveva carattere soltanto generico.
Ti sarò grato se mi indicherai dove posso trovare, per disteso, quanto dice circa la deduzione relativistica del coefficiente di trascinamento di Fresnel.

Affettuosi saluti.

**201**
Gian Antonio Maggi a Carlo Somigliana

Pisa 29 Aprile 1923

Carissimo Somigliana,

Grazie della cara tua e dell'articolo in bozze (credo poterle trattenere), che ho letto con interesse, come ogni cosa tua. Convengo, in buona parte, coi tuoi giudizii, conformemente alla chiusa della mia lettera, che non presta alla teoria della relatività che una adesione relativa. L'Armellini, diventato apostolo della Relatività, ti potrebbe dire come gli tenetti simile discorso.

Mantengo però la mia opinione, per le ragioni della mia Noticina, che deve essersi incrociata alle tue bozze, sulla incomparabilità delle interpretazioni della trasformazione di Lorentz. E per quanto alla velocità di propagazione della luce non posso altrimenti che insistere sulla questione di fatto che non è mai stato contestato, ch'io sappia, che uno stesso valore appartiene alla velocità di propagazione assoluta della luce emessa da una sorgente fissa e da una sorgente animata da una certa velocità, e che la dipendenza, nel secondo caso, dalla velocità della sorgente appartiene alla velocità di propagazione relativa ad uno spazio che partecipa al movimento della sorgente, per modo che, rispetto ad esso, la sorgente è fissa. La qual dipendenza scaturita dalla suddetta ipotesi della costanza della veloc. di propagaz. assoluta, senza invocare punto l'assimilazione della velocità di propagazione a quella di un projettile, ma semplicemente in base alla distinzione intrinseca di velocità di propagazione, e alla più diretta intuizione, alla quale la teoria della relatività tira un frego, spazzando via il mezzo di propagazione, e postulando l'equivalenza di due spazii mobili, al noto modo, l'uno relativamente all'altro. Né la deduzione di questa dipendenza è da addebitarsi alla teoria einsteiniana, perché di molto anteriore al sorgere di questa, tanto che le prime esperienze di Michelson, a cui diede luogo, risalgono al 1881.

[238] Probabilmente si tratta di: “Sulle varie interpretazioni della trasformazione di Lorentz”, *Rendiconti della R. Accademia Nazionale dei Lincei*, v. XXXII, 1923, pp. 196−197.

Di nuovo, i migliori saluti e una cordiale stretta di mano dal

Tuo aff.[mo]
G.A. Maggi.

**202**

Gian Antonio Maggi a Carlo Somigliana
[carta intestata: R. Università di Milano - Facoltà di Scienze - Il Preside]

Milano 22 Febbraio 1927

Carissimo Somigliana,

Grazie della cara tua, e dei risultati che, non occorre dire, ho comunicato alla riunione del Seminario. La quale, per tuo principale merito, è riuscita molto bene, perché la mia esposizione fu seguita con vivo interesse dall'uditorio numeroso (con intervento di Palatini, Gerbaldi, Amerio di Pavia, e dei Dottori Giotti e Martin venuti da Merate), e diede per luogo ad un'interessante discussione, alla quale presero parte parecchi dei presenti, cominciando dal Bianchi, che, pel primo, mi manifestò il desiderio di quell'argomento. Concordemente lodati i pregi del tuo concetto. Per quanto al risultato $\frac{1}{2\varepsilon}$ per lo schiacciamento, questo è parso alquanto troppo superiore a quello ordinariamente accettato per non prestarsi a riserve, come discutibile è parsa la ragione dei coefficienti *[*adottati*]* nella formola pel calcolo di $g_l, g_q, g_p$. Il Bianchi mi ha domandato di mandargli le formole della tua lettera, per farvi nuovi calcoli, e non mancherò di soddisfare il suo desiderio, col trascrivergliele, perché non sarebbe il caso di mandargli anche il giudizio sul Segretario Accademico, che ha smarrito la traccia della tua Nota.

Ti faccio ora presente d'aver trovato, nella tua Nota lincea, da correggere un errore materiale di calcolo, che venialissmo per sé stesso, non permette però di ricondurre l'equivalenza delle tre equazioni (le due (7) e la (9)) a cui deve soddisfare la *[*i*]*. Le (6') vanno scritte (come esattamente è scritta la (6),

$$\left(\frac{g_3}{z_3}-\frac{g_1}{z_1}\right)\frac{z_1^2 z_3^2}{z_3^2-z_1^2}=\cdots\left(\frac{g_1}{z_1}-\frac{g_2}{z_2}\right)\frac{z_2^2 z_1^2}{z_1^2-z_2^2}=\cdots$$

Con questo, invece delle (7), si scriverà

$$\frac{g_2 z_3-g_3 z_2}{z_1}\,\frac{1}{z_1^2-z_3^2}=\frac{g_3 z_1-g_1 z_3}{z_2}\,\frac{1}{z_3^2-z_1^2}=\frac{g_1 z_2-g_2 z_1}{z_3}\,\frac{1}{z_1^2-z_2^2},$$

dalle quali, componendo, si ha

$$\frac{\dfrac{g_2 z_3-g_3 z_2}{z_1}+\dfrac{g_3 z_1-g_1 z_3}{z_2}+\dfrac{g_1 z_2-g_2 z_1}{z_2}}{0}=(z_1 z_2 z_3)^{-1}cB,$$

e per conseguenza la (9), che alla sua volta, permette di ricondurre, colla composizione, l'equivalenza delle due suddette equazioni.

**203**
Gian Antonio Maggi a Carlo Somigliana

Lettera a Somigliana[239]

Lanzo d'Intelvi, 18 Settembre 1927.

Carissimo Somigliana,

A completamento del noto nostro discorso di Como, ecco in breve la mia analisi, quale si trova nella mia Dinamica dei Sistemi alla quale (2[a] ediz.) si riferiscono le citazioni. Per maggior comodità di confronto sostituisco qua il sistema di punti al sistema di pezzi rigidi.

Indichiamo i vincoli del sistema composto di n punti P di coordinate x, y, z, con $\mu<\nu$ equazioni indipendenti

(1) $\sum_P \left( L_i \frac{dx}{dt} + M_i \frac{dy}{dt} + N_i \frac{dz}{dt} \right) = T_i$ (i= 1,2,…μ),

dove le $L_i$, $M_i$, $N_i$, $T_i$ indicano funzioni dalle coordinate di tutti i punti P e del tempo t (pag, 129, §51, equ. II e pag. 132 equ. III).

Ne seguono

(2) $\sum_P (L_i \delta x + M_i \delta y + N_i \delta z) = 0$ (i= 1,2,….μ)

per definizione dell'atto di movimento dal sistema, indicando con $\delta x, \delta y, \delta z$ le componenti della velocità del punto P (pag. 149, §56, equ. (2)) (se vuoi, per definizione dello spostamento virtuale intendendo che $\delta x, \delta y, \delta z$ rappresentino le componenti dello spostamento del punto P).

Enuncio poi i due postulati delle pressioni vincolari (pag. 160, §62), i quali semplicemente si traducono nelle equazioni

(3) 1° postulato: $m \frac{d^2x}{dt^2} = X + X'$; $m \frac{d^2y}{dt^2} = Y + Y'$, $m \frac{d^2z}{dt^2} = Z + Z'$,

dove X, Y, Z e X', Y', Z' rappresentano le componenti della forza motrice impressa, funzioni note delle coordinate di tutti i punti, dalle loro derivate rispetto a t, e di t, e dalla pressione vincolare applicata al punto generico P: e

2° postulato

(4) $\sum_P (X' \delta x + V' \delta y + Z' \delta z) = 0$.

Si aggiungono così alle 3ν nominate X, Y, Z, allora 3 ν che sono le X', Y', Z'.

Ora, si hanno, per formare un numero eguale di equazioni, le (3) che sono 3 ν, le (1) che sono μ, e 3ν–μ equazioni che si ricavano da (4) introducendosi pur $\delta x, \delta y, \delta z$, quello ch'io chiamo un sistema fondamentale di soluzioni della (2): cioè $3\,\delta x, \delta y, \delta z - \mu$ (numero dei gradi di libertà del sistema) soluzioni tra loro indipendenti, (dei cui singoli parametri sono poi una funzione lineare omogenea i parametri omologhi di ogni altra soluzione) (pag. 167-168, §62).

[239] Questa minuta si trovava nelle *Rusticationes*.

Io poi introduco le equazioni pure, dalle quali sono eliminate le X', Y', Z' (§63 e segu.), salvo ritornare e insistere; coll'ajuto della 1[a] forma delle equazioni di Lagrange, sulla loro determinazione in funzione di t, delle x, y, z, e delle $\frac{dx}{dt}, \frac{dy}{dt}, \frac{dz}{dt}$ di tutti i punti (§82).

Ma quanto precede non è, in sostanza, l'analisi del Poli? Al qual, per cortesia, ho parlato solo di analogia, di cui mi compiacevo. Pare invece che l'amico Levi-Civita gli abbia dato una tiratura d'orecchi. E a cotesta rispose con una lettera cortese, mentre della mia - che non poteva avere infine troppo diverso significato - si direbbe che non abbia trovato far di mejo che prenderne atto
...

(Nella lettera mandata, sostituito a x, y, z, la sola x, come una qualunque delle tre coordinate del punto P, analogamente X a X, Y, Z e X' a X', Y', Z').

**204**
Gian Antonio Maggi a Carlo Somigliana
[incompleta nella parte iniziale; s.l. e s. d.]

Lasciando da parte la teoria del Voigt, per discorrere in generale, la divergenza fra la teoria classica e la teoria della relatività fa capo alla distinzione sostanziale fra spazio fisso, a cui è connesso il mezzo di propagazione, e spazio mobile, che fa la prima, la qual distinzione è soppressa, insieme col mezzo, dalla teoria della relatività. Per la teoria classica sono eguali le velocità di propagazione della luce emessa da una sorgente fissa o da una sorgente mobile, intese velocità assolute, cioè rispetto allo spazio fisso: dipende invece dal senso e dalla direzione della velocità della sorgente la velocità di propagazione della luce relativa ad uno spazio che possiede la stessa velocità, pel modo che, rispetto ad esso, la sorgente è fissa. Invece, per la teoria della relatività, questa stessa velocità relativa, che non ha più ragione di essere distinta da una velocità di propagazione assoluta, è indipendente dalla velocità della sorgente. D'accordo che questa concezione ripugna all'intuizione. Ed è questa una delle ragioni per cui io non presto alla teoria della relatività che un'adesione relativa.
Ti lascio su questo punto, sul quale sono sicuro di avere il tuo consenso. Arrivederci al più tardi ai Lincei. Intanto una cordiale stretta di mano, colla quale mi confermo

Tuo aff[mo]
G.A. Maggi.

**205**
**Francesco Tavani** a Gian Antonio Maggi

92 Loughborough Rd
London S.W. 9
17.V.29

Ill$^{imo}$ Sig. Prof. Gian A. Maggi
della R. Università di Milano

Illustrissimo Signor Professore

Mi prendo la libertàdi rivolgermi a Lei per un esame e, se la crede degna, per presentazione della nota qui acclusa al R. Istituto Lombardo. Il Sig. Prof. Levi Civita ebbe la bontà di suggerirmi il suo nome. In attesa e domando scusa del disturbo ch'io possa darle.
Mi creda con la Massima Stima

Dev. Suo
F. Tavani

**206**
Gian Antonio Maggi a Francesco Tavani
[carta intestata: R. Università di Milano - Istituto Matematico
Via C. Saldini, 50 (Città degli Studi) - Il Direttore]

Al Sig. Francesco Tavani
92 Longborrough Rd. *[sic!]*
London S.W. 9.

Milano 7 Giugno 1929

Pregiatissimo Signore,

Il R. Istituto Lombardo non ammette comunicazioni di persone non appartenenti all'Accademia alla lettura e all'inserzione nei Rendiconti che in seguito ad approvazione della Sezione competente, che, in questo caso è la sezione di Scienze Matematiche, su conforme parere di un Membro Relatore. Mi spiace di non poter proporre l'approvazione del Suo lavoro, perché, consultato anche qualche collega della Sezione, non apparisce abbastanza chiaramente l'applicazione delle promesse trasformazioni agli accennati problemi, non rilevandosi così, come sarebbe desiderabile, lo scopo a cui sono destinate.
Non posso quindi far altro che rimandarLe il manoscritto, spiacente di non aver potuto soddisfare il Suo desiderio.
Aggradisca intanto i miei distinti saluti e mi creda

Dev.$^{mo}$ Suo
Gian Antonio Maggi.

**207**
**Paolo Ubaldi** a Gian Antonio Maggi
[busta intestata: Ill.[mo] Sig. Prof. G.A. Maggi della R. Univ.
Corso Plebisciti, 3 - Milano]

Ill.[mo] Sig. Professore,

Le scrivo, dietro consiglio dell'amico mio prof. Ettore Bignone dell'Univ. di Firenze, benché il mio nome Le sia del tutto sconosciuto. Fui incaricato dal prof. E. Codignola, pure di Firenze, di curare la edizione dell'opera postuma del mio venerato e indimenticabile maestro Giuseppe Fraccaroli, la quale consiste nella versione italiana della "Repubblica" di Platone. Il lavoro è quasi ultimato; ma c'è un punto che io non riesco a chiarire bene, perché non sono un matematico. Si tratta del passo platonico a pag. 546 B-D, che il Fraccaroli stesso dice (in una postilla) di non aver ancor bene dichiarato. E io so, dal Brignone, che il Fr. voleva discutere con Lei di questo passo. Vuole Ella esser tanto gentile di fare ciò che il mio Maestro non ebbe tempo di fare? Quello che soprattutto non capisco è il tratto che va da pag. 546 C a D. Qui ci vorrebbe una nota esplicativa ch'io non riesco a combinare! Tanto più che io temo che il testo sia non del tutto sicuro. Sono pronto a venire da Lei, qualora Ella lo credesse bene.

Le mando la trascrizione del passo tradotto dal Fraccaroli. E grazie fin d'ora di tutto quello che Ella farà in proposito.

Accolga i miei umili ossequi.

della S.V. Ill[ma]
Milano 23.XI.31
Via Copernico, 9

Dev[mo] D. Paolo Ubaldi
dell'Univ Cattolica

**208**
Gian Antonio Maggi a Paolo Ubaldi

Casa 4 Dicembre 1931.

Pregiatissimo Professore,

Le invio con questa la traduzione del noto passo, ch'io farei, ne' filosofo, ne' ellenista, traduzione organizzata come quella del Fraccaroli, sulla decifrazione del Jowett, dal quale però mi discosto in principio, coll'interpretazione di τρὶς αὐξηθείς con triplicato, che mi fornisce con $\frac{4}{3}\times 3 = 4$, il numero 4, e di πεμπάδι συζυγείς, con $4\times5^2$=100, moltiplicando il quale (ἑκατὸν τοσαυτάκις) per 100, ottengo più semplicemente il numero 10000 indicato da Jowett come prima ἁρμονία. Del τρίς αὐξησείς il Jowett non si vede bene quale uso faccia. Il Fraccaroli se ne scosta pure, interpretandolo come elevato della terza potenza, donde segue il risultato indicato nella Nota (3), notevole, ma non saprei come conciliabile col δύο ἁρμονίας παρέχεται. Credo poi indispensabile di accompagnare la traduzione con una Nota che ne spieghi il senso, come quella che io faccio seguire alla mia.

Sempre a Sua disposizione, conto, ad ogni modo, sul piacere di vederLa all'Inaugurazione degli Studii, l'8 corrente, alla quale sono cortesemente invitato. Intanto La prego di aggradire i miei più distinti saluti, e di credermi sempre

Suo Dev.$^{mo}$ Collega
G.A.M.

Traduzione (Platone Πολιτεία, 546)[240]

Di questi il fondamento epitrito, triplicato, combinato col numero cinque, fornisce due armonie, l'una [rappresentata da] un quadrato, cento per cento, l'altra, eguale, per una dimensione, a questa, ma rettangolare,[formata] da cento numeri [ricavati col quadrato] da diametri razionali del numero cinque, diminuito ciascuno di un'unità, [oppure] irrazionali [diminuito ciascuno] di due, e di cento cubi del numero 3.

Nota

Il fondamento epitrito, triplicato, combinato col numero cinque: $\left(\frac{4}{3}\times 3\right)\times 5^2 = 100$,

100×100, la prima armonia. Per la seconda armonia, va tenuto presente che il diametro del numero cinque, cioè l'ipotenusa del triangolo rettangolo, isoscele di cateto 5, è rappresentata da $\sqrt{50}$, e quindi, sotto questa forma, irrazionale, mentre 7= $\sqrt{49}$ ne fornisce la forma (approssimativa) razionale. Ciò posto la seconda armonia è data da

$$(7^2-1)\times 100+3^3\times 100=75000=100\times 75$$

oppure da

$$\left(\left(\sqrt{50}\right)^2-2\right)\times 100+3^3\times 100=75000=100\times 75.$$

Si ritrovano così i numeri che il Jowett, da cui è presa, per la seconda parte, questa decifrazione, propone come armonie

$$10000=160\times 100=10\times 25\times 40$$
$$7500=100\times 75=15\times 20\times 25,$$

Il primo, osserva il Jowett, riducibile alla semplice progressione aritmetica 2, 5, 8, il secondo riducibile ai lati del triangolo pitagorico 3, 4, 5.

*[Appunti]*

Πλάτων – Πολιτεία – η – 546.

Di questi il fondamento epitrito (ἐπίτριτος πυθμήν), combinato col numero cinque, fornisce, moltiplicato per tre (τρὶς αὐξηθείς), due armonie, l'una rappresentata da un quadrato, cento per cento, l'altra da un rettangolo, con un lato comune con quello, [formato] da cento numeri ricavati dai diametri razionali del numero cinque, diminuito ciascuno di uno, [oppure] irrazionali [diminuito ciascuno] di due, e di cento cubi del numero tre.

fondamento epirito $\equiv \frac{4}{3}$: moltiplicato per tre $\equiv \frac{4}{3}\times 3=4$: combinato con cinque $\equiv 4\times 5^2=100$: prima armonia $\equiv 100\times 100=10.000$: seconda armonia

$(7^2-1)\times 100+3^3\times 100=7500=100\times 75$ (4800+2700=7500)

$[(7{,}0204)^2=50]-2]\times 100+3^3\times 100=75000=100\times 75.$

Si ottengono così i numeri supposti da Jowett.

[240] *Repubblica*, VIII, 546c2-7. Si tratta di un passo noto per il suo carattere enigmatico che ha dato origine a molti tentativi di interpretazione. Ad esempio, si veda G. de Callataÿ, "Il numero geometrico", in: M. Vegetti (a cura di), *Platone. La Repubblica.Traduzione e commento. Vol. VI. Libro VIII-IX,* Napoli, 2005, pp. 169-187. Insieme a questo carteggio si trova un foglio con l'elenco delle opere consultate da Maggi per preparare la risposta.

**209**
Paolo Ubaldi a Gian Antonio Maggi
[busta intestata: Ill.mo Chiar.mo Sig. Prof. G.A. Maggi della R. Univ.
Corso Plebisciti, 3 - Milano]

Torino 6.XII.31
Via Cottolengo, 32

Chiar.mo e Carissimo Sig. Professore,

Ricevo qui a Torino, in questo momento, la sua gentilissima lettera, che mi toglie davvero da un grave imbroglio. Per quanto poco matematico io sia, mi pare che la sua interpretazione sia chiara e convincente! Vivissime grazie. Scrivo subito all'amico Bignone che il rebus (per noi due un vero rebus) è stato da Lei risolto.

Martedì sera non potrò, con mio vivo rincrescimento, trovarmi presente all'inaugurazione dell'anno accademico dell'Univ. Cattolica. Una grave disgrazia ha colpito me e tutti i miei soci. Ieri è spirato improvvisamente il Rettor Maggiore della Pia Società Salesiana, D. Filippo Rinaldi, della quale Società faccio parte fina dall'età di quattordici anni. Rimarrò a Torino per essere presente alle esequie che si faranno l'8 c.m. verso le 11. A Milano tornerò la settimana prossima: e allora Le telefonerò per sapere dove e quando avrò l'onore di fare la sua conoscenza personale.

Auguri a Lei e a tutta la sua famiglia di vita lunga e, per quanto è possibile in questo povero mondo, felice! Dio mi esaudisca.

Di Lei con profonda stima e sincero affetto

Devmo D. Paolo Ubaldi

**210**
**Gustavo Uzielli** a Gian Antonio Maggi

Firenze, 8 Giugno 1897
Viale Michelangelo 3 bis - Villa Nobili

Illustre Prof. G.A. Maggi,

L'Egregio collega Bertini mi ha fatto sapere che potevo scriverle direttamente circa la questione del Tunnel del Sempione. Ciò mi ha fatto vivo piacere.

Ho tardato a scriverle perché, oppresso da occupazioni, non potevo scriverle una breve lettera e perché infine speravo venirla a trovare a Pisa. La cosa infatti a me preme moltissimo tanto più che non può tardare la questione del Sempione a sorgere in Parlamento. Quindi, occorrendo, verrò costà, perché in dieci parole si chiarisce sovente un argomento più che con venti pagine.

Le sono in primo luogo gratissimo di avermi riassunto con parole chiare e semplici il modo di trattare il problema teorico in generale.

Per passare dal campo teorico al pratico credo doverLe sottomettere alcune notizie e osservazioni sull'argomento:

1° I calcoli eseguiti da vari fatti per conoscere la temperatura interna al centro della Galleria del Sempione concordano nel calcolarla di 40° centigradi e il capitolato obbliga l'impresa a ridurla a 25° cent. al massimo.

2° I mezzi proposti dall'impresa per ridurre in modo costante la temperatura da 40° a 25° sono:
L'introduzione dell'aria mediante un getto d'acqua con apparecchi simili ad un iniettore, e se occorre, anche impiegando dell'acqua polverizzata. Della relazione di uno dei periti Italo-Svizzeri risulta che questi "hanno potuto vedere nelle officine Sulzer in Winterthur un apparecchio nel quale l'aria riscaldata, trascinata da un altro tubo con un getto centrale d'acqua fredda a 12° e sotto la pressione di 5 atmosfera usciva molto raffreddata dal tubo: non poterono fare delle determinazioni precise della temperatura, ma l'effetto refrigerante d'acqua era evidente."
Il sottolineato è posto da me. A me pare che la mancanza d'esperienza rende illusorio il ragionamento che segue dopo quelle parole. Quelle esperienze erano la cosa fondamentale. Non ho possibilità di fare, ne' probabilità di vederle fare in modo scientifico. Quindi ho pensato che forse si poteva risolvere il problema, a priori con qualche approssimazione.
Ecco cosa aggiunge il Prof. Colombo: *[nota a piè di pagina:* (1) Il Politecnico, Anno XLIV, Aprile 1896, p. 227-228. Le do questa indicazione perché sarà bene che Ella veda tutto l'articolo del prof. Colombo*]* uno dei periti: "Basandosi sul potere assorbente dell'acqua, se si suppone che basti di raffreddare di 10°, in media, uno spessore di roccia di 3 metri ad una temperatura interna di 40° a 3 metri dalla parete della galleria, e ad una temperatura di 20° presso la stessa parete, se si calcola anche il calore che si deve togliere alla roccia staccata dalle mine ed il calore continuamente irradiante dalla parete della galleria attraverso lo spessore supposto di 3 metri, si trova che una quantità d'acqua di 50 litri per secondo basterebbe per abbassare a 20° la temperatura di ogni parte della galleria. Siccome nel progetto venne previsto una circolazione di 50 metri cubi d'aria al secondo con 84 litri d'acqua fredda dal lato nord e di 75 litri dal lato sud, l'impresa si crede perfettamente sicura di soddisfare alla condizione imposta dal contratto".
È esatto questo calcolo? Perché su questo punto fondamentale non si trova in nessuna relazione concernente al Sempione maggiori schiarimenti?
Come scrissi al Prof. Bertini al Gottardo non poterono abbassare la temperatura che di un grado. È vero che ciò dipendeva ancora da questo; che le macchine per la perforazione dovevano servire a procurare il raffreddamento e la ventilazione; ora erano obbligate a rimanere attive nel momento essenziale, cioè alquanto dopo lo scoppio delle mine che empivano la Galleria di gas nocivi. Perciò appunto, è indispensabile che le macchine per la ventilazione e il raffreddamento siano affatto separate da quelle per la perforazione, e perciò è indispensabile che vi sia al Sempione, come è progettato, una seconda galleria più piccola della principale, ove sarà disposto tutto il macchinismo igienico!
Ecco ora alcune osservazioni.
1°. Ho motivo di credere che l'esperienza del raffreddamento non sia stata fatta in condizioni, identiche alle naturali, cioè che la parete dell'ambiente raffreddato non siano state mantenute a temperatura costante.
2°. L'aria trascinata dall'acqua, arrivando sotto pressione produrrà un raffreddamento per la sua espansione; ma in un punto della Galleria alla distanza d dallo sbocco, l'acqua e l'aria inclusa si saranno riscaldate per aver camminato per la lunghezza d in un tubo metallico (certamente di ferro) e che avrà circa la temperatura della roccia.
3°. Non ho potuto trovare notizie esatte circa la disposizione per raffreddare la piccola Galleria, il che è certo cosa non insolubile, ma certo con aumento di spesa.
4°. Io ho identificato un poco inesattatamente il problema pratico col teorico, considerando in questo una sezione determinata della Galleria.

La temperatura di questa sezione è certo influenzata da quella dell'aria inclusa nelle due parti della galleria fra le quali è compresa la sezione considerata.
5°. Ella certo avrà notizia di vari lavori sulla conduttività delle roccie. Non ho qui i miei libri sull'argomento ma mi par certo che vi siano in proposito esperienze del Jannetaz.[241]

Quel che sarebbe molto importante, a vantaggio del problema che mi sono posto, si è di arrivare a trovare i limiti numerici, in gradi termometrici, fra i quali potrà trovarsi la temperatura della galleria, al punto dove sarà massima, cioè al centro ossia a km 10 (esattamente $\frac{19731}{2}$) dagli imbocchi.

Le unisco uno schizzo della sezione della Galleria e Le ricopio un frammento che già ho inviato al Prof.re Bertini, estratto da una relazione di un Ingegnere del Gottardo, e che mostra la terribile mortalità, ivi avvenuta, in causa della temperatura ecc.

Sarebbe certo un gran servizio d'umanità, se opportune previsioni scientifiche, potessero contribuire a diminuire la grande mortalità che avverrà nell'eseguire la Galleria del Sempione. Quanto ho letto in proposito è oscuro e scritto da chi è evidentemente interessato finanziariamente alla riuscita dell'impresa.
Questa, una volta concepita e senza stare a vedere se sarà più utile a Milano o a Genova (cosa ora in discussione) sarà certo vantaggiosa all'Italia.

Con la massima stima mi creda.

Suo Devotissimo
Gustavo Uzielli

P.S. Qualunque altro schiarimento desiderasse prego scrivermi e se gli fosse più comodo può anche farlo sapere a mio figlio Vittorio studente al primo anno di matematiche e latore del presente.

[241] Si veda la lettera #19, nota 55.

Grande Galleria per il transito

Piccola Galleria per il macchinismo di ventilazione e raffreddamento ecc.

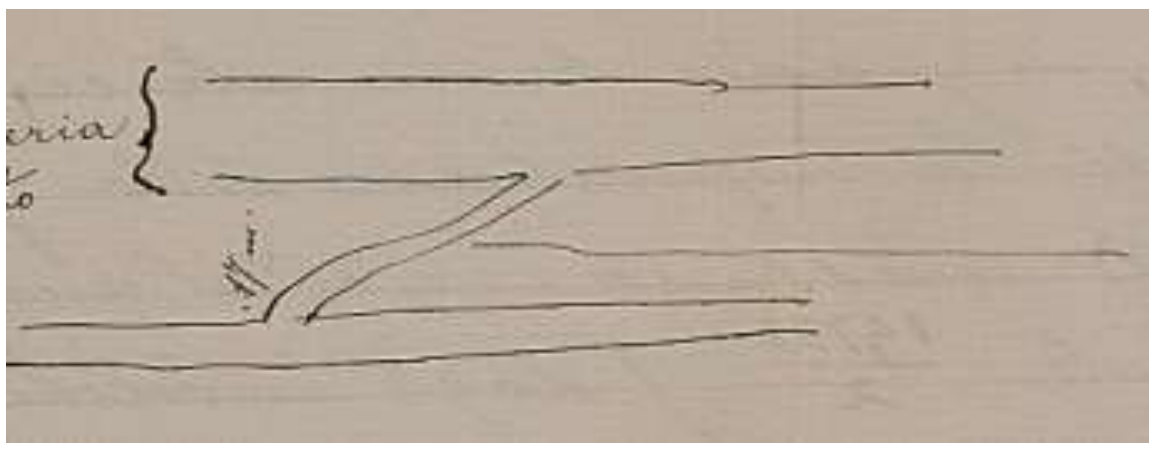

Queste piccole Gallerie di comunicazione fra le due principali si trovano a 200 metri

Sezione della gran Galleria

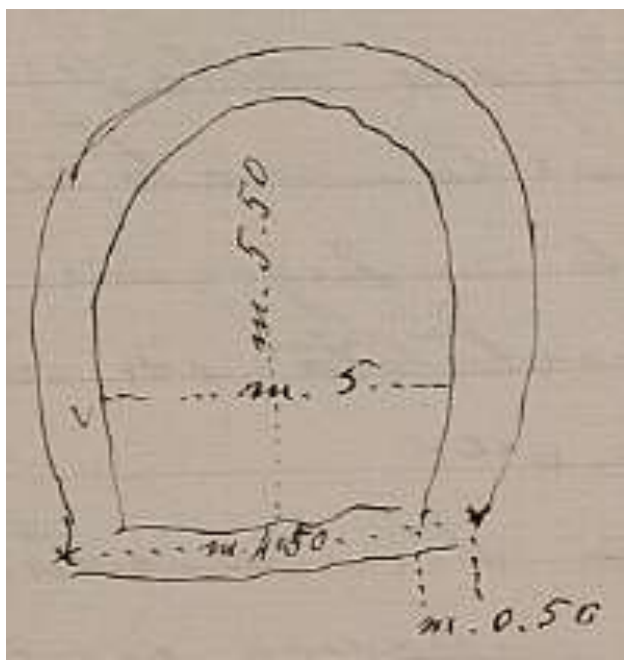

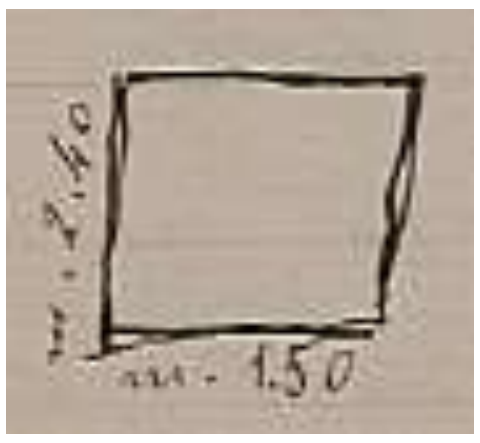

Sezione della piccola Galleria

Rivestimento in mattoni
Roccia del monte nella parte centrale della Galleria: Gneiss.

N.B. La forma del contorno della piccola Galleria non è perfettamente rettangolare, ma l'area della sezione rimane quasi invariata.

[Uzielli Lezioni 1887-88 p. 459]
Parole dell'Ing.re Mancy direttore durante 9 anni della Galleria Sud di Airolo]

"Questa temperatura (30°) non è tale da non permettere il lavoro all'esterno, ma in un ambiente saturo di umidità essa annienta completamente le forze umane. Gli uomini lavoravano del tutto nudi, e ciò nonostante erano incapaci di ottenere dei profittevoli risultati; il minimo movimento, la stessa parola erano una fatica ed il lavoro prodotto in queste condizioni era quasi nullo. Tutti gli operai divennero poco alla volta anemici e furono obbligati a lasciare il cantiere per essere sostituiti da altri che non producevano molto di più dei primi e mancavano della loro esperienza. Ci si farà un'idea di ciò che era la vita, nel tunnel l'ultimo inverno avanti la perforazione, quando dirò che il cuore arrivava a battere da 155 a 160 pulsazioni e la temperatura interna del corpo umano sorpassava i 39°."

**211**
Gustavo Uzielli a Eugenio Bertini
spedita poi a Gian Antonio Maggi
[busta indirizzata a: Ill$^{mo}$ Signor Cav. Prof. Eugenio Bertini
Università di Pisa; reindirizzata a: Firenze Viale Michelangiolo 3bis]

Firenze 16 Giugno 1897.

Egregio collega,

Sono lieto che il Prof. Maggi continuerà ad occuparsi della questione. Mi farebbe piacere sapere per quando potrò avere qualche risposta. Per i motivi che Le dissi non potrei tardar tanto (cioè non oltre a Luglio) a scrivere l'articolo che andrà pubblicato in due numeri del Giornale dei Lavori Pubblici, che si stampa qui in Firenze.

Si potrebbe far così. Del primo articolo, in cui accennerei al problema da me posto al Prof. Maggi, manderei le bozze a questo, ed egli modificherebbe quanto riguarda il detto annuncio, o anche, (per quanto certo non ne sarei lieto) toglierebbe l'indicazione del suo nome. Nel secondo articolo sarebbe contenuta, fra altre cose, la risposta del Prof. Maggi.

Prego quindi, per non aver la noia di rispondermi, di far avere, per ogni buon fine, a mio figlio Vittorio, l'indirizzo estivo del Prof. Maggi e l'accenno del tempo in cui esso lascerà Pisa.

Mi pregio accludere per la sottoscrizione Enrico Betti un vaglia di £20 spiacente di non poter mandare di più.

I miei più cordiali saluti

Devotissimo
Gustavo Uzielli

**212**
Gian Antonio Maggi a Gustavo Uzielli

Pisa 19 Giugno 97

Chiarissimo Signore

Le chiedo innanzi tutto scusa d'aver lasciato qualche giorno senza risposta la pregiatissima Sua dell'8 corrente, avendo avuto, nella scorsa settimana, il tempo completamente assorbito da altri e prevalenti impegni. E tanto più mi rincresce di non aver potuto occuparmene subito, poiché mi è bastatto rendermi conto della natura delle questioni proposte, per riconoscere che Ella mi fa l'onore d'attribuirmi una competenza, che mi manca affatto; per modo che, col miglior desiderio di servirLa, questa volta non lo posso fare in alcun modo. Le osservazioni contenute nella mia risposta al prof. Bertini riguardano, com'Ella vede, i principii più generali della teoria della propagazione del calore: movendo dai quali, ci s'inoltra, da una parte, nelle equazioni alle derivate parziali, e, dall'altra, nelle diverse applicazioni, che formano oggetto della tecnologia del calore. Questo è ciò che occorre pel Suo problema: ed io non ho tralasciato di dichiarare al collega Bertini che non era affar mio. Che se anche avessi maggiori cognizioni della questione dal lato fisico mi mancherebbero poi sempre, nel modo più completo, quei documenti della pratica, che certamente, in materia d'applicazione, sono gli argomenti che hanno veramente peso.

Perciò aggiungerei anche, avendomi il Bertini comunicato la Sua ultima lettera, che, restandole gratissimo pe’ suoi pensieri a mio riguardo, non mi parrebbe il caso né di nominarmi come autore di quelle osservazioni, ché non posso riconoscervi sufficiente importanza, né di mandarmi le bozze de’ suoi articoli, che volentieri leggerò stampati, essendo quel poco tutto quello che, a proposito, io posso dire.

M’è grata intanto l’occasione per inviarLe i miei rispetti, e professarmi

Dev.mo Suo
G.A. Maggi.

**213**
Gian Antonio Maggi a **Karl Theodor Vahlen**
[tramite Ulisse Dini]

Sul Teorema di Brioschi degli 8 quadrati[242]
Per K. Th. Vahlen.

Scopo dell’A. è di fornire un procedimento per dedurre i teoremi dei 4 e degli 8 quadrati analogo a quello che consiste nel fare nell’identità

$$ax \cdot a'x' = aa' \cdot xx'$$

per a e a’, x e x’ numeri complessi conjugati, per dedurre il teorema dei 2 quadrati, o nel porre nell’identità

$$\begin{vmatrix} a & -b \\ b' & a' \end{vmatrix} \begin{vmatrix} x & -y \\ y' & x' \end{vmatrix} = \begin{vmatrix} ax+by & ay'-bx' \\ b'x-a'y & b'y'+a'x' \end{vmatrix}$$

a=a’, b=b’, x=x’, y=y’ per dedurne il suddetto teorema, e per a e a’, b e b’, x e x’, y e y’ numeri complessi conjugati per dedurne il teorema dei 4 quadrati.

Egli si vale per ciò di un’identità fornita della moltiplicazione di due matrici di due righe e 4 colonne (pel teorema dei 4 e degli 8 quadrati) e di 8 colonne (pel teorema degli 8 quadrati). Questi teoremi si ottengono in tal modo prendendo gli elementi o a due a due eguali, o a due a due complessi conjugati.

L’A. rileva che il prof Studnicka, applica la ( ) all’ipotesi che a e a’, x e x’ rappresentino quaternioni conjugati, senza avvertire che la stessa identità implica la legge commutativa.

Noto che la stesssa objezione si può movere all’espressione della ( ) supponendosi le lettere simboli di quaternioni o di biquaternioni, dalla quale l’A. prende le mosse, pur restando il suo procedimento indipendente da essi

Al prof. Dini
2 Aprile 900.

[242] K.T. Vahlen, “Sul teorema di Brioschi degli 8 quadrati”, *Giornale di matematiche di Battaglini*, v. 39, Napoli, 1901, pp. 181-184.

**214**
**Quirino Valente** a Gian Antonio Maggi[243]

Chiarissimo Sig. Professore,

Il mio nome non Le sarà mai certamente nuovo, essendo stato un suo studente. Perciò mi prendo la libertà di sottoporre alla Sua saggezza il lavoro qui accluso, che vorrei esporre per tesi di laurea in Matematica in cotesta Regia Università.

Esso consta di due parti essenziali; nella 1° ho esposto ed ampliati i risultati ottenuti da Paul Suchar, trattando la *trasformazione reciproca di alcuni movimenti* (Bulletin de la Societé Mathématique, tom. XXXIII)[244]; nella 2° mi son proposto di vedere fino a che punto si possono spingere tali risultati. Ella vedrà adesso, se sono riuscito nell'intento.

Memore della Sua squisita cortesia, son certo che mi sarà largo di benigno compatimento, pel disturbo che Le procuro.

Se Le riuscisse d'incomodo a farmi sapere qui la Sua opinione, potrebbe comunicarla al Prof. Di Vestea.

Mentre Le chiedo tante scuse, ben distintamente La riverisco.[245]

Dev$^{mo}$
Roma, 15 gennaio 1907
Via Domenico Fontana, 12
Quirino Valente

**215**
Quirino Valente a Gian Antonio Maggi

Ill.$^{mo}$ Sig. Professore,

Fo seguito alla mia del 15 corr., prima di tutto per ringraziarla della cortese sollecitudine, con cui ha preso a cuore la faccenda che mi riguarda e poi per aggiungere le qui accluse osservazioni al capitolo che tratta la trasformazione dei movimenti nello spazio.

Con perfetta stima ed osservanza La riverisco.

Dev$^{mo}$
Quirino Valente

Roma, 22-1-907.

Osservazioni
da inserire nel Capitolo che tratta la trasformazione dei mov.$^{ti}$ nello spazio

---

Per dimostrare come dalle (22) si passi alle (16), basta prendere in considerazione il sistema lineare nelle l, m, n ..., che si ottiene paragonando i gruppi (I') (II') (III') con le equazioni che si hanno dai gruppi (I), (II), (III), moltiplicandoli opportunamente (come abbiamo già visto) per a, b, c … e il sistema delle (13), (14) e (15) ottenuto dagli stessi gruppi fondamentali ma paragonati in modo differente. Questi due sistemi, di cui la 1ª e l'ultima equazione sono le seguenti:

$$(w)\ \begin{cases} a'm + a''n = bl' + cl'' \\ cl'' + c'm'' = a''n + b''n' \end{cases} \qquad (w')\ \begin{cases} \Delta_{21}k_{12} + \Delta_{31}k_{13} = \Delta_{12}k_{21} + \Delta_{13}k_{31} \\ \Delta_{13}k_{31} + \Delta_{23}k_{32} = \Delta_{31}k_{13} + \Delta_{32}k_{23} \end{cases}$$

---

[243] Questo carteggio è raccolto in un fascicolo denominato *Corrispondenza con Valente*.
[244] "Sur une transformation réciproque en mécanique", 1905, pp. 210-224.
[245] Tale tesi verrà pubblicata: Q. Valente, *Intorno alla trasformazione reciproca dei movimenti*, Mareggiani, Bologna 1910.

Devono essere contemporaneamente soddisfatti e se si hanno le (22), le (w') possiamo scriverle

$$\begin{cases} a'\Delta_{21} + a''\Delta_{31} = b\Delta_{12} + c\Delta_{13} \\ c\Delta_{13} + c'\Delta_{23} = a''\Delta_{31} + b''\Delta_{32} \end{cases}$$

le quali, confrontate con le (w) ci danno precisamente le (16).

Infine giova osservare che le a, b, c …, le $K_{ir}$ e le l, m, n … sono dei numeri. Ora, in tutti i problemi di meccanica bisogna guardare a che l'omogeneità delle formule sia rispettata e nel nostro caso, un semplice sguardo alle formule (8) (9) di trasformazione, ci fa accorgere che le (26), le quali ci servono per l'interpretazione geometrica, possono avere un significato logico, solo quando i simboli che vi compariscono, siano assunti come misure, di quantità lineari e non come quantità vere e proprie.

## 216

Quirino Valente a Gian Antonio Maggi

Chiarissimo Professore,

Dopo una più matura riflessione, mi sono accorto che nel trattare la trasformazione nello spazio, sono caduto in un errore d'interpretazione.

La risoluzione del problema proposto, dipende, come ho detto, dalla possibilità di soddisfare contemporaneamente i gruppi fondamentali (I) (II) (III); (I') (II') (III'), oppure dalla possibilità di soddisfare il sistema delle (13), (14), (15), contenente le $\Delta_{ir}$, e le $K_{\rho s}$, insieme alle equazioni (α), (β), (γ), ottenute dalle (I), (II), (III).

Ora il confronto delle due equazioni:

$$\begin{cases} \Delta_{21}k_{12} + \Delta_{31}k_{13} = \Delta_{12}k_{21} + \Delta_{13}k_{31} \\ \Delta_{13}k_{31} + \Delta_{23}k_{32} = \Delta_{31}k_{13} + \Delta_{32}k_{23} \end{cases}$$

con la (12), ci dà la (16) qualunque sieno le $K_{\rho s}$ ($\rho \neq s$) e per conseguenza qualunque sieno le a, b… purché soddisf.[ti] a k=1. Le (16) e le (16') ci portano alle (21) e (21'), che risolvono completamente il problema, perché il nostro scopo è quello di esprimere le l, m, n… in funzione delle a, b, c… Il λ è arbitrario, perché è arbitrario ρ; ma ne vogliamo una trasformazione reale, bisogna escludere i sei valori di ρ che annullano il trinomio $\rho^4 + \frac{1}{\rho^2} - \rho$. Se fosse λ=0 con che ρ=1, ricadremmo nel caso del moto centrale, il che si verifica anche direttamente dalle formule (21) e (21').

La considerazione delle (22), le quali si ottengono confrontando le due equazioni:

$$\begin{cases} a'm + a''n = bl' + cl'' \\ cl'' + c'm'' = a''n + b''n' \end{cases}$$

con le (12), qualunque sieno le m, n, l', n', l'', m'', interviene qualora si vogliano esprimere le a, b, c..., soddisfacenti ai gruppi fondamentali, in funzione delle l', n', m, n... prese ad arbitrio.

Tanto le (16), quanto le (22) si possono ottenere indipendemente le une dalle altre, com'è facile verificare. Quindi non ho detto bene, affermando che le (22) si possono scrivere mercè le (16) e viceversa.

Dopo di che, l'interpretazione geometrica, cade.

Resta però ferma la conclusione e cioè che anche nel caso più generale di un movimento, la cui legge è espressa dalle formule (1), la trasformazione di Paul Suchar è estensibile allo spazio. Sostituendo nelle (8)-(9) alle l, m, n..., i loro valori che si ricavano dalle (21) (21'), si ottengono le formule che ci fanno passare dal movimento diretto al trasformato e viceversa.

Son certo che mi perdonerà il fastiodio che Le arreco, perché comprende che nelle questioni che presentano una certa difficoltà, non è difficile inoltrarsi in circoletti viziosi.

Rispettosamente La riverisco.

Dev$^{mo}$
Quirino Valente

Roma, 27 gennaio 1907.

**217**

Quirino Valente a Gian Antonio Maggi
[cartolina postale indirizzata: All'Ill$^{mo}$ Cav. G.A. Maggi
Professore di Meccanica Razionale nell'Università di Pisa]

Ill$^{mo}$ Sig. Professore,

Se mi è lecito di farLe un'altra preghiera, non esito a raccomandarLe di affrettare, per quanto Le è possibile, il Suo autorevole giudizio sull'accettabilità o meno della mia tesi, che, per non perdere il diritto alla mia regolare iscrizione in questa Scuola di Applicazione, dovrei discutere prima del 6 marzo p.v., giorno improrogabile che mi è stato concesso per la presentazione del Diploma di laurea.

Colgo l'occasione per informarLa che ho già riformata l'ultima parte della tesi nel senso della mia lettera ultima, dimostrando inoltre che le condizioni (16) per le l, m, n... sono necessarie e bastano per effettuare la trasformazione del Suchar nello spazio.

Con perfetta stima La riverisco

Dev$^{mo}$
Quirino Valente

Roma, 3 febbraio 1907.

**218**
Quirino Valente a Gian Antonio Maggi

Ill.[mo] Sig. Professore

Le accuso ricevuta del mio manoscritto, ringraziandoLa sentitamente delle indicazioni e dei consigli, che mi saranno certamente di utile indirizzo per qualche ricerca. Mi permetta però di dare alcune delucidazioni intorno alla 2.[a] parte del problema.

Ricordo di avere affermato che dalle relazioni

$$(1)\quad \begin{cases} l' = \Delta_{12} & m = \Delta_{21} & n' = \Delta_{32} \\ l'' = \Delta_{13} & n = \Delta_{31} & m'' = \Delta_{23} \end{cases}$$

si ricavano senz'altro:

$$(2)\quad \begin{cases} -\Delta_{11} = 1 + m' + n'' \\ -\Delta_{22} = 1 + l + n'' \\ -\Delta_{33} = 1 + l + m' \end{cases}$$

*[scritto a lato in matita:* <u>non in generale</u>
Es. l'=m=n'=l''=n=m''=0
l=4, m'=n''=$\frac{1}{2}$
$\Delta_{11} = \frac{1}{4}$, $\Delta_{22}$=$\Delta_{33}$=2 non valgono né la (2) né la (2*).*]*

E qui il Prof. Niccoletti osserva che questo non è affatto vero e che le (2) sono invece nuove condizioni da imporre. Ma basta mettere la (1) sotto la forma:

$l'(1+n'') = l''n'$

$l''(1+m') = l'm''$

$m(1+n'') = nm''$

$n(1+m') = mn'$

$n'(1+l) = nl'$

$m''(1+l) = ml''$

*[scritto a lato in matita:* Conviene guardare quando il processo è illegittimo. Questo è certo quando l'=l''=m=m''=n=n'=0.*]*

e dividere in croce per es. la 1[a] e la 5[a] di queste per avere $-\Delta_{11}$ e così di seguito. Pertanto la mia affermazione è esatta. Con le stesse (2) e ricordando che $\Delta$=1 abbiamo poi facilmente
(2*) $\Delta_{11}$=l; $\Delta_{22}$=m'; $\Delta_{33}$=n''
Un procedimento analogo, ci conduce dalle

$\rho l' = \Delta_{12}$

$\rho l'' = \Delta_{13}$ etc.

$\rho n' = \Delta_{32}$

alle altre

$$\Delta_{11} = \rho^2 + \rho l - \frac{1}{\rho}$$

$$\Delta_{22} = \rho^2 + \rho m' - \frac{1}{\rho}$$

$$\Delta_{33} = \rho^2 + \rho n'' - \frac{1}{\rho}$$

*[scritto a lato in matita:* Uguale os$^{ne}$ che disopra.*]*

per cui i risultati susseguenti, non sono privi di significato e ci conducono direttamente ad esprimere le l, m, n, l',… in funzione delle costanti a, b ,c,… che in fondo è il nostro scopo, prescindendo da qualunque via che si tiene per giungervi.

Ad ogni modo, appena corretto e ampliato il lavoro nel senso da Lei indicatomi per le ricerche nel piano, spero di fare una scappatina a Pisa per poterLe parlare a voce.

Con perfetta stima ed osservanza

La riverisco

Dev$^{mo}$

Q Valente

Roma, 9 febbraio 1907.

**219**

**Antonio Della Valle** a Gian Antonio Maggi

[busta indirizzata a: Ill.$^{mo}$ Prof. Gian Antonio Maggi della R. Università di Milano (Bergamo) Valnegra]

Napoli (Vomero)
Via Aniello Falcone, n. 212

18 Sett. 1930

Carissimo Amico.

La tua lettera mi ha fatto molto piacere e molto dolore!

Povero mio fratello! Come penosamente traeva il fiato con quella sua mancanza di elementi respiratorii per la scarsezza di globuli rossi. Quando l'anno passato si fermò per alcune settimane qui sul *[*tepido*]* Vomero, vi furono sì alcuni giorni mediocri; ma, purtroppo, avendo noi saputo della diagnosi del Prof Zoja, noi che lo conoscevamo sempre pieno di brio, lo consideravamo in penosa *[*malinconia*]*. Ahimè!

Gaetano (il primogenito del caro babbo) aveva avuto la premura di scrivere a mio nipote Mario (figlio dell'altro fratello nostro Romualdo morto pochi anni fa) intorno all'ultima giornata della penosa esistenza. Come era penosa quella lettera!

E difatti mio figlio Guido, che pure era partito da Milano immediatamente, non giunse a tempo a rivederlo in vita!

Tu mi ripeti che con Michele spesso parlavate anche di me. Quanto sei stato e quanto sempre sei ora, buono col vecchio amico tuo. Come sempre ricordo quei bei giorni di

Modena passati insieme,[246] e specialmente quella ricorrenza del nostro Santantonio, quando la tua affettuosa cortesia ti ti spinse a venire nel Laboratorio di Zoologia a propormi di passare insieme a casa tua la ricorrenza del nostro onomastico. Oh la lieta accoglienza della gentilissima tua Signora!

L'ho ricordata più volte alle mie figliuole!

Siamo rimasti soli, carissimo Maggi.

Che la Divina Provvidenza ci protegga!

Ti abbraccio col più vario affetto

AD Valle

**220**

Gian Antonio Maggi a **Alessandro Visconti**[247]

Lanz D'Intelvi, 17 Settember 1935

Car el mè sur Alex,

Minga, Dio ne libra! per pedantaria, ma per amor, s'el voeur per gelosia, de la puritaa del noster meneghin, me permetti de osservà che i milanes del temp de Maria Teresa, ben pussee fedel al meneghin che i noster present concittadin desambrosianaa, avarissen ditt "sur Giusepp", e minga "sciur Giusepp"; come poeu el Porta "Ma saal el mè sur Lella…" "Come fala mo' a dì sure Lenin..." ecc; e come mi el ciami lu "sur Alex" e minga "sciur", augurandagh che "scior" el poda ciammass, cont el pusse giust significaa de la parola.

Credi minga da sbagliamm, in quanto al milanes. Ma forsi me sbagli cont el fagh perd el tempo con sti ciaccer. Allora, car el mè sur Alex, ch'el me scusa, ma ch'el veda però come me lassi minga scàpà i sò bei articol, e come ghe metta attenzion, magari anca troppo!

A rivedess prest a Milan! Intant una bonna strengiuda de man del

sò

Gianantoni Maggi.

Ch.mo Prof. Alessandro Visconti
della R. Università di Ferrara
Via Montenero 66
Milano.

[246] Nel 1885 Maggi era professore di Analisi infinitesimale presso l'Università di Modena, mentre Della Valle vi insegnava Anatomia comparata.

[247] In due copie di minuta che si trovavano nelle *Rusticationes*.

**221**

[Giulio Vivanti] a Gian Antonio Maggi[248]

*[scritto da Maggi:* Portatomi da Vivanti il 18 Aprile 1935*]*

Caso semplice

Se $(\xi,y)$ è un estremo della funzione $f(x,y)$, e se è $\varphi(\xi,y)=0$, $(\xi,y)$ è anche un estremo condizionato di $f(x,y)$ colla condizione $\varphi(x,y)=0$.

Per le ipotesi fatte $(\xi,y)$ soddisfa alle condizioni:

(1) $\frac{\delta f}{\delta x}=0, \quad \frac{\delta f}{\delta y}=0, \quad \varphi=0.$

D'altra parte le condizioni per un estremo di f sotto la condizione $\varphi=0$ sono:

(2) $\begin{vmatrix} \frac{\delta f}{\delta x} & \frac{\delta f}{\delta y} \\ \frac{\delta \varphi}{\delta x} & \frac{\delta \varphi}{\delta y} \end{vmatrix}=0, \quad \varphi=0;$

ora le (2) sono evidentemente conseguenze delle (1).

Caso generale

Siano date le funzioni f, $\varphi_1, \varphi_2, \ldots, \varphi_m$, delle variabili $x_1, x_2, \ldots, x_n$, essendo $m<n$. Se, posto $r<m$, $(\xi_1, \xi_2, ..., \xi_n)$ è un estremo di f sotto le condizioni $\varphi_1, \ldots, \varphi_r=0$, e se soddisfa inoltre alle condizioni:

(1) $\varphi_{r+1}=0, \ldots, \varphi_m=0$,

allora $(\xi_1, \xi_2, ..., \xi_n)$ è un estremo di f sotto le condizioni

(2) $\varphi_1=0, ..., \varphi_m=0$.

Per ipotesi $(\xi_1, \xi_2, ..., \xi_n)$ soddisfa alle condizioni:

(3) $\begin{vmatrix} \frac{\delta f}{\delta x_1} \cdots \frac{\delta f}{\delta x_n} \\ \frac{\delta \varphi_1}{\delta x_1} \cdots \frac{\delta \varphi_1}{\delta x_n} \\ \cdots\cdots\cdots \\ \frac{\delta \varphi_r}{\delta x_1} \cdots \frac{\delta \varphi_r}{\delta x_n} \end{vmatrix}=0,$

[Intendo con questa notazione che sono nulli tutti i minori d'ordine massimo della matrice]

(4) $\varphi_1=0, \ldots, \varphi_r=0$.

D'altra parte un estremo di f sotto le condizioni (2) deve soddisfare alle equazioni:

[248] Si veda la nota alla lettera #90.

$$(5)\quad \begin{vmatrix} \frac{\delta f}{\delta x_1} \cdots \frac{\delta f}{\delta x_n} \\ \frac{\delta \varphi_1}{\delta x_1} \cdots \frac{\delta \varphi_1}{\delta x_n} \\ \cdots\cdots\cdots \\ \frac{\delta \varphi_m}{\delta x_1} \cdots \frac{\delta \varphi_m}{\delta x_n} \end{vmatrix},$$

(6) $\varphi_1=0, \ldots, \varphi_m=0$;
ora le (6) coincidono colle (1), (4) e la (5) è evidentemente conseguenza della (3).

Geometricamente la cosa è evidente.
Prendo il caso semplice. Se ($\xi$, y) è un massimo sul piano xy per la funzione f(x,y), e se per quel punto passa la curva $\varphi(x,y)=0$, il punto è un massimo sulla curva.

Analogamente il caso generale.

*[aggiunto da Maggi:*

Questa dimostrazione si presta alla seguente forma analitica: valori delle variabili che estremizzano f incondizionatamente, soddisfanno df=0 confrontati con valori qualisivogliano delle variabili medesime, e per conseguenza anche confrontati con valori che soddisfanno particolari condizioni, che è l'ipotesi della estremizzazione condizionata.*]*

**222**

**Vito Volterra** a Gian Antonio Maggi

Via in Lucina, 17
Roma 27 Novembre 1907

Carissimo Amico

Dovendo riguardare per un concorso una Nota del Dott. Chella[249] pubblicata nel Luglio 1905 nei Rendiconti dell'Accademia dei Lincei io vi ho trovato varie difficoltà.

Il Battelli mi disse che questa Nota è un estratto della tesi. Mi faresti cosa grata se tu potessi dirmi qualche cosa intorno ad essa se tu avessi avuto già occasione di leggerla.
Scusami del disturbo e gradisci i più affettuosi saluti

Dal tuo aff$^{mo}$ amico
Vito Volterra

[249] S. Chella, "Su di una misura assoluta del coefficiente di attrito interno dei gas", *Rend. della R. Accad. dei Lincei*, v. XIV, 1905, pp. 23-30. Nota presentata dal Corrispondente A. Battelli. Lavoro eseguito nell'Istituto di Fisica della R. Università di Pisa.

**223**
Gian Antonio Maggi a Vito Volterra

Pisa 2 Dicembre 1907

Carissimo amico,

Ho tardato un pajo di giorni a risponderti, perché mi è occorso di cercare gli scritti del Meyer, e volevo anche rivedere la tesi di laurea del D.r Chella, che non fu trovata nell'Archivio dell'Università. Oggi poi è venuto da me, mandato dal Battelli, lo stesso D.r Chella, e ho pur creduto di dirgli che mi ero occupato del suo lavoro, poiché tu incidentalmente, in una tua recente lettera, mi avrai scritto di averci riscontrato delle "difficoltà". Tu non mi dicevi quali fossero: ma ecco quelle che ho incontrato io, e che ho fatto presenti al Chella.

1) Un $=\Delta_2\psi$, che non sta.

2) Un'equazione che sarebbe, se stesse la precedente eguaglianza, quella delle funzioni cilindriche, ed è affermato, ma non so se possa essere dimostrato che s'integra per mezzo di queste.

Queste due affermazioni non danno però luogo ad alcuna conseguenza.

3) Una sommatoria introdotta in principio, per poi ridursi tacitamente a considerare un solo termine.

Qui è questione d'intendersi meglio fin dal principio.

4) Ipotesi relative al contorno diverse al contorno diverse dalle generali dell'idrodinamica, e, in particolare, da quelle di Meyer, mentre poi il momento delle pressioni d'attrito è calcolato con la formola di Meyer.

Non ho cercato altro ma quest'ultima objezione è capitale. Vi ho specialmente insistito col Chella, il quale mi ha detto che proverebbe di render ragione del suo procedimento. Purché ci riesca!

La tesi di laurea comprendeva anche una parte sperimentale, sulla quale riferì particolarmente il Battelli e la Commissione si attenne alla sua relazione. Mi spiace che, non essendomi occupato dell'argomento dell'attrito interno dei fluidi in modo speciale, non feci, a suo tempo, sufficiente attenzione alla parte teorica, per rilevarvi gli scogli in discorso.

Mi farà molto piacere di conoscere il tuo giudizio. Intanto t'invio, da parte nostra, i migliori saluti, e ti prego di presentare i miei ossequi alla tua Signora,[250] e a tua Madre,[251] che spero sarà a quest'ora completamente ristabilita.

In particolare, aggradisci una cordiale stretta di mano dal

Tuo affmo amico
G.A. Maggi

P.S. Sono rimasto esterrefatto e addoloratissimo per la morte del povero Sella![252]

Al Ch.mo Sig. Prof. Vito Volterra
Senatore del Regno
Via Lucina 17 Roma.

---

[250] Almagià Virginia.

[251] Almagià Angelica.

[252] Si tratta di Alfonso, figlio di Quintino, che morì il 25 novembre. Egli rivestiva un ruolo chiave nella Società Italiana per il progresso delle Scienze, fondata da Volterra l'anno precedente.

**224**
Gian Antonio Maggi a Vito Volterra

Pisa 4 Dicembre 1907

Carissimo Volterra,

Aggiungo alla mia risposta, che avrai ricevuto jeri, che, essendomi inoltrato nella lettura degli articoli del Meyer, ho riconosciuto che questi finisce poi per fare ipotesi simili a quelle del Chella il quale perciò né s'ingannava, né traeva altri in inganno, affermando di essersi attenuto al Meyer. Resta a vedersi come conciliarle con quelle che pone a fondamento della teoria, ed egli si accorge bene di non procedere a perfetto rigore, tanto che dice (Crelle, Vol. 59, pag. 247[253]) "Um zu vermeiden, dass durch die eingeführte Annäherung das mathematische Interesse an der Aufgabe geschemälert werde etc" E cita l'esperienza… che però non salvò mai una incongruenza logica! Per esempio, è supposta la condizione alla faccia del disco $\eta = \frac{d\Psi}{dx} = E(\Psi - \Psi_1)$ e poi e poi *[sic!]* fatto $\Psi = \Psi_1$, conservando a $\frac{d\Psi}{dx}$ un valore diverso da 0. Io trovo da comporre la contraddizione (non so se il Meyer in qualche luogo ragioni allo stesso modo), osservando che $\Psi = \Psi_1$ vuol dire mancanza di scorrimento fra disco e fluido, quindi attrito esterno estremamente grande, allora sensibilmente $E = \infty$.

Di nuovo i miglior i saluti e una cordiale stretta di mano dal

Tuo aff. amico
G.A. Maggi

Al Ch. Sig. Prof. V. Volterrra

[253] O.E. Meyer, "Ueber die Reibung der Flüssigkeiten", *Journal für die reine und angewandte Mathematik*, v. 59, 1861, pp. 229-303.

## 225
Gian Antonio Maggi a **Francesco Zambaldi**[254]

Gentilissimo Professore,

Il discorso di Platone mi sembra essere il seguente: Data un'area sotto una forma qualsivoglia, cioè rappresentata da una figura (piana) limitata da un contorno qualsivoglia, distendiamola entro un dato cerchio, conformente al dato contorno παρὰ τὴν δοθεῖσαν αὐτοῦ γραμμὴν; il risultato, rispetto all'iscrizione di un triangolo la cui area sia la supposta, nel cerchio, sarà diversa secondo che la parte mancante del cerchio (cioè non coperta) sarà eguale, o no, alla parte coperta.

L'affermazione è indiscutibile, ma, se altro non aggiungo, poco significativa. Non so se si possa conferire significato col seguente discorso. Teniamo presente che l'area del massimo triangolo che può essere iscritto in un dato cerchio risulta press'a poco la metà dell'area del cerchio. Allora se accettiamo questo rapporto come esatto, e intendiamo che il triangolo debba essere regolare, alla quale ipotesi compete appunto la suddetta area massima, quest'area dovrà essere la metà dell'area del cerchio, poiché il triangolo sia iscrivibile nel cerchio. E per conseguenza, coperta una parte del cerchio colla figura data, avente l'area supposta, il triangolo [regolare] avente la stessa area. sarà iscrivibile nel cerchio, oppur no, secondo che la parte mancante sarà, o non sarà, eguale alla parte coperta.

La tacita ipotesi della regolarità del triangolo non sarebbe disforme da simili limitazioni che si trovano, per esempio, in Archimede.

Tolta la limitazione della regolarità, e ammesso sempre il suddetto rapporto, il criterio per l'iscrivibilità di un triangolo avente la supposta area del cerchio, diventa che, distesa entro il cerchio la figura data, la parte mancante (non coperta) sia eguale almeno alla parte coperta. Ma questo almeno non si trova nel testo.

Con questo, salvo consegnarLe il foglietto brevi manu, e sempre a Sua disposizione per ogni schiarimento, La prego di aggradire i miei più cordiali saluti e di credermi sempre

Suo aff.mo collega
Gian Antonio Maggi.

Casa, 16 Marzo 1916

Al prof. Zambaldi

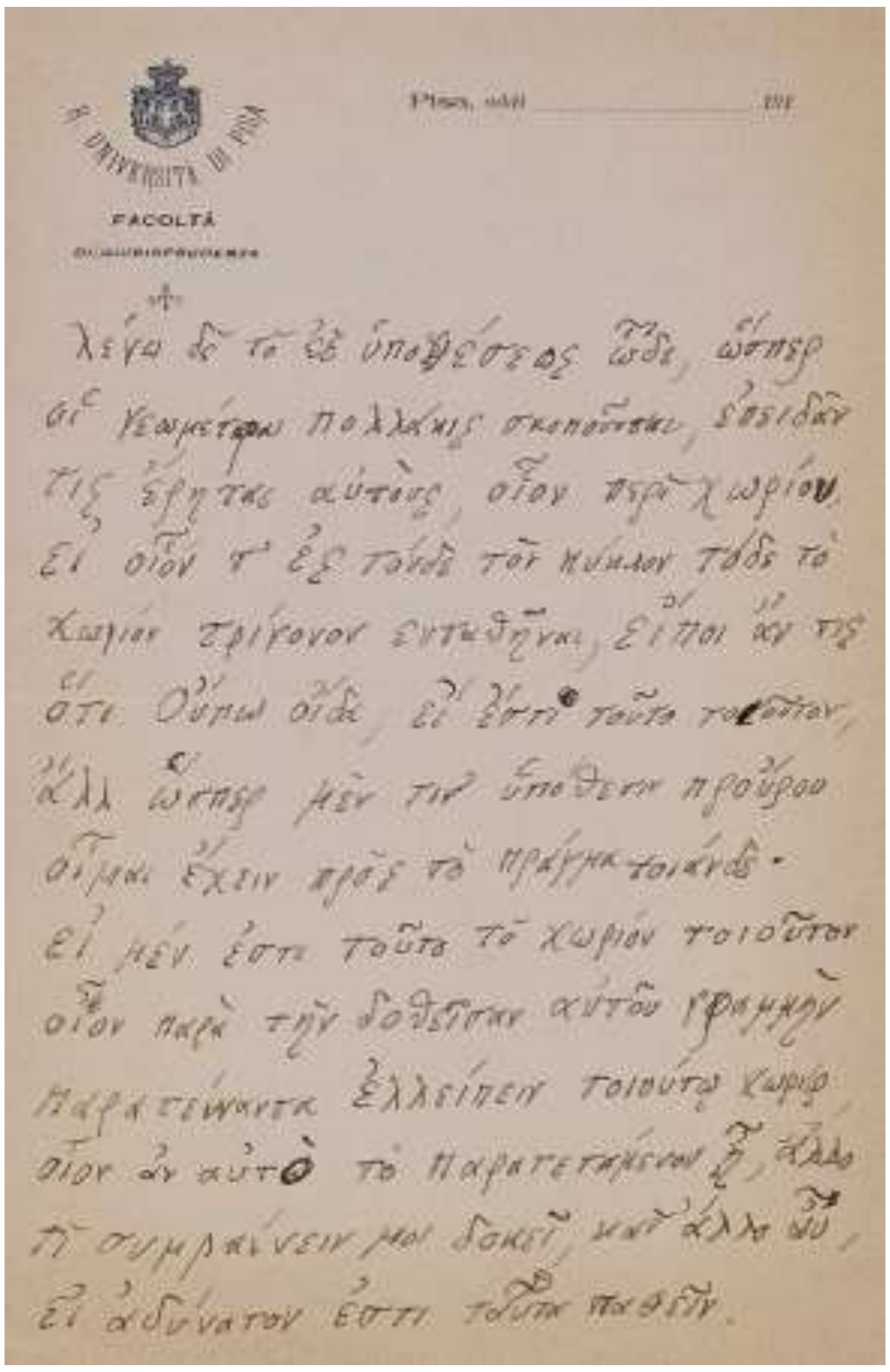
R. Università di Pisa
Facoltà di Giurisprudenza
Pisa, addì ... 191

λέγω δὲ τὸ ἐξ ὑποθέσεως ὧδε, ὥσπερ
οἱ γεωμέτραι πολλάκις σκοποῦνται, ἐπειδάν
τις ἔρηται αὐτούς, οἷον περὶ χωρίου,
εἰ οἷόν τ' ἐς τόνδε τὸν κύκλον τόδε τὸ
χωρίον τρίγωνον ἐνταθῆναι, εἴποι ἄν τις
ὅτι Οὔπω οἶδα εἰ ἔστι τοῦτο τοιοῦτον,
ἀλλ' ὥσπερ μέν τιν ὑπόθεσιν προὔργου
οἶμαι ἔχειν πρὸς τὸ πρᾶγμα τοιάνδε·
εἰ μέν ἐστι τοῦτο τὸ χωρίον τοιοῦτον
οἷον παρὰ τὴν δοθεῖσαν αὐτοῦ γραμμὴν
παρατείναντα ἐλλείπειν τοιούτῳ χωρίῳ
οἷον ἂν αὐτὸ τὸ παρατεταμένον ᾖ, ἄλλο
τι συμβαίνειν μοι δοκεῖ, καὶ ἄλλο αὖ,
εἰ ἀδύνατόν ἐστι ταῦτα παθεῖν.

[254] Questa lettera e gli scritti preparatori (tra cui la trascrizione in greco del passo e la sua traduzione in latino) sono contenuti nel fascicolo: *Dal Menone di Platone. Passo dubbio dell'art. 22, sul quale mi consulta il collega prof. Zambaldi*. Gli appunti in greco, dei quali è inserita la scansione, vanno dalla riga 86e4 alla 87a8, secondo la numerazione basata sull'*editio princeps*. Sul retro Maggi scrive: *possitne in hoc circulo triquetrum hoc spatium extendi, doxerit aliquid, nondum novi, an hoc hujusmodi sit, sed ad hanc rem hypothesin arbitror talem conducere: si enim est hoc spatium tale, ut quod ad positam ejus lineam protendis ab eo tanto spatio deficiat, quantum ipsum est protensum, aliud quiddam evento puto... (Ediz. Didot).*

*[Frammento]*
Gentilissimo Professore,

Non so se la seguente possa essere accettata come un'interpretazione di quel passo di Menone, che mi sembra veramente di colore oscuro.

Sta il fatto che l'area del massimo triangolo, il triangolo equilatero, che può essere iscritto in un cerchio dato, sta all'area del cerchio nel rapporto di $\frac{3\sqrt{3}}{4}$ a π, che è, all'ingrosso, $\frac{1}{2}$*.

Allora, se ammettiamo questo rapporto, e supponiamo, per un momento che il χωρίον τρίγωνον debba intendersi regolare, cioè il triangolo equilatero, il discorso di Platone potrebbe essere codesto: Distendiamo nel cerchio la supposta area conformemente al contorno dato (παρὰ τὴν δοθεῖσαν αὐτοῦ γραμμὴν), la stessa area in forma di triangolo [regolare] sarà iscrivibile nel cerchio dato, o no, secondo che la parte del cerchio che non risulta coperta è eguale, oppure no, alla parte coperta.

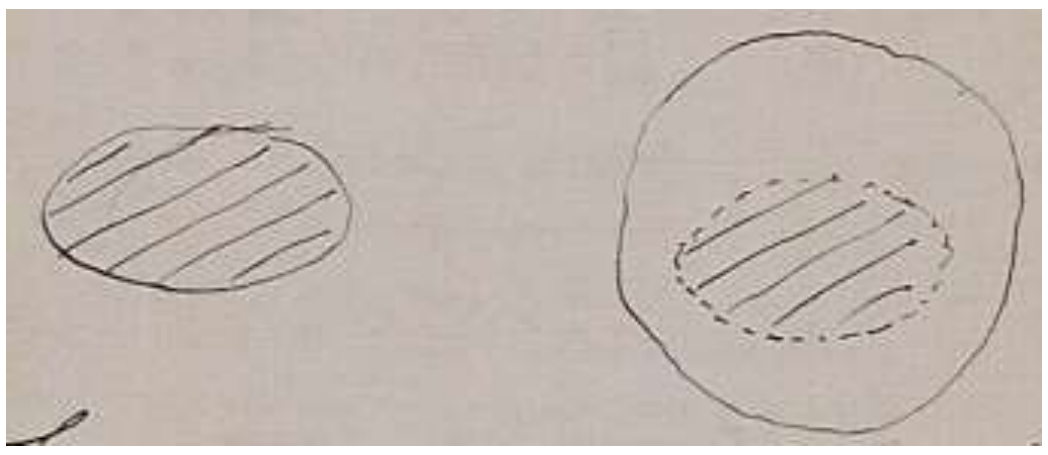

La tacita ipotesi della regolarità del triangolo non sarebbe difforme da simili limitazioni che si trovano implicite, per esempio, in Archimede.

A voler prescindere da quella dimostrazione, il discorso dovrebbe essere il seguente: l'area, in forma di triangolo, sarà iscrivibile nel cerchio dato, o no, secondo che la parte del cerchio che non risulta coperta non è minore, o è minore, della parte coperta

* 0,4136

## 226
## [Giuseppe] Zuppone-Strani a Gian Antonio Maggi[255]

Chiarissimo Signor Professore,

Ecco quell'innocentino d'un problema. Per tranquillità della mia conoscenza la prego di mostrarmi la via per risolverlo o di dichiararmelo di impossibile geometrica soluzione.
Scusi e gradisca molti ossequi
Messina 19/1/93

Suo dev[mo]
G. Zuppone-Strani

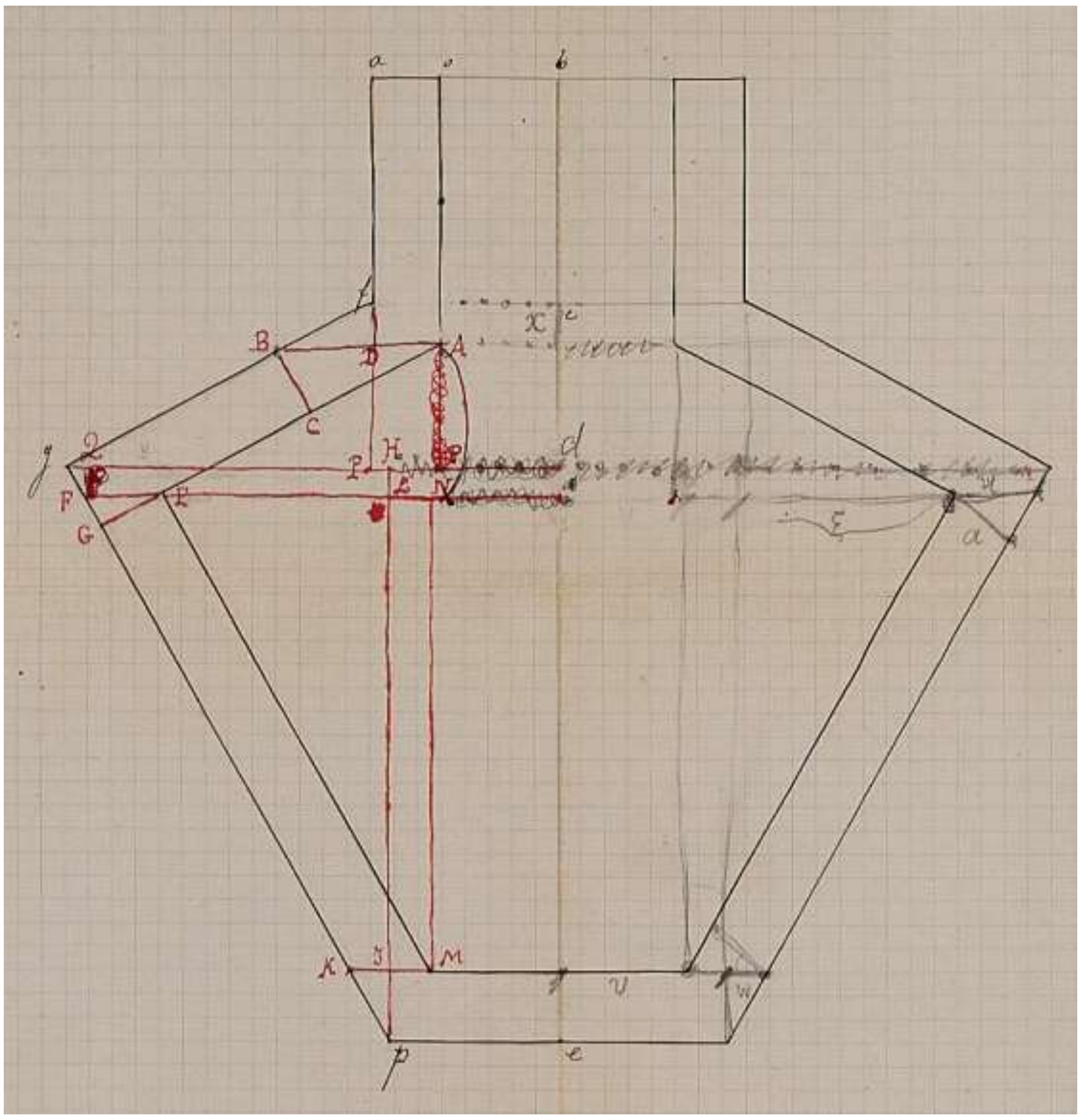

Conosco <u>tutti</u> gli elementi per trovare il volume totale dell'anfora. Conosco oltre a ciò lo spessore che è da pertutto costante. Chiedo la capacità senza fare uso di trigonometria.

Conosco ab=r<br>
gd=R,<br>
pe=r<br>
bc=a<br>
cd=h<br>
de=h'<br>
as=s

[255] Si veda la nota alla lettera #13.

*Appunti di Maggi per la soluzione del problema proposto:*

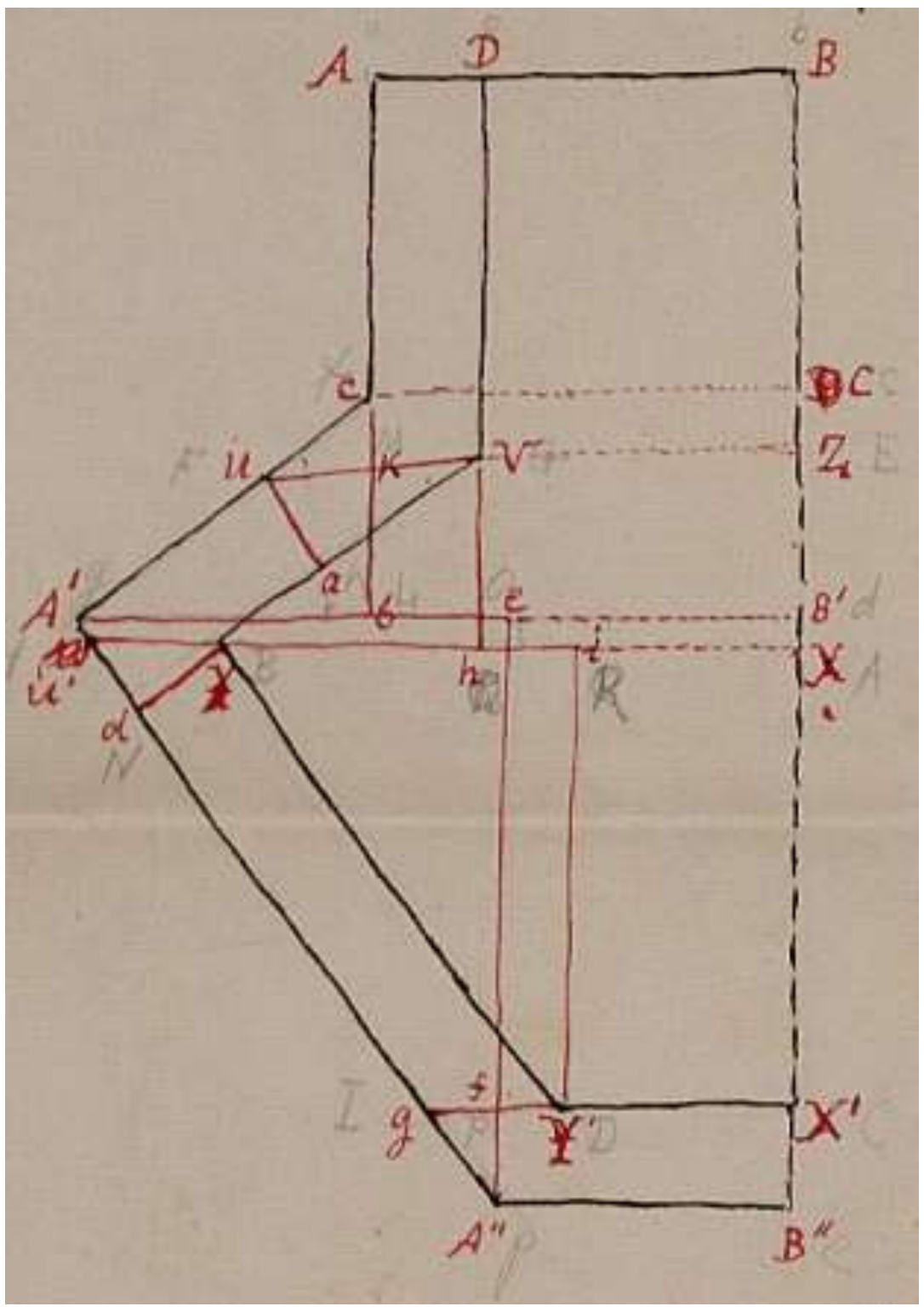

Dati AB=r, A'B'=r', A''B''=r'', BC=a, CB'=a', B'B''=a'', AD=S.
Incognite principali: BD=x, XY=x', X'Y'=x'', CZ=y, B'X'=y'.
Incognite ausiliare: UV=u, U'Y=gY'=u'.

1°. Si ha BD=AB–AD. Quindi x=r–S.

2°. Da UVa e A'bc: $\frac{u}{S}=\frac{\sqrt{a'^2+(r'-r)^2}}{a'}$. E si ha u.

3°. Da UVa e ckU: $\frac{y}{u-S}=\frac{S}{\sqrt{u^2-S^2}}$. E si ha y.

4°. Da u'Yd e A'A''e : $\frac{u'}{S}=\frac{\sqrt{a''^2+(r'-r'')^2}}{a''}$.E si ha u'

5°. Da A''fg=u'Yd: gf=U'd= $\sqrt{u'^2-S^2}$ .Quindi: x''=gf+fX'–gY'=$\sqrt{u'^2-S^2}$ +r''–u'.

6°. Le due coppie di triangoli A'bc e VhY, A'A''e e YY'i forniscono ciascuna una relazione lineare fra x' e y'.

Dalla prima: $\frac{a'-y+y'}{x'-x}=\frac{a''}{r'-r}$. $\qquad \frac{Vh}{Yh}=\frac{bc}{A'b}$

Dalla seconda: $\frac{a''-y'-S}{x'-x''}=\frac{a''}{r'-r''}$. $\qquad \frac{Y'i}{Yi}=\frac{A''e}{A'e}$

Posto $\frac{a''}{r'-r}=q$, $\frac{a''}{r'-r''}=q'$, le due relazioni si possono scrivere:

y'–qx'=y–a'–qx, $\qquad$ y'+q'x'=a''–S+q'x''. E di qui

$x'=\frac{a''-S-y+a'+ax+q'x''}{q+q'}$ $\qquad$ $y'=\frac{(a''-S)q-(y-a')q'-qq'(x''-x)}{q+q'}$.

**227**
Gian Antonio Maggi a destinatario sconosciuto[256]
[biglietto]

Milano 13 Febbrajo 1930

Carissimo amico,

Ho ricevuto jeri l'altro dallo Zanichelli il vostro bellissimo volume, che ho finito ora di tagliare, procurandomi il matematico piacere, di scorrerlo per intero, e te ne rendo vivissimi ringraziamenti, colla preghiera di farne parte al collega Amaldi.

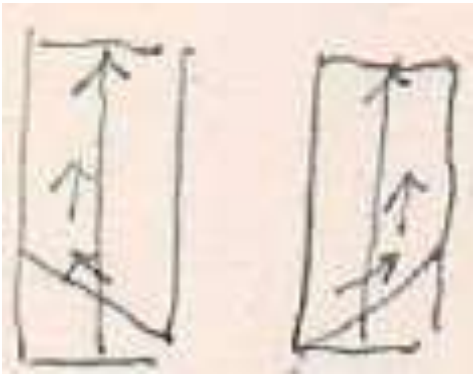

Ho rilevato lo scrupolo etimologico del destrorso e sinistrorso. Comunque sia, vite destrorsa e sinistrorsa sono termini di uso comune, che mi sembrano giustificabili dal salire dell'elica a destra e a sinistra di un osservatore posto davanti alla vite, ipotesi più alla mano di quella dell'osservatore posto nell'asse della vite, riservata agli iniziati. E con questo crederei che si possa intendere destrorso e sinistrorso anche pel sistema cartesiano. Questione, del resto, di parole. Non ti sarà sfuggita nel Maxwell (§23, Nota) la proposta di chiamare i due sistemi "of the vine and the hop" che torna sistema del vino e della birra.[257]

I miei più cordiali ossequi etc.

**228**
Gian Antonio Maggi a destinatario sconosciuto[258]
[s.l.e s.d.; "Non mandata"]

...
Non colgo bene la differenza, fondata sui termini infinitesimali di ordine superiore, che fai tra i due casi relativi a

(1) $$dq = C_p \frac{\partial T}{\partial v} + C_v \frac{\partial T}{\partial p} dp:$$

e cioè il caso che il secondo membro sia un differenziale esatto e quello che non sia. Poiché, se, nel primo caso, interpretiamo dq semplicemente come simbolo di differenziale, con che non vi ha luogo a ricordare termini di ordine superiore, si può, conformemente, nel secondo caso interpretare dq come $\frac{dq}{dt}dt$; dove t è un parametro, di cui v e p sono funzioni, relative alla particolare serie di trasformazioni supposta, e

$$dq = C_p \frac{\partial T}{\partial v}\frac{dv}{dt} + C_v \frac{\partial T}{\partial p}\frac{dp}{dt} = \lim_{\Delta t = 0} \frac{\Delta q}{\Delta t},$$

e neppure in questo caso vi ha luogo a ricordare termini di ordine superiore. Se poi, nel secondo caso, interpretiamo dq (da te indicato con $\partial q$) come incremento, con che la (1) risulta ottenuta col trascurare l'aggiunta di termini di ordine inferiore, egualmente nel primo caso, dq non si può inter*[pretare]*[259] come incremento altrimenti che intendendo di trascurare i termini di ordine superiore, conformemente al noto sviluppo del Δq.

---

[256] Potrebbe trattarsi di T. Levi-Civita che scrisse diversi libri di Meccanica razionale con Ugo Amaldi.
[257] J.C. Maxwell, in una nota a p. 24 del testo *A Treatise on Electricity and Magnetism* (1873), scrive che il prof. W.H. Miller gli aveva suggerito l'idea che i viticci del vino e del luppolo sono rispettivamente destrorsi e sinistrorsi, così nominò i due sistemi di relazioni "of the vine and the hop".
[258] Potrebbe trattarsi di Moisè Ascoli (si veda la lettera #2) per l'argomento trattato.
[259] L'angolo è strappato.

Per quanto sl significato della (1), dove dq si interpreti come incremento infinitesimale, direi poi che, se si ammette il principio dell'equivalenza, va intesa dedotta da

$$\Delta q - (\Delta q)_{p=cost.} - (\Delta q)_{v=cost.} = \frac{1}{E} ABC$$

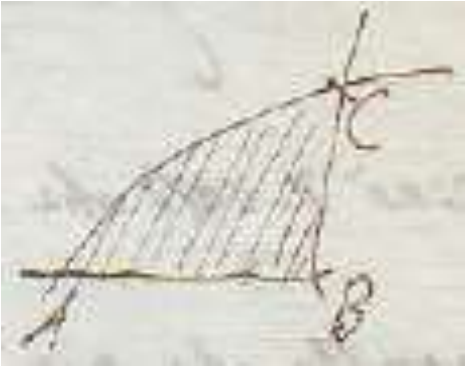

trasformando gli infinitesimali di ordine superiore.
E se si prescinde da quel principio, va intesa formata, applicando il principio della indipendente coesistenza delle variazioni infinitesimali, che si verifica appunto, trascurandol'aggiunta degli infinitesimali di ordine superiore, per cui, a mio modo di vedere, non ne segueil principio della conservazione del calore: allo stesso modo che l'espressione Xdx+Ydy+Zdz del lavoro infinitesimale non segue l'ipotesi del potenziale.

Non mandata

**229**

Gian Antonio Maggi a destinatario sconosciuto[260]

Pregiatissimo Signore,

Parecchie circostanze m'impedirono di redigere prima d'ora una risposta al Suo pregiato scritto in data 31 Ottobre 1885,[261] e alle posteriori Note, che, mi consegnò, la scorsa estate, a Milano. In quel primo lavoro, Ella usa benevolmente parole assai lusinghiere a proposito del mio scritto in data 5 Ottobre 1885,[262] delle quali La ringrazio vivamente. Esso non ebbe tuttavia troppa fortuna: poiché Ella mantiene in sostanza tutte le Sue obbjezioni: e invariabilmente insiste perché a 274+T sia sostituito T, dove fu mia prima preoccupazione di convincerLa che non si può, perché il numero 274, che, nel caso in discussione, non occorre rammentare che, preso col segno − rappresenta il così detto 0 assoluto di temperatura, si aggiunge spontaneamente e necessariamente a T, e non già perché noi leviamo da T −274, per riferirci allo 0 assoluto, piuttosto che allo 0 ordinario (§§ 5 e 6).[263] Confesso che, in difesa della termodinamica, non saprei trovare migliori ragioni di quelle che già Le esposi. Ma gioverà forse movere da un diverso punto di partenza.

1. Questo sia il seguente fatto:
Suppongasi che l'unità di massa di un gas perfetto qualsivoglia assorba la quantità di calore Q, e che T la temperatura (misurata nel modo ordinario) del gas, nell'istante che finisce d'assorbirla. Noi potremo fare infinite ipotesi sulle condizioni presentate dal gas, durante quell'assorbimento di calore: e a ciascuna di queste ipotesi corrisponderà una diversa dipendenza di Q dalla temperatura finale T.

In particolare, possiamo supporre che mentre il gas assorbe Q, il suo volume si mantenga costante, oppure che si mantenga costante la pressione a cui è assoggettato. In ambedue questi casi, l'assorbimento della quantità di calore Q, dev'essere accompagnato da un aumento della temperatura del gas, il quale, trovandosi in fine alla temperatura T, dovrà quindi aver cominciato ad assorbire la quantità di calore Q ad una temperatura inferiore, che indicherò con $t_0$: e, indicando con c un numero (che nei due casi anzidetti riceve valori diversi), sarà

$$Q=c(T-t_0) \qquad (1).$$

[260] Potrebbe trattarsi di Luigi Sala, per l'argomento trattato (si vedano le lettere #162 e #163).
[261] Non si trova nel Fondo.
[262] La lettera #163 è datata 1/10/1885 ed è piuttosto consistente, mentre non ne esiste una datata 5/10: Maggi potrebbe aver riportato il giorno in modo errato oppure potrebbe averla copiata qualche giorno dopo.
[263] Si riferisce ai diversi punti nei quali è divisa la lettera #163.

Supponiamo invece che il gas assorba la quantità di calore Q, mentre la sua temperatura si mantiene costantemente eguale a T. In questo caso, il gas, in conseguenza dell'assorbimento di calore, dovrà espandersi e passare da un volume $v_0$ ad un volume V. Dinotando J l'equivalente meccanica del calore, e R la costante dei gas, in questo caso si troverà

$$Q=\frac{R}{J}\log\frac{V}{v_0}(T+274). \qquad (2)$$

Nell'ipotesi che il gas compia un ciclo di Carnot, esso assorbirà una quantità di calore Q alla temperatura costante T della sorgente, e cederà al refrigerante una quantità di calore q, ad una temperatura costante t. Si avrà quindi per la (2), diventando $v_0'$, V' il volume del gas alla fine e al principio della sua contrazione in contatto del refrigerante

$$Q=\frac{R}{J}\log\frac{V}{v_0}(T+274)$$

$$q=\frac{R}{J}\log\frac{V'}{v_0'}(t+274).$$

Si dimostra che dovrà essere

$$\frac{V}{v_0}=\frac{V'}{v_0'},$$

e da ciò segue che, indicando con k una quantità costante per un medesimo ciclo di Carnot, sarà

$$Q=k(T+274) \qquad (3).$$

Così la (1) e la (3) si riferiscono ad ipotesi, che si deducono a vicenda: e cioè la (1) all'ipotesi che il gas, mentre assorbisce Q, mantenga costante il volume, o la pressione: la (3), al contrario, all'ipotesi che, mentre assorbisce Q, il gas mantenga costante la temperatura.

V.S. invece suppone sempre implicitamente che la (3) sia una conseguenza della (1), dedotta da essa col fare $t_0=-274$, e riferendosi ad un altro corpo, ch'Ella chiama maggiore o minore di quello a cui si riferisce (3), e ch'io supponga intenda di massa maggiore o minore. Ora questo, mi permetta di esprimermi francamente, è il Suo errore fondamentale. La (1) e la (3), lungi dall'essere la seconda una conseguenza della prima, si escludono a vicenda: e sono riferite entrambe all'unità di massa d'uno stesso gas: sol che sta la (1), se si suppone che questo corpo assorba la quantità di calore Q, a volume o a pressione costante, e sta invece (3) se si suppone che assorba quella quantità di calore a temperatura costante.

Ammesso che la quantità di calore Q si debba esprimere con (1), non v'ha dubbio che, ove Q riesca proporzionale a T+274, bisognerà supporre $t_0=-274$, e quindi che il gas assorba Q a partire da −274. Ma, se si suppone che il gas assorba Q a temperatura costante, Q non si deve esprimere colla (1), bensì colla (3): e, non essendo questa formola una conseguenza della (1), dalla circostanza che, per (3), Q è proporzionale a T+274 non si può inferire in alcun modo, com'Ella fa, che Q debba essere attinto a partire da −274.

Dai Suoi scritti mi sembra di poter desumere che Ella suppone che, quando un gas assorbe una quantità di calore Q, e T indica la temperaturafinale, la (1) debba necessariamente sussistere (posto per $t_0$ il debito valore). Invece non è così. A seconda delle ipotesi varrà la (1), o la (3): come potranno valere altre relazioni, tutte egualmente deducibili dalla relazione
Q=aumento di calore interno del gas
+quantità di calore equivalente al lavoro d'espansione,
che abbraccia tutti casi possibili.

2. Come da questa relazione generale, nell'ipotesi che il gas assorba Q, mentre la sua temperatura si mantiene costantemente eguale a T (nel qual caso l'aumento di calore interno è nullo) si deduca la (3), lo esposi diffusamente nel citato mio scritto (§ 5): e forse, colla persuasione che la formola stessa si possa ricavare assai più facilmente dalla (1), su questa deduzione Ella non ha riflesso abbastanza. Io debbo pregarLa di volerla rivedere attentamente: e non potrà far a meno di persuadersi come la (3) non implichi in alcun modo che il corpo considerato, o un altro corpo, che occorra tirare in scena, assorba Q a partire da −274: come il numero 274 non si aggiunga a T, nell'espressione di Q, perché noi sostituiamo le temperature assolute alle temperature contate dallo 0 ordinario, ma spontaneamente e necessariamente, introdotto dalla formola dei gas

$$pv=R(T+274),$$

per mezzo della quale, allo scopo di eseguire l'integrazione, esprimiamo la pressione in funzione del volume: come infine la (3) sia una diretta conseguenza della formola dei gas e del principio di Mayer, per modo che o si accetta la (3), o bisogna negare tutta quanta la termodinamica.

3. Un attento esame di quella deduzione La convincerà della necessità di ammettere la (3): e, ove Ella si persuada che la (1), col ciclo di Carnot, non ha nulla a che vedere, è probabile che la (3) non Le parrà neppure un inevitabile paradosso: né Le parrà che qualche paradosso contenga la relazione

$$\frac{Q-q}{Q}=\frac{H-h}{H},$$

che ne consegue; poiché, esclusa la necessità che Q sia attinto a −274 (necessità puramente condizionata alla credenza che la (3) discenda dalla (1)) parmi che non avrà più difficoltà ad ammettere le considerazioni che, in difesa di quella formola, da Lei attaccata, si trovano nel § 7 del mio scritto. Ciò mi lascia supporre quanto Ella afferma in proposito nel § 4 del Suo scritto in data 31 Ottobre 85.

A questo proposito osserverò che la dimostrazione, data nel § 10 di quello scritto, che indicando con d, A e P rispettivamente il calore utilizzato in un ciclo, la quantità attinta dalla sorgente, e la quantità preesistente nel corpo, il coefficiente economico non sarà giustamente rappresentato da $\frac{d}{A+P}$, se non nel caso che P sia 0, ossia che inizialmente il corpo si trovi allo 0 assoluto di temperatura, caso che non potrà mai verificarsi nelle macchine industriali, non potrebbe peccare in nulla. Il guajo è che V.S. vuol farsene un'arma, per negare che il coefficiente economico in discorso possa essere rappresentato da $\frac{H-h}{H}$; e a ciò io non posso che objettare, come già nel mio precedente scritto (§ 7) che, non essendo H−h il calore, d, utilizzato dal ciclo, vale a dire, per effetto del ciclo, convertito in lavoro disponibile, nulla toglie che, dividendo questa quantità diversa da d per H, ossia per A+P, invece che per A, si ottenga per quoziente lo stesso numero come dividendo d per A. Perché quella dimostrazione avesse il valore d'un'obbjezione, bisognerebbe che, non H−h, ma d, si dividesse per H=A+P, e si pretendesse di ottenere in tal modo il coefficiente economico del ciclo.

Finalmente, per quanto all'argomentazione ab absurdo, di cui nella IV[a] delle Sue "Obbjezioni fondamentali" mi limiterò ad osservare che la macchina ch'Ella suppone utilizzare costantemente $\frac{2}{3}$ del calore attinto dalla sorgente, funzionando fra limiti di temperatura diversi e tali che il rapporto $\frac{T-t}{T+274}$ non conserva lo stesso valore, si compone di un gas, che assorbe C.150 calorie a volume costante, per trasformarne in lavoro C.100 con un'espansione adiabatica, che ne abbassa la temperatura di 100 gradi. Con ciò esso utilizza realmente $\frac{2}{3}$ del calore assorbito: ma quelle due operazioni non costituiscono un ciclo chiuso qualsiasi, nonché un ciclo di Carnot: e perciò nessuna meraviglia che la relazione

$$\frac{Q-q}{Q}=\frac{T-t}{T+274},$$

che concerne il rapporto del calore utilizzato al calore attinto, in un ciclo di Carnot, applicata alla successione di quelle due operazioni conduce a risultati assurdi.

Alla successione di quelle operazioni che, non riconducendo il gas allo stato iniziale, ripeto, non costituiscono un ciclo chiuso qualsiasi, non si potrà neppure applicare la proporzione che il relativo coefficiente economico non possa superare quello di un ciclo di Carnot dove le temperature estreme siano le medesime. Feci già questa osservazione nel mio scritto più volte ricordato (§ 8), dove dissi anche come, immaginando una successione di operazioni, che, in seguito all'espansione adiabatica, riconducessero il gas allo stato iniziale, o formassero così un ciclo chiuso, trovai il coefficiente economico di questo ciclo 0, 136, mentre quello di un ciclo di Carnot, dove la temperatura della sorgente si suppone 150, e 50 quella del refrigerante, è 0,236. Terminerò colle seguenti considerazioni sulle osservazioni, che, a questo proposito, si trovano nel Suo scritto 31/8bre 85.

Ella dice (§ 9) che non sa se quel ciclo, da Lei ideato, sia un vero ciclo chiuso, ma sa che si uniforma al modo di agire delle macchine industriali a fuoco, e che vi è in esso il ritorno ad uno stato iniziale, o per lo meno un complesso di operazioni equivalenti a questo ritorno. Ora che susseguentemente al riscaldamento e all'espansione adiabatica del gas la macchina possa riprendere lo stato iniziale, sta bene: ma è patente che non lo riprende per opera di quelle sole due operazioni, che sono le sole ch'Ella considera per conchiudere che il coefficiente economico della Sua macchina sarà $\frac{2}{3}$. Al principio noi abbiamo una massa di gas a 0°: immaginiamola chiusa in un cilindro munito di stantuffo. Si comincia col riscaldare questa massa di gas da 0° a 150°, mantenendone costante il volume, e perciò tenendo fisso lo stantuffo: e con tale operazione si attingono dalla sorgente C.150 calorie (per ogni unità di massa): poi, liberato lo stantuffo, si lascia che il gas si espanda finché (supposto che durante l'espansione ne *[sic!]* riceva né ceda calore) la sua temperatura da 150° scende a 50°: operazione, colla quale C100 calorie (per ogni unità di massa) si trasformano in lavoro. A questo punto Ella non si occupa più della Sua macchina: della quale non occorrerebbe di occuparsi ulteriormente, se non premesse di riaverla periodicamente nello stato iniziale. Se questo si vuole, allora bisogna badare che in seguito alla 2[a] delle precedenti operazioni lo stantuffo si trova spostato d'un tratto corrispondente alla corsa compiuta per l'espansione del gas, e il cilindro pieno di una massa di gas a 50°. Per riavere la macchina nello stato iniziale, bisognerà che lo stantuffo riprenda la primitiva posizione, e, in quella posizione, chiuda nel cilindro un'egual massa di gas alla temperatura di 0°. Per ottenere questo risultato, è chiaro che bisognerà che alle 2 operazioni da Lei considerate ne seguano altre: e sono queste, che, richiedendo una spesa di lavoro e la dispersione dell'equivalente quantità di calore, che se ne ottiene per frutto, nella mia ipotesi, riducono il coefficiente economico della medesima a 0,136.

È con queste ultime operazioni che alla produzione di un certo lavoro a spese di un'equivalente quantità di calore attinta dalla sorgente, si contrappone la produzione di una certa quantità di calore che va disperso nel refrigerante, a spese di un equivalente lavoro. Alla fine dell'espansione adiabatica del gas del calore, da esso attinto dalla sorgente, non è stata ancora dispersa alcuna parte; sta bene che di C.150 calorie ne abbiamo utilizzato solamente C.100; ma le rimanenti C.50 non sono perciò andate disperse, e noi potremo utilizzarle come le precedenti lasciando espandere ulteriormente il gas, finché la sua temperatura riprenda il valore iniziale 0. Ella dice che arresta a quel punto l'espansione del gas, allo scopo di riavere la macchina nelle condizioni iniziali. Orbene, vediamo quali operazioni si devono fare succedere a tal fine all'espansione in tal modo limitata: e troveremo che questa operazioni richiederanno che si disperda nel refrigerante una quantità di calore ben superiore a C.50 calorie (io ho trovato C.129,6)

Solo coll'aggiunta di queste ultime operazioni, la macchina in discorso riprodurrà le condizioni delle macchine industriali, per esempio, della macchina a vapore. Difatti, come agisce una macchina a vapore, che, per fissare le idee, funziona fra le temperature (normalmente contate) di 150° e di 50°? In sostanza, e nelle migliori condizioni, nel seguente modo.

1°. Una massa d'acqua, dal condensatore mantenuto a 50°, è spinta, per mezzo della tromba alimentatrice, nella caldaja, mantenuta a 150°, e , attingendo calore dall'acqua che vi bolle, da 50° passa a 150°. La caldaja è in comunicazione col cilindro, munito di stantuffo: e questo si sposta in modo da limitare colla sezione inizialmente occupata un volume pari a quello dell'acqua introdotta. 2°. Una massa d'acqua bollente a 150°, eguale a quella dell'acqua introdotta (che, senza pericolo d'arrivare perciò a conclusioni diverse possiamo supporre la stessa massa d'acqua introdotta, portata a 150°) si trasforma in vapore saturo alla temperatura di 150° e quindi alla pressione di 5 atmosfere, e, passando nel cilindro, fa spostare lo stantuffo. 3° Tolta la comunicazione fra il cilindro e la caldaja, quella massa di vapore si espande senza comunicazione né cessione di calore, finché si riduce ad una miscela di vapore saturo alla temperatura di 50°, e quindi alla pressione di   atmosfera *[sic! c'è solo uno spazio, senza un valore]*, e di goccioline d'acqua a 50°: in seguito alla quale operazione lo stantuffo si avanza di un nuovo tratto. 4°. Posto in comunicazione il cilindro col condensatore, dove regna la pressione di   atmosfera *[sic! come sopra]*, propria del vapor d'acqua saturo a 50°, lo stantuffo è obbligato a retrocedere superando quella pressione e a comprimere quel residuo di vapore nel condensatore, dove, condensandosi in acqua a 50°, formatasi per l'espansione. Con ciò si ritorna alle condizioni iniziali: mentre nel condensatore ritorna la stessa massa d'acqua a 50°, che originariamente fu spinta nella caldaja, lo stantuffo riprende la posizione primitiva. Ora la 4ª operazione, necessaria a tal fine, non ha riscontro nel processo da Lei ideato, (almeno per quanto se ne occupa V.S.): in seguito al quale lo stantuffo resta appunto alla sezione, che raggiunge per l'espansione del gas.

La retrocessione dello stantuffo (che, in pratica, viene sospinto dal vapore, passato, pel giuoco del cassetto di distribuzione, dalla parte opposta del cilindro) assorbe un lavoro, misurato dalla quantità di calore, che genera, e si disperde nel condensatore. Supposto che la massa d'acqua vaporizzata fosse un chilogrammo, e ammesso (come si credeva prima dei lavori di Clausius) che, durante l'espansione, non si condensasse parzialmente in acqua, per modo che alla fine dell'espansione (3[a] operazione) si trovasse nel cilindro un chilogrammo di vapore saturo alla temperatura di 50°, quella quantità di calore ammonterebbe a 572 calorie, mentre il calore attinto dalla sorgente, nelle nostre ipotesi, sono 590 calorie. Il coefficiente economico della macchina sarebbe in questo caso $\frac{590\text{-}572}{590}=\frac{28}{590}=0{,}047$: mentre il limite massimo assegnato dalla termodinamica è la nota frazione $\frac{100}{424}=0{,}236$. Fortunatamente, ha luogo la condensazione, in virtù della quale il coefficiente economico si eleva a 0,218, di soli 0,018 minore del limite anzidetto.

S'Ella ben riflette su questi esempi, si convincerà anche facilmente che la macchina si obbliga a compiere un ciclo, con assorbimento di lavoro, non già per mantener sempre in giuoco la stessa massa d'agente intermediario, e far così economia dell'agente medesimo, ma per riavere periodicamente lo stato iniziale. L'abbondanza dell'acqua è una buona ragione perché poco c'importi di ciò che avviene del vapor acqueo, che lo stantuffo, per retrocedere, e riprendere la posizione iniziale, caccia dal cilindro. Ma per cacciarlo dal cilindro e riprendere la posiz. iniz. lo stantuffo deve superarne la pressione per un tratto corrispondente alla sua corsa e questa operazione assorbe un lavoro, il cui frutto, rappresentato dal calore prodotto dalla liquefazione di quel vapore, va buttato via. Nella locomotiva il vapore cacciato dal cilindro si disperde nell'atmosfera, e la caldaja è alimentata con quantità d'acqua sempre nuova: lo stesso avviene nelle macchine ad alta pressione. Orbene, il modo d'agire di queste macchine è in sostanza lo stesso come se vi fosse un condensatore mantenuto a 100°, e dove per conseguenza regnasse la pressione d'un'atmosfera: nel qual caso sarà indifferente supporre che nella macchina circoli sempre la stessa massa d'acqua, oppure soltanto un'egual quantità, senza che l'acqua spinta nella caldaja al principio di un ciclo sia precisamente quella stessa, che alla fine del ciclo precedente si è raccolta nel condensatore.

Finalmente, poiché ho detto che il coefficiente economico della macchina a vapore precedentemente considerata, tenendo conto della condensazione, risulta 0,218, osserverò che, se si confrontasse il lavoro effettivamente utilizzabile, svolto da un ciclo, col calore simultaneamente generato dal combustibile, che arde nel fornello, si troverebbe un numero certamente assai minore: poiché quel lavoro non sarà che una parte di quello che svolge il ciclo, e di quel calore, soltanto una parte servirà a vaporizzare l'acqua della caldaja, il resto andando disperso per l'irradiazione, la conduttività dei materiali ecc. Ma la circostanza che quel rapporto, che, per intenderci, potremmo chiamare il coefficiente economico, che potremmo chiamare <u>pratico</u>, è per le singole macchine, inferiore al corrispondente coefficiente economico teorico, non ha evidentemente altro risultato che di rendere applicabile a fortiori la proposiz. che una macchina non può trasformare in lavoro, per effetto d'ogni ciclo, una fraz. del calore attinto dalla sorgente, superiore a quella che la teoria assegna per limite massimo del coeff. economico d'una macchina che funz. tra le stesse temp. estreme, non toglie nulla alla proposizione che il coeff. economico di una macchina funzionante fra gli stessi limiti di temperatura, secondo un ciclo di Carnot, non possa essere superato, poiché, se non può esserlo dal coefficiente teorico, per tal ragione, <u>a fortiori</u>, non potrà esserlo dal coefficiente pratico.

*[Una seconda versione della* <u>*prima parte*</u> *di questa lettera:*
Pregiatissimo Signore,

Parecchie circostanze m'impedirono di redigere prima d'ora una risposta al Suo pregiato scritto in data 31 Ottobre 1885, e alle posteriori Note, che, mi consegnò la scorsa estate a Milano. In quel primo lavoro Ella usa benevolmente parole assai lusinghiere a proposito del mio scritto in data 5 Ottobre 1885, delle quali La ringrazio vivamente. Esso non ebbe tuttavia troppa fortuna: poiché Ella mantiene in sostanza tutte le Sue obbjezioni, e invariabilmente insiste perché, nella formola, che fornisce il limite massimo del coefficiente economico di una macchina termica, funzionante fra i limiti di temperatura T e t, si sostituisca T a T+274, dove fu mio principale intento di convincerLa che non si deve fare, perché il numero 274, che ivi non occorre rammentare che , preso col segno −, rappresenta il così detto 0 assoluto di temperatura, si aggiunge spontaneamente e necessariamente a T, e non già perché noi leviamo da T −274, per riferirci allo 0 assoluto, piuttosto che allo 0 ordinario (§§ 5,6). Confesso che, in difesa della Termodinamica, non saprei trovare migliori ragioni di quelle che già Le esposi. Ma gioverà farle movere da un diverso punto di partenza?*]*

*[Una seconda versione della* <u>*parte finale*</u> *di questa lettera:*
non sarà che una parte di quello che svolge il ciclo, e di quel calore soltanto una parte servirà a mantenere in ebollizione l'acqua della caldaja, il resto andando disperso per l'irradiazione, la conduttività dei materiali ecc. Ma la circostanza che quel rapporto, che potremmo chiamare, per intenderci, il coefficiente economico <u>pratico</u> è per le singole macchine, inferiore al corrispondente coefficiente economico teorico, non ha evidentemente altro risultato che di rendere applicabile <u>a</u> <u>fortiori</u> dal punto di vista dell'utilizzazione pratica la proposizione che una macchina non può utilizzare, <u>per</u> <u>effetto</u> <u>d'ogni</u> <u>ciclo</u>, una frazione del calore attinto dalla sorgente, superiore a quella che misura il coeff. economico d'una macchina che fra le stesse temp. estreme funziona secondo un ciclo di Carnot.

Termino con ciò le mie considerazioni, né ho bisogno di dirLe che, se Le resterà in proposito qualche difficoltà, mi sarà sempre grato ch'Ella mi comunichi le osservazioni che Le parrà di fare.

Aggradisca, Pregiatissimo Signore, i segni della mia perfetta stima, e mi creda sempre

Dev.mo Aff. mo Suo
Dr. Gian Antonio Maggi.

Modena 15 Dicembre 1885.*]*

**230**
Gian Antonio Maggi a destinatario sconosciuto[264]
[s.l. e s.d.; pagine numerate da 1 a 3]

La definizione di massa secondo il metodo di Mach da Lei modificato e reso rigoroso poggia sul seguente attributo delle figure materiali:
Due figure materiali, concepite isolate in partenza l'una dall'altra, si imprimono accelerazioni il cui rapporto è invariabile con ogni circostanza modificatrice del movimento
È questo rapporto che dà la misura della massa d'una fig. mat. rispetto all'altra assunta come campione.
L'impossibilità di ottenere due fig. mat. isolate l'una con l'altra toglie all'attributo in discorso ogni carattere di principio sperimentale.
Ho quindi pensato se si possa introdurre la massa per una via leggermente diversa, considerando cioè le accelerazioni mutue che due corpi s'imprimono non già concepiti isolati in presenza, ma collegati fra di loro e sotto l'influenza di azioni esterne.

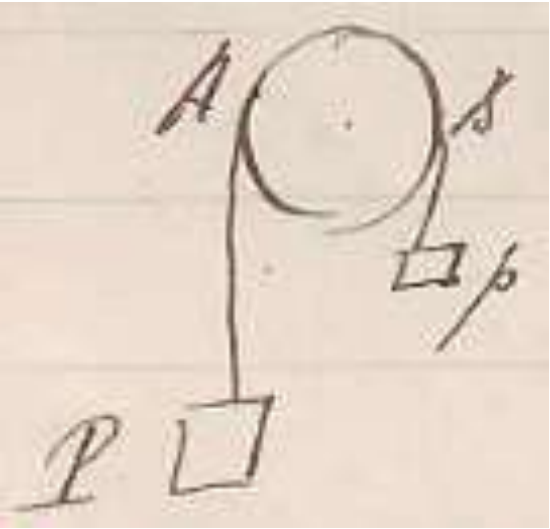

Imaginiamo *[sic!]* i corpi P e p legati da un filo che passa sulla gola di una carrucola. Se faccio cadere liberamente P e p, tagliando cioè i fili PA e pB, è noto che essi acquistano entrambi la stessa accelerazione, l'accelerazione di gravità g. Lasciando cadere invece il peso P collegato con p così come s'è detto si osserva che P (e anche tutto il sistema) acquista un'accelerazione g'<g. Ne deduciamo che il peso p imprime a P in questo caso (cioè sotto l'azione della terra e con questo speciale legame) l'accelerazione g−g'. Facendo lo stesso ragionamento per il corpo p si deduce che P imprime a p l'accelerazione g+g'. Il rapporto $\frac{g-g'}{g+g'}$ diremo che misura la massa del corpo p quando si assuma P per campione. La costanza di questo rapporto col tempo e con la velocità è evidente, giacché tali sono g e g'. Inoltre esso, per il modo stesso come l'abbiamo dedotto, dà l'indice della maggiore attitudine che ha P a muoversi rispetto a p quando entrambi siano, come nel caso, tratti a muoversi egualmente (adopro vocaboli impropri perché voglio solo giustificare la denominazione di rapporto da masse *[sic!]* data a questo rapporto), cioè con la stessa accelerazione per effetto delle cause esterne (in questo caso della Terra).
Ma in questo caso l'azione esterna (della Terra) si manifesta sempre con la stessa accelerazione g che impone a tutti i corpi. Si potrebbe supporre che al mutare di questa accelerazione anche il rapporto che ci dà la misura delle masse cambiasse. Per vedere la costanza del rapporto anche in questi casi ho pensato a due piani inclinati egualmente su un piano orizzontale. Se p e P cadessero sul piano senza essere fra loro collegati si muoverebbero con l'accelerazione g' entrambi. Si osserva che con il legame P si muove con l'accelerazione g''.

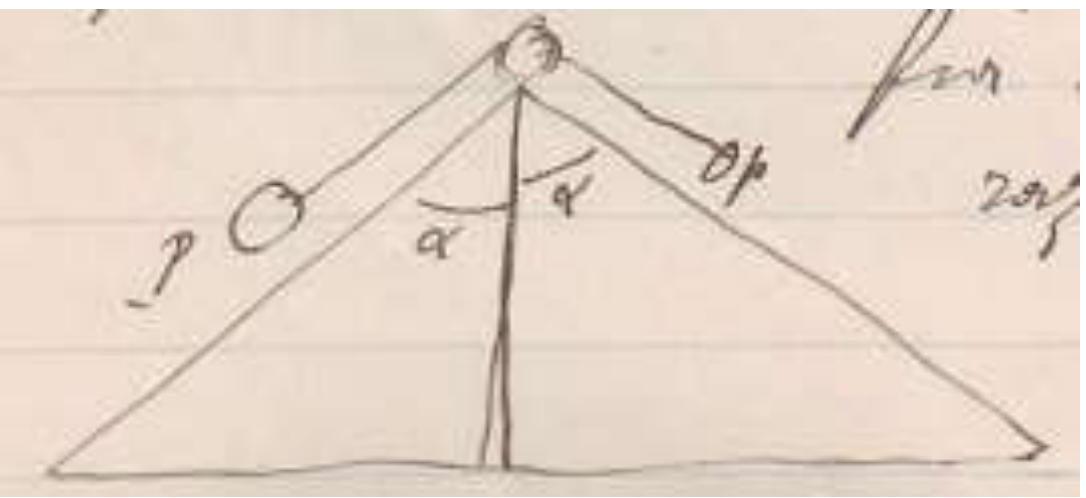

Il rapporto $\frac{g'-g''}{g'+g''}$ è l'analogo a quello di prima. Se ora si fa mutare α si può osservare che questo rapporto ha sempre lo stesso valore, malgrado cambi g'.
Si potrà infine far variare la specie di azione esterna, considerando l'apparecchio a forza centrifuga che il Mach ha proposto, ed osservare che il rapporto fra le accelerazioni che le due sfere collegate mediante la molla si *[*imprimono*]* sotto l’azione del momento rotatorio che noi imprimiamo sono in rapporto costante.
Da tutto questo si vede, secondo me, che la definizione di massa non ha bisogno assoluto dell’attributo dei corpi in presenza, anzi credo che per la via indicata il concetto di massa scaturisca più evidente, come indice di inerzia al moto.
Il caso dei corpi in presenza isolati non è che un caso più particolare quanto quelli accennati *[sic!]*, giacché si riduce al legame etere (l’etere non è che una molla) fra i due corpi e alla mancanza di azioni esterne. Per giustificare la definizione di massa dopo aver introdotto l’attributo di cui ho fatto cenno mi pare che si ricorra implicitamente al principio d’inerzia che a priori non è affatto evidente. Perché infatti il rapporto fra le accelerazioni che due corpi in presenza isolati s’imprimono deve dare l’indice del rapporto delle masse dei due corpi? Solo perché si suppone implicitamente che se i due corpi non si imprimessero accelerazioni mutue con la loro presenza essi non avrebbero accelerazioni.
Per quell’altra via invece la giustificazione della definizione di massa poggia sul fatto sperimentale che P e p acquistano, se fatti cadere liberamente, la stessa accelerazione.
Concludendo *[*non*]* pare che l’attributo delle figure materiali potrebbe subire per questa via qualche modificazione, e in ogni modo introdurre la massa nel modo da me indicato si presti benissimo almeno per una trattazione elementare.

[264] Si veda il carteggio con Polara, del quale questa minuta, insieme alla successiva, potrebbe far parte.

**231**

Gian Antonio Maggi a destinatario sconosciuto[265]

[s.l. e s.d. ; pagine numerate da 1 a 3]

Massa e forza motrice

Nel modo da me proposto di introdurre la massa in un insegnamento elementare è supposto che alla Dinamica preceda una breve trattazione della Statica: si suppone quindi che si abbia una prima idea di forza e di tensione.
Per una trattazione elementare io credo che non si possa spogliare il concetto di massa dell’altro comune della maggiore o minore inerzia che presentano i diversi corpi a cambiare di accelerazione sotto una stessa azione. Del resto, quando si assume per misura della massa d’un corpo rispetto ad un altro scelto come campione il rapporto inverso delle accelerazioni che i due corpi si imprimono quando sono isolati, si ha di mira il fatto che la loro presenza determina sui due corpi azioni eguali e contrarie, mentre le loro accelerazioni sono in rapporto costante: questo rapporto dà così l’indice della maggiore o minore inerzia al moto dell’un corpo rispetto all’altro sotto la stessa azione.
Ora a me pare che per dare questa definizione di massa che rispecchi questo criterio non sia necessario partire dal fatto che due corpi isolati si imprimono accelerazioni mutue in rapporto costante, ma si possa partire dalle variazioni di accelerazioni che subiscono due corpi in movimento qualsiasi (e quindi soggetti anche ad azioni esterne) quando si colleghino opportunamente.
Adoprando la carrucola di cui Le ho fatto cenno, si vede che, legando i due pesi P e p in quel modo, il corpo P subisce una variazione di accelerazione g−g’, mentre il corpo p subisce una variazione di accelerazione g+g’. Si riscontra che il rapporto $\frac{g-g'}{g+g'}$ resta costante.

[265] Si veda la nota precedente e la lettera precedente per il contenuto.

Dobbiamo quindi argomentare che il legame introduce qualcosa (tensione) che modifica le accelerazioni dei due corpi nel rapporto $\frac{g-g'}{g+g'}$. E poiché questa tensione agisce ugualmente sui due corpi. tanto che il filo resta teso (nel caso della molla essa assume una configurazione d'equilibrio), ne segue che la diversa variazione di accelerazione che han subito i due corpi è dovuta alla diversa inerzia dei due corpi sotto la stessa influenza. (La variazione di accelerazione è quindi intimamente collegata con le masse dei due corpi)
In vista delle considerazioni fatte in principio diremo che questo rapporto misura la massa d'un corpo quando si assuma l'altro come campione.
Potrebbe pensarsi che questo rapporto non assumesse lo stesso valore col variare delle azioni esterne che agiscono sui due corpi. Le esperienze con i due piani inclinati o con l'apparecchio a forza centrifuga provano che ciò non è.
Quale è ora la funzione del legame? A me pare che essa sia solo di rendere sensibile il mutuo comando del movimento dei due corpi, anche in questo caso più generale di corpi soggetti ad azioni esterne. È, in certo modo, un espediente per isolare i due corpi, per modo che essi siano soggetti da parte del legame a tensioni eguali e contrarie (restando sempre per le azioni esterne che agiscono sui due corpi).
E da questo vediamo di che natura deve essere il legame. Esso deve essere tale che il movimento di uno dei due corpi (che assumerebbe se non ci fosse il legame) possa influire su quello dell'altro: in tal caso il filo di legame sarà manifestamente teso e le tensioni che esso determina sui due corpi saranno eguali e contrarie.
Mi pare quindi che si possa così concludere:
Dati due corpi P e p in movimento e soggetti ad azioni esterne qualisivogliano, se si collegano in modo che muovendo P nella direzione nella quale si muoverebbe se non ci fosse il legame anche p si muova nella direzione in cui esso (p) si muoverebbe se non ci fosse il legame, essi subiscono variazioni nelle loro accelerazioni il cui rapporto non varia col tempo, e col variare delle azioni esterne. Tale rapporto è uguale al rapporto inverso delle masse.
Naturalmente queste condizioni per il legame sono soddisfatte nel caso della carrucola, dei due piani inclinati e dell'apparecchio a forza centrifuga.
Mi pare inoltre che il fatto che il rapporto delle accelerazioni che due corpi isolati si imprimono è inverso a quello delle loro masse sia un caso particolare di quanto ora accennato. Il legame infatti che in generale introduciamo per rendere sensibile il mutuo comando dei due corpi dev'essere tale, come ho già detto, che le tensioni da esso esercitate sui due corpi siano eguali e contrarie. Ora il legame etereo interposto fra i due corpi è proprio tale, giacché le azioni mutue dei due corpi, che si possono attribuire a questo legame, sono eguali e contrarie.
Infine, la questione di dedurre dal metodo da me indicato per definire la massa le leggi fondamentali della dinamica, mi pare molto semplice. Intanto, poiché queste riguardano le forze bisognerà introdurre la definizione di forza, dicendo che ad un corpo di massa m che possieda l'accelerazione c sia impressa la forza (1) f=mc. E si può far vedere come questa definizione dinamica di forza collimi con quella statica. Se infatti il peso P ha subito una variazione di accelerazione (g−g') vuol dire che il legame esercita su di lui una forza, secondo la (1) M(g−g'), e analogamente se p ha subito una variazione di accelerazione (g+g') il legame esercita su di lui una forza m(g+g'). Per la definizione di massa si ha M(g−g')=m(g+g'), ossia questa due forze sono eguali.

Questo concorda perfettamente con il fatto delle due tensioni eguali e contrarie del filo che abbiamo dedotto in principio dalla circostanza che il filo è teso (cioè per via Statica).
Del resto poi, ammessa la (1) si può seguire il metodo generalmente adottato, giacché tutte le deduzioni si fanno formalmente partendo dalla (1).

**232**

Gian Antonio Maggi a destinatario sconosciuto
[s.l. e s.d., ma estate 1936, perché si trovava nelle *Rusticationes* di quell'anno]

Lanzo d'Intelvi (Como)

Carissimo amico,

Non posso trattenermi dal far seguire al telegramma queste righe, colle quali mi sembra di poter venir Loro più da vicino. Ma che altro posso dire, a Lei, mio carissimo, alla Mamma, alla Sorella della cara scomparsa, se non è che partecipo, con tutto l'animo e con tutto il cuore, ad uno schianto, di cui comprendo intera l'acerbità? Per qualche particolare di così dolorosa fine a cui nulla poteva preparare, per notizia di Loro, ho subito scritto al nostro G*[...]* che supponevo a Milano. Ancora non ne ho ricevuto risposta. Ma non l'aspetto oltre per rendermi Loro quest'altra volta presente. Perché, senza timore di esacerbare un estremo dolore, so pure che, se, in questi momenti della vita, in cui le nostre forze a stento bastano per reggerci, è possibile un conforto umano, questo è il sincero compianto di chi sappiamo che ci vuole veramente bene.

Col più affettuoso abbraccio

il Suo ...

## Trascrizione dei *Giudizii*

**I**

[s.l. e s.d.]

L'argomento della Dissertazione, assegnato al candidato, a termini dell'Art. 131 del Regolamento, 21 Agosto 1905, era "L'Esperienza nelle Teoria di Maxwell e di Lorentz o l'interpretazione meccanica dei fenomeni elettrici". Il candidato lo ha svolto in forma di un opuscolo di 169 pagine (Catania, Galatolo), nel quale esprime succintamente i principali argomenti attenenti alle più recenti teorie del campo elettro-magnetico, sotto il doppio aspetto matematico e sperimentale, li coordina, come occorre per metterne in rilievo la mutua connessione, e li correda di proprie osservazioni, le quali particolarmente riguardano le relazioni della Meccanica Classica colle conclusioni della suddetta teoria. La teoria di Maxwell e le sue conclusioni, verificate dall'esperienza, sulla condizione fisica dei fenomeni luminosi, formano il primo capitolo. Accennato che alcuni fatti, superati in questi tempi, non rientrano in tale teoria, segue un capitolo delle ipotesi degli elettroni fondata sui fenomeni dei raggi catodici, della ionizzazione dei gas, della radioattività. E da questi passa ai fondamenti della teoria elettronica di Lorentz, e alle sue più immediate conseguenze, quali la propagazione delle azioni elettro-magnetiche, le tensioni nell'etere, la massa elettromagnetica, la formazione dei raggi X. Il cap. IV è dedicato all'emissione e assorbimento del calore, in quanto le nuove teorie immaginate per questo fenomeno si connettono con quella del campo elettro-magnetico. Il cap. V ha per oggetto il fenomeno di Zeemann *[sic!]*, sotto le diverse forme, e la relativa spiegazione fondata sull'ipotesi degli elettroni. Il cap. VI confronta le due ipotesi dell'etere mobile coi corpi e dell'etere fisso, e, accennate le ragioni che stanno per la seconda, vi è dedotta, in base alla celebre esperienza di Michelson e Morley, la relativa ipotesi della contrazione lorentziana. Il principio della relatività di Einstein, esposto e corredato di delucidazioni sul significato delle relative conclusioni, concernenti le misure della lunghezza e del tempo forma l'oggetto del cap. VII. Nel cap. VIII sono rilevate le divergenze fra la meccanica elettromagnetica e la meccanica classica, ed è reso conto della teoria di Thomson, che, al movimento dell'etere fornito di massa, procura di eliminarle, e della teoria di Lewis,[266] che si propone di cambiarle, abbandonando solo la costanza della massa, con che la meccanica classica riesce una prima approssimazione della teoria del movimento. Il cap. IX concerne l'esistenza dell'etere, i rapporti tra l'etere e la materia, e l'immagine elettromagnetica dell'universo fisico. Infine il cap. X contiene le osservazioni di Ritz[267] intorno alla teoria di Lorentz.

Tale la materia dell'opuscolo. Era forse desiderabile una più particolareggiata trattazione degli argomenti sperimentali, per modo che questi entrassero in alquanto maggior proporzione, come pure, nella parte critica, una maggior sintesi, che rappresentasse più chiaramente, fra concetti diversi, quello dell'autore. Tuttavia è un lavoro pregevole, per abbondanza di materia, ampiezza di disegno, chiara ed originale esposizione, il quale attesta nel candidato soda coltura e elevate doti di studioso, oltre di che costituisce una Monografia, che si presta utilmente ad iniziare allo studio del campo elettromagnetico, col suo più ampio significato.

[266] Potrebbe trattarsi di Gilbert Newton, il quale però pubblica su questi argomenti nel 1908 e 1909.

[267] Potrebbe trattarsi di Walter, il quale pubblica su questi argomenti a partire dal 1908.

## II

Dr. A. Signorini. Sulla teoria analitica dei fenomeni luminosi nei mezzi cristallini uniassici.[268]

Caduta la soluzione delle equazioni di Lamé, che Sofia Kovalevskij aveva creduto trovare, colla critica di Volterra, il quale rilevava inoltre che gl'integrali di Lamé, rivelando infiniti sugli assi ottici, e polidromi secondo i circuiti che circondano un asse ottico, non si prestano a rappresentare il movimento luminoso, emanante da un unico centro, si riapriva la questione, se ed in quanto il movimento per superficie d'onda, e la costruzione di Huyghens, fondata su di esso, nei mezzi cristallini, si conciliano coi principii della teoria elastica, o, se si vuole, della teoria elettromagnetica della luce. La qual questione è tanto più interessante, perché la costruzione di Huyghens fornisce risultati conformi all'esperienza. Essa forma l'oggetto del presente lavoro, almeno per quanto concerne i mezzi cristallini uniassici, che sono pure particolarmente importanti, come i meglio studiati sperimentalmente.

Il S. comincia perciò con una ricerca delle soluzioni elementari, corrispondenti al movimento per superficie d'onda, delle equazioni di Lamé, ridotte al caso dei mezzi uniassici, e adattando a questo caso il procedimento generale di Volterra, ritrova le formole di Brill, che comprendono gli integrali elementari di Lamé.

Applica poi a questi integrali formalmente il principio di Huyghens, adattando al presente caso il procedimento trovato da Poincaré, per dedurre la formola di Kirchoff nei mezzi isotropi, e giunge, per questa via, alla formola della Kovalevskij. Ma, introducendo tali espressioni al posto delle componenti dello spostamento luminoso, nelle equazioni di Lamé, ridotta al caso dei cristalli uniassici, constata che, in generale, esse non vi soddisfanno, né soddisfanno la condizione di incompressibilità. È da notarsi che questa conclusione non scaturisce senz'altro da quella di Volterra, perché le singolarità degli integrali elementari si riducono col coincidere dei due assi ottici. Indaga poi direttamente e rileva le circostanze analitiche, che, nel presente caso, infirmano tale processo di integrazione.

La precedente verifica conduce il S. al risultato che le formole della Kovalevskij, anzi che alle equazioni di Lamé, soddisfanno

$$\frac{\delta^2\xi}{\delta t^2} - a^2\Delta_2\xi - (a^2 - c^2)\frac{\delta}{\delta y}\left(\frac{\delta\eta}{\delta x} - \frac{\delta\xi}{\delta y}\right) = X, \quad \frac{\delta^2\eta}{\delta t^2} - a^2\Delta_2\eta + (a^2 - c^2)\frac{\delta}{\delta x}\left(\frac{\delta\eta}{\delta x} - \frac{\delta\xi}{\delta y}\right) = Y$$

$$\frac{\delta^2\zeta}{\delta t^2} - a^2\Delta_2\zeta = Z, \qquad (1)$$

con particolari forme di X, Y, Z (Z=0).

Ora, queste equazioni si riducono alle equazioni di Lamé, supponendo X=Y=Z=0, e stabilendo la condizione d'incompressibilità

$$\frac{\delta\xi}{\delta x} + \frac{\delta\eta}{\delta y} + \frac{\delta\zeta}{\delta z} = 0:$$

al quale scopo è necessario e sufficiente (in conseguenza delle stesse equazioni) che soddisfacciano ad una simile condizione, oltre i valori iniziali di $\xi$, $\eta$, $\zeta$, anche quelli di $\frac{\delta\xi}{\delta t}$, $\frac{\delta\eta}{\delta t}$, $\frac{\delta\zeta}{\delta t}$.

[268] *Annali della Scuola Normale Superiore di Pisa. Classe di Scienze*, v. 12, 1912, pp. 1-133.

Perciò, il S. si propone e conduce a fornire l'integrazione delle (1), con X, Y, Z e valori iniziali delle $\xi$, $\eta$, $\zeta$, $\frac{\delta\xi}{\delta t}$, $\frac{\delta\eta}{\delta t}$, $\frac{\delta\zeta}{\delta t}$ qualisivogliano. Facendo poi l'ipotesi che detti valori sono nulli per $t=-\infty$, trova le formole note, nella stessa ipotesi, del Grünwald, come soluzione del problema del movimento luminoso, determinato da date forze perturbatrici.

La precedente ricerca fornisce al S. degli integrali semplici dalle equazioni (1) - con X=Y=Z=0 - le quali, pel complesso delle loro proprietà, rappresentano l'estensione al presente caso della soluzione $\frac{f(r-at)}{r}$ relativa ai mezzi isotropi.

Egli se ne vale per dedurne le formole corrispondenti, sotto il medesimo punto di vista, alla Formola di Kirchoff, con procedimento analogo a quello di Kirchoff. Questa formole riescono poi verificate da formole più generali, dedotte con altro procedimento, indicato al S. da Grünwald.

Un'importante applicazione di queste formole fa il S. colla determinazione del movimento luminoso, emanante da un unico centro. Supposto che per $t=\infty$ il movimento in discorso sia nullo, ne risulta che il movimento per superficie d'onda rappresenta, per approssimazione, il movimento luminoso emanante da un unico centro, dopo un tempo sufficiente *[sic!]* lungo. A differenza delle soluzioni di Lamé, risulta connesso lo stato delle due falde oltre di che la singolarità è ridotta al centro d'emanazione. Questo risultato collima con uno analogo ottenuto da Grünwald, ma, sotto molteplici aspetti, ne resta distinto.

Le formole precedenti non possono assumersi puramente come espressione del principio di Huyghens, perché il movimento corrispondente ad ogni elemento del contorno del campo considerato non soddisfa la condizione di incompressibilità. Mediante una trasformazione analitica il S. ne deduce delle formole che soddisfanno anche a questa condizione. Esse coincidono con quelle trovate da Conway per via del tutto diversa.

Inoltre il S. fa un'applicazione delle precedenti formole, che rappresentano l'estensione della Formula di Kirchoff, al problema specifico della doppia rifrazione uniassica di un pennello di luce, valendosi di un procedimento analogo a quello usato da Kirchoff pei mezzi isotropi, e ritrova le note circostanze, verificate dall'esperienza.

L'analisi, occorrente per la trattazione di codesto argomento riesce, di sua natura, piuttosto complicata. È da rilevarsi la cura, che ripetutamente pone il S. di verificare con rigore, e coll'applicazione a casi particolari, la correttezza de' suoi risultati.

Il S., con questo lavoro, in parte abilmente coordina, con uno schema originale, alcuni risultati ottenuti da altri, sotto varii punti di vista, in parte vi aggiunge il proprio contributo di qualche notevole risultato, e d'interessanti applicazioni. Egli fornisce con esso non dubbia prova di ampia conoscenza degli strumenti dell'analisi, di perizia nel maneggiarli, di perspicacia nella discussione dei risultati, infine di quella attitudine alla sintesi di elementi diversi, che procura un'opera efficace per la trattazione dei problemi della Fisica Matematica.

G.A. Maggi

Con questo lavoro il S. in parte esamina e abilmente coordina, con uno schema originale, risultati poco diffusi, ottenuti da altri, sotto diversi punti di vista, in parte vi aggiunge il proprio contributo di alcuni risultati importanti, di una nuova analisi, di interessanti applicazioni.

**III**
[su carta intestata: Regia Università di Pisa
Facoltà di Scienze Fisiche Matemat. e Naturali]

Il Principio di Relatività dell'Elettrodinamica
di Otto Berg.[269]

Il Principio di Relatività dell'Elettrodinamica consiste in questo che, conseguite due terne d'assi in movimento traslatorio, rettilineo, uniforme, l'una rispetto all'altra, e due osservatori collegati con esso, in modo da giudicare fissa la rispettiva terna, questi due osservatori giudicheranno, ciascuno, che la luce emanata da un punto, ad un istante, si propaga per onde sferiche, colla stessa velocità di propagazione, e con quel punto - fisso - come centro. Questo "postulato" di Einstein concilia l'ipotesi dell'etere immobile, reclamata dal presente assetto della teoria del campo elettro-magnetico, col risultato dell'esperienza di Michelson che, pur raggiungendo un estremo grado di sensibilità, i fenomeni luminosi fra corpi terrestri, non rivelano il movimento annuo della Terra. Ne conseguono, per altro, conclusioni paradossali sul giudizio che i due osservatori rendono della posizione e del tempo: due avvenimenti non sono contemporanei per ambedue, ognuno di essi giudica più breve un segmento appartenente allo spazio con cui è collegato l'altro, e, in un punto di questo, più lento il decorrere del tempo. Einstein ha presentato così sotto l'aspetto di un diverso giudizio relativo della posizione e del tempo la "trasformazione di Lorentz", che lo stesso Lorentz fondava sul concetto del "tempo locale" e sull'ipotesi di una effettiva contrazione dei corpi in movimento. Minkowski ha procurato, alla sua volta, di presentare sotto un nuovo aspetto questa trasformazione, e il principio di relatività, col concetto del "punto universale" (Weltpunkt). Posizione e tempo sono indissolubilmente connessi, l'una con l'altro, ma, per la descrizione dei fenomeni naturali, sono indifferenti diverse determinazioni della posizione e del tempo, insieme connesse, purché ne risulti invariata la rappresentazione analitica di una certa quadrica iperspaziale, dove la velocità di propagazione della luce acquista il significato di limite superiore delle velocità dei corpi naturali. La suddetta diversità dei giudizii trova in questo la sua ragione.

L'A. fa una lucida ed abbastanza ampia esposizione di questa questione, particolarmente importante per la Fisica Matematica, considerata in sé stessa e pelle sue conseguenze. Esamina in fine la possibilità che l'esperienza di Michelson possa trovare altre spiegazioni. Poiché il principio di relatività "non si deduce dall'esperienza, ma risponde piuttosto ad una esigenza metafisica" e la sua efficacia è subordinata al mantenimento dei concetti informatori della teoria di Maxwell, per amor della quale si è disposti, se occorre, a far getto della Meccanica Classica e degli antichi concetti dello Spazio e del Tempo. "Singolare destino" (conclude l'A. con parole di Ritz) "di una teoria, di cui, durante la vita dell'autore, nessuno voleva sapere".

Per il prof. A. Forte, al quale mandata il 29 Luglio 1913.

L'Articolo si trova in un Numero del 1910 dell'Abhandlungen der Fries'schen Schule etc.

[269] Si tratta della traduzione dal tedesco di: O. Berg, "Das Relativitätsprinzip der Elektrodynamik", *Abhandlungen der Fries'schen Schule*, Vandenhoeck & Ruprecht, 1910.

## IV

Relaz. del Concorso di Pavia Bollettino etc.
Anno 1913, 1° Sem. pag. 1480 [5 Giugno]
Giudizio Complessivo

Daniele[270] I lavori del Daniele sono frutto di elaborazione coscienziosa. Egli ha [recato] tendenza (fra) buoni contributi in più campi: equilibrio delle reti, determinazioni infinitesime delle superficie, (tendenza) dei sistemi materiali a sfuggire l'attrito, centri di vibrazione, armoniche ellissoidali.

Per quanto abbia coltivato di preferenza la Meccanica Analitica, appare manifestamente maturo anche per la Fis. Matem. che ne è la naturale continuazione oltre l'ambito dei fenomeni immediatamente avvertiti quale movimento dei corpi ponderabili.

N.B. Presenta 31 pubblicazioni.

Relaz. del Concorso di Pavia. 5 Giugno 1913
Boll. Uff. etc (Vol. I pag. 1480)

Laura[271] Il Laura si rivela conoscitore profondo della Meccanica Analitica, e di quelle teorie della Fisica Matematica che più direttamente si collegano colla Meccanica. Presenta interessanti contributi, ben maturati, di Meccanica e di Acustica, e uno studio ampio e sistematico delle vibrazioni smorzate. A lui spetta indubbiamente un posto onorevole nel presente Concorso.

N.B. Le pubblicazioni arrivano alla 12.

## V

Tre pubblicazioni fatte dal prof. Daniele dopo il concorso di Pavia sono di argomento puramente analitico. Per cui il giudizio resta, in via assoluta, lo stesso come pel suddetto concorso.

[270] Ermenegildo Daniele: è libero docente di Meccanica Razionale a Pavia dall'a.a. 1902-03 al 1914-15.
[271] Ernesto Laura: diventerà professore ordinario a Pavia sulla cattedra di Meccanica Razionale solo per l'a.a. 1921-22.

# VI

## Concorso di Meccanica Razionale pel Politecnico di Torino 1915.

## Prof. Ernesto Laura

Presenta il prof. Laura, in aggiunta alle pubblicazioni giudicate nel Concorso di Pavia (v. Boll. Uff. 5 Giugno 1913), altre cinque, numerate XI, e poi da XIII a XVI. La XIII è uno studio, a base di calcoli, dell'equilibrio elastico di un solido isotropo vincolato, con ipotesi informate al concetto di Saint Venant:[272] del quale interessante risultato è il collegamento del problema in discorso con quello delle distorsioni di Volterra. I quattro rimanenti lavori formano un gruppo, concernente il movimento ondulatorio provocato da un contorno agente, e precisamente il problema esterno, secondo il generale concetto di Love del fronte mobile, parallelo al contorno, che, completando il teorema di Huygens-Kirchhoff, ha dato forma matematica al principio delle onde inviluppo, inclusa la circostanza della inattività della coda. Il Laura postula l'esistenza e le circostanze caratteristiche del fronte di Love, ammettendo però la velocità c (di propagazione delle onde sferiche libere) per la velocità $\omega$ del fronte secondo la normale, donde, nelle applicazioni, la identificazione di certe condizioni cinematiche e dinamiche. Con questo mi sembra che il Laura non attribuisca sufficiente importanza alla circostanza, che il Love dimostra, della identità di $\omega$ con c, la quale invece ha tanta parte nella traduzione del principio in forma matematica.

La Nota XI contiene una dimostrazione, col metodo dei coni caratteristici nello spazio a quattro dimensioni (x, y, z, t) (Volterra-Tedone), della formola di Huygens-Kirchhoff, per lo spazio esterno al fronte. Il Love ne dà una dimostrazione che trovo preferibile sia dal lato della semplicità, sia perché mette in rilievo la proprietà importante che l'integrale esteso al fronte (nell'ipotesi consueta della continuità della funzione) è quello, che s'interpreta colla suddetta inattività della coda. Inoltre non è esatta l'identificazione di questo problema col problema esterno nel primitivo significato: per questo sono indispensabili le proprietà asintotiche della funzione $\varphi$. E se non regge il discorso di Kirchhoff, conforme all'intuizione è l'ipotesi che, all'infinito, $\varphi$ si comporti come $\frac{f\left(t-\frac{r}{c}\right)}{r}$, sufficiente per rappresentare $\varphi$ nel campo esterno.

I lavori XIV, XV, XVI si possono considerare come parti, sovrapposte in alcuni punti, di una medesima ricerca: la rappresentazione dello spostamento elastico sullo spazio esterno all'onda di Love. Nella XVI è completamente risoluto il problema per le onde sferiche, e dimostrata la sovrabbondanza dei dati dimostrati (nella XIV) atti a determinare univocamente la soluzione. Questi lavori consistono precisamente in calcoli laboriosi, i quali (verifica a parte) dimostrano la perizia e la coltura del concorrente: e i risultati rappresentano certamente un positivo contributo ad un importante problema della Fisica Matematica.

Il favorevole giudizio espresso sul Laura nel Concorso di Pavia non potrebbe essere altrimenti che avvalorato e accresciuto dai nuovi lavori, sia dal lato dell'attività che da quello del valore.

[272] Adhémar-Jean-Claude Barré de Saint-Venant.

Prof. Ing. Giuseppe Armellini.

Laureato in Ingegneria nel 1910, e due anni dopo in Matematica, presenta quindici pubblicazioni che si riferiscono tutte alla tesi dell'attrazione newtoniana, considerata sotto diversi aspetti, i quali trovano la loro applicazione in corrispondenti questioni di Meccanica Celeste. La maggior parte di questi lavori si può distribuire in tre gruppi, che fanno capo alle due maggiori Memorie, pubblicate negli Atti dei XL, le quali ne riassumono e completano i risultati. Il primo gruppo riguarda il movimento newtoniano di due punti materiali di massa variabile, in moto prestabilito, col tempo τ e la bella Memoria a cui fa capo può dirsi veramente esauriente, contenendo un copioso esame delle circostanze del movimento, dedotte direttamente dalle equazioni differenziali: la riduzione dell'integrazione al principio della stella di Mittag Leffler: la deduzione di una soluzione approssimante, atta alle applicazioni alla Meccanica Celeste. Il secondo gruppo riguarda il movimento newtoniano di più punti, di massa data, in presenza l'un dell'altro, e particolarmente il moto di un punto attratto da più punti fissi, problema a cui si riferisce l'accennata Memoria. Due Note riguardano invece il problema degli n corpi propriamente detto. L'autore riconduce la soluzione al principio della stella di Mittag Leffler, ispirandosi al concetto enunciato, per la prima volta dal Volterra, nel caso speciale di punti fissi in linea retta, con costante delle aree diversa da 0, ma generalizzando col valersi di un tempo ridotto, mediante il quale la regolarità delle funzioni incognite si dimostra estesa a *[*tutta una*]* lista comprendente l'asse reale dei tempi. Generalizza ancora l'analisi del Sundman coll'estendere i risultati da punti materiali a sferette infinitesimali, guadagnando così un significato concreto, di grande interesse, con ingegnoso ragionamento. Se la soluzione generale del proposto problema è implicita nella serie di Mittag Leffler, soddisfanno abbondantemente le esigenze della descrizione del movimento i numerosi e interessanti risultati sulle circostanze dell'urto. Escono dai limiti degli indicati gruppi una Nota sullo schiacciamento di Giove, dedotto dal moto del V satellite, con risultato pienamente conforme all'esperienza, e questo lavoro, è da contarsi come uno dei più notevoli: una Nota preliminare sulle perturbazioni del satellite di Nettuno: o la Nota sulle comete iperboliche, che è un saggio elegantissimo di soluzione di un complicato problema di Calcolo delle Probabilità.

La produzione dell'Armellini già notevole per l'attività del giovane autore, è pregevolissima per la scelta d'importanti problemi, per originalità e genialità di ragionamento, per perizia di condotta, per la deduzione rigorosa e completa d'interessanti risultati.

Spetta indubbiamente all'Armellini uno dei primi posti, nel presente Concorso.[273]

Gian Antonio Maggi

[273] Armellini vince il Concorso con voto unanime della Commissione; alcuni lavori commentati qui da Maggi sono disponibili in rete sul sito dell'Accademia Nazionale dei Lincei.

# VII

## Prof. Antonio Signorini.

Meritevole del maggior rilievo è l'opera scientifica del prof. Signorini, che, nel periodo relativamente breve di un settennio,[274] ha pubblicato venti lavori, sopra svariati argomenti di Analisi, di Meccanica e di Fisica Matematica. Iniziata quest'opera con ricerche di Geometria Differenziale, informate a recenti lavori del Bianchi (1, 2), dove però il Signorini dimostra la propria distinta iniziativa, i seguenti scritti di Analisi recano notevoli contributi a questioni - come il teorema di Stokes sulla forma del Geoide (4,5) e il Teorema di Wittaker *[sic!]* sulle orbite periodiche che interessano del pari l'applicazione alla Meccanica e alla Fisica Matematica. I primi lavori di Fisica Matematica - una lunga Memoria (3) e una Nota (7) esauriscono la questione del significato del Principio di Huyghens nei mezzi cristallini uniassici. Una Memoria (10), completata da una Nota (11) confermano il risultato di tradurre in equazioni, di cui sono anche sviluppate le prime conseguenze, il problema del campo elettrico, senza introdurre i vincoli cinematici, oggetto di molteplici discussioni. Nei lavori (12, 13-14, 17) sul moto di un punto, e di un sistema olonomo, soggetto a resistenza idraulica e forza elastica, l'autore, da un semplice enunciato, sviluppa, con singolare perizia, una copia di risultati, che rappresentano una completa soluzione del proposto problema. I più recenti scritti (15, 18-20) trattano a fondo, sulla base delle equazioni di Hertz, il problema del movimento dell'elettricità in un conduttore cilindrico e toroidale, con manifesto vantaggio della teoria della corrente elettrica nel reoforo concreto. Le questioni, di varia specie, che imposta il Signorini, attraggono l'attenzione per la loro importanza intrinseca; e la condotta attesta una vera maturità di spirito matematico: i risultati ne rappresentano un'effettiva soluzione. Spetta al Signorini, nel presente Concorso, un posto emergente.

[274] Signorini si laureò nel 1909 e iniziò subito a pubblicare, quindi si può dedurre che il presente scritto sia databile al 1916. Inoltre vinse la cattedra a Palermo proprio in quell'anno ed è probabile che questa relazione sia relativa a quel concorso.

## VIII

Dr. Lucio Silla.
(Cfr. Elenco dei Documenti e delle Pubblicazioni[275]).

1). Monografia, che, per la condotta, e la copia della materia, ben risponde al proposto scopo di iniziare a orientare il lettore nell'importante questione del "Principio di Dirichlet".
2). È pure, per la massima parte, un ben condotto e bene esposto lavoro di compilazione, dedicato all'esposizione dei principali saggi di riduzione del duplice teorema delle velocità virtuali e postulati che soddisfacciano maggiormente l'intuizione. Il "nuovo metodo" dell'autore si trova a pag. 31, e si fonda principalmente sulla ipotesi che i vincoli siano rappresentati da azioni rettilinee tra i punti. L'A. distingue poi gli spostamenti virtuali in normali per cui $\sum Rdr = 0$, e non normali, per cui trova $\sum Rdr > 0$.
3). Appartiene ai saggi di subordinare alla definizione statica la definizione cinetica della forza.

4). Analisi di un caso di equilibrio in cui il teorema di Dirichlet è in difetto per la ragione che il massimo del potenziale non è effettivo. Non saprei come questo caso - segnalato da Hadamard . sia "assai interessante" e il metodo non mi sembra "molto semplice" come afferma l'autore.

5). Contiene una parte teorica, e una parte - la seconda - che, esaminando varie applicazioni del giroscopio, interessa più particolarmente la Tecnica. La parte teorica riproduce ragionamenti elementari intesi a render ragione delle circostanze del movimento del giroscopio, mancanti del desiderabile rigore, e, in qualche punto, di esattezza [nel movimento rotatorio di un piano, v e $\frac{dv}{dt}$ perpendicolare al piano], i quali dovrebbero essere oramai abbandonati per applicare al problema del giroscopio le equazioni generali della Dinamica. La teoria della Bilancia Giroscopica si trova a §74 della mia Stereodinamica,[276] e quella del Giroscopio di Foucault al §116.

6). 7) La prima 6) e un riassunto della seconda 7). La ricerca concerne il movimento della "stella doppia elettronica", introdotte da Righi, per spiegare i raggi magnetici e l'A. generalizza l'ipotesi del Righi, e arriva alla deduzione generale delle equazioni del movimento, per mezzo delle funzioni ellittiche. Sotto l'aspetto meccanico, è un esempio di movimento piano di un punto, con tempo dipendente dalla posizione e dalla velocità. Ma sono notevoli, per la Fisica, le osservazioni di carattere concreto, e, per l'Analisi, l'abile discussione dei risultati.

8). È risoluto, coi principii della Meccanica, senza ricorrere al Calcolo delle Variazioni, il problema: Trovare una catenaria tangente ad un cerchio, tale che la somma dell'arco di catenaria e dell'arco di cerchio, limitati ai punti di raccordo, abbia un valore assegnato. La ricerca si riduce alla determinazione delle costanti relative al tratto catenario (introdotte dall'integrazione delle relative equazioni) e della lunghezza di uno dei due tratti. L'A. parte da considerazioni generali: e coll'esempio in questione intende, se non erro, indicare il metodo di risoluzione, che, d'altronde, apparisce abbastanza ovvio e naturale.

Nei rimanenti lavori si rileva l'ispirazione dell'opera del compianto Lauricella.[277]

---

[275] Non è presente nel Fondo.
[276] Il volume è del 1903.
[277] Giuseppe Lauricella morì nel 1913.

9). Estende un risultato del Lauricella ad ipotesi più generali, valendosi di un teorema di Weyl sulle serie che procedono per Funzioni Ortogonali.

10) Ricerca preliminare, e complementare, della seguente: condurre la dimostrazione di formole di risoluzione di un'equazione integrale e di un sistema di equazioni integrali di 1° specie, per mezzo delle autofunzioni corrispondenti al nucleo.

11) I, II, III. Riduzione del problema dell'equilibrio elastico con dati spostamenti in superficie alle equazioni integrali, con principio o metodo conforme a quello usato dal Lauricella pel problema di Dirichlet sul piano.

12) Nota di Analisi, concernente lo sviluppo in serie *[precedente]* di funzioni ortogonali, fondata sui più recenti lavori del Lauricella, sullo stesso argomento delle fz. ortogonali.

I, II. Buone dispense, ad uso didattico. Salvo le mie riserve sulla terminologia, e sull'introduzione del vettore col punto di applicazione.

L'opera scientifica del Dr. Silla, alquanto lenta e modesta sul principio, s'intensifica in fine, e sorge a più ardue questioni, coi più recenti lavori. Questi però sono ricerche di indole sostanzialmente analitica, dove si rileva manifesta la traccia dell'opera del compianto Lauricella. Delle pubblicazioni appartenenti più strettamente alla Meccanica, alcune hanno piuttosto carattere di esercizii, altre di monografia. Se poi da codest'opera non si potrebbe dire che emerga pregio di originalità, è notevole, in compenso, la felice condotta e la perspicuità dell'esposizione, qualità che si ritrovano nei corsi litografati, ad uso degli studenti, ed avvalorano la fama di chiaro ed efficace insegnante, che, col suo non breve esercizio, nell'Università di Roma, si è acquistata il concorrente.

## IX

Prof. Lucio Silla.

L'opera scientifica del Dr. Silla, alquanto lenta e modesta sul principio, s'intensifica in fine, e sorge a più ardue questioni, coi più recenti lavori. Questi però sono ricerche di indole sostanzialmente analitica, dove si rileva manifesta la traccia dell'opera del compianto Lauricella. Delle pubblicazioni appartenenti più strettamente alla Meccanica, alcune hanno piuttosto carattere di esercizii, altri di monografie. Particolarmente buona, per le considerazioni di carattere concreto, sotto l'aspetto fisico, e per l'abile condotta e completa discussione dei risultati, sotto l'aspetto dell'Analisi, è la Nota sui raggi magnetici, e, se vogliamo, sulla stella doppia elettronica, secondo il concetto di Righi. Non si potrebbe dire, a mio giudizio, che dall'opera del Silla emerga pregio di originalità, ma è notevole, in compenso, la felice condotta e la perspicuità dell'esposizione: qualità che si ritrovano nei corsi litografati, ad uso degli studenti, ed avvalorano la fama di chiaro ed efficace insegnante, che, col suo non breve esercizio nell'Università di Roma, si è acquistata il concorrente. Posto onorevole egli tiene certamente nel presente Concorso.

Gian Antonio Maggi

Prof. Matteo Bottasso.

L’opera del concorrente Prof. Bottasso rivela essenzialmente un geometra. Oltre che alla Geometria, nel senso più esteso della parola, appartiene gran parte de’ suoi lavori, è per la via di questa dottrina, cioè del Calcolo Geometrico, ch’egli è passato nel campo della Meccanica. Ma ivi ha concentrato, quasi esclusivamente, la sua attività sulla Statica, trattandovi e risolvendovi complesse e delicate questioni di indole prettamente geometrica, e mostrando il vantaggio che ritrae un abile uso dei nuovi metodi vettoriali. Escono da questo campo le recenti note sopra un problema di valori al contorno, e sull’equilibrio delle piastre elastiche; ed è da segnalarsi, come pure informato a diverso concetto, oltre che per gl’intrinseci pregi, la nota sul movimento del filo flessibile e inestensibile. Non potrò esimermi dal dichiarare che sarebbe stato desiderabile, per gli effetti del presente Concorso, che il prof. Bottasso aggiungesse maggiori saggi di ricerche nel campo più genuino della Meccanica Razionale. Ma la perizia e l’ingegno che rivela la sua produzione, coll’ampio affidamento che forniscono, gli assegnano indubbiamente una onorevole posizione.

Gian Antonio Maggi

Prof. Ing. Goffredo Mancini

La Memoria “Studii preliminari per una turbina a gas” e i due articoli nel giornale “L’Industria” sono, il primo, di Termotecnica, e i due seguenti di Meccanica Applicata. Sono ricerche elementari che non so quanto possano interessare la Scienza Applicata: certo né interessano la Meccanica Razionale, né servono a fornire indizii della coltura del concorrente in questa materia. Nella Memoria l’uso della formola $T_2=T_1\sqrt[4]{\frac{p_2}{p_1}}$ invece di $T_2=T_1\left(\frac{p_2}{p_1}\right)^{1-\frac{1}{k}}$ non saprei come giudicarlo. Il libro “Meccanica Costruttiva” è un semplice trattato scolastico elementare di materia applicata. La Memoria sugli effetti giroscopici, il solo lavoro che, per l’argomento appartiene alla dottrina della cattedra posta a Concorso, non potrebbe, più del resto, sotto alcun aspetto, fornire prova dell’idoneità del candidato.
Le poche pubblicazioni dell’Ing. G. Mancini, concernenti, salvo l’ultima, questioni tecniche, sembrerebbero più indicate per un concorso di Meccanica Applicata. Non saprei, per conto mio, quanto interesse potrebbero offrire sotto codesto aspetto. Certo, non si saprebbe trarne prova della preparazione del concorrente, per la cattedra posta a concorso, e piuttosto ci sarebbe motivo per concludere l’opposto

Gian Antonio Maggi

## X

## Decomposizione intrinseca del vortice e sue applicazioni
Nota dell'Ing. Dr. Arnaldo Masotti

Il Dr. Masotti trova, in primo luogo, quelle che chiama componenti intrinseche del vortice relativo ad un punto qualsivoglia di un fluido in movimento, e cioè le sue componenti secondo le direzioni principali, nel punto, della linea di flusso passante pel punto medesimo le quali vanno concepite, nel campo occupato dal fluido, come funzioni del punto medesimo. Dalle espressioni così trovate deduce poi relazioni interessanti tra gli elementi geometrici della linea di flusso ed elementi cinetici del supposto movimento, donde, in varie ipotesi, che sogliono essere considerate nell'Idrodinamica, si deducono notevoli conseguenze. Reputo la Nota meritevole di essere inserita nei Rendiconti dell'Istituto.[278]

Gian Antonio Maggi.

## XI

## Aldo Pontremoli. Sulla conducibilità elettrica e termica nei metalli.

Il prof. A. Pontremoli si propone la deduzione di espressioni teoriche della conducibilità elettrica e termica di un metallo, nell'ipotesi che l'elemento conduttore abbia la struttura rivelata dai raggi X, attribuendo la corrente elettrica e termica al moto di elettroni, paragonabili a quello di molecole in un gas ultrararefatto, entro tubetti in cui s'immagina di poter dividere un cubetto unitario del considerato conduttore. Arriva con questo a formole formalmente coincidenti con quelle del Lorentz, salvo che al posto del cammino libero dell'elettrone di conduzione interviene il diametro del tubetto elementare. Dalle trovate formole deduce, valendosi anche di darne spiegazioni, diverse considerazioni, tra cui la variabilità della conducibilità colla direzione, che confronta con relativi risultati sperimentali. Accenna infine a un nuovo metodo per le distanze interatomiche che se ne potrà ricavare, per questo sia, se le previsioni teoriche saranno verificate dall'esperienza.

La sezione si pregia proporre l'inserzione di questo saggio del valoroso coltore della nuova teoria fisica nei Rendiconti dell'Istituto.[279]

---

[278] *Rendiconti Istituto Lombardo*, v. 60, 1927, pp. 869-874.
[279] *ibidem*, pp. 273-280.

## XII

Per l'Istituto Lombardo
Febbrajo 1928.

Sull'impiego delle isogamme.
Nota di Arnaldo Belluigi.

L'A. si propone la determinazione della variazione della grandezza g dell'accelerazione di gravità, mediante la misura, per mezzo della bilancia di Eötvös, della componente, secondo un asse orizzontale, del gradiente G della stessa g, in una serie di punti; formando i vertici di una poligonale, in base a $\Delta g=\int_{x_1}^{x_2} G_x dx$, che se $x_2-x_1$ è abbastanza piccolo, fornisce sensibilmente $\Delta g=G_{x_1}(x_2-x_1)$. Da questa ricerca potranno essere dedotte le linee, lungo le quali la variazione di g, nel campo orizzontale esplorato, è nulla, che si chiamano isogramme. Premessa l'applicazione ad una poligonale chiusa, dove dalla somma delle differenze osservabili $\Delta g$ relative ai singoli lati si può dedurre una correzione da ripartire sulle supposte osservazioni, sviluppa le formole nell'ipotesi di una variazione lineare e quadratica della $G_x$, e conclude con una formola svincolata da simili ipotesi, fondata, con un processo di interpolazione dalle eseguite osservazioni.

Reputo che la Nota del Sig.r Arnaldo Belluigi possa essere ammessa alla Lettura, e accolta nei Rendiconti dell'Istituto Lombardo.[280]

Milano 15 Febbrajo 1928

N.B. L'A. potrà tener conto sulle bozze di alcuni desiderati emendamenti di forma accennati con matite colorate.

Gian Antonio Maggi.

## XIII

Concorso Cagnola. a 1928[281]
Aude audenda.
L'azione dei campi elettrici intermolecolari sulla emissione delle righe spettrali.

L'A., premesse alcune considerazioni sulle varie circostanze atte a modificare le modalità delle righe spettrali, tra cui l'effetto Stark-Lo Surdo corrispondente a campi elettrici intermolecolari, si propone, nella presente Memoria, di indagare le leggi che possono governare l'allargamento e lo spostamento delle righe spettrali, per azione dei campi elettrici intermolecolari, come pure le mutue relazioni tra questi fenomeni, allo scopo di proporne convenienti basi teoriche alle ricerche sperimentali.

A tal fine, prende le mosse dalla formola che, secondo la teoria del Sommerfeld dell'effetto l'effetto Stark-Lo Surdo, esprime la frequenza di un componente emesso sotto l'azione di un campo elettrico costante, e, al campo costante, vi sostituisce un campo elettrico variabile colla legge di probabilità del Maxwell: con che ne trae, in conformità a questa legge, un'espressione della probabilità che l'inserimento dovuto al campo elettrico intemolecolare, della frequenza della radiazione $i^{ma}$, corrispondente ad un certo salto della quantizzazione, sia $t_i$ a meno di $dt_i$. Da questo risultato deduce poi un'espressione di $I_\sigma d\sigma$ ($I_\sigma$ intensità della riga di frequenza $\sigma$ a meno di $d\sigma$), come somma delle quantità omologhe appartenenti alle 2n+1 componenti, in cui la radiazione considerata è decomposta dall'effetto Stark-Lo Surdo, computata ciascuna componente in ragione della suddetta probabilità i dove però per le $\sigma_i d\sigma_i$, appartenenti alle singole componenti, è sostituito $\sigma d\sigma$, mentre rimangono i coefficienti dipendenti da i, cioè dai diversi salti di quantizzazione. Di questa sostituzione non mi risulta che l'A. renda ragione. Deduce in seguito alcune conseguenze dell'espressione ridotta, all'effetto limitato al 1° grado, cioè alla formola fondamentale arrestata al termine lineare nella forza elettrica del campo. Riproponendo poi l'ipotesi dell'effetto di 2° grado, deduce del suddetto risultato un'espressione dell'allargamento e dello spostamento di una riga, per azione di un campo elettrico variabile, come il supposto, di cui sia data la media quadratica, la quale non è altro che la nota costante della originaria legge di probabilità del Maxwell. Notevole conseguenza di quelle espressioni è che il rapporto tra lo spostamento della riga e il quadrato della sua semilarghezza è indipendente dell'intensità del supposto campo intermolecolare.

[280] A. Belluigi, "Su l'impiego delle isogamme", *Rendiconti Istituto Lombardo*, v. 61, 1928, pp. 420-426.

[281] Premio istituito dalla *Fondazione scientifica Cagnola* (ente nato nel 1848 dal medico milanese Antonio Cagnola come atto testamentario) che consisteva in L. 2500 e una medaglia d'oro del valore di L. 500. La presente relazione, scritta dalla Commissione formata da Murani, Maggi e Somigliana, riguarda il concorso scaduto il 31 dicembre 1927 (pubblicata sui *Rendiconti* del 1928) il cui tema era: "L'atomo del Rutherford: ipotesi del Bohr su la emissione e sull'assorbimento dell'energia per *quanta*. Portare qualche contributo allo studio della questione dal lato teorico e sperimentale".

Sopra l'applicazione dei precedenti calcoli alla riga $H_v$, studiata da Takamine e Kokeibo, con risultati numerici, a proposito dei quali è da notare che l'unità a cui si riferiscono le misure delle due suddette quantità risulta quella stessa della frequenza. Infine l'A., dopo brevi notizie su modificazioni del fenomeno prodotto da variazioni di pressione, accenna a esperienze preliminari da lui eseguite sulla stessa $H_v$, per pressioni varianti fino ad un'atmosfera, le quali gli fornirono un andamento che qualitativamente corrisponde alle previsioni teoriche, e allega una fotometria, che mostra la $H_v$ - tra le righe di riferimento - per una pressione di 10 cm di mercurio, della quale figura alcune indicazioni esplicative non mi parrebbero superflue.

L'A. rivela conoscenza del recente e delicato argomento, e interessante certo è il problema, alla cui soluzione si è proposto di recare contributo. Ma il lavoro mi sembra lasciar abbastanza desiderare in maggior grado consistenza di costruzione, chiarezza di esposizione, e conferma sperimentale.[282]

Milano 15 Febbrajo 1928

Gian Antonio Maggi.

[282] Infatti la relazione conclude: "Un solo concorrente. La relazione propone e l'Istituto approva che il premio non venga conferito". Non compare il nome del candidato.

## XIV

Prof. Bruno Finzi. Integrazione per successive approssimazioni delle equazioni di un liquido viscoso in moto stazionario.

Il Prof. Finzi introduce nell'equazione vettoriale che traduce le equazioni euleriane del movimento stazionario di un liquido viscoso l'ipotesi che la velocità v relativa al punto generico sia sviluppabile in serie di potenze intere di V, coefficiente di viscosità, e postulando, assai plausibilmente (a malgrado di una accennata obbiezione di C.W. Oseen) che v tenda a $v_0$, velocità del fluido perfetto, col tendere di V a 0, ne ricava una successione di equazioni differenziali atte a determinare successivamente i coefficienti $v_0$, $v_1$, $v_2$ ... della serie, inteso di saper determinare successivamente ciascuno di essi, in basa alla cognizione del precedente, da $v_0$ in poi. I coefficienti così calcolati danno una serie convergente entro il debito cerchio di convergenza. Applica poi il procedimento all'ipotesi del moto piano, in cui comparisce lo sviluppo in serie della funzione di corrente, e al moto per piani di un fascio, in cui presenta lo sviluppo in serie della funzione di Stokes, questa, come la precedente, relativa a v. Introduce in seguito l'ipotesi che il moto del fluido perfetto corrispondente a V=0, sia irrotazionale ($\mathrm{rot}v_0=0$), senza che sia, in generale, quello del fluido viscoso ($\mathrm{rot}v\neq0$), con che la determinazione dei coefficienti $\Psi_n$ delle precedenti serie, si collega colla teoria dell'equazione $\Delta_2\Psi_n=\Phi_n$ (nel caso del fascio di piani, riferendosi però ad uno spazio non euclideo). Fa seguito la considerazione delle condizioni al contorno, con cui univocizzare la soluzione delle accennate equazioni differenziali.

L'ipotesi che il contorno sia formato da linea di flusso fornisce, in ambedue i problemi, $\Psi_n=0$, $n>0$, e allora, pel caso del moto piano, si assegnano le $\Psi_n$ nei punti interni coll'uso della funzione di Green. L'ipotesi del fluido aderente al contorno risulta da escludersi, come impossibile. Determinato il moto, una formola generale fornisce anche lo sforzo applicato ad un elemento. L'A. vi si trattiene brevemente nell'ipotesi del fluido debolmente viscoso - $V^2$ trascurabile di fronte a V. Chiudono il lavoro alcuni esempii relativi al moto piano, nel quale, come si è accennato, si può condurre più oltre la determinazione delle incognite del problema.

Reputo il lavoro del prof. Finzi perfettamente meritevole di essere ammesso alla lettura all'Istituto e accolto nei Rendiconti.[283]

Gian Antonio Maggi.

Milano 30 Maggio 1928

[283] *Rendiconti Istituto Lombardo*, v. 61, 1928, pp. 695-715.

**XV**

[su carta intestata R. Università di Milano - Istituto Matematico
Via C. Saldini, 50 (Città degli Studi)]

Antonio Signorini. n. Arezzo, 1888: dott. in matem. (Pisa) 1909: abilitato all'insegn. della matem. (Pisa, Sc. Norm. Sup.) 1911: ingegnere civile (Palermo) 1921: medaglia d'oro dei XL (1920): prof. di Fis. Matem. nella R. Univ. di Napoli.

Esordì, come allievo di Luigi Bianchi, con lavori di geometria differenziale (La trasformaz. $B_k$ delle superf. applic. sulle quadriche nello spazio ellittico, Ann. della Sc. Norm. di Pisa, 1909, e Sulla permutabilità della trasform. H colla trasform. $B_k$ ecc., Annali di Matem. 1910. *[sic! non chiude la parentesi]*

Tosto dopo poi cominciò a dedicarsi a questioni affatto diverse, di Fisica Matem. colla rilevante Memoria Sulla teoria analitica dei fenomeni luminosi nei mezzi cristallini uniassici, Ann. della Sc. Normale, 1911 che si può dire esaurisca il denso problema, nella convenuta ipotesi. Della quale è un conciso estratto la Nota dei Lincei Sulle vibrazioni luminose di un mezzo uniassico, 1911.

Seguono, tra il 1914 e il 1916 Note pure di Fisica Matematica. Una a parte riguarda la Dinamica: Moto di un punto soggetto a resistenza idraulica e forza di richiamo, Ist. Veneto 1914. Caratterizzazione energetica dei moti soggetti a resistenza viscosa o idraulica, Lincei 1914: e, nell'indirizzo idrodinamico piano del Levi-Civita: Sull'inizio dell'efflusso dei liquidi Circ. Matem. di Palermo, 1916. E sull'idrodinamica, sotto diverso punto di vista, ritorna con le Note ai Lincei del 1928: Espressione asintotica di una formola del Levi-Civita e Sul teorema di Kutta-Joukowski. Tutt'altra parte riguarda invece la teoria affatto distinta dalla precedente del campo elettromagnetico: o, anche in questa parte, alcune concernono la dinamica degli elettroni, nella fase prerelativistica, Sulla dinamica dell'elettrone svincolata da arbitrarie ipotesi cinematiche, N. Cimento 1912, e Sui moti dilatatorii di un elettrone sferico Circ. Mat. di Palermo 1914: mentre altre concernono la propagazione delle onde elettromagnetiche nei conduttori corporei; tre Note dei Lincei: Sulla propagaz. di onde elettromagnetiche in un conduttore cilindrico 1914, Resistenza effettiva e resistenza ohmica, 1915, Sulla propagaz. di onde magnetiche in un conduttore toroidale, 1915.

Colla guerra, alla quale prese attiva parte come ufficiale di Artiglieria, e dopo, il Signorini riprende le pubblicazioni con ricerche di balistica: Sul moto dei projettili di bombarda, Ist. Veneto 1917 dove considera, oltre il consueto moto del baricentro, anche il moto riferito al medesimo, e, problema affatto diverso, Un teorema di confronto in balistica esterna ed alcune sue applicazioni Circ. Mat. di Palermo 1919, e Sull'integraz. appross. delle equaz. classiche della balistica esterna, Memoria dei Lincei 1928, dove è tracciato e sviluppato un metodo originale di approssimazione consistente nel rinserrare la trajettoria proposta tra due trajettorie di cui si assegnano le circostanze.

Da qualche anno l'attività del Signorini si trova dedicata alla Statica del cemento armato. Un teorema di esistenza ed unicità nella statica dei materiali poco resistenti a trazione, Lincei 1925, Sulla statica del cemento armato Seminario di Roma 1926, Sulla pressoflessione del cemento armato Ann. di Matem. 1928. Memoria, quest'ultima, dove il Signorini svolge per intero la sua teoria informata a concetti originali sulla costruzione e le condizioni d'equilibrio di sistemi atti a rappresentare il cemento armato.

Infine appartengono a svariate questioni di analisi le Note Sul criterio di Stephanos Lincei 1911, Esistenza di un'estremale chiusa dentro un contorno di Whittaker, Palermo 1912, Sopra un problema al contorno nella teoria di fz di variab. complessa Ann. di Matem. 1916.

La produzione del Signorini, consistente in una trentina fra Note e Memorie, se non s'impone per mole, spicca, in compenso, per varietà di argomenti, originalità di concetti, genialità di condotta e significato di risultati, pregi che mi sembra non si possa mancare di riconoscere in nessuno dei suoi accennati lavori.

## XVI

Relazione della Commissione Giudicatrice
della Libera Docenza in Fisica Teorica del Dr. Gleb Wataghin.

Il Dr. Gleb Wataghin, nato a Birzula, Circondario di Ananiev, Vescovado di Cherson (Russia) nel 1899, compiuti gli studii in patria, a Kiev, li riprese all'Università di Torino, dopo un'interruzione di due anni, per la rivoluzione e la guerra civile, che determinano l'emigrazione della sua famiglia. Laureato in quell'Università in Matematica e in Fisica, nel 1922, con pieni voti, vi occupò successivamente diversi posti di Assistente in Fisica e in materie matematiche, e attualmente vi tiene quello di Incaricato della Meccanica Razionale, oltre di che dall'anno scolastico 1925-26 fu nominato Incaricato di Calcolo Infinitesimale nella R. Accademia di Artiglieria e Genio, e dal 1928-29 Incaricato dell'insegnamento della Fisica Sperimentale, nella stessa Accademia.

Presenta sedici pubblicazioni attinenti particolarmente alla Fisica Teorica.

Una Nota "Ueber einige Periodizitätseigenschaften von mechanischen Systemen und Quantentheorie" (Annalen d. Physik 1925)[284] fa parte a sé. L'A. vi dimostra che per i sistemi multiploperiodici che posseggono u periodi, vale la relazione ρ=u, indicando ρ la dimensione della varietà densamente riempita da ogni trajettoria del punto (p,q) nello spazio delle relative fasi.

Un gruppo di tre lavori è dedicato ad una discussione della teoria balistica della propagazione dell'energia luminosa. (Sulla dipendenza della velocità della luce dal movimento della sorgente, Lincei, 1925. Sull'ipotesi balistica e l'effetto Doppler, Lincei, 1926. Ueber eine experimentelle Prüfung der ballistischen Hypothese, Zeitschrift f. Physik 1926). La discussione comprende un'ampia parte teorica sugli argomenti recati a favore e in opposizione alla teoria, e l'esperienza ideata e eseguita dall'autore, valendosi dei raggi canali. Le conclusioni dell'autore sono per l'esclusione dell'ipotesi balistica, per mancata conferma delle sue conseguenze da parte dell'esperienza.

La Nota Ueber die Quantenbedingungen (Zeitschrift f. Physik 1925) ha per oggetto la dimostrazione che le condizioni di quantizzazione di Bohr-Sommerfeld, per una classe estesa di sistemi meccanici, che comprende tutti i casi considerati da Sommerfeld, possono essere ridotte ad una condizione sola che vi equivale.

Nella Nota Sull'aberrazione della luce e la teoria della relatività (Lincei, 1926) l'autore presenta l'aberrazione annua come un semplice effetto del movimento della terra intorno al Sole, liberando questo fenomeno dalle note discussioni, con intervento della Relatività, suscitate da altre spiegazioni.

[284] Questa Nota si trova online sul sito onlinelibrary.wiley.com

La Nota Teoria della Diffrazione svolta in base alla meccanica ondulatoria (Nuovo Cimento, 1925) e Beitrag zu einer wellenmechanischen Theorie der Fraunhoferschen Beugungserscheinungen (Zeitschrift f. Physik 1927) contengono una teoria dei fenomeni di diffrazione di Fraunhofer, fondata sulla meccanica ondulatoria. Considerato perciò un reticolo spaziale atomico, si procura di spiegare la distribuzione dell'intensità nello spettro di diffrazione ottenuto coll'incidenza sul reticolo di un'onda piana mediante la teoria di Schrödinger.

Un gruppo di quattro Note (Sulla possibilità di conciliare la teoria ondulatoria delle interferenze luminose coll'ipotesi dei quanti di luce (Nuovo Cimento, 1927) Gl'integrali generali di alcune equazioni differenziali della fisica matematica (Lincei 1928). Versuch einer Korpuskularen Theorie der Interferenz und Beugung (Zeitschrift f. Physik 1928) Sulla teoria dei quanti di luce (Nuovo Cimento 1929) sono dedicate ad un tentativo di dare una teoria corpuscolare delle interferenze luminose, mediante un modello, informato ai principii della meccanica ondulatoria, che rende ragione della circostanza dell'interferenza senza violare i principii della Dinamica e della Fisica.

Specialmente la Nota Sulla teoria dei quanti di luce (Nuovo Cimento 1929) ha per oggetto la dimostrazione della possibilità di far derivare tutte le equazioni della Dinamica Quantistica da un unico principio variazionale, del tipo del teorema di Hamilton.

La Commissione riconosce nei lavori del Dr. Wataghin il pregio della scelta di argomenti importanti, concernenti punti svariati della Fisica Teorica, della lucida interpretazione del proposto problema, dell'ingegnosa ideazione della soluzione, della correttezza e purezza della condotta per ottenerla. Se qualche riserva può essere fatta sulla consistenza di alcune conclusioni, di ognuna deve pur riconoscersi che reca alla trattata questione un interessante contributo.

Svolse la lezione sul tema proposto "I fondamenti della meccanica ondulatoria" dimostrando cultura profonda nel nuovo e delicato argomento, originalità nella sintetica composizione dello svolgimento, perfetta chiarezza e sicurezza d'esposizione.

La Commissione pertanto reputa unanimente che al Dottor Gleb Wataghin possa essere conferita la Libera Docenza in Fisica Teorica. [285]

[285] Wataghin ottiene la Libera Docenza in Fisica Teorica presso l'Università di Torino nel 1929, che è quindi la presumibile data di questo scritto.

**XVII**

[su carta intestata R. Università di Milano - Istituto Matematico
Via C. Saldini, 50 (Città degli Studi) - Il Direttore]

Relazione sulla Nota del Sig.r Guido Facciotti
"Sulla corrispondenza fra spostamenti e deformazioni in un solido elastico"

Sotto l'indicato titolo, l'A. si applica in sostanza, a dedurre le componenti dello spostamento infinitesimale dall'ipotesi che le componenti del tensore di deformazione siano rappresentati da polinomii di grado n, con che le componenti dello spostamento risultano rappresentate da polinomi di grado n+1, e un ovvia *[sic!]* generalizzazione permette anche di concludere che, nell'ipotesi che le componenti del tensore siano sviluppabili in serie di polinomi convergenti integrabili, risulteranno pure sviluppate in serie di polinomii le componenti dello spostamento.

La Sezione Matematica rileva, oltre la correttezza una certa abilità di condotta, e propone che la Nota del Sig.r Guido Facciotti sia ammessa alla lettura e inserita nei Rendiconti.[286]

Milano 22 Aprile 1930

Gian Antonio Maggi - Relatore.

**XVIII**

[su carta intestata R. Università di Milano - Facoltà di Scienze;
l'angolo sinistro in alto, dove probabilmente compariva la data, è strappato]

Relazione della Nota della Dott. Maria Pacifico
"Sopra alcuni problemi al contorno per le funzioni armoniche".

La Sig.na Dott. Maria Pacifico, allieva del prof. M. Picone, riprende i problemi di Neumann nella corona circolare, e misto di Dirichlet e di Neumann nel cerchio, oggetto di classiche soluzioni del Dini, per rimuovere la restrizione, comunque tenue, della sviluppabilità in serie di Fourier delle funzioni date sul contorno, supposta dal Dini, e valersi di un metodo simile a quello usato, con vantaggio, dal prof. Picone pel problema di Dirichlet nella corona circolare.

La Sezione Matematica, avuto anche riguardo all'invariabile significato di quei problemi, interessanti per la Fisica Matematica, non meno che per l'Analisi, reputa la Nota della Dott. M. Pacifico meritevole di essere *[angolo strappato]*, e inserita nei Rendiconti.[287]

---

[286] *Rendiconti Istituto Lombardo*, v. 63, 1930, pp. 873-879.
[287] *ibidem*, pp. 541-555.

## XIX

Bruno Finzi. Nato a Gardone Val Trompia il 12 Febbrajo 1899, ex combattente, come Tenente del Genio, è attualmente Professore Straordinario di Meccanica Razionale in questa R. Università di Milano, designatovi dalla Facoltà di Scienze, in base al suo onorevole risultato nel Concorso, per la stessa cattedra nella R. Università di Genova, dello scorso anno. Già attivo nel nostro Politecnico, vi ottenne nel 1920 il Diploma di Ingegnere, e sul principio ne esercitò la professione. Ma, attratto dalla Scienza pura, tanto che, già nel 1921, conseguiva con lode la Laurea in Matematica Pura nella R. Università di Pavia, abbandonò ben presto la professione, per dedicarsi alle discipline matematiche e all'insegnamento, con tanto successo da giungere poco più che trentenne, alla cattedra universitaria.

Copiosa, varia e ricca di pregi è la produzione scientifica del Finzi, in non piccola parte contenuta in Note dei Rendiconti del nostro Istituto, il resto nei Rendiconti dei Lincei, e in altri tra i maggiori Periodici Scientifici. Un notevole gruppo concerne l'Idrodinamica, sotto i diversi aspetti della Cinematica, dell'Energetica, dell'Integrazione delle equazioni del movimento, della critica di questioni discusse, e un contributo particolarmente originale va riconosciuto nelle sue Note sui veli liquidi. Un altro pur notevole gruppo concerne quel recente ramo dell'Analisi, che ha pure tanta attinenza colla Meccanica e colla Fisica Matematica, chiamato Calcolo Tensoriale. E qui le ampie generalizzazioni del Finzi, coi suoi sistemi multipli, gli valgono un posto tra i presenti indagatori dello spazio geometrico, sotto l'aspetto della sua rappresentazione coi principii del Calcolo Differenziale Assoluto. All'infuori di questi gruppi, si trovano poi ricerche di Meccanica Razionale, tra cui la brillante teoria del Bumerang, contributi alla teoria della Relatività, un pregevole trattato, in collaborazione colla Dr. Pastori, "Lezioni di Matematica per Chimici" che riassume questo suo insegnamento, del quale fu incaricato per tre anni dalla nostra Facoltà di Scienze.

Tutta l'opera del Finzi rivela vivo e acuto ingegno, spiccata personalità, costante attività, ampia e soda coltura. Per cui, colla sua nomina a Socio Corrispondente,[288] *[al]*l'Istituto si aggregherà un giovane scienziato, che onora la Scuola lombarda, in cui si è formato, che immancabilmente gli assicura un attivo e valoroso contributo.

[288] Da questa informazione, si può dedurre che questo scritto sia del 1933, anno nel quale Finzi fu nominato Socio Corrispondente dell'Istituto Lombardo (Adunanza del 19/1).

## XX

*[in penna rossa]* Copia del M.S. consegnato a Chisini pel prossimo numero del Periodico il 19 Nov. 1933.[289]

In quanto tempo un pianeta, fermato, cadrebbe nel Sole?

È ben da aspettarsi che, supponendo l'orbita ellittica, anziché circolare, il calcolo fornisca tempi diversi, per la caduta da punti diversi dell'orbita, invece del mio unico tempo. Ma non parrebbe da aspettarsi, allo stesso modo, l'unico valore del rapporto della durata della rivoluzione del pianeta a quest'unico tempo, nell'ipotesi dell'orbita ellittica, che, per ogni pianeta, rappresenta valor medio. Questo è il risultato, che mi è parso meritevole di qualche attenzione. Il prof. Federici, nella sua Nota dello stesso titolo della mia (*),[290] riportando, come "risultato noto", la pura formola, senza commenti, sembra non tener conto di quel significato, che io considero pure come unica ragione della mia piccola Nota.

G.A. Maggi

(*) Periodico di Matematica, 1 Novembre 1933[291]

## XXI

Osservazioni sulla Nota di G. Vitali "Sulla forza centrifuga"[292]

I due casi del punto che si libera dall'appoggio sopra una piattaforma girevole e del punto che rompe il filo della fionda sono identici, e non possono servire a distinguere le due specie di forza centrifuga. Avevo osservato, a suo tempo, come una svista, questa distinzione degli autori della recente pubblicazione. Questo pel primo capoverso.

Per secondo, il Maggi[293] non ha mai inteso che le due specie di forza centrifuga corrispondano a fatti identici, e in essi varii soltanto il punto di vista dell'osservatore. Il Maggi si è ripetutamente ingegnato a "mettere in rilievo" precisamente il contrario.

La seguente deduzione è esatta, ma alle stesse conclusioni si arriva semplicemente invocando, dalla teoria del movimento relativo, il teorema del parallelogrammo delle velocità o il teorema di Coriolis. Per il primo, è nulla, nella posizione iniziale, la velocità del movimento riferita alla coppia girante colla piattaforma. Pel secondo, l'accelerazione di questo movimento è la somma dell'accelerazione centrifuga e della centrifuga composta. Questa è nulla, colla velocità del movimento relativo, nella posizione iniziale; e ne risulta l'accelerazione del movimento riferito alla coppia girante, o alla piattaforma, in quella posizione, l'accelerazione centrifuga. Moltiplicando l'accelerazione per la massa del punto, se ne ottiene la forza centrifuga.

---

[289] Non c'è di questa comunicazione nel verbale dell'adunanza del 19 novembre.

[290] Nota di Maggi pubblicata nel 1927 (*Periodico di matematiche*, pp. 329-330) e ristampata in *Selecta*, 1932, pp. 91-92.

[291] C. Federici, v. XIII, pp. 289-291.

[292] *Periodico di matematiche*, v. XIV, 1934, pp. 190-193. In tale Nota, Vitali critica una posizione presentata nella pubblicazione: A. Palatini, R. Serini, *Elementi di Fisica*, v. 1 alla p. 127, come si legge nella nota a piè di pagina.

[293] Maggi pubblica due Note sulla forza centrifuga: la prima nel 1926 su *Il Nuovo Cimento* e la seconda nel 1930 sul *Periodico di Matematiche*.

Questa forza centrifuga, a cui l'autore si arresta, è la specie applicata al punto la quale non esiste che per l'osservatore connesso colla piattaforma girante, e, in generale, connesso cogli assi mobili.

Ma esiste l'altra specie, applicata dal punto mobile ad un altro punto, appartenente al corpo con cui è realizzato il vincolo, come reazione rispetto alla pressione vincolare. È la forza che strappa il filo della fionda, che sfonda la guida che obbliga il punto al movimento circolare, che logora l'appoggio del punto alla piattaforma girevole (che gli accennati autori, non si saprebbe perché, esimono da simile logoramento).

Le due specie di forza centrifuga hanno comuni gli elementi dinamici, grandezza, direzione e verso, ma diverso il punto d'applicazione, e sono di natura interamente diversa.

Il caso del punto appoggiato alla piattaforma girevole, e liberato ad un certo momento, dal *[*sistema*]* che ve lo manteneva connesso, va considerato, nei riguardi della forza centrifuga, nell'ipotesi che il movimento del punto liberato non sia soggetto a sensibile attrito da parte della piattaforma. È a questo modo che questo caso s'identifica con quello della fionda.

Infine va osservato che il movimento del punto riferito alla piattaforma, salvo l'accennato caso che si annulli la velocità del movimento medesimo, risulta soggetto oltre che alla forza centrifuga - la sola che l'Autore sembra considerare - anche alla forza centrifuga composta.

N.B. *[in penna rossa]* Copia di relazione consegnata a Chisini su lettura del M.S. presentato pel Period. di Matematica.[294]

## XXII

Analyse Mathématique sur un problème mixte
Comptes Rendus 1930[295]

Sollecitazioni iperastatiche. Ist. Lomb. 1932[296]
"Viene proposta e motivata la qualifica di 'sollecitaz. iperastatiche' per una particolare categoria di sollecitazioni astatiche; coll'aggiunta di qualche osservaz. sul teorema di Da Silva e sulla definiz. degli assi d'equilibrio di Möbius".

Alcune proprietà di media nella Elastostatica ordinaria - Lincei, 1932.[297]

Sopra alcune questioni di Statica dei sistemi continui - Sc. Norm. di Pisa, 1933.[298]
"Uno degli scopi di questa Memoria è la ricerca di tutti i tipi di sollecitazione totale che - almeno nei limiti delle Elastostatica ordinaria - danno luogo, in ogni sistema elastico omogeneo ad uno stress indipendente dalla natura del materiale".

---

[294] Non vi è traccia di questa comunicazione nelle adunanze successive. Si trova invece una risposta degli autori del volume oggetto della critica alle pp. 317-322.
[295] A. Signorini, p. 712.
[296] A. Signorini, v. 65, pp. 1144-1150.
[297] A. Signorini, v. 15, pp. 151-156.
[298] A. Signorini, pp. 231-151.

## XXIII

Il discorso del Bertrand (§§ 189 e segn.) è dominato dalla inconsistente identificazione di quantità e relativa misura. La lunghezza è una quantità inerente al segmento rettilineo, le cui proprietà geometriche sono sufficienti per i criterii d'eguaglianza, e di maggiore e minore, e per la definizione della somma di più lunghezze (somma di quantità); dal cui complesso si ricava la misura di una lunghezza, rispetto ad una convenuta unità di misura (pure lunghezza): questo un numero variabile, in modo assegnabile, coll'unità.

Colla legge della caduta libera di un punto grave, la misura della lunghezza del tratto di verticale descritto dal grave alla fine d'ogni tempo riesce proporzionale al quadrato della misura di questo intervallo di tempo. Si potrà quindi valersi di questi quadrati di misure di tempo per formare la misura delle corrispondenti lunghezze di segmenti di verticale. Ma ciò non conferisce alcun significato alla domanda se queste lunghezze sono per avventura quadrati di tempi. Basta esaminare di che cose realmente si tratta.

E così per tutto il resto.

J. Bertrand, Leçons sur la théorie mathématique de l'électricité (Paris, 1890) pag. 267.

Milano 27 Aprile 1936.

## XXIV

Mi è ben grato attestare, a favore della Sig.na Maria Pastori, perché di questa mia attestazione essa possa pienamente valersi, che la stessa Sig.na Maria Pastori, Dottore in Matematica, e Libera docente di questa Facoltà di Scienze, dall'anno accademico 1925-26 al 1930-31, col titolo di Assistente di Analisi, ma prestando effettivamente opera di Assistente di Meccanica Razionale, impartì, sotto la mia direzione, il regolare corso di Esercizi di Meccanica Razionale, colla più lodevole attività, non solo, ma dimostrando, coll'efficace esposizione e l'illuminata scelta e trattazione degli argomenti, veramente eccellenti qualità didattiche, così da doversi, nei successivi esami dei giovani, riconoscere il frutto e renderle cospicuo merito del suo insegnamento.

Milano, 27 Maggio 1936, XIV

In fede

Prof. Gian Antonio Maggi.

# XXV

## Statistica di Fermi.

| | | |
|---|---|---|
| Valore dell'energia | $E_1$, | $E_2$ |
| Numero delle particelle a cui appartiene | $n_1=2$, | $n_2=3$, |
| Celle disponibili | $g_1=3$, | $g_2=4$, |
| Combinazioni relative ai due valori | $\frac{1\cdot2\cdot3}{1\cdot2\cdot1}=3$, | $\frac{1\cdot2\cdot3\cdot4}{1\cdot2\cdot3\cdot1}=4$. |
| Combinazioni composte (con ciascuna del primo gruppo una del secondo) | $3\times4=12$ | |
| Combinazioni possibili in ambedue i casi | $\left\{\begin{matrix} n=n_1+n_2=5 \\ g=g_1+g_2=7 \end{matrix}\right\}$ | |

$$\frac{1\cdot2\cdot3\cdot4\cdot5\cdot6\cdot7}{1\cdot2\cdot3\cdot4\cdot5\cdot1\cdot2}=\frac{6\cdot7}{1\cdot2}=21.$$

La probabilità della supposta distribuzione risulta

$$\frac{12}{21}=\frac{4}{7}.$$

Milano, 9.VI.1937.

## Indice dei corrispondenti e dei nomi citati nelle lettere e nei *Giudizii*

Nell'elenco che segue si trovano sia i corrispondenti che i nomi citati nei documenti in ordine lessicografico. Di tutti, ove possibile, sono riportate le date di nascita, di morte e l'attività svolta. Quasi tutti i personaggi dell'elenco possono essere trovati facilmente in rete; di alcuni, meno noti, si sono inserite le informazioni in possesso; i rimanenti restano, al momento, ignoti al curatore.

I numeri a destra rimandano alle trascrizioni dei documenti nei quali si trova il nome considerato; in corsivo sono indicati i corrispondenti e i numeri scritti in grassetto corsivo rimandano alle lettere nelle quali il corrispondente è mittente o destinatario.

Accanto al numero della lettera si trova l'eventuale errata grafia con la quale l'autore ha scritto il citato.

## Materiali su Gian Antonio Maggi

L'elenco completo delle pubblicazioni scientifiche di G.A. Maggi si trova in: U. Cisotti, "Gli scritti scientifici di Gian Antonio Maggi", *Rendiconti del Seminario Matematico e Fisico di Milano,* 1938, pp. 183-189.

Di seguito l'elenco delle commemorazioni:

1937- U. Cisotti, "G.A. Maggi", *Rendiconti del Seminario Matematico e Fisico di Milano*, pp. XIX-XX
Qui si trovano la fotografia e la firma di Maggi.
1937- U. Cisotti, "Gian Antonio Maggi", *Rendiconti del R. Istituto Lombardo*, pp. 169-171
1937- "Necrologio", *Bollettino della Unione Matematica Italiana*, p. 168
1937- A.M. Pizzagalli, "Necrologio", *Atene e Roma*, s.III, pp. 214-215
1937- B. Finzi, "Necrologio di G.A. Maggi", *Periodico di matematiche,* pp. 184-188
1937/38- C. Somigliana, "Gian Antonio Maggi", *Atti della R. Accademia delle Scienze di Torino*, pp. 519-528
1938- A. Signorini, "Commemorazione del Socio Gian Antonio Maggi", *Rendiconti dell'Accademia dei Lincei*, pp. 455-466
1938- U. Cisotti, "Gli scritti scientifici di Gian Antonio Maggi", *Rendiconti del Seminario Matematico e Fisico di Milano,* pp. 167-189
Il 10/6/1938 presso la sede del *Seminario Matematico* si tenne una commemorazione di Maggi - il cui verbale si trova alle pp. IX-X dei *Rendiconti* - alla quale parteciparono, oltre a numerosi colleghi, i parenti: Clara Maggi, Bianca Ferratini Maggi, Alberto Ferratini, Gian Emilio Maggi, Gianvincenzo Maggi, Angelo Maria Pizzagalli, Bianca Franchi Maggi Casorati, Maria Luisa Pizzini Maggi, Maria Luisa Maggi Lolli, Miotty Maggi Gherini.
1957- B. Finzi, "Gian Antonio Maggi", *Rendiconti del Seminario Matematico e Fisico di Milano,* pp. X-XIV
2001- Giannantonio 'Popi' Maggi, *Note biografiche di Gian Antonio Maggi*, Comerio
Si tratta di una raccolta scritta dal bisnipote e conservata presso l'archivio privato Andrea Ferratini; contiene la trascrizione di gran parte delle commemorazioni del 1937/38 qui indicate, della Memoria "In quanto tempo un pianeta, fermato, cadrebbe nel sole?", di alcuni discorsi tenuti da Maggi e la presentazione - scritta con taglio personale - della figura di Maggi, arricchita da aneddoti famigliari.
2006- S. D'Agostino, F. Della Torre, S. Leva, A. P. Morando, A. Rossi, "Gian Antonio Maggi", XXVI Congresso Nazionale di Storia della Fisica e dell'Astronomia, Roma, 15-17 giugno 2006; in attesa di pubblicazione sugli Atti del Convegno

## APPENDICE
## Fonti correlate e trascrizioni di due lettere a Ernesto Cesàro e a Max Abraham

In altri fondi archivistici accademici si trovano documenti relativi a Maggi.
Ad oggi ho reperito:

- 1 lettera del 1898 indirizzata a Alessandro Paoli, conservata presso la Biblioteca universitaria dell'Università degli Studi di Pisa
- un fascicolo conservato presso l'Accademia Nazionale dei Lincei, Fondo "Tullio Levi Civita" (cfr. P. Cagiano de Azevedo, "Fondo Tullio Levi Civita-Inventario della corrispondenza"):
  Busta 11, fascicolo 517, Maggi Gian Antonio
  1900 set. 9: biglietto con busta;
  1901 ago. 1: lettera con busta;
  1902 set. 3: lettera con busta;
  1903 giu. 4, nov. 1: 2 lettere con buste;
  1904 mar. 8, ago. 1: 2 lettere con buste;
  1905 mag. 1: lettera con busta;
  1906 set. 18, dic. 31: lettera e cartolina illustrata;
  1908 mag. 29, giu. 6, lug. 26 e 30: 4 lettere con 3 buste;
  1910 giu. 4, dic. 29; lettera e biglietto con buste;
  1912 mar. 17, nov. 28; 2 lettere con buste;
  1914 feb. 5, dic. 20: 2 lettere con buste;
  1915 mar. 31; lettera con busta;
  1917 lug. 18 cartolina postale;
  1918 dic. 31: lettera con busta;
  1921 giu. 19, dic. 19: biglietto con busta, cartolina illustrata;
  1922 apr. 17, giu. 24, set. 19: 2 lettere e biglietto con buste;
  1923 mar. 2, ago. 21, set. 3: 2 lettere con buste, cartolina illustrata;
  1925 gen. 2, apr. 29, giu. 20, nov. 22: lettera e biglietto con buste, 2 cartoline postali;
  1926 lug. 2, 9, ott. 2: 2 lettere e 1 biglietto con buste;
  1927 gen. 28, ago. 20, nov. 21: 2 lettere con buste, cartolina postale;
  1928 apr. 5, ago.25, dic. 6: cartolina postale, 2 cartoline illustrate;
  1929 mag. 2, ott. 22, nov. 9, 14: 3 lettere con buste, cartolina postale;
  1930 apr. 21, set. 3: 2 lettere con buste;
  1932 feb. 8: cartolina illustrata;
  1933 set. 24, ott. 7, 11, 23, 27: 4 lettere e biglietto con 4 buste;
  1934 gen. 4, apr. 28, giu. 15, ago. 3, set 28: 3 lettere e 1 biglietto con buste, cartolina illustrata;
  1935 gen. 1, 4, [apr. 23], ago. 6: 2 lettere con buste, 2 cartoline illustrate;
  1936 giu. 12, lug. 15, ago. 2: 2 lettere e 1 biglietto con buste;
  1937 apr. 10, giu. 2: lettera con busta, cartolina illustrata.
- 1 lettera del 12/2/1895 indirizzata a Ernesto Cesàro conservata presso il Dipartimento di Matematica e Applicazioni "R. Cacciopoli" dell'Università degli Studi "Federico II" di Napoli (disponibile in rete e trascritta qui di seguito).
- 1 lettera del 20/3/1913 indirizzata a Giuseppe Colombo conservata presso l' Archivio del Politecnico di Milano, fasc. personale di Max Abraham (la trascrizione che compare qui di seguito è disponibile in rete). È probabile che presso questo Archivio siano presenti ulteriori documenti relativi all'attività di Maggi per il Seminario matematico.

- Materiale diverso presso l'archivio privato del nipote Andrea Ferratini:
  - quaderni con accurati disegni di uccelli realizzati in età scolare;
  - quaderni di appunti di lezioni di scuola superiore e universitaria;
  - diplomi di laurea con relativi libretti;
  - diplomi di assegnazione dei premi e delle onorificenze;
  - comunicazioni di elezione a membro di diverse Società scientifiche;
  - documenti di assegnazione delle diverse cattedre e di dimissione dalle università;
  - raccolta di note biografiche (95 pagine) a cura del bisnipote Giannantonio 'Popi' Maggi;
  - raccolta di poesie dedicate a Maggi nel suo giorno onomastico 1889 da parte di un poeta messinese non identificato;
  - minute di lettere: 3 al Direttore del Corriere della Sera (sulla corretta grafia dei nomi slavi); stralcio alla mamma (Pisa 29/10/1908); 1 in cirillico (destinatario sconosciuto, 14/1/1909);
  - 1 lettera di A.V. Vassilief (1898) che annuncia l'elezione di Maggi a membro della Société Physico-Mathématique de l'Université Impériale de Kasan;
  - 1 cartolina postale (mittente sconosciuto, Palermo 20/4/1889) di carattere personale;
  - 1 lettera (mittente sconosciuto, Genova 22/11/1889) su una questione di pensionamento.

Lettera di Gian Antonio Maggi a Ernesto Cesàro[299]

Messina 12 Febbrajo 95

Egregio Collega

Ella ha forse giustificato il mio lungo silenzio col supporre che avendo ricevuto il Suo libro[300] colla scritta "omaggio dell'Editore" attribuissi il bel dono ai signori Bocca. Ma in tal caso è stato troppo benevolo; poiché, conoscendo la Sua bontà, ne sono subito restato obbligato a Lei, anche prima di sapere come stava veramente la cosa, come poi seppi dal prof. Brambilla. Quando parlai con lui nella mia fermata di poche ore a Napoli, al ritorno da Milano, fermata troppo breve per aver il tempo di venire a vederLa, aveva appena ricevuto il libro: e mi proponevo di scriverLe e ringraziarneLa, appena arrivato. Ma all'indomani ci sorprese il trambusto dei tremuoti;[301] e questa sia la mia prima scusa d'aver tanto tardato e mancato. Intanto ho avuto occasione di vedere diverse parti del Suo interessantissimo lavoro, e di apprezzare la misura, la precisione e la chiarezza a cui è informato. Né forse Ella pensa tutto il piacere che mi hanno fatto i cenni relativi al Clifford, la cui traduzione anonima è un lavoro mio, che per l'indole del libro, non fu dei più lievi.[302]

[299] Questa lettera è consultabile in rete a partire dall'indirizzo http://pzmath.unibas.it/PRIN/STMAT/
[300] Potrebbe trattarsi di: *Teoria dell'elasticità*, Torino, Bocca, 1894.
[301] Si riferisce al sisma che ebbe luogo il 16 novembre 1894 e alle scosse dei mesi successivi.
[302] *Il senso comune nelle scienze esatte. Traduzione, con note, dell'opera: W.K. Clifford, The common Sense of the Exact Sciences*, Dumolard, Milano, 1886.

Forse Le scriverò ancora fra non molto per chiedere il Suo consiglio intorno una questione che concerne appunto il movimento elastico. E questa fu anche una ragione per cui ho aspettato. Ma poiché al presente mi sto occupando d'altro, non starò oltre ad attendere il momento, e a rimandare il pagamento d'un debito già troppo antico.

Coi miei vivi ringraziamenti Ella aggradisca anche i migliori saluti, e mi creda sempre

Dev.mo Suo
Gian Antonio Maggi

P.S. Martinetti mi incarica di mandarLe i suoi saluti, e tutt'e due La preghiamo di ricordarci ai colleghi di costì.

Lettera di Gian Antonio Maggi a Giuseppe Colombo[303]

Roma 20 marzo 1913

Gentilissimo Senatore,

Ho il piacere d'informarLa che è stato testé presentata al Ministero la Relazione della Commissione, che conclude all'unanime proposta della promozione del Prof. Abraham.

Portate in campo le voci della imperfetta conoscenza della lingua italiana e le relative conseguenze, nei rapporti colla Scuola, io ho dovuto contarmi tra quelli a cui tali voci erano arrivate, perché è la verità. Del resto, esattezza a parte, esse sono oramai ampiamente diffuse, e sarebbe stato opportuno che un'esplicita dichiarazione del Consiglio dei Professori, a completamento del voto favorevole alla domanda di promozione, ne prevenisse gli effetti. La questione, già sollevata nell'adunanza di Domenica, pareva eliminata. Diversamente, Le avrei chiesto un breve colloquio, dopo il Teatro. Cercai poi di Lei, il giorno dopo, all'albergo Milano, ma mi fu confermato c'era *[sic!]* partito, quella mattina.

La prego intanto di aggradire i miei più distinti e cordiali rispetti, e di credermi sempre

Devmo Suo
Gian Antonio Maggi

[303] La presente trascrizione è contenuta in C. Citrini, "Matematica e vita civile nel Politecnico di cento anni fa: la vicenda di Max Abraham", *Annali di Storia delle Università italiane*, volume 12, 2008 che è consultabile in rete a partire dall'indirizzo www.cisui.unibo.it

## Ringraziamenti

È stato nel 1992, con l'inizio della tesi, che ho iniziato ad apprezzare il lavoro sugli archivi di corrispondenze di matematici e qui voglio ringraziare, pur a distanza di anni, la mia relatrice, prof.ssa Simonetta Di Sieno, per avermi introdotto in questo mondo e per aver fatto in modo che io potessi rimanerci fino ad oggi. Ringrazio anche il signor Giuliano Moreschi, Direttore della Biblioteca "G. Ricci" del Dipartimento di Matematica dell'Università di Milano - al quale si deve il ritrovamento del Fondo Maggi - che mi ha affidato le lettere sia durante la stesura della tesi che in questi mesi.

Per quanto riguarda le ricerche da me effettuate per la presente pubblicazione, devo ringraziare molte persone.

Grazie al prof. Massimo Galuzzi per le pazienti spiegazioni che mi ha dato sui passi di Platone (rimproverandomi di non aver studiato il greco!) e per la sua costante attenzione nei confronti delle mie attività sul fronte della Storia della Matematica. Devo a lui il contatto con il prof. Ferruccio Franco Repellini, che ringrazio sia per avermi fatto dono della sua competenza e disponibilità nel correggere le parti in greco e in latino delle corrispondenze Paolo Ubaldi e Francesco Zambaldi che per avermi fornito utili notizie.
La dott.ssa Silvia Donghi, Bibliotecaria della Famiglia Meneghina di Milano, con la sua competenza e disponibilità, mi ha permesso di trovare utili informazioni negli archivi di quella Società in merito alle corrispondenze Luigi Medici e Direttore de "la Meneghina" nell'ambiente davvero particolare di via San Paolo.

Ringrazio inoltre: il dott. Jacopo De Tullio che ha rivisto la lettera #42; il prof. Francesco Della Torre per l'invio immediato del testo del suo intervento del 2006 su Maggi; il personale degli Archivi Storici del Campus Bovisa Durando del Politecnico di Milano e la signora Lucia Baroni della Biblioteca Universitaria di Pisa per aver cortesemente effettuato alcune ricerche nei rispettivi archivi; il signor Pasquale Abiusi della Direzione Servizi Civici, Partecipazione e Sport per aver ricercato le sepolture di Clara e Bianca Maggi presso il Cimitero Monumentale; la dott.ssa Tina Paternoster per alcune traduzioni dal greco.

Ho incontrato, durante la stesura di questo lavoro - grazie all'interessamento della prof.ssa Marina Bertolini - dapprima una pronipote di Maggi, la prof.ssa Giovanna Pugno Vanoni e in seguito i nipoti Cecilia e Andrea Ferratini che, con la moglie Lina, hanno dimostrato interesse nei confronti di questo lavoro, mi hanno fornito molte informazioni sulla famiglia e mi hanno dato accesso ai documenti del nonno. Mi è caro qui ringraziarli per la calorosa accoglienza.

Un grazie particolare alla dott.ssa Giovanna Dimitolo, amica e collega, che mi ha aiutato a diversi livelli anche in questo lavoro.